\documentclass[openany]{memo-l}

\usepackage{amsmath} 
\usepackage{amssymb}

\newcommand{\forces}{\Vdash}

\newcommand{\lbv}{[\![} 
\newcommand{\rbv}{]\!]}
\newcommand{\V}{{\bf V}} 
\newcommand{\lesdot}{\mathrel{\mathord{<}\!\!\raise 
0.8 pt\hbox{$\scriptstyle\circ$}}} 
\newcommand{\comp}{\circ} 
  
\newcommand{\add}{{\rm {\bf add}}\/} 
\newcommand{\unif}{{\rm {\bf non}}\/} 
\newcommand{\cof}{{\rm {\bf cof}}\/} 
\newcommand{\non}{{\rm {\bf non}}\/} 
\newcommand{\cov}{{\rm {\bf cov}}\/}
\newcommand{\con}{{\mathfrak c}} 
\newcommand{\dominating}{{\mathfrak d\/}} 
\newcommand{\unbounded}{{\mathfrak b}}

\newcommand{\can}{2^{\textstyle \omega}} 
\newcommand{\fs}{2^{\textstyle <\!\omega}} 
\newcommand{\baire}{\omega^{\textstyle \omega}} 
\newcommand{\iso}{[\omega]^{\textstyle \omega}} 
\newcommand{\fsuo}{[\omega]^{\textstyle <\!\omega}} 
\newcommand{\fseo}{\omega^{\textstyle <\!\omega}} 
\newcommand{\pfs}{\omega^{\mathunderaccent\smile-3 \omega}}
\newcommand{\conc}{{}^\frown\!}
\newcommand{\lh}{\ell g\/} 
\newcommand{\rest}{{\restriction}}
\newcommand{\mrot}{{\rm root}\/} 
\newcommand{\suc}{{\rm succ}} 
\newcommand{\spliting}{{\rm split}}
\newcommand{\dom}{{\rm dom}} 
\newcommand{\rng}{{\rm rng}}
\newcommand{\nor}{{\rm {\bf nor}}\/} 
\newcommand{\NOR}{{\rm NOR}} 

\newcount\skewfactor
\def\mathunderaccent#1#2 {\let\theaccent#1\skewfactor#2
\mathpalette\putaccentunder}
\def\putaccentunder#1#2{\oalign{$#1#2$\crcr\hidewidth
\vbox to.2ex{\hbox{$#1\skew\skewfactor\theaccent{}$}\vss}\hidewidth}}

\newcommand{\random}{{\mathbb B}}
\newcommand{\C}{{\mathcal C}}
\newcommand{\D}{{\mathcal D}}
\newcommand{\F}{{\mathcal F}} 
\newcommand{\cH}{{\mathcal H}}  
\newcommand{\J}{{\mathcal J}}  
\newcommand{\K}{{\mathcal K}}

\newcommand{\M}{{\mathcal M}}  
\newcommand{\N}{{\mathcal N}}  
\newcommand{\cP}{{\mathcal P}}
\newcommand{\p}{{\mathbb P}}

\newcommand{\q}{{\mathbb Q}}  

\newcommand{\cS}{{\mathcal S}}

\newcommand{\mbR}{{\mathbb R}}

\newcommand{\apr}{{\rm apr}}
\newcommand{\dn}{{\rm dn}}
\newcommand{\up}{{\rm up}}
\newcommand{\val}{{\bf val}}
\newcommand{\dis}{{\bf dis}}
\newcommand{\pos}{{\rm pos}}

\newcommand{\CR}{{\rm CR}}
\newcommand{\WCR}{{\rm WCR}}
\newcommand{\bH}{{\bf H}}
\newcommand{\uhalf}{{\rm half}}
\newcommand{\basis}{{\bf basis}}
\newcommand{\Wtil}{{\dot{W}}}
\newcommand{\FC}{{\rm FC}}
\newcommand{\PFC}{{\rm PFC}}
\newcommand{\PC}{{\rm PC}}

\newcommand{\dcl}{{\rm dcl}}
\newcommand{\Fin}{{\rm Fin}}
\newcommand{\TCR}{{\rm TCR}}

\newcommand{\bsum}{{\rm sum}}
\newcommand{\tsum}{{\rm tsum}}

\newcommand{\tree}{{\rm tree}}
\newcommand{\fH}{\varphi_\bH}
\newcommand{\NMP}{{\bf NMP}}
\newcommand{\NNP}{{\bf NNP}}
\newcommand{\UP}{{\bf UP}}
\newcommand{\sUP}{{\bf sUP}}
\newcommand{\AB}{{\bf AB}}
\newcommand{\dpt}{{\rm dp}}

\newtheorem{theorem}{Theorem}[section] 
\newtheorem{claim}{Claim}[theorem]
\newtheorem{lemma}[theorem]{Lemma} 
\newtheorem{proposition}[theorem]{Proposition} 
\newtheorem{corollary}[theorem]{Corollary} 
\newtheorem{SSthesis}[theorem]{Thesis}

\theoremstyle{definition}

\newtheorem{definition}[theorem]{Definition}
\newtheorem{gendef}[theorem]{General Definitions}
\newtheorem{example}[theorem]{Example}

\theoremstyle{remark}
\newtheorem{notation}[theorem]{Notation}
\newtheorem{conclusion}[theorem]{Conclusion}
\newtheorem{remark}[theorem]{Remark}

\numberwithin{section}{chapter}
\numberwithin{equation}{chapter}

\begin{document}
\frontmatter
\title{Norms on possibilities I: forcing with trees and creatures}

\author{Andrzej Ros{\l}anowski}
\address{Institute of Mathematics\\
 The Hebrew University of Jerusalem\\
 91904 Jerusalem, Israel\\
 and Mathematical Institute of Wroclaw University\\
 50384 Wroclaw, Poland} 
\curraddr{Department of Mathematics and Computer Science\\
 Boise State University\\
 Boise ID 83725, USA}
\email{roslanow@math.idbsu.edu}
\urladdr{http://math.idbsu.edu/$\sim$roslanow}
\thanks{The first author thanks the Hebrew University of Jerusalem and the
Lady Davis Foundation for the Golda Meir Postdoctoral Fellowship, and the
KBN (Polish Committee of Scientific Research) for partial support through
grant 2P03A01109.} 
\author{Saharon Shelah}
\address{Institute of Mathematics\\
 The Hebrew University of Jerusalem\\
 91904 Jerusalem, Israel\\
 and  Department of Mathematics\\
 Rutgers University\\
 New Brunswick, NJ 08854, USA}
\email{shelah@math.huji.ac.il}
\urladdr{http://www.math.rutgers.edu/$\sim$shelah}
\thanks{The research of the second author was partially supported by ``Basic
Research Foundation'' of the Israel Academy of Sciences and Humanities.
Publication 470} 

\date{December 01, 1997}


\maketitle

\tableofcontents

\setcounter{page}{5}

\chapter*{Annotated Content}

\noindent {\sc Chapter 0:\quad Introduction}
\medskip

\noindent {\sc Chapter 1:\quad Basic definitions}\quad
We introduce a general method of building forcing notions with use of {\em
norms on possibilities} and we specify the two cases we are interested in.
\begin{enumerate}
\item[1.0] {\bf Prologue}
\item[1.1] {\bf Weak creatures and related forcing notions}
[We define weak creatures, weak creating pairs and forcing notions determined
by them.]
\item[1.2] {\bf Creatures}
[We introduce the first specific case of the general schema: creating pairs
and forcing notions of the type $\q^*_{\C(\nor)}(K,\Sigma)$.]
\item [1.3] {\bf Tree creatures and tree--like forcing notions}
[The second case of the general method: forcing notions $\q^{\tree}_e(K,
\Sigma)$ in which conditions are trees with norms; tree creatures and
tree--creating pairs.] 
\item[1.4] {\bf Non proper examples}
[We show several examples justifying our work in the next section: the method
may result in forcing notions collapsing $\aleph_1$, so special care is
needed to ensure properness.]
\end{enumerate}

\noindent {\sc Chapter 2:\quad Properness and the reading of names}\quad
We define properties of weak creating pairs which guarantee that the forcing
notions determined by them are proper. Typically we get a stronger property
than properness: names for ordinals can be read continuously. 
\begin{enumerate}
\item[2.1] {\bf Forcing notions ${\bf Q}^*_{{\rm s}\infty}(K,\Sigma)$, ${\bf
Q}^*_{{\rm w}\infty}(K,\Sigma)$}  
[We show that the respective forcing notions are proper if $(K,\Sigma)$ is
finitary and either growing or captures singletons.]
\item[2.2] {\bf Forcing notion ${\bf Q}^*_{f}(K,\Sigma)$: bigness and halving}
[We introduce an important property of creatures: bigness. We note that it is
useful for deciding ``bounded'' names without changing the finite part of a
condition in forcing notions discussed in 2.1. Next we get properness of
$\q^*_f(K,\Sigma)$ when the creating pair $(K,\Sigma)$ is big and has the
Halving Property.] 
\item[2.3] {\bf Tree--creating $(K,\Sigma)$}
[We show that properness is natural for forcing notions
$\q^{\tree}_e(K,\Sigma)$ determined by tree--creating pairs (with our norm
conditions). With more assumptions on $(K,\Sigma)$ we can decide names on
fronts.] 
\item[2.4] {\bf Examples}
[We recall some old examples of forcing notions with norms putting them in our
setting and we build more of them.]
\end{enumerate}

\noindent {\sc Chapter 3:\quad More properties}\quad
We formulate conditions on weak creating pairs which imply that the
corresponding forcing notions: do not add unbounded reals, preserve non-null
sets or preserve non-meager sets. 
\begin{enumerate}
\item[3.1] {\bf Old reals are dominating}
[From the results of section 2 we conclude that various forcing notions are
$\baire$--bounding.]
\item[3.2] {\bf Preserving non-meager sets}
[We deal with preservation of being a non-meager set. We show that if a
tree--creating pair $(K,\Sigma)$ is $T$-omittory then the forcing notion
$\q^{\tree}_1(K,\Sigma)$ preserves non-meager sets. We formulate a weaker
property (being of the $\NMP$--type) which in the finitary case implies that
forcing notions $\q^{\tree}_1(K,\Sigma)$, $\q^*_f(K,\Sigma)$ preserve
non-meager sets. We get a similar conclusion for $\q^*_{{\rm
w}\infty}(K,\Sigma)$ when $(K,\Sigma)$ is a finitary creating pair which
captures singletons.]
\item[3.3] {\bf Preserving non-null sets}
[We formulate a property of tree--creating pairs which implies that the
forcing notion $\q^{\tree}_e(K,\Sigma)$ preserves non-null sets.]
\item[3.4] {\bf (No) Sacks Property}
[An easy condition ensuring ``no Sacks property'' for forcing notions of our
type.] 
\item[3.5] {\bf Examples}
[We build a tree--creating pair $(K,\Sigma)$ such that the forcing notion
$\q^{\tree}_1(K,\Sigma)$ is proper, $\baire$-bounding, preserves the outer
measure, preserves non-meager sets but does not have the Sacks property.]
\end{enumerate}

\noindent {\sc Chapter 4:\quad Omittory with Halving}\quad
We explain how omittory creating pairs with the weak Halving Property produce
almost $\baire$-bounding forcing notions. 
\begin{enumerate}
\item[4.1] {\bf What omittory may easily do}
[We show why natural examples of forcing notions $\q^*_{{\rm s}\infty }(K,
\Sigma)$ (for an omittory creating pair $(K,\Sigma)$) add a Cohen real and
make ground model reals meager.]
\item[4.2] {\bf More operations on weak creatures}
[Just what the title says: we present more ways to put weak creatures
together.] 
\item [4.3] {\bf Old reals are unbounded}
[We say when a creating pair $(K,\Sigma)$ is of the $\AB$--type and we show
that this property may be concluded from easier--to--check properties. We show
that $\q^*_{{\rm s}\infty}(K,\Sigma)$ is almost $\baire$-bounding if
$(K,\Sigma)$ is growing condensed and of the $\AB$--type. For omittory
creating pairs we do not have to assume ``condensed'' but then we require a
stronger variant of the $\AB$, $\AB^+$.]
\item[4.4] {\bf Examples}
[We generalize the forcing notions from \cite{Sh:207}, \cite{RoSh:501} building
examples for properties investigated before.]
\end{enumerate}

\noindent {\sc Chapter 5:\quad Around not adding Cohen reals}\quad
We try to ensure that the forcing notions built according to our schema do not
add Cohen reals even if iterated. We generalize ``$(f,g)$--bounding'' and
further we arrive to a more general iterable condition implying ``no Cohen
reals''. 
\begin{enumerate}
\item[5.1] {\bf $(f,g)$--bounding}
[We present easy ways to make sure that our method results in $(f,g)$--bounding
forcing notions.]
\item[5.2] {\bf $(\bar{t},\bar{\F})$--bounding} [We introduce a
natural (in our context) generalization of $(f,g)$--bounding property. For the
sake of completeness we show that the new property is preserved in CS
iterations.] 
\item[5.3] {\bf Quasi--generic $\Gamma$ and preserving them}
[We formulate a reasonably weak but still iterable condition for not adding
Cohen reals. We define $\bar{t}$--systems, we say when $\Gamma$ is
quasi-$W$-generic and when a forcing notion is $\Gamma$--genericity
preserving. These notions will be crucial in the next section too.]
\item[5.4] {\bf Examples}
[We construct a sequence $\langle W^k_{\ell,h}: h\in \F_\ell, k,\ell\rangle$
such that $W^k_{\ell,h}$ are $\bar{t}$--systems and various forcing notions
(including the random algebra) are $(\Gamma,W^k_{\ell,h})$--genericity
preserving (for quasi-generic $\Gamma$). We build a forcing notion $\q^*_{{\rm
w}\infty}(K,\Sigma)$ which is proper, $\baire$-bounding, $(f,g)$-bounding,
makes ground model reals null and we use the technology of
``$\Gamma$--genericity'' to conclude that its CS iterations with Miller's
forcing, Laver's forcing and random algebra do not add Cohen reals.] 
\end{enumerate}

\noindent {\sc Chapter 6:\quad Playing with ultrafilters}\quad Our aim here is
to build a model in which there is a $p$-point generated by $\aleph_1$
elements which is not a $q$-point and ${\mathfrak m}_1=\aleph_2$.
\begin{enumerate}
\item[6.1] {\bf Generating an ultrafilter} 
[We say when and how quasi-$W$-generic $\Gamma$ determines an
ultrafilter on $\omega$.] 
\item[6.2] {\bf Between Ramsey and $p$-points}
[We define semi--Ramsey and almost--Ramsey ultrafilters and we have a short
look at them.]
\item[6.3] {\bf Preserving ultrafilters}
[We give conditions on a tree--creating pair $(K,\Sigma)$ which imply that the
forcing notion $\q^{\tree}_1(K,\Sigma)$ preserves ``$\D$ is an ultrafilter''
for $\D$ which is Ramsey, almost Ramsey. We say when the filter generated by
$\D$ in the extension is almost Ramsey.]
\item[6.4] {\bf Examples}
[We construct $\bar{t}$--systems $W^n_L$ such that if a
quasi-$W^N_L$-generic $\Gamma$ generates a semi--Ramsey ultrafilter then it
generates an almost--Ramsey ultrafilter, and we build a suitable quasi-generic
$\Gamma$. For a function $\psi\in\baire$ we give a tree--creating pair
$(K,\Sigma)$ such that the forcing notion $\q^{\tree}_1(K,\Sigma)$ preserves
``$\D$ is an almost--Ramsey ultrafilter'' and it adds a function $\dot{W}$
(with $\dot{W}(m)\in[\psi(m)]^{\textstyle m+1}$) such that for each partial
function $h\in\prod\limits_{m\in\dom(h)}\psi(m)$ infinitely often
$h(m)\in\dot{W}(m)$. Next we apply it to get an answer to Matet's problem.]
\end{enumerate} 

\noindent {\sc Chapter 7:\quad Friends and relatives of PP}\quad We deal with
Balcerzak--Plewik number and various properties resembling PP--property.
\begin{enumerate}
\item[7.1] {\bf Balcerzak--Plewik number}
[We recall the definition of $\kappa_{\rm BP}$  and we show that it is bounded
by the dominating number of the relation determined by the strong
PP--property.] 
\item[7.2] {\bf An iterable friend of strong PP--property}
[We introduce a property slightly stronger than the strong PP--property but
which can be easily handled in CS iterations. We show that this property
is natural for forcings $\q^{\tree}_0(K,\Sigma)$, $\q^*_{{\rm w}\infty}(K,
\Sigma)$ (in finitary cases).]
\item[7.3] {\bf Bounded relatives of PP}
[We define various PP-like properties for localizing functions below a given
function. We say how one gets them for our forcing notions and how we may
handle them in iterations.]
\item[7.4] {\bf Weakly non-reducible $p$-filters in iterations}
[We show that a property of filters, crucial for getting PP-like properties
for our forcing notions, is easy to preserve in CS iterations.]
\item[7.5] {\bf Examples}
[For a perfect set $P\subseteq\can$ we build a creating pair $(K,\Sigma)$ such
that the forcing notion $\q^*_f(K,\Sigma)$ is proper, $\baire$-bounding and
adds a perfect subset $Q$ of $P$ whith property that $(\forall K\in\iso)(Q
\rest K\neq 2^{\textstyle K})$. We use this forcing notion to get consistency
of $\dominating<\kappa_{\rm BP}$. We show how forcing notions from other parts
of the paper may be used to distinguish PP-like properties (and the
corresponding cardinal invariants). We build an example of a forcing notion
which is $\baire$--bounding and preserves non-meager sets but which does not
have the strong PP-property.]
\end{enumerate}

\mainmatter
\setcounter{chapter}{-1}
\chapter{Introduction} 
Set Theory\footnote{A large part of the beginning of this introduction is
based on notes of Haim Judah. I really think that they fit to the present
paper, though Saharon Shelah is not convinced --- Andrzej Ros{\l}anowski.}
began with Georg Cantor's work when he was studying some special sets of reals
in connection with the theory of trigonometric series. This study led Cantor
to the following fundamental question: does there exist a bijection between
the natural numbers and the set of real numbers? He answered this question
negatively by showing that there is no such function. Cantor's work did not
stop here and with his sharp intuition he discovered new concepts like the
aleph's scale:
\[0,1,\ldots,\aleph_0,\aleph_1,\ldots,\aleph_\omega,\aleph_{\omega+1},\ldots.\]
Thus Cantor's theorem says that $\aleph_0<2^{\aleph_0}$ and Cantor's question 
was:\ \ \ is $2^{\aleph_0}$ equal to $\aleph_1$?

A real advance on Cantor's question was given by Kurt G\"odel when he showed
that it is (relatively) consistent that $2^{\aleph_0}=\aleph_1$. In 1963 Paul
Cohen showed that if the ZF--axioms for Set Theory are consistent then there
is a model for Set Theory where the continuum is bigger than $\aleph_1$.
Cohen's work is the end of classical set theory and is beginning of a new era. 

When the cardinality of the continuum is $\aleph_1$ (i.e.~CH holds) most of
the combinatorial problems are solved. When the continuum is at least
$\aleph_3$ then most of the known technology fails and we meet very strong
limitations and barriers. 

When the continuum is $\aleph_2$ there are many independence results; moreover
there are reasonably well developed techniques for getting them by sewing the
countable support iterations of proper forcing notions together with theorems
on preservation of various properties.

The aim of this paper is to present some tools applicable in the last case.
We present here a technique of constructing of (proper) forcing notions that
was introduced by Shelah for solving problems related to cardinal invariants
like the unbounded number or the splitting number as well as questions of
existence of special kinds of P-points (see Blass Shelah \cite{BsSh:242} and
Shelah \cite{Sh:207}, \cite{Sh:326}). That method was successfully applied in
Fremlin Shelah \cite{FrSh:406}, Ros{\l}anowski Shelah \cite{RoSh:501},
Ciesielski Shelah \cite{CiSh:653} and other papers. The first attempt to
present a systematic study of the technique was done in the late eighties when
the second author started work on preparation of a new edition of \cite{Sh:b}.
For a long time the new book, \cite{Sh:f}, was supposed to contain 19
chapters. The last chapter, {\em Norms on possibilities}, contained a series
of general definitions and statements of some basic results. However, there
was no new application (or: a good question to solve) and the author of the
book decided to put this chapter aside. Several years later, when the first
author started his cooperation with Shelah some new applications of {\em Norms
  on possibilities} appeared. But the real shape was given to the work due to
questions of Tomek Bartoszy\'nski and Pierre Matet. The answers were very
stimulating for the development of the general method.

This paper is meant as the first one in a series of works presenting
applicability of the method of {\em norms on possibilities}. In
\cite{RoSh:670} we will present more applications of this technique -- for
example we develop the ideas of Ciesielski Shelah \cite{CiSh:653} to
build models without magic sets and their relatives.  Though one can get an
impression here that our method results in non-ccc forcing notions, we managed
to generalize it slightly and get a tool for constructing ccc forcing notions.
That was successfully applied in \cite{RoSh:628} to answer a problem of Kunen
by constructing a ccc Borel ideal on $\can$ which is translation--invariant
index--invariant and is distinct from the null ideal, the meager ideal and
their intersection. It should be pointed out here, that already in Judah
Ros{\l}anowski Shelah \cite{JRSh:373} an example of a ccc forcing notion built
with the use of {\em norms on possibilities} was given (the forcing notion
there can be presented as some $\q^{\tree}_1(K,\Sigma)$ in the terminology
here). Investigations of ccc forcing notions constructed according to our
schema are continued in \cite{RoSh:672}. There are serious hopes that the
technique presented here might be used to deal with problems of large
continuum due to special products. This would continue Goldstern Shelah
\cite{GoSh:448}. Another direction is study of $\sigma$-ideals related to
forcing notions built according to the schema.

Let us note that most of the forcing notions constructed here fall into the
category of snep--forcing notions of \cite{Sh:630}. Consequently, the general
machinery of definable forcing notions is applicable here. We may use it to
improve some of our results, and also to get more tools for handling
iterations (see \cite{Sh:630} and \cite{Sh:669} for more details).  

We want to emphasize that though the aim of the paper is strongly related
to independence proofs it should have some value for those firmly committed to
unembellished ZFC, too. This is nicely expressed by the following: 
\medskip

\begin{SSthesis}
We cannot discover the (candidates for) Theorems of ZFC without having
good forcing techniques to show they are hard nuts.
\end{SSthesis}

\section{The content of the paper} 
Most of the results of the paper originated in answering a particular question
by constructing an example of a forcing notion. However, the general idea of
the paper is to extract those properties of the example which are responsible
for the fact that it works, with the hope that it may help in further
applications of the method. That led us to separation of ``the general
theory'' from its applications, and caused us to introduce a large number of
definitions specifying various properties of weak creating pairs. Each chapter
ends with a section presenting examples and applications of the tools
developed in previous sections. Moreover, at the end of the paper we present
the list of all definitions which appeared in it. This is not a real index,
but should be helpful. The first chapter introduces basic definitions and the
general scheme. In the next chapter we deal with the fundamental question of
when our forcing notions are proper. The first two chapters are a basis for the
rest of the paper. After reading them one can jump to any of the following
chapters. 

In the third chapter we show how we may control some basic properties of
forcing notions built according to our scheme. The properties we deal with
here are related to measure and category and they lead us to example
\ref{tomek1} of a proper forcing notion which is $\baire$--bounding, preserves
non-meager sets and the outer measure, but does not have the Sacks
Property. This answers Bartoszy\'nski's request \cite[Problem 5]{Ba94}.

The fourth part continues \cite{RoSh:501}, dealing with localizations of
subsets of $\omega$. Though the example constructed here is a minor
modification of the one built there, it is presented according to our general
setting. We show explicitly how the weak Halving Property works in this type
of examples. 

A serious problem in getting models of ZFC with given properties of measure
and category is that of not adding Cohen reals. What is disturbing here is
that we do not have any good (meaning: sufficiently weak but iterable)
conditions for this. In the fifth chapter we show how one can ensure that our
general scheme results in forcing notions not adding Cohen reals. A new
iterable condition for this appears here and quite general tools are developed
(see \ref{Gampres}, \ref{superconc}). Finally, in \ref{tomek2},
\ref{tomek2conc}, we fully answer another request of Bartoszy\'nski formulated
in \cite[Problem 4]{Ba94}. We build a proper $\baire$--bounding forcing notion
which preserves non-meager sets, makes ground model reals null, is
$(f,g)$--bounding and such that countable support iterations of this forcing
with Laver forcing, Miller forcing and random real forcing do not add Cohen
reals.

The next chapter leads to answering a question of Matet and Pawlikowski. In
\ref{conmatet} we show that it is consistent that there exists a $p$-point
generated by $\aleph_1$ elements which is not a $q$-point and that for every
$\psi\in\baire$ and a family $\F$ of $\aleph_1$ partial infinite functions
$h:\dom(h)\longrightarrow \omega$ such that $h(n)<\psi(n)$ for $n\in\dom(h)
\subseteq\omega$ there is $W\in\prod\limits_{n\in\omega}[\psi(n)]^{\textstyle
n+1}$ with $(\forall h\in\F)(\exists^\infty n\in\dom(h))(h(n)\in W(n))$. 
Several general results on preserving special properties of ultrafilters are
presented on the way to this solution.

A starting point for chapter 7 was a problem of Balcerzak and Plewik. We show
that the Balcerzak--Plewik number $\kappa_{\rm BP}$ (see \ref{BPnumber}) is
bounded by a cardinal invariant related to the strong PP--property (in
\ref{kbelPP}). Next we show the consistency of both ``$\dominating<\kappa_{\rm
BP}$'' (in \ref{BPabod}) and ``$\kappa_{\rm BP}<\con$'' (in \ref{BPbelcon}). 
We treat our solution as a good opportunity to look at various properties of
forcing notions related to the PP--property (and corresponding cardinal
invariants).  

\section{Notation} Most of our notation is standard and compatible with
that of classical textbooks on Set Theory (like Bartoszy\'nski Judah
\cite{BaJu95} or Jech \cite{J}). However in forcing we keep the convention
that {\em a stronger condition is the larger one}. 

\begin{notation}
\label{notacja}
\begin{enumerate}
\item $\mbR^{{\geq}0}$ stands for the set of non-negative reals. The integer
part of a real $r\in\mbR^{{\geq}0}$ is denoted by $\lfloor r\rfloor$.
\item For two sequences $\eta,\nu$ we write $\nu\vartriangleleft\eta$ whenever
$\nu$ is a proper initial segment of $\eta$, and $\nu\trianglelefteq\eta$ when
either $\nu\vartriangleleft\eta$ or $\nu=\eta$. The length of a sequence
$\eta$ is denoted by $\lh(\eta)$.
\item A {\em tree} is a family of finite sequences closed under initial
segments. (In \ref{quasitree} we will define more general objects.) For a tree
$T$ the family of all $\omega$--branches through $T$ is denoted by $[T]$. We
may use the notation $\lim(T)$ for this object too (see \ref{quasitree}, note
that a tree is a quasi tree).
\item The quantifiers $(\forall^\infty n)$ and $(\exists^\infty n)$ are
abbreviations for  
\[(\exists m\in\omega)(\forall n>m)\quad\mbox{ and }\quad(\forall m\in\omega)
(\exists n>m),\] 
respectively.
\item For a function $h:X\longrightarrow X$ and an integer $k$ we define
$h^{(k)}$ as the $k^{\rm th}$--iteration of $h$: $h^{(1)}=h$, $h^{(k+1)}=
h\circ h^{(k)}$.
\item For a set $X$,\ \ \ $[X]^{\textstyle{\leq}\omega}$,
$[X]^{\textstyle{<}\omega}$ and ${\mathcal P}(X)$ will stand for families of
countable, finite and all, respectively, subsets of the set $X$. The family of
$k$-element subsets of $X$ will be denoted by $[X]^{\textstyle k}$. The set of
all finite sequences with values in $X$ is called $X^{\textstyle {<}\omega}$
(so domains of elements of $X^{\textstyle {<}\omega}$ are integers). The
collection of all {\em finite partial} functions from $\omega$ to $X$ is
$X^{\mathunderaccent\smile-3 \omega}$. 
\item The Cantor space $\can$ and the Baire space $\baire$ are the spaces of
all functions from $\omega$ to $2$, $\omega$, respectively, equipped with
natural (Polish) topology. 
\item For a forcing notion $\p$, $\Gamma_\p$ stands for the canonical
$\p$--name for the generic filter in $\p$. With this one exception, all
$\p$--names for objects in the extension via $\p$ will be denoted with a dot
above (e.g.~$\dot{\tau}$, $\dot{X}$).
\item $\con$ stands for the cardinality of the continuum. The dominating
number (the minimal size of a dominating family in $\baire$ in the ordering of
eventual dominance) is denoted by $\dominating$ and the unbounded number (the
minimal size of an unbounded family in that order) is called
$\unbounded$. $\M$, $\N$ stand for the $\sigma$--ideals of meager and null
sets on the real line, respectively. 
\end{enumerate}
\end{notation}

\begin{gendef}
\label{gendef}
\begin{enumerate}
\item For an ideal ${\mathcal J}$ of subsets of a space $X$ we define its
cardinal characteristics (called {\em additivity}, {\em covering number}, {\em
uniformity} and {\em cofinality}, respectively):
\begin{itemize}
\item $\add({\mathcal J})=\min\{|{\mathcal A}|: {\mathcal A} \subseteq
  {\mathcal J} \ \&\ \bigcup {\mathcal A} \notin {\mathcal J}\}$,
\item $\cov({\mathcal J})=\min\{|{\mathcal A}|: {\mathcal A} \subseteq
  {\mathcal J} \ \&\ \bigcup {\mathcal A} =X\}$, 
\item $\unif({\mathcal J})=\min\{|Y|: Y \subseteq X \ \& \ Y\notin {\mathcal
    J}\}$, 
\item $\cof({\mathcal J})=\min\{|{\mathcal A}|: {\mathcal A} \subseteq
  {\mathcal J}\ \&\ (\forall A\in {\mathcal J})(\exists B\in {\mathcal A})(A
  \subseteq B)\}$. 
\end{itemize}
\item Assume that $X,Y$ are Polish spaces and $R\subseteq X\times Y$ is a Borel
relation. Suppose that $\V\subseteq\V'$ are models of ZFC and that all
parameters we need are in $\V$. We say that the extension $(\V,\V')$ has {\em
the $R$-localization property} if  
\[(\forall x\in X\cap\V')(\exists y\in Y\cap\V)((x,y)\in R).\] 
If $x\in X\cap\V'$, $y\in Y\cap\V$ and $(x,y)\in R$ then we say that {\em $y$
$R$-localizes $x$}. 

We say that a forcing notion $\p$ has the $R$--localization property if every
generic extension of $\V$ via $\p$ has this property.
\item For a relation $R\subseteq X\times Y$ we define two cardinal numbers
(the unbounded and the dominating number for $R$):
\[\unbounded(R)=\min\{|B|:(\forall y\in Y)(\exists x\in B)((x,y)\notin R)\}\]
\[\dominating(R)=\min\{|D|:(\forall x\in X)(\exists y\in D)((x,y)\in R)\}.\]
\end{enumerate}
\end{gendef}

\chapter{Basic definitions} 
In this chapter we introduce our heroes: forcing notions built of weak
creatures. The Prologue is intended to give the reader some intuitions
needed to get through a long list of definitions. A general scheme is
presented in the second part, where we define weak creatures, sub--composition
operations, weak creating pairs $(K,\Sigma)$ and corresponding forcing notions
$\q_{\C(\nor)}(K,\Sigma)$. However, in practice (at least in this paper) we
will be interested in two special cases of the scheme. The first main family
of weak creating pairs (and related forcing notions) are creating pairs
determined by composition operations on creatures. The second family consists
of tree--creating pairs coming from tree compositions on tree--creatures. 
These two options are introduced in the following two parts of the chapter. It
should be underlined here that the rest of the paper will deal with these two
(essentially disjoint and parallel) cases of the general scheme. In the last
section we give some justifications for our work in the next chapter, showing
that without extra care our schema may result in forcing notions collapsing
$\aleph_1$.   
\medskip

\noindent Note: Our terminology (weak creatures, creatures, tree--creatures
etc) might be slightly confusing, but it was developed during a long period of
time (see introduction) and large parts of it are established in literature
already. 
\medskip

\noindent {\bf Basic Notation:} In this paper $\bH$ will stand for a function
with domain $\omega$ such that $(\forall m\in\omega)(|\bH(m)|\geq 2)$. We
usually assume that $0\in \bH(m)$ (for all $m\in\omega$); if it is not the
case then we fix an element of $\bH(m)$ and we use it whenever appropriate
notions refer to $0$. Moreover we fix ``a sufficiently large'' uncountable
regular cardinal $\chi$ and we assume that at least $\bH\in {\mathcal
H}(\chi)$ (the family of sets of cardinality hereditarily less than $\chi$)
or, what is more natural, even $\bH\in {\mathcal H}(\aleph_1)$. 

\setcounter{section}{-1}

\section{Prologue}
If one looks at forcing notions appearing naturally in the Set Theory of Reals
(i.e.~the forcing notions adding a real with certain properties and preserving
various properties of the ground model reals) then one realizes that they
often have a common pattern. A condition in such a forcing notion determines
an initial segment of the real we want to add and it puts some restrictions on
possible further extensions of the initial segment. When we pass to a stronger
condition we extend the determined part of the generic real and we put more
restrictions on possible extensions. But we usually demand that the amount of
freedom which is left by the restrictions still goes to infinity in a
sense. The basic part of the definition of such a forcing notion is to
describe the way a condition puts a restriction on possible initial segments
of the generic real. Typically the restriction can be described locally by use
of ``atoms'' or ``black boxes'', in our terminology called {\em weak
creatures}. So a weak creature $t$ has a domain (contained in finite
sequences) and for each sequence $w$ from the domain it gives a family of  
extensions of $w$ (this is described by a relation $\val[t]$: if $\langle
u,v\rangle\in\val[t]$ then $v$ is an allowed extension of $u$). Moreover, such
a $t$ has a norm $\nor[t]$ which measures the amount of freedom it
leaves. Further, we are told what we are allowed to do with weak creatures:
typically we may shrink them, glue together or just forget about them
(i.e.~{\em omit} them). The results of permitted operations on a family $\cS$
of weak creatures are elements of $\Sigma(\cS)$ in our notation (where
$\Sigma$ is a {\em sub-composition operation} on the considered family $K$ of
weak creatures), see 1.1.4. Now a condition in our forcing notion can be
viewed as $(w,\cS)$, where $w$ is a finite sequence (the determined part of
the generic real) and $\cS$ is a countable family of weak creatures satisfying
some demands on its structure and requirements on $\nor[t]$ for $t\in
\cS$. When we want to build a stronger condition then we take $t\in \cS$ such
that $w$ is in the domain of $t$ and we pick up one of the possible extensions
of $w$ allowed by $t$. We may repeat this procedure finitely many times and we
get a sequence $w^*$ extending $w$. Next we choose a family $\cS^*$ of weak
creatures such that each $s\in \cS^*$ is obtained by permitted procedures from
some $\cS_s\subseteq\cS$ (i.e.~$s\in\Sigma(\cS_s)$). The pair $(w^*,\cS^*)$ is
an extension of $(w,\cS)$ provided $\cS^*$ satisfies the structure demands and
norm requirements.

However, this general schema breaks to two cases, which, though very similar,
are of different flavors. In the first case we demand that the family $\cS$
in a condition $(w,\cS)$ has a linear structure. Then we usually represent the
condition as $(w,t_0,t_1,t_2,\ldots)$, where for some sequence $0\leq m_0<m_1<
m_2<\ldots<\omega$
\begin{quotation}
\noindent $w$ is a sequence of length $m_0$ and for each $i<\omega$:

\noindent $t_i$ is a weak creature saying in which way sequences of length
$m_i$ may be extended to sequences of length $m_{i+1}$.
\end{quotation}
So it is natural in this context to consider only weak creatures $t$ such that
for some integers $m^t_{\dn}<m^t_{\up}$ ($\dn$ stands for ``down''), the
domain of $t$ is contained in sequences of length $m^t_{\dn}$ and every
extension of a sequence from the domain allowed by $t$ is of length 
$m^t_{\up}$. In other words we require that if $\langle u,v\rangle\in\val[t]$
then $\lh(u)=m^t_{\dn}$ and $\lh(v)=m^t_{\up}$. In applications the domain of
the relation $\val[t]$ consist of {\em all legal} sequences of length
$m^t_{\dn}$. Let us describe a simple example of this kind. Consider the
Silver forcing $\q$ ``below $2^n$'': a condition in $\q$ is a function $p:
\dom(p)\longrightarrow\omega$ such that
\[\dom(p)\subseteq\omega\quad\&\quad |\omega\setminus\dom(p)|=\omega\quad\&
\quad (\forall n\in\dom(p))(p(n)<2^n).\]
Let us look at this forcing in a different way. 

Let $K$ consist of all triples $t=(\nor[t],\val[t],\dis[t])$ such that
\begin{itemize}
\item $\dis[t]=(m^t,k^t)$ where $m^t<\omega$ and $k^t\in\{*\}\cup 2^{m^t}$,
\item $\nor[t]=m^t$ if $k^t=*$ and $\nor[t]=0$ otherwise,
\item $\val[t]=\{\langle u,v\rangle\in \prod\limits_{i<m^t} 2^i\times
\prod\limits_{i\leq m^t} 2^i: u\vartriangleleft v\ \ \&\ \ (k^t\neq *\
\Rightarrow\ v(m^t)=k^t)\}$.
\end{itemize}
For $\cS\subseteq K$ we let $\Sigma(\cS)=\emptyset$ if $|\cS|\neq 1$ and

if $t\in K$, $k^t\neq *$ then $\Sigma(\{t\})=\{t\}$,

if $t\in K$, $k^t= *$ then $\Sigma(\{t\})$ consists of all $s\in K$ with $m^s=
m^t$. 

\noindent Now, a condition $p$ in $\q$ may be represented as a sequence $(w^p,
t^p_0,t^p_1,\ldots)$, where
\begin{enumerate}
\item[(a)] each $t^p_i$ is from $K$,
\item[(b)] $w^p$ is a finite sequence of length $m^{t^p_0}$ such that $w^p(n)
<2^n$ for $n<m^{t^p_0}$,
\item[(c)] $m^{t^p_i}=m^{t^p_0}+i$,
\item[(d)] $\lim\sup\limits_{i\to\omega}\nor[t^p_i]=\infty$.
\end{enumerate}
The order $\leq$ of $\q$ is defined by $(w^p,t^p_0,t^p_1,t^p_2,\ldots)\leq
(w^q,t^q_0,t^q_1,t^q_2,\ldots)$ if and only if for some $N<\omega$:
\begin{quotation}
\noindent $w^p\trianglelefteq w^q$, $\lh(w^q)=\lh(w^p)+N$,

\noindent $\langle w^q\rest\lh(w^p)+i,w^q\rest (w^p)+i+1\rangle\in\val[t^p_i]$
for each $i<N$, and

\noindent $t^q_j\in \Sigma(t^p_{N+j})$ for every $j<\omega$.
\end{quotation}
The pair $(K,\Sigma)$ is an example of a {\em creating pair} and the forcing
notion $\q$ (represented as above) is $\q^*_{{\rm w}\infty}(K,\Sigma)$ (see
1.2.2 and 1.2.4). 

On the other pole of possible weak creatures we have those which provide
possible extensions for only one sequence $\eta$ (i.e.~those $t$ for which
$|\dom(\val[t])|=1$). Weak creatures of this type are called {\em
tree--creatures} and they say to us simply:
\begin{quotation}
{\em \noindent I know what the restrictions on extensions of a single
sequence $\eta$ are, and I do not look at other sequences at all.}
\end{quotation}
Tree--creatures are fundamental for building forcing notions in which
conditions are trees of a special kind. In these forcing notions a condition
$p=(w,\cS)$ is such that $w$ is a root (stem) of a tree $T^p$ and each
$t\in\cS$ is a part of the tree $T^p$; usually such a $t$ describes how one
passes from an element $\eta\in T^p$ to its extensions in $T^p$ (not
necessarily immediate successors). It is natural to put some requirements on
sub--composition operations $\Sigma$ when the weak creatures we consider are
tree--creatures, and this leads to the definition of {\em tree composition}
and {\em tree--creating pair}, see 1.3.3. Moreover, it turns out to be very
practical to consider special demands on the norms $\nor[t]$ to take an
advantage of the tree--form of a condition, see 1.3.5. (Note that in further
definitions we do not require that $T^p$ is a tree but we demand that it is a
{\em quasi tree} only. This will simplify the notation a little bit.) Let us
illustrate this by a suitable representation of the Laver forcing ${\mathbb
L}$. Recall that a condition in ${\mathbb L}$ is a tree $T\subseteq\fseo$ such
that if $\eta\in T$, $\mrot(T)\trianglelefteq\eta$ then $|\suc_T(\eta)|=
\omega$.

Let $K$ consist of all triples $t=(\nor[t],\val[t],\dis[t])$ such that
\begin{itemize}
\item $\dis[t]=(\eta^t,A^t)$, where $A^t\in\iso$ and $\eta^t\in\fseo$,
\item $\nor[t]=\lh(\eta^t)$,
\item $\val[t]=\{\langle \eta^t,\nu\rangle: \eta^t\vartriangleleft\nu\ \&\
\lh(\nu)=\lh(\eta^t)+1\ \&\ \nu(\lh(\eta^t))\in A^t\}$.
\end{itemize}
For $t\in K$ we let $\Sigma(\{t\})=\{s\in K:\eta^s=\eta^t\ \&\ A^s\subseteq
A^t\}$ and for $\cS\subseteq K$ with $|\cS|\neq 1$ we declare $\Sigma(\cS)=
\emptyset$. Now, a condition $p$ in ${\mathbb L}$ can be represented as
$(\eta^p,\langle t^p_\nu: \nu\in T^p\rangle)$, where $T^p\subseteq\fseo$ is a
tree such that $\mrot(T^p)=\eta^p$ and for each $\nu\in T^p$,
$\mrot(T^p)\trianglelefteq\nu$ we have $\suc_{T^p}(\nu)=\{\rho: \langle
\nu,\rho\rangle\in\val[t^p_\nu]\}$ (so $\cS$ is $\{t^p_\nu: \nu\in T^p\}$
here). Moreover, we demand that for each infinite branch $\eta$ through $T^p$
the norms $\nor[t^p_{\eta\rest i}]$ go to infinity. The order $\leq$ of
${\mathbb L}$ is given by $(\eta^p,\langle t^p_\nu:\nu\in T^p\rangle)\leq
(\eta^q,\langle t^q_\nu:\nu\in T^q\rangle)$ if and only if $\eta^q\in T^p$ and
\[(\forall\nu\in T^q)(\nu\in T^p\ \&\ t^q_\nu\in\Sigma(t^p_\nu)).\]
The pair $(K,\Sigma)$ defined above is a tree creating pair and the forcing
notion ${\mathbb L}$ is $\q^{\tree}_1(K,\Sigma)$.   

\section{Weak creatures and related forcing notions}

\begin{definition}
\label{weakcreat}
\begin{enumerate}
\item A triple $t=(\nor,\val,\dis)$ is {\em a weak creature for $\bH$} if: 
\begin{enumerate}
\item[(a)] $\nor\in\mbR^{{\geq}0}$,
\item[(b)] $\val$ is a non-empty subset of
\[\{\langle x,y\rangle\in\bigcup_{m_0<m_1<\omega}[\prod_{i<m_0}\bH(i)\times
\prod_{i<m_1}\bH(i)]: x\vartriangleleft y\},\] 
\item[(c)] $\dis\in{\mathcal H}(\chi)$.
\end{enumerate}
The family of all weak creatures for $\bH$ is denoted by $\WCR[\bH]$.
\item In the above definition we write $\nor=\nor[t]$, $\val=\val[t]$ and
$\dis=\dis[t]$. 
\end{enumerate}
[$\val$ is for value, $\nor$ is for norm, $\dis$ is for distinguish.]
\end{definition}

\begin{remark}
The $\dis[t]$ in a weak creature $t$ plays the role of an additional parameter
which allows as to have distinct creatures with the same values of $\val$ and
$\nor$. This may be sometimes important in defining sub-composition operations
on $K$ (see \ref{subcom} below): we will be able to have distinct values of
$\Sigma(t_0)$, $\Sigma(t_1)$  though $\val[t_0]=\val[t_1]$ and $\nor[t_0]=
\nor[t_1]$. One may think that this additional parameter describes the way the
weak creature $t$ was constructed (while $\val[t]$, $\nor[t]$ give the final
effect of the construction). We may sometimes ``forget'' to mention $\dis[t]$
explicitly -- in most of the results and applications $\dis[t]$ might be
arbitrary. In the examples we construct, if we do not mention $\dis[t]$ we
mean that either it is 0 or its form is clear.
\end{remark}

\begin{definition}
\label{morecreat}
\begin{enumerate}
\item If we omit $\bH$ we mean for some $\bH$ or the $\bH$ is clear from the
context, etc.  
\item We say that $\bH$ is {\em finitary} (or {\em of a countable character},
respectively) if $\bH(n)$ is finite (countable, resp.) for each $n\in\omega$.
We say that $K\subseteq\WCR[\bH]$ is {\em finitary} if $\bH$ is finitary and
$\val[t]$ is finite for each $t\in K$. 
\end{enumerate}
\end{definition}

\begin{definition}
\label{subcom}
Let $K\subseteq \WCR[\bH]$.
\begin{enumerate}
\item A function $\Sigma:[K]^{\textstyle{\leq}\omega}\longrightarrow\cP(K)$ is
a {\em sub-composition operation on $K$} if: 
\begin{enumerate}
\item[(a)] {\em (transitivity)}\ \ \ if $\cS\in [K]^{\textstyle{\leq}\omega}$
and for each $s\in\cS$ we have $s\in\Sigma(\cS_s)$ then $\Sigma(\cS)\subseteq
\Sigma(\bigcup\limits_{s\in\cS}\cS_s)$,
\item[(b)]  $r\in\Sigma(r)$ for each $r\in K$ and $\Sigma(\emptyset)=
\emptyset$. 
\end{enumerate}
[Note that $\Sigma(\cS)$ may be empty for non-empty $\cS$; in future defining
$\Sigma$ we will describe it only for the cases it provides a non-empty
result, in all other cases we will assume that $\Sigma(\cS)=\emptyset$.]
\item In the situation described above the pair $(K,\Sigma)$ is called {\em a
weak creating pair for $\bH$}. 
\item Suppose that $(K,\Sigma)$ is a weak creating pair, $t_0,t_1\in K$. We
say that $t_0,t_1$ are $\Sigma$--equivalent (and we write then
$t_0\sim_{\Sigma} t_1$) if $\nor[t_0]=\nor[t_1]$, $\val[t_0]=\val[t_1]$ and
for each $\cS\subseteq [K\setminus\{t_0,t_1\}]^{\textstyle \leq\omega}$ we
have $\Sigma(\cS\cup\{t_0\})=\Sigma(\cS\cup\{t_1\})$.
\end{enumerate}
\end{definition}

\begin{remark}
Note that the relation $\sim_\Sigma$ as defined in \ref{subcom}(3) does
not have to be transitive in a general case (so, perhaps, we should not use
the name {\em $\Sigma$--equivalent}). However, if $(K,\Sigma)$ is either a
creating pair (see \ref{compos}) or a tree--creating pair (see
\ref{treecreature}) then $\sim_\Sigma$ is an equivalence relation. Then, if
additionally $\Sigma(\cS)$ is non-empty for finite $\cS$ only, the value of
$\Sigma(\cS)$ depends on the $\sim_\Sigma$--equivalence classes of elements of
$\cS$ only. Therefore we will tend to think in these situations that we
identify all $\sim_\Sigma$--equivalent elements of $K$ (or just consider a
selector $K^*\subseteq K$ of $K/{\sim}_\Sigma$). If $\Sigma(\cS)$ may be
non-empty for an infinite $\cS\subseteq K$ (which may happen for
tree--creating pairs), then we have to be more careful before we consider this
identification: we should check that the values of $\Sigma$ depend on
$\sim_\Sigma$--equivalence classes only.
\end{remark}

\begin{definition}
\label{basis}
Let $(K,\Sigma)$ be a weak creating pair for $\bH$.
\begin{enumerate}
\item  For a weak creature $t\in K$ we define its basis (with respect to
$(K,\Sigma)$) as 
\[\basis(t)=\{w\in\bigcup_{m<\omega}\prod_{i<m}\bH(i):\; (\exists
s\in\Sigma(t))(\exists u)(\langle w,u\rangle\in\val[s])\}.\] 
\item For $w\in\bigcup\limits_{m<\omega}\prod\limits_{i<m}\bH(i)$ and
$\cS\in [K]^{\textstyle {\leq}\omega}$ we define the set $\pos(w,\cS)$ of
possible extensions of $w$ from the point of view of $\cS$ (with respect to
$(K,\Sigma)$) as: 
\[\pos^*(w,\cS)=\{u: (\exists s\in\Sigma(\cS))(\langle
w,u\rangle\in\val[s])\},\]     
\[\hspace{-0.5cm}
\begin{array}{ll}
\pos(w,\cS)=\{u:&\!\!\!\!\mbox{there are disjoint sets }\cS_i\mbox{ (for 
$i<m<\omega$) with }\bigcup\limits_{i<m}\cS_i=\cS\\
\ &\mbox{and a sequence }0<\ell_0<\ldots<\ell_{m-1}<\lh(u)\mbox{ such that}\\
\ & u{\restriction} \ell_0\in\pos^*(w,\cS_0)\ \mbox{ and}\\
\ & u{\restriction} \ell_1\in\pos^*(u{\restriction}\ell_0,\cS_1)\ \&\ \ldots\
\&\ u\in\pos^*(u{\restriction} \ell_{m-1},\cS_{m-1})\}.\\  
\end{array}
\]
\item Whenever we use $\basis$ or $\pos$ we assume that the weak creating pair
$(K,\Sigma)$ with respect to which these notions are defined is understood.
\end{enumerate}
\end{definition}

\begin{definition}
\label{maindef}
Suppose $(K,\Sigma)$ is a weak creating pair for $\bH$ and $\C(\nor)$ is a
property of $\omega$-sequences of weak creatures from $K$ (i.e.~$\C(\nor)$ can
be thought of as a subset of $K^\omega$). We define the forcing notion
$\q_{\C(\nor)}(K,\Sigma)$ by  
\medskip

\noindent{\bf conditions} are pairs $(w,T)$ such that for some $k_0<\omega$: 
\begin{enumerate}
\item[(a)] $w\in\prod\limits_{i<k_0}\bH(i)$
\item[(b)] $T=\langle t_i: i<\omega\rangle$ where:
\begin{enumerate}
\item[(i)]    $t_i\in K$  for each $i$, 
\item[(ii)]   $w\in\basis(t_i)$ for some $i<\omega$ and for each
$u\in\pos(w,\{t_i: i\in I_0\})$, $I_0\subseteq\omega$ there is
$i\in\omega\setminus I_0$ such that $u\in\basis(t_i)$,
\end{enumerate}
\item[(c)] the sequence $\langle t_i: i<\omega\rangle$ satisfies the condition
$\C(\nor)$;
\end{enumerate}
\smallskip

\noindent{\bf the order} is given by: $(w_1,T^1)\leq (w_2,T^2)$ if and only if 

\noindent for some disjoint sets $\cS_0,\cS_1,\cS_2,\ldots\subseteq\omega$ we
have:  
\[w_2\in\pos(w_1,\{t^1_\ell:\ell\in\cS_0\})\quad\mbox{ and }\quad t^2_i\in
\Sigma(\{t^1_{\ell}: \ell\in\cS_{i+1}\})\mbox{ for each }i<\omega\] 
(where $T^\ell=\langle t^\ell_i: i<\omega\rangle$).

If $p=(w,T)$ we let $w^p=w$, $T^p=T$ and if $T^p=\langle t_i: i<\omega\rangle$
then we let $t^p_i=t_i$. We may write $(w,t_0,t_1,\ldots)$ instead of $(w,T)$
(when $T=\langle t_i: i<\omega\rangle$).
\end{definition}

\begin{proposition}
If $(K,\Sigma)$ is a weak creating pair and $\C(\nor)$ is a property of
sequences of elements of $K$ then $\q_{\C(\nor)}(K,\Sigma)$ is a forcing
notion. 
\end{proposition}

\begin{remark}
The reason for our notation $\C(\nor)$ for the property relevant for {\bf (c)}
of \ref{maindef} is that in the applications this conditions will say that the
norms $\nor[t_i]$ go to the infinity in some sense. Some of the possibilities
here are listed in \ref{conditions} below.
\end{remark}

\begin{definition} 
\label{conditions}
For a weak creature $t$ let us denote 
\[m_{\dn}(t)=\min\{\lh(u): u\in \dom(\val[t])\}.\] 
We introduce the following (basic) properties of sequences of weak creatures
which may serve as $\C(\nor)$ in \ref{maindef}: 
\begin{enumerate}
\item[({\rm s}$\infty$)]\ \ \ A sequence $\langle t_i: i<\omega\rangle$
satisfies $\C^{{\rm s}\infty}(\nor)$ if and only if
\[(\forall i<\omega)(\nor[t_i]>\max\{i,m_{\dn}(t_i)\}).\]

\item[($\infty$)\ ]\ \ \ \ A sequence $\langle t_i: i<\omega\rangle$ satisfies 
$\C^\infty(\nor)$ if and only if
\[\lim\limits_{i\rightarrow\omega}\nor[t_i]=\infty.\]

\item[({\rm w$\infty$})]\ \ \ A sequence $\langle t_i: i<\omega\rangle$
satisfies $\C^{{\rm w}\infty}(\nor)$ if and only if
\[\lim\sup\limits_{i\rightarrow\omega}\nor[t_i]=\infty.\]
\end{enumerate}
Let $f:\omega\times\omega\longrightarrow\omega$. We define the property
introduced by $f$ by
\begin{enumerate}
\item[($f$)\ ]\ \ \ \ A sequence $\langle t_i: i<\omega\rangle$ satisfies
$\C^f(\nor)$ if and only if 
\[(\forall k<\omega)(\forall^\infty i)(\nor[t_i]> f(k,m_{\dn}(t_i))).\]
\end{enumerate}
For notational convenience we will sometimes use the empty norm condition:
\begin{enumerate}
\item[($\emptyset$)\ ]\ \ \ \ Each sequence $\langle t_i: i<\omega\rangle$
satisfies $\C^\emptyset(\nor)$.  
\end{enumerate}
\smallskip 

The forcing notions corresponding to the above properties (for a weak creating
pair $(K,\Sigma)$) will be denoted by $\q_{{\rm s}\infty}(K,\Sigma)$, 
$\q_{\infty}(K,\Sigma)$, $\q_{{\rm w}\infty}(K,\Sigma)$, $\q_{f}(K,\Sigma)$
and $\q_{\emptyset}(K,\Sigma)$, respectively.
\end{definition}

\begin{remark}
1)\ \ \ Note that the second component of a pair $(w,T)\in\q_{\C(\nor)}
(K,\Sigma)$ is {\em a sequence} of weak creatures, and in the most general
case the order of its members may be important. For example the property
$\C^{{\rm s}\infty}(\nor)$ introduced in \ref{conditions} is not permutation
invariant and some changes of the order in the sequence $\langle t_i:
i<\omega\rangle$ may produce a pair $(w,T')$ which is not a legal
condition. This is not what we would like to have here, so in applications in
which this kind of problems appears we will restrict ourselves to suborders
$\q^*_{\C(\nor)}(K,\Sigma)$ of $\q_{\C(\nor)}(K,\Sigma)$ in which we put
additional structure demands on the sequences $\langle t_i:i<\omega\rangle$
(see \ref{forcingstar}). Moreover, to get properness for forcing notions
$\q^*_{{\rm s}\infty}(K,\Sigma)$ we will have to put some demands on
$(K,\Sigma)$ (see \ref{omitoryproper}). These demands will cause that various
variants of the norm condition $\C^{{\rm s}\infty}(\nor)$ result in equivalent
forcing notions (see \ref{sinfty}). So, from the point of view of
applications, the main reason for introducing $\C^{{\rm s}\infty}(\nor)$ is a
notational convenience.\\
2)\ \ \ Note that 
\[\q_{{\rm s}\infty}(K,\Sigma)\subseteq\q_{\infty}(K,\Sigma)\subseteq \q_{{\rm
w}\infty}(K,\Sigma)\subseteq\q_{\emptyset}(K,\Sigma)\]
where the inclusions mean ``suborder'' (but often not ``complete suborder'').
If we put some conditions on $f$ (e.g.~$f$ is fast, see \ref{fast}) then we
may easily have $\q_{f}(K,\Sigma)\subseteq\q_{\infty}(K,\Sigma)$.\\
3)\ \ \ In our applications we will consider the forcing notions
$\q_{f}(K,\Sigma)$ only for functions $f:\omega\times\omega\longrightarrow 
\omega$ which are growing fast enough (see \ref{fast} below).
\end{remark}

\begin{definition}
\label{fast}
A function $f:\omega\times\omega\longrightarrow\omega$ is {\em fast} if 
\[(\forall k\in\omega)\big(\forall \ell\in\omega)(f(k,\ell)\leq f(k,\ell+1)\ \
\&\ \ 2\cdot f(k,\ell)<f(k+1,\ell)\big).\]
The function $f$ is {\em $\bH$-fast} if additionally  ($\bH$ is finitary and)
for each $k,\ell\in\omega$:
\[2^{\fH(\ell)}\cdot(f(k,\ell)+\fH(\ell)+2) < f(k+1,\ell),\]  
where $\fH(\ell)=|\prod\limits_{i<\ell}\bH(i)|$.
\end{definition}

\begin{definition}
\label{thereal}
Suppose that $(K,\Sigma)$ is a weak creating pair and $\C(\nor)$ is a property
of sequences of elements of $K$. Let $\dot{W}$ be a
$\q_{\C(\nor)}(K,\Sigma)$-name such that 
\[\forces_{\q_{\C(\nor)}(K,\Sigma)}\Wtil=\bigcup\{w^p: p\in
\Gamma_{\q_{\C(\nor)}(K,\Sigma)}\}.\]  
\end{definition}

\begin{proposition}
Suppose that $(K,\Sigma)$ is a weak creating pair and $\C(\nor)$ is a property
such that the forcing notion $\q_{\C(\nor)}(K,\Sigma)$ is non-empty. Then:
\begin{enumerate}
\item  $\forces_{\q_{\C(\nor)}(K,\Sigma)}$ ``$\Wtil$ is a member of
$\prod\limits_{i<\omega}\bH(i)$''. 
\item  If $(\forall i\in\omega)(\bH(i)=2)$ then
$\forces_{\q_{\C(\nor)}(K,\Sigma)}$ ``$\Wtil$ is a real''.  
\item If for every $t\in K$, $u\in\basis(t)$ the set $\pos(u,t)$ has at least
two elements 

then $\forces_{\q_{\C(\nor)}(K,\Sigma)}$``$\Wtil\notin \V$''. 
\end{enumerate}
\end{proposition}

\begin{remark}
1)\ \ \ We will always assume that the considered weak creating pairs
$(K,\Sigma)$ (and norm conditions $\C(\nor)$) are such that
$\q_{\C(\nor)}(K,\Sigma)\neq\emptyset$. Usually, it will be enough that $K$
contains enough creatures with large norms and in each particular example this
requirement will be easy to verify.\\
2)\ \ \ In general, the $\dot{W}$ defined in \ref{thereal} does not have to
encode the generic filter. We may formulate a condition ensuring this. Let
$(K,\Sigma)$ be a weak creating pair and $\C(\nor)$ be a norm condition such
that $\q_{\C(\nor)}(K,\Sigma)$ is not empty. For $p\in\q_{\C(\nor)}(K,\Sigma)$
define
\[S(p)\stackrel{\rm def}{=}\{w\in\bigcup_{n\in\omega}\prod_{i<n}\bH(i):
(\exists q\geq p)(w\trianglelefteq w^q)\}\]
Clearly $S(p)$ is a subtree of $\bigcup\limits_{n\in\omega}\prod\limits_{i<n}
\bH(i)$. Moreover, for each $w\in\bigcup\limits_{n\in\omega}\prod\limits_{i<n} 
\bH(i)$ and $p,q\in \q_{\C(\nor)}(K,\Sigma)$:
\[\begin{array}{lll}
p\not\forces_{\q_{\C(\nor)}(K,\Sigma)}\mbox{``}w\ntriangleleft\dot{W}\mbox{''}
&\mbox{ if and only if }&w\in S(p)\quad\mbox{ and}\\
p\forces_{\q_{\C(\nor)}(K,\Sigma)}\mbox{``}\dot{W}\in [S(q)]\mbox{''}
&\mbox{ if and only if }&S(p)\subseteq S(q).\\
  \end{array}\]
Now we may define a $\q_{\C(\nor)}(K,\Sigma)$--name $\dot{H}$ by
\[\forces_{\q_{\C(\nor)(K,\Sigma)}}\dot{H}=\{p\in\q_{\C(\nor)}(K,\Sigma):
\dot{W}\in [S(p)]\}\]
and we may want to claim that $\forces\dot{H}=\Gamma_{\q_{\C(\nor)}(K,
\Sigma)}$. But for this we need to know that any two conditions in
$\dot{H}$ are compatible. A sufficient and necessary requirement for this is:
\begin{enumerate}
\item[($\boxdot$)] {\em if} $p,q\in \q_{\C(\nor)}(K,\Sigma)$ and
$S(p)\subseteq S(q)$

{\em then} $p\forces q\in\Gamma_{\q_{\C(\nor)}(K,\Sigma)}$ (or in other words
$p\geq q$ modulo the equivalence of conditions).
\end{enumerate}
In most of our examples and applications, the condition ($\boxdot$) will be
easy to check. We will not mention it in future as we will not use its
consequences. 

\noindent Note however that it is very easy to build examples of weak creating
pairs $(K,\Sigma)$ (even creating pairs or tree creating pairs) for which
($\boxdot$) fails. Some of these examples might appear naturally.
\end{remark}

\section{Creatures}
Now we will deal with the first specific case of the general scheme: creating
pairs and forcing notions $\q^*_{\C(\nor)}(K,\Sigma)$. Notation and
definitions introduced here are applicable to this case only and should not be
confused with that for tree--creating pairs. 

\begin{definition}
\label{creatures}
Let $t$ be a weak creature for $\bH$.
\begin{enumerate}
\item If there is $m<\omega$ such that $(\forall\langle u,v\rangle\in\val[t])
(\lh(u)=m)$\\ 
then this unique $m$ is called $m^t_{\dn}$.
\item If there is $m<\omega$ such that $(\forall\langle u,v\rangle\in\val[t])
(\lh(v)=m)$\\
then this unique $m$ is called $m^t_{\up}$.
\item If both $m_{\dn}^t$ and $m_{\up}^t$ are defined then $t$ is called {\em
an $(m^t_\dn, m^t_\up)$--creature} (or just {\em a creature}).  
\item $\CR_{m_{\dn},m_{\up}}[\bH]=\{t\in\WCR[\bH]: m^t_{\dn}=m_{\dn}$ and
$m^t_{\up}=m_{\up}\}$,

$\CR[\bH]=\bigcup\limits_{m_{\dn}<m_{\up}<\omega}\CR_{m_{\dn},m_{\up}}[\bH]$.
\end{enumerate}
\end{definition}

\begin{definition}
\label{compos}
Suppose that $K\subseteq\CR[\bH]$ and $\Sigma$ is a sub-composition operation
on $K$. We say that $\Sigma$ is {\em a composition} on $K$ (and we say that
$(K,\Sigma)$ is {\em a creating pair for $\bH$}) if:
\begin{enumerate}
\item if $\cS\in [K]^{\textstyle{\leq\omega}}$ and $\Sigma(\cS)\neq\emptyset$
then $\cS$ is finite and for some enumeration $\cS=\{t_0,\ldots,t_{m-1}\}$ we
have $m^{t_i}_{\up}= m^{t_{i+1}}_{\dn}$ for all $i<m-1$, and 
\item for each $s\in\Sigma(t_0,\ldots,t_{m-1})$ we have $m^s_{\dn}=
m^{t_0}_{\dn}$ and $m^s_{\up}=m^{t_{m-1}}_{\up}$. 
\end{enumerate}
In this paper we will always assume that the creating pair under
considerations is additionally {\em nice} and {\em smooth} (see
\ref{niceandsmo} below) and we will not repeat this demand later.
\end{definition}

\begin{remark}
\label{simpler}
Sets of creatures with pairwise distinct $m^t_{\dn}$'s might be naturally
ordered according to this value and therefore in similar situations we
identify sets of creatures with the corresponding sequences of creatures.
\end{remark}

\begin{definition}
\label{candidates}
\begin{enumerate}
\item For $K\subseteq\CR[\bH]$ and a composition operation $\Sigma$ on $K$
we define {\em finite candidates} ($\FC$) and {\em pure finite candidates}
($\PFC$) with respect to $(K,\Sigma)$:
\[\hspace{-0.4cm}\begin{array}{ll}
\FC(K,\Sigma)=&\{(w,t_0,\ldots,t_n):
w\in\basis(t_0)\mbox{ and for each }i\leq n\\
\ & \quad t_i\in K, m^{t_i}_{\up} = m^{t_{i + 1}}_{\dn}\mbox{ and }
\pos(w,t_0,\ldots,t_{i})\subseteq\basis(t_{i+1})\},\\
\end{array}\]
\[\hspace{-0.3cm}\PFC(K,\Sigma)=\{(t_0,\ldots,t_n):(\exists w\in\basis(t_0))
((w,t_0, \ldots,t_n)\in \FC(K))\}.\]  

\item We have a natural partial order $\leq$ on $\FC(K,\Sigma)$ (like in
\ref{maindef}). The partial order $\leq$ on $\PFC$ is defined by 

\noindent $(t_0,\ldots,t_{n-1})\leq (s_0,\ldots,s_{m-1})$\quad if and only if
\quad $m^{t_{n-1}}_{\up}=m^{s_{m-1}}_{\up}$, and 
\[(\forall w\in\basis(t_0))\big(w\in\basis(s_0)\mbox{ and }(w,t_0,\ldots,
t_{n-1})\leq (w,s_0,\ldots,s_{m-1})\big)\]
(so $(t_0)\leq (s_0)$ means that $s_0\in\Sigma(t_0)$ {\em and} $\basis(t_0)
\subseteq\basis(s_0)$).  
\item A sequence $\langle t_0,t_1,t_2,\ldots\rangle$ of creatures from $K$ is 
{\em a pure candidate} with respect to a creating pair $(K,\Sigma)$ if  
\[(\forall i<\omega)(m^{t_i}_{\up}=m^{t_{i+1}}_{\dn})\quad\mbox{ and}\]
\[(\exists w\in\basis(t_0))(\forall i<\omega)(\pos(w,t_0,\ldots,t_i)\subseteq
\basis(t_{i+1})).\] 
The set of pure candidates with respect to $(K,\Sigma)$ is denoted by
$\PC(K,\Sigma)$. The partial order $\leq$ on $\PC(K,\Sigma)$ is defined
naturally.
\item For a norm condition $\C(\nor)$ the family of {\em $\C(\nor)$-normed
pure candidates} is
\[\PC_{\C(\nor)}(K,\Sigma)\stackrel{\rm def}{=}\{\langle t_0,t_1,\ldots\rangle
\in \PC(K,\Sigma): \langle t_0,t_1,\ldots,\rangle \mbox{ satisfies
}\C(\nor)\}. \]
\end{enumerate}
\end{definition}

\begin{definition}
\label{niceandsmo}
Let $(K,\Sigma)$ be a creating pair for $\bH$. We say that
\begin{enumerate}
\item $(K,\Sigma)$ is {\em nice} if for all $t_0,\ldots,t_{n-1}\in K$ and
$s\in\Sigma(t_0,\ldots,t_{n-1})$ we have $\basis(t_0)\subseteq\basis(s)$.
\item $(K,\Sigma)$ is {\em smooth} provided that: 

\noindent{\em if} $(w,t_0,\ldots,t_{n-1})\in\FC(K,\Sigma)$, $m<n$ and
$u\in\pos(w,t_0,\ldots,t_{n-1})$ 

\noindent{\em then}
$u\restriction m^{t_m}_{\dn}\in\pos(w,t_{0},\ldots,t_{m-1})$ and
$u\in\pos(u\restriction m^{t_m}_{\dn},t_m,\ldots,t_{n-1})$. 
\item $K$ is {\em forgetful} if for every creature $t\in K$ we have:
\[[\langle w,u\rangle\in\val[t]\ \& \ w'\in\!\prod_{n<m^t_{\dn}}\!\bH(n)]\ \ 
\Rightarrow\ \ \langle w',w'\conc u\restriction [m^t_{\dn},m^t_{\up})\rangle
\in\val[t].\]   
\item $K$ is {\em full} if $\dom(\val[t])=\prod\limits_{n<m^t_\dn}\bH(n)$ for
every $t\in K$.
\end{enumerate}
\end{definition}
As we said in \ref{compos}, we will always demand that a creating pair is nice
and smooth (but these properties occur naturally in applications). The main
reason for the first assumption is to have the effect presented in
\ref{implications}(2) below and the second demand is to get the conclusion of
\ref{smoessdec}. Before we state these observations let us modify a little bit
the forcing notions we are interested in.

\begin{definition}
\label{forcingstar}
Let $(K,\Sigma)$ be a creating pair and $\C(\nor)$ be a property of
$\omega$-sequences of creatures. The forcing notion
$\q^*_{\C(\nor)}(K,\Sigma)$ is a suborder of $\q_{\C(\nor)}(K,\Sigma)$
consisting of these conditions $(w,t_0,t_1,\dots)$ for which additionally
\begin{enumerate}
\item[$(\boxtimes_{\ref{forcingstar}})$] $\quad\quad\quad (\forall
i<\omega)(m^{t_i}_{\up}=m^{t_{i+1}}_{\dn})$. 
\end{enumerate}
\end{definition}

\begin{remark}
1)\ \ \ The forcing notions introduced in \ref{forcingstar} fit better to the
idea of creatures and compositions on them. Moreover in most of the
applications the forcing notions $\q^*_{\C(\nor)}(K,\Sigma)$ and
$\q_{\C(\nor)}(K,\Sigma)$ will be equivalent. Even in the most general case 
they are not so far from each other; note that if $p\in\q_{\C(\nor)}^*(K,
\Sigma)$ and $q\in\q_{\C(\nor)}(K,\Sigma)$, $p\leq q$ (in $\q_{\C(\nor)}(K,
\Sigma)$) then $q\in\q^*_{\C(\nor)}(K,\Sigma)$. Of course it may happen that
$\q^*_{\C(\nor)}(K,\Sigma)$ is trivial -- this usually suggests that the
tree--approach is more suitable (see \ref{treecreature}).\\
2)\ \ \ Several notions simplify for the forcing notions $\q^*_{\C(\nor)}(K,
\Sigma)$. For example if $t_0,\ldots,t_{n-1}\in K$ are such that $m=
m^{t_0}_{\dn}$, $m^{t_i}_{\up}= m^{t_{i+1}}_{\dn}$ and $w\in\prod\limits_{i<m}
\bH(i)$ then 
\[\begin{array}{ll}
\pos(w,t_0,\ldots,t_{n-1})=\{u: &\mbox{for some }0\!\leq\!
k_1\!<\!\ldots\!<\! k_\ell\!<\! n-1\mbox{ we have}\\ 
\ & u\restriction m^{t_{k_1}}_\up\in\pos^*(w,t_0,\ldots,t_{k_1})\ \&\ \\
\ & u\restriction m^{t_{k_2}}_\up\in\pos^*(u\restriction
m^{t_{k_1}}_\up,t_{k_1+1},\ldots,t_{k_2})\ \&\ \ldots\ \&\\ 
\ & u\in\pos^*(u\restriction m^{t_{k_\ell}}_\up,t_{k_\ell+1},\ldots,
t_{n-1})\}.\\ 
\end{array}\]   
3)\ \ \ The norm condition (${\bf s}\infty$) (see \ref{conditions}) can be
presented slightly simpler for $\q^*_{{\rm s}\infty}(K,\Sigma)$. For 
$(w,t_0,t_1,\ldots)\in\q^*_{{\rm s}\infty}(K,\Sigma)$ it says just that 
\[(\forall i<\omega)(\nor[t_i]>m^{t_i}_{\dn}).\]
\end{remark}

\begin{proposition}
\label{implications}
Suppose $(K,\Sigma)$ is a creating pair for $\bH$. 
\begin{enumerate}
\item Assume that $K$ is full. Then $(K,\Sigma)$ is nice and if  $\C(\nor)$ is
one of the conditions $({\rm s}\infty)$, $(\infty)$, or $(f)$ (where $f$ is
any fast function), then the forcing notion $\q^*_{\C(\nor)}(K,\Sigma)$ is a
dense subset of $\q_{\C(\nor)}(K,\Sigma)$. 
\item If $(K,\Sigma)$ is nice, $(w,t_0,t_1,t_2,\ldots)$ is a condition in
$\q^*_\emptyset(K,\Sigma)$ and  $\langle s_n: n\in\omega\rangle$ is such that
for some $0=k_0<k_1<\ldots<\omega$, $s_n\in\Sigma(t_i: k_n\leq i<k_{n+1})$
(for all $n\in\omega$) then $(w,s_0,s_1,s_2,\ldots)$ is a condition in
$\q^*_\emptyset(K,\Sigma)$ (stronger than $(w,t_0,t_1,t_2,\ldots)$).
\item If $(K,\Sigma)$ is forgetful then it is full.
\end{enumerate}
\end{proposition}

\begin{definition}
\label{essapprox}
Let $(K,\Sigma)$ be a creating pair, $\C(\nor)$ be a norm condition,
$p\in\q^*_{\C(\nor)}(K,\Sigma)$ and $\dot{\tau}$ be a
$\q^*_{\C(\nor)}(K,\Sigma)$-name for an ordinal. We say that 
\begin{enumerate}
\item $p$ {\em essentially decides the name} $\dot{\tau}$ if
\[(\exists m\in\omega)(\forall u\in\pos(w^p,t^p_0,\ldots,t^p_{m-1}))((u,
t^p_m, t^p_{m+1},\ldots)\mbox{ decides the value of }\dot{\tau}),\]
\item $p$ {\em approximates $\dot{\tau}$ at $n$} (or at $t^p_n$) whenever:

\noindent for each $w_1\in\pos(w^p,t^p_0,\ldots,t^p_{n-1})$, if there is a
condition $r\in\q^*_{\C(\nor)}(K,\Sigma)$ stronger than $p$ and such that
$w^r=w_1$ and $r$ decides the value of $\dot{\tau}$ then the condition
$(w_1,t^p_n,t^p_{n+1},\ldots)$ decides the value of $\dot{\tau}$. 
\end{enumerate}
\end{definition}

\begin{lemma}
\label{smoessdec}
Suppose that $(K,\Sigma)$ is a smooth creating pair, $\C(\nor)$ is a norm
condition and $\dot{\tau}$ is a $\q^*_{\C(\nor)}(K,\Sigma)$--name for an
ordinal. Assume that a condition $p\in\q^*_{\C(\nor)}(K,\Sigma)$ essentially
decides $\dot{\tau}$ (approximates $\dot{\tau}$ at each $n$, respectively). 
Then each $q\geq p$ essentially decides $\dot{\tau}$ (approximates
$\dot{\tau}$ at each $n$, respectively).  
\end{lemma}

\begin{proof}
Immediate by smoothness. 
\end{proof}

\begin{definition} 
\label{orders}
Let $(K,\Sigma)$ be a creating pair for $\bH$.
\begin{enumerate}
\item For a property $\C(\nor)$ of $\omega$-sequences of creatures from $K$
and conditions $p,q\in\q_{\C(\nor)}^*(K,\Sigma)$ we define

$p\leq_{\apr} q$ (in $\q^*_{\C(\nor)}(K,\Sigma)$)\quad if and only if

$p\leq q$ and for some $k$ we have $(\forall i<\omega)(t^p_{i+k}=t^q_i)$ 

(so then $w^q\in\pos(w^p,t^p_0,\ldots, t^p_{k-1})$ too).

\item We define relations $\leq^{{\rm s}\infty}_n$ (for $n<\omega$) on
$\q^*_{{\rm s}\infty}(K,\Sigma)$ by:
\begin{enumerate}
\item[$(\alpha)$] $p\leq_{0}^{{\rm s}\infty} q$ (in $\q^*_{{\rm
s}\infty}(K,\Sigma)$) if $p\leq q$ and $w^p=w^q$,
\item[$(\beta)$] $p\leq_{n+1}^{{\rm s}\infty} q$ (in $\q^*_{{\rm
s}\infty}(K,\Sigma)$) if $p\leq_{0}^{{\rm s}\infty} q$ and $t^p_i=t^q_i$ for
$i<n+1$.  
\end{enumerate}

\item Relations $\leq^\infty_n$ on $\q^*_{\infty}(K,\Sigma)$ (for $n<\omega$)
are defined by:
\begin{enumerate}
\item[$(\alpha)$\ ] $p\leq_{0}^{\infty} q$ (in $\q^*_{\infty}(K,\Sigma)$) if
$p\leq q$ and $w^p=w^q$,
\item[$(\beta^*)$] $p\leq_{n+1}^{\infty} q$ (in $\q^*_{\infty}(K,\Sigma)$) if
$p\leq_{0}^{\infty} q$ and 
\[t^p_j=t^q_j\mbox{ for all }j\leq\min\{i<\omega: \nor[t^p_i]>n+1\}\mbox{
and}\] 
\[\{t^q_i: i<\omega\ \&\ \nor[t^q_i]\leq n+1\}\subseteq\{t^p_i: i<\omega\}.\]
\end{enumerate}
\item Relations $\leq^{{\rm w}\infty}_n$ on $\q^*_{{\rm w}\infty}(K,\Sigma)$
are defined by  
\begin{enumerate}
\item[$(\alpha)$\ ] $p\leq_{0}^{{\rm w}\infty} q$ (in $\q^*_{{\rm
w}\infty}(K,\Sigma)$) if $p\leq q$ and $w^p=w^q$,
\item[$(\beta^+)$] $p\leq_{n+1}^{{\rm w}\infty} q$ (in $\q^*_{{\rm
w}\infty}(K,\Sigma)$) if $p\leq_{0}^{{\rm w}\infty} q$ and 
\[t^p_j=t^q_j\mbox{ for all }j\leq\min\{i<\omega:\nor[t^p_i]>n+1\}.\]
\end{enumerate}

\item Let $f:\omega\times\omega\longrightarrow\omega$ be a fast function.
Relations $\leq^f_n$ on $\q^*_f(K,\Sigma)$ are defined by: 
\begin{enumerate}
\item[$(\alpha)$\ ] $p\leq_{0}^{f} q$ (in $\q^*_{\infty}(K,\Sigma)$) if
$p\leq q$ and $w^p=w^q$,
\item[$(\beta^f)$] $p\leq_{n+1}^{f} q$ (in $\q^*_{f}(K,\Sigma)$) if
$p\leq_{0}^f q$ and 
\[t^p_j=t^q_j\mbox{ for all }j\leq\min\{i<\omega: \nor[t^p_i]>
f(n+1,m^{t^p_i}_{\dn})\}\quad\quad\mbox{ and}\]  
\[\{t^q_i\!: i<\omega\ \&\ \nor[t^q_i]\leq f(n+1,m^{t^q_i}_{\dn})\}\subseteq
\{t^p_i\!: i<\omega\}.\]   
\end{enumerate}

\item We may omit superscripts in $\leq^{{\rm s}\infty}_n$, $\leq^\infty_n$,
$\leq^{{\rm w}\infty}_n$ and $\leq^f_n$ if it is clear from the context in
which forcing notion we are working (i.e.~what is the norm condition we deal
with).  
\end{enumerate}
\end{definition}

\begin{remark}
The difference between e.g.\ (3) and (4) is in the last condition of (3),
of course.  
\end{remark}

\begin{proposition}
\label{propord}
Suppose $(K,\Sigma)$ is a creating pair for $\bH$.  Let $\C(\nor)$ be one of
the following properties of $\omega$-sequences: $\C^{{\rm s}\infty}(\nor)$,
$\C^{\infty}(\nor)$, $\C^{{\rm w}\infty}(\nor)$ or $\C^f(\nor)$ for some fast
function $f$ (see \ref{conditions}) and let $\leq_n$ be the corresponding
relations (defined in \ref{orders}). Then 
\begin{enumerate}
\item $\leq_{\apr}$ is a partial order (stronger than $\leq$) on
$\q^*_{\C(\nor)}(K,\Sigma)$.  
\item $\leq_n$ (with superscripts) are partial orders
(stronger than $\leq$) on the respective $\q^*_{\C(\nor)}(K,\Sigma)$ and
$p\leq_{n+1}q$ implies $p\leq_n q$.  
\item Suppose that $p_n\in\q^*_{\C(\nor)}(K,\Sigma)$ (for $n\in\omega$) are
such that 
\[(\forall n\in\omega)(p_n\leq_{n+1} p_{n+1}).\]
Then the naturally defined limit condition $p=\lim\limits_n p_n$ satisfies:
\[p\in\q^*_{\C(\nor)}(K,\Sigma)\quad\mbox{ and }\quad(\forall n<\omega)(p_n
\leq_{n+1}p).\] 
\end{enumerate}
\end{proposition}

\begin{remark}
A natural property one could ask for in the context of creating pairs is some
kind of monotonicity:
\[\basis(t)=\dom(\val[t])\quad\mbox{ and }\quad \pos(u,t)=\{v: \langle
u,v\rangle\in \val[t]\},\]
for $t\in K$ and $u\in\basis(t)$. However, there is no real need for it, as
all our demands and assumptions on creating pairs will refer to $\pos$ (and
not $\val$). But for tree--creating pairs we will postulate the respective
demand, mainly to simplify notation (and have explicit tree--representations
of conditions), see \ref{treecreature}(3).
\end{remark}

\section{Tree creatures and tree--like forcing notions}
Here we introduce the second option for our general scheme: forcing notions in
which conditions are trees with norms. This case, though parallel to the one
of creating pairs, is of different character and therefore we reformulate all
general definitions for this particular context.

\begin{definition}
\label{quasitree}
\begin{enumerate}
\item {\em A quasi tree} is a set $T$ of finite sequences with the
$\vartriangleleft$-smallest element denoted by $\mrot(T)$. 

\item A quasi tree  $T$ is {\em a tree} if it is closed under initial
segments. If $T$ is a quasi tree then $\dcl(T)$ is the smallest tree
containing $T$ (the downward closure of $T$).

\item For a quasi tree $T$ and $\eta\in T$ we define {\em the successors of
$\eta$ in $T$}, {\em the restriction of $T$ to $\eta$}, {\em the splitting
points of $T$} and {\em maximal points of $T$} by:
\[\suc_T(\eta)=\{\nu\in T: \eta\vartriangleleft\nu\ \&\ \neg(\exists\rho\in
T)(\eta\vartriangleleft\rho\vartriangleleft\nu)\},\]
\[T^{[\eta]}=\{\nu\in T: \eta\trianglelefteq \nu\},\]
\[\spliting(T)=\{\eta\in T: |\suc_T(\eta)|\geq 2\},\]
\[\max(T)=\{\nu\in T:\mbox{ there is no }\rho\in T\mbox{ such that }
\nu\vartriangleleft\rho\}.\]   
We put $\hat{T}= T\setminus\max(T)$.

\item The set of all limit infinite branches through a quasi tree $T$ is
\[\lim(T)\stackrel{\rm def}{=}\{\eta:\eta\mbox{ is an $\omega$--sequence }\
\mbox{ and }\ (\exists^\infty n)(\eta\rest n\in T)\}.\]  
The quasi tree $T$ is {\em well founded} if $\lim(T)=\emptyset$.

\item A subset $F$ of a quasi tree $T$ is {\em a front} in $T$ if no two
distinct members of $F$ are $\vartriangleleft$-comparable and 
\[(\forall\eta\in\lim(T)\cup\max(T))(\exists n\in\omega)(\eta\restriction n\in
F).\] 
\end{enumerate}
\end{definition}

\begin{remark}
Note the difference between $\lim(T)$ and $\lim(\dcl(T))$ for a quasi tree
$T$. In particular, it is possible that a quasi tree $T$ is well-founded but
there is an infinite branch through $\dcl(T)$. Moreover, a front in $T$ does
not have to be a front in $\dcl(T)$.
\end{remark}

\begin{definition}
\label{treecreature}
\begin{enumerate}
\item A weak creature $t\in\WCR[\bH]$ is a {\em tree--creature} if $\dom(\val[
t])$ is a singleton $\{\eta\}$ and no two distinct elements of $\rng(\val[t])$
are $\vartriangleleft$--comparable;  

$\TCR[\bH]$ is the family of all tree--creatures for $\bH$. 

\item $\TCR_\eta[\bH]=\{t\in\TCR[\bH]: \dom(\val[t])=\{\eta\}\}$.

\item A sub-composition operation $\Sigma$ on $K\subseteq\TCR[\bH]$ is {\em a
tree composition} (and then $(K,\Sigma)$ is called {\em a tree--creating pair}
(for $\bH$)) if: 
\begin{enumerate}
\item[(a)] if $\cS\in [K]^{\textstyle{\leq}\omega}$, $\Sigma(\cS)\neq\emptyset$
then $\cS=\{s_\nu:\nu\in\hat{T}\}$ for some well founded quasi tree
$T\subseteq\bigcup\limits_{n<\omega}\prod\limits_{i<n}\bH(i)$ and a system
$\langle s_\nu: \nu\in\hat{T}\rangle\subseteq K$ such that for each finite
sequence $\nu\in\hat{T}$ 
\[s_\nu\in\TCR_\nu[\bH]\quad\mbox{ and }\quad\rng(\val[s_\nu])=\suc_T(\nu),\]
and
\item[(b)] if $t\in\Sigma(s_\nu: \nu\in\hat{T})$ then
$t\in\TCR_{\mrot(T)}[\bH]$ and $\rng(\val[t])\subseteq\max(T)$.
\end{enumerate}
If $\hat{T}=\{\mrot(T)\}$, $t=t_{\mrot(T)}\in\TCR_{\mrot(T)}[\bH]$ and
$\rng(\val[t])=\max(T)$ then we will write $\Sigma(t)$ instead of
$\Sigma(t_\nu:\nu\in\hat{T})$. 

\item A tree-composition $\Sigma$ on $K$ is {\em bounded} if for each $t\in
\Sigma(s_\nu:\nu\in\hat{T})$ we have
\[\nor[t]\leq\max\{\nor[s_\nu]: (\exists \eta\in\rng(\val[t]))(\nu
\vartriangleleft\eta)\}.\] 
\end{enumerate}
\end{definition}

\begin{remark}
\label{treerem}
1)\ \ \ Note that sets of tree creatures relevant for tree compositions have a
natural structure: we identify here $\cS$ with $\{s_{\nu(s)}: s\in\cS\}$ where
$\nu(s)$ is such that $s\in\TCR_{\nu(s)}$ and $s_{\nu(s)}=s$.\\
2)\ \ \ To check consistency of our notation for tree creatures with that of
\ref{maindef} note that in \ref{treecreature}(3), if $s_\nu\in
\Sigma(s_{\nu,\eta}:\eta\in \hat{T}_\nu)$ for each $\nu\in\hat{T}$, $T$ is a
well founded quasi tree as in (3){\bf (a)} of \ref{treecreature} then
$T^*\stackrel{\rm def}{=}\bigcup\limits_{\nu\in\hat{T}} T_\nu$ is a well
founded quasi tree, $\hat{T}^*=\bigcup\limits_{\nu\in\hat{T}} \hat{T}_\nu$ and
$\langle s_{\nu,\eta}: \nu\in\hat{T}, \eta\in\hat{T}_\nu\rangle$ is a system
for which $\Sigma$ may be non-empty, i.e.~it satisfies the requirements of
\ref{treecreature}(3){\bf (a)}.\\
3)\ \ \ Note that if $(K,\Sigma)$ is a tree--creating pair for $\bH$, $t\in
\TCR_\eta[\bH]$ then $\basis(t)=\{\eta\}$ and $\pos(\eta,t)=\rng(\val[t])$ (see
\ref{basis}). For this reason we will write $\pos(t)$ for $\pos(\eta,t)$ and
$\rng(\val[t])$ in the context of tree--creating pairs.\\ 
4)\ \ \ Tree--creating pairs have the properties corresponding to the niceness
and smoothness of creating pairs (see \ref{niceandsmo}, compare with
\ref{treesmooth}). 
\end{remark}

When dealing with tree--creating pairs it seems to be more natural to consider
both very special norm conditions and some restrictions on conditions of the
forcing notions we consider. The second is not very serious: the forcing
notions $\q^{\tree}_e(K,\Sigma)$ (for $e<5$) introduced in \ref{treeforcing}
below are dense subsets of the general forcing notions
$\q_{\C(\nor)}(K,\Sigma)$ (for suitable conditions $\C(\nor)$). We write the
definition of $\q^{\tree}_e(K,\Sigma)$ fully, not referring the reader to
\ref{maindef}, to show explicitly the way tree creating pairs work.

\begin{definition} 
\label{treeforcing} 
Let $(K,\Sigma)$ be a tree--creating pair for $\bH$.
\begin{enumerate}
\item We define the forcing notion $\q^{\tree}_1(K,\Sigma)$ by letting:

\noindent {\bf conditions } be sequences $p=\langle t_\eta: \eta\in T\rangle$
such that 
\begin{enumerate}
\item[(a)] $T\subseteq\bigcup\limits_{n\in\omega}\prod\limits_{i<n}\bH(i)$ is
a non-empty quasi tree with $\max(T)=\emptyset$,
\item[(b)] $t_\eta\in\TCR_\eta[\bH]\cap K$ and $\pos(t_\eta)=\suc_T(\eta)$
(see \ref{treerem}(3)),
\item[(c)${}_1$] for every $\eta\in\lim(T)$ we have:
\[\mbox{the sequence }\langle\nor[t_{\eta \rest k}]:k<\omega, \eta\rest k\in
T\rangle\mbox{ diverges to infinity;}\]
\end{enumerate}
\noindent{\bf the order} be given by:

\noindent $\langle t^1_\eta: \eta\in T^1\rangle\leq\langle t^2_\eta: \eta\in
T^2\rangle$ if and only if 

\noindent $T^2\subseteq T^1$ and for each $\eta\in T^2$ there is a well
founded quasi tree $T_{0,\eta}\subseteq (T^1)^{[\eta]}$ such that
$t^2_\eta\in\Sigma (t^1_\nu: \nu\in \hat{T}_{0,\eta})$.

If $p=\langle t_\eta:\eta\in T\rangle$ then we write $\mrot(p)=\mrot(T)$,
$T^p= T$, $t^p_\eta = t_\eta$ etc.
  
\item Similarly we define forcing notions $\q^{\tree}_e(K,\Sigma)$ for
$e=0,2,3,4$ replacing the condition {\bf (c)${}_1$} by {\bf
(c)${}_e$} respectively, where:
\begin{enumerate}
\item[(c)${}_0$]  for every $\eta\in\lim(T)$:
\[\lim\sup\langle\nor[t_{\eta \rest k}]:k<\omega, \eta\rest k\in T\rangle
=\infty,\] 
\item[(c)${}_2$]  for every $\eta\in T$ and $n<\omega$ there is $\nu$ such
that  $\eta\vartriangleleft\nu\in T$ and $\nor[t_\nu]\geq n$,
\item[(c)${}_3$] for every $\eta\in T$ and $n<\omega$ there is $\nu$ such
that $\eta\vartriangleleft\nu\in T$ and 
\[(\forall\rho\in T)(\nu\vartriangleleft\rho\ \ \Rightarrow\ \
\nor[t_\rho]\geq n),\]  
\item[(c)${}_4$] for every $n<\omega$, the set 
\[\{\nu\in T: (\forall \rho\in T)(\nu\vartriangleleft\rho\quad\Rightarrow
\quad\nor[t_\rho]\geq n)\}\]
contains a front of the quasi tree $T$.
\end{enumerate}
\item If $p\in\q^{\tree}_e(K,\Sigma)$ then we let $p^{[\eta]}=\langle t^p_\nu:
\nu\in (T^p)^{[\eta]}\rangle$ for $\eta\in T^p$.
\item For the sake of notational convenience we define partial order
$\q^{\tree}_{\emptyset}(K,\Sigma)$ in the same manner as
$\q^{\tree}_e(K,\Sigma)$ above but we omit the requirement {\bf (c)$_e$} (like
in \ref{conditions}; so this is essentially $\q_\emptyset(K,\Sigma)$).
\end{enumerate}
\end{definition}

\begin{remark}
1)\ \ \ In the definition above we do not follow exactly the notation of
\ref{maindef}: we omit the first part $w^p$ of a condition $p$ as it can be
clearly read from the rest of the condition. Of course the missing item is
$\mrot(p)$. In this new notation the name $\dot{W}$ of \ref{thereal} may be
defined by 
\[\forces_{\q^{\tree}_e(K,\Sigma)}\Wtil=\bigcup\{\mrot(p): p\in
\Gamma_{\q^{\tree}_e(K,\Sigma)}\}.\] 
2)\ \ \ Note that 
\[\begin{array}{rr}
\q^{\tree}_4(K,\Sigma)\subseteq\q^{\tree}_1(K,\Sigma)\subseteq \q^{\tree}_0
(K,\Sigma)\subseteq\q^{\tree}_2(K,\Sigma)&\quad\mbox{ and}\\
\q^{\tree}_1(K,\Sigma)\subseteq\q^{\tree}_3(K,\Sigma)\subseteq
\q^{\tree}_2(K,\Sigma)\\ 
  \end{array}\]
but in general these inclusions do not mean ``complete suborders''. If the
tree--creating pair is t-omittory (see \ref{tomit}) then $\q^{\tree}_4(K,
\Sigma)$ is dense in $\q^{\tree}_2(K,\Sigma)$ and thus all these forcing
notions are equivalent. If $(K,\Sigma)$ is $\bar{2}$--big (see \ref{kbig})
then $\q^{\tree}_4(K,\Sigma)$ is dense in $\q^{\tree}_1(K,\Sigma)$ (see
\ref{fronor2}). 
\end{remark}
Let us give two simple examples of tree--creating pairs.

Let $\bH(i)=\omega$ (for $i\in\omega$).

Let $K_0\subseteq\TCR[\bH]$ consists of these tree--creatures $s$ that if
$s\in\TCR_\eta[\bH]$ then $\rng(\val[s])\subseteq\{\eta\conc\langle k\rangle:
k\in\omega\}$ and
\[\nor[s]=\left\{
\begin{array}{ll}
\lh(\eta)&\mbox{if }\val[s]\mbox{ is infinite,}\\
0        &\mbox{otherwise.}
\end{array}\right.\]
The operation $\Sigma_0$ gives non-empty values for singletons only; for $s\in
K_0$ we let $\Sigma_0(s)=\{t\in K_0:\val[t]\subseteq\val[s]\}$ (an operation
$\Sigma$ defined in this manner will be further called {\em trivial\/}). 
Clearly $(K_0,\Sigma_0)$ is a tree--creating pair. Note that:
\begin{enumerate}
\item[(a)] the forcing notions $\q^{\tree}_2(K_0,\Sigma_0)$ and
$\q^{\tree}_0(K_0,\Sigma_0)$ are equivalent to Miller's Rational Perfect Set
Forcing; 
\item[(b)] the forcing notions $\q^{\tree}_1(K_0,\Sigma_0)$,
$\q^{\tree}_3(K_0,\Sigma_0)$, $\q^{\tree}_4(K_0,\Sigma_0)$ are equivalent to
the Laver forcing
\end{enumerate}
(thus $\q^{\tree}_1(K_0,\Sigma_0)$ is not a complete suborder of
$\q^{\tree}_0(K_0,\Sigma_0)$, and $\q^{\tree}_3(K_0,\Sigma_0)$ is not a
complete suborder of $\q^{\tree}_2(K_0,\Sigma_0)$).

Let us modify the norms on the tree--creatures a little. For this we define a
function $f:\fseo\longrightarrow\omega$ by
\[f(\langle\rangle)=0,\quad f(\eta\conc\langle k\rangle)=\left\{
\begin{array}{ll}
f(\eta)+1&\mbox{if }k=0\\
f(\eta)  &\mbox{otherwise.}
\end{array}\right.\]
Now, let $K_1$ consist of tree creatures $s\in\TCR[\bH]$ such that
$\rng(\val[s])\subseteq\{\eta\conc\langle k\rangle: k\in\omega\}$ (where $s\in
\TCR_\eta[\bH]$) and  
\[\nor[s]=\left\{
\begin{array}{ll}
f(\eta)&\mbox{if }\val[s]\mbox{ is infinite,}\\
0      &\mbox{otherwise.}
\end{array}\right.\]
Let $\Sigma_1$ be the trivial tree-composition on $K_1$, so it is nonempty for
singletons only and then $\Sigma_1(s)=\{t\in K_1: \val[t]\subseteq\val[s]\}$
(clearly $(K_1,\Sigma_1)$ is a tree creating pair). Then 
\begin{enumerate}
\item[(a)] the forcing notion $\q^{\tree}_0(K_1,\Sigma_1)$ is a dense suborder
of $\q^{\tree}_2(K_1,\Sigma_1)$, 
\item[(b)] the partial orders $\q^{\tree}_1(K_1,\Sigma_1)$ and
$\q^{\tree}_4(K_1,\Sigma_1)$ are empty, but
\item[(c)] $\q^{\tree}_3(K_1,\Sigma_1)$ is not--trivial (it adds a new real)
and it is not a complete suborder of $\q^{\tree}_2(K_1,\Sigma_1)$ (e.g.
incompatibility is not preserved) and it is disjoint from
$\q^{\tree}_0(K_1,\Sigma_1)$.
\end{enumerate}

\begin{definition}
\label{thick}
Let $(K,\Sigma)$ be a tree creating pair, $e<5$, $p\in\q^{\tree}_e(K,\Sigma)$.
A set $A\subseteq T^p$ is called {\em an $e$-thick antichain} (or just {\em a
thick antichain}) if it is an antichain in $(T^p,\vartriangleleft)$ and for
every condition $q\in\q^{\tree}_e(K,\Sigma)$ stronger than $p$ the
intersection $A\cap\dcl(T^q)$ is non-empty. 
\end{definition}

\begin{proposition}
\label{frontetc}
Suppose that $(K,\Sigma)$ is a tree--creating pair for $\bH$, $e<5$,
$p\in\q^{\tree}_e(K,\Sigma)$ and $\eta\in T^p$. Then:
\begin{enumerate}
\item $\q_e^{\tree}(K,\Sigma)$ is a partial order.
\item Each $e$-thick antichain in $T^p$ is a maximal antichain. Every front of
$T^p$ is an $e$-thick antichain in $T^p$.
\item If $e\in\{1,3,4\}$, $n<\omega$ then the set 
\[\begin{array}{ll}
B_n(p)\stackrel{\rm def}{=}\{\eta\in T^p: & \mbox{(i) } \ (\forall \nu\in
T^p)(\eta\trianglelefteq\nu\ \ \Rightarrow\ \ \nor[t_\nu]> n)\mbox{ but }\\
\ & \mbox{(ii) } \mbox{ no }\eta'\vartriangleleft\eta, \eta'\in T^p \mbox{
satisfies (i)}\}\\ 
\end{array}\]
is a maximal $\vartriangleleft$-antichain in $T^p$. If $e=4$ then $B_n(p)$ is
a front of $T^p$. 
\item For every $m,n<\omega$ the set
\[\begin{array}{ll}
F_n^m(p)\stackrel{\rm def}{=}\{\eta\in T^p: & \mbox{(i) } \ \nor[t_\eta]>n
\mbox{ and} \\   
\ & \mbox{(ii)}\ |\{\eta'\in T^p: \eta'\vartriangleleft\eta\ \&\
\nor[t_{\eta'}]>n\}|=m\}\\ 
\end{array}\]
is a maximal $\vartriangleleft$-antichain of $T^p$. If $e\in\{0,1,4\}$ then
$F^m_n(p)$ is a front of $T^p$.  
\item If $K$ is finitary (so $|\val[t]|<\omega$ for $t\in K$, see
\ref{morecreat}) then every front of $T^p$ is a front of $\dcl(T^p)$ and hence
it is finite.  
\item If $\Sigma$ is bounded then each $F^m_n(p)$ is a thick antichain of
$T^p$.  
\item $p\leq p^{[\eta]}\in \q_e^{\tree}(K,\Sigma)$\quad and\quad $\mrot(p^{[
\eta]})=\eta$. 
\end{enumerate}
\end{proposition}

\begin{remark}
\label{treesmooth}
One of the useful properties of tree--creating pairs $(K,\Sigma)$ and forcing
notions $\q^{\tree}_e(K,\Sigma)$ is the following:
\begin{enumerate}
\item[$(*)_{\ref{treesmooth}}$] {\em Suppose that $p,q\in\q^{\tree}_e(K,
\Sigma)$, $p\leq q$ (so in particular $T^q\subseteq T^p$), $\eta\in T^q$ and
$\nu\vartriangleleft\eta$, $\nu\in T^p$.\\
Then $p^{[\nu]}\leq q^{[\eta]}$.}
\end{enumerate}
\end{remark}

\begin{definition}
\label{AxA}
Let $p,q\in\q^{\tree}_e(K,\Sigma)$, $e<3$ (and $(K,\Sigma)$ a tree--creating
pair). We define relations $\leq^e_n$ for $n\in\omega$ by:
\begin{enumerate}
\item If $e\in\{0,2\}$ then:

\noindent $p\leq^e_0 q$ (in $\q^{\tree}_e(K,\Sigma)$) if $p\leq q$ and
$\mrot(p)=\mrot(q)$,

\noindent $p\leq^e_{n+1} q$ (in $\q^{\tree}_e(K,\Sigma)$) if $p\leq^e_0 q$ and

if $\eta\in F^0_n(p)$ (see \ref{frontetc}(4)) and $\nu\in T^p$,
$\nu\trianglelefteq\eta$ then $\nu\in T^q$ and $t^q_\nu=t^p_\nu$. 

\item The relations $\leq^1_{n}$ (on $\q^{\tree}_1(K,\Sigma)$) are defined by:

\noindent  $p\leq^1_0 q$ (in $\q^{\tree}_1(K,\Sigma)$) if and only if $p\leq
q$ and $\mrot(p)=\mrot(q)$,

\noindent $p\leq_{n+1}^1 q$ (in $\q^{\tree}_1(K,\Sigma)$) if $p\leq^1_0 q$ and

if $\eta\in F^0_n(p)$ (see \ref{frontetc}(4)) and $\nu\in T^p$,
$\nu\trianglelefteq\eta$ then $\nu\in T^q$, $t^p_\nu=t^q_\nu$, and
\[\{t^q_\eta:\eta\in T^q\ \&\ \nor[t^q_\eta]\leq n\}\subseteq\{t^p_\eta:\eta
\in T^p\}.\]  
\item We may omit the superscript $e$ in $\leq^e_n$ if it is clear in which
of the forcing notions $\q^{\tree}_e(K,\Sigma)$ we are working.
\end{enumerate}
\end{definition}

\begin{proposition}
\label{fusAxA}
Let $(K,\Sigma)$ be a tree--creating pair for $\bH$, $e<3$.
\begin{enumerate}
\item The relations $\leq^e_n$ are partial orders on $\q^{\tree}_e(K,\Sigma)$
stronger than $\leq$. The partial order $\leq^e_{n+1}$ is stronger than
$\leq^e_n$. 
\item Suppose that conditions $p_n\in\q^{\tree}_e(K,\Sigma)$ are such that
$p_n\leq_{n+1}^e p_{n+1}$. 

Then $\lim\limits_{n\in\omega}p_n\in\q^{\tree}_e(K,\Sigma)$ and $(\forall 
n\in\omega)(p_n\leq^e_{n+1}\lim\limits_{n\in\omega} p_n)$ (where the limit
condition $p=\lim\limits_{n\in\omega}p_n$ is defined naturally; $T^p
=\bigcap\limits_{n\in\omega} T^{p_n}$). 
\end{enumerate}
\end{proposition}

\section{Non proper examples}
In the next chapter we will see that if one combines a norm condition with
suitable properties of a weak creating pair then the resulting forcing notion
is proper. In particular we will see that (with the norm conditions defined in
\ref{treeforcing}) getting properness in the case of tree--creating pairs is
relatively easy. Here, however, we show that one cannot expect a general
theorem like ``$\q^{\tree}_{\C(\nor)}(K,\Sigma)$ is always proper'' and that
we should be always careful a little bit. The forcing notions resulting from
our general schema may collapse $\aleph_1$! For example, looking at the norm
conditions introduced in \ref{treeforcing}(2) one could try to consider the
following condition 
\begin{enumerate}
\item[${\bf (c)}_5$]  $(\forall k\in\omega)(\forall^\infty n)(\forall\eta\in
T^p)(\lh(\eta)\geq n\ \ \Rightarrow\ \ \nor[t_\nu]\geq k)$.
\end{enumerate}
If a creating pair $(K,\Sigma)$ is finitary then, clearly, the forcing
notions $\q^{\tree}_5(K,\Sigma)$ and $\q^{\tree}_4(K,\Sigma)$ are the same.

The forcing notion $\q^{\tree}_5(K,\Sigma)$ might be even not proper. The
following example shows this bad phenomenon which may be made quite general.

\begin{example}
\label{nonpro1}
Let $\bH(i)=\omega$ for $i\in\omega$.

\noindent There is a tree creating pair
$(K_{\ref{nonpro1}},\Sigma_{\ref{nonpro1}})$ for $\bH$ which is simple (see
\ref{simpglui}) and the forcing notion $\q^{\tree}_5(K_{\ref{nonpro1}},
\Sigma_{\ref{nonpro1}})$ is not proper (and collapses $\con$ onto $\omega$).
\end{example}

\begin{proof}[Construction]
To define the family $K_{\ref{nonpro1}}$ of tree creatures for $\bH$
choose families $\emptyset\neq S^\ell_\eta\subseteq R_\eta\subseteq\iso$ and
functions $h_\eta: R_\eta\longrightarrow \omega$ (for $\eta\in\fseo$,
$0<\ell\leq\lh(\eta)$) such that for every $\eta,\ell$ we have:
\begin{enumerate}
\item[$(\alpha)$] $\omega\in R_\eta$, $h_\eta(\omega)=\lh(\eta)+1$,
\item[$(\beta)$]  if $F\in S^\ell_\eta$ then $h_\eta(F)=\ell$, each
$S^\ell_\eta$ is infinite,
\item[$(\gamma)$] if $F_0,F_1\in S^\ell_\eta$, $F_0\neq F_1$ then $F_0\cap
F_1=\emptyset$,
\item[$(\delta)$] if $A\in R_\eta$, $h_\eta(A)\geq\ell+1$ then for each $F\in
S^\ell_\eta$ 
\[A\cap F\in R_\eta\quad\mbox{ and }\quad h_\eta(A\cap F)=\ell,\]
\item[$(\varepsilon)$] if $A_0,A_1\in R_\eta$, $A_0\subseteq A_1$ then
$h_\eta(A_0)\leq h_\eta(A_1)$.
\end{enumerate}
\medskip

\noindent There are several possibilities to construct $S^\ell_\eta,R_\eta,
h_\eta$ as above. One can do it for example in the following way. Fix $\eta
\in\omega^{\textstyle n}$. Take a system $\{K_\sigma:\sigma\in \omega^{
\textstyle {\leq}n}\}$ of infinite subsets of $\omega$ such that
\begin{enumerate}
\item[a)] $K_{\langle\rangle}=\omega$,
\item[b)] $\sigma_0\vartriangleleft\sigma_1\in\omega^{\textstyle{\leq}n}\ \
\Rightarrow\ \ K_{\sigma_1}\subseteq K_{\sigma_0}$,
\item[c)] $\sigma_0,\sigma_1\in \omega^{\textstyle \ell}\ \&\ \ell\leq n\ \&\
\sigma_0\neq\sigma_1\ \ \Rightarrow\ \ K_{\sigma_0}\cap K_{\sigma_1}=
\emptyset$, 
\item[d)] $\sigma\in\omega^{\textstyle{<}n}\ \ \Rightarrow\ \ K_\sigma=
\bigcup\limits_{m\in\omega}K_{\sigma\conc\langle m\rangle}$. 
\end{enumerate}
Now put $R_\eta=\{\bigcup\limits_{\sigma\in I}K_\sigma:\emptyset\neq I
\subseteq\omega^{\textstyle n}\}$. For $A\in R_\eta$ we declare that
\begin{enumerate}
\item[e)] $h_\eta(A)\geq 1$, $h_\eta(\omega)=h_\eta(K_{\langle\rangle})=n+1$,
\item[f)] if for some $\sigma\in\omega^{\textstyle n-\ell}$, $\ell\leq n$ the
set $A$ contains the set $K_\sigma$ then $h_\eta(A)\geq\ell+1$,
\item[g)] if $A_0\subseteq A_1$, $A_0,A_1\in R_\eta$ then $h_\eta(A_0)\leq
h_\eta(A_1)$. 
\end{enumerate}
Next we put 
\begin{quotation}
\noindent $F^{\ell,m}_\eta=\bigcup\{K_{\sigma\conc\langle m\rangle}:
\sigma\in\omega^{\textstyle n-\ell}\}$, \ \ \ \ for $0<\ell\leq n$,
$m\in\omega$\quad and 

\noindent $S^\ell_\eta=\{F^{\ell,m}_\eta: m\in\omega\}$.
\end{quotation}
It should be easy to check that $S_\eta^\ell$, $R_\eta$, $h_\eta$ (for
$\eta\in\fseo$, $0<\ell\leq\lh(\eta)$) defined in this way satisfy the
requirements $(\alpha)$--$(\varepsilon)$.  
\medskip

A tree creature $t\in\TCR_\eta[\bH]$ is in $K_{\ref{nonpro1}}$ if $\eta\in
\fseo$, $\lh(\eta)>0$ and 
\[\dis[t]\in R_\eta,\quad\val[t]=\{\langle\eta,\eta\conc\langle m\rangle
\rangle: m\in \dis[t]\}\quad\mbox{ and }\quad \nor[t]=h_\eta(\dis[t]).\]
If $A\in R_\eta$ then the unique tree--creature $t\in\TCR_\eta[\bH]\cap
K_{\ref{nonpro1}}$ such that $\dis[t]=A$ will be denoted by $t^{\eta,A}$. 

\noindent The operation $\Sigma_{\ref{nonpro1}}$ is trivial: it gives a
non-empty result for singletons only and then $\Sigma_{\ref{nonpro1}}(t)=\{s
\in K_{\ref{nonpro1}}: \val[s]\subseteq\val[t]\}$.  

\begin{claim}
\label{cl6}
The forcing notion $\q^{\tree}_5(K_{\ref{nonpro1}},\Sigma_{\ref{nonpro1}})$
collapses $\con$ onto $\omega$. 
\end{claim}

\begin{proof}[{\em Proof of the claim}]
Fix $\eta\in\fseo\setminus\{\langle\rangle\}$ for a moment.

\noindent Elements of $S^\ell_\eta$ are pairwise disjoint so we may naturally
order them according to the smallest element. Say $S^\ell_\eta=\{
F^{\ell,m}_\eta:m<\omega\}$. Let $f:[\lh(\eta),\omega)\longrightarrow\omega$.
We define a condition $p^{f,\eta}\in\q^{\tree}_5(K_{\ref{nonpro1}},
\Sigma_{\ref{nonpro1}})$ putting (we keep the notation as for the forcing
notions $\q^{\tree}_e(K,\Sigma)$):    
\begin{quotation}
$\mrot(p^{f,\eta})=\eta$;

let $k_0=\lh(\eta)$, $k_{\ell+1}=f(k_\ell)+k_\ell+1$ (for $\ell<\omega$);

if $\nu\in T^{p^{f,\eta}}$, $k_{\ell-1}\leq\lh(\nu)<k_\ell$ then
$t^{p^{f,\eta}}_\nu=t^{\nu,F^{\ell,f(\lh(\nu))}_\nu}$ and

$\suc_{T^{p^{f,\eta}}}(\nu)=\{\nu\conc\langle m\rangle: m\in
F^{\ell,f(\lh(\nu))}_\nu\}$. 
\end{quotation}
Clearly this defines $p^{f,\eta}\in\q^{\tree}_5(K_{\ref{nonpro1}},
\Sigma_{\ref{nonpro1}})$. Note that 

if $f,g:[\lh(\eta),\omega) \longrightarrow\omega$ are distinct  

then the conditions $p^{f,\eta}, p^{g,\eta}$ are incompatible in
$\q^{\tree}_5(K_{\ref{nonpro1}}, \Sigma_{\ref{nonpro1}})$ 

\noindent (by the requirement $(\gamma)$). Let $\dot{\tau}$ be a
$\q^{\tree}_5(K_{\ref{nonpro1}},\Sigma_{\ref{nonpro1}})$--name for a function
defined on $\fseo$ such that $\forces$``$\dot{\tau}(\eta)\in\V\ \&\ \dot{\tau}
(\eta):[\lh(\eta),\omega)\longrightarrow\omega$'' for $\eta\in\fseo$ and for
$f:[\lh(\eta),\omega) \longrightarrow\omega$ we have  
\[p^{f,\eta}\forces_{\q^{\tree}_5(K_{\ref{nonpro1}},\Sigma_{\ref{nonpro1}})}
\dot{\tau}(\eta)=f.\] 
This definition is correct as $\{p^{f,\eta}\!: f\!:[\lh(\eta),\omega)
\rightarrow\omega\}$ is an antichain (of course it is not necessarily maximal
in $\q^{\tree}_5(K_{\ref{nonpro1}},\Sigma_{\ref{nonpro1}})$). The claim will
be shown if we prove that  
\[\forces_{\q^{\tree}_5(K_{\ref{nonpro1}},\Sigma_{\ref{nonpro1}})} (\forall
g\in\baire\cap\V)(\exists\eta\in\fseo)(\forall
n\geq\lh(\eta))(\dot{\tau}(\eta)(2n)=g(n)).\] 
For this suppose that $g\in\baire$, $p\in\q^{\tree}_5(K_{\ref{nonpro1}},
\Sigma_{\ref{nonpro1}})$. Choose an increasing sequence $\lh(\mrot(p))<
k_0<k_1<\ldots$ of odd integers such that for each $\ell<\omega$
\[(\forall\nu\in T^p)(k_\ell\leq\lh(\nu)\ \ \Rightarrow\ \ \nor[t^p_\nu]\geq
\ell +2).\]
Let $f:[k_0,\omega)\longrightarrow\omega$ be a function such
that:
\begin{enumerate}
\item $f(k_\ell)=k_{\ell+1}-k_\ell-1$
\item if $n\geq k_0$ then $f(2n)=g(n)$.
\end{enumerate}
(Note that these clauses are compatible as the $k_\ell$'s are odd. Of course
there is still much freedom left in defining $f$.)

Choose $\eta\in T^p\cap\omega^{\textstyle k_0}$ and look at the condition
$p^{f,\eta}$. Due to the requirement $(\delta)$ this condition is compatible
with $p$: \ \ \ define
$r\in\q^{\tree}_5(K_{\ref{nonpro1}},\Sigma_{\ref{nonpro1}})$ by 
\begin{quotation}
$\mrot(r)=\eta$, $T^r\subseteq T^p$\ \ \ and\\
if $\nu\in T^r$, $k_{\ell-1}\leq\lh(\nu)<k_\ell$ then
\[t^r_\nu=t^{\nu,A_\nu}\quad\mbox{ and }\quad\suc_{T^r}(\nu)=\{\nu\conc\langle
m\rangle: m\in A_\nu\}\] 
where $A_\nu=\dis[t^p_\nu]\cap F^{\ell,f(\lh(\nu))}_\nu$.
\end{quotation}
Note that by the choice of $k_\ell$ and the requirement $(\delta)$ we have
\[t^{\nu,A_\nu}\in\Sigma_{\ref{nonpro1}}(t^p_\nu)\cap\Sigma_{\ref{nonpro1}}
(t^{p^{f,\eta}}_\nu)\quad\mbox{ and }\quad\nor[t^{\nu,A_\nu}]=\ell.\]
Consequently the definition of $r$ is correct. Clearly $r$ is stronger than
both $p$ and $p^{f,\eta}$. Thus 
\[r\forces_{\q^{\tree}_5(K_{\ref{nonpro1}},\Sigma_{\ref{nonpro1}})}(\forall
n\geq k_0)(\dot{\tau}(\eta)(n)=f(n))\]
which together with 
\[(\forall n\geq k_0)(f(2n)=g(n))\quad\mbox{ and }\quad k_0=\lh(\eta)\]
finishes the proof of the claim. 
\end{proof}
\end{proof}

Our next example shows that the assumption that $(K,\Sigma)$ is finitary in
\ref{omitoryproper} is crucial.

\begin{example}
\label{nonpro2}
Let $\bH(i)=\omega$ for $i\in\omega$.

\noindent There is a creating pair $(K_{\ref{nonpro2}},
\Sigma_{\ref{nonpro2}})$ for $\bH$ which is forgetful and growing (see
\ref{omitoryetc}(3)) but the forcing notion $\q^*_\infty(K_{\ref{nonpro2}},
\Sigma_{\ref{nonpro2}})$ is not proper (and collapses $\con$ onto $\omega$).

\noindent (By \ref{sinfty} we may replace $\q^*_\infty$ by either $\q^*_{{\rm
w}\infty}$ or $\q^*_{{\rm s}\infty}$ or $\q^*_f$ (for a fast function $f$).)
\end{example}

\begin{proof}[Construction]
This is similar to \ref{nonpro1}: for $0<\ell\leq i<\omega$ choose
$\emptyset\neq S^\ell_i\subseteq R_i\subseteq\iso$ and functions $h_i:R_i
\longrightarrow\omega$ satisfying the requirements ($\alpha$)--($\varepsilon$)
of the construction of \ref{nonpro1} (with $i$ instead of $\eta$ and
$\lh(\eta)$) and 
\begin{enumerate}
\item[$(\zeta)$] $\bigcup S^\ell_i=\omega$ for each $0<\ell\leq i<\omega$
\end{enumerate}
(this additional condition is satisfied by the example constructed there). Fix
an enumeration $S^\ell_i=\{F^{\ell,m}_i: m\in\omega\}$.  

A creature $t\in K_{\ref{nonpro2}}$ may be described in the following
way. For each $i\in [m^t_{\dn},m^t_{\up})$ we have a set $A_i\in R_i$. Now:
\[\begin{array}{l}
\dis[t]=\langle A_i:m^t_{\dn}\leq i<m^t_{\up}\rangle\\
\val[t]=\{\langle u,v\rangle\!\in\!\prod\limits_{i<m^t_{\dn}}\!\!\bH(i)\!
\times\!\prod\limits_{i<m^t_{\up}}\!\!\bH(i)\!: u\vartriangleleft v\ \&\
(\forall i\!\in\![m^t_{\dn},m^t_{\up}))(v(i)\in A_i)\},\\
\nor[t]=\max\{h_i(A_i): i\in [m^t_{\dn},m^t_{\up})\}.
\end{array}\] 
If creatures $t_0,\ldots,t_{n-1}\in K_{\ref{nonpro2}}$ are determined by sets
$A^j_i\in R_i$ (for $j<n$, $i\in [m^{t_j}_{\dn},m^{t_j}_{\up})$) in the way
described above and $m^{t_{j+1}}_{\dn}=m^{t_j}_{\up}$ (for $j<n-1$) then
$\Sigma_{\ref{nonpro2}}(t_0,\ldots,t_{n-1})$ consists of all creatures $t\in
K_{\ref{nonpro2}}$ which are determined (in the way described above) by some
sets $A_i\in R_i$ (for $i\in [m^{t_0}_{\dn},m^{t_{n-1}}_{\up})$) such that
$A_i\subseteq A^j_i$ whenever $m^{t_j}_{\dn}\leq i<m^{t_j}_{\up}$, $j<n$. 

It is easy to check that $\Sigma_{\ref{nonpro2}}$ is a composition operation
on $K_{\ref{nonpro2}}$. The creating pair
$(K_{\ref{nonpro2}},\Sigma_{\ref{nonpro2}})$ is forgetful and growing. 

\begin{claim}
\label{cl18}
The forcing notion $\q^*_\infty(K_{\ref{nonpro2}},
\Sigma_{\ref{nonpro2}})$ collapses $\con$ onto $\omega$.
\end{claim}

\begin{proof}[{\em Proof of the claim}]
We proceed like in \ref{cl6} (with small modifications however). Let
$\pi:\omega\longrightarrow\omega\times\omega:n\mapsto(\pi_0(n),\pi_1(n))$ be
a bijection. Let $\pi^*_n:\omega^{\textstyle [n,\omega)}\longrightarrow\baire$
(for $n\in\omega$) be mappings defined in the following  manner. Let
$f:[n,\omega)\longrightarrow\omega$; 
\begin{quotation}
\noindent inductively define $n_0=n$, $n_{\ell+1}=n_\ell+\pi_0(f(n_\ell))+
\pi_1(f(n_\ell))+2$ (for $\ell<\omega$) and then $m_\ell=n_\ell+
\pi_0(f(n_\ell))+1$ (so $n_\ell<m_\ell<n_{\ell+1}$);

\noindent now put $\pi^*_n(f)(\ell)=f(m_\ell)$ (for $\ell\in\omega$).
\end{quotation}
For $0<n<\omega$, $u\in\prod\limits_{i<n}\bH(i)$ and a function $f:[n,\omega)
\longrightarrow\omega$ define a condition $p^{u,f}\in\q^*_\infty
(K_{\ref{nonpro2}},\Sigma_{\ref{nonpro2}})$:
\begin{quotation}
\noindent $k_0=n$, $k_{2\ell+1}=k_{2\ell}+\pi_0(f(k_{2\ell}))+1$, $k_{2\ell+2}
=k_{2\ell+1}+\pi_1(f(k_{2\ell}))+1$,

\noindent $w^{p^{u,f}}=u$,

\noindent if $k_{2\ell}\leq i+n<k_{2\ell+2}$ then $t^{p^{u,f}}_i\in
K_{\ref{nonpro2}}$ is such that $m^{t^{p^{u,f}}_i}_{\dn}=n+i$,
$m^{t^{p^{u,f}}_i}_{\up}=n+i+1$, $\dis[t^{p^{u,f}}_i]=\langle F^{\ell+1,
f(i+n)}_{i+n}\rangle$.
\end{quotation}
As in \ref{cl6}, if $f,g:[n,\omega)\longrightarrow\omega$ are distinct, $u\in
\prod\limits_{i<n}\bH(i)$ then the conditions $p^{u,f},p^{u,g}$ are
incompatible. Consequently we may choose a $\q^*_\infty(K_{\ref{nonpro2}},
\Sigma_{\ref{nonpro2}})$--name $\dot{\tau}$ for a function on $\omega$ such
that $\forces(\forall n\in\omega)(\dot{\tau}(n):[n,\omega)\rightarrow\omega)$
and $p^{u,f}\forces\dot{\tau}(n)=f$. To finish it is enough to show that  
\[\forces_{\q^*_\infty(K_{\ref{nonpro2}},\Sigma_{\ref{nonpro2}})} (\forall
g\in\baire\cap\V)(\exists n\in\omega)(\pi^*_n(\dot{\tau}(n))=g).\]
Suppose that $p\in\q^*_\infty(K_{\ref{nonpro2}},\Sigma_{\ref{nonpro2}})$,
$g\in\baire$. Choose $2<i_0<i_1<\ldots<\omega$ such that $\nor[t^p_{i_\ell}]
\geq \ell+2$ and next choose $k_0<k_1<k_2<\ldots<\omega$ such that for each
$\ell\in\omega$: 
\begin{quotation}
\noindent $m^{t^p_{i_\ell}}_{\dn}\leq k_\ell<m^{t^p_{i_\ell}}_{\up}$ and for
some set $A_\ell\in R_{k_\ell}$ we have
\[h_{k_\ell}(A_\ell)\geq \ell+2\quad\mbox{ and }\quad (\forall n\in A_\ell)
(\exists \langle u,v\rangle\in \val[t^p_{i_\ell}])(v(k_\ell)=n)\]
\end{quotation}
(possible by the way we defined $(K_{\ref{nonpro2}},\Sigma_{\ref{nonpro2}})$).
Choose any $v\in\pos(w^p,t^p_0,\ldots,t^p_{i_0})$ and let $u=v\rest k_0$.
Next choose $f:[k_0,\omega)\longrightarrow\omega$ such that for each
$\ell\in\omega$:
\[\pi_0(f(k_{2\ell}))=k_{2\ell+1}-k_{2\ell}-1,\quad \pi_1(f(k_{2\ell}))=
k_{2\ell+2}-k_{2\ell+1}-1,\quad f(k_{2\ell+1})=g(\ell),\]
and if $k\in (k_{2\ell},k_{2\ell+2})\setminus\{k_{2\ell+1}\}$,
$m^{t^p_i}_{\dn}\leq k<m^{t^p_i}_{\up}$ ($\ell<\omega$, $i<\omega$) then 
\[(\exists \langle u,v\rangle\in\val[t^p_i])(v(k)\in F^{\ell+1,f(k)}_k).\]
One easily checks that the choice of $f$ is possible (remember the additional
requirement $(\zeta)$) and that the conditions $p^{u,f}$ and $p$ are
compatible in $\q^*_\infty(K_{\ref{nonpro2}},\Sigma_{\ref{nonpro2}})$. As 
\[p^{u,f}\forces_{\q^*_\infty(K_{\ref{nonpro2}},\Sigma_{\ref{nonpro2}})}
\pi^*_n(\dot{\tau}(n))=\pi^*_n(f)=g,\]
we finish the proof of the claim. 
\end{proof}
\end{proof}

One could expect that the main reason for collapsing $\con$ in the two
examples constructed above is that the $(K,\Sigma)$'s there are not
finitary. But this is not the case. Using similar ideas we may build a
finitary creating pair $(K,\Sigma)$ for which the forcing notion
$\q^*_\infty(K,\Sigma)$ collapses $\con$ onto $\omega$ as well. This is the
reason why we have to use forcing notions $\q^*_f(K,\Sigma)$ with $(K,\Sigma)$
satisfying extra demands (including Halving and bigness, see
\ref{halbigprop}) and why $\q^*_\infty(K,\Sigma)$ is used only for growing
$(K,\Sigma)$ (so then $\q^*_{{\rm s}\infty}(K,\Sigma)$ is dense in
$\q^*_\infty(K,\Sigma)$, see \ref{2.1}). This bad effect can be made quite
general and we will present it in this way, trying to show the heart of the
matter. One could try to cover the previous examples by our ``negative
theory'' too, but this would involve much more complications.   

\begin{definition}
\label{local}
We say that a weak creating pair $(K,\Sigma)$ is {\em local} if for every
$t\in K$, $w\in\basis(t)$ and $u\in\pos(w,t)$ we have $\lh(u)=\lh(w)+1$.  
\end{definition}

\begin{definition}
\label{bad}
Let $(K,\Sigma)$ be a (nice and smooth) creating pair for $\bH$ which is local
(so $t\in K\ \Rightarrow\ m^t_{\up}=m^t_{\dn}+1$) and simple (which means that
$\Sigma(\cS)\neq\emptyset\ \Rightarrow\ |\cS|=1$; see \ref{simpglui}(1)).\\
We say that $(K,\Sigma)$ is {\em definitely bad} if there are a perfect tree
$T\subseteq\fseo$ and mappings $F_0,F_1$ such that
\begin{enumerate}
\item $T\cap\omega^{\textstyle m}$ is finite for each $m\in\omega$,\quad
$\dom(F_0)=\dom(F_1)=\bigcup\limits_{m<\omega}\prod\limits_{i<m}\bH(i)$, 
\item if $v\in\prod\limits_{i\leq m}\bH(i)$, $m<\omega$ then $F_1(v): T\cap
\omega^{\textstyle m}\longrightarrow 2$ and $F_0(v):T\cap\omega^{\textstyle m}
\longrightarrow T$ is such that $(\forall\eta\in T\cap\omega^{\textstyle m})(
F_0(v)(\eta)\in\suc_T(\eta))$,
\item {\em if\/} $t\in K$, $\nor[t]\geq 2$, $m=m^t_{\dn}>0$, $i<2$ and
$F^*:T\cap\omega^{\textstyle m}\longrightarrow T$ is such that $F^*(\nu)\in
\suc_T(\nu)$ for $\nu\in T\cap\omega^{\textstyle m}$ {\em then\/} there is
$s\in\Sigma(t)$ such that $\nor[s]\geq\nor[t]-1$ and for each $\eta\in T\cap
\omega^{\textstyle m-1}$ there is $\nu\in\suc_T(\eta)$ with 
\[(\forall u\in\basis(s))(\forall v\in\pos(u,s))(F_0(v)(\nu)=F^*(\nu)\ \ \&\ \
F_1(v)(\nu)=i).\]
\end{enumerate}
\end{definition}

\begin{proposition}
Suppose that $(K,\Sigma)$ is a local, simple and definitely bad creating pair
for $\bH$ such that $\Sigma(t)$ is finite for each $t\in K$. Then
\begin{enumerate}
\item[(a)] the forcing notion $\q^*_\infty(K,\Sigma)$ collapses $\con$ onto
$\omega$,
\item[(b)] if $f:\omega\times\omega\longrightarrow\omega$ is a fast function
then the forcing notion $\q^*_f(K,\Sigma)$ collapses $\con$ onto $\omega$. 
\end{enumerate}
\end{proposition}

\begin{proof}
In both cases the proof is exactly the same, so let us deal with {\bf (a)}
only. So suppose that a finitely branching perfect tree $T\subseteq\fseo$ and
functions $F_0,F_1$ witness that $(K,\Sigma)$ is definitely bad. Let $G^0:
T\times\prod\limits_{i<\omega}\bH(i)\longrightarrow [T]$ and $G^1:[T]\times
\prod\limits_{i<\omega}\bH(i)\longrightarrow\can$ be defined by
\[\begin{array}{l}
G^0(\eta,W)\rest\lh(\eta)=\eta,\qquad \mbox{ and}\\
G^0(\eta,W)\rest (m+1)=F_0(W\rest(m+1))(G^0(\eta,W)\rest m)\quad\mbox{for }m
\geq\lh(\eta),\\
G^1(\rho,W)(n)=F_1(W\rest (n+1))(\rho\rest n)\qquad\mbox{ for }n\in\omega.
  \end{array}\]
We are going to show that
\[\forces_{\q^*_\infty(K,\Sigma)}(\forall r\in\can\cap\V)(\exists\eta\in T)(
\forall^\infty n\in\omega)(r(n)=G^1(G^0(\eta,\dot{W}),\dot{W})(n)),\]
where $\dot{W}$ is the $\q^*_\infty(K,\Sigma)$--name defined in
\ref{thereal}. To this end suppose that $r\in\can$ and $p\in\q^*_\infty(K,
\Sigma)$. We may assume that $\lh(w^p)>0$ and $(\forall i\in\omega)(\nor[
t^p_i]\geq 3)$. 

Fix $i_0\in\omega$ for a moment. By downward induction on $i\leq i_0$ we
choose $s^{i_0}_i\in \Sigma(t^p_i)$ and $F^*_{i_0,i}: T\cap\omega^{\textstyle
\lh(w^p)+i-1}\longrightarrow T$ such that
\begin{enumerate}
\item[($\alpha$)] $\nor[s_{i_0,i}]\geq\nor[t^p_i]-1$,
\item[($\beta$)]  $F^*_{i_0,i}(\nu)\in\suc_T(\nu)$ for $\nu\in T$, $\lh(\nu)=
\lh(w^p)+i-1$,
\item[($\gamma$)] for all sequences $u\in\basis(s_{i_0,i})$ and $v\in\pos(u,
s_{i_0,i})$ and every $\nu\in T$ of length $\lh(\nu)=\lh(w^p)+i-1$ we have
\[F_0(v)(F^*_{i_0,i}(\nu))=F^*_{i_0,i+1}(\nu)\quad\mbox{ and }\quad F_1(v)(
F^*_{i_0,i}(\nu))=r(\lh(w^p)+i)\]
[for $i=i_0$ we omit the first part of the above demand].
\end{enumerate}
It should be clear that the choice of the $s_{i_0,i}$'s and $F^*_{i_0,i}$'s as
above is possible by \ref{bad}(3). All levels of the tree $T$ are finite, so
for each $i\in\omega$ there are finitely many mappings $F^*: T\cap\omega^{
\textstyle \lh(w^p)+i-1}\longrightarrow T\cap\omega^{\textstyle\lh(w^p)+i}$. 
Moreover, each $\Sigma(t^p_i)$ is finite (for $i\in\omega$). Hence, by
K\"onig's lemma, we find a sequence $\langle t^q_i: i\in\omega\rangle$ and a
mapping ${\bf F}^*:T\longrightarrow T$ such that
\begin{enumerate}
\item[(i)] \ $(\forall\nu\in T)({\bf F}^*(\nu)\in\suc_T(\nu))$ and
\item[(ii)]  for each $i\in\omega$ there is $j(i)>i$ such that for every
$j\leq i$
\[s^{j(i)}_j=t^q_j\quad\mbox{ and }\quad (\forall\nu\in
T\cap\omega^{\textstyle \lh(w^p)+j-1})({\bf F}^*(\nu)=F^*_{j(i),j}(\nu)).\]
\end{enumerate}
By $(\alpha)$ we have that $q=(w^p,t^q_0,t^q_1,\ldots)\in\q^*_\infty(K,
\Sigma)$ and it is stronger than $p$. Take any $\eta\in T$ with $\lh(\eta)=
\lh(w^p)-1$ and let $\eta_0={\bf F}^*(\eta)$, $\eta_{i+1}={\bf F}^*(\eta_i)$,
and $\eta^+=\lim(\eta_i)\in [T]$. It follows from $(\gamma)$ and (ii) (e.g.\
inductively using smoothness) that for each $n\in\omega$ and $v\in\pos(w^p,
t^q_0,\ldots,t^q_{n-1})$ we have 
\[\begin{array}{ll}
(v,t^q_n,t^q_{n+1},\ldots)\forces&\mbox{`` }G^0(\eta_0,\dot{W})\rest \lh(v)=
\eta^+\rest \lh(v)\ \mbox{ and}\\
\ &\ \ G^1(G^0(\eta_0,\dot{W}),\dot{W})\rest[\lh(\eta_0),\lh(v))= r\rest
[\lh(\eta_0),\lh(v))\mbox{''.}
  \end{array}\] 
Hence we conclude
\[q\forces_{\q^*_\infty(K,\Sigma)}(\forall n\geq\lh(\eta_0))(G^1(G^0(\eta_0,
\dot{W}),\dot{W})(n)=r(n)),\]
finishing the proof. 
\end{proof}

\begin{example}
\label{nonpro3}
Let $f:\omega\times\omega\longrightarrow\omega$ be a fast function (for
example $f(k,\ell)=2^{2k} (\ell+1)$). There are a finitary function $\bH$ and
a creating pair $(K_{\ref{nonpro3}},\Sigma_{\ref{nonpro3}})$ for $\bH$ such
that  
\begin{enumerate}
\item[(a)] $(K_{\ref{nonpro3}},\Sigma_{\ref{nonpro3}})$ is local, simple,
forgetful and definitely bad (and smooth),
\item[(b)] $\Sigma_{\ref{nonpro3}}(t)$ is finite for each $t\in
K_{\ref{nonpro3}}$,
\item[(c)] the forcing notions $\q^*_\infty(K_{\ref{nonpro3}},\Sigma_{
\ref{nonpro3}})$ and $\q^*_f(K_{\ref{nonpro3}},\Sigma_{\ref{nonpro3}})$ are not
trivial and thus collapse $\con$ onto $\omega$.
\end{enumerate}
\end{example}

\begin{proof}[Construction] 
Let $f:\omega\times\omega\longrightarrow\omega$ be fast. For $n\in\omega$ let
$k_n=2^{f(n+1,n+1)}$. Next, for $n\in\omega$, let $\bH(n)$ consist of all
pairs $\langle z_0,z_1\rangle$ such that 
\[z_0:\prod_{i<n}k_i\longrightarrow k_n\quad\mbox{ and }\quad z_1:\prod_{i<n}
k_i\longrightarrow 2.\]
Immediately by the definition, one sees that $\bH$ is finitary. Now we define
the creating pair $(K_{\ref{nonpro3}},\Sigma_{\ref{nonpro3}})$ for $\bH$. A
creature $t\in\CR[\bH]$ with $m^t_{\dn}>0$ is in $K_{\ref{nonpro3}}$ if:
\begin{itemize}
\item $\dis[t]=\langle m^t_{\dn},\ \langle A^t_\nu: \nu\in\prod\limits_{i<
m^t_{\dn}-1} k_i\rangle,\ F^t_0,F^t_1\rangle$, where

$A^t_\nu\subseteq k_{m^t_{\dn}-1}$ for $\nu\in\prod\limits_{i<m^t_{\dn}-1}
k_i$,

$F^t_0:\{\nu\in\prod\limits_{i<m^t_{\dn}} k_i: \nu(m^t_{\dn}-1)\in
A^t_{\nu\rest (m^t_{\dn}-1)}\}\longrightarrow k_{m^t_{\dn}}$,

$F^t_1:\{\nu\in\prod\limits_{i<m^t_{\dn}} k_i: \nu(m^t_{\dn}-1)\in
A^t_{\nu\rest (m^t_{\dn}-1)}\}\longrightarrow 2$,

\item $\val[t]$ consists of all pairs $\langle w,u\rangle\in\prod\limits_{
i<m^t_{\dn}}\bH(i)\times\prod\limits_{i\leq m^t_{\dn}}\bH(i)$ such that
$w\vartriangleleft u$ and if $u(m^t_{\dn})=\langle z_0,z_1\rangle$ then
$z_0\supseteq F^t_0$ and $z_1\supseteq F^t_1$.  
\item $\nor[t]=k_{m^t_{\dn}}-\max\{|A^t_\nu|:\nu\in\prod\limits_{i<m^t_{\dn}
-1}k_i\}$.
\end{itemize}
If $t\in\CR[H]$, $m^t_{\dn}=0$ then we take $t$ to $K_{\ref{nonpro3}}$ if:\\
$\nor[t]=0$,\quad $\dis[t]=\langle 0,x^t,i^t\rangle$, where $x^t\in k_0$,
$i^t<2$,\quad $\val[t]=\{\langle\;\langle\rangle,u\;\rangle\}$, where
$u\in\prod\limits_{i\leq 0}\bH(i)$, $u(0)=\langle x^t,i^t\rangle$ (both $x^t$
and $i^t$ are treated here as functions from $\{\langle\rangle\}$).

The composition operation $\Sigma_{\ref{nonpro3}}$ is the trivial one (so
$\Sigma_{\ref{nonpro3}}(\cS)$ is non-empty for singletons only and $\Sigma_{
\ref{nonpro3}}(t)=\{s\in K_{\ref{nonpro3}}:\val[s]\subseteq\val[t]\}$). 
Easily $(K_{\ref{nonpro3}},\Sigma_{\ref{nonpro3}})$ is a local, simple and
forgetful creating pair. Note that if $n>0$ and $t\in K_{\ref{nonpro3}}$ is
such that $m^t_{\dn}=n$ and $A^t_\nu=\emptyset$ for each $\nu\in\prod\limits_{
i<n-1}k_i$ then $\nor[t]=k_{n-1}>f(n,n)$, so the forcing notions $\q^*_\infty(
K_{\ref{nonpro3}},\Sigma_{\ref{nonpro3}})$, $\q^*_f(K_{\ref{nonpro3}},\Sigma_{
\ref{nonpro3}})$ are non-trivial. Finally let $T=\bigcup\limits_{n\in\omega}
\prod\limits_{i<n}k_i$  and for $v\in\prod\limits_{i\leq m}\bH(i)$ let
$F_0(v)$, $F_1(v)$ be such that $v(m)=\langle F_0(v),F_1(v)\rangle$. It should
be clear that $T,F_0,F_1$ witness that $(K_{\ref{nonpro3}},\Sigma_{
\ref{nonpro3}})$ is definitely bad. 
\end{proof}

\chapter{Properness and the reading of names} 
This chapter is devoted to getting the basic property: properness. The
first two sections deal with forcing notions determined by creating
pairs. We define properties of creating pairs implying that
appropriate forcing notions are proper. Some of these properties may
look artificial, but in applications they appear naturally. The third
part deals with forcing notions $\q^{\tree}_e(K,\Sigma)$ (determined
by a tree creating pair). Here, properness is an almost immediate
consequence of our choice of norm conditions. In most cases, proving
properness of a forcing notion we get much stronger property:
continuous reading of names for ordinals. This property will be
intensively used in the rest of the paper.  Finally, in the last part
of the chapter we give several examples for properties introduced and
studied before.

\section{Forcing notions ${\bf Q}^*_{{\rm s}\infty}(K,\Sigma)$, ${\bf
Q}^*_{{\rm w}\infty}(K,\Sigma)$}  
\label{2.1}

\begin{definition}
\label{omitoryetc}
Let $(K,\Sigma)$ be a creating pair for $\bH$.
\begin{enumerate}
\item For $t\in K$, $m_0\leq m^t_{\dn}$, $m^t_{\up}\leq m_1$ we define
the creature $s\stackrel{\rm def}{=}t\Rsh [m_0, m_1)$ by: 
\[\begin{array}{lll}
\nor[s]&=&\nor[t],\\
\dis[s]&=&\langle 4,m_0,m_1\rangle\conc\langle\dis[t]\rangle,\\
\val[s]&=&\{\langle w,u\rangle\in\prod\limits_{i<m_0}\bH(i)\times
\prod\limits_{i<m_1}\bH(i): \langle u\rest m^t_{\dn},u\rest m^t_{\up}\rangle
\in\val[t]\ \&\\  
\ &\ &\qquad w\vartriangleleft u\ \&\ (\forall i\in [m_0,m^t_{\dn})\cup
[m^t_{\up},m_1))(u(i)=0)\}.\\   
\end{array}\]
[Note that $t\Rsh [m_0,m_1)$ is well defined only if $\val[s]\neq\emptyset$
above and then $m^s_{\dn}=m_0$, $m^s_{\up}=m_1$.]
\item The creating pair $(K,\Sigma)$ is {\em omittory} if:
\begin{enumerate}
\item[$(\boxtimes_0)$] if $t\in K$ and $u\in\basis(t)$ then $u\conc {\bf
0}_{[m^t_{\dn},m^t_{\up})}\in\pos(u,t)$ but there is $v\in\pos(u,t)$ such that
$v\rest [m^t_{\dn},m^t_{\up})\neq {\bf 0}_{[m^t_{\dn},m^t_{\up})}$,
\item[$(\boxtimes_1)$] for every $(t_0,\ldots,t_{n-1})\in\PFC(K,\Sigma)$ and 
$i<n$: 
\[t_i\Rsh[m^{t_0}_{\dn},m^{t_{n-1}}_{\up})\in\Sigma(t_0,\ldots,t_{n-1}),\] 
\item[$(\boxtimes_2)$] if $t, t\Rsh [m_0,m_1)\in K$ then for every $u\in\basis
(t\Rsh [m_0,m_1))$ and $v\in\pos(u,t\Rsh [m_0,m_1))$ we have 
\[v(n)\neq 0\ \&\ n\in [\lh(u),\lh(v))\quad\Rightarrow\quad n\in [m^t_{\dn},
m^t_{\up}).\] 
\end{enumerate}
[Note that $(\boxtimes_0)$ implies that in the cases relevant for
$(\boxtimes_1)$, $t_i\Rsh [m^{t_0}_{\dn},m^{t_{n-1}}_{\up})$ is well defined.]
\item $(K,\Sigma)$ is {\em growing} if for any $(t_0,\ldots,t_{n-1})\in
\PFC(K,\Sigma)$ there is a creature $t\in\Sigma(t_0,\ldots,t_{n-1})$ such that
$\nor[t]\geq\max\limits_{i<n}\nor[t_i]$.  
\end{enumerate}
\end{definition}

\begin{proposition}
If $(K,\Sigma)$ is omittory then it is growing.
\end{proposition}

\begin{proposition}
\label{sinfty}
Suppose that a creating pair $(K,\Sigma)$ is growing. 
\begin{enumerate}
\item Then $\q^*_{{\rm s}\infty}(K,\Sigma)$ is a dense subset of both
$\q^*_\infty(K,\Sigma)$ and $\q^*_{{\rm w}\infty}(K,\Sigma)$ and
$\q^*_f(K,\Sigma)$ (for every fast (see \ref{fast}) function $f:\omega\times
\omega\longrightarrow\omega$). Consequently whenever we work with growing
creating pairs we may interchange the respective forcing notions as they are
equivalent.  
\item Moreover, if $g:\omega\times\omega\longrightarrow\omega$,
$p\in\q^*_{{\rm w}\infty}(K,\Sigma)$ then there is $q\in\q^*_{{\rm
s}\infty}(K,\Sigma)$ such that $p\leq_0^{{\rm w}\infty} q$ and
\[(\forall n\in\omega)(\nor[t^q_n]> g(n,m^{t^q_n}_{\dn})).\]
\end{enumerate}
\end{proposition}

\begin{proof} Suppose that $g:\omega\times\omega\longrightarrow\omega$ and
$(w,t_0,t_1,\ldots)\in\q^*_{{\rm w}\infty}(K,\Sigma)$. Choose an increasing
sequence $k_0<k_1<k_2<\ldots$ such that 
\[\nor[t_{k_0}]>g(0,m^{t_0}_{\dn})\quad\mbox{ and }\quad\nor[t_{k_{n+1}}]>
g(n+1,m^{t_{k_n}}_{\up})\]
(exists by \ref{conditions}(${\rm w}\infty$)) and choose
$s_0\in\Sigma(t_0,\ldots,t_{k_0})$, $s_{n+1}\in\Sigma(t_{k_n+1},\ldots,
t_{k_{n+1}})$ such that $\nor[s_n]\geq\nor[t_{k_n}]$ 
(exist by \ref{omitoryetc}(3)). Hence (by \ref{implications}(2); remember that
we assume $(K,\Sigma)$ is nice)  
\[q\stackrel{\rm def}{=}(w,s_0,s_1,s_2,\ldots)\in\q^*_{{\rm
s}\infty}(K,\Sigma)\]   
and clearly $(w,t_0,t_1,t_2,\ldots)\leq_0 q$. 
\end{proof}

\begin{theorem}
\label{deciding}
Assume $(K,\Sigma)$ is a finitary creating pair. Further assume that $p\in
\q^*_{{\rm s}\infty}(K,\Sigma)$ and for $n<\omega$ we have a $\q^*_{{\rm s}
\infty}(K,\Sigma)$-name $\dot{\tau}_n$ such that $\forces_{\q^*_{{\rm s}
\infty}(K,\Sigma)}$``$\dot{\tau}_n$ is an ordinal'' and $\ell<\omega$.
{\em Then} there is  $q=(w^p,s_0,s_1,s_2,\ldots)$ such that:
\begin{enumerate}
\item[(a)]  $p\leq_\ell^{{\rm s}\infty} q\in \q^*_{{\rm s}\infty}(K,\Sigma)$
and 
\item[(b)]  {\em if} $\ell\leq n<\omega$, $m\leq m^{s_{n-1}}_{\up}$ then the
condition $q$ approximates $\dot{\tau}_m$ at $s_n$ (see \ref{essapprox}(2)).
\end{enumerate}
\end{theorem}

\begin{proof} Let $p=\langle w^p,t^p_0,t^p_1,t^p_2,\ldots\rangle$. Let $s_i=
t^p_i$ for $i<\ell$. Now, by induction on $n\geq\ell$ we define $q_n,s_n,
t^n_{n+1},t^n_{n+2},\ldots$ such that:
\begin{enumerate}
\item[(i)]    $q_\ell=p$,
\item[(ii)]   $q_{n+1}=(w^p,s_0,\ldots,s_n,t^n_{n+1},t^n_{n+2},\ldots)\in
\q^*_{{\rm s}\infty}(K,\Sigma)$ 
\item[(iii)]  $q_n\leq_{n}^{{\rm s}\infty} q_{n+1}$
\item[(iv)]   if $w_1\in\pos(w^p,s_0,\ldots,s_{n -1})$, $m\leq
m^{s_{n-1}}_{\up}$ and there is a condition $r\in\q^*_{{\rm s}\infty}(K,
\Sigma)$, $\leq_0^{{\rm s}\infty}$--stronger than $(w_1, s_n,t^n_{n+1},t^n_{
n+2},\ldots)$ which decides the value of $\dot{\tau}_m$ 

then the condition $(w_1,s_n,t^n_{n+1},t^n_{n+2},\ldots)$ does it.  
\end{enumerate}
Arriving at the stage $n+1>\ell$ we have defined 
\[q_n=(w^p,s_0,\ldots,s_{n-1},t^{n-1}_n,t^{n-1}_{n+1},\ldots).\] 
Fix an enumeration $\langle (w^n_i, m^n_i):
i<k_n\rangle$ of
\[\pos(w^p,s_0,\ldots,s_{n-1})\times(m^{s_{n-1}}_{\up}+1)\]
(since each $\bH(m)$ is finite, $k_n$ is finite). Next choose by induction on
$k\leq k_n$ conditions $q_{n,k}\in\q^*_{{\rm s}\infty}(K,\Sigma)$ such that:
\begin{enumerate}
\item[$(\alpha)$] $q_{n, 0}=q_n$ 
\item[$(\beta)$]  $q_{n,k}$ is of the form $(w^p,s_0,\ldots,s_{n-1},t_n^{n,k}, 
t^{n, k}_{n+1}, t^{n, k}_{n + 1},\ldots)$
\item[$(\gamma)$] $q_{n,k}\leq_n q_{n,k+1}$
\item[$(\delta)$] if, in $\q^*_{{\rm s}\infty}(K,\Sigma)$, there is a
condition $r\geq_0 (w^n_k,t^{n,k}_n,t^{n,k}_{n+1},t^{n,k}_{n+2},\ldots)$ which
decides (in $\q^*_{{\rm s}\infty}(K,\Sigma)$) the value of $\dot{\tau}_{
m^n_k}$, then
\[(w^n_k,t^{n,k+1}_n,t^{n,k+1}_{n+1},t^{n,k+1}_{n+2},\ldots)\in \q^*_{{\rm s}
\infty}(K,\Sigma)\]
is a condition which forces a value to $\dot{\tau}_{{m^n_k}}$.    

(Note: $(w^n_k,t^{n,k}_n,t^{n,k}_{n+1},t^{n,k}_{n+2},\ldots)\in \q^*_{{\rm
s}\infty}(K,\Sigma)$.) 
\end{enumerate}
For this part of the construction we need our standard assumption that
$(K,\Sigma)$ is nice. Note that choosing
$(w^n_k,t^{n,k+1}_n,t^{n,k+1}_{n+1},t^{n,k+1}_{n+2},\ldots)$ we want to be 
sure that 
\[(w^p,s_0,\ldots,s_{n-1},t^{n,k+1}_n,t^{n,k+1}_{n+1},t^{n,k+1}_{n+2},\ldots)
\in\q^*_{{\rm s}\infty}(K,\Sigma)\]
(remember that \ref{maindef}(b)(ii) might fail). But by \ref{implications}(2)
it is not a problem. Next, the condition $q_{n+1}\stackrel{\rm def}{=}
q_{n,k_n}\in \q^*_{{\rm s}\infty}(K,\Sigma)$ satisfies (iv): the keys are the
clause ($\delta$) and the fact that 
\[(w^n_k,t^{n,k+1}_n,t^{n,k+1}_{n+1},t^{n,k+1}_{n+2},\ldots)\leq^{{\rm s}
\infty}_0 (w^n_k,t^{n,k_n}_n, t^{n,k_n}_{n+1},t^{n,k_n}_{n+2},\ldots)\in
\q^*_{{\rm s}\infty}(K,\Sigma).\]
Thus $s_n\stackrel{\rm def}{=} t^{n,k_n}_n$, $q_{n+1}$ and
$t^{n}_{n+k}\stackrel{\rm def}{=} t^{n,k_n}_{n+k}$ are as required. 

Now, by \ref{propord}:
\[q\stackrel{\rm def}{=}(w^p,s_0,s_1,\ldots,s_l,s_{l+1},\ldots)=\lim_n q_n\in
\q^*_{{\rm s}\infty}(K,\Sigma).\] 
Easily it satisfies the assertions of the theorem.
\end{proof}

A small modification of the proof of \ref{deciding} shows the corresponding
result for the forcing notion $\q^*_{{\rm w}\infty}(K,\Sigma)$:

\begin{theorem}
\label{winfty}
Assume $(K,\Sigma)$ is a finitary creating pair and $p\in\q^*_{{\rm w}\infty}(
K,\Sigma)$. Let $\dot{\tau}_n$ be $\q^*_{{\rm w}\infty}(K,\Sigma)$-names
for ordinals (for $n<\omega$), $\ell<\omega$. Then there is a condition
$q=(w^p,s_0,s_1,\ldots)\in\q^*_{{\rm w}\infty}(K,\Sigma)$ 
such that  
\begin{enumerate}
\item[(a)] $p\leq_\ell q$ and
\item[(b)] there is an increasing sequence $\ell=k_0<k_1<k_2<\ldots<\omega$
such that if $n<\omega$ and $m\leq m^{s_{k_n -1}}_{\up}$ then the condition
$q$ approximates $\dot{\tau}_m$ at $s_{k_n}$.
\end{enumerate}
\end{theorem}

As an immediate corollary to theorem \ref{deciding} we get the following.

\begin{corollary}
\label{omitoryproper}
Assume that $(K,\Sigma)$ is a finitary creating pair.  
\begin{enumerate}
\item[(a)] Suppose that $\dot{\tau}_n$ are $\q^*_{{\rm s}\infty}(K,
\Sigma)$-names for ordinals and $q\in\q^*_{{\rm s}\infty}(K,\Sigma)$ is a
condition satisfying {\bf (b)} of \ref{deciding} (for $\ell=0$). Further
assume that $q\leq r\in\q^*_{{\rm s}\infty}(K,\Sigma)$, $n<\lh(w^r)$ and
$r\forces$``$\dot{\tau}_n=\alpha$'' (for some ordinal $\alpha$).\\ 
{\em Then} for some $q'\in\q^*_{{\rm s}\infty}(K,\Sigma)$, $q\leq_{\apr}q'
\leq_{0} r$, we have $q'\forces$``$\dot{\tau}_n=\alpha$''.
\end{enumerate}
(Note: $\{q'\!\in\!\q^*_{{\rm s}\infty}(K,\Sigma)\!: q\leq_{\apr} q'\}$ is
countable provided $\bigcup\limits_{i<\omega}\bH(i)$ is countable.)   
\begin{enumerate}
\item[(b)] The forcing notion $\q^*_{{\rm s}\infty}(K,\Sigma)$ is proper (and
$\alpha$-proper for $\alpha<\omega_1$).
\end{enumerate}
\end{corollary}

It should be underlined here that \ref{omitoryproper} applies to forcing
notions $\q^*_\infty(K,\Sigma)$ for finitary growing creating pairs (remember
\ref{sinfty}). To get the respective conclusion for $\q^*_{{\rm w}\infty}(K,
\Sigma)$ we need to assume more.

\begin{definition}
\label{simpglui}
We say that:
\begin{enumerate}
\item A weak creating pair $(K,\Sigma)$ is {\em simple} if $\Sigma(\cS)$ is
non-empty for singletons only.
\item A creating pair $(K,\Sigma)$ is {\em gluing} if it is full and for every
$k<\omega$ there is $n_0<\omega$ such that for every $n_0\leq n<\omega$,
$(t_0,\ldots,t_n)\in\PFC(K,\Sigma)$, for some $s\in \Sigma(t_0,\ldots,t_n)$ we
have 
\[\nor[s]\geq \min\{k,\nor[t_0],\ldots,\nor[t_n]\}.\]
In this situation the integer $n_0$ is called {\em the gluing witness for
  $k$}.  
\end{enumerate}
\end{definition}

The two properties defined above are, in a sense, two extremal situations
under which we may say something on $(K,\Sigma)$. The demand ``either simple
or gluing'' (like in \ref{halbigdec}) should not be surprising if one realizes
that then we know what may happen when $\Sigma$ is applied, at least in terms
of $m^t_{\dn}$, $m^t_{\up}$.

\begin{corollary}
\label{winprop}
Assume that $(K,\Sigma)$ is a finitary and simple creating pair. 
\begin{enumerate}
\item[(a)] Suppose that $\dot{\tau}_n$ are $\q^*_{{\rm
w}\infty}(K,\Sigma)$--names for ordinals, $k_0<k_1<\ldots<\omega$ and $q\in
\q^*_{{\rm w}\infty}(K,\Sigma)$ are as in \ref{winfty}(b). Suppose that $q\leq
r\in\q^*_{{\rm w}\infty}(K,\Sigma)$, and $r$ decides the value of one of the
names $\dot{\tau}_n$, say $r\forces$``$\dot{\tau}_n=\alpha$''.

{\em Then} for some $q'\in\q^*_{{\rm w}\infty}(K,\Sigma)$ we have 
\[q\leq_{\apr}q',\quad q'\forces\mbox{``}\dot{\tau}_n=\alpha\mbox{''\quad and
} q',r\mbox{ are compatible}.\]
\item[(b)] The forcing notion $\q^*_{{\rm w}\infty}(K,\Sigma)$ is proper (and
even more).
\end{enumerate}
\end{corollary}

\begin{remark}
Note the presence of ``simple'' in the assumptions of \ref{winprop}. In
practical applications of forcing notions of the type $\q^*_{{\rm w}\infty}$
we can get more, see \ref{sinwin} below.
\end{remark}

\begin{definition}
\label{singleton}
Let $(K,\Sigma)$ be a creating pair for $\bH$. We say that $(K,\Sigma)$ {\em
captures singletons} if $(K,\Sigma)$ is forgetful and for every
$(t_0,\ldots,t_n)\in\PFC(K,\Sigma)$ and for each $u\in\basis(t_0)(=
\prod\limits_{m<m^{t_0}_{\dn}}\bH(m))$ and $v\in\pos(u,t_0,\ldots,t_n)$ {\em
there is} $(s_0,\ldots,s_k)\in\PFC(K,\Sigma)$ such that $(t_0,\ldots,t_n)\leq
(s_0,\ldots,s_k)$ (see \ref{candidates}) and
\[\pos(u,s_0,\ldots,s_k)=\{v\},\quad m^{t_0}_{\dn}=m^{s_0}_{\dn},\quad
m^{s_k}_{\up}=m^{t_n}_{\up}.\]
[Note that we put no demands on the norms of the $s_i$'s.] 
\end{definition}

\begin{proposition}
\label{sinimpzero}
Suppose that $(K,\Sigma)$ is a creating pair which captures singletons (so in
particular it is forgetful), $p\in\q^*_{{\rm w}\infty}(K,\Sigma)$ and
$\dot{\tau}$ is a $\q^*_{{\rm w}\infty}(K,\Sigma)$--name for an ordinal. Then 
there is $q\in\q^*_{{\rm w}\infty}(K,\Sigma)$ such that 
\[p\leq^{{\rm w}\infty}_0 q\qquad\mbox{ and\qquad $q$ decides $\dot{\tau}$}.\]
\end{proposition}

\begin{proof} Take $r\in \q^*_{{\rm w}\infty}(K,\Sigma)$ such that $p\leq r$
and $r\forces\dot{\tau}=\alpha$ (for some $\alpha$). Look at $w^r$: for some
$n\in\omega$ we have $w^r\in\pos(w^p,t^p_0,\ldots,t^p_{n-1})$. By
\ref{singleton} we find $s_0,\ldots,s_k$ such that $\pos(w^p,s_0,\ldots,
s_k)=\{w^r\}$ and 
\[(w^p,t^p_0,t^p_1,\ldots)\leq (w^p,s_0,\ldots,s_k,t^r_0,t^r_1,\ldots)
\stackrel{\rm def}{=}q\in \q^*_{{\rm w}\infty}(K,\Sigma).\]
Clearly $p\leq^{{\rm w}\infty}_0 q$. To show that $q\forces\dot{\tau}=\alpha$
we use our standard assumption that $(K,\Sigma)$ is smooth. Suppose that
$q'\geq q$ is such that $q'\forces\dot{\tau}\neq\alpha$. We may assume that
$\lh(w^{q'})>m^{s_k}_{\up}$. By the smoothness we have $w^{q'}\rest m^{
s_k}_{\up}\in\pos(w^p,s_0,\ldots,s_k)$, and so $w^{q'}\rest m^{s_k}_{\up}=
w^r$, and $w^{q'}\in\pos(w^r,t^r_0,\ldots,t^r_\ell)$ (for a suitable
$\ell<\omega$). Consequently $q'\geq r$ and this is a contradiction. 
\end{proof}

\begin{corollary}
\label{sinwin}
Assume $(K,\Sigma)$ is a finitary creating pair which captures singletons.
\begin{enumerate}
\item Let $p\in\q^*_{{\rm w}\infty}(K,\Sigma)$, $\dot{\tau}_n$ be $\q^*_{{\rm
w}\infty}(K,\Sigma)$-names for ordinals (for $n<\omega$) and $\ell<\omega$.
Then there is a condition $q\in\q^*_{{\rm w}\infty}(K,\Sigma)$ such that  
\begin{enumerate}
\item[(a)] $p\leq_\ell q$ and
\item[(b)] the condition $q$ essentially decides (see \ref{essapprox}(1)) each
name $\dot{\tau}_n$.
\end{enumerate}
\item The forcing notion $\q^*_{{\rm w}\infty}(K,\Sigma)$ is proper.
\end{enumerate}
\end{corollary}

\begin{proof} 1)\quad Follows from \ref{winfty} and \ref{sinimpzero}.\\
2)\quad Follows from 1). 
\end{proof}
 
\section{Forcing notion ${\bf Q}^*_{f}(K,\Sigma)$: bigness and halving} 

\begin{definition}
\label{bigetc}
Let $\bar{r}=\langle r_m: m<\omega\rangle$ be a non-decreasing sequence of
integers $\geq 2$. For a creating pair $(K,\Sigma)$ for $\bH$ we say that:
\begin{enumerate}
\item $(K,\Sigma)$ is {\em big} if for every $k<\omega$ there is
$m<\omega$ such that:
 
\noindent{\em if} $t\in K$, $\nor[t]\geq m$, $u\in\basis(t)$ and $c:\pos(u,t)
\longrightarrow\{0,1\}$  

\noindent{\em then} there is $s\in\Sigma(t)$ such that $\nor[s]\geq k$, and
$c\rest\pos(u,s)$ is constant.  

\noindent In this situation we call $m$ {\em a bigness witness for $k$}.

\item $(K,\Sigma)$ is {\em $\bar{r}$--big} if for each $t\in K$ such that
$\nor[t]>1$ and $u\in\basis(t)$ and $c:\pos(u,t)\longrightarrow r_{m^t_{\dn}}$
there is $s\in \Sigma(t)$ such that $\nor[s]\geq\nor[t]-1$ and $c\restriction
\pos(u,s)$ is constant.
\end{enumerate}
\end{definition}

\begin{remark}
Clearly, for a creating pair $(K,\Sigma)$, $\bar{r}$--big implies
implies big. 
\end{remark}

To show how the notions introduced in definition~\ref{bigetc} work we start
with proving an application of \ref{deciding} to the case when the creating
pair is additionally big and growing. 

\begin{proposition}
\label{decbel}
Assume $(K, \Sigma)$ is a finitary, growing and big creating pair. If $p\in
\q^*_{{\rm s}\infty}(K,\Sigma)$, $p\forces$``$\dot{\tau}<m$'', $m<\omega$ {\em
then} there is $q$, $p\leq_0 q\in\q^*_{{\rm s}\infty}(K,\Sigma)$ such that
$q\forces$``$\dot{\tau} = m_0$'' for some $m_0$.   
\end{proposition}

\begin{proof} Let $h\in\baire$ be such that $h(k)$ is a bigness witness for
$k$ (remember that $(K,\Sigma)$ is big, see \ref{bigetc}(1)). Note that 
\ref{deciding} $+$ \ref{sinfty}(2) give us a condition
$p'=(w^p,s_0,s_1,s_2,\ldots)\in\q^*_{{\rm s}\infty}(K,\Sigma)$ such that
$p\leq^{{\rm s}\infty}_0 p'$ and
\begin{enumerate}
\item[$(\alpha)$] $\nor[s_\ell]\geq h^{(m\cdot|\pos(w^p,s_0,\ldots,
s_{\ell-1})|)}(m^{s_\ell}_{\dn}+1)$ for all $\ell<\omega$ and 
\item[$(\beta)$] $p'$ approximates the name $\dot{\tau}$ at each $n<\omega$.
\end{enumerate}
Using iteratively the choice of $h(k)$ we will have then
\begin{enumerate}
\item[$(\gamma)$] for every $\ell<\omega$ and each function 
\[d:\{\langle u,v\rangle: u\in\pos(w^p,s_0,\ldots,s_{\ell-1})\ \&\
v\in\pos(u,s_\ell)\} \longrightarrow m+1\]
there is a creature $s\in\Sigma(s_\ell)$ such that $\nor[s]>m^{s_\ell}_{\dn}$
and  
\[d\restriction\{\langle u,v\rangle\in\dom(d): v\in\pos(u,s)\}\]
depends on the first coordinate only. 
\end{enumerate}
(Since, as usual, we assume that $(K,\Sigma)$ is nice, we have in ($\gamma$)
above that $\basis(s)\supseteq\pos(w^p,s_0,\ldots,s_{\ell-1})$.) 

Now apply $(\gamma)$ to find $s^\prime_\ell\in\Sigma(s_\ell)$ (for $\ell<
\omega$) such that $\nor[s^\prime_\ell]>m^{s_\ell}_{\dn}=m^{s_\ell^\prime}_{
\dn}$ and for every $u\in\pos(w^p,s^\prime_0,\ldots,s^\prime_{\ell-1})$ we
have  
\begin{enumerate}
\item[$(\delta)$] for each $v_0,v_1\in\pos(u,s^\prime_\ell)$, $i<m$
\[(v_0,s^\prime_{\ell+1},s^\prime_{\ell+2},\ldots)\forces
\dot{\tau}=i\quad\mbox{ iff }\quad
(v_1,s^\prime_{\ell+1},s^\prime_{\ell+2},\ldots)\forces \dot{\tau}=i.\]
\end{enumerate}
Look at $q\stackrel{\rm def}{=}(w^p,s^\prime_0,s^\prime_1,s^\prime_2,\ldots)$.
First it is a condition in $\q^*_{{\rm s}\infty}(K,\Sigma)$ as $(K,\Sigma)$ is
nice and $s^\prime_\ell\in\Sigma(s_\ell)$, $\nor[s_\ell^\prime]>
m^{s_\ell^\prime}_{\dn}$. To show that 
\[(\exists m_0<m)(q\forces \dot{\tau}=m_0)\]
take a condition $r=(w_1,t_0,t_1,\ldots)\geq q$ such that $r$ decides the
value of $\dot{\tau}$, $w_1\in\pos(w^p,s_0^\prime,\ldots,s_{\ell-1}^\prime)$
and $\ell$ is the smallest possible. By $(\beta)$ we know that the condition
$(w_1,s^\prime_\ell,s^\prime_{\ell+1},\ldots)$ forces a value to $\dot{\tau}$,
say $m_0$. We claim that $\ell=0$, i.e.~$w_1=w^p$ (which is enough as then $q$
decides the value of $\dot{\tau}$). Why? Suppose that $\ell>0$ and look at the
requirement $(\delta)$ for $\ell-1$, $u=w_1\restriction m^{s^\prime_{\ell-
1}}_{\dn}$. By the smoothness $u\in\pos(w^p,s^\prime_0,\ldots,s^\prime_{\ell
-2})$ and consequently, by $(\delta)$, for each
$v\in\pos(u,s^\prime_{\ell-1})$ we have 
\[(v,s^\prime_\ell,s^\prime_{\ell+1},\ldots)\forces\dot{\tau}=m_0.\]
Applying smoothness once again we note that for each $w\in\pos(u,
s^\prime_{\ell-1},s^\prime_{\ell},\ldots,s^\prime_k)$ 
\[w\restriction m^{s^\prime_\ell}_{\dn}\in\pos(u,s^\prime_{\ell-1})\quad\mbox{
and }\quad w\in\pos(w\restriction m^{s^\prime_\ell}_{\dn},
s^\prime_\ell,\ldots, s^\prime_{k}).\]
Hence for each such $w$ we have
\[(w\restriction m^{s^\prime_\ell}_{\dn}, s^\prime_\ell, s^\prime_{\ell+1},
\ldots,)\forces\dot{\tau}=m_0\quad\mbox{ and}\]
\[(w, s^\prime_{k+1},s^\prime_{k+2},\ldots)\geq (w\restriction
m^{s^\prime_\ell}_{\dn}, s^\prime_\ell, s^\prime_{\ell+1},\ldots)\]
and so 
\[(w, s^\prime_{k+1}, s^\prime_{k+2},\ldots,)\forces\dot{\tau}=m_0.\]
Hence we may conclude that $(u,s^\prime_{\ell-1},s^\prime_{\ell},\dots)
\forces\dot{\tau}=m_0$ which contradicts the choice of $\ell$.
\end{proof}

\begin{remark}
One may notice that the assumptions of \ref{decbel} are difficult to satisfy
in most natural cases. First examples of growing creating pairs one has in
mind are omittory creating pairs. However, if we demand that an omittory
creating pair $(K,\Sigma)$ is smooth then we get to 
\[t\in K\ \ \&\ \ u\in\basis(t)\quad\ \Rightarrow\ \quad u\conc{\bf
0}_{[m^t_{\dn},m^t_{\up})}\in\pos(u,t).\]
This excludes bigness as defined in \ref{bigetc}. Thus it is desirable to
consider in this case a weaker condition, which more fits to specific
properties of omittory creating pairs.
\end{remark}

\begin{definition}
\label{omitbig}
An omittory creating pair $(K,\Sigma)$ is {\em omittory--big} if for every
$k<\omega$ there is  $m<\omega$ such that:
 
\noindent{\em if} $t\in K$, $\nor[t]\geq m$, $u\in\basis(t)$ and $c:\pos(u,t)
\longrightarrow\{0,1\}$ {\em then} there is $s\in\Sigma(t)$ such that $\nor[s]
\geq k$ and $c\rest(\pos(u,s)\setminus\{{\bf 0}_{[m^t_{\dn},m^t_{\up})}\})$ is
constant. 

\noindent We may call $m$ {\em an omittory-bigness witness for $k$}.
\end{definition}

\begin{proposition}
\label{omitdecbel}
Assume $(K,\Sigma)$ is a finitary, omittory and omittory--big creating pair. 
Suppose that $p\in\q^*_{{\rm s}\infty}(K,\Sigma)$, $p\forces$``$\dot{\tau}<
m$'', $m<\omega$. Then there is a condition $q\in\q^*_{{\rm s}\infty}(K,
\Sigma)$ such that $p\leq_0^{{\rm s}\infty} q$ and $q$ decides the value of
$\dot{\tau}$. 
\end{proposition}

\begin{proof} We start as in the proof of \ref{decbel}, but in $(\gamma)$
there we say that 
\[d\restriction\{\langle u,v\rangle\in\dom(d): v\in\pos(u,s)\ \&\ (\exists
k\in [\lh(u),\lh(v)))(v(h)\neq 0)\}\]
depends on the first coordinate only, and therefore we get $s'_\ell\in\Sigma
(s_\ell)$ as there but with $(\delta)$ replaced by
\begin{enumerate}
\item[$(\delta)^-$] for each $v_0,v_1\in\pos(u,s^\prime_\ell)\setminus\{u\conc
{\bf 0}_{[m^{s_\ell'}_{\dn},m^{s'_\ell}_{\up})}\}$, $i<m$
\[(v_0,s^\prime_{\ell+1},s^\prime_{\ell+2},\ldots)\forces
\dot{\tau}=i\quad\mbox{ iff }\quad
(v_1,s^\prime_{\ell+1},s^\prime_{\ell+2},\ldots)\forces \dot{\tau}=i.\]
\end{enumerate}
Now comes the main modification of the proof of \ref{decbel}. We choose an
infinite set $I=\{i_0,i_1,i_2,\ldots\}\subseteq\omega$ such that for every
$i<m$ we have 
\begin{description}
\item[if] $k<\ell\leq\ell'<\omega$, $w\in\pos(w^p,s'_0,\ldots,s'_{i_k})$,
$v_0\in\pos(w,s'_{i_{\ell}}\Rsh [m^{s'_{i_k}}_{\up}, m^{s'_{i_\ell}}_{\up}))$,
$v_1\in\pos(w,s'_{i_{\ell'}}\Rsh [m^{s'_{i_k}}_{\up}, m^{s'_{i_{\ell'}}}_{
\up}))$, and 
\[(\exists m\in [m^{s'_{i_\ell}}_{\dn},m^{s'_{i_\ell}}_{\up}))(v_0(m)\neq 0)
\quad\mbox{and}\quad(\exists m\in [m^{s'_{i_{\ell'}}}_{\dn},m^{s'_{i_{
\ell'}}}_{\up}))(v_1(m)\neq 0)\]
\item[then]\quad $(v_0,s^\prime_{i_\ell+1},s^\prime_{i_\ell+2},\ldots)\forces
\dot{\tau}=i\quad\mbox{ iff }\quad (v_1,s^\prime_{i_{\ell'}+1},s^\prime_{i_{
\ell'}+2},\ldots)\forces \dot{\tau}=i$
\end{description}
and a similar condition for the case of $w=w^p$, $0\leq\ell\leq\ell'<\omega$.
The construction of the set $I$ is rather standard (by induction) and it goes
like in the proof of the suitable property for the Mathias forcing (see e.g.
\cite[7.4.6]{BaJu95}). Next we look at 
\[q\stackrel{\rm def}{=}(w^p,s'_{i_0}\Rsh [m^{s'_0}_{\dn},m^{s'_{i_0}}_{\up}),
s'_{i_1}\Rsh [m^{s'_{i_0}}_{\up},m^{s'_{i_1}}_{\up}), s'_{i_2}\Rsh
[m^{s'_{i_1}}_{\up},m^{s'_{i_2}}_{\up}),\ldots).\]
It should be clear that it is a condition in $\q^*_{{\rm s}\infty}(K,\Sigma)$
which is $\leq^{{\rm s}\infty}_0$--stronger than $p$. Note that, as
$(K,\Sigma)$ is omittory (remember the demand $(\boxtimes_0)$ of
\ref{omitoryetc}(2)), by the choice of the set $I$ we have  
\begin{quotation}
\noindent if $w_1\in\pos(w^q,t^q_0,\ldots,t^q_k)$, $k<\omega$, $w_1\rest
[m^{t^q_k}_{\dn}, m^{t^q_k}_{\up})={\bf 0}_{[m^{t^q_k}_{\dn},
m^{t^q_k}_{\up})}$ and the condition $(w_1, t^q_{k+1},t^q_{k+2},\ldots)$
decides the value of $\dot{\tau}$ 

\noindent then $(w_1\rest m^{t^q_k}_{\dn},t^q_k,t^q_{k+1},t^q_{k+2},\ldots)$
does so.
\end{quotation}
Now, like in \ref{decbel} we show that the condition $q$ decides the value of
$\dot{\tau}$, using the remark above and the choice of the set $I$.
\end{proof}

\begin{definition}
\label{halving}
Let $(K,\Sigma)$ be a creating pair.
\begin{enumerate}
\item We say that $(K,\Sigma)$ has the {\em Halving Property} if there is a
mapping 
\[\uhalf:K\longrightarrow K\]
such that
\begin{enumerate}
\item[(a)] for each $t\in K$,  $\uhalf(t)\in\Sigma(t)$ and $\nor[\uhalf(t)]
\geq\frac{1}{2}\nor[t]$, 
\item[(b)] if $t_0,\ldots,t_n\in K$, $\min\{\nor[t_i]: i\leq n\}\geq 2$ and
a creature $t\in\Sigma(\uhalf(t_0),\ldots,\uhalf(t_n))$ is such that
$\nor[t]>0$ then there is $s\in\Sigma(t_0,\ldots,t_n)$ such that 
\[\nor[s]\geq\min\{\frac{1}{2}\nor[t_i]: i\leq n\}\quad\mbox{and}\quad(\forall
u\!\in\!\basis(t_0))(\pos(u,s)\subseteq\pos(u,t)).\]  
\end{enumerate}
\item We say that $(K,\Sigma)$ has the {\em weak Halving Property} if there is
a mapping $\uhalf:K\longrightarrow K$ which satisfies {\bf (a)} above and
\begin{enumerate}
\item[(b)$^-$] if $t_0\in K$, $\nor[t_0]\geq 2$ and $t\in\Sigma(\uhalf(t_0))$
is such that $\nor[t]>0$

then there is a creature $s\in\Sigma(t_0)$ such that 
\[\nor[s]\geq\frac{1}{2}\nor[t_0]\quad\mbox{and}\quad(\forall u\in
\basis(t_0))(\pos(u,s)\subseteq\pos(u,t)).\]  
\end{enumerate}
\item Whenever we say that $(K,\Sigma)$ has the (weak) Halving Property we
assume that the function $\uhalf:K\longrightarrow K$ witnessing this is fixed.
\end{enumerate}
\end{definition}

\begin{remark}
Remember that we standardly assume that creating pairs are nice, so in
\ref{halving}(1b), \ref{halving}(2b$^-$) we have $\basis(t_0)\subseteq
\basis(t)$ and $\basis(t_0)\subseteq\basis(s)$. Of course, the Halving
Property implies the weak Halving Property. Moreover, the two notions agree
for simple creating pairs.  
\end{remark}

Our next lemma shows how we are going apply the Halving Property.

\begin{lemma}
\label{halback}
Assume that $(K,\Sigma)$ is a  creating pair with the Halving Property
(witnessed by a mapping $\uhalf:K\longrightarrow K$). Suppose that $f:\omega
\times\omega\longrightarrow\omega$ is a fast function, $\dot{\tau}$ is a
$\q^*_f(K,\Sigma)$-name for an ordinal, $n<\omega$, $0<\varepsilon\leq 1$ and
$p\in\q^*_f(K,\Sigma)$ is a condition such that
\[(\forall i\in\omega)(\nor[t^p_i]\geq\varepsilon\cdot
f(n,m^{t^p_i}_{\dn})).\]
Further assume that there is a condition $r\in\q^*_f(K,\Sigma)$ such that
\[(\forall i\in\omega)(\nor[t^r_i]>0)\quad\mbox{ and }\quad (w^p,\uhalf(
t^p_0),\uhalf(t^p_1),\ldots)\leq_0^f r\quad\mbox{ and}\]
$r$ essentially decides $\dot{\tau}$ (see \ref{essapprox}). {\em Then} there
is a condition $q\in\q^*_f(K,\Sigma)$ such that  
\[(\forall i\in\omega)(\nor[t^q_i]\geq\frac{\varepsilon}{2}\cdot
f(n,m^{t^q_i}_{\dn})),\quad p\leq_0^f q\quad\mbox{ and $q$ essentially decides
}\dot{\tau}.\]  
\end{lemma}

\begin{proof} First note that the niceness implies that $(w^p,\uhalf(t^p_0),
\uhalf(t^p_1),\ldots)$ is a condition in $\q^*_f(K,\Sigma)$ (by
\ref{halving}(1)(a) and \ref{implications}(2); remember that $f$ is fast). Now
suppose that $r\in\q^*_f(K,\Sigma)$ is as in the assumptions of the
lemma. Take $m<\omega$ so large that  
\[(\forall u\in\pos(w^p,t^r_0,\ldots,t^r_{m-1}))((u,t^r_m,t^r_{m+1},\ldots)
\mbox{ decides the value of }\dot{\tau})\quad\mbox{ and}\]
\[(\forall i\geq m)(\nor[t^r_i]\geq\frac{\varepsilon}{2}
f(n,m^{t^r_i}_{\dn}))\] 
(for the first requirement remember that $(K,\Sigma)$ is smooth; the second is
possible since $\varepsilon\leq 1$). Next choose integers 
$0=i_0<i_1<\ldots<i_{m-1}<i_m$ such that 
\[(\forall \ell<m)(t^r_\ell\in\Sigma(\uhalf(t^p_{i_\ell}),\ldots,
\uhalf(t^p_{i_{\ell+1}-1}))).\] 
Applying the Halving Property (see \ref{halving}(1b); remember that we have
assumed $\nor[t^r_\ell]>0$ for each $\ell\in\omega$) we find $s_\ell\in 
\Sigma(t^p_{i_\ell},\ldots,t^p_{i_{\ell+1}-1})$ (for $\ell<m$) such that
\[\begin{array}{l}
\nor[s_\ell]\geq\min\{\frac{1}{2}\nor[t^p_i]: i_\ell\leq i<i_{\ell+1}\}\quad
\mbox{ and}\\
\quad (\forall u\!\in\!\basis(t^p_{i_\ell}))(\pos(u,s_\ell)\subseteq\pos(u,
t^r_\ell)).
  \end{array}\]
Then easily 
\[q\stackrel{\rm def}{=}(w^p,s_0,\ldots,s_{m-1},t^r_m,t^r_{m+1},\ldots)\in
\q^*_f(K,\Sigma),\quad\ p\leq_0^f q\quad\mbox{ and}\]
\[(\forall i\in\omega)(\nor[t^q_i]\geq\frac{\varepsilon}{2}
f(n,m^{t^q_i}_{\dn}))\]
(for the last statement remember that $f$ is fast so $f(n,\cdot)$ is
non-decreasing). Moreover $q$ essentially decides the value of $\dot{\tau}$ as 
\[\pos(w^p,s_0,\ldots,s_{m-1})\subseteq\pos(w^p,t^r_0,\ldots,t^r_{m-1}).\]
\end{proof}

\begin{remark}
One could ask why we cannot in the conclusion of \ref{halback} require that
simply $\nor[t^q_i]\geq \frac{1}{2}\nor[t^p_j]$ (for $j$ chosen somehow
suitably, e.g.\ such that $m^{t^q_i}_{\dn}=m^{t^p_j}_{\dn}$). The reason is
that ``the upgrading procedure'' given by \ref{halving}(1b) takes care of
possibilities only: no other relation between $s,t$ there is required. In
particular we do not know if $\Sigma(t)\cap\Sigma(s)\neq\emptyset$. 
Consequently, if we apply this procedure to {\em all} $m$ (replacing each
$t^r_m$ by suitable $s_m$) then we may get a condition incompatible with $r$.
\end{remark}

\begin{theorem}
\label{halbigdec}
Assume that a creating pair $(K,\Sigma)$ for $\bH$ is finitary, $\bar{2}$-big
and has the Halving Property. Further suppose that a function $f:\omega\times
\omega\longrightarrow\omega$ is $\bH$-fast and that $(K,\Sigma)$ is either
simple or gluing (see \ref{simpglui}). Let $\dot{\tau}_m$ be $\q^*_f(K,
\Sigma)$-names for ordinals (for $m\in\omega$), $p\in\q^*_f(K,\Sigma)$ and
$n<\omega$.  

\noindent Then there is a condition $q\in\q^*_f(K,\Sigma)$ such that
$p\leq_n^f q$ and $q$ essentially decides (see \ref{essapprox}) all the names
$\dot{\tau}_m$ (for $m\in\omega$). 
\end{theorem}

\begin{proof} Let $(K,\Sigma)$ and $f$ be as in the assumptions of the
theorem. 

\begin{claim}
\label{cl7}
Suppose that $\dot{\tau}$ is a $\q^*_f(K,\Sigma)$-name for an ordinal,
$n<\omega$, a condition $p\in\q^*_f(K,\Sigma)$ and a real $\varepsilon$
are such that
\[2^{|\pos(w^p,t^p_0)|}\cdot 2^{-\fH(m^{t^p_1}_{\dn})}\leq\varepsilon\leq 1\] 
(where $\fH(\ell)=|\prod\limits_{i<\ell}\bH(i)|$, see \ref{fast}), and
\[\nor[t^p_0]>1\quad\mbox{ and }\quad(\forall i>0)(\nor[t^p_i]\geq\varepsilon
\cdot f(n+1, m^{t^p_i}_{\dn})).\] 
Then there is a condition $q\in\q^*_f(K,\Sigma)$ such that $p\leq_0^f q$,
$t^q_0\in\Sigma(t^p_0)$, $\nor[t^q_0]\geq\nor[t^p_0]-1$, $q$ essentially
decides the value of $\dot{\tau}$ and 
\[(\forall i>0)(\nor[t^q_i]\geq\frac{\varepsilon}{2^{|\pos(w^p,t^p_0)|}}\cdot
f(n+1, m^{t^q_i}_{\dn})).\] 
\end{claim}

\noindent{\em Proof of the claim:}\ \ \ First note that our assumptions on
$\varepsilon$ (and the fact that $f$ is $\bH$-fast) imply that if
$t\in K$ is such that $m^t_{\dn}\geq m^{t^p_1}_{\dn}$ and
\[\nor[t]\geq\frac{\varepsilon}{2^{|\pos(w^p,t^p_0)|}}\cdot f(n+1,m^t_{\dn})\]
then 
\[\nor[t]>f(n,m^t_{\dn})+\fH(m^t_{\dn})+2.\]
Let $\{w^0_m: m< m_0\}$ enumerate $\pos(w^p,t^p_0)$. We inductively choose
conditions $p^0_m$ for $m\leq m_0$ such that 
\begin{enumerate}
\item[$(\alpha)$] $p^0_0=(w^0_0,t^p_1,t^p_2,\ldots)$,
\item[$(\beta)$]  if there is a condition $r\in\q^*_f(K,\Sigma)$ such that
$p^0_m\leq_0 r$ and
\[(\forall i\in\omega)(\nor[t^r_i]\geq\frac{\varepsilon}{2^{m+1}}\cdot
f(n+1,m^{t^r_i}_{\dn}))\]
and $r$ essentially decides $\dot{\tau}$ then we choose such an $r$ and we put
\[p^0_{m+1}=(w^0_{m+1}, t^r_0,t^r_1,\ldots),\]
\item[$(\gamma)$] if we cannot apply the clause $(\beta)$ (i.e.~there is no
$r$ as above) then 
\[p^0_{m+1}=(w^0_{m+1},\uhalf(t^{p^0_m}_0),\uhalf(t^{p^0_m}_1),\ldots).\]
\end{enumerate}
[For the sake of the uniformity of the inductive definition we let $w^0_{m_0}$
to be e.g.~$w^0_0$.] 

\noindent Note that by the niceness there are no problems in the above
construction (i.e.\ $p^0_m\in\q^*_f(K,\Sigma)$). Let $c:\pos(w^p,t^p_0)
\longrightarrow 2$ be such that 
\[c(w^0_m)=\left\{
\begin{array}{ll}
0 &\mbox{ if the clause $(\beta)$ was applied to choose } p^0_{m+1},\\
1 &\mbox{ otherwise.}
\end{array}
\right.\]
Due to the bigness of $(K,\Sigma)$ we find a creature $s_0\in\Sigma(t^p_0)$
such that $c$ is constant on $\pos(w^p,s_0)$ and $\nor[s_0]\geq
\nor[t^p_0]-1$. Note that 
\[q\stackrel{\rm def}{=}(w^p,s_0,t^{p^0_{m_0}}_0,t^{p^0_{m_0}}_1,
t^{p^0_{m_0}}_2,\ldots)\in\q^*_f(K,\Sigma)\]
(the niceness applies here once again, see \ref{implications}; the norm
condition should be obvious as $p^0_{m_0}\in\q^*_f(K,\Sigma)$). Moreover 
\[\nor[t^{p^0_{m_0}}_i]\geq\frac{\varepsilon}{2^{|\pos(w^p,t^p_0)|}}\cdot
f(n+1,m^{t^{p^0_{m_0}}_i}_{\dn})\] 
and hence, in particular,
\[\nor[t^{p^0_{m_0}}_i]> f(n,m^{t^{p^0_{m_0}}_i}_{\dn})+
\fH(m^{t^{p^0_{m_0}}_i}_{\dn}) +2.\]  
If the constant value of $c\restriction\pos(w^p,s_0)$ is $0$ then easily the
condition $q$ satisfies the requirements of the claim (use \ref{smoessdec}).

So we want to exclude the possibility that the constant value is $1$. For
this assume that it is the case. First note that then, due to the way we
constructed $p^0_{m_0}$, we may apply lemma~\ref{halback} and conclude that
there are no $u\in\pos(w^p,s_0)$ and $r\in\q^*_f(K,\Sigma)$ such that
$\nor[t^r_i]>0$ for all $i\in\omega$ and  
\[(u,t^{p^0_{m_0}}_0,t^{p^0_{m_0}}_1,\ldots)\leq_0 r\quad\mbox{ and }\quad
r\mbox{ essentially decides the value of }\dot{\tau}.\]
Now we inductively choose an increasing sequence $\ell_0<\ell_1<\ldots$ of
integers, creatures $s_1,s_2,\ldots\in K$ and conditions $p_0,p_1,\ldots\in
\q^*_f(K,\Sigma)$ such that 
\begin{enumerate}
\item $p_0=q$ (defined above), $\ell_0=0$,
\item $\ell_{k+1}$ is such that $\ell_{k+1}>\ell_k$ and
\[(\forall i\geq\ell_{k+1})(\nor[t^{p_k}_i]\geq f(n+k+2,m^{t_i^{p_k}}_{
\dn})),\] 
\item $s_i=t^{p_{k}}_i$ for $i\leq\ell_k$ and $\nor[t^{p_k}_i]>f(n+k,
m^{t^{p_k}_i}_{\dn})+\fH(m^{t^{p_k}_i}_{\dn})+2$ for $i>\ell_k$,
\item there are no $u\in\pos(w^p,s_0,\ldots,s_{\ell_k})$ and $r\in\q^*_f(K,
\Sigma)$ such that 
\[(\forall i\in\omega)(\nor[t^r_i]>0)\]
and $(u,t^{p_k}_{\ell_k+1},t^{p_k}_{\ell_k+2},\ldots)\leq_0 r$ and $r$
essentially decides $\dot{\tau}$,
\item for $i\in (\ell_k,\ell_{k+1}]$, $s_i\in \Sigma(t^{p_k}_i)$ and
$\nor[s_i]>f(n+k,m^{s_i}_{\dn})$,
\item $\nor[s_{\ell_{k+1}}]> f(n+k+1,m^{s_{\ell_{k+1}}}_{\dn})$ and
$p_k\leq_{k}^f p_{k+1}$. 
\end{enumerate}
Suppose we have defined $\ell_k$, $p_k$ and $s_i$ for $i\leq\ell_k$. Choose
$\ell_{k+1}$ according to the requirement (2) above. Fix an enumeration 
\[\{w^k_m: m<m_k\}=\pos(w^p,s_0,\ldots,s_{\ell_k},t^{p_k}_{\ell_k+1},\ldots,
t^{p_k}_{\ell_{k+1}}).\]
We inductively choose conditions $p^k_m\in\q^*_f(K,\Sigma)$ (in the way
analogous to the construction of $p^0_m$'s):
\begin{enumerate}
\item[$(\alpha)_k$] $p^k_0=(w^k_0,t^{p_k}_{\ell_{k+1}+1},
t^{p_k}_{\ell_{k+1}+2},\ldots)$,
\item[$(\beta)_k$]  if there is a condition $r\in\q^*_f(K,\Sigma)$ such that
$p^k_m\leq_0 r$ and
\[(\forall i\in\omega)(\nor[t^r_i]\geq\frac{1}{2^{m+1}}\cdot
f(n+k+2,m^{t^r_i}_{\dn}))\]
and $r$ essentially decides $\dot{\tau}$ then we choose such an $r$ and we put
\[p^k_{m+1}=(w^k_{m+1},t^r_0,t^r_1,\ldots),\]
\item[$(\gamma)_k$] if we cannot apply the clause $(\beta)_k$ then 
\[p^k_{m+1}=(w^k_{m+1},\uhalf(t^{p^k_m}_0),\uhalf(t^{p^k_m}_1),\ldots).\]
\end{enumerate}
As previously, due to the assumptions about $(K,\Sigma)$, we can carry out the
construction. Let $c_k:\pos(w^p,s_0,\ldots,s_{\ell_k},t^{p_k}_{\ell_k+1},
\ldots,t^{p_k}_{\ell_{k+1}})\longrightarrow 2$ be a function such that 
\[c_k(w^k_m)=\left\{
\begin{array}{ll}
0 &\mbox{ if the clause $(\beta)_k$ was applied to choose } p^k_{m+1}\\
1 &\mbox{ otherwise.}
\end{array}
\right.\]
Applying successively $\bar{2}$-bigness (to each $t^{p_k}_i$ for $i=
\ell_{k+1},\ell_{k+1}-1,\ldots,\ell_k+1$) we find creatures $s_{\ell_k+1},
\ldots,s_{\ell_{k+1}}$ such that for each $i\in [\ell_k+1,\ell_{k+1}]$  
\[s_i\in\Sigma(t^{p_k}_i)\quad\mbox{ and }\quad\nor[s_i]\geq\nor[t^{p_k}_i]
-|\pos(w^p,s_0,\ldots,s_{\ell_k},t^{p_k}_{\ell_k+1},\ldots,t^{p_k}_{i-1})|\]
and for each $u\in\pos(w^p,s_0,\ldots,s_{\ell_k},s_{\ell_k+1},\ldots,
s_{\ell_{k+1}})$ the value of $c_k(u)$ depends on $u\restriction m^{s_{
\ell_k}}_{\up}$ only. If the constant value of $c_k\restriction\pos(v,
s_{\ell_k+1},\ldots,s_{\ell_{k+1}})$ (for some sequence $v\in\pos(w^p,s_0,
\ldots,s_{\ell_k})$) is $0$ then easily 
\[(v,s_{\ell_k+1},\ldots,s_{\ell_{k+1}},t^{p^k_{m_k}}_0,t^{p^k_{m_k}}_1,
\ldots)\in \q^*_f(K,\Sigma)\]
(remember that $(K,\Sigma)$ is nice), it is $\leq_0$-stronger than
$(v,t^{p^k}_{\ell_k+1},t^{p^k}_{\ell_k+2},\ldots)$ and it essentially decides
$\dot{\tau}$ (by \ref{smoessdec}). This contradicts the inductive assumption
(4). So the constant value is always $1$ and consequently $c_k\restriction
\pos(w^p,s_0,\ldots,s_{\ell_k},\ldots,s_{\ell_{k+1}})$ is constantly $1$. Put 
\[p_{k+1}\stackrel{\rm def}{=}(w^p,s_0,\ldots,s_{\ell_{k+1}},t^{p^k_{m_k}}_0,
t^{p^k_{m_k}}_1,\ldots)\in\q^*_f(K,\Sigma).\]
Note that for all $i\in\omega$ 
\[\nor[t^{p^k_{m_k}}_i]>f(n+k+1, m^{t^{p^k_{m_k}}_i}_{\dn})+
\fH(m^{t^{p^k_{m_k}}_i}_{\dn})+2\] 
and thus $p_{k+1}$ satisfies (3). By the construction (and \ref{halback}) the
condition $p_{k+1}$ satisfies the inductive requirement (4) (for $k+1$). Since
$f$ is $\bH$-fast we have that for $i\in [\ell_k+1,\ell_{k+1}]$ (by the
inductive assumption (3))
\[\nor[t^{p_{k+1}}_i]\geq\nor[t^{p_k}_i]-|\pos(w^p,s_0,\ldots,s_{\ell_k},
t^{p_k}_{\ell_k+1},\ldots t^{p_k}_{i-1})|>f(n+k,m^{t^{p_{k+1}}_i}_{\dn}),\]
and for $i>\ell_{k+1}$
\[\nor[t^{p_{k+1}}_i]>f(n+k+1,m^{t^{p_{k+1}}_i}_{\dn})\]
and (if $k>0$)
\[\nor[t^{p_k}_{\ell_k}]=\nor[s_{\ell_k}]>f(n+k,m^{t^{p_k}_{\ell_k}}_{\dn}).\]
Hence $p_k\leq_{k}^f p_{k+1}$. Moreover
$\nor[s_{\ell_{k+1}}]>f(n+k+1,m^{s_{\ell_{k+1}}}_{\dn})$. Thus the
requirements (5) and (6) are satisfied.  

Finally look at the limit condition 
\[p^*=(w^p,s_0,s_1,\ldots)=\lim\limits_{k\in\omega} p_k\in\q^*_f(K,\Sigma).\]
Now we have to distinguish the two cases: $(K,\Sigma)$ is simple and
$(K,\Sigma)$ is gluing. 

If our creating pair is simple then we take a condition $r\geq p^*$ which
decides the value of $\dot{\tau}$ and such that $\nor[t^r_i]>0$ for all $i\in
\omega$. We may assume that for some $k\in\omega$ we have $w^r\in\pos(w^p,s_0,
\ldots,s_{\ell_k})$. Then $(w^r, t^{p_k}_{\ell_k+1},t^{p_k}_{\ell_k+2},
\ldots)\leq_0 r$ and this contradicts the assumption (4) about $p_k$. 

Now suppose that $(K,\Sigma)$ is gluing. Note that choosing $\ell_{k+1}$ we
may take it arbitrarily large. So we may assume that (additionally)
$\ell_{k+1}-\ell_k$ is larger then the gluing witness for $f(k,
m^{s_{\ell_k}}_{\up})$. Then we find $s^*_{k+1}\in\Sigma(s_{\ell_k+1},\ldots,
s_{\ell_{k+1}})$ such that
\[\nor[s^*_{k+1}]\geq\min\{f(k,m^{s_{\ell_k}}_{\up}),\nor[s_{\ell_k+1}],
\ldots,\nor[s_{\ell_{k+1}}]\}=f(k,m^{s_{\ell_k}}_{\up}).\]
Put $s^*_0=s_0$ and consider the condition 
\[p^{**}=(w^p,s^*_0,s^*_1,s^*_2,\ldots)\geq p^*.\]
Now we finish choosing the $r$ as earlier above $p^{**}$. The claim is proved. 

\begin{claim}
\label{cl8}
Suppose that $\dot{\tau}$ is a $\q^*_f(K,\Sigma)$-name for an ordinal,
$n<\omega$ and $p\in\q^*_f(K,\Sigma)$. Then there is $q\in\q^*_f(K,\Sigma)$
such that $p\leq^f_n q$ and $q$ essentially decides $\dot{\tau}$. 
\end{claim}

\noindent{\em Proof of the claim}:\ \ \ Take $i_0<\omega$ so large that 
\[(\forall i\geq i_0)(\nor[t^p_i]>f(n+1, m^{t^p_i}_{\dn})).\]
Let $\{w_m:m<m^*\}$ enumerate $\pos(w^p,t^p_0,\dots,t^p_{i_0})$. Choose
inductively $q_m$, $\varepsilon_m$ (for $m<m^*$) such that
\begin{quotation}
\noindent $q_0$ is given by \ref{cl7} for $(w_0,t^p_{i_0+1},t^p_{i_0+2},
\ldots)$,\ \ $\varepsilon_0=1$  

\noindent $\varepsilon_1=\frac{1}{2^{|\pos(w_0,t^p_{i_0+1})|}}$,\ \
$\varepsilon_{m+1}=\frac{\varepsilon_m}{2^{|\pos(w_{m},t^{q_{m-1}}_0)|}}$ (for
$m>0$), 

\noindent $q_{m+1}$ is given by claim \ref{cl7} for $(w_{m+1},t^{q_m}_0,
t^{q_m}_1,\ldots)$, $\varepsilon_{m+1}$.
\end{quotation}
Note that arriving at the stage $m+1<m^*$ of this construction we have 
\[|\pos(w_0,t^p_{i_0+1})|+\sum_{\ell\leq m}|\pos(w_{\ell+1},t^{q_\ell}_0)|\leq
\fH(m^{t^{q_m}_1}_{\dn}),\]
so $\varepsilon_{m+1}$ satisfies the respective demand. Moreover,
\[\nor[t^{q_m}_i]\geq\varepsilon_{m+1}\cdot f(n+1,m^{t^{q_m}_i}_{\dn})>
f(n,m^{t^{q_m}_i}_{\dn})\quad\quad\mbox{ for } i>0,\ m<m^*,\]
and for each $m<m^*$
\[\nor[t^{q_m}_0]\geq\nor[t^p_{i_0+1}]-m> f(n,m^{t^{q_m}_0}_{\dn})+2.\]
Consequently we may carry out the construction and finally letting 
\[q\stackrel{\rm def}{=}(w^p,t^p_0,\ldots,t^p_{i_0},t^{q_{m^*-1}}_0,
t^{q_{m^*-1}}_1,\ldots)\]
we will clearly have a condition as required in the claim.
\medskip

\noindent Applying inductively claim \ref{cl8} to $\dot{\tau}_m$ and $n+m$ we
finish the theorem (using \ref{propord} and \ref{smoessdec}). 
\end{proof}

\begin{corollary}
\label{halbigprop}
Suppose that a creating pair $(K,\Sigma)$ is finitary, $\bar{2}$-big and has
the Halving Property. Further assume that it is either simple or gluing. Let
$f:\omega\times\omega\longrightarrow\omega$ be an $\bH$-fast function. Then
the forcing notion $\q^*_f(K,\Sigma)$ is proper.
\end{corollary}

\section{Tree--creating $(K,\Sigma)$}

\begin{proposition}
\label{fusfront}
Let $e<3$, $n<\omega$, $p\in\q^{\tree}_e(K,\Sigma)$ and let $A\subseteq T^p$ 
be an antichain in $T^p$ such that $(\forall\eta\in A)(\exists\nu\in F^0_n(p)
)(\nu\vartriangleleft\eta)$. Assume that for each $\eta\in A$ we have a
condition $q_\eta\in\q^{\tree}_e(K,\Sigma)$ such that $p^{[\eta]}\leq_0
q_\eta$ and  
\begin{quotation}
\noindent if $e=1$ then $\{t^{q_\eta}_\nu: \nu\in T^{q_\eta}\ \&\ \nor[t^{
q_\eta}_\nu]\leq n\}\subseteq\{t^{p}_\nu:\nu\in T^p\}$.
\end{quotation}
{\em Then} there exists $q\in\q^{\tree}_e(K,\Sigma)$ such that $p\leq^e_{n+1}
q$, $A\subseteq T^q$, $q^{[\eta]}=q_\eta$ for $\eta\in A$ and if $\nu\in T^p$
is such that there is no $\eta\in A$ with $\eta \trianglelefteq\nu$ then
$\nu\in T^q$ and $t^p_\nu=t^q_\nu$.
\end{proposition}

\begin{definition}
\label{kbig}
Let $(K,\Sigma)$ be a tree-creating pair for $\bH$, $k<\omega$.
\begin{enumerate}
\item A tree creature $t\in K$ is called {\em $k$-big} if $\nor[t]>1$ and for
every function $h:\pos(t)\longrightarrow k$ there is $s\in\Sigma(t)$ such that
$h\rest\pos(s)$ is constant and $\nor[s]\geq\nor[t]-1$. 
\item We say that $(K,\Sigma)$ is {\em $k$-big} if every $t\in K$ with
$\nor[t]>1$ is $k$-big.  
\end{enumerate}
\end{definition}

\begin{remark}
The difference with \ref{bigetc} is not serious - we could have interfering.
\end{remark}

\begin{definition}
\label{tomit}
A tree creating pair $(K,\Sigma)$ is t-omittory if for each system $\langle
s_\nu: \nu\in\hat{T}\rangle\subseteq K$ such that $\hat{T}$ is a well founded
quasi tree, $\mrot(s_\nu)=\nu$, $\pos(s_\nu)=\suc_T(\nu)$ (for
$\nu\in\hat{T}$) and for every $\nu_0\in\hat{T}$ such that $\pos(s_{\nu_0})
\subseteq\max(T)$ there is $s\in\Sigma(s_\nu:\nu\in\hat{T})$ such that
\[\nor[s]\geq\nor[s_{\nu_0}]-1\quad\mbox{ and }\quad\pos(s)\subseteq
\pos(s_{\nu_0}).\] 
\end{definition}

\begin{remark}
\label{remtomit}
The name ``t-omittory'' comes from ``tree--omittory'': it is a natural notion
corresponding to omittory creating pairs for the case of tree--creating
pairs. The main point of being t-omittory is that if $p,q\in\q^{\tree}_e(K,
\Sigma)$, $p\leq q$ then we have a condition $r\in\q^{\tree}_e(K,\Sigma)$ such
that $p\leq^e_0 r$ and $\dcl(T^r)\subseteq\dcl(T^q)$ and $t^r_\nu=t^q_\nu$ for
each $\nu\in T^r$, $\mrot(r)\vartriangleleft\nu$. [Why? Let $\eta=\mrot(q)$
and let $T^*\subseteq T^p$ be a well founded quasi tree such that 
\[(\forall\nu\in\hat{T}^*)(\suc_{T^*}(\nu)=\pos(t^p_\nu)),\quad
\mbox{and}\quad\mrot(T^*)=\eta,\quad\mbox{and}\quad t^q_\eta\in\Sigma(
t^p_\nu:\nu\in\hat{T}^*).\]
Let $T^-=\{\mrot(p)\}\cup\bigcup\{\pos(t^p_\nu): \nu\vartriangleleft\eta\ \&\
\nu\in T^p\}\cup \pos(t^q_\eta)$. Clearly $T^-$ is a well founded quasi tree
and we may apply \ref{tomit} to $\langle t^q_\eta,t^p_\nu: \nu\vartriangleleft
\eta\ \&\ \nu\in T^p\rangle$ and $\eta$. Thus we get $t^r_{\mrot(p)}\in
\Sigma(t^q_\eta,t^p_\nu:\nu\vartriangleleft\eta\ \&\ \nu\in T^p)$ such that
$\pos(t^r_{\mrot(p)})\subseteq\pos(t^q_\eta)$. Note that, by transitivity of
$\Sigma$, $t^r_{\mrot(p)}\in\Sigma(t^p_\nu:\nu\in\hat{T}^-\cup\hat{T}^*)$. For
$\nu\in\pos(t^r_{\mrot(p)})$ let $t^r_\nu=t^q_\nu$ and so on. Easily, this
defines a condition $r$ as required.]

Moreover this property implies that $\q^{\tree}_4(K,\Sigma)$ is dense in
$\q^{\tree}_2(K,\Sigma)$. 
\end{remark}

\begin{lemma}
\label{bigfront}
Suppose $(K,\Sigma)$ is a tree--creating pair, $e<3$, $p\in\q^{\tree}_e(K,
\Sigma)$, $n<\omega$ and $\dot{\tau}$ is a $\q^{\tree}_e(K,\Sigma)$-name for
an ordinal. Further assume that if $e=2$ then $(K,\Sigma)$ is bounded (see
\ref{treecreature}(4)). Then: 
\begin{enumerate}
\item There exist a condition $q\in\q^{\tree}_e(K,\Sigma)$ and a maximal
antichain $A\subseteq T^q$ of $T^q$ such that:    
\begin{enumerate}
\item[$(\alpha)$]  $p\leq_{n}^e q$,
\item[$(\beta)$]   for every $\eta\in A$ the condition $q^{[\eta]}$ decides
the value of $\dot{\tau}$, 
\item[$(\gamma)$]  $A$ is an $e$-thick antichain of $T^p$ (see definition
\ref{thick}). 
\end{enumerate}
\item Assume additionally that either $e=0$ and $(K,\Sigma)$ is t-omittory or
$e=1$ and $(K,\Sigma)$ is 2-big. Then there are $q\in\q^{\tree}_e(K,\Sigma)$
and a front $F$ of $T^q$ satisfying clauses $(\alpha)$ and $(\beta)$ of (1)
above. 
\end{enumerate}
\end{lemma}

\begin{proof} 1)\quad Put 
\[\begin{array}{ll}
A_0=\{\nu\in T^p\!: &(\exists \eta\in F^0_n(p))(\eta\vartriangleleft\nu)
\mbox{ and there is }q\in\q^{\tree}_e(K,\Sigma)\mbox{ such that}\\
\ &p^{[\nu]}\leq^e_0 q,\ q\mbox{ decides the value of }\dot{\tau}\ \mbox{
and}\\ 
\ &\mbox{if $e=1$ then } (\forall \eta\in T^q)(\nor[t^q_\eta]>n)\}.\\
  \end{array}\]

\begin{claim}
\label{cl1}
$(\forall r\in\q^{\tree}_e(K,\Sigma))(p\leq r\ \ \Rightarrow\ \ A_0\cap
T^r\neq\emptyset)$.
\end{claim}

\noindent{\em Proof of the claim}: \ \ \ Suppose $r\in\q^{\tree}_e(K,\Sigma)$
is such that $p\leq r$. We may assume that $r$ decides the value of
$\dot{\tau}$ and if $e=1$ then additionally we have $(\forall\eta\in
T^r)(\nor[t^r_\eta]>n)$ (by \ref{frontetc}(3)). Now, as $F^0_n(p)$ is an
$e$-thick antichain of $T^p$ (see \ref{frontetc}(4), (6); remember that if
$e=2$ then we assume that $(K,\Sigma)$ is bounded), we find $\nu\in T^r$ such
that $(\exists\eta\in F^0_n(p))(\eta\vartriangleleft\nu)$. Look at the
condition $r^{[\nu]}$: it witnesses that $\nu\in A_0$. 
\medskip

\noindent Thus 
\[A\stackrel{\rm def}{=}\{\nu\in A_0:\mbox{ there is no }\nu'\vartriangleleft
\nu\mbox{ which is in }A_0\}\]
is an $e$-thick antichain of $T^p$ and for each $\eta\in A$ we may take a
witness $q_\eta$ for $\eta\in A_0$. Now apply \ref{fusfront} to find
$q\in\q^{\tree}_e(K,\Sigma)$ such that $p\leq^e_n q$ and $q^{[\eta]}=q_\eta$
for $\eta\in A$.
\medskip

\noindent 2)\quad If $e=0$ and $(K,\Sigma)$ is t-omittory then for each
$\nu\in F^1_n(p)$ we may choose a condition $q_\nu\geq_0 p^{[\nu]}$
deciding the value of $\dot{\tau}$ (see remark \ref{remtomit}). Now apply
\ref{fusfront}. 
\medskip

\noindent Assume now that $(K,\Sigma)$ is 2-big, $e=1$. Let
\[\begin{array}{ll}
A_1\stackrel{\rm def}{=}\{\nu\in T^p: &(\exists\eta\in F^1_n(p))(\eta
\trianglelefteq\nu)\mbox{ and there is }q\geq_0 p^{[\nu]}\mbox{ such that}\\
\ &\{t^q_\rho:\rho\in T^q\ \&\ \nor[t^q_\rho]\leq n\}\subseteq\{t^p_\rho:\rho
\in T^p\}\mbox{ and}\\
\ &\mbox{there is a front }F\mbox{ of }T^q\mbox{ with }(\forall \rho\in
F)(q^{[\rho]}\mbox{ decides }\dot{\tau})\}.\\ 
  \end{array}\]
Our aim is to show that $F^1_n(p)\subseteq A_1$ which will finish the proof
(applying \ref{fusfront} remember that $F^1_n(p)$ is a front of $T^p$
``above'' $F^0_n(p)$).

Fix $\eta_0\in F^1_n(p)$. For each $\eta_0\trianglelefteq\eta\in T^p$ such
that $\nor[t^p_\eta]>1$ the creature $t^p_\eta$ is 2-big so there is
$s_\eta\in\Sigma(t^p_\eta)$ such that $\nor[s_\eta]\geq\nor[t^p_\eta]-1$ and 
\[\mbox{either }\quad \pos(s_\eta)\cap A_1=\emptyset\quad\mbox{ or }\quad
\pos(s_\eta)\subseteq A_1. \]

\begin{claim}
\label{cl2}
If $\eta_0\trianglelefteq\eta\in T^p$, $\nor[t^p_\eta]\geq n+2$ and
$\pos(s_\eta)\cap A_1\neq\emptyset$ then $\eta\in A_1$. 
\end{claim}

\noindent{\em Proof of the claim}: \ \ \ By the the choice of $s_\eta$ we know
that then $\pos(s_\eta)\subseteq A_1$ so for $\rho\in\pos(s_\eta)$ we may
choose a condition $q_\rho$ and a front $F^\rho\subseteq T^{q_\rho}$
witnessing $\rho\in A_1$. Look at the quasi tree
\[T^{r}\stackrel{\rm def}{=}\{\eta\}\cup\pos(s_\eta)\cup
\bigcup_{\rho\in\pos(s_\eta)}T^{q_\rho}.\]
It determines a condition $r\in\q^{\tree}_1(K,\Sigma)$. It follows from the
assumption $\nor[t^p_\eta]\geq n+2$ that $\nor[s_\eta]>n$ and therefore
\[\{t^r_\nu: \nu\in T^r\ \&\ \nor[t^r_\nu]\leq n\}\subseteq\{t^p_\nu:\nu\in
T^{p^{[\eta]}}\}.\]
Hence the condition $r$ together with $F^\eta=\bigcup\limits_{\rho\in\pos(
s_\eta)} F^\rho$ (which is clearly a front of $T^r$) witness that $\eta\in
A_1$, finishing the proof of the claim. 
\medskip

Now we construct inductively a condition $q$:

$\mrot(T^q)=\eta_0$, $T^q\subseteq T^p$,

if $\nu\in T^q$ and $\nor[t^p_\nu]<n+2$  then $t^q_\nu=t^p_\nu$,
$\suc_{T^q}(\nu)=\pos(t^p_\nu)$, 

if $\nu\in T^q$ and $\nor[t^p_\nu]\geq n+2$  then $t^q_\nu=s_\nu$, and
$\suc_{T^q}(\nu)=\pos(s_\nu)$. 

\noindent It should be clear that $q\in\q^{\tree}_1(K,\Sigma)$ and 
\[\{t^q_\nu:\nu\in T^q\ \&\ \nor[t^q_\nu]\leq n\}\subseteq\{t^p_\nu: \nu\in
T^{p^{[\eta_0]}}\}.\] 

\begin{claim}
\label{cl3}
$B\stackrel{\rm def}{=}\{\eta\in T^q: (\forall\nu\in T^q)(\eta\trianglelefteq
\nu\ \ \Rightarrow\ \ \nor[t^q_\nu]\geq n+2)\}\subseteq A_1$.
\end{claim}

\noindent{\em Proof of the claim}: \ \ \ Suppose that $\eta\in B$. Take a
condition $r\in\q^{\tree}_1(K,\Sigma)$, $q^{[\eta]}\leq r$ which decides
$\dot{\tau}$. We may assume that $(\forall\nu\in T^r)(\nor[t^r_\nu]\geq n+2)$
(by \ref{frontetc}(3)). Consequently $\emptyset=\{t^r_\nu: \nu\in T^r\ \&\
\nor[t^r_\nu]\leq n\}\subseteq\{t^p_\nu: \nu\in T^p\}$ and $\mrot(r)\in A_1$.
But now we note that for each $\nu\in T^q$, if $\eta\trianglelefteq\nu
\vartriangleleft\mrot(r)$ then $\nor[t^p_\nu]\geq n+2$ (as $\eta\in B$) and
$t^q_\nu=s_\nu$. Thus we may apply \ref{cl2} inductively to conclude that all
these $\nu$, including $\eta$, are in $A_1$, finishing the proof of the claim. 

\begin{claim}
\label{cl4}
If $\eta\in T^q$, $\eta\notin A_1$ then there is $\nu\in\pos(t^q_\eta)$ such
that $\nu\notin A_1$.
\end{claim}

\noindent{\em Proof of the claim}: \ \ \ Should be clear.

\begin{claim}
\label{cl5}
$\eta_0\in A_1$.
\end{claim}

\noindent{\em Proof of the claim}: \ \ \ Assume not. Then we inductively
choose a sequence
\[\eta_0\vartriangleleft\eta_1\vartriangleleft\eta_2\vartriangleleft\ldots\in
T^q\] 
such that 
\[(\forall i\in\omega)(\eta_i\notin A_1\ \&\ \nor[t^q_{\eta_{i+1}}]< n+2).\]
For this suppose that we have defined $\eta_i\notin A_1$. Take
$\eta^*\in\pos(t^q_{\eta_i})\setminus A_1$ (possible by \ref{cl4}). By claim
\ref{cl3} we know that $\eta^*\notin B$ so there is $\nu\in (T^q)^{[\eta^*]}$
such that $\nor[t^q_\nu]<n+2$. Let $\eta_{i+1}$ be the shortest such $\nu$,
i.e.~$\eta_{i+1}$ is such that
\[[\eta\in T^q\ \&\ \eta^*\trianglelefteq\eta\vartriangleleft\eta_{i+1}]\
\Rightarrow\ \nor[t^q_\eta]\geq n+2\quad\quad\mbox{ and }\
\nor[t^q_{\eta_{i+1}}]<n+2.\]
By repeating applications of \ref{cl2} we conclude that $\eta_{i+1}\notin
A_1$, as otherwise $\eta^*\in A_1$. 

\noindent Now look at the branch through $T^q$ determined by $\langle\eta_i:
i<\omega\rangle$ -- it contradicts $q\in\q^{\tree}_1(K,\Sigma)$. This finishes
the proof of the claim and the lemma. 
\end{proof}

\begin{theorem}
\label{treedec}
Suppose $(K,\Sigma)$ is a tree--creating pair, $e<3$, $p\in
\q^{\tree}_e(K,\Sigma)$. Further suppose that if $e=2$ then $\Sigma$ is
bounded. Let $n<\omega$ and let $\dot{\tau}_k$ be
$\q^{\tree}_e(K,\Sigma)$-names for ordinals (for $k\in\omega$). Then:
\begin{enumerate}
\item There exist a condition $q\in\q^{\tree}_e(K,\Sigma)$ and $e$-thick
antichains $A_k\subseteq T^q$ of $T^q$ such that for each $k\in\omega$:   
\begin{enumerate}
\item[$(\alpha)$]  $p\leq_{n}^e q$,
\item[$(\beta)$]   for every $\eta\in A_k$ the condition $q^{[\eta]}$ decides
the value of $\dot{\tau}_k$, 
\item[$(\gamma)$]  $(\forall\nu\in A_{k+1})(\exists\eta\in A_k)(\eta
\vartriangleleft\nu)$.
\end{enumerate}
\item If {\em either} $e=0$ and $(K,\Sigma)$ is t-omittory {\em or} $e=1$ and
$(K,\Sigma)$ is 2-big then there are $q\in\q^{\tree}_e(K,\Sigma)$ and fronts
$F_k$ of $T^q$ such that for each $k\in\omega$ the conditions
$(\alpha)$--$(\gamma)$ of (1) above are satisfied. 
\end{enumerate}
\end{theorem}

\begin{proof} This is an inductive application of \ref{bigfront} and
\ref{fusfront} (and \ref{frontetc} + \ref{fusAxA}). 
\end{proof}

\begin{lemma}
\label{4dec}
Let $(K,\Sigma)$ be a tree--creating pair, $\dot{\tau}$ be a $\q^{\tree}_3(K,
\Sigma)$--name for an ordinal, $p\in\q^{\tree}_3(K,\Sigma)$, $n<\omega$. Then
there are a condition $q\in\q^{\tree}_3(K,\Sigma)$ and a $3$--thick antichain
$A\subseteq T^p$ of $T^p$ such that $q\geq p$ and
\begin{enumerate}
\item[$(\alpha)$] $(\forall\eta\in A)(\forall\nu\in T^p)(\nu\vartriangleleft
\eta\ \ \Rightarrow\ \ \nu\in T^q\ \&\ t^q_\nu=t^p_\nu)$,
\item[$(\beta)$]  for every $\eta\in A$ the condition $q^{[\eta]}$ decides the
value of $\dot{\tau}$,
\item[$(\gamma)$] $(\forall\eta\in A)(\forall\nu\in T^q)(\eta\trianglelefteq
\nu\ \ \Rightarrow\ \ \nor[t^q_\nu]>n)$.
\end{enumerate}
\end{lemma}

\begin{proof} Look at the set 
\[\begin{array}{ll}
B\stackrel{\rm def}{=}\{\eta\in T^p:&\mbox{there is a condition }q\geq
p^{[\eta]}\mbox{ such that }\ \mrot(q)=\eta,\\
\ &(\forall\nu\in T^q)(\nor[t^q_\nu]>n)\ \mbox{ and }\ q\mbox{ decides the
name }\dot{\tau}\}. 
  \end{array}\]
Easily, $B\cap T^r\neq\emptyset$ for every condition $r\geq p$. Hence the set 
\[A\stackrel{\rm def}{=}\{\eta\in B: \neg(\exists\nu\in B)(\nu\vartriangleleft
\eta)\}\]
is a $3$--thick antichain of $T^p$. Now we finish in a standard way. 
\end{proof}

\begin{theorem}
\label{4prop}
Let $(K,\Sigma)$ be a tree--creating pair, $\dot{\tau}_k$ be $\q^{\tree}_3(K,
\Sigma)$--names for ordinals (for $k<\omega$) and $p\in\q^{\tree}_3(K,
\Sigma)$. Then there are a condition $q\in\q^{\tree}_3(K,\Sigma)$ and
$3$--thick antichains $A_k$ of $T^p$ such that $q\geq p$ and for every
$k\in\omega$: 
\begin{enumerate}
\item[$(\alpha)$] $A_k\subseteq T^q$,
\item[$(\beta)$]  if $\eta\in A_k$ then $q^{[\eta]}$ decides $\dot{\tau}_k$,
\item[$(\gamma)$] $(\forall\nu\in A_{k+1})(\exists \eta\in
A_k)(\eta\vartriangleleft\nu)$.
\end{enumerate}
\end{theorem}

\begin{proof}
Build the condition $q$ by induction using \ref{4dec}.
\end{proof}

\begin{corollary}
\label{treeproper}
If $e<4$, $(K,\Sigma)$ is a tree creating pair which is bounded if $e=2$, and
$\bigcup\limits_{i\in\omega}\bH(i)$ is countable then the forcing notion
$\q^{\tree}_e(K,\Sigma)$ is proper (and even more). 
\end{corollary}

\begin{proof} By \ref{treedec}, \ref{4prop} (or rather the proofs of them) and
the definition of thick antichains (remember \ref{treesmooth}).
\end{proof}

\begin{lemma}
\label{fronor1}
Assume that $(K,\Sigma)$ is a $2$--big tree--creating pair, $n<\omega$, and
$p\in\q^{\tree}_1(K,\Sigma)$. Then there are $q\in\q^{\tree}_1(K,\Sigma)$ and
a front $F$ of $T^q$ such that $p\leq^1_n q$ and
\[(\forall\nu\in F)(\forall\eta\in T^q)(\nu\trianglelefteq\eta\ \ \
\Rightarrow\ \ \ \nor[t^q_\eta]>n+1).\]
\end{lemma}

\begin{proof}
It is like \ref{bigfront}(2). We consider the set
\[\begin{array}{ll}
A_1^*\stackrel{\rm def}{=}\{\nu\in T^p: &(\exists\eta\in F^1_n(p))(\eta
\trianglelefteq\nu)\mbox{ and there is }q\geq_0 p^{[\nu]}\mbox{ such that}\\
\ &\{t^q_\rho: \rho\in T^q\ \&\ \nor[t^q_\rho]\leq n\}\subseteq\{t^p_\rho:
\rho\in T^p\}\ \mbox{ and}\\
\ &\mbox{there is a front }F\mbox{ of }T^q\mbox{ with }\\
\ &(\forall \rho\in F)(\forall\eta\in T^q)(\rho\trianglelefteq\eta\ \ \ 
\Rightarrow\ \ \ \nor[t^q_\eta] >n+1)\}.\\  
  \end{array}\]
We proceed exactly as in \ref{bigfront}(2) to show that $F^1_n(p)\subseteq
A^*_1$. 
\end{proof}

\begin{corollary}
\label{fronor2}
Suppose that $(K,\Sigma)$ is a $2$-big tree creating pair, $n<\omega$,
$p\in\q^{\tree}_1(K,\Sigma)$. Then there are $q\in\q^{\tree}_1(K,\Sigma)$ and
fronts $F_m$ of $T^q$ (for $m\in\omega$) such that 
\[p\leq^1_n q\quad\mbox{ and }\quad (\forall\eta\in T^q)(\forall\nu\in F_m)
(\nu\trianglelefteq\eta\ \ \ \Rightarrow\ \ \ \nor[t^q_\eta]\geq
m).\]
Hence, in particular, $\q^{\tree}_4(K,\Sigma)$ is dense in $\q^{\tree}_1(K,
\Sigma)$.
\end{corollary}

\section{Examples}
In this part we will give several examples of weak creating pairs, putting
some of the known forcing notions into our setting. It seems that the main
ingredient of any application of our technique is, next to an appropriate
choice of the function $\bH$, the definition of the norm we use to measure
possibilities. Often such a norm is an application of a particular type of
pre--norms. 

\begin{definition}
\label{prenorms}
A function $H:{\mathcal P}(A)\longrightarrow\mbR^{{\geq}0}$ is {\em a pre-norm
on the set $A$} (or rather ${\mathcal P}(A)$) if
\begin{enumerate}
\item[(a)] $B\subseteq C\subseteq A$ implies $H(B)\leq H(C)$
\item[(b)] $H(A)>0$ and if $a\in A$ then $H(\{a\})\leq 1$.
\end{enumerate}
A pre-norm $H$ on $A$ is {\em nice} if additionally
\begin{enumerate}
\item[(c)] if $B\subseteq C\subseteq A$, $H(C)>1$ then either $H(B)\geq
H(C)-1$ or $H(C\setminus B)\geq H(C)-1$.
\end{enumerate}
\end{definition}

\begin{definition}
\label{exprno}
\begin{enumerate}
\item For a non-empty finite set $A$ we let $\dpt^0(A)=|A|$.
\item For a finite family $A\subseteq\fsuo$ such that $(\forall a\in
A)(|a|>1)$ we define $\dpt^1(A)\in\omega$ by the following induction
\[\begin{array}{ll}
\dpt^1(A)\geq 0&\mbox{always,}\\
\dpt^1(A)\geq 1&\mbox{if }A\neq\emptyset,\\
\dpt^1(A)\geq n+2&\mbox{if for every set $X\subseteq\omega$ one of the
following conditions holds:}\\
\ &\dpt^1(\{a\in A: a\subseteq X\})\geq n+1\ \mbox{ or }\\
\ &\dpt^1(\{a\in A:a\subseteq \omega\setminus X\})\geq n+1.
  \end{array}\]
\item For a non-empty finite family $A$ of non-empty subsets of $\omega$ we
let  
\[\dpt^2(A)=\min\{|I|: (\forall a\in A)(a\cap I\neq\emptyset)\}.\]
\item For $n\in\omega$, $i<3$ and $A$ in the domain of $\dpt^i$ we let 
\[\dpt^i_n=\log_{2+n}(\dpt^i(A)).\]
\end{enumerate}
\end{definition}

\begin{proposition}
\label{exareok}
\begin{enumerate}
\item Let $i<3$, $n\in\omega$. Suppose that $A$ is a finite set in the domain
of $\dpt^i$ such that $\dpt^i_n(A)>0$. Then $\dpt^i_n\rest {\mathcal P}(A)$ is
a nice pre-norm on ${\mathcal P}(A)$. 
\item If $H$ is a nice pre-norm on ${\mathcal P}(A)$, $r<H(A)$ is a positive
real number and $H^r:{\mathcal P}(A)\longrightarrow\mbR^{{\geq}0}$ is defined
by 
\[H^r(B)=\max\{0,H(B)-r\}\qquad\mbox{ for }\ B\subseteq A,\]
then $H^r$ is a nice pre-norm on ${\mathcal P}(A)$.
\end{enumerate}
\end{proposition}

\begin{proof} 
1)\quad Note that (in all cases), if $B,C\subseteq A$ then
\[\dpt^i(B\cup C)\leq \dpt^i(B)+\dpt^i(C).\]
The only unclear instance here might be $i=1$, but note that if $\bigcup
A\subseteq [m_0,m_1)$, $B\subseteq A$ then
\begin{quotation}
\noindent $\dpt^1(B)\geq k+2$ \quad if and only if\\
for every partition $\langle I_\ell:\ell<\ell_0\rangle$ of $[m_0,m_1)$,
$\ell_0\leq 2^{k+1}$, there are $b\in B$ and $\ell<\ell_0$ such that $b
\subseteq I_\ell$.
\end{quotation}
2)\quad Check.
\end{proof}

\begin{remark}
\ref{exareok}(2) will be of special importance when defining creating pairs
with the Halving Property. Then we will use $H^r$ for $r=\lfloor\frac{1}{2}
H(A)\rfloor$ (see \ref{bighalex}, \ref{locex} and \ref{BPex}).
\end{remark}

Our first example of a creating pair recalls Blass--Shelah forcing notion
applied in \cite{BsSh:242} to show the consistency of the following statement:
\begin{quotation}
{\em 
\noindent if $\D_1,\D_2$ are non-principal ultrafilters on $\omega$

\noindent then there is a finite-to-one function $f\in\baire$ such that
$f(\D_1)=f(D_2)$.
}
\end{quotation}
(The suitable model was obtained there by a countable support iteration of
forcing notions close to $\q^*_{{\rm s}\infty}(K_{\ref{blsh}}^*,\Sigma^*_{
\ref{blsh}})$ over a model of CH.)

\begin{example}
\label{blsh}
Let $\bH(m)=2$ for $m\in\omega$. We build a full, omittory and omittory--big
(and smooth) creating pair $(K_{\ref{blsh}},\Sigma_{\ref{blsh}})$ for $\bH$.
\end{example}

\begin{proof}[Construction]
A creature $t\in\CR[\bH]$ is in $K_{\ref{blsh}}$ if $m^t_{\dn}+2<m^t_{\up}$
and there is a sequence $\langle A^t_u: u\in\prod\limits_{i<m^t_{\dn}}\bH(i)
\rangle$ such that for every $u\in\prod\limits_{i<m^t_{\dn}}\bH(i)$:
\begin{enumerate}
\item[($\alpha$)] $A^t_u$ is a non-empty family of subsets of $[m^t_{\dn},
m^t_{\up})$, each member of $A^t_u$ has at least $2$ elements,
\item[($\beta$)]  $\langle u,v\rangle\in\val[t]$ \quad if and only if 

$u\vartriangleleft v\in\prod\limits_{i<m^t_{\up}}\bH(i)$ and $\{i\in
[m^t_{\dn}, m^t_{\up}): v(i)=1\}\in A^t_u\cup\{\emptyset\}$,

\item[($\gamma$)] $\nor[t]=\min\{\dpt^1_0(A^t_u):u\in\prod\limits_{i<m^t_{
\dn}}\bH(i)\}$.
\end{enumerate}
[Note that we do not specify here what are the $\dis[t]$ for $t\in
K_{\ref{blsh}}$. We have a total freedom in this, we may allow all possible
values of $\dis[t]$ to appear.]

The composition operation $\Sigma_{\ref{blsh}}$ on $K_{\ref{blsh}}$ is defined
as follows. Suppose that $t_0,\ldots,t_n\in K_{\ref{blsh}}$ are such that
$m^{t_i}_{\up}=m^{t_{i+1}}_{\dn}$ for $i<n$. Then
\begin{quotation}
\noindent $s\in\Sigma_{\ref{blsh}}(t_0,\ldots,t_n)$\quad if and only if

\noindent $s\in K_{\ref{blsh}}$, $m^s_{\dn}=m^{t_0}_{\dn}$, $m^s_{\up}=
m^{t_n}_{\up}$ and for every $\langle u,v\rangle\in\val[s]$ for each $i\leq n$
we have $\langle v\rest m^{t_i}_{\dn},v\rest^{t_i}_{\up}\rangle\in\val[t_i]$.
\end{quotation}
It is an easy exercise to check that $(K_{\ref{blsh}},\Sigma_{\ref{blsh}})$ is
a full, omittory, smooth and omittory--big creating pair (for the last
property use \ref{exareok}). Note that the forcing notion $\q^*_{{\rm s}
\infty}(K_{\ref{blsh}},\Sigma_{\ref{blsh}})$ is non-trivial as for each $m_0<
m_0+2^{n+1}<m_1$ there is $t\in K_{\ref{blsh}}\cap\CR_{m_0,m_1}[\bH]$ such
that $\nor[t]=\log_2(n+1)$. 

One may consider a modification of $(K_{\ref{blsh}},\Sigma_{\ref{blsh}})$
making it forgetful. For this we let $K^*_{\ref{blsh}}=\{t\in
K_{\ref{blsh}}:(\forall u_0,u_1\in\prod\limits_{i<m^t_{\dn}} \bH(i))(A^t_{u_0}
=A^t_{u_1})\}$ and $\Sigma^*_{\ref{blsh}}(t_0,\ldots,t_n)=\Sigma_{
\ref{blsh}}(t_0,\ldots,t_n)\cap K^*_{\ref{blsh}}$ (for suitable $t_0,\ldots,
t_1 \in K^*_{\ref{blsh}}$). Check that
$(K^*_{\ref{blsh}},\Sigma^*_{\ref{blsh}})$ is a forgetful, omittory and
omittory--big creating pair.
\end{proof}

\begin{example}
\label{bighalex}
We define functions $\bH:\omega\longrightarrow\omega$ and $f:\omega\times
\omega\longrightarrow\omega$ and a creating pair $(K_{\ref{bighalex}},\Sigma_{
\ref{bighalex}})$ for $\bH$ such that:
\begin{itemize}
\item $f$ is $\bH$--fast,
\item $(K_{\ref{bighalex}},\Sigma_{\ref{bighalex}})$ is $\bar{2}$--big,
forgetful, simple and has the Halving Property, 
\item the forcing notion $\q^*_f(K_{\ref{bighalex}},\Sigma_{\ref{bighalex}})$
is non--trivial.
\end{itemize}
\end{example}

\begin{proof}[Construction]
Let $F\in\baire$ be an increasing function. Define inductively functions
$\bH=\bH^F$ and $f=f^F$ such that for each $n,k,\ell\in\omega$: 
\begin{enumerate}
\item[(i)] \ $\bH(n)=F(\fH(n))\cdot 2^{f(n,n)}$ (remember $\fH(n)=
\prod\limits_{i<n}|\bH(i)|$,\ $\fH(0)=1$),
\item[(ii)] $f(0,\ell)=\ell+1$,\ \ $f(k+1,\ell)=2^{\fH(\ell)}\cdot
(f(k,\ell)+F(\fH(\ell))\cdot\fH(\ell)+2)$
\end{enumerate}
(note that {\bf (i)}+{\bf (ii)} uniquely determine $\bH$ and $f$ and $f$ is
$\bH$--fast). 

A creature $t\in\CR[\bH]$ belongs to $K_{\ref{bighalex}}$ if $m^t_{\up}=
m^t_{\dn}+1$ and 
\begin{itemize}
\item $\dis[t]=\langle m^t_{\dn}, A_t, H_t\rangle$, where $A_t$ is a subset of
$\bH(m^t_{\dn})$ and $H_t:{\mathcal P}(A_t)\longrightarrow\omega$ is a nice
pre-norm, 
\item $\val[t]=\{\langle u,v\rangle\in\prod\limits_{i<m^t_{\dn}}\bH(i)\times
\prod\limits_{i\leq m^t_{\dn}}\bH(i): u\vartriangleleft v\ \&\ v(m^t_{\dn})\in
A_t\}$,
\item $\nor[t]=H_t(A_t)$.
\end{itemize}
For $t\in K_{\ref{bighalex}}$ let $\Sigma_{\ref{bighalex}}(t)$ consist of all
$s\in K_{\ref{bighalex}}$ such that $m^s_{\dn}=m^t_{\dn}$, $A_s\subseteq A_t$
and $(\forall B\subseteq A_s)(H_s(B)\leq H_t(B))$. (And if $\cS\subseteq
K_{\ref{bighalex}}$, $|\cS|\neq 1$ then we let $\Sigma_{\ref{bighalex}}(\cS)=
\emptyset$.) Next, for $t\in K_{\ref{bighalex}}$ we define $\uhalf(t)\in
K_{\ref{bighalex}}$ as follows:

if $\nor[t]<2$ then $\uhalf(t)=t$,

if $\nor[t]\geq 2$ then $\uhalf(t)\in K_{\ref{bighalex}}$ is (the unique
creature) such that 
\[m^{\uhalf(t)}_{\dn}=m^t_{\dn},\quad A_{\uhalf(t)}=A_t,\quad \mbox{ and}\quad
H_{\uhalf(t)}=(H_t)^r \quad\mbox{ (see \ref{exareok}(2))},\]
where $r=\lfloor\frac{1}{2}\nor[t]\rfloor$. 

\noindent It should be clear that $(K_{\ref{bighalex}},\Sigma_{\ref{bighalex}
})$ is a forgetful, simple and $\bar{2}$--big creating pair (for the last
remember the definition of nice pre-norms). Moreover, the function $\uhalf$
witnesses that $(K_{\ref{bighalex}},\Sigma_{\ref{bighalex}})$ has the weak
Halving Property (and so the Halving Property). [Why? Note that if $s\in
\Sigma_{\ref{bighalex}} (\uhalf(t))$, $\nor[s]>0$, $\nor[t]\geq 2$ then $A_s
\subseteq A_t$ and
\[1\leq H_s(A_s)\leq H_{\uhalf(t)}(A_s)=H_t(A_s)-\lfloor\frac{1}{2}\nor[t]
\rfloor.\]
Thus $H_t(A_s)\geq\lfloor\frac{1}{2}\nor[t]\rfloor +1>\frac{1}{2}\nor[t]\geq
1$. Now look at a creature $t'\in K_{\ref{bighalex}}$ such that $\dis[t']=
\langle m^t_{\dn},A_s,H_t\rest{\mathcal P}(A_s)\rangle$.] Finally note that if
$m<\omega$ and $t\in K_{\ref{bighalex}}$ is such that $\dis[t]=\langle
m,\bH(m),\dpt^0_0\rest {\mathcal P}(\bH(m))\rangle$ then $\nor[t]>f(m,m)$.
\end{proof}

Note that the creating pair $(K,\Sigma)$ described in the Prologue to represent
the Silver forcing $\q$ ``below $2^n$'' is an example of a finitary creating
pair which captures singletons. 
\medskip

The first serious application of tree--creating pairs appeared in
\cite{Sh:326}. The forcing notion ${\mathcal LT}^f_d$ constructed there was
later modified in various ways and several variants of it found their
applications (see e.g.~\cite{BJSh:368}, \cite{FrSh:406} and \ref{frsh}
below). This forcing notion is essentially the $\q^{\tree}_1(K^0_{\ref{326}},
\Sigma^0_{\ref{326}})$. (One should note similarities with the forcing notion
$\q^*_{{\rm s}\infty}(K_{\ref{blsh}},\Sigma_{\ref{blsh}})$.) 

\begin{example}
\label{326}
Let $f\in\baire$ be a strictly increasing function, $\bH(m)=f(m)+1$ for
$m\in\omega$. We construct finitary tree--creating pairs $(K^\ell_{\ref{326}},
\Sigma^\ell_{\ref{326}})$, $\ell<4$, for $\bH$ such that
\begin{enumerate}
\item $K^0_{\ref{326}}=K^1_{\ref{326}}$, $K^2_{\ref{326}}=K^3_{\ref{326}}$,
\item $(K^0_{\ref{326}},\Sigma^0_{\ref{326}})$, $(K^1_{\ref{326}},
\Sigma^1_{\ref{326}})$ are $2$--big local tree--creating pairs,
\item $(K^2_{\ref{326}},\Sigma^2_{\ref{326}})$, $(K^3_{\ref{326}},
\Sigma^3_{\ref{326}})$ are $2$--big t-omittory tree--creating pairs.
\end{enumerate}
\end{example}

\begin{proof}[Construction]
First we define $(K^0_{\ref{326}},\Sigma^0_{\ref{326}})$. A tree creature
$t\in \TCR_\eta[\bH]$ (where $\eta\in\bigcup\limits_{n<\omega}\prod\limits_{i<
n}\bH(i)$) is in $K^0_{\ref{326}}$ if 
\begin{itemize}
\item $\dis[t]=\langle\eta,A_t,H_t\rangle$, where $A_t\subseteq\bH(\lh(\eta))$
and $H_t$ is a nice pre-norm on ${\mathcal P}(A_t)$,
\item $\val[t]=\{\langle\eta,\nu\rangle:\eta\vartriangleleft\nu\ \ \&\ \
\lh(\nu)=\lh(\eta)+1\ \ \&\ \ \nu(\lh(\eta))\in A_t\}$,
\item $\nor[t]=H_t(A_t)$.
\end{itemize}
The operation $\Sigma^0_{\ref{326}}$ is the trivial one and for $t\in\TCR_\eta
[\bH]\cap K^0_{\ref{326}}$:
\[\Sigma^0_{\ref{326}}(t)=\{s\in \TCR_\eta[\bH]\cap K^0_{\ref{326}}: A_s
\subseteq A_t\ \ \&\ \ H_s=H_t\rest {\mathcal P}(A_s)\}.\]
Easily, $(K^0_{\ref{326}},\Sigma^0_{\ref{326}})$ is a finitary $2$--big simple
tree--creating pair.

To define $\Sigma^1_{\ref{326}}$ on $K^1_{\ref{326}}=K^0_{\ref{326}}$ we let
for $t\in\TCR_\eta[\bH]$: 
\[\Sigma^1_{\ref{326}}(t)=\{t\in \TCR_\eta[\bH]\cap K^1_{\ref{326}}: A_s
\subseteq A_t\ \ \&\ \ (\forall B\subseteq A_s)(H_s(B)\leq H_t(B))\}.\]
Plainly, $(K^1_{\ref{326}},\Sigma^1_{\ref{326}})$ is a finitary $2$--big
simple tree--creating pair too.

To have t-omittory variants of the tree-creating pairs defined above we
declare that a tree--creature $t\in\TCR_\eta[\bH]$ is in $K^2_{\ref{326}}=
K^3_{\ref{326}}$ if 
\begin{itemize}
\item $\dis[t]=\langle\eta,\eta^*_t,A_t,H_t\rangle$, where $\eta
\trianglelefteq\eta^*_t\in \bigcup\limits_{n<\omega}\prod\limits_{i<n}\bH(i)$,
$A_t\subseteq\bH(\lh(\eta^*_t))$ and $H_t$ is a nice pre-norm on ${\mathcal
P}(A_t)$, 
\item $\val[t]=\{\langle\eta,\nu\rangle:\eta^*_t\vartriangleleft\nu\ \ \&\ \
\lh(\nu)=\lh(\eta^*_t)+1\ \ \&\ \ \nu(\lh(\eta^*_t))\in A_t\}$,
\item $\nor[t]=H_t(A_t)$.
\end{itemize}
The operations $\Sigma^2_{\ref{326}}$, $\Sigma^3_{\ref{326}}$ are such that
if $T$ is a well founded quasi tree, $\langle s_\nu: \nu\in\hat{T}\rangle$ is
a system of tree--creatures from $K^2_{\ref{326}}$ such that 
\[(\forall \nu\in\hat{T})(s_\nu\in\TCR_\nu[\bH]\ \&\ \rng(\val[s_\nu])=
\suc_T(\nu))\]
then $\Sigma^2_{\ref{326}}(s_\nu:\nu\in \hat{T})$ consists of all
tree-creatures $s\in K^2_{\ref{326}}\cap\TCR_{\mrot(T)}[\bH]$ such that for
some $\nu_0\in\hat{T}$ we have
\[\rng(\val[s])\subseteq\rng(\val[s_{\nu_0}])\subseteq\max(T)\quad\mbox{ and
}\quad H_s=H_{s_{\nu_0}}\rest {\mathcal P}(A_s),\]
and $\Sigma^3_{\ref{326}}(s_\nu:\nu\in \hat{T})$ is defined in a similar manner
but we replace the last demand (on $H_s$) by ``$(\forall B\subseteq A_s)(H_s(B)
\leq H_{s_{\nu_0}}(B))$''.  Now check that $\Sigma^2_{\ref{326}}$,
$\Sigma^3_{\ref{326}}$ are t-omittory $2$--big tree compositions on
$K^2_{\ref{326}}=K^3_{\ref{326}}$.
\end{proof}

Let us finish our overview of ``classical'' examples recalling Fremlin--Shelah
forcing notion. This forcing notion is essentially $\q^{\tree}_1(K_{\ref{frsh}
},\Sigma_{\ref{frsh}})$, and it is a relative of $\q^{\tree}_1(K^1_{\ref{326}
},\Sigma^1_{\ref{326}})$. It was applied in \cite{FrSh:406} to construct a
model in which there is a countable relatively pointwise compact set of
Lebesgue measurable functions which is not stable.   

\begin{example}
\label{frsh}
We build a function $\bH$ and a finitary, local $2$--big tree--creating pair
$(K_{\ref{frsh}},\Sigma_{\ref{frsh}})$.
\end{example}

\begin{proof}[Construction]
Choose inductively increasing sequences $\langle n_k:k<\omega\rangle$ and 
$\langle m_k: k<\omega\rangle$ such that $n_0=m_0=4$, $m_{k+1}>m_k\cdot
2^{n_k}$ and  
\[(n_{k+1})^{(m_{k+1})^{-(k+6)}}\cdot(m_{k+1})^{-(k+6)}>\log_2(n_{k+1}),\quad
n_{k+1}>2^{(m_{k+1})^{k+6}}\cdot (m_{k+1})^{k+6}.\]
For $i\in\omega$ let $\bH(i)=\{a\subseteq n_i: \frac{|a|}{n_i}\geq 1-\frac{1}{
2^{i+2}}\}$.\\
A tree--creature $t\in\TCR_\eta[\bH]$ is taken to $K_{\ref{frsh}}$ if
\begin{itemize}
\item $\dis[t]=\langle\eta,A_t\rangle$, where $A_t\subseteq \bH(\lh(\eta))$ is
non-empty, 
\item $\val[t]=\{\langle\eta,\nu\rangle:\eta\vartriangleleft\nu\ \ \&\ \
\lh(\nu)=\lh(\eta)+1\ \ \&\ \ \nu(\lh(\eta))\in A_t\}$,
\item $\nor[t]=\frac{\lh(\eta)+1}{\log_2(n_{\lh(\eta)})-(\lh(\eta)+2)}\cdot
\dpt^2_0(A_t)$.
\end{itemize}
[Note that $\dpt^2_0(\bH(i))\geq\log_2(n_i)-(i+2)$.] The operation
$\Sigma_{\ref{frsh}}$ is trivial and for $t\in\TCR_\eta[\bH]\cap
K_{\ref{frsh}}$ 
\[\Sigma_{\ref{frsh}}(t)=\{s\in K_{\ref{frsh}}\cap\TCR_\eta[\bH]:A_s\subseteq
A_t\}.\]
Check that $(K_{\ref{frsh}},\Sigma_{\ref{frsh}})$ is a local 2--big
tree--creating pair.
\end{proof}

\begin{remark}
\begin{enumerate}
\item Note that if $\dot{W}$ is the $\q^{\tree}_1(K_{\ref{frsh}},\Sigma_{
\ref{frsh}})$--name for the generic real (see \ref{thereal}), then
\[\begin{array}{ll}
\forces_{\q^{\tree}_1(K_{\ref{frsh}},\Sigma_{\ref{frsh}})}&\mbox{`` }(\forall
i\in\omega)(\dot{W}(i)\subseteq n_i\quad\&\quad \frac{|\dot{W}(i)|}{
n_i}\geq 1-\frac{1}{2^{i+2}})\qquad\mbox{ and}\\
\ &\ \ (\forall\eta\in\V\cap\prod\limits_{i<\omega} n_i)(\forall^\infty i\in
\omega)(\eta(i)\notin \dot{W}(i))\mbox{ ''}.
  \end{array}\]
Thus, after forcing with $\q^{\tree}_1(K_{\ref{frsh}},\Sigma_{\ref{frsh}})$,
the ground model reals are of measure zero.
\item One can define a t-omittory variant of $(K_{\ref{frsh}},\Sigma_{
\ref{frsh}})$ (similarly to the definition of the pair
$(K^2_{\ref{326}},\Sigma^2_{\ref{326}})$; in forcing this would correspond to
considering $\q^{\tree}_0(K_{\ref{frsh}}, \Sigma_{\ref{frsh}})$). 
\item In practical applications, forcing notions of the type $\q^{\tree}_0(K,
\Sigma)$ can be represented in an equivalent form as $\q^{\tree}_1(K^*,
\Sigma^*)$ for some t-omittory pair $(K^*,\Sigma^*)$.
\end{enumerate}
\end{remark}

In the next two examples we want to show that the choice of the type of
forcing notion or the norm condition may be very crucial. Even if we use the
same or very similar weak creating pairs, different approaches may result in
forcing notions with extremely different properties.

\begin{example}
\label{rh3vs1}
There exists a finitary, local and $2$--big tree--creating pair
$(K_{\ref{rh3vs1}},\Sigma_{\ref{rh3vs1}})$ such that the forcing notion
$\q^{\tree}_1(K_{\ref{rh3vs1}},\Sigma_{\ref{rh3vs1}})$ is $\baire$--bounding
but the forcing notion $\q^{\tree}_3(K_{\ref{rh3vs1}},\Sigma_{\ref{rh3vs1}})$
adds an unbounded real.  
\end{example}

\begin{proof}[Construction]
Let $\bH(i)=2^{i+1}$ for $i\in\omega$. A tree--creature $t\in\TCR_\eta[\bH]$
is taken to $K_{\ref{rh3vs1}}$ if 
\begin{itemize}
\item $\dis[t]=\langle\eta,A_t\rangle$ for some $A_t\subseteq\bH(\lh(\eta))$,
\item $\val[t]=\{\langle\eta,\nu\rangle: \eta\vartriangleleft\nu\ \&\ \lh(\nu)
=\lh(\eta)+1\ \&\ \nu(\lh(\eta))\in A_t\}$,
\item $\nor[t]=\dpt^0_0(A_t)$.
\end{itemize}
The operation $\Sigma_{\ref{rh3vs1}}$ is trivial and $\Sigma_{\ref{rh3vs1}}(t)
=\{s\in K_{\ref{rh3vs1}}:\val[s]\subseteq\val[t]\}$.

Plainly, $(K_{\ref{rh3vs1}},\Sigma_{\ref{rh3vs1}})$ is a finitary, local and
$2$--big tree creating pair. By \ref{treedec}(2) we conclude that
$\q^{\tree}_1(K_{\ref{rh3vs1}},\Sigma_{\ref{rh3vs1}})$ is $\baire$--bounding
(compare \ref{treebound}). Note that, by \ref{fronor2}, $\q^{\tree}_4(K_{
\ref{rh3vs1}},\Sigma_{\ref{rh3vs1}})$ is dense in $\q^{\tree}_1(K_{
\ref{rh3vs1}},\Sigma_{\ref{rh3vs1}})$. 

Suppose now that
$p\in\q^{\tree}_3(K_{\ref{rh3vs1}},\Sigma_{\ref{rh3vs1}})$. By induction on
$i<\omega$ we build an increasing sequence $\langle n_i: i<\omega\rangle
\subseteq\omega$, a condition $q\in\q^{\tree}_3(K_{\ref{rh3vs1}},\Sigma_{
\ref{rh3vs1}})$ and a function $f:\{\eta\in T^q: (\exists i<\omega)(\lh(\eta)
=n_i)\}\longrightarrow\omega$ such that 
\begin{enumerate}
\item[($\alpha$)] $q\geq p$ and $\mrot(q)=\mrot(p)$,
\item[($\beta$)]  $n_0=\lh(\mrot(q))$, $f(\mrot(q))=0$,
\item[($\gamma$)] if $\nu,\eta\in T^q$, $\lh(\eta)=n_i$, $\lh(\nu)=n_{i+1}$
and $\eta\vartriangleleft\nu$ then $f(\nu)\in\{f(\eta),f(\eta)+1\}$,
\item[($\delta$)] for each $\eta\in T^q$ such that $\lh(\eta)=n_i$ there is
exactly one $\nu\in T^q$ such that $\eta\vartriangleleft\nu$, $\lh(\nu)=n_{i+
1}$ and $f(\nu)=f(\eta)+1$,
\item[($\varepsilon$)] if $\nu\in T^q$, $n_i\leq\lh(\nu)<n_{i+1}$ then
$\nor[t^q_\nu]=f(\nu\rest n_i)$. 
\end{enumerate}
The construction is quite straightforward. Suppose we have defined $n_i$,
$T^q\cap\omega^{\textstyle {\leq}n_i}\subseteq T^p\cap\omega^{\textstyle
{\leq}n_i}$ and $f\rest T^q\cap\omega^{\textstyle {\leq}n_i}$ in such a way
that 
\[(\forall \eta\in T^q\cap\omega^{\textstyle n_i})(\forall\nu\in T^p)(\eta
\trianglelefteq\nu\ \ \Rightarrow\ \ \nor[t^p_\nu]\geq f(\eta)).\]
By the definition of $\q^{\tree}_3(K_{\ref{rh3vs1}},\Sigma_{\ref{rh3vs1}})$ we
find $n_{i+1}>n_i$ such that for each $\eta\in T^q\cap\omega^{\textstyle n_i}$
there is $\eta^*\in T^p\cap\omega^{\textstyle n_{i+1}}$ such that $\eta
\vartriangleleft\eta^*$ and
\[(\forall \nu\in T^p)(\eta^*\trianglelefteq\nu\ \ \Rightarrow\ \ \nor[
t^p_\nu] \geq f(\eta)+1)\]
(remember $T^q\cap\omega^{\textstyle n_i}$ is finite). Now, for each $\eta\in
T^q\cap\omega^{\textstyle n_i}$ we continue building the condition $q$ above
$\eta$ in such a manner that each $t^q_\nu$ (for $n_i\leq \lh(\nu)<n_{i+1}$)
has norm $f(\eta)$ and the $\eta^*$ is taken to $T^q$. Declare $f(\eta^*)=
f(\eta)+1$ and $f(\eta')=f(\eta)$ for all $\eta'\in T^q\cap \omega^{\textstyle
n_{i+1}}$ extending $\eta$ but different from $\eta^*$.\\
It is easy to check that $q$ built in this manner is a condition in $\q^{
\tree}_3(K_{\ref{rh3vs1}},\Sigma_{\ref{rh3vs1}})$ stronger than $p$.

Note that for each $m\in\omega$ the set $\{\eta\in T^q:\eta\in\dom(f)\ \&\
f(\eta)=m+1\}$ is the $B_m(q)$ (and thus it is a $3$--thick antichain of
$T^q$). We will be done if we show the following claim.

\begin{claim}
\label{cl37}
If $q\leq r\in \q^{\tree}_3(K_{\ref{rh3vs1}},\Sigma_{\ref{rh3vs1}})$ then for
some $m\in\omega$ the set $\{\eta\in T^r: \eta\in\dom(f)\ \&\ f(\eta)=m\}$ is
infinite. 
\end{claim}

\noindent{\em Proof of the claim}: \ \ \ Choose $\eta\in T^r$ such that
\[(\forall\nu\in T^r)(\eta\trianglelefteq\nu\ \ \Rightarrow\ \ \nor[t^r_\nu]
\geq 2)\quad\mbox{ and }\ \ \lh(\eta)=n_i\mbox{ (for some $i\in\omega$)}.\]
Let $m=f(\eta)+1$. Note that if $\nu\in T^r$, $\lh(\nu)=n_j\geq n_i$, $\eta
\trianglelefteq\nu$ and $f(\nu)=f(\eta)$ then:
\begin{enumerate}
\item $|\{\nu^*\in T^r: \nu\vartriangleleft\nu^*\ \&\ \lh(\nu^*)=n_{j+1}\}|
\geq 4$,
\item there is at most one $\nu^*\in T^r$ such that
\[\lh(\nu^*)=n_{j+1},\quad \nu\vartriangleleft\nu^*,\quad\mbox{ and }\quad
f(\nu^*)=m,\]
\item there are $j^*>j$ and $\nu^*\in T^r$ such that 
\[\nu\vartriangleleft\nu^*,\quad \lh(\nu^*)=n_{j^*}\quad\mbox{ and }\quad
f(\nu^*)=m.\]
\end{enumerate}
Hence the set $\{\nu\in T^r:\eta\trianglelefteq\nu\ \&\ \nu\in\dom(f)\ \&\
f(\nu)=m\}$ is infinite. 
\end{proof}

\begin{remark}
Note that the proof that $\q^{\tree}_3(K_{\ref{rh3vs1}},\Sigma_{\ref{rh3vs1}}
)$ adds an unbounded real does not use the specific form of $(K_{\ref{rh3vs1}
},\Sigma_{\ref{rh3vs1}})$. With not much changes we may repeat it for any
local tree creating pair $(K,\Sigma)$ such that 
\begin{enumerate}
\item if $\nu\in\pos(t)$, $t\in K$ and $m\leq\nor[t]$, $m\in\omega$

then there is $s\in\Sigma(t)$ such that $\nor[s]=m$ and $\nu\in\pos(s)$, and

\item if $\nor[t]>2$ then $|\pos(t)|>2$. 
\end{enumerate}
\end{remark}
In the last example of this section we try to show the difference between the
use of tree creating pairs and that of creating pairs.

\begin{example}
\label{two}
Let $\langle n_i:i<\omega\rangle\subseteq\omega$, $f:\omega\times\omega
\longrightarrow\omega$ and $\bH:\omega\longrightarrow\fsuo$ be such that
\begin{enumerate}
\item $f(0,\ell)=\ell+1$,
$f(k+1,\ell)=2^{\fH(\ell)+1}\cdot(f(k,\ell)+\fH(\ell)+2)$, 
\item $n_0=0$, $n_{i+1}\geq (n_i+1)\cdot 2^{2^{2\cdot f(i,i)\cdot\fH(i)^2+2}}
+n_i$,
\item $\bH(i)=[[n_i,n_{i+1})]^{\textstyle n_i+1}$
\end{enumerate}
(so $f$ is $\bH$--fast).\\
We construct weak creating pairs $(K^0_{\ref{two}},\Sigma^0_{\ref{two}})$ and
$(K^1_{\ref{two}},\Sigma^1_{\ref{two}})$ for $\bH$ such that
\begin{enumerate}
\item[(a)] $(K^0_{\ref{two}},\Sigma^0_{\ref{two}})$ is a $2$--big, local and
finitary tree creating pair,
\item[(b)] $(K^1_{\ref{two}},\Sigma^1_{\ref{two}})$ is a simple,
$\bar{2}$--big, finitary and forgetful creating pair with the Halving
Property,
\item[(c)] if $\p$ is either $\q^{\tree}_1(K^0_{\ref{two}},\Sigma^0_{
\ref{two}})$ or $\q^*_f(K^1_{\ref{two}},\Sigma^1_{\ref{two}})$ and $\dot{W}$
is the corresponding name for the generic real (see \ref{thereal}) interpreted
as an infinite subset of $\omega$, then
\[\begin{array}{l}
\forces_{\p}\mbox{`` }(\forall i\in\omega)(|\dot{W}\cap [n_i,n_{i+1})|=n_i+1)
\ \mbox{ and}\\
\ \ \ \ (\forall X\!\in\!\iso\cap\V)(\forall^\infty i\!\in\!\omega)(\dot{W}
\cap [n_i,n_{i+1})\subseteq X\mbox{ or }\dot{W}\cap [n_i,n_{i+1})\cap X=
\emptyset)\mbox{''}.
\end{array}\]
\end{enumerate}
\end{example}

\begin{proof}[Construction]
We try to define minimal forcing notions adding a set $\dot{W}
\subseteq\iso$ with the property stated in the clause {\bf (c)}. The most
natural way is to use weak creatures giving approximations to $\dot{W}$ with
norms related to $\dpt^1$. 

Defining $(K^0_{\ref{two}},\Sigma^0_{\ref{two}})$ we may follow the simplest
possible pattern presented already in \ref{frsh} and \ref{rh3vs1}. So a
tree--creature $t\in\TCR_\eta[\bH]$ is in $K^0_{\ref{two}}$ if
\begin{itemize}
\item $\dis[t]=\langle\eta,A_t\rangle$, where $A_t\subseteq\bH(\lh(\eta))$,
\item $\val[t]=\{\langle\eta,\nu\rangle:\eta\vartriangleleft\nu\ \&\
\lh(\nu)=\lh(\eta)+1\ \&\ \nu(\lh(\eta))\in A_t\}$,
\item $\nor[t]=\frac{\dpt^1_0(A_t)}{2\cdot\fH(\lh(\eta))^2}$.
\end{itemize}
The operation $\Sigma^0_{\ref{two}}$ is trivial (and so
$s\in\Sigma^0_{\ref{two}}(t)$ if and only if $\val[s]\subseteq\val[t]$). One
easily checks that $(K^0_{\ref{two}},\Sigma^0_{\ref{two}})$ is a local 2--big
and finitary tree creating pair. Note that (see the proof of \ref{exareok})
\[\dpt^1(\bH(i))\geq 2^{2\cdot\fH(i)^2\cdot f(i,i)}+1\quad\mbox{ and thus }\
\dpt^1_0(\bH(i))>f(i,i)\cdot 2\cdot\fH(i)^2.\]
Consequently, $\q^{\tree}_1(K^0_{\ref{two}},\Sigma^0_{\ref{two}})$ is a
non-trivial forcing notion. Checking that it satisfies the demand {\bf (c)} is
easy if you remember the definition of $\dpt^1$. 

Now we want to define a creating pair $(K^1_{\ref{two}},\Sigma^1_{\ref{two}})$
in a similar way as $(K^0_{\ref{two}},\Sigma^0_{\ref{two}})$. However, we
cannot just copy the previous case (making suitable adjustments) as we have to
get a new quality: the Halving Property. But we use \ref{exareok}(2) for this. 
Thus a creature $t\in\CR[\bH]$ is taken to $K^1_{\ref{two}}$ if $m^t_{\up}=
m^t_{\dn}+1$ and
\begin{itemize}
\item $\dis[t]=\langle m^t_{\dn},B_t,r_t\rangle$, where $B_t\subseteq\bH(
m^t_{\dn})$ and $r_t$ is a non-negative real,
\item $\val[t]=\{(u,v)\in\prod\limits_{i<m^t_{\dn}}\bH(i)\times\prod\limits_{
i\leq m^t_{\dn}}: u\vartriangleleft v\ \&\ v(m^t_{\dn})\in B_t\}$,
\item $\nor[t]=\max\{0,\frac{\dpt^1_0(B_t)}{2\cdot\fH(m^t_{\dn})^2}-r_t\}$.
\end{itemize}
The operation $\Sigma^1_{\ref{two}}$ is defined by:
\[\Sigma^1_{\ref{two}}(t)=\{s\in K^1_{\ref{two}}: m^t_{\dn}=m^s_{\dn}\ \&\
B_s\subseteq B_t\ \&\ r_s\geq r_t\}.\]
It is not difficult to verify that $(K^1_{\ref{two}},\Sigma^1_{\ref{two}})$ is
a finitary, forgetful, simple and $\bar{2}$--big creating pair (remember
\ref{exareok}(2)). Let $\uhalf:K^1_{\ref{two}}\longrightarrow K^1_{\ref{two}}$
be such that $s=\uhalf(t)$ if and only if $m^s_{\dn}=m^t_{\dn}$, $B_s=B_t$ and
$r_s=r_t+\frac{1}{2}\nor[t]$.\\
We claim that the function $\uhalf$ witnesses the Halving Property for
$(K^1_{\ref{two}},\Sigma^1_{\ref{two}})$. Clearly $\uhalf(t)\in\Sigma(t)$ and
\[\begin{array}{l}
\nor[\uhalf(t)]=\max\{0,\frac{\dpt^1_0(B_t)}{2\cdot\fH(m^t_{\dn})^2}-r_t-
\frac{1}{2}\nor[t]\}=\\
\max\{0,\frac{\dpt^1_0(B_t)}{2\cdot\fH(m^t_{\dn})^2}-r_t-\frac{1}{2}\cdot
\frac{\dpt^1_0(B_t)}{2\cdot\fH(m^t_{\dn})^2}+\frac{1}{2}r_t\}=\\
\frac{1}{2}\max\{0,\frac{\dpt^1_0(B_t)}{2\cdot\fH(m^t_{\dn})^2}-r_t\}=
\frac{1}{2}\nor[t].
  \end{array}\]
Suppose now that $t_0\in K^1_{\ref{two}}$, $\nor[t_0]\geq 2$ and $t\in
\Sigma^1_{\ref{two}}(\uhalf(t_0))$ is such that $\nor[t]>0$. Then $m^t_{\dn}=
m^{t_0}_{\dn}$, $B_t\subseteq B_{t_0}$ and $r_t\geq r_{t_0}+\frac{1}{2}\nor[
t_0]$. Let $s\in\Sigma^1_{\ref{two}}(t_0)$ be such that $B_s=B_t$ and $r_s=
r_{t_0}$. Clearly $\val[s]=\val[t]$ and $\nor[s]=\max\{0,\frac{\dpt^1_0(B_s)}{
2\cdot\fH(m^s_{\dn})^2}-r_{t_0}\}$. But we know that
\[\begin{array}{l}
0<\nor[t]=\max\{0,\frac{\dpt^1_0(B_t)}{2\cdot\fH(m^t_{\dn})^2}-r_t\}=\\
\frac{\dpt^1_0(B_s)}{2\cdot\fH(m^s_{\dn})^2}-r_t\leq\frac{\dpt^1_0(B_s)}{2
\cdot\fH(m^s_{\dn})^2}-r_{t_0}-\frac{1}{2}\nor[t_0]
  \end{array}\]
and hence $\frac{1}{2}\nor[t_0]\leq\nor[s]$.\\
Moreover, by standard arguments, the forcing notion $\q^*_f(K^1_{\ref{two}},
\Sigma^1_{\ref{two}})$ is not trivial and satisfies the demand {\bf (c)}. 

Let us try to show what may distinguish the two forcing notions. We do not
have a clear property of the extensions, but we will present a technical hint
that they may work differently. Let us start with noting the following
property of the pre-norm $\dpt^1$.

\begin{claim}
\label{cl38}
Suppose that $A_0,\ldots,A_{k-1}\subseteq\fsuo$ are finite families of sets
with at least 2 elements, $m>\frac{k(k+3)}{2}$ and $\dpt^1(A_i)\geq m$ for
each $i<k$. Then there is a set $X\subseteq\omega$ such that for each $i<k$ 
\[\begin{array}{l}
\mbox{both}\quad\dpt^1(\{a\in A_i: a\subseteq X\})\geq m-\frac{k(k+3)}{2}\\
\mbox{and }\quad\dpt^1(\{a\in A_i: a\cap X=\emptyset\})\geq m-\frac{k(k+
3)}{2}.
  \end{array}\]
\end{claim}
 
\noindent{\em Proof of the claim}: \ \ \ We prove the claim by induction on
$k$. 

\noindent {\sc Step}\quad $k=1$.\\
We have $A_0\subseteq{\mathcal P}([m_0,m_1))$ such that $\dpt^1(A_0)\geq
m>2$. Take a set $X\subseteq [m_0,m_1)$ of the smallest possible size such
that $\dpt^1(\{a\in A_0: a\subseteq X\})\geq m-1$. Pick any point $n\in X$ and
let $Y=X\setminus\{n\}$. Then
\[\dpt^1(\{a\in A_0: a\subseteq Y\})<m-1\quad\mbox{ and }\quad\dpt^1(\{a\in
A_0:a\subseteq\{n\}\})=0.\]
Hence we get $\dpt^1(\{a\in A_0: a\cap X=\emptyset\})\geq m-2$ (remember the
characterization of $\dpt^1$ from the proof of \ref{exareok}(1)). 
Consequently, the set $X$ is as required. 
\smallskip

\noindent {\sc Step}\quad $k+1$.\\
Suppose $A_0,\ldots,A_{k-1},A_k\subseteq{\mathcal P}([m_0,m_1))$ are such that
$\dpt^1(A_i)\geq m>\frac{(k+1)(k+4)}{2}$. Let $X_0\subseteq [m_0,m_1)$ be such
that for each $i<k$
\[\begin{array}{l}
\dpt^1(\{a\in A_i: a\subseteq X_0\})\geq m-\frac{k(k+3)}{2}\quad\mbox{ and}\\
\dpt^1(\{a\in A_i: a\cap X_0=\emptyset\})\geq m-\frac{k(k+3)}{2}
  \end{array}\]
(exists by the inductive hypothesis). Since $\dpt^1(A_k)\geq m$, one of the
following holds:
\[\dpt^1(\{a\in A_k: a\subseteq X_0\})\geq m-1\quad \mbox{ or }\quad
\dpt^1(\{a\in A_k: a\cap X_0=\emptyset\})\geq m-1.\]
We may assume that the first takes place. Take $X_1\subseteq X_0$ such that
both
\[\dpt^1(\{a\in A_k: a\subseteq X_1\})\geq m-3\ \mbox{ and }\
\dpt^1(\{a\in A_k: a\subseteq X_0\setminus X_1\})\geq m-3.\]
Let $I_0=\{i<k:\dpt^1(\{a\in A_i: a\subseteq X_0\setminus X_1\})\geq m-\frac{k
(k+3)}{2}-1\}$. If $I_0=k$ then we finish this procedure, otherwise we fix
$i_0\in k\setminus I_0$ and we choose a set $X_2\subseteq X_1$ such that 
\[\begin{array}{l}
\dpt^1(\{a\in A_k: a\subseteq X_2\})\geq m-5,\\
\dpt^1(\{a\in A_k: a\subseteq X_1\setminus X_2\})\geq m-5,\quad\mbox{ and}\\
\dpt^1(\{a\in A_{i_0}: a\subseteq X_1\setminus X_2\})\geq m-\frac{k(k+3)}{2}
-2.
  \end{array}\]
We let $I_1=\{i<k:i\in I_0$ or $\dpt^1(\{a\in A_i: a\subseteq X_1\setminus
X_2\})\geq m-\frac{k(k+3)}{2}-2\}$. Note that $I_0\subsetneq I_1$. We continue
in this fashion till we get $I_\ell=k$. Note that this has to happen for some
$\ell\leq k$. Look at the set $X_{\ell+1}$ constructed at this stage. It has
the property that   
\[\begin{array}{l}
\dpt^1(\{a\in A_k: a\subseteq X_{\ell+1}\})\geq m-(2k+1),\\
\dpt^1(\{a\in A_k: a\subseteq X_0\setminus X_{\ell+1}\})\geq m-(2k+1),\quad
\mbox{ and for }i<k\\
\dpt^1(\{a\in A_i: a\subseteq X_0\setminus X_{\ell+1}\})\geq m-\frac{k(k+3)}{2}
-(k+1)\geq m-\frac{(k+1)(k+4)}{2}.
\end{array}\]
So let $X=X_0\setminus X_{\ell+1}$ and check that it is as required for
$A_0,\ldots, A_k$ (and $k+1$). 
\medskip

Now we may show an extra property of $\dot{W}$ which we may get in the case of
$\q^{\tree}_1(K^0_{\ref{two}},\Sigma^0_{\ref{two}})$.
\begin{claim}
\label{cl39}
The following holds in $\V^{\q^{\tree}_1(K^0_{\ref{two}},\Sigma^0_{
\ref{two}})}$:

there are sequences $\langle i_k: k<\omega\rangle$, $\langle
X_i:i<\omega\rangle$ from $\V$  such that
\begin{enumerate}
\item $i_0<i_1<\ldots<\omega$, $X_i\subseteq [n_i,n_{i+1})$,
\item for each $k\in\omega$, for exactly one $i\in [i_k,i_{k+1})$ we have
$\dot{W}\cap [n_i,n_{i+1})\subseteq X_i$,
\item the set $\{i\in\omega:\dot{W}\cap [n_i,n_{i+1})\subseteq X_i\}$ is {\em
not} in $\V$.
\end{enumerate}
[The last demand is to avoid a triviality like $X_i\in\{\emptyset,
[n_i,n_{i+1})\}$.]
\end{claim}

\noindent{\em Proof of the claim}: \ \ \ Let $p\in\q^{\tree}_1(K^0_{
\ref{two}},\Sigma^0_{\ref{two}})$. We may assume that $(\forall\nu\in T^p)(
\nor[t^p_\nu]>4)$. Let $i\geq \lh(\mrot(p))$ and let $\eta\in T^p$,
$\lh(\eta)=i$. Then
\[\dpt^1(A_{t^p_\eta})>2^{8\cdot\fH(i)^2}>\frac{\fH(i)(\fH(i)+3)}{2}.\]
Since $|\{\eta\in T^p:\lh(\eta)=i\}|\leq\fH(i)$, we may use \ref{cl38} to find
a set $X_i\subseteq [n_i,n_{i+1})$ such that for every $\eta\in T^p$ with
$\lh(\eta)=i$ we have both
\[\begin{array}{l}
\dpt^1(\{a\in A_{t^p_\eta}: a\subseteq X_i\})\geq \dpt^1(A_{t^p_\eta})-
\frac{\fH(i)(\fH(i)+3)}{2}\qquad\mbox{ and }\\ 
\dpt^1(\{a\in A_{t^p_\eta}: a\cap X_i=\emptyset\})\geq \dpt^1(A_{t^p_\eta})-
\frac{\fH(i)(\fH(i)+3)}{2}.
  \end{array}\]
Note that
\[\frac{\log_2(\dpt^1(A_{t^p_\eta})-\frac{\fH(i)(\fH(i)+3)}{2})}{2\cdot\fH(
i)^2}\geq \nor[t^p_\eta]-1.\]
Therefore we may inductively build a condition $q\in \q^{\tree}_1(K^0_{
\ref{two}},\Sigma^0_{\ref{two}})$ and a sequence $\langle i_k:
k<\omega\rangle$ such that
\begin{enumerate}
\item $p\leq q$, $\mrot(q)=\mrot(p)$, $\lh(\mrot(q))=i_0<i_1<i_2<\ldots<
\omega$,  
\item for each $\eta\in T^q$, if $\lh(\eta)=i$ then either 
\[A_{t^q_\eta}=\{a\in A_{t^p_\eta}: a\cap X_i=\emptyset\}\ \mbox{ or }\
A_{t^q_\eta}=\{a\in A_{t^p_\eta}: a\subseteq X_i\}\]
(and so $\nor[t^q_\eta]\geq \nor[t^p_\eta]-1$),
\item for each $\eta\in T^q$ with $\lh(\eta)=i_k$ there is exactly one
$i=i(\eta)\in [i_k,i_{k+1})$ such that 
\begin{quotation} 
\noindent if $\eta\vartriangleleft\nu\in T^q$, $\lh(\nu)=i_{k+1}$ and
$\nu_0=\nu\rest i$ 

\noindent then $A_{t^q_{\nu_0}}=\{a\in A_{t^p_{\nu_0}}:a\subseteq X_i\}$,
\end{quotation}
and for distinct $\eta$ as above the values of $i(\eta)$ are distinct.
\end{enumerate}
Now check that the condition $q$ forces that $\langle X_i: i<\omega\rangle$
and $\langle i_k: k<\omega\rangle$ are as required. 
\medskip

\noindent Finally look at \ref{cl39} in the context of the {\em forgetful\/}
creating pair $(K^1_{\ref{two}},\Sigma^1_{\ref{two}})$. 
\end{proof}

\chapter{More properties}
While the properness is the first property we usually ask for when building a
forcing notion, the next request is preserving some properties of ground model
reals. In this chapter we start investigations in this direction dealing with
three properties of this kind. We formulate conditions on weak creating pairs
which imply that the corresponding forcing notions: do not add unbounded
reals, preserve non--null sets or preserve non--meager sets. Applying the
methods developed here we answer Bartoszy\'nski's request (see \cite[Problem
5]{Ba94}), building a proper forcing notion $\p$ which 
\begin{enumerate}
\item preserves non--meager sets, and
\item preserves non--null sets, and
\item is $\baire$-bounding, and
\item does not have the Sacks property.
\end{enumerate}
A forcing notion with these properties is associated with the cofinality of
the null ideal (see \cite{BaJu95}). The construction is done in \ref{tomek1},
\ref{tomek1conc} and it fulfills promise of \cite[7.3A]{BaJu95}. 

\section{Old reals are dominating}
Recall that a forcing notion $\p$ is {\em $\baire$--bounding} if it does not
add unbounded reals, i.e.
\[\forces_{\p}(\forall x\in\baire)(\exists y\in\baire\cap\V)(\forall^\infty
n)(x(n)<y(n)).\]
Any countable support iteration of proper $\baire$--bounding forcing notions
is $\baire$--bounding (see \cite[Ch VI, 2.8A--C, 2.3]{Sh:f}).

\begin{conclusion}
\label{treebound}
Suppose that $(K,\Sigma)$ is a finitary tree-creating pair. 
\begin{enumerate}
\item If $(K,\Sigma)$ is 2-big then the forcing notion
$\q^{\tree}_1(K,\Sigma)$ is $\baire$-bounding.  
\item If $(K,\Sigma)$ is t-omittory then the forcing notion
$\q^{\tree}_0(K,\Sigma)$ is $\baire$-bounding.  
\end{enumerate}
\end{conclusion}

\begin{proof}  
By \ref{treedec}(2) and \ref{frontetc}(5). 
\end{proof}

\begin{conclusion}
\label{gfbound}
Let $(K,\Sigma)$ be a finitary creating pair for $\bH$, and let $f:\omega
\times\omega\longrightarrow\omega$ be an $\bH$-fast function. Suppose that
$(K,\Sigma)$ is $\bar{2}$-big, has the Halving Property and is either
simple or gluing. Then the forcing notion $\q^*_f(K,\Sigma)$ is
$\baire$-bounding.  
\end{conclusion}

\begin{proof}
By \ref{halbigdec}.
\end{proof}

\begin{conclusion}
\label{winbound}
Let $(K,\Sigma)$ be a finitary creating pair which captures singletons.
Then the forcing notion $\q^*_{{\rm w}\infty}(K,\Sigma)$ is $\baire$-bounding.
\end{conclusion}

\begin{proof}
By \ref{sinwin}.
\end{proof}

\section{Preserving non-meager sets}
An important question concerning forcing notions is if ``large'' sets of reals
from the ground model remain ``large'' after the forcing. Here we interpret
``large'' as ``non-meager''. Preserving this property in countable support
iteration is relatively easy. Any countable support iteration of proper
$\baire$--bounding forcing notions which preserve non-meager sets is of the
same type (see \cite[6.3.21, 6.3.22]{BaJu95}). If we omit
``$\baire$--bounding'' then we may consider a condition slightly stronger than
``preserving non-meager sets'':

\begin{definition}
\label{meapres}
Let $\p$ be a proper forcing notion. We say that $\p$ is {\em
Cohen--preserving} if
\begin{enumerate}
\item[$\otimes^{\rm meager}_{\p}$] for every countable elementary submodel $N$
of $(\cH(\chi),{\in},{<^*_\chi})$, a condition $p\in\p$ and a real $x\in\can$
such that\ \ $p,\p,\ldots\in N$ and $x$ is a Cohen real over $N$, {\em there
is} an $(N,\p)$--generic  condition $q\in\p$ stronger than $p$ such that 
\[q\forces_{\p}\mbox{``$x$ is a Cohen real over $N[\Gamma_{\p}]$''}.\]
\end{enumerate}
\end{definition}
In practice, forcing notions preserving non-meagerness of sets from the ground
model are Cohen--preserving. Now, to deal with iterations we may use \cite[Ch
XVIII, 3.10]{Sh:f} (considering $(\bar{R},S,{\bf g})$ as there with ${\bf
g}_a$ being a Cohen real over $a$).

\begin{theorem}
\label{tominonmea}
Suppose that $\bigcup\limits_{i\in\omega}\bH(i)$ is countable and $(K,\Sigma)$
is a t-omittory tree--creating pair. Then the forcing notion
$\q^{\tree}_0(K,\Sigma)$ is Cohen--preserving.
\end{theorem}

\begin{proof} Suppose that $N\prec (\cH(\chi),{\in},{<^*_\chi})$ is countable,
$\bH,K,\Sigma,p,\ldots\in N$, $p\in\q^{\tree}_0(K,\Sigma)$ and $x\!\in\!\can$
is a Cohen real over $N$. Let $\langle\dot{\tau}_n:n<\omega\rangle$, $\langle
\langle\dot{k}^n_i: i<\omega\rangle\!: n<\omega\rangle$, $\langle\langle
\dot{\sigma}^n_i: i<\omega\rangle\!: n<\omega\rangle$ list all $\q^{\tree}_0 
(K,\Sigma)$--names from $N$ for ordinals and sequences of integers and
sequences of finite functions, respectively, such that for each $n<\omega$:
\begin{enumerate}
\item[(a)] $\forces_{\q^{\tree}_0(K,\Sigma)}$``the sequence $\langle
\dot{k}^n_i: i<\omega\rangle$ is strictly increasing'',
\item[(b)] $\forces_{\q^{\tree}_0(K,\Sigma)}$``$(\forall i<\omega)(
\dot{\sigma}^n_i:[\dot{k}^n_i,\dot{k}^n_{i+1})\longrightarrow 2)$''.
\end{enumerate}
Thus each $\langle\dot{k}^n_i,\dot{\sigma}^n_i: i<\omega\rangle$ is essentially
a name for a canonical co-meager set 
\[\{y\in\can: (\exists^\infty i)(y\rest [\dot{k}^n_i, \dot{k}^n_{i+1})=\dot{
\sigma}^n_i)\}.\]  
Of course, the enumerations are not in $N$, but their initial segments are
there.

\begin{claim}
\label{cl20}
Suppose that $q\in\q^{\tree}_0(K,\Sigma)\cap N$, $\nu\in T^q$, $n\in\omega$.
Then there exists a condition $q^n_\nu\in\q^{\tree}_0(K,\Sigma)\cap N$ such
that 
\begin{enumerate}
\item $q^{[\nu]}\leq^0_0 q^n_\nu$, 
\item $\nor[t^{q^n_\nu}]>n$,
\item $q^n_\nu$ decides the values of $\dot{\tau}_m$ for $m\leq n$,
\item for every $m\leq n$:
\[q^n_\nu\forces (\exists i_0<\ldots<i_n<\omega)\big(x\rest [\dot{k}^m_{i_0},
\dot{k}^m_{i_0+1})=\dot{\sigma}^m_{i_0}\ \&\ \ldots\ \&\ x\rest
[\dot{k}^m_{i_n},\dot{k}^m_{i_n+1})=\dot{\sigma}^m_{i_n}\big).\] 
\end{enumerate}
\end{claim}

\noindent{\em Proof of the claim:}\ \ \ First, using \ref{treedec}(2), choose
a condition $q^*\in\q^{\tree}_0(K,\Sigma)\cap N$ such that $q^{[\nu]}\leq
q^*$, $q^*$ decides the values of $\dot{\tau}_m$ for $m\leq n$ and for some
fronts $F^*_j$ of $T^{q^*}$ (for $j<\omega$) we have
\begin{enumerate}
\item[($\alpha$)] $(\forall j\in\omega)(\forall\eta_0\in F^*_{j+1})(\exists
\eta_1\in F^*_j)(\eta_1\vartriangleleft \eta_0)$,
\item[($\beta$)]  for each $m\leq n$, $j\in\omega$, $\eta\in F^*_j$ and $i\leq
j$, the condition $(q^*)^{[\eta]}$ decides the values of $\dot{k}^m_i$,
$\dot{\sigma}^m_i$,
\item[($\gamma$)] $\langle F^*_j: j<\omega\rangle\in N$
\end{enumerate}
(remember that $N$ is an elementary submodel of $(\cH(\chi),{\in},
{<^*_\chi})$). Now take an infinite branch $\rho\in N\cap\lim(T^{q^*})$
through the quasi tree $T^{q^*}$. Take an increasing sequence $\langle \ell_j:
j<\omega\rangle\in N$ of integers such that $\rho\rest\ell_j\in F^*_j$. Then,
by ($\beta$), we have two sequences $\langle k^m_j: m\leq n, j<\omega\rangle$
and $\langle \sigma^m_j: m\leq n, j<\omega\rangle$, both in $N$, such that for
$m\leq n$, $j<\omega$:
\[(q^*)^{[\rho\rest\ell_j]}\forces_{\q^{\tree}_0(K,\Sigma)}\mbox{``} 
\dot{k}^m_j=k^m_j\ \ \&\ \ \dot{\sigma}^m_j=\sigma^m_j\mbox{''.}\]
Look at the set 
\[\{y\in\can: (\forall m\leq n)(\exists^\infty j)(y\rest [k^m_j,k^m_{j+1})
=\sigma^m_j)\}.\]
It is a dense $\Pi^0_2$--set (coded) in $N$ and therefore $x$ belongs to it.
Hence we find $j^*<\omega$ such that
\[(\forall m\leq n)(|\{j<j^*: x\rest [k^m_j,k^m_{j+1})=\sigma^m_j\}|>n).\]
Now take $\eta^*\in T^{(q^*)^{[\rho\rest\ell_{j^*}]}}$ such that
$\nor[t^{q^*}_{\eta^*}]>n+1$. Since $(K,\Sigma)$ is t-omittory we find a
condition $q^n_\nu\in\q^{\tree}_0(K,\Sigma)\cap N$ such that 
\[q^{[\nu]}\leq^0_0 q^n_\nu,\quad \pos(t^{q^n_\nu}_\nu)\subseteq
\pos(t^{q^*}_{\eta^*}),\quad \nor[t^{q^n_\nu}_\nu]>n,\quad\mbox{ and}\]
$t^{q^n_\nu}_\eta=t^{q^*}_\eta$ for all $\eta\in T^{q^n_\nu}$,
$\nu\vartriangleleft\eta$ (compare \ref{remtomit}). Now one easily checks that
$q^n_\nu$ is as required in the claim.
\medskip

We inductively build a sequence $\langle q_n: n<\omega\rangle$, a condition
$q\in\q^{\tree}_0(K,\Sigma)$ and an enumeration $\langle\nu_n: n<\omega
\rangle$ of $T^q$ such that for all $m,n\in\omega$:
\begin{enumerate}
\item $\nu_n\vartriangleleft \nu_m\quad\Rightarrow\quad n<m$,
\item $\langle\nu_n:n<\omega\rangle\subseteq T^p$,
\item $q_n\in\q^{\tree}_0(K,\Sigma)\cap N$, $p^{[\nu_n]}\leq^0_0 q_n$ and if
$\nu_m\vartriangleleft\nu_n$ (so $m<n$) and $\nu_n\in\pos(t^{q_m}_{\nu_m})$
then $(q_m)^{[\nu_n]}\leq^0_0 q_n$,
\item $\nor[t^{q_n}_{\nu_n}]>n$,
\item the condition $q_n$ decides the values of $\dot{\tau}_m$ for $m\leq n$,
\item for every $m\leq n$:
\[q_n\forces(\exists i_0<\ldots<i_n<\omega)\big(x\rest [\dot{k}^m_{i_0},
\dot{k}^m_{i_0+1})=\dot{\sigma}^m_{i_0}\ \&\ \ldots\ \&\ x\rest
[\dot{k}^m_{i_n},\dot{k}^m_{i_n+1})=\dot{\sigma}^m_{i_n}\big),\] 
\item $t^q_{\nu_n}=t^{q_n}_{\nu_n}$. 
\end{enumerate}
The construction is actually described by the conditions above: with a
suitable bookkeeping we build sequences $\langle\nu_n: n\in\omega\rangle$ and
$\langle A_n: n<\omega\rangle\subseteq{\mathcal P}(\omega)$. Arriving at the
stage $n$ of the construction we know $\nu_m$ for $m\leq n$ and $q_m$ for
$m<n$. Applying \ref{cl20} we find $q_n\in\q^{\tree}_0(K,\Sigma)\cap N$ such
that the requirements (3)--(6) above are satisfied. Next we choose
$A_n\subseteq\omega\setminus\bigcup\limits_{m<n} A_m$ of size $|\pos(t^{q_n}_{
\nu_n})|$ and we assign numbers from $A_n$ to elements of $\pos(t^{q_n}_{
\nu_n})$ in such a way that $\pos(t^{q_n}_{\nu_n})=\{\nu_k: k\in A_n\}$, the
set $\omega\setminus\bigcup\limits_{m\leq n} A_m$ is infinite and $\min(\omega
\setminus\bigcup\limits_{m<n} A_m)\in A_n$.

\noindent The $q$ constructed above is a condition in $\q^{\tree}_0(K,\Sigma)$
due to (4), and it is stronger than $p$ by (3) and (7). Clearly $q$ is not in
$N$, but as every finite step of the construction takes place in $N$, the
condition $q$ is $(N,\q^{\tree}_0(K,\Sigma))$--generic (by (5)). Moreover, by
(6),  
\[q\forces_{\q^{\tree}_0(K,\Sigma)}(\forall m<\omega)(\exists^\infty i)(x\rest
[\dot{k}^m_i,\dot{k}^m_{i+1})=\dot{\sigma}^m_i),\]
what implies that
\[q\forces_{\q^{\tree}_0(K,\Sigma)}\mbox{``$x$ is a Cohen real over
$N[\Gamma_{\q^{\tree}_0(K,\Sigma)}]$''.}\]
This finishes the proof.
\end{proof}

The definition \ref{nmp} below was inspired by \ref{tominonmea} and its proof.
We distinguish here the two cases: ``$(K,\Sigma)$ is a creating pair'' and
``$(K,\Sigma)$ is a tree--creating pair'', but in both of them the flavor of
being of the $\NMP$--type is the same: it generalizes somehow the notion of
t-omittory tree--creating pairs. (Note that if $(K,\Sigma)$ is a t-omittory
tree--creating pair then it is of the $\NMP$--type.) One could formulate a
uniform condition here, but that would result in unnecessary complications in
formulation.   

\begin{definition}
\label{nmp}
\begin{enumerate}
\item A finitary creating pair $(K,\Sigma)$ is {\em of the \NMP--type} if the
following condition is satisfied: 
\begin{enumerate}
\item[$(\circledast)_{\NMP}$] \hspace{0.15in} Suppose that $(w,t_0,t_1,\ldots)
\in\q^*_\emptyset(K,\Sigma)$ is such that
\[(\forall k\in\omega)(\nor[t_k]>\fH(m^{t_k}_{\dn}))\]
and let $n_0<n_1<n_2<\ldots<\omega$. Further, assume that 
\[g:\bigcup_{i\in\omega}\pos(w,t_0,\ldots,t_{n_i-1})\longrightarrow
\bigcup_{i\in\omega}\pos(w,t_0,\ldots,t_{n_{i+1}-1})\]
is such that $g(v)\in\pos(v,t_{n_i},\ldots,t_{n_{i+1}})$ for $v\in\pos(w,t_0,
\ldots,t_{n_i-1})$ (so $v\vartriangleleft g(v)$). {\em Then} there are
$0<i<\omega$ and a creature $s\in\Sigma(t_{n_0},\ldots,t_{n_i-1})$ such that  
\begin{enumerate}
\item[$(\alpha)$] $\nor[s]\geq\min\{\nor[t_m]-\fH(m^{t_m}_{\dn}): n_0\leq m<
n_i\}$,    
\item[$(\beta)$]  for each $v\in\pos(w,t_0,\ldots,t_{n_0-1},s)$ there is $j<i$
such that  
\[v\rest m^{t_{n_{j+1}-1}}_{\up}=g(v\rest m^{t_{n_j-1}}_{\up}).\]
\end{enumerate}
\end{enumerate}
\item A tree--creating pair $(K,\Sigma)$ is {\em of the $\NMP$--type} if the
following condition is satisfied:
\begin{enumerate}
\item[$(\circledast)^{\tree}_{\NMP}$] Suppose that $\langle t_\eta:\eta\in
T\rangle\in \q^{\tree}_{\emptyset}(K,\Sigma)$ (see \ref{treeforcing}(4)) is
such that 
\[(\forall \eta\in T)(\nor[t_\eta]>1)\]
and let $F_0,F_1,F_2,\ldots$ be fronts of the quasi tree $T$ such that for
$i<\omega$: 
\[(\forall \nu\in F_{i+1})(\exists\nu'\in F_i)(\nu'\vartriangleleft\nu).\]
Assume that a function $g:\bigcup\limits_{i\in\omega} F_i\longrightarrow
\bigcup\limits_{i\in\omega} F_{i+1}$ is such that $\nu\vartriangleleft g(\nu)
\in F_{i+1}$ provided $\nu\in F_i$. {\em Then} there are $0<i<\omega$ and a
tree--creature $s\in\Sigma(t_\eta:(\exists\nu\in F_i)(\eta\vartriangleleft
\nu))$ such that 
\begin{enumerate}
\item[$(\alpha)^{\tree}$] $\nor[s]\geq\inf\{\nor[t_\eta]-1: \eta\in T\}$,
\item[$(\beta)^{\tree}$]  for each $\nu\in\pos(s)$ there are $j<i$ and $k<
\lh(\nu)$ such that 
\[\nu\rest k\in F_j\quad\mbox{ and }\quad g(\nu\rest k)\trianglelefteq\nu.\]
\end{enumerate}
\end{enumerate}
\end{enumerate}
\end{definition}

\begin{theorem}
\label{nonmeager}
Assume $(K,\Sigma)$ is a finitary, gluing and $\bar{2}$--big creating pair. 
Suppose that $(K,\Sigma)$ is of the $\NMP$--type and has the Halving
Property. Let $f:\omega\times\omega\longrightarrow\omega$ be an $\bH$-fast
function. Then the forcing notion $\q^*_f(K,\Sigma)$ is Cohen--preserving.
\end{theorem}

\begin{proof} By \ref{gfbound} we know that the forcing notion $\q^*_f(K,
\Sigma)$ is $\baire$-bounding. Consequently it is enough to show that if
$A\subseteq\can$ is a non-meager set then
\[\forces_{\q^*_f(K,\Sigma)}\mbox{``}A\mbox{ is not meager''}\]
(see \cite[6.3.21]{BaJu95}). So suppose that $A\subseteq\can$ is not meager
but some condition in $\q^*_f(K,\Sigma)$ forces that this set is meager. Thus
we find a condition $p_0\in\q^*_f(K,\Sigma)$ and $\q^*_f(K,\Sigma)$--names
$\dot{k}_n, \dot{\sigma}_n$ such that 
\[\begin{array}{ll}
p_0\forces_{\q^*_f(K,\Sigma)} &\mbox{`` }\dot{k}_0<\dot{k}_1<\ldots<\omega\ 
\mbox{ and }\ \dot{\sigma}_n: [\dot{k}_n,\dot{k}_{n+1})\longrightarrow 2\
\mbox{ (for $n<\omega$) }\ \ \mbox{ and}\\ 
\ &\ \ (\forall x\in A)(\forall^\infty n)(x\rest [\dot{k}_n,\dot{k}_{n+1})\neq
\dot{\sigma}_n)\mbox{''}.
  \end{array}\]
As $\q^*_f(K,\Sigma)$ is $\baire$-bounding, we find a condition $p_1\geq p_0$,
a sequence $0=k_0<k_1<k_2<\ldots<\omega$ and names $\dot{\rho}_n$ such that
\[\begin{array}{ll}
p_1\forces_{\q^*_f(K,\Sigma)} &\mbox{``}\dot{\rho}_n: [k_n,k_{n+1})
\longrightarrow 2\ \mbox{ (for $n<\omega$)\ \ and }\\
\ &\ \ (\forall x\in A)(\forall^\infty n)(x\rest [k_n,k_{n+1})\neq
\dot{\rho}_n)\mbox{''}. 
  \end{array}\]
Further, applying \ref{halbigdec} we find $p_2\geq p_1$ such that $p_2$
essentially decides all the names $\dot{\rho}_n$. Clearly we may assume that
$\nor[t^{p_2}_n]>f(2,m^{t^{p_2}_n}_{\dn})$ for all $n<\omega$ (and thus
$\nor[t^{p_2}_n]>\fH(m^{t^{p_2}_n}_{\dn})$). Choose $0=n_0<n_1<\ldots<\omega$
and $\ell_0<\ell_1<\ell_2<\ldots<\omega$ such that 
\[(\forall m\in\omega)(\fH(m^{t^{p_2}_{n_m}}_{\dn})<\ell_{m+1}-\ell_m)\] 
and for each $m<\omega$ and every sequence $w\in\pos(w^{p_2},t^{p_2}_0,\ldots,
t^{p_2}_{n_{m+1}-1})$ the condition $(w,t^{p_2}_{n_{m+1}},t^{p_2}_{n_{m+1}+1},
\ldots)$ decides all the names $\dot{\rho}_j$ for $j\in [\ell_m,\ell_{m+1})$
(remember the choice of $p_2$). Let 
\[g:\bigcup_{m\in\omega}\pos(w^{p_2},t^{p_2}_0,\ldots,t^{p_2}_{n_m-1})
\longrightarrow\bigcup_{m\in\omega}\pos(w^{p_2},t^{p_2}_0,\ldots,t^{p_2}_{n_{
m+1}-1})\]
be such that $v\vartriangleleft g(v)\in\pos(w^{p_2},t^{p_2}_0,\ldots,t^{p_2}_{
n_{m+1}-1})$ for $v\in\pos(w^{p_2},t^{p_2}_0,\ldots,t^{p_2}_{n_m-1})$. Next,
for each $v\in\pos(w^{p_2},t^{p_2}_0,\ldots,t^{p_2}_{n_m-1})$, $m<\omega$ fix
$\ell(v)\in [\ell_m,\ell_{m+1})$ and $\rho(v)$ such that
\begin{itemize}
\item there are no repetitions in $\langle \ell(v)\!: v\in\pos(w^{p_2},
t^{p_2}_0,\ldots,t^{p_2}_{n_m-1})\rangle$,
\item $\rho(v): [k_{\ell(v)},k_{\ell(v)+1})\longrightarrow 2$ is such that 
\[(g(v),t^{p_2}_{n_{m+1}},t^{p_2}_{n_{m+1}+1},\ldots)\forces_{\q^*_f(K,\Sigma)}
\mbox{`` }\dot{\rho}_{\ell(v)}=\rho(v)\mbox{ ''.}\]
\end{itemize}
Now we apply successively \ref{nmp}(1) to the condition $(w^{p_2},t^{p_2}_0,
t^{p_2}_1,\ldots)$, the sequence $\langle n_i:i<\omega\rangle$ and the mapping
$g$. As a result we construct an increasing sequence $0=i_0<i_1<i_2<\ldots<
\omega$ of integers and creatures $s_j\in\Sigma(t^{p_2}_{n_{i_j}},\ldots,
t^{p_2}_{n_{i_{j+1}}-1})$ such that for all $j<\omega$:
\begin{enumerate}
\item $\nor[s_j]\geq\min\{\nor[t^{p_2}_{n_{i_j}}]-\fH(m^{t^{p_2}_{n_{i_j}}}_{
\dn}), \ldots,\nor[t^{p_2}_{n_{i_{j+1}}-1}]-\fH(m^{t^{p_2}_{n_{i_{j+1}}-1}}_{
\dn})\}$,
\item for each $v\in\pos(w^{p_2},t^{p_2}_0,\ldots,t^{p_2}_{n_{i_j-1}},s_j)$
there is $i^*\in [i_j,i_{j+1})$ such that $v\rest m^{t^{p_2}_{n_{i^*+1}}}_{
\dn}=g(v\rest m^{t^{p_2}_{n_{i^*}}}_{\dn})$.  
\end{enumerate}
It should be clear that $(w^{p_2},s_0,s_1,\ldots)\in\q^*_f(K,\Sigma)$ (note
that if $\nor[t^{p_2}_n]>f(k+1,m^{t^{p_2}_n}_{\dn})$ for all $n\in
[n_{i_j},n_{i_{j+1}})$ then $\nor[s_j]>f(k,m^{s_j}_{\dn})$). Look at the set 
\[\begin{array}{ll}
A^*\stackrel{\rm def}{=}\{x\in\can:& \mbox{for infinitely many }j\in\omega,
\mbox{ for every }i^*\in [i_j,i_{j+1})\\
\ &(\forall v\in\pos(w^{p_2},t^{p_2}_0,\ldots,t^{p_2}_{n_{i^*}-1}))(x\rest
[k_{\ell(v)},k_{\ell(v)+1})=\rho(v))\}. 
  \end{array}\]
It is a dense $\Pi^0_2$--subset of $\can$ and hence $A^*\cap A\neq\emptyset$.
Take $x\in A^*\cap A$. Note that the choice of the $s_j$'s (see clause $2.$
above) implies that for each $j<\omega$ and $v\in\pos(w^{p_2},s_0,\ldots,s_j)$
there is $i^*\in [i_j,i_{j+1})$ such that letting $v'=v\rest m^{t^{p_2}_{
n_{i^*}}}_{\dn}$ we have
\[(v,s_{j+1},s_{j+2},\ldots)\forces_{\q^*_f(K,\Sigma)}\mbox{`` }\dot{\rho}_{
\ell(v')}=\rho(v')\mbox{ ''.}\]
Hence $(w^{p_2},s_0,s_1,\ldots)\forces_{\q^*_f(K,\Sigma)}$`` $(\exists^\infty
n\in\omega)(x\rest [k_n,k_{n+1}) = \dot{\rho}_n)$ '', a contradiction. 
\end{proof}

\begin{theorem}
\label{treenonmea}
Assume that $(K,\Sigma)$ is a finitary $2$-big tree--creating pair of the
$\NMP$--type. Then the forcing notion $\q^{\tree}_1(K,\Sigma)$ is
Cohen--preserving. 
\end{theorem}

\begin{proof} Like \ref{nonmeager} but using \ref{treebound}, \ref{treedec}(2),
\ref{fronor2}, and \ref{nmp}(2): we choose $p_0$, $\dot{k}_n$, $\dot{
\sigma}_n$, $p_1$, $k_n$ and $\dot{\rho}_n$ as there. Further we inductively
build an increasing sequence $\langle\ell_m: m\in\omega\rangle$ of integers, a
condition $p_2\geq p_1$ and fronts $F_m$ of $T^{p_2}$ (for $m<\omega$) such
that $|F_m|<\ell_{m+1}-\ell_m$ and for every $\eta\in F_{m+1}$
\[\eta\trianglelefteq\nu\in T^{p_2}\ \Rightarrow\ \nor[t^{p_2}_\nu]\geq m+1\ \
\mbox{ and }\ \ p^{[\eta]}_2\mbox{ decides all }\dot{\rho}_j\mbox{ for }j\in
[\ell_m,\ell_{m+1}),\big)\]  
and the front $F_{m+1}$ is above $F_m$. For $\nu\in F_m$ we define $g(\nu)\in
F_{m+1}$, $\ell(\nu)\in [\ell_m,\ell_{m+1})$ and $\rho(\nu):[k_{\ell(\nu)},
k_{\ell(\nu)+1})\longrightarrow 2$ in a manner parallel to that in the proof
of \ref{nonmeager}. Next we build a condition $q\geq p_2$ and an increasing
sequence $\langle m_i: i\in\omega\rangle$ such that each $F_{m_i}\cap T^q$ 
is a front of $T^q$ and
\begin{quotation}
\noindent if $\eta\in F_{m_i}\cap T^q$, $i\in\omega$ then $\pos(t^q_\eta)
\subseteq \bigcup\{F_m: m_i<m\leq m_{i+1}\}$ and for every $\nu\in\pos(
t^q_\eta)$ there are $j$ and $k<\lh(\nu)$ such that $\nu\rest k\in F_j$ and
$g(\nu\rest k)\trianglelefteq\nu$. 
\end{quotation}
Finally we let
\[A^*\stackrel{\rm def}{=}\{x\in\can: (\exists^\infty i\in\omega)(\forall m\in
[m_i,m_{i+1}))(\forall\nu\in F_m)(x\rest [k_{\ell(\nu)},k_{\ell(\nu)+})=\rho
(\nu))\}\]
and we finish as in \ref{nonmeager}. 
\end{proof}

\begin{theorem}
\label{winnonmea}
Suppose that $(K,\Sigma)$ is a finitary creating pair which captures
singletons. Then the forcing notion $\q^*_{{\rm w}\infty}(K,\Sigma)$ is 
Cohen--preserving. 
\end{theorem}

\begin{proof} Like \ref{nonmeager}, but using \ref{winbound}, \ref{sinwin} and
\ref{singleton}. Note that the last implies that the pair $(K,\Sigma)$ has the
following property:
\begin{enumerate}
\item[($\circledcirc$)] {\em If} $(t_0,\ldots,t_N)\in\PFC(K,\Sigma)$, $\fH(
m^{t_0}_{\dn})\leq k$, $0=n_0<n_1<\ldots<n_k<N$ and $u\in\pos(u\rest m^{t_0}_{
\dn},t_0,\ldots,t_N)$\quad {\em then} there are $s_0,\ldots,s_\ell\in K$ such
that $\pos(u\rest m^{s_0}_{\dn}, s_0,\ldots,s_\ell)=\{u\}$, $m^{s_0}_{\dn}=
m^{t_0}_{\dn}$, $m^{s_\ell}_{\up}=m^{t_N}_{\up}$ and $\langle t_0,\ldots,t_N
\rangle\leq\langle s_0,\ldots,s_\ell\rangle$. \quad {\em Consequently,}
choosing an enumeration $\{w_j: j<\fH(m^{t_0}_{\dn})\}$ of $\basis(t_0)$ 
and letting $u_j=w_j\conc u\rest [m^{t_0}_{\dn}, m^{t_{n_{j+1}-1}}_{\up})$ we
will have (remember $(K,\Sigma)$ is forgetful) 
\begin{enumerate}
\item[($\alpha$)] $u_j\in\pos(w_j,t_0,\ldots,t_{n_{j+1}-1})$,
\item[($\beta$)]  $(t_0,\ldots,t_N)\leq(s_0,\ldots,s_\ell)$, 
\item[($\gamma$)] for each $w\in\basis(s_0)$ and $v\in\pos(w,s_0,\ldots,
s_\ell)$ there is $j\!<\!\fH(m^{t_0}_{\dn})$ such that $v\rest m^{t_{n_{j+1}
-1}}_{\up}=u_j$. 
\end{enumerate}
\end{enumerate}
Thus we may repeat the proof of \ref{nonmeager} with not many changes.
\end{proof}

Let us note that it is not an accident that we have the results on preserving
of non-meager sets only for forcing notions which are very much like the ones
coming from t-omittory tree--creating pairs. If we look at the opposite pole:
local weak creating pairs, then we notice that they easily produce forcing
notions making the ground reals meager.

\begin{definition}
\label{trmea}
Let $\bH$ be of countable character. We say that a weak creating pair
$(K,\Sigma)$ for $\bH$ is {\em trivially meagering} if for every $t\in K$ with
$\nor[t]>1$ and each $u\in\basis(t)$ and $v\in\pos(u,t)$ there is $s\in
\Sigma(t)$ such that $\nor[s]\geq\nor[t]-1$ and $v\notin\pos(u,s)$.
\end{definition}

\begin{proposition}
\label{trmamea}
\begin{enumerate}
\item If $(K,\Sigma)$ is a local trivially meagering tree--creating pair for
$\bH$, $\bH$ is of countable character and $e=1,3$ then 
\[\forces_{\q^{\tree}_e(K,\Sigma)}\mbox{`` }\baire\cap\V\mbox{ is meager
''.}\] 
\item If $(K,\Sigma)$ is a simple finitary and trivially meagering creating
pair for $\bH$ and $f:\omega\times\omega\longrightarrow\omega$ is $\bH$--fast
then 
\[\forces_{\q^*_f(K,\Sigma)}\mbox{`` }\baire\cap\V\mbox{ is meager ''.}\]
\end{enumerate}
\end{proposition}

\begin{proof} In the first case remember that if $p\in\q^{\tree}_e(K,\Sigma)$,
$e=1,3$ then for some $\nu\in T^p$ we have $(\forall\eta\in
T^p)(\nu\trianglelefteq\eta\ \Rightarrow\ \nor[t^p_\eta]\geq 2)$. Hence, as
$(K,\Sigma)$ is local and trivially meagering we easily get
\[\forces_{\q^{\tree}_e(K,\Sigma)}\mbox{`` }(\forall x\in\prod_{m\in\omega}
\bH(m)\cap\V)(\forall^\infty m\in\omega)(\dot{W}(m)\neq x(m))\mbox{ ''.}\]
For the second case suppose that $p\in\q^*_f(K,\Sigma)$. Since $f(n+1,\ell)>
\fH(\ell)+f(n,\ell)$  using the assumptions that $(K,\Sigma)$ is simple and
trivially meagering we immediately see that
\[p\forces_{\q^*_f(K,\Sigma)}\mbox{`` }(\forall x\in\prod_{m\in\omega}
\bH(m)\cap\V)(\forall^\infty n\in\omega)(\dot{W}\rest [m^{t^p_n}_{\dn},
m^{t^p_n}_{\up})\neq x\rest [m^{t^p_n}_{\dn},m^{t^p_n}_{\up}))\mbox{
'',}\] 
finishing the proof. 
\end{proof}

Later, in \ref{makemeager}, we will see that the forcing notions $\q^*_{{\rm
s}\infty}(K,\Sigma)$ may make the ground model reals meager too. This suggests
that if one wants to build a forcing notion preserving non-meagerness then the
most natural approach is $\q^{\tree}_e$ for $e=0,2$ or $\q^*_{{\rm w}\infty}$.

\section{Preserving non-null sets}
In this section we introduce a property of tree--creating pairs which implies
that forcing notions $\q^{\tree}_e(K,\Sigma)$ preserve non-null sets. Though
preserving non-null sets alone is not enough to use the preservation theorem
of \cite[Ch. XVIII, \S3]{Sh:f}, one may apply the methods of \cite{Sh:630},
\cite{Sh:669} when dealing with countable support iterations of forcing
notions of the type presented here.

\begin{definition}
\label{tnnp}
We say that a weak creating pair $(K,\Sigma)$ is of the $\NNP$--type if the
following condition is satisfied: 
\begin{enumerate}
\item[$(\circledast)_{\NNP}$]\hspace{1.5cm} there are increasing sequences
$\bar{a}=\langle a_n: n<\omega\rangle\subseteq (0,1)$ and $\bar{k}=\langle
k_n: n\in\omega\rangle\subseteq\omega$ such that $\lim\limits_{n\to\infty}
a_n<1$ and:

\noindent{\em if} $t\in K$, $\nor[t]>2$, $v\in\basis(t)$, $\lh(v)\geq k_n$,
$N<\omega$ and a function
\[g:\pos(v,t)\longrightarrow {\mathcal P}(N)\]
is such that $(\forall u\in\pos(v,t))(a_{n+1}\leq\frac{|g(u)|}{N})$ {\em then}
the set 
\[\begin{array}{ll}
\{n<N: &\mbox{there is }s\in\Sigma(t)\mbox{ such that }\nor[s]\geq\nor[t]-1
\mbox{ and}\\ 
\ &\basis(t)\subseteq\basis(s)\ \mbox{ and }\ (\forall u\in\pos(v,s))(n\in
g(u))\} 
\end{array}\]
has not less than $N\cdot a_n$ elements.
\end{enumerate}
(The sequences $\bar{a},\bar{k}$ from $(\circledast)_{\NNP}$ are said to
witness $\NNP$.) 
\end{definition}

\begin{definition}
\label{gluing}
We say that a tree--creating pair $(K,\Sigma)$ for $\bH$ is {\em gluing}
(respectively: {\em weakly gluing}) if for each well founded quasi tree 
$T\subseteq\bigcup\limits_{n<\omega}\prod\limits_{i<n}\bH(i)$ and a system
$\langle s_\nu:\nu\in\hat{T}\rangle\subseteq K$ such that 
\[(\forall\nu\in\hat{T})(s_\nu\in\TCR_\nu[\bH]\quad\&\quad\pos(s_\nu)=
\suc_T(\nu))\] 
there is $s\in\Sigma(s_\nu:\nu\in\hat{T})$ such that
\[\nor[s]\geq\sup\{\nor[s_\nu]-1: \nu\in\hat{T}\}\]
($\nor[s]\geq\inf\{\nor[s_\nu]-1: \nu\in\hat{T}\}$, respectively).
\end{definition}

\begin{remark}
The above definition, though different from \ref{simpglui}(2) (for creating
pairs),  has actually the same meaning: we may glue together creatures without
loosing too much on norms.
\end{remark}

\begin{definition}
\label{strfin}
We say that a weak creating pair $(K,\Sigma)$ is {\em strongly finitary} if
$K$ is finitary (see \ref{morecreat}(2)) and  $\Sigma(\cS)$ is finite for each
$\cS\subseteq K$. If $\sim_\Sigma$ (see \ref{subcom}(3)) is an equivalence
relation on $K$ and $\Sigma$ depends on $\sim_\Sigma$--equivalence classes
only, then what we actually require is that $\Sigma(\cS)/{\sim_\Sigma}$ is
finite. 
\end{definition}

\begin{theorem}
\label{treenonnull}
Suppose $(K,\Sigma)$ is a strongly finitary tree--creating pair of the
$\NNP$--type. Further suppose that:
\begin{quotation}
\noindent either $(K,\Sigma)$ is t-omittory and $e=0$

\noindent or $(K,\Sigma)$ is $2$-big and weakly gluing and $e=1$.
\end{quotation}
{\em If} $A\subseteq\can$ is a non-null set {\em then} $\forces_{\q^{\tree}_e
(K,\Sigma)}\mbox{``}A$ is not null''.
\end{theorem}

\begin{proof} In both cases (i.e.\ $e=0$ and $e=1$) the proof is actually the
same so let us deal with the case $e=1$ only (and thus we assume that
$(K,\Sigma)$ is $2$-big and weakly gluing). Suppose that $A\subseteq\can$ is a
set which is not null but for some $p_0\in\q^{\tree}_1(K,\Sigma)$ 
\[p_0\forces_{\q^{\tree}_1(K,\Sigma)}\mbox{``$A$ is null''.}\]
We may assume that $A$ is of outer measure 1 -- just consider the set
\[\{x\in\can: (\exists y\in A)(\forall^\infty n)(x(n)=y(n))\}.\]
Let $\bar{a}=\langle a_n: n<\omega\rangle$, $\bar{k}=\langle k_n:n<\omega
\rangle$ witness that $(K,\Sigma)$ is of the $\NNP$--type. Let $\dot{T}$ be a
$\q^{\tree}_1(K,\Sigma)$-name for a subtree of $\fs$ such that
\[p_0\forces_{\q^{\tree}_1(K,\Sigma)}\mbox{``}\mu([\dot{T}])>
\lim_{n\to\infty} a_n\ \ \&\ \ [\dot{T}]\cap A=\emptyset\mbox{''}.\]
By \ref{treedec}(2) we find a condition $p_1\geq p_0$ and fronts $F_n$
of $T^{p_1}$ such that for each $n<\omega$:
\[(\forall \eta\in F_n)(\mbox{the condition }p_1^{[\eta]}\mbox{ decides
}\dot{T}\cap 2^{\textstyle n}).\] 
Clearly we may assume that
\[F_0=\{\mrot(p_1)\}\quad\mbox{and}\quad(\forall n\!\in\!\omega)(\forall\nu\!
\in\! F_{n+1})(\exists\nu'\!\in\! F_n)(\nu'\vartriangleleft\nu\ \&\ 
k_{n+1}\leq\lh(\nu))\]
and $(\forall n\in\omega)(\forall\nu\in F_n)(\forall\eta\in T^{p_1})(\nu
\trianglelefteq\eta\ \Rightarrow\ \nor[t^{p_1}_\eta]>4+n)$ (remember
\ref{fronor2}). As $(K,\Sigma)$ is weakly gluing and finitary we find a
condition $p_2\geq p_1$ such that 
\[T^{p_2}=\bigcup_{n\in\omega}F_n \quad\mbox{ and }\quad (\forall \eta\in
T^{p_2})(\nor[t^{p_2}_\eta]>3).\]
For $\eta\in F_n$, $n<\omega$ let $g_{\eta}:\pos(t^{p_2}_{\eta})
\longrightarrow {\mathcal P}(2^{\textstyle n+1})$ be a function such that  
\[(\forall\nu\in\pos(t^{p_2}_{\eta}))(p_2^{[\nu]}\forces_{\q^{\tree}_1(K,
\Sigma)}\dot{T}\cap 2^{\textstyle n+1} =g_{\eta}(\nu)).\]
Clearly, if $\eta\in F_n$, $n<\omega$ and $\nu\in\pos(t^{p_2}_\eta)$ then
$a_{n+1}\leq\frac{|g_{\eta}(\nu)|}{2^{n+1}}$. Let (for $\eta\in F_n$,
$n<\omega$) 
\[\begin{array}{ll}
X^n_{\eta}\stackrel{\rm def}{=}\{\sigma\in 2^{\textstyle n+1}\!: &\mbox{there
is }s\in \Sigma(t^{p_2}_{\eta})\mbox{ such that}\\ 
\ &\nor[s]\geq\nor[t^{p_2}_{\eta}]-1,\ \mbox{ and }\ (\forall \nu\in\pos(s))
(\sigma\in g_{\eta}(\nu))\}.\\   
\end{array}\]
Due to $(\circledast)_{\NNP}$ we know that $2^{n+1}\cdot a_n\leq |X^n_\eta|$
(for all $\eta\in F_n$, $n\in\omega$). Fix $n<\omega$. By downward induction
on $m< n$ we define sets $X^n_\eta$ for $\eta\in F_m$: 

if $\eta\in F_m$, $m<n$ then
\[\begin{array}{ll}
X^n_\eta=\{\sigma\in 2^{n+1}: &\mbox{there is }s\in\Sigma(t^{p_2}_\eta)
\mbox{ such that }\\
\ &\nor[s]\geq\nor[t^{p_2}_\eta]-1\ \mbox{ and }\ (\forall\nu\in
\pos(t^{p_2}_\eta))(\sigma\in X^n_\nu))\}.\\
\end{array}\]
Now we may apply the choice of $\bar{a},\bar{k}$ (remembering that $\lh(\eta)
\geq k_m$ for each $\eta\in F_m$) to conclude (by the downward induction on
$m\leq n$) that 
\[(\forall m\leq n)(\forall\eta\in F_m)(|X^n_\eta|\geq 2^{n+1}\cdot a_m).\]
Hence, in particular, $2^{n+1}\cdot a_0\leq |X^n_{\mrot(p_2)}|$ for each
$n<\omega$. So look at the set
\[F\stackrel{\rm def}{=}\{x\in\can: (\exists^\infty n)(x\rest (n+1)\in
X^n_{\mrot(p_2)})\}.\]
Necessarily $\mu(F)\geq a_0>0$ and thus we may take $x\in F\cap A$. For each
$n<\omega$ such that $x\rest (n+1)\in X^n_{\mrot(p_2)}$ we may choose a well
founded quasi tree $S_n$ and a system $\langle s^n_\nu:\nu\in\hat{S}_n\rangle$
of creatures from $K$ such that:  
\[\hat{S}_n\subseteq\bigcup\limits_{m\leq n} F_m,\ \max(S_n)\subseteq
F_{n+1},\ \mrot(S_n)=\mrot(p_2)\ \mbox{ and for all $\nu\in\hat{S}_n$ we
have}\] 
\[\pos(s^n_\nu)=\suc_{S_n}(\nu),\ \  s^n_\nu\in\Sigma(t^{p_2}_\nu),\ \  
\nor[s^n_\nu]\geq\nor[t^{p_2}_\nu]-1\quad\mbox{and } x\rest (n+1)\in
X^n_\nu,\]
and if $\nu\in S_n\cap F_n$, $\nu^*\in\pos(s^n_\nu)$ then $x\rest (n+1)\in
g_\nu(\nu^*)$. Note that if one constructs a condition $q_n$ such that
\[\begin{array}{l}
S_n\subseteq T^{q_n}\subseteq\bigcup\limits_{m<\omega} F_m,\\
(\forall m\leq n)(\forall\nu\in T^{q_n}\cap F_m)(\nu\in \hat{S}_n\ \ \&\ \
t^{q_n}_\nu= s^n_\nu)\quad\mbox{ and}\\
(\forall m>n)(\forall\nu\in T^{q_n}\cap F_m)(t^{q_n}_\nu= t^{p_2}_\nu)\\
\end{array}\]
then $q_n\forces_{\q^{\tree}_1(K,\Sigma)} x\rest (n+1)\in \dot{T}$. Hence, in
particular, 
\begin{enumerate}
\item[$(\oplus)$] {\em if} $n$ is such that $x\rest (n+1)\in
X^n_{\mrot(p_2)}$, $m\leq n+1$ and $\nu\in S_n\cap F_m$\\
{\em then} $p_2^{[\nu]}\forces x\rest m\in \dot{T}$.
\end{enumerate}
Applying the K\"onig Lemma (remember that $(K,\Sigma)$ is strongly finitary!)
we find a quasi tree $S\subseteq\bigcup\limits_{n\in\omega} F_n$ and a system
$\langle s_\nu: \nu\in S\rangle$ such that
\begin{enumerate}
\item $\max(S)=\emptyset$,
\item for all $\nu\in S$:
\[\suc_S(\nu)=\pos(s_\nu)\ \mbox{ and }\ \nor[s_\nu]\geq\nor[t^{p_2}_\nu]-1\
\mbox{ and }\ s_\nu\in\Sigma(t^{p_2}_\nu),\]
\item for some increasing sequence $0<n_0<n_1<n_2<\ldots<\omega$ we have:
\[\begin{array}{l}
(\forall i\in\omega)(x\rest (n_i+1)\in X^{n_i}_{\mrot(p_2)})\quad\mbox{ and}\\
(\forall i\in\omega)(\forall\nu\in S\cap F_i)(\forall j\geq i)(\nu\in S_{n_j}\
\&\ s_\nu=s^{n_j}_\nu).
\end{array}\] 
\end{enumerate}
The quasi tree $S$ determines a condition $q\in\q^{\tree}_1(K,\Sigma)$
stronger than $p_2$. We claim that\ \ \ $q\forces_{\q^{\tree}_1(K,\Sigma)}
(\forall i<\omega)(x\rest i\in\dot{T})$. Why? Let $i\in\omega$ and $\nu\in
S\cap F_i$. Then 
\[s_\nu=s^{n_i}_\nu,\mbox{ and }\nu\in S_{n_i}\cap F_i,\ i\leq n_i\mbox{
and }x\rest(n_i+1)\in X^{n_i}_{\mrot(p_2)}.\]
But now we may use $(\oplus)$ to conclude that for each such $\nu$
\[p_2^{[\nu]}\forces_{\q^{\tree}_1(K,\Sigma)} x\rest i\in\dot{T},\]
and hence $q\forces x\rest i\in\dot{T}$. Consequently $q\forces x\in
[\dot{T}]$, contradicting $q\geq p_0$. 
\end{proof}

We may get a variant of \ref{treenonnull} for tree creating pairs which are
not gluing (e.g.\ for local $(K,\Sigma)$, see \ref{local}). Then, however, we
have to require more from the witnesses for the $\NNP$--type. 

\begin{theorem}
Suppose that $(K,\Sigma)$ is a strongly finitary $2$--big tree creating pair
of the $\NNP$--type with witnesses $\bar{a},\bar{k}$ such that $k_n=n$. Then 
\[\forces_{\q^{\tree}_1(K,\Sigma)}\mbox{`` }A\mbox{ is non-null ''}\]
whenever $A\subseteq \can$ is a set of positive outer measure.
\end{theorem}

\begin{proof} It is similar to \ref{treenonnull}. We start exactly like there
choosing $A$, $\dot{T}$, $p_1$, fronts $F_n$ of $T^{p_1}$ and functions
$g_n:F_n\longrightarrow{\mathcal P}(2^{\textstyle n})$ such that $p_1^{[
\eta]}\forces\dot{T}\cap 2^{\textstyle n}=g_n(\eta)$ for each $n\in\omega$ and
$\eta\in F_n$. Next we define $g_n(\nu)$ for $\nu\in T^{p_1}$ below $F_n$ by
downward induction, in such a way that:
\begin{quotation}
\noindent if $\nu\in T^{p_1}$ and there is $\eta_0\in F_n$ such that $\nu
\vartriangleleft\eta_0$ and $g_n(\eta)$ has been defined already for all
$\eta\in\pos(t^{p_1}_\nu)$ and $(\forall\eta\in\pos(t^{p_1}_\nu))(a_{\lh(
\eta)}\leq\frac{|g_n(\eta)|}{2^n})$ then 
\[\begin{array}{ll}
g_n(\nu)=\{\sigma\in 2^{\textstyle n}\!: &\mbox{there is }s\in
\Sigma(t^{p_1}_{\nu})\mbox{ such that}\\  
\ &\nor[s]\geq\nor[t^{p_1}_{\nu}]-1,\ \mbox{ and }\ (\forall \eta\in\pos(s))
(\sigma\in g_n(\eta))\}.\\   
\end{array}\]
\end{quotation}
We continue as in \ref{treenonnull} getting suitable $S_n$ for $n\in\omega$
and applying the K\"onig lemma we get $S\subseteq T^{p_1}$ with the
corresponding properties.\\
Note that the main difference is that, in the above construction, we may keep
the demand $\frac{|g_n(\nu)|}{2^n}\geq a_{\lh(\nu)}$ (and this is the
replacement for ``weakly gluing''). 
\end{proof}

\section{(No) Sacks Property}
Recall that a forcing notion $\p$ has {\em the Sacks property} is for every
$\p$--name $\dot{x}$ for a real in $\baire$ we have
\[\forces_{\p} (\exists F\in\prod_{n\in\omega} [\omega]^{\textstyle
n+1}\cap\V)(\forall n\in\omega)(\dot{x}(n)\in F(n)).\]
The Sacks property is equivalent to preserving the basis of the null ideal:
every Lebesgue null set in the extension may be covered by a null set (coded)
in the ground model. Here we are interested in refusing this property, i.e.\
getting forcing notions which do not preserve the basis of the null ideal.
 
\begin{definition}
\label{srsp}
We say that a weak--creating pair $(K,\Sigma)$ {\em strongly violates the
Sacks property} if
\begin{enumerate}
\item[$(\otimes)_{\ref{srsp}}$] for some nondecreasing unbounded function
$f\in\baire$ we have

{\em if}\ \ $t\in K$, $\nor[t]>1$

{\em then} for each $u\in\basis(t)$ there is $n\geq \lh(u)$ such that
\[f(n)<|\{w(n): w\in\pos(u,t)\ \&\ n<\lh(w)\}|.\]
\end{enumerate}
\end{definition}

\begin{theorem}
\label{nosacks}
Let $\bH$ be of countable character and let $(K,\Sigma)$ be a weak creating
pair for $\bH$ which strongly violates the Sacks property. Assume that

{\em either} $(K,\Sigma)$ is a creating pair and $\p$ is one of $\q^*_{{\rm
s}\infty}(K,\Sigma)$, $\q^*_\infty(K,\Sigma)$, 

$\q^*_{{\rm w}\infty}(K,\Sigma)$, $\q^*_f(K,\Sigma)$  

{\em or} $(K,\Sigma)$ is a tree--creating pair and then $\p$ is one of
$\q^{\tree}_e(K,\Sigma)$ ($e<5$). 

\noindent Then the forcing notion $\p$ fails the Sacks property.
\end{theorem}

\begin{proof} 
For simplicity we may assume that $\bH(i)=\omega$ (for all $i\in\omega$). Take
an increasing sequence $\langle n_k: k<\omega\rangle$ of positive integers
such that $k+1<f(n_k)$ for all $k\in\omega$. Let $\dot{W}$ be the $\p$--name
for the generic real (see \ref{thereal}) and let $\dot{x}$ be the $\p$-name
for an element of $\baire$ such that 
\[\forces_{\p} (\forall k\in\omega)(\dot{x}(k)=\pi(\dot{W}\rest n_k)).\]
where $\pi:\fseo\longrightarrow\omega$ is the canonical bijection. Now we
claim that 
\[\forces_{\p} (\forall F\in\prod_{n\in\omega} [\omega]^{\textstyle
n+1}\cap\V)(\exists^\infty n)(\dot{x}(n)\notin F(n)).\]
Why? Suppose that $p\in\p$, $F\in\prod\limits_{n\in\omega}[\omega]^{\textstyle
n}$, $N\in\omega$. By \ref{srsp}, in all relevant cases, we find $k>N$ and
$n\in [n_k,n_{k+1})$ such that $p$ ``allows'' more than $f(n)$ values for
$\dot{W}(n)$. But $k+2\leq f(n_k)\leq f(n)$ and thus the condition $p$
``allows'' more than $k+1$ values for $\dot{x}(k+1)$. 
\end{proof}

\section{Examples}

\begin{example}
\label{tomek1}
We build a tree--creating pair $(K_{\ref{tomek1}},\Sigma_{\ref{tomek1}})$
which is: strongly finitary, $2$-big, t-omittory, gluing, of the $\NNP$--type,
and which strongly violates the Sacks property.  
\end{example}

\begin{proof}[Construction] 
Before we define $(K_{\ref{tomek1}},\Sigma_{\ref{tomek1}})$ let us note some
basic properties of the nice pre-norm $\dpt^0_k$ defined in
\ref{exprno}(1),(4).
\begin{claim}
\label{cl9}
Let $M<\omega$, $k<\omega$.
\begin{enumerate}
\item If $A\subseteq M$ then either $\dpt^0_k(A)\geq\dpt^0_k(M)-\log_{2+k}(2)$

\noindent or $\dpt^0_k(M\setminus A)\geq\dpt^0_k(M)-\log_{2+k}(2)$.
\item Suppose that $a\in (0,1)$, $N<\omega$ and a function
$g:M\longrightarrow{\mathcal P}(N)$ is such that
\[(\forall m<M)(N\cdot\frac{a(k+1)+1}{k+2}\leq |g(m)|).\]
Then $a\cdot N\leq |\{n<N: \dpt^0_k(\{m<M: n\in g(m)\})\geq\dpt^0_k(M)-1\}|$.
\end{enumerate}
\end{claim}

\noindent{\em Proof of the claim:}\ \ \ {\em 2)}\ \ Let $u(n)=\{m<M:n\in
g(m)\}$ and 
\[X=\{n<N: \dpt^0_k(u(n))\geq \dpt^0_k(M)-1\}.\]
Look at the set $Y\stackrel{\rm def}{=}\{(m,n)\in M\times N: n\in g(m)\}$ and
note that 
\[|Y|=\sum_{m<M} |g(m)| \geq N\cdot\frac{a(k+1)+1}{k+2}\cdot M.\]
On the other hand, noticing that \ \ $n\in X\quad\Leftrightarrow\quad
|u(n)|\geq\frac{M}{k+2}$,\ \ we have:
\[|Y|=\sum_{n<N} |u(n)|\leq\sum_{n\in X}|u(n)|+\sum_{n\notin X}\frac{M}{k+2}
\leq |X|\cdot M + (N-|X|)\cdot \frac{M}{k+2}.\] 
Consequently $\frac{a(k+1)+1}{k+2}\cdot N\leq |X|+\frac{N-|X|}{k+2}$ and hence
$a\cdot N\leq |X|$, proving the claim. 
\medskip

Let $\bH(n)=(n+2)^n$.

Let $K_{\ref{tomek1}}$ be the collection of all tree--creatures $t$ for $\bH$
such that  
\begin{enumerate}
\item $\dis[t]$ is a pair $(d_0(t),d_1(t))$ such that $d_0(t)\subseteq\bigcup 
\limits_{n\in\omega}\prod\limits_{m<n}\bH(m)$ is a finite tree, $|d_0(t)|>2$
and $d_1(t)\trianglelefteq\mrot(d_0(t))$,
\item $\nor[t]=\min\{\dpt^0_k(\suc_{d_0(t)}(\eta)): \eta\in\spliting(d_0(t))\
\ \&\ \ \lh(\eta)=k\}$,
\item $\val[t]=\{\langle d_1(t),\nu\rangle: \nu\in\max(d_0(t))\}$.
\end{enumerate}
For a well founded quasi tree $T$ and a system $\langle s_\nu: \nu\in\hat{T}
\rangle$ of tree--creatures from $K_{\ref{tomek1}}$ such that the requirement
{\bf (a)} of \ref{treecreature} is satisfied we let
\[\Sigma_{\ref{tomek1}}(s_\nu:\nu\in\hat{T})=\{s\in K_{\ref{tomek1}}:
d_1(s)=\mrot(T)\ \&\ \rng(\val[s])\subseteq\max(T)\}.\]
It should be clear that $\Sigma_{\ref{tomek1}}$ is a tree-composition on
$K_{\ref{tomek1}}$. Now, the tree--creating pair $(K_{\ref{tomek1}},\Sigma_{
\ref{tomek1}})$ is strongly finitary, t-omittory and gluing. For the last two
properties we apply the procedure similar to the one below.

Note that if $t\in K_{\ref{tomek1}}$ and $\nu\in d_0(t)$ is a splitting point
of $d_0(t)$ then choosing $\eta_\rho\in\max(d_0(t))$ (for $\rho\in
\suc_{d_0(t)}(\nu)$) such that $\rho\trianglelefteq\eta_\rho$ we may build a
tree creature $s\in\Sigma_{\ref{tomek1}}(t)$ such that $\pos(s)=\{\eta_\rho:
\rho\in\suc_{d_0(t)}(\nu)\}$. Then we will have $\nor[s]\geq\nor[t]$. If
additionally $\nu\in d_0(t)$ is a splitting point such that $\nor[t]=
\dpt^0_{\lh(\nu)}(\suc_{d_0(t)}(\nu))$ then $\nor[s]=\nor[t]$ (see the
definition of $K_{\ref{tomek1}}$). In this case, let us call the respective
tree--creature $s(t)$ (here we just fix one such).   

Considering suitable $s(t)$'s and using \ref{cl9}(1) one can easily show that
the creating pair $(K_{\ref{tomek1}},\Sigma_{\ref{tomek1}})$ is $2$-big.

Let $k_n=2^{n+4}-2$, $a_0=\frac{1}{2}$, $a_{n+1}=\frac{a_n\cdot (k_n+1)+1}{k_n
+2}=a_n-\frac{a_n}{2^{n+4}}+\frac{1}{2^{n+4}}$. We are going to show that the
sequences $\bar{k}=\langle k_n: n<\omega\rangle$, $\bar{a}=\langle a_n: n<
\omega\rangle$ witness that $(K_{\ref{tomek1}},\Sigma_{\ref{tomek1}})$ is of
the $\NNP$--type.
 
First note that $\bar{a},\bar{k}$ are strictly increasing, $\bar{a}\subseteq
(0,1)$ and $\lim\limits_{n\to\infty}a_n\leq\frac{5}{8}$. Now suppose that
$t\in K_{\ref{tomek1}}$, $k_n\leq m^t_{\dn}$, $N<\omega$ and $g:\pos(t)
\longrightarrow {\mathcal P}(N)$ is such that 
\[(\forall\nu\in\pos(t))(a_{n+1}\leq\frac{|g(\nu)|}{N}).\]
Take $s(t)$ (defined above) and look at $h=g\rest \pos(s(t))$. We may think
that actually $h:\suc_{d_0(s(t))}(\nu)\longrightarrow{\mathcal P}(N)$, where
$\nu$ is the unique splitting point of $s(t)$ (note that $\lh(\nu)\geq
m^{s(t)}_{\dn}\geq k_n$). Applying claim \ref{cl9}(2) we get  
\[a_n\cdot N\leq |\{n<N: \dpt^0_{k_n}(\{m: n\in h(m)\})\geq
\dpt^0_{k_n}(\suc_{d_0(s(t))}(\nu))-1\}|.\]
For each $n<N$ from the set on the right hand side of the inequality above
choose $s_n\in\Sigma_{\ref{tomek1}}(s(t))$ such that $\pos(s_n)=\{\eta\in 
\pos(s(t)): n\in h(\eta(\lh(\nu)))\}$. By the definition of $\dpt^0_k$,
$\dpt^0_k(w_0)\geq \dpt^0_k(w_1)-1$ implies that $\dpt^0_{k+1}(w_0)\geq
\dpt^0_{k+1}(w_1)-1$, and therefore $\dpt^0_{k_n}(\{m:n\in h(m)\})\geq
\dpt^0_{k_n}(\suc_{d_0(s(t))}(\nu))-1$ implies that $\nor[s_n]\geq \nor[s(t)]
-1$. Now we may conclude that the set 
\[\begin{array}{ll}
\{n<N:&\mbox{there is }s\in\Sigma_{\ref{tomek1}}(t)\mbox{ such that
}\nor[s]\geq\nor[t]-1\mbox{ and}\\
\ &(\forall\nu\in\pos(s))(n\in g(\nu))\}.\\
\end{array}\]
has not less than $a_n\cdot N$ elements. 

Finally note that $(K_{\ref{tomek1}},\Sigma_{\ref{tomek1}})$ satisfies the
condition $(\oplus_{\ref{srsp}})$ for $f(n)=n$ (so it strongly violates the
Sacks property). 
\end{proof}

\begin{conclusion}
\label{tomek1conc}
The forcing notions $\q^{\tree}_e(K_{\ref{tomek1}},\Sigma_{\ref{tomek1}})$ for
$e\!<\!5$ are equivalent. They are proper, preserve the outer measure, preserve
non-meager sets, are $\baire$-bounding, but do not have the Sacks property.
\end{conclusion}

\begin{example}
\label{loctom}
We construct a finitary, $2$--big and  local tree--creating pair
$(K_{\ref{loctom}},\Sigma_{\ref{loctom}})$ which is trivially meagering and of
the $\NNP$--type with the sequence $\bar{k}=\langle n:n<\omega\rangle$
witnessing it.  
\end{example}

\begin{proof}[Construction] This is similar to \ref{tomek1}. We define
$\bar{k}^*$, $\bar{k}$ and $\bar{a}$ letting $k^*_n=2^{n+4}-2$, $k_n=n$,
$a_0=\frac{1}{2}$, $a_{n+1}=\frac{a_n(k^*_n+1)+1}{k^*_n+2}$. Let $\bH(n)=
2^{(n+4)^2}$. 

The family $K_{\ref{loctom}}$ consists of these $t\in\TCR[\bH]$ that:
\begin{itemize}
\item $\dis[t]=\langle\eta,A_t\rangle$, where $\eta$ is such that $t\in
\TCR_\eta[\bH]$ and $A_t\subseteq\bH(\lh(\eta))$,
\item $\val[t]=\{\langle\eta,\nu\rangle: \eta\vartriangleleft\nu\ \&\
\lh(\nu)=\lh(\eta)+1\ \&\ \nu(\lh(\eta))\in A_t\}$,
\item $\nor[t]=\dpt^0_{k^*_{\lh(\eta)}}(A_t)$.
\end{itemize}
The operation $\Sigma_{\ref{loctom}}$ is trivial and $s\in
\Sigma_{\ref{loctom}}(t)$ if and only if $\val[s]\subseteq\val[t]$. 

Plainly, $(K_{\ref{loctom}},\Sigma_{\ref{loctom}})$ is a $2$--big local and
trivially meagering tree--creating pair. Checking that it is of the
$\NNP$--type (with witnesses $\bar{a}$ and $\bar{k}$) is exactly like in
\ref{tomek1} (just apply \ref{cl9}). 
\end{proof}

\begin{conclusion}
$\q^{\tree}_1(K_{\ref{loctom}},\Sigma_{\ref{loctom}})$ is a proper
$\baire$--bounding forcing notion which preserves outer measure but makes the
ground model reals meager.
\end{conclusion}

\begin{example}
\label{ponic}
We define functions $\bH$, $f$ and a creating pair $(K_{\ref{ponic}},\Sigma_{
\ref{ponic}})$ such that
\begin{enumerate}
\item $\bH$ is finitary, $f:\omega\times\omega\longrightarrow\omega$ is
$\bH$--fast, 
\item $(K_{\ref{ponic}},\Sigma_{\ref{ponic}})$ is gluing, forgetful,
$\bar{2}$--big, has the Halving Property and is of the $\NMP$--type,
\item the forcing notion $\q^*_f(K_{\ref{ponic}},\Sigma_{\ref{ponic}})$ is not
trivial. 
\end{enumerate}
\end{example}

\begin{proof}[Construction] 
Define inductively $\bH$ and $f$ such that for all $n,k,\ell\in\omega$:
\begin{enumerate}
\item[(i)] \ $\bH(0)=8$,\ \ $\bH(n)=2^{\fH(n)+f(n,n)}$,
\item[(ii)] $f(0,\ell)=\ell+1$,\ \ $f(k+1,\ell)=2^{\fH(\ell)+1}\cdot
(f(k,\ell)+\fH(\ell)+2)$
\end{enumerate}
(compare with \ref{bighalex}; note that the above conditions uniquely
determine $\bH$ and $f$). By their definition $\bH$ is finitary and $f$ is
$\bH$--fast. 

A creature $t\in\CR_{m_0,m_1}[\bH]$ belongs to $K_{\ref{ponic}}$ if:
\begin{itemize}
\item $\dis[t]=\langle m_0,m_1,X_t,H_t,\langle A^k_t,H^k_t,u^k_t:k\in X_t
\rangle\rangle$, \quad where $X_t\subseteq [m_0,m_1)$, $H_t$ is a nice
pre-norm on ${\mathcal P}(X_t)$, and for each $k\in X_t$:\\
$A^k_t\subseteq\bH(k)$, $H^k_t$ is a nice pre-norm on ${\mathcal
P}(A^k_t)$ and $u^k_t\in\prod\limits_{i\in [m_0,m_1)\setminus\{k\}}\bH(i)$,
\item $\val[t]=\{\langle u,v\rangle\!\in\!\prod\limits_{i<m_0}\bH(i)\!\times\!
\prod\limits_{i<m_1}\bH(i)\!: u\vartriangleleft v\ \ \&\ \ (\exists k\in X_t)
(u^k_t\subseteq v\ \&\ v(k)\in A^k_t)\}$,
\item $\nor[t]=\min\{H_t(X_t),H^k_t(A^k_t): k\in X_t\}$.
\end{itemize}
Now we describe the operation $\Sigma_{\ref{ponic}}$. Suppose that $t_0,
\ldots,t_n\in K_{\ref{ponic}}$ are such that $m^{t_i}_{\up}=m^{t_{i+1}}_{\dn}$
for $i<n$. Let $\Sigma_{\ref{ponic}}(t_0,\ldots,t_n)$ consist of all creatures
$s\in K_{\ref{ponic}}$ such that $m^s_{\dn}=m^{t_0}_{\dn}$, $m^s_{\up}=
m^{t_n}_{\up}$ and for some $i\leq n$
\begin{enumerate}
\item[$(\alpha)$] $X_s\subseteq X_{t_i}$,\quad $(\forall B\subseteq X_s)(H_s(
B)\leq H_{t_i}(B))$ and for every $k\in X_s$:
\item[$(\beta)$]  $A^k_s\subseteq A^k_{t_i}$ and $(\forall A\subseteq A^k_s)(
H^k_s(A)\leq H^k_{t_i}(A))$, \quad and
\item[$(\gamma)$] $u^k_{t_i}\subseteq u^k_s$ and for every $j\in (n+1)
\setminus \{i\}$ for some $\ell\in X_{t_j}$ we have
\[u^\ell_{t_j}\subseteq u^k_s\quad\mbox{ and }\quad u^k_s(\ell)\in A^\ell_{
t_j}.\]
\end{enumerate}
It should be clear that $(K_{\ref{ponic}},\Sigma_{\ref{ponic}})$ is a gluing
and forgetful creating pair for $\bH$. It is $\bar{2}$--big as for each $t\in
K_{\ref{ponic}}$ both $H_t$ and $H^k_t$ (for $k\in X_t$) are nice pre-norms. 
We may use similar arguments as in \ref{winnonmea} to show that
$(K_{\ref{ponic}},\Sigma_{\ref{ponic}})$ is of the $\NMP$--type. Now define
function $\uhalf:K_{\ref{ponic}}\longrightarrow K_{\ref{ponic}}$ by
\begin{quotation}
\noindent $\uhalf(t)=s$ \quad if and only if

\noindent $m^s_{\dn}=m^t_{\dn}$,\ \ $m^s_{\up}=m^t_{\up}$,\ \ $X_s=X_t$,\ \
$A^k_s=A^k_t$ for $k\in X_s$ and 
\[H_s=(H_t)^{\frac{1}{2}\nor[t]},\quad H^k_s=(H^k_t)^{\frac{1}{2}\nor[t]}\ \
\mbox{ for }k\in X_s\]
(see \ref{exareok}(2)).
\end{quotation}
It is not difficult to check that the function $\uhalf$ witnesses that
$(K_{\ref{ponic}},\Sigma_{\ref{ponic}})$ has the Halving Property.

Note that $(K_{\ref{ponic}},\Sigma_{\ref{ponic}})$ resembles an omittory
creating pair. 
\end{proof}

\chapter{Omittory with Halving}
In \ref{halbigdec}, \ref{gfbound} we saw how the Halving Property and the
bigness apply to forcing notions $\q^*_f(K,\Sigma)$. In this chapter we will
look at another combination: omittory creating pairs with the weak Halving
Property. Since an omittory creating pair cannot be big, it is natural that we
consider in this context the (s$\infty$) norm condition. The first example of
a forcing notion of the type $\q^*_{{\rm s}\infty}(K,\Sigma)$ for an omittory
creating pair $(K,\Sigma)$ with the weak Halving Property appeared in 
\cite{Sh:207} (but in the real application there a different norm condition
was used). A direct application of a forcing notion of this type was presented
in \cite{RoSh:501}. In the last part of this chapter we will develop the
example from that paper. Before, in the first section, we show that the
forcing notions $\q^*_{{\rm s}\infty}(K,\Sigma)$ with $(K,\Sigma)$ omittory
tend to add Cohen reals and make ground reals meager. Next we introduce some
general operations on creating pairs and, in the third section, we explain how
the weak Halving Property may prevent them from adding dominating reals.

\section{What omittory may easily do}
Natural examples of omittory creating pairs with the weak Halving Property
are meagering and anti-big (see \ref{meagering} below). We will show how these
properties cause that forcing notions $\q^*_{{\rm s}\infty}(K,\Sigma)$ do some
harm to the old reals. Examples and applications are presented in the last
part of this chapter.  

First note the following easy observation.

\begin{proposition}
\label{adunbreal}
If $(K,\Sigma)$ is an omittory creating pair such that for each $t\in K$,
$u\in\basis(t)$
\[\nor[t]>0\quad\Rightarrow\quad |\pos(u,t)|>2\]
then $\forces_{\q^*_{{\rm s}\infty}(K,\Sigma)}$ ``there is an unbounded real
over $\V$''.
\end{proposition}

\begin{definition}
\label{meagering}
Let $(K,\Sigma)$ be a creating pair.
\begin{enumerate}
\item We say that $(K,\Sigma)$ is {\em meagering} if for every $(t_0,\ldots,
t_{n-1})\in\PFC(K,\Sigma)$, $t\in\Sigma(t_0,\ldots,t_{n-1})$ and $\langle k_i:
i<n\rangle$ such that for each $i<n$: 
\[\nor[t_i]> 2\quad\mbox{ and }\quad m^{t_i}_{\dn}\leq k_i <m^{t_i}_{\up}
\quad\mbox{ and }\quad\nor[t]> 2\]
there is $s\in\Sigma(t)$ satisfying 
\[\begin{array}{l}
\nor[s]\geq\nor[t]-1\qquad\mbox{ and}\\
(\forall u\in\basis(t_0))(\exists v\in\pos(u,s))(\exists k\in [\lh(u),\lh(v)))
(v(k)\neq 0)\qquad\mbox{ and}\\
(\forall u\in\basis(t_0))(\forall v\in\pos(u,s))(\forall i<n)(v(k_i)=0).
  \end{array}\]
\item The creating pair $(K,\Sigma)$ is called {\em anti-big} if there are
colourings 
\[c_t:\bigcup\limits_{u\in\basis(t)}\pos(u,t)\longrightarrow 3\qquad \mbox{
for } t\in K\]
such that:\quad if $(t_0,\ldots,t_{n-1})\in\PFC(K,\Sigma)$, $\nor[t_i]>1$ (for
$i<n$) and $t\in\Sigma(t_0,\ldots,t_{n-1})$, $\nor[t]>1$ then for each
$u\in\basis(t_0)$ there are $v_0,v_1\in\pos(u,t)$ and $\ell<n$ satisfying
\[\begin{array}{l}
v_0\rest m^{t_\ell}_{\dn}=v_1\rest m^{t_\ell}_{\dn},\quad c_{t_\ell}(v_0\rest
m^{t_\ell}_{\up})=0,\quad c_{t_\ell}(v_1\rest m^{t_\ell}_{\up})=1,\quad\mbox{
and}\\
(\forall i\in n\setminus\{\ell\})(c_{t_i}(v_0\rest m^{t_i}_{\up})=c_{t_i}(v_1
\rest m^{t_i}_{\up})=2).
  \end{array}\]
\end{enumerate}
\end{definition}

\begin{theorem}
\label{makemeager}
Let $(K,\Sigma)$ be a growing creating pair. 
\begin{enumerate}
\item If $(K,\Sigma)$ is meagering then
\[\forces_{\q^*_{{\rm s}\infty}(K,\Sigma)}\mbox{ ``}\baire\cap\V\mbox{ is
meager''.}\]
\item If $(K,\Sigma)$ is anti-big then
\[\forces_{\q^*_{{\rm s}\infty}(K,\Sigma)}\mbox{ ``there is a Cohen real over
}\V\mbox{''.} \]
\end{enumerate}
\end{theorem}

\begin{proof} 
1)\ \ \ Let $p=(w^p,t^p_0,t^p_1,\ldots)\in\q^*_{{\rm s}\infty}(K,\Sigma)$ and
for $n\in\omega$ let $f(n)={\mathcal P}([m^{t^p_n}_{\dn}, m^{t^p_n}_{\up}))$. 
The space $\prod\limits_{n\in\omega}f(n)$ equipped with the product topology
(of the discrete $f(n)$'s) is a perfect Polish space. Thus it is enough to
show that  
\[p\forces_{\q^*_{{\rm s}\infty}(K,\Sigma)}\mbox{`` }\prod_{n\in\omega} f(n)
\cap\V\mbox{ is a meager subset of }\prod_{n\in\omega} f(n)\mbox{ ''.}\]
Note that if $X\in\prod\limits_{n\in\omega} f(n)$ is such that
$(\exists^\infty n\in\omega)(X(n)\neq\emptyset)$ then the set 
\[\{Y\in\prod_{n\in\omega} f(n): (\forall^\infty n\in\omega)(Y(n)=\emptyset
\mbox{ or }Y(n)\neq X(n))\}\]
is meager in $\prod\limits_{n\in\omega} f(n)$. Let $\dot{X}$ be a $\q^*_{{\rm
s}\infty}(K,\Sigma)$--name for an element of $\prod\limits_{n\in\omega} f(n)$
such that 
\[p\forces_{\q^*_{{\rm s}\infty}(K,\Sigma)}(\forall n\in\omega)(\dot{X}(n)=
\{k\in [m^{t^p_n}_{\dn},m^{t^p_n}_{\up}): \dot{W}(k)\neq 0\}),\]
where $\dot{W}$ is the $\q^*_{{\rm s}\infty}(K,\Sigma)$--name for the generic
real (see \ref{thereal}). It follows from the remarks above that it is enough
to show that
\begin{enumerate}
\item[($\alpha$)] $p\forces_{\q^*_{{\rm s}\infty}(K,\Sigma)}(\exists^\infty
n\in\omega)(\dot{X}(n)\neq \emptyset)$\quad and
\item[($\beta$)]  $p\forces_{\q^*_{{\rm s}\infty}(K,\Sigma)}(\forall Y\in
\prod\limits_{n\in\omega} f(n)\cap \V)(\forall^\infty n\in\omega)(Y(n)=
\emptyset\mbox{ or }Y(n)\neq\dot{X}(n))$.
\end{enumerate}
To this end suppose that $p\leq q=(w^q,t^q_0,t^q_1,\ldots)\in\q^*_{{\rm s}
\infty}(K,\Sigma)$ is such that $\lh(w^q)>1$ and $\nor[t^q_i]>m^{t^q_i}_{\dn}
+2$ for $i\in\omega$ and let $Y\in\prod\limits_{n\in\omega} f(n)$. For each
$n\in\omega$ choose $k_n\in [m^{t^p_n}_{\dn},m^{t^p_n}_{\up})$ such that
$Y(n)\neq \emptyset\ \Rightarrow\ k_n\in Y(n)$. Let $0\leq n_0<n_1<n_2<\ldots<
\omega$ be such that $w^q\in\pos(w^p,t^p_0,\ldots,t^p_{n_0-1})$ and $t^q_i\in
\Sigma(t^p_{n_i},\ldots,t^p_{n_{i+1}-1})$ for $i\in\omega$. Note that
necessarily $m^{t^p_{n_0}}_{\dn}\geq 2$ and thus $\nor[t^p_n]>2$ for each
$n\geq n_0$. Applying \ref{meagering}(1) we find $s_i\in\Sigma(t^q_i)$ such
that for each $i\in\omega$:
\begin{enumerate}
\item[$(*)^1_i$] $\nor[s_i]\geq\nor[t^q_i]-1>m^{t^q_i}_{\dn}$, and
\item[$(*)^2_i$] $(\forall u\!\in\!\pos(w^p,t^p_0,\ldots,t^p_{n_i-1}))
(\exists v\!\in\!\pos(u,s_i))(\exists k\!\in\! [\lh(u),\lh(v)))(v(k)\neq 0)$,
\item[$(*)^3_i$] $(\forall u\in\pos(w^p,t^p_0,\ldots,t^p_{n_i-1}))(\forall
v\in\pos(u,s_i))(\forall n\in [n_i,n_{i+1}))(v(k_n)=0)$.
\end{enumerate}
Look at $r\stackrel{\rm def}{=}(w^q,s_0,s_1,s_2,\ldots)$. Clearly $r\in \q^*_{
{\rm s}\infty}(K,\Sigma)$ is a condition stronger than $q$. Moreover, by the
choice of the $s_i$'s we have
\[r\forces_{\q^*_{{\rm s}\infty}(K,\Sigma)} (\forall n\geq n_0)(\dot{W}(k_n)=
0)\]
and therefore, by the choice of the $k_n$'s, we get 
\[r\forces_{\q^*_{{\rm s}\infty}(K,\Sigma)} (\forall n\geq n_0)(Y(n)\neq
\emptyset\quad\Rightarrow\quad Y(n)\neq\dot{X}(n)).\]
Further note that if $v\in\pos(w^p,s_0)$ is given by $(*)^2_0$ above for $w^q$
then 
\[(v,s_1,s_2,\ldots)\forces (\exists n\in [n_0,n_1))(\dot{X}(n)\neq
\emptyset).\] 
We finish by density arguments.
\medskip

\noindent 2)\ \ \ Let $p=(w^p,t^p_0,t^p_1,\ldots)\in\q^*_{{\rm s}\infty}(K,
\Sigma)$ and let $\dot{Z}$ be a $\q^*_{{\rm s}\infty}(K,\Sigma)$--name for a
subset of $\omega$ such that
\[p\forces_{\q^*_{{\rm s}\infty}(K,\Sigma)}\dot{Z}=\{n\in\omega: c_{t^p_n}
(\dot{W}\rest m^{t^p_n}_{\up})<2\}.\]
Note that $p\forces$``$\dot{Z}$ is infinite''. Why? Suppose $p\leq q\in
\q^*_{{\rm s}\infty}(K,\Sigma)$, $\lh(w^q)>1$. Let $n_0<n_1<\omega$ be such
that $w^q\in\pos(w^p,t^p_0,\ldots,t^p_{n_0-1})$, $t^q_0\in\Sigma(t^p_{n_0},
\ldots,t^p_{n_1-1})$. Necessarily $\nor[t^p_i]>1$ for $i\in [n_0,n_1)$ and
$\nor[t^q_0]>1$. So we find $\ell\in [n_0,n_1)$ and $v_0,v_1\in\pos(w^q,
t^q_0)$ as in \ref{meagering}(2). Now look at the condition
$r=(v_0,t^q_1,t^q_2, \ldots)\in\q^*_{{\rm s}\infty}(K,\Sigma)$. It is stronger
than $q$ and forces that $\ell\in \dot{Z}$. 

Now let $\dot{c}$ be a $\q^*_{{\rm s}\infty}(K,\Sigma)$--name for a real in
$\can$ such that 
\[p\forces_{\q^*_{{\rm s}\infty}(K,\Sigma)}\mbox{`` if $k\in\omega$ and $n$ is
the $k^{\rm th}$ member of $\dot{Z}$ then }\dot{c}(k)=c_{t^p_n}(\dot{W}\rest
m^{t^p_n}_{\up})\mbox{ ''.}\]
We claim that $p\forces$ ``$\dot{c}$ is a Cohen real over $\V$''. So suppose
that $p\leq q\in \q^*_{{\rm s}\infty}(K,\Sigma)$, $\lh(w^q)>1$ and ${\mathcal
U}\subseteq \can$ is an open dense set. Let $0\leq n_0<n_1<\ldots<\omega$ be
such that 
\[w^q\in\pos(w^p,t^p_0,\ldots,t^p_{n_0-1})\quad\mbox{ and }\quad t^q_i\in
\Sigma(t^p_{n_i},\ldots,t^p_{n_{i+1}-1})\quad\mbox{for }i\in\omega.\]
Let $m=|\{n<n_0: c_{t^p_n}(w^q\rest m^{t^p_n}_{\up})<2\}|$ and let $\nu\in
2^{\textstyle m}$ be such that $\nu(k)=c_{t^p_n}(w^q\rest m^{t^p_n}_{\up})$ if
$k<m$ and $n<n_0$ is the $k^{\rm th}$ member of the set $\{n<n_0:
c_{t^p_n}(w^q\rest m^{t^p_n}_{\up})<2\}$. Choose $\eta\in\fs$ such that
$\nu\vartriangleleft\eta$ and 
\[(\forall x\in\can)(\eta\vartriangleleft x\ \Rightarrow\ x\in {\mathcal
U}).\]
Let $j=\lh(\eta)-m$. Use \ref{meagering}(2) to define inductively $u\in
\pos(w^q,t^q_0,\ldots,t^q_{j-1})$ such that for each $i<j$, for some $\ell\in
[n_i,n_{i+1})$ we have 
\[c_{t^p_\ell}(u\rest m^{t^p_\ell}_{\up})=\eta(\ell)\quad\mbox{ and }\quad
(\forall k\in [n_i,n_{i+1})\setminus\{\ell\})(c_{t^p_k}(u\rest
m^{t^p_k}_{\up})=2).\]
Look at the condition $r=(u,t^q_j,t^q_{j+1},\ldots)\in\q^*_{{\rm s}\infty}(K,
\Sigma)$: it is stronger than $q$ and it forces that
$\eta\vartriangleleft\dot{c}$. We finish by density argument.
\end{proof}

\section{More operations on weak creatures}
Below we define some operations on creatures and tree--creatures which provide
for (some) systems of weak creatures a new weak creature (of the same
type). These operations may be used to define sub--composition operations.

\begin{definition}
\label{sum}
Suppose $0<m<\omega$ and for $i<m$ we have $t_i\in\CR[\bH]$ such that
$m^{t_i}_{\up}\leq m^{t_{i+1}}_{\dn}$. Then we define {\em the sum of the
creatures $t_i$} as a creature $t=\Sigma^{\bsum}(t_i:i<m)$ such that (if well
defined then):
\begin{enumerate}
\item[(a)] $m^t_{\dn}=m^{t_{0}}_{\dn}$, $m^t_{\up}=m^{t_{m-1}}_{\up}$,
\item[(b)] $\val[t]$ is the set of all pairs $\langle h_1, h_2\rangle$ such
that:  
\begin{quotation}
\noindent $\lh(h_1)=m^t_{\dn}$, $\lh(h_2)=m^t_{\up}$, $h_1\vartriangleleft
h_2$,  

\noindent and $\langle h_2\rest m^{t_i}_{\dn}, h_2\rest
m^{t_i}_{\up}\rangle\in\val[t_i]$ for $i<m$, 

\noindent and $h_2\rest [m^{t_i}_{\up}, m^{t_{i+1}}_{\dn})$ is identically
zero for $i<m-1$,
\end{quotation}
\item[(c)] $\nor[t]=\min\{\nor[t_i]:i<m\}$,
\item[(d)] $\dis[t]=\langle 0\rangle\conc\langle\dis[t_i]: i<m\rangle$.
\end{enumerate}
If for all $i<m-1$ we have $m^{t_i}_{\up}=m^{t_{i+1}}_{\dn}$ then we call the
sum {\em tight}.  
\end{definition}

\begin{remark}
Note that the sum $\Sigma^{\bsum}(t_i: i<m)$ is defined only for these
sequences $\langle t_i: i<m\rangle \subseteq \CR[\bH]$ for which
$m^{t_i}_{\up}\leq m^{t_{i+1}}_{\dn}$ and part (b) of the definition gives a
nonempty value of $\val[t]$. 
\end{remark}

\begin{definition}
\label{dsum}
If $m<\omega$, $u\subseteq m$, $d\in{\mathcal H}(\chi)$ is a function such
that $\dom(d)\supseteq (\mbR^{{\geq}0})^{\textstyle u}$ and $\rng(d)\subseteq
\mbR^{{\geq}0}$, and for $i<m$ creatures $t_i\in\CR[\bH]$ are such that
$m^{t_i}_{\up}\leq m^{t_{i+1}}_{\dn}$ then we define {\em the $(d,u)$-sum
$t=\Sigma_{d,u}^{\bsum}(t_i: i<m)$ of the $t_i$'s} by: 
\begin{enumerate}
\item[(a)] $m^t_{\dn}=m^{t_{0}}_{\dn}$, $m^t_{\up}=m^{t_{m-1}}_{\up}$,
\item[(b)] $\val[t]$ is the set of pairs $\langle h_1,h_2\rangle$ such that:
\begin{quotation}
\noindent $\lh(h_1)=m^t_{\dn}$, $\lh(h_2)= m^t_{\up}$, $h_1\vartriangleleft
h_2$ and 

\noindent $\langle h_2\rest m^{t_i}_{\dn},h_2\rest m^{t_i}_{\up}\rangle\in
\val[t_i]$ for $i\in u$,

\noindent $h_2\rest [m^{t_i}_{\dn},m^{t_i}_{\up})$ is identically zero for
$i\notin u$ and

\noindent $h_2 \rest [m^{t_i}_{\up}, m^{t_{i+1}}_{\dn})$ is identically zero
for $i<m-1$. 
\end{quotation}
\item[(c)] $\nor[t]=d(\langle\nor[t_i]: i\in u\rangle)$,
\item[(d)] $\dis[t]=\langle 1,d,u\rangle\conc\langle\dis[t_i]:i<m\rangle$. 
\end{enumerate}
[Note: the $(d,u)$-sum is defined only if clause (b) gives a nonempty value
for $\val[t]$.] 
\end{definition}

\begin{definition}
\label{saturated}
\begin{enumerate}
\item For a pre-norm $H$ on $\omega$ (see \ref{prenorms}) let $D_H$ be the
family of all functions $d$ such that for some finite set
$u_d\subseteq\omega$, $H(u_d)>0$ and 
\[d:(\mbR^{{\geq}0})^{\textstyle u_d}\longrightarrow \mbR^{{\geq}0}:\langle
r_i:i\in u_d\rangle\mapsto\min\{H(u_d),r_i: i\in u_d\}.\]
\item We say that a creating pair $(K,\Sigma)$ is {\em saturated with respect
to a pre-norm $H$ on $\omega$} if for each $d\in D_H$ and $(t_i:m_0\leq i<m_1)
\in\PFC(K,\Sigma)$ such that $u_d\subseteq [m_0,m_1)$ and $\nor[t_i]>0$ for
$i\in u_d$: 
\[\Sigma^{\bsum}_{d,u_d}(t_i:m_0\leq i<m_1)\mbox{ is well defined and
belongs to }\Sigma(t_i: m_0\leq i<m_1),\]
and if $t\in \Sigma(\Sigma^{\bsum}_{d,u_d}(t_i: m_0\leq i<m_1))$, $\nor[t]>0$
then for some $d^*\in D_H$ and $s_i\in \Sigma(t_i)$ (for $m_0\leq i<m_1$) we
have
\[u_{d^*}\subseteq u_d,\quad \val[\Sigma^{\bsum}_{d^*,u_{d^*}}(s_i: m_0\leq
i<m_1)] \subseteq\val[t],\quad\mbox{and }\nor[s_i]>0\mbox{ for }i\in
u_{d^*}.\]
We say that $(K,\Sigma)$ is {\em saturated with respect to (nice) pre--norms}
if for each (nice) pre-norm $H$ on $\omega$, $(K,\Sigma)$ is saturated with
respect to $H$. Similarly for other classes of pre--norms.
\end{enumerate}
\end{definition}

\begin{remark}
Note that in practical realizations of \ref{saturated}(2) the additional
parameter $\dis$ may play a crucial role. Looking at a creature $t$ we may
immediately recognize if it comes from the operation $\Sigma^{\bsum}_{d,u_d}$
and we do not have to worry that the last demand gives a contradiction. It
may happen that for distinct $d$'s from $D_H$ we get (as a result of
$\Sigma^{\bsum}_{d,u_d}$) creatures with the same values of $\val$, $\nor$,
however they are distinguished by $\dis$. Moreover, the same effect appears
for distinct pre--norms $H$: we can read from $\dis$ the function $d$ and
consequently the function $H$ restricted to subsets of $u_d$. In applications
we may redefine $\dis[\Sigma^{\bsum}_{d,u_d}(t_i: i<m)]$, but we should keep
this coding property.  
\end{remark}

\begin{definition}
\label{treesum}
Let $T\subseteq\bigcup\limits_{n<\omega}\prod\limits_{i<n} \bH(i)$ be a well
founded quasi tree and let $\langle s_\nu:\nu\in \hat{T}\rangle\subseteq\TCR
[\bH]$ be a system of tree creatures such that for each $\nu\in\hat{T}$:
\[\dom(\val[s_\nu])=\{\nu\}\quad\mbox{ and }\quad \pos(s_\nu)=\suc_T(\nu).\]
\begin{enumerate}
\item {\em The tree--sum $t=\Sigma^{\tsum}(s_\nu:\nu\in \hat{T})$ of tree
creatures $s_\nu$} (for $\nu\in\hat{T}$) is defined by:
\begin{enumerate}
\item[$(\alpha)$] $\val[t]=\{\langle\mrot(T),\eta\rangle:\eta\in\max(T)\}$,
\item[$(\beta)$]  $\nor[t]=\inf\{\nor[s_\nu]:\nu\in \hat{T}\ \&\ \nor[s_\nu]
\geq 1\}$, if $\nor[s_\nu]<1$ for all $\nu\in \hat{T}$ then we let $\nor[t]
=0$,
\item[$(\gamma)$] $\dis[t]=\langle 2\rangle\conc\langle s_\nu:\nu\in\hat{T}
\rangle$. 
\end{enumerate}
\item For a function $g\in\baire$, {\em the special tree--sum $t=
\Sigma^{\tsum}_g(s_\nu: \nu\in\hat{T})$ of tree creatures $s_\nu$ (for
$\nu\in\hat{T}$) with respect to $g$} is defined in a similar manner as
$\Sigma^{\tsum}$ but the conditions $(\beta)$, $(\gamma)$ introducing the norm
and $\dis$ are replaced by 
\begin{enumerate}
\item[$(\beta)^*_g$] $\nor[t]=\max\{k<\omega:(\forall\eta\in\max(T))(|\{\ell
<\lh(\eta)\!:\eta\rest\ell\in \hat{T}$ and $\nor[s_{\eta\rest\ell}]\geq k\}|
\geq g(k))\}$,
\item[$(\gamma)^*_g$] $\dis[t]=\langle 3,g\rangle\conc\langle \dis[s_\nu]:
\nu\in\hat{T}\rangle$. 
\end{enumerate}
\end{enumerate}
\end{definition}

\section{Old reals are unbounded}
Recall that a forcing notion $\p$ is {\em almost $\baire$--bounding} if for
every $\p$--name $\dot{f}$ for an element of $\baire$ and any $p\in\p$ we have
\[(\exists g\in\baire)(\forall A\in\iso)(\exists q\geq p)(q\forces_{\p}
\mbox{``}(\exists^\infty n\in A)(\dot{f}(n)<g(n))\mbox{''}).\]
Almost $\baire$--bounding forcing notions do not add dominating reals (i.e.
they force that ``$(\forall x\in\baire)(\exists y\in\baire\cap\V)
(\exists^\infty n)(x(n)<y(n))$''). If $\q$ is a forcing notion not adding
dominating reals and
\[\forces_{\q}\mbox{``}\dot{\p}\mbox{ is almost $\baire$--bounding''}\]
then the composition $\q*\dot{\p}$ does not add a dominating real (see
\cite[Ch VI, 3.6]{Sh:f}). Thus the notion of ``being almost
$\baire$--bounding'' is very useful from the point of view of iterations: in a
countable support iteration of proper forcing notions no dominating reals are
added at limit stages (see \cite[Ch VI, 3.17]{Sh:f}). (Note that ``not adding
dominating reals'' is not preserved by compositions.)

In the definition \ref{ab}(1) below one can think about the following
situation (explaining the name ``decision function''). Suppose that
$(K,\Sigma)$ is a creating pair, $\dot{\tau}$ is a $\q^*_{{\rm s}\infty}
(K,\Sigma)$--name for an ordinal and $p\in\q^*_{{\rm s}\infty}(K,\Sigma)$,
$N_0<\omega$ are such that $p$ approximates $\dot{\tau}$ at each $n\geq N_0$
(see \ref{essapprox}, remember \ref{deciding}). Let us define a function 
\[z:\pos(w^p,t^p_0,\ldots,t^p_{N_0-1})\times\PC(K,\Sigma)\longrightarrow
\bigcup_{k\geq N_0}\pos(w^p,t^p_0,\ldots,t^p_{k-1})\]
such that for every $v\in\pos(w^p,t^p_0,\ldots,t^p_{N_0-1})$ and $\langle
t_0',t_1',\ldots\rangle\in\PC(K,\Sigma)$ satisfying $\langle t^p_{N_0},
t^p_{N_0+1},\ldots\rangle\leq \langle t_0',t_1',\ldots\rangle$:
\begin{enumerate}
\item If $(v,t_0',t_1',\ldots)\in\q^*_{{\rm s}\infty}(K,\Sigma)$ then $z(v,
\langle t_0',t_1',\ldots\rangle)$ is the first (in a fixed ordering of
$\bigcup\limits_{n<\omega}\prod\limits_{m<n}\bH(m)$) of the shortest $v^*$
such that $v^*\!\in\!\pos(v,t_0',\ldots,t_{k-1}')$ (for some $k<\omega$) and
$(v^*,t^p_m,t^p_{m+1},\ldots)$ decides the value of $\dot{\tau}$ ($m$ is just
suitable: $m^{t^p_m}_{\dn}=m^{t_{k-1}'}_{\up}$).
\item If $(v,t_0',t_1',\ldots)\notin\q^*_{{\rm s}\infty}(K,\Sigma)$ then we
take the first $k$ such that $\nor[t_k']\leq m^{t_k'}_{\dn}$ and we ask if
there is $\langle t_0'',t_1'',\ldots\rangle\in\PC(K,\Sigma)$ such that 
\[(v,t^p_{N_0},t^p_{N_0+1},\ldots)\leq (v,t_0'',t_1'',\ldots)\in\q^*_{{\rm s}
\infty}(K,\Sigma),\quad\mbox{ and for some }\ell\leq k\]
\[z(v,\langle t_0'',t_1'',\ldots\rangle)\in\pos(v,t_0'',\ldots,t_{\ell-1}'')
\quad \mbox{ and }\quad\langle t_0',\ldots,t_{\ell-1}'\rangle=\langle t_0'',
\ldots,t_{\ell-1}''\rangle.\] 
If the answer is ``yes'' then we choose such a sequence $\langle t_0'',t_1'',
\dots\rangle$ and we let  
\[z(v,\langle t_0',t_1',\dots\rangle)=z(v,\langle t_0'',t_1'',\dots\rangle)\]
(note that this does not depend on the choice of the particular $\langle
t_0'',t_1'',\dots\rangle$; see the previous case). If the answer is ``no''
then $z(v,\langle t_0',t_1',\ldots\rangle)$ is the first element of
$\pos(v,t_0',\ldots,t_k')$.  
\end{enumerate}
This $z$ is a canonical example of a decision function for $p$, $N_0$,
$(K,\Sigma)$; we will call it $z(p,N_0,\dot{\tau})$ (assuming that
$(K,\Sigma)$ is understood).

\begin{definition}
\label{ab}
Let $(K,\Sigma)$ be a creating pair.
\begin{enumerate}
\item Let $p\in\q^*_\emptyset(K,\Sigma)$, $N_0\in\omega$. We say that a
function 
\[z:\pos(w^p,t^p_0,\ldots,t^p_{N_0-1})\times\PC(K,\Sigma)\longrightarrow
\bigcup_{k\geq N_0}\pos(w^p,t^p_0,\ldots,t^p_{k-1})\]
is {\em a decision function for $p$, $N_0$, $(K,\Sigma)$} if:
\begin{enumerate}
\item[$(*)_{\ref{ab}}$] for every $v\in\pos(w^p,t^p_0,\ldots,t^p_{N_0-1})$
and $\langle t^\prime_0,t^\prime_1,\ldots\rangle\in\PC(K,\Sigma)$ such that 
$\langle t^p_{N_0},t^p_{N_0+1},\ldots\rangle\leq\langle
t^\prime_0,t^\prime_1,\ldots\rangle$, there is $k\in\omega$ such that: 
\[z(v,\langle t^\prime_0,t^\prime_1,\ldots\rangle)\in\pos(v,t^\prime_0,
\ldots,t^\prime_{k-1})\]
and if $\langle t^p_{N_0},t^p_{N_0+1},\ldots\rangle\leq\langle t''_0,t''_1,
\ldots\rangle\in\PC(K,\Sigma)$ is such that $t''_i\sim_\Sigma t'_i$ for all
$i<k$, then  
\[z(v,\langle t''_0,t''_1,\ldots\rangle)=z(v,\langle t'_0,t'_1,\ldots
\rangle).\] 
\end{enumerate}
\item We say that $(K,\Sigma)$ is {\em of the $\AB$-type} whenever the
following two conditions are satisfied: 
\begin{enumerate}
\item[$(\circledast)_{\AB}^0$] {\em if} $(t_0,\ldots,t_{n-1})\in
\PFC(K,\Sigma)$, $k<n$ {\em then} there is $t\in\Sigma(t_0,\ldots,t_{n-1})$
such that 
\[\nor[t]\geq\min\{\nor[t_\ell]: \ell<n\}\]
and if $(w,t_0,\ldots,t_{n-1})\in\FC(K,\Sigma)$, $t'\in\Sigma(t)$,
$\nor[t']>0$, then there is $t''\in\Sigma(t_k)$ such that $\nor[t'']>0$ and 
\[(\exists u'\in\pos(w,t_0,\ldots,t_{k-1}))(\forall u''\in\pos(u,t''))
(\exists v\in\pos(w,t'))(u''\trianglelefteq v);\]

\item[$(\circledast)_{\AB}^1$] {\em if} $p\in\q^*_\emptyset(K,\Sigma)$,
$N_0\in\omega$, $\nor[t^p_i]>2$ for $i\geq N_0$ and 
\[z:\pos(w^p,t^p_0,\ldots,t^p_{N_0-1})\times\PC(K,\Sigma)\longrightarrow
\bigcup_{k\geq N_0}\pos(w^p,t^p_0,\ldots,t^p_{k-1})\]
is a decision function for $p, N_0, (K,\Sigma)$

{\em then} there are $N_1>N_0$ and $t^*\in\Sigma(t^p_{N_0},\ldots,
t^p_{N_1-1})$ such that
\[\nor[t^*]\geq\frac{1}{2}\min\{\nor[t^p_{N_0}],\ldots,\nor[t^p_{N_1-1}]\}
\qquad\mbox{ and}\]
for each $v\in\pos(w^p,t^p_0,\ldots,t^p_{N_0-1})$ and $t\in\Sigma(t^*)$ with
$\nor[t]>0$ there is $\langle t^\prime_0,t^\prime_1,\ldots\rangle\in\PC(K, 
\Sigma)$ such that $\langle t^p_{N_0},t^p_{N_0+1},\ldots\rangle\leq\langle
t^\prime_0,t^\prime_1,\ldots\rangle$ and
\begin{enumerate}
\item[($\alpha$)] if $t'_i\in\Sigma(t^p_{N_0+k},\ldots,t^p_{N_0+k+\ell})$
($i,k,\ell <\omega$) 

then $\nor[t'_i]\geq\frac{1}{2}\min\{\nor[t^p_{N_0+k}],\ldots, 
\nor[t^p_{N_0+k+\ell}]\}$, and
\item[($\beta$)] $(\exists w\in\pos(v,t))(z(v, \langle t^\prime_0,t^\prime_1,
t^\prime_2,\ldots\rangle)\trianglelefteq w)$.
\end{enumerate}
\end{enumerate}
\item We say that $(K,\Sigma)$ is {\em condensed} if for every $(w,t_0,\ldots,
t_{n-1})\in \FC(K,\Sigma)$ with $\nor[t_i]>0$ for $i<n$, and $t\in\Sigma(t_0,
\ldots,t_{n-1})$, $\nor[t]>0$, there exist $k<n$, a creature $s\in\Sigma(t_k)$
and $v\in\pos(w,t_0,\ldots,t_{k-1})$ such that 
\[\nor[s]>0\quad\mbox{ and }\quad(\forall u\in\pos(v,s))
(\exists u^*\in\pos(w,t))(u\trianglelefteq u^*).\]
\end{enumerate}
\end{definition}

\begin{remark}
\label{getzero}
Note that the condition $(\circledast)^0_{\AB}$ is easy to satisfy: e.g.~if
$(K,\Sigma)$ is omittory and has the property that for every $t\in K$:
\[\mbox{if }m_0\leq m^t_{\dn}<m^t_{\up}\leq m_1\mbox{ then }\Sigma(t\Rsh
[m_0,m_1))=\{s\Rsh [m_0,m_1): s\in\Sigma(t)\}\]
then it satisfies this requirement (the $(\circledast)^0_{\AB}$ for $(t_0, 
\ldots,t_{n-1})$, $k<n$ is witnessed by $t_k\Rsh [m^{t_0}_{\dn},
m^{t_{n-1}}_{\up})$ ). 
\end{remark}

\begin{theorem}
\label{almbound}
Suppose that $(K,\Sigma)$ is a finitary, growing and condensed creating pair
of the $\AB$--type. Then the forcing notion $\q^*_{{\rm s}\infty}(K,\Sigma)$
is almost $\baire$-bounding.  
\end{theorem}

\begin{proof}
Let us start with the following claim which will be used later too.
\begin{claim}
\label{cl35}
Let $(K,\Sigma)$ be as in the assumptions of \ref{almbound}. Suppose that
$\dot{\tau}$ is a $\q^*_{{\rm s}\infty}(K,\Sigma)$--name for a function in
$\baire$, $q\in \q^*_{{\rm s}\infty}(K,\Sigma)$ and $n\in\omega$. Then there
are a condition $p=(w^p,t^p_0,t^p_1,\ldots)\in\q^*_{{\rm s}\infty}(K,\Sigma)$
and a strictly increasing function $g\in\baire$ such that $q\leq^{{\rm s}
\infty}_n p$ and for every $\ell\in\omega$
\begin{enumerate}
\item[($\boxplus^*_\ell$)] for each $v\in\pos(w^p,t^p_0,\ldots,t^p_{n-1},
\ldots,t^p_{n+\ell-1})$ and $t\in\Sigma(t^p_{n+\ell})$ with $\nor[t]>0$ there
is $w\in\pos(v,t)$ such that the condition $(w,t^p_{n+\ell+1},t^p_{n+\ell+2},
\ldots)$ decides the value of $\dot{\tau}(g(\ell))$ and the decision is
smaller than $g(\ell+1)$.
\end{enumerate}
\end{claim}

\noindent{\em Proof of the claim:}\ \ \  We define inductively conditions
$p_\ell\in \q^*_{{\rm s}\infty}(K,\Sigma)$ and the values $g(\ell)$ for
$\ell\in\omega$ such that $g(0)=0$, $q=p_0\leq^{{\rm s}\infty}_n p_1\leq^{{\rm
s}\infty}_{n+1}\ldots\leq^{{\rm s}\infty}_{n+\ell}p_{n+\ell}\leq^{{\rm
s}\infty}_{n+\ell+1}\ldots$ and $p_{\ell+1}$, $g(\ell)$, $g(\ell+1)$ have the
property stated in ($\boxplus^*_\ell$). It should be clear that then the limit
condition $p=\lim\limits_{\ell}p_\ell$ (see \ref{propord}(3)) is as required
in the claim.

Suppose we have defined $p_\ell,g(\ell)$. Using \ref{sinfty} and
\ref{deciding} we find a condition $p^*_\ell\in\q^*_{{\rm s}\infty}(K,\Sigma)$
such that

$p_\ell\leq_{n+\ell} p^*_\ell$, $\nor[t^{p^*_\ell}_k]>2\cdot
m^{t^{p^*_\ell}_k}+2$ for all $k\geq n+\ell$ and 

$p^*_\ell$ approximates $\dot{\tau}(g(\ell))$ at each $t^{p^*_\ell}_k$ for
$k\geq n+\ell$.

\noindent Let $z_\ell=z(p^*_\ell,n+\ell,\dot{\tau}(g(\ell)))$ be the canonical
decision function as defined before \ref{ab} (remember the choice of
$p^*_\ell$). Thus
\[z_\ell:\pos(w^q,t^{p^*_\ell}_0,\ldots,t^{p^*_\ell}_{n+\ell-1})\times\PC(K,
\Sigma)\longrightarrow\bigcup_{m\geq n+\ell}\pos(w^q,t^{p^*_\ell}_0,\ldots,
t^{p^*_\ell}_{m-1})\]
is a decision function such that
\begin{quotation}
\noindent{\em if\/} $w\in\pos(w^q,t^{p^*_\ell}_0,\ldots,t^{p^*_\ell}_{n+\ell
-1})$, $(w,t^\prime_0,t^\prime_1,\ldots)\in\q^*_{{\rm s}\infty}(K,\Sigma)$,\\
$(w,t^{p^*_\ell}_{n+\ell},t^{p^*_\ell}_{n+\ell+1},\ldots)\leq (w,t^\prime_0,
t^\prime_1,\ldots)$\\
{\em then\/} $z_\ell(w,\langle t^\prime_0,t^\prime_1,\ldots\rangle)\in
\pos(w,t^\prime_0,t^\prime_1,\ldots,t^\prime_{m-1})$ (for some $m\in\omega$)
is such that the condition $(z_\ell(w,\langle t^\prime_0,t^\prime_1,\ldots
\rangle),t^{p^*_\ell}_M,t^{p^*_\ell}_{M+1},\ldots)$ gives a value to
$\dot{\tau}(g(\ell))$, where $M$ is such that $m^{t^\prime_{m-1}}_{\up}=
m^{t^{p^*_\ell}_{M-1}}_{\up}$.
\end{quotation}
Now apply $(\circledast)_{\AB}^1$ to $p^*_\ell$, $n+\ell$ and $z_\ell$ to find
$N>n+\ell$ and $t^*\in\Sigma(t^{p^*_\ell}_{n+\ell},\ldots,t^{p^*_\ell}_{N-1})$
such that  
\[\nor[t^*]\geq\frac{1}{2}\min\{\nor[t^{p^*_\ell}_{n+\ell}],\ldots,\nor[
t^{p^*_\ell}_{N-1}]\}>m^{t^{p^*_{n+\ell}}}_{\dn}+1\]
and for each $v\in\pos(w^q,t^{p^*_\ell}_0,\ldots,t^{p^*_\ell}_{n+\ell-1})$ and
$t\in\Sigma(t^*)$ with $\nor[t]>0$ there is $w\in\pos(v,t)$ for which $(w,
t^{p^*_\ell}_N,t^{p^*_\ell}_{N+1},\ldots)$ decides $\dot{\tau}(g(\ell))$. For
this note that if $\langle t'_0,t'_1,\ldots\rangle\in \PC(K,\Sigma)$ is given
by $(\circledast)^1_{\AB}$ for $z_\ell,t,v$ then, as for some $k_i\leq\ell_i<
\omega$
\[\begin{array}{l}
t'_i\in\Sigma(t^{p^*_\ell}_{n+\ell+{k_i}},\ldots,t^{p^*_\ell}_{n+\ell+
\ell_i}),\qquad\mbox{ and}\\
\nor[t'_i]\geq\frac{1}{2}\min\{\nor[t^{p^*_\ell}_{n+\ell+k_i}],\ldots,
\nor[t^{p^*_\ell}_{n+\ell+\ell_i}]\}>m^{t^\prime_i}_{\dn},
  \end{array}\] 
the condition $(v,t'_0,t'_1,\ldots)\in \q^*_{{\rm s}\infty}(K,\Sigma)$ is
stronger than $(v,t^{p^*_\ell}_{n+\ell},t^{p^*_\ell}_{n+\ell+1},\ldots)$. Thus
our requirements on $z_k$ apply. Finally we define 

$p^*_{\ell+1}=(w^q,t^{p^*_\ell}_0,\ldots,t^{p^*_\ell}_{n+\ell-1},t^*,
t^{p^*_\ell}_N,t^{p^*_\ell}_{N+1},\ldots)$\qquad and
\[\begin{array}{l}
g(\ell+1)=1+g(\ell)+\\
+\max\{i\!<\!\omega: (\exists v\!\in\!\pos(w^q,t^{p^*_\ell}_0,
\ldots t^{p^*_\ell}_{n+\ell-1},t^*))((v,t^{p^*_\ell}_N,t^{p^*_\ell}_{N+1},
\ldots) \forces \dot{\tau}(g(\ell))=i)\}.
  \end{array}\]
Clearly they are as required. This finishes the inductive construction and the
proof of the claim. 
\medskip

Now we are going to show that $\q^*_{{\rm s}\infty}(K,\Sigma)$ is almost
$\baire$--bounding. For this suppose that $\dot{\tau}$ is a name for a
strictly increasing function in $\baire$ and $q\in \q^*_{{\rm
s}\infty}(K,\Sigma)$. Applying claim \ref{cl35} to $\dot{\tau},q$ and $n=0$ we
get a condition $p\geq q$ and an increasing function $g\in\baire$ as there (so
they satisfy ($\boxplus^*_\ell$) for $\ell\in\omega$). Note that, as
$\dot{\tau}$ is (forced to be) increasing, for every $\ell\in\omega$ we have
\begin{quotation}
\noindent if $v\in\pos(w^p,t^p_0,\ldots,t^p_{\ell-1})$ and $t\in\Sigma(
t^p_\ell)$ is such that $\nor[t]>0$ 

\noindent then $(w,t^p_{\ell+1},t^p_{\ell+2},\ldots)\forces$``$\dot{\tau}
(\ell)<g(\ell+1)$'', for some $w\in\pos(v,t)$. 
\end{quotation}
We will be done when we show the following claim.
\begin{claim}
\label{cl11}
For each $A\in\iso$ there is $p'\geq p$ such that
\[p'\forces_{\q^*_{{\rm s}\infty}(K,\Sigma)} (\exists^\infty k\in A)
(\dot{\tau}(k)<g(k+1)).\] 
\end{claim}

\noindent{\em Proof of the claim:}\ \ \  Let $A\in\iso$. Choose $0=n_0<n_1<
\ldots<\omega$ and creatures $t_i\in\Sigma(t^p_{n_i},\ldots,t^p_{n_{i+1}-1})$
and $k_i\in A$ such that
\begin{enumerate}
\item $n_i\leq k_i<n_{i+1}$, $\nor[t_i]\geq\min\{\nor[t^p_k]: n_i\leq
k<n_{i+1}\}$,
\item if $w\in\pos(w^p,t^p_0,\ldots,t^p_{n_i-1})$, $t'\in\Sigma(t_i)$,
$\nor[t']>0$ then there is $t''\in\Sigma(t^p_{k_i})$ such that $\nor[t'']>0$
and  
\[(\exists u'\in\pos(w,t^p_{n_i},\ldots,t^p_{k_i-1}))(\forall u''\in
\pos(u',t''))(\exists v\in\pos(w,t'))(u''\trianglelefteq v)\]
\end{enumerate}
(possible by $(\circledast)^0_{\AB}$). Now let $p'=(w^p,t_0,t_1,\ldots)$.
Plainly, $p'\in\q^*_{{\rm s}\infty}(K,\Sigma)$, $p\leq p'$. We want to show
that 
\[p'\forces_{\q^*_{{\rm s}\infty}(K,\Sigma)} (\exists^\infty k\in
A)(\dot{\tau}(k)<g(k+1)).\]
So assume not. Thus we find a condition $p^+\geq p'$ and $k^*\in\omega$ such
that 
\[p^+\forces_{\q^*_{{\rm s}\infty}(K,\Sigma)}(\forall k\in A)(k^*\leq k\ \ \
\Rightarrow\ \ \ \dot{\tau}(k)\geq g(k+1)).\]
Let $i\in\omega$ be such that $w^{p^+}\in\pos(w^p,t_0,t_1,\ldots,t_{i-1})$. As
we may pass to an extension of $p^+$ we may assume that $k^*<k_i$ and
$\lh(w^{p^+})>1$. Let $j>i$ be such that $t^{p^+}_0\in\Sigma(t_i,\ldots,
t_{j-1})$. Since $(K,\Sigma)$ is condensed we find $\ell\in [i,j)$, $s\in
\Sigma(t_\ell)$ and $v\in\pos(w^{p^+},t_i,\ldots,t_{\ell-1})$ such that 
\[\nor[s]>0\quad\mbox{ and }\quad(\forall u\in\pos(v,s))(\exists
u^*\in\pos(w^{p^+},t^{p^+}_0))(u\trianglelefteq u^*).\]
Now look at our choice of $t_\ell$: since $v\in\pos(w^p,t^p_0,\ldots,
t^p_{n_\ell-1})$ and $s\in\Sigma(t_\ell)$, $\nor[s]>0$, therefore we find
$t''\in\Sigma(t^p_{k_\ell})$ and $u'\in\pos(v,t^p_{n_\ell},\ldots,
t^p_{k_\ell-1})$ such that $\nor[t'']>0$ and
\[(\forall u''\in\pos(u',t''))(\exists u^+\in\pos(v,s))(u''\trianglelefteq
u^+).\] 
But now look at the choice of $p$: since $u'\in\pos(w^p,t^p_0,\ldots,t^p_{
k_\ell-1})$, $t''\in\Sigma(t^p_{k_\ell})$ and $\nor[t'']>0$, we find
$w\in\pos(u',t'')$ such that the condition $(w,t^p_{k_\ell+1},t^p_{k_\ell+2},
\ldots)$ forces that $\dot{\tau}(k_\ell)<g(k_\ell+1)$. But now, going
back, we know that there is $u^+\in\pos(v,s)$ such that $w\trianglelefteq
u^+$. Further we find $u^*\in\pos(w^{p^+}, t^{p^+}_0)$ such that
$u^+\trianglelefteq u^*$. So look at the condition $(u^*,t^{p^+}_1,t^{p^+}_2,
\ldots)$. It is stronger than $p^+$ and it forces that $\dot{\tau}(k_\ell)<
g(k_\ell+1)$, contradicting the choice of $p^+$ and $k^*<k_i\leq k_\ell$. This
finishes the proof of the claim and the theorem. 
\end{proof}

\begin{lemma}
\label{getnorm}
Suppose that $(K,\Sigma)$ is a strongly finitary and omittory creating pair
with the weak Halving Property which is saturated with respect to nice 
pre--norms with values in $\omega$ (see \ref{saturated}(2)). Further suppose
that for each $t\in K$ 
\begin{enumerate}
\item[($\otimes_1$)] $(\forall s\in\Sigma(t))(\val[s]\subseteq\val[t])$ and
\item[($\otimes_2$)] if $m_0\leq m^t_{\dn}<m^t_{\up}\leq m_1$, $s\in\Sigma(
\uhalf(t))$ 

\noindent then $s\Rsh [m_0,m_1)\in\Sigma(\uhalf(t\Rsh [m_0,m_1)))$.
\end{enumerate}
Assume that $p\in\q^*_\emptyset(K,\Sigma)$, $N_0\in\omega$, $\nor[t^p_i]>2$
for $i\geq N_0$, $m\geq 1$ and
\[z:\pos(w^p,t^p_0,\ldots,t^p_{N_0-1})\times\PC(K,\Sigma)\longrightarrow
\bigcup_{n\geq N}\pos(w^p,t^p_0,\ldots,t^p_{n-1})\]
is a decision function for $p,N_0,(K,\Sigma)$.\\
{\em Then} there are a nice pre--norm $H:\fsuo\longrightarrow\omega$ (so
$(K,\Sigma)$ is saturated with respect to $H$) and $d\in D_H$ (see
\ref{saturated}) such that $u_d=[N_0,N_1)$, $H([N_0,N_1))\geq m$ and
\begin{enumerate}
\item[if] $t\in\Sigma\big(\Sigma^{\bsum}_{d,u_d}(\uhalf(t_{N_0}^p),\ldots,
\uhalf(t_{N_1-1}^p))\big)$, $v\in\pos(w^p,t^p_0,\ldots,t^p_{N_0-1})$, and
$\nor[t]>0$
\item[then] there is $\langle t_0',t_1',\ldots\rangle\in\PC(K,\Sigma)$
such that for each $i\in\omega$, for some $k\leq\ell<\omega$ 
\[t'_i\in\Sigma(t^p_{N_0+k},\ldots,t^p_{N_0+\ell})\quad\&\quad \nor[t'_i]\geq
\frac{1}{2}\min\{\nor[t^p_{N_0+k}],\ldots,\nor[t^p_{N_0+\ell}]\}\] 
\ and $(\exists w\in\pos(v,t))(z(v,\langle t_0',t_1',\ldots\rangle)
\trianglelefteq w)$. 
\end{enumerate}
\end{lemma}

\begin{proof}
This is essentially \cite[2.14]{Sh:207}.\\
First note that if $i_0<\ldots<i_k$, $j\leq j_0<\ldots<j_\ell<\omega$, $\{i_0,
\ldots,i_k\}\subseteq \{j_0,\ldots,j_\ell\}$ (and $k\leq\ell<\omega$), $w\in
\pos(w^p,t^p_0,\ldots,t^p_{j-1})$ and $s_{j_n}\in\Sigma(t^p_{j_n})$ (for
$n\leq\ell$) then  
\[\langle w,v\rangle\in\val[\Sigma^{\bsum}(s_{i_0}\Rsh[m^{t^p_j}_{\dn},
m^{s_{i_0}}_{\up}),s_{i_1},\ldots,s_{i_k})]\]
implies
\[\langle w,v\conc{\bf 0}_{[m^{s_{i_k}}_{\up},m^{s_{i_\ell}}_{\up})}\rangle
\in\val[\Sigma^{\bsum}(s_{j_0}\Rsh[m^{t^p_j}_{\dn},
m^{s_{j_0}}_{\up}),s_{j_1},\ldots,s_{j_\ell})].\] 
Why? Suppose that $\langle w,v\rangle\in\val[\Sigma^{\bsum}(s_{i_0}
\Rsh[m^{t^p_j}_{\dn},m^{s_{i_0}}_{\up}),s_{i_1},\ldots,s_{i_k})]$. If
$j_0=i_0$ this immediately implies $\langle w,v\rest m^{s_{j_0}}_{\up}\rangle
\in\val[s_{j_0}\Rsh[m^{t^p_j}_{\dn},m^{s_{j_0}}_{\up})]$. Otherwise
necessarily $j_0\!<\!i_0$ and $v\rest [m^{t^p_j}_{\dn},m^{s_{j_0}}_{\up})$ is
constantly zero. Now, $w\in\pos(w^p,t^p_0,\ldots,t^p_{j-1})$, so using
the assumptions that $(K,\Sigma)$ is omittory and ($\otimes_1$) we get 
$\langle w,v\rest m^{s_{j_0}}_{\up}\rangle\in\val[s_{j_0}\Rsh[m^{t^p_j}_{\dn},
m^{s_{j_0}}_{\up})]$. Proceeding in this fashion further we get the desired
conclusion. 

Note that above we use ``$\val[\Sigma^{\bsum}(\ldots)]$'' and not
``$\pos(w,\Sigma^{\bsum}(\ldots))$'' as we do not claim that $\Sigma^{\bsum}
(\ldots)$ is in $K$. 

Let us define the function $H:\fsuo\longrightarrow\omega$ by:  
\begin{description}
\item[$H(u)\geq 0$] always,
\item[$H(u)\geq 1$] if $|u\setminus N_0|>1$, $u\setminus N_0=\{i_0,\ldots,
i_{k-1}\}$ (the increasing enumeration), and\\
{\em if} $s_\ell\in\Sigma(\uhalf(t^p_{i_\ell}))$, $\nor[s_\ell]>0$ (for $\ell<
k$) and $v\in\pos(w^p,t^p_0,\ldots,t^p_{N_0-1})$\\
{\em then} there exists $\langle t'_0,t'_1,\ldots\rangle\in\PC(K,\Sigma)$ such
that for each $i\in\omega$, for some $k\leq\ell<\omega$
\[t'_i\in\Sigma(t^p_{N_0+k},\ldots,t^p_{N_0+\ell})\quad\&\quad \nor[t'_i]\geq
\frac{1}{2}\min\{\nor[t^p_{N_0+k}],\ldots,\nor[t^p_{N_0+\ell}]\}\] 
and for some $w$
\[\langle v,w\rangle\in\val[\Sigma^{\bsum}(s_0\Rsh[m^{t^p_{N_0}}_{\dn},
m^{t^p_{i_0}}_{\up}),s_1,\ldots,s_{k-1})]\quad\&\quad z(v,\langle t'_0,t'_1,
t'_2,\ldots\rangle)\trianglelefteq w,\]  
\item[$H(u)\geq n+1$] if for every $u'\subseteq u$ either $H(u')\geq n$ or
$H(u\setminus u')\geq n$ (for $n>0$).
\end{description}
Note that this defines correctly a nice pre--norm on $\fsuo$; for monotonicity
use the remark we started with. Thus $(K,\Sigma)$ is saturated with respect to
$H$.  

Now it is enough to find $N_1>N_0$ such that $H([N_0,N_1))\geq m$ and then
take $d\in D_H$ with $u_d=[N_0,N_1)$. Why? Suppose that we have such an $N_1$
(and the respective $d$) and let $t\in\Sigma(\Sigma^{\bsum}_{d,u_d}(
\uhalf(t^p_{N_0}),\ldots,\uhalf(t^p_{N_1-1})))$, $\nor[t]>0$. By
\ref{saturated}(2) (the second demand) we have that there are $d^*\in D_H$ and
$s_i\in\Sigma(\uhalf(t^p_i))$ (for $i\in [N_0,N_1)$) such that
$u_{d^*}\subseteq u_d$ and $\nor[s_i]>0$ for $i\in u_{d^*}$, and $\val[
\Sigma^{\bsum}_{d^*,u_{d^*}}(s_i: N_0\leq i<N_1)]\subseteq\val[t]$. But now
look at the definition of $H$ (and remember the definitions of
$\Sigma^{\bsum}$, $\Sigma^{\bsum}_{d^*,u_{d^*}}$; note $H(u_{d^*})>0$).  

As $H$ is monotonic, it is enough to find a set $u\in [\omega\setminus
N_0]^{\textstyle {<}\omega}$ with $H(u)\geq m$. We will do this by induction
on $m$ for all $p,N_0,z$.
\medskip

\noindent{\sc Case 1:}\ \ \ $m=1$\\
For $t\in\Sigma(\uhalf(t^p_{N_0+i}))$ with $\nor[t]>0$, $i\in\omega$ fix
$s(t)\in \Sigma(t^p_{N_0+i})$ such that
\[\nor[s(t)]\geq\frac{1}{2}\nor[t]\quad\mbox{ and }\quad(\forall w\in
\basis(t))(\pos(w,s(t))\subseteq\pos(w,t))\]
(possible by \ref{halving}(2)(b$^-$)). We may additionally require that if
$t\sim_\Sigma t'$ then $s(t)=s(t')$. Let 
\[{\mathcal X}\stackrel{\rm def}{=}\{\langle t'_0,t'_1,\ldots\rangle\in
\PC(K,\Sigma): (\forall i\in\omega)(t'_i\in\Sigma(\uhalf(t^p_{N_0+i}))\ \ \&\
\ \nor[t^\prime_i]>0)\}.\]
As $(K,\Sigma)$ is strongly finitary each $\Sigma(\uhalf(t^p_{N_0+i}))$ is
finite (up to $\sim_\Sigma$--equivalence, but we may consider representatives
only) and thus the space $\mathcal X$ equipped with the product topology (of
discrete $\Sigma(\uhalf(t^p_{N_0+i}))$'s) is compact.\\
For $w\in\pos(w^p,t^p_0,\ldots,t^p_{N_0-1})$ and $\langle t'_0,t'_1,\ldots
\rangle\in {\mathcal X}$ let $M(w,\langle t'_0,t'_1,\ldots\rangle)$ be the
unique $M<\omega$ such that   
\[z(w,\langle s(t'_0),s(t'_1),\ldots\rangle)\in \pos(w,s(t'_0),\ldots,
s(t'_{M-1})).\] 
Note that the function $M:\pos(w^p,t^p_0,\ldots,t^p_{N_0-1})\times{\mathcal
X}\longrightarrow\omega$ is continuous. Why? Look at the definition of
decision functions: if $\langle t''_0,t''_1,\ldots\rangle\in{\mathcal X}$ is
such that $t''_i\sim_\Sigma t'_i$ for all $i<M(w,\langle t'_0,t'_1,\ldots
\rangle)$ then
\[z(w,\langle s(t'_0),s(t'_1),\ldots\rangle)=z(w,\langle s(t''_0),s(t''_1),
\ldots\rangle).\]
Hence, by compactness of $\mathcal X$, the function $M$ is bounded. Let 
$N'_1>1$ be such that $\rng(M)\subseteq N'_1$ and let $N_1=N_0+N'_1$. We want
to show $H([N_0,N_1))\geq 1$. So suppose that $s_\ell\in\Sigma(\uhalf(t^p_{
N_0+\ell}))$, $\nor[s_\ell]>0$ (for $\ell<N'_1$) and $v\in\pos(w^p,t^p_0,
\ldots,t^p_{N_0-1})$. Look at $\langle s_0,\ldots,s_{N_1'-1},\uhalf(t^p_{
N_1}),\uhalf(t^p_{N_1+1}),\ldots\rangle\in {\mathcal X}$, and let 
\[\bar{\bf s}=\langle s(s_0),\ldots,s(s_{N_1'-1}), s(\uhalf(t^p_{N_1})),
s(\uhalf(t^p_{N_1+1})),\ldots\rangle\in\PC(K,\Sigma).\] 
By the choice of $N'_1$ we know that for some $k<N'_1$:
\[z(v,\bar{\bf s})\in\pos(v,s(s_0),\ldots,s(s_k))\subseteq\pos(v,s_0,\ldots,
s_k)\] 
(by the choice of $s(t)$'s). Take $w^*\in\pos(v,s_0,\ldots,s_{N_1'-1})$ such
that $z(v,\bar{\bf s})\trianglelefteq w^*$. Applying ($\otimes_1$) we get
$\langle v,w^*\rangle\in\val[\Sigma^{\bsum}(s_0,\ldots,s_{N_1'-1})]$. To
finish this case note that $\nor[s(s_\ell)]\geq \frac{1}{2}\nor[t^p_{N_0+
\ell}]$ for $\ell<N'_1$ and $\nor[s(\uhalf(t^p_k))]\geq \frac{1}{2}\nor[
t^p_k]$ for $k\geq N_1$.
\smallskip

\noindent{\sc Case 2:}\ \ \ $m'=m+1\geq 2$\\
Now suppose that we always can find a finite subset of $\omega$ of the
pre--norm $H$ at least $m$. Thus we find an increasing sequence
$N_0=\ell_0<\ell_1<\ldots<\omega$ such that $H([\ell_i,\ell_{i+1}))\geq m$ for
each $i$. Consider the space of all increasing functions $\psi\in\baire$ such
that $\psi\rest N_0$ is the identity and $(\forall i\in\omega)(\psi(N_0+i)\in
[\ell_i,\ell_{i+1}))$ - it is a compact space. For each $\psi$ from the space
we may consider a condition
\[(w^p,t^p_0,\ldots,t^p_{N_0-1},t^p_{\psi(N_0)}\Rsh [m^{t^p_{\ell_0}}_{\dn},
m^{t^p_{\psi(N_0)}}_{\up}), t^p_{\psi(N_0+1)}\Rsh [m^{t^p_{\psi(N_0)}}_{\up},
m^{t^p_{\psi(N_0+1)}}_{\up}),\ldots)\]
(and call it $p_\psi$) and the respective pre--norm $H_\psi$ defined like $H$
but for $p_\psi$, $N_0$, $z$ (note that $p\leq p_\psi\in\q^*_\emptyset(K,
\Sigma)$, so $z$ may be interpreted as a decision function for $p_\psi$). 

\begin{claim}
\label{cl19}
For each finite set $u\subseteq\omega$: 
\[H_\psi(u)\leq H(\{\psi(k): k\in u\}).\] 
\end{claim}

\noindent{\em Proof of the claim:}\ \ \ Suppose $H_\psi(u)\geq 1$. We may
assume that $u\subseteq\omega\setminus N_0$. Let $s_k\in\Sigma(\uhalf(
t^p_{\psi(k)}))$, $\nor[s_k]>0$ for $k\in u$ and $v\in\pos(w^p,t^p_0,\ldots,
t^p_{N_0-1})$. By the assumption ($\otimes_2$) we know that  
\[s^*_k\stackrel{\rm def}{=}s_k\Rsh [m^{t^{p_\psi}_k}_{\dn},
m^{t^{p_\psi}_k}_{\up})\in\Sigma\big(\uhalf(t^p_{\psi(k)}\Rsh
[m^{t^{p_\psi}_k}_{\dn},m^{t^{p_\psi}_k}_{\up}))\big)=\Sigma(\uhalf(
t^{p_\psi}_k)).\] 
So, as $H_\psi(u)\geq 1$, we find $\langle t'_0,t'_1,\ldots\rangle\in\PC(K,
\Sigma)$ such that 
\begin{enumerate}
\item $t'_i\in\Sigma(t^{p_\psi}_{k_i},\ldots,t^{p_\psi}_{n_i})$
(for some $N_0\leq k_i\leq n_i<\omega$),  
\item $\nor[t'_i]\geq\frac{1}{2}\min\{\nor[t^{p_\psi}_{k_i}],\ldots,\nor[
t^{p_\psi}_{n_i}]\}$, and
\item for some $w$
\[\langle v,w\rangle\in\val\big[\Sigma^{\bsum}(s_k^*:k\in u)\Rsh
[m^{t^{p_\psi}_{N_0}}_{\dn}, m^{t^{p_\psi}_{\max(u)}}_{\up})\big]\quad \&\quad
z(v,\langle t'_0,t'_1,\ldots\rangle)\trianglelefteq w.\]
\end{enumerate}
But then also $\langle v,w\rangle\in\val\big[\Sigma^{\bsum}(s_k:k\in u)\Rsh
[m^{t^p_{N_0}}_{\dn}, m^{t^p_{\psi(\max(u))}}_{\up})\big]$. As  
\[\min\{\nor[t^{p_\psi}_{k_i}],\ldots,\nor[t^{p_\psi}_{n_i}]\}\geq 
\min\{\nor[t^p_{\psi(k_i-1)+1}],\ldots,\nor[t^p_{\psi(n_i)}]\}\]
and $t'_i\in\Sigma(t^p_{\psi(k_i-1)+1},\ldots,t^p_{\psi(n_i)})$, we may
conclude that $1\leq H(\{\psi(k): k\in u\})$.  

Next we easily proceed by induction, finishing the proof of the claim.
\medskip

By the induction hypothesis for each suitable function $\psi$ we find $N_\psi
>N_0$ such that $H_\psi([N_0,N_\psi))\geq m$. By the compactness of the space
of all these functions we find one $n$ such that $H_\psi([N_0,N_0+n))\geq m$
(for each $\psi$). Look at the interval $[N_0,\ell_n)$ -- we claim that its
$H$-norm is greater or equal than $m+1$. Why? By the choice of $\ell_i$'s we
have that $H([N_0,\ell_n))\geq m$. Suppose that $u\subseteq [N_0,\ell_n)$. If
$u\cap [\ell_k,\ell_{k+1})\neq\emptyset$ for each $k<n$ then we may take a
function $\psi$ from our space such that $\psi[[N_0,N_0+n)]\subseteq u$. But
$H_\psi([N_0,N_0+n))\geq m$ and by \ref{cl19} 
\[m\leq H(\{\psi(k): N_0\leq k<N_0+n\})\leq H(u).\]
If $u\cap[\ell_k,\ell_{k+1})=\emptyset$ for some $k<n$ then necessarily 
\[m\leq H([\ell_k,\ell_{k+1}))\leq H([N_0,\ell_n)\setminus u).\]
This finishes the induction and the proof of the lemma. 
\end{proof}

\begin{theorem}
\label{thmalmbound}
Assume that $(K,\Sigma)$ is a strongly finitary and omittory creating pair
with the weak Halving Property. Further suppose that $(K,\Sigma)$ is saturated
with respect to nice pre--norms with values in $\omega$ and for each $t\in K$:
\begin{enumerate}
\item[($\otimes_1$)] $(\forall s\in\Sigma(t))(\val[s]\subseteq\val[t])$,
\item[($\otimes_2$)] if $m_0\leq m^t_{\dn}<m^t_{\up}\leq m_1$, $s\in\Sigma(
\uhalf(t))$\\ 
then $s\Rsh [m_0,m_1)\in\Sigma(\uhalf(t\Rsh [m_0,m_1)))$,
\item[($\oplus_3$)]  $\Sigma(t\Rsh [m_0,m_1))=\{s\Rsh [m_0,m_1): s\in\Sigma(t)
\}$. 
\end{enumerate}
Then the creating pair $(K,\Sigma)$ is of the $\AB$--type. Consequently if
$(K,\Sigma)$ is additionally  condensed then the forcing notion $\q^*_{{\rm
s}\infty}(K,\Sigma)$ is almost $\baire$-bounding. 
\end{theorem}

\begin{proof}
By \ref{getzero}, \ref{getnorm} and \ref{almbound}.
\end{proof}

\begin{remark}
The assumptions of \ref{thmalmbound} may look very complicated, but in the
real examples they are relatively easy to check and appear naturally. 
Sometimes it is easy to check directly that a creating pair is of the
$\AB$--type, but then it may happen that it is not condensed (this happens
e.g.~for $(K_{\ref{blsh}},\Sigma_{\ref{blsh}})$; see \ref{blshrev}). To get
that $\q^*_{{\rm s}\infty}(K,\Sigma)$ is almost $\baire$--bounding we do not
have to require that $(K,\Sigma)$ is condensed, but then we should strengthen
the demands of \ref{ab}(2) a little bit.
\end{remark}

\begin{definition}
\label{abplus}
We say that a creating pair $(K,\Sigma)$ is {\em of the $\AB^+$--type} if it
satisfies the demand $(\circledast)_{\AB}^0$ of \ref{ab} and the following
strengthening of $(\circledast)_{\AB}^1$:
\begin{enumerate}
\item[$(\circledast^+)_{\AB}^1$] {\em if} $p\in\q^*_{{\rm s}\infty}(K,\Sigma)$,
$N_0\in\omega$, and 
\[z:\pos(w^p,t^p_0,\ldots,t^p_{N_0-1})\times\PC(K,\Sigma)\longrightarrow
\bigcup_{k\geq N_0}\pos(w^p,t^p_0,\ldots,t^p_{k-1})\]
is a decision function for $p, N_0, (K,\Sigma)$ {\em then} there are $N_1>N_0$
and $t^*\in\Sigma(t^p_{N_0},\ldots, t^p_{N_1-1})$ such that
\[\nor[t^*]\geq\frac{1}{2}\min\{\nor[t^p_{N_0}],\ldots,\nor[t^p_{N_1-1}]\}\]
and for each $v\in\pos(w^p,t^p_0,\ldots,t^p_{N_0-1},t^*)$ such that 
\[(\exists k\in [m^{t^*}_{\dn},m^{t^*}_{\up}))(v(k)\neq 0)\]
there is $\langle t^\prime_0,t^\prime_1,\ldots\rangle\in\PC_{{\rm s}\infty}(
K,\Sigma)$ such that $\langle t^p_{N_0},t^p_{N_0+1},\ldots\rangle\leq\langle
t^\prime_0,t^\prime_1,\ldots\rangle$ and $z(v\rest m^{t^*}_{\dn}, \langle
t^\prime_0,t^\prime_1, t^\prime_2,\ldots\rangle)\trianglelefteq v$.
\end{enumerate}
\end{definition}

\begin{theorem}
\label{analbo}
Suppose that $(K,\Sigma)$ is a finitary and omittory creating pair of the
$\AB^+$--type such that for each $t\in K$
\begin{enumerate}
\item[($\oplus_0$)]  if $\nor[t]>1$ and $u\in\basis(t)$ then $|\pos(u,t)|>2$
and 
\item[($\oplus_3$)]  $\Sigma(t\Rsh [m_0,m_1))=\{s\Rsh [m_0,m_1): s\in\Sigma(t)
\}$. 
\end{enumerate}
Then the forcing notion $\q^*_{{\rm s}\infty}(K,\Sigma)$ is almost
$\baire$--bounding. 
\end{theorem}

\begin{proof}
It is fully parallel to \ref{almbound}. First one proves that 
\begin{claim}
\label{cl40}
If $\dot{\tau}$ is a $\q^*_{{\rm s}\infty}(K,\Sigma)$--name for an element of 
$\baire$, $q\in \q^*_{{\rm s}\infty}(K,\Sigma)$ and $n\in\omega$, then there
are a condition $p\in\q^*_{{\rm s}\infty}(K,\Sigma)$ and an increasing $g\in
\baire$ such that $q\leq^{{\rm s}\infty}_n p$ and for every $\ell\in\omega$
\begin{enumerate}
\item[($\boxplus^+_\ell$)] if $w\in\pos(w^p,t^p_0,\ldots,t^p_{n-1},\ldots,
t^p_{n+\ell})$ is such that $w\rest [m^{t^p_{n+\ell}}_{\dn},m^{t^p_{n+
\ell}}_{\up})\neq{\bf 0}$\\
then the condition $(w,t^p_{n+\ell+1},t^p_{n+\ell+2},\ldots)$ decides the
value of $\dot{\tau}(g(\ell))$ and the decision is smaller than $g(\ell+1)$.
\end{enumerate}
\end{claim}

\noindent{\em Proof of the claim:}\ \ \  Repeat the proof of \ref{cl35}.
\medskip

Next, assuming that $\dot{\tau}$ is a name for a strictly increasing function
in $\baire$, $n=0$, and $q\in\q^*_{{\rm s}\infty}(K,\Sigma)$, we take the
condition $p\geq_0 q$ and the function $g\in\baire$ given by \ref{cl40}. They
have the property that for each $\ell\in\omega$
\begin{quotation}
\noindent if $v\in\pos(w^p,t^p_0,\ldots,t^p_\ell)$ is such that $(\exists k\in
[m^{t^p_\ell}_{\dn},m^{t^p_\ell}_{\up}))(v(k)\neq 0)$

\noindent then $(v,t^p_{\ell+1},t^p_{\ell+2},\ldots)\forces$``$\dot{\tau}
(\ell)<g(\ell+1)$''. 
\end{quotation}
To show that for every $A\in\iso$ there is $p'\geq p$ such that
\[p'\forces_{\q^*_{{\rm s}\infty}(K,\Sigma)} (\exists^\infty k\in A)
(\dot{\tau}(k)<g(k+1))\] 
we slightly modify the proof of \ref{cl11}. So suppose $A\in\iso$. Choose
$0=n_0<n_1<\ldots<\omega$ and $k_i\in A$ such that $n_i\leq k_i<n_{i+1}$ and
let $t^{p'}_i=t^p_{k_i}\Rsh [m^{t^p_{n_i}}_{\dn},m^{t^p_{n_{i+1}}}_{\dn})$. 
Look at the condition $p'=(w^p,t^{p'}_0,t^{p'}_1,\ldots)$. Assume that
$p^+\geq p'$, $k^*\in\omega$, $\lh(w^{p^+})>1$ and $k_i>k^*$, where $i$ is
such that $w^{p^+}\in\pos(w^p,t^{p'}_0,\ldots,t^{p'}_{i-1})$. Take $j>i$ such
that $t^{p^+}_0\in\Sigma(t^{p'}_i,\ldots,t^{p'}_{j-1})$. Choose $v\in\pos(
w^{p^+},t^{p^+}_0)$ such that $v(\ell)\neq 0$ for some $\ell\in [m^{t^{
p^+}_0}_{\dn},m^{t^{p^+}_0}_{\up})$ (exists by $(\oplus_0)$). Let $i^*\in
[i,j)$ be such that $\ell\in [m^{t^{p'}_{i^*}}_{\dn},m^{t^{p'}_{i^*}}_{\up})$. 
By smoothness we know that $v\rest m^{t^{p'}_{i^*}}_{\up}\in\pos(v\rest m^{t^{
p'}_{i^*}}_{\dn}, t^{p'}_{i^*})$, and therefore, by $(\oplus_3)$ and the
choice of $t^{p'}_{i^*}$ we get 
\[v\rest m^{t^p_{k_{i^*}}}_{\dn}\in\pos(w^p,t^p_0,\ldots,t^p_{k_{i^*}})\quad
\mbox{ and }\quad m^{t^p_{k_{i^*}}}_{\dn}\leq \ell<m^{t^p_{k_{i^*}}}_{\up}.\]
Hence $(v\rest m^{t^p_{k_{i^*}}}_{\up},t^p_{k_{i^*}+1},t^p_{k_{i^*}+2},\ldots)
\forces\dot{\tau}(k_{i^*})<g(k_{i^*}+1)$, so we are done. 
\end{proof}

Let us finish this section by proving a parallel of \ref{thmalmbound} for the
tree--like forcing notions.

\begin{theorem}
\label{omitalmbound}
Suppose that $\bigcup\limits_{i\in\omega}\bH(i)$ is countable and $(K,\Sigma)$
is a t-omittory tree creating pair for $\bH$. Then the forcing notion
$\q^{\tree}_0(K,\Sigma)$ is almost $\baire$-bounding.
\end{theorem}

\begin{proof}
Let $p\in\q^{\tree}_0(K,\Sigma)$ and let $\dot{f}$ be a $\q^{\tree}_0
(K,\Sigma)$--name for a function in $\baire$. For simplicity we may assume
that for every $t\in T^p$ we have $|\pos(t)|=\omega$ or at least that above
each $\nu\in T^p$ we find may an infinite front of $T^{p^{[\nu]}}$ (compare
\ref{treebound}(2)).  

Like in \ref{bigfront}(2) we may construct a condition $q\in\q^{\tree}_0(K,
\Sigma)$ stronger than $p$ and fronts $F_n$ of $T^q$ such that for all
$n\in\omega$:
\begin{enumerate}
\item $F_n=\{\nu^s: s\in\omega^{\textstyle n+1}\}$ (just a fixed enumeration),
and for each $s\in\omega^{\textstyle n\!+\!1}$:
\item if $m\in\omega$ then $\nu^s\vartriangleleft\nu^{s\conc\langle m
\rangle}$,
\item $\{\nu^{s\conc\langle m\rangle}:m\in\omega\}$ is a front of $T^{q^{[
\nu^s]}}$ and $\nor[t^q_{\nu^s}]\geq n+1$,
\item the condition $q^{[\nu^s]}$ decides the value of $\dot{f}\rest\big(n+1+
\sum\limits_{k\leq n}s(k)\big)$.
\end{enumerate}
For $m\in\omega$ let $g(m)$ be
\[1+\max\{\ell<\omega: (\exists s\in\omega^{\textstyle {\leq}m{+}1})
\big(\sum_{k<\lh(s)} s(k)\leq 4m\ \ \&\ \ q^{[\nu^s]}\forces\mbox{`` }\dot{f}
(m)=\ell\mbox{ ''}\big)\}.\]
Let $A\in\iso$. For $s\in\fseo$ choose $m(s)\in A$ such that $\lh(s)+\sum
\limits_{k<\lh(s)} s(k)<m(s)$ and let $c(s)=s\conc\langle m(s)\rangle$. Note
that $q^{[\nu^{c(s)}]}\forces\dot{f}(m(s))<g(m(s))$. Now build inductively a
condition $p'\in\q^{\tree}_0(K,\Sigma)$ such that $q\leq^0_0 p'$, $F_0
\subseteq T^{p'}$ and for each $n\in\omega$:
\begin{quotation}
\noindent if $s\in\omega^{\textstyle 2n+1}$, $\nu^s\in F_{2n}\cap T^{p'}$ then
\[t^{p'}_{\nu^s}\in \Sigma\big(t^p_\rho: (\exists\eta\in F_{2n+1})(\nu^s
\trianglelefteq\rho\trianglelefteq\eta)\big),\quad\pos(t^{p'}_{\nu^s})
\subseteq\pos(t^q_{\nu^{c(s)}})\quad\mbox{ and}\]
$\nor[t^{p'}_{\nu^s}]\geq 2n+1$ 
\end{quotation}
(possible as $(K,\Sigma)$ is t-omittory and by the third requirement on
$F_n$'s). We claim that  
\[p'\forces_{\q^{\tree}_0(K,\Sigma)}(\exists^\infty m\in
A)(\dot{f}(m)<g(m)).\] 
To see this suppose that $p''\geq p'$, $N\in\omega$. Choose
$s\in\omega^{\textstyle 2N+1}$ such that $\nu^s\in F_{2N}\cap\dcl(T^{p''})$.
Then necessarily $\pos(t^{p'}_{\nu^s})\cap\dcl(T^{p''})\neq\emptyset$ so we
may choose $\eta\in T^{p''}$  such that some initial segment of $\eta$ is in
$\pos(t^{p'}_{\nu^s})\subseteq\pos(t^q_{\nu^{c(s)}})$ (see the construction of
$p'$). But now we conclude
\[(p'')^{[\eta]}\forces_{\q^{\tree}_0(K,\Sigma)}\dot{f}(m(s))<g(m(s))\]
what finishes the proof as $N<m(s)\in A$.
\end{proof}

\section{Examples}
Let us start with noting that the Blass--Shelah forcing notion $\q^*_{{\rm s}
\infty}(K_{\ref{blsh}}^*,\Sigma^*_{\ref{blsh}})$ is a good application of the
notions introduced in this section.

\begin{proposition}
\label{blshrev}
The creating pair $(K_{\ref{blsh}}^*,\Sigma_{\ref{blsh}}^*)$ (see the end of
the construction for \ref{blsh}) is meagering, of the $\AB^+$--type and
satisfies the demands $(\oplus_0)$, $(\oplus_3)$ of \ref{analbo}.
Consequently, the forcing notion $\q^*_{{\rm s}\infty}(K^*_{\ref{blsh}},
\Sigma^*_{\ref{blsh}})$:
\begin{enumerate}
\item[$(\alpha)$] makes ground model reals meager,
\item[$(\beta)$]  adds an unbounded real,
\item[$(\gamma)$] is almost $\baire$--bounding,
\item[$(\delta)$] does not add Cohen reals.
\end{enumerate}
\end{proposition}

\begin{proof}
To show that $(K^*_{\ref{blsh}},\Sigma^*_{\ref{blsh}})$ is meagering
assume that $(t_0,\ldots,t_{n-1})\in\PFC(K^*_{\ref{blsh}},\Sigma^*_{
\ref{blsh}})$, $t\in\Sigma^*_{\ref{blsh}}(t_0,\ldots,t_{n-1})$ and $\langle
k_i:i<n\rangle$ are such that $\nor[t_i]\geq 2$, $m^{t_i}_{\dn}\leq k_i<
m^{t_i}_{\up}$ (for $i<n$) and $\nor[t]\geq 2$. Fix $u\in\prod\limits_{i<
m^{t_0}_{\dn}}\bH(i)$. By the definition of $K^*_{\ref{blsh}}$ (see clause
$(\gamma)$ there) we know that $2\leq \nor[t]\leq\dpt^1_0(A^t_u)$. Moreover,
by the definition of $\Sigma^*_{\ref{blsh}}$ we know that no element of
$A^t_u$ is included in $\{k_i:i<n\}$ (by clause $(\alpha)$ of \ref{blsh};
remember $a\cap [m^{t_i}_{\dn},m^{t_i}_{\up})\in A^{t_i}_{u\cup (a\cap
m^{t_i}_{\dn})}\cup\{\emptyset\}$ for all $a\in A^t_u$). Consequently, if
$A^s_u=\{a\in A^t_u: a\cap\{k_i:i<n\}=\emptyset\}$, then $\dpt^1(A^s_u)\geq
\dpt^1(A^t_u)-1$ and $\dpt^1_0(A^s_u)\geq \nor[t]-1$. This determines a
condition $s\in\Sigma^*_{\ref{blsh}}(t)$ which is as required in
\ref{meagering}(1).  

It should be clear that $(K^*_{\ref{blsh}},\Sigma^*_{\ref{blsh}})$ satisfies
the conditions $(\oplus_0)$, $(\oplus_3)$ of \ref{analbo} (actually,
$(\oplus_3)$ is satisfied if interpreted ``modulo $\sim_{\Sigma^*_{
\ref{blsh}}}$'', but this makes no problems). The proof that $(K^*_{
\ref{blsh}},\Sigma^*_{\ref{blsh}})$ is of the $\AB^+$--type follows exactly
the lines of \cite[2.6]{BsSh:242} (see \cite[7.4.20]{BaJu95} too) and is left
to the reader. 

Consequently, the assertion $(\alpha)$ follows from \ref{makemeager}(1),
clause $(\beta)$ is a consequence of \ref{adunbreal} and $(\gamma)$ follows
from \ref{analbo}. To show $(\delta)$ one uses \ref{omitdecbel}, or see
\ref{presGpoint}. 

Note, that if $\dot{W}$ is interpreted as a name for an infinite subset of
$\omega$, then 
\[\forces_{\q^*_{{\rm s}\infty}(K^*_{\ref{blsh}},\Sigma^*_{\ref{blsh}})}(
\forall X\in\iso\cap\V)(|\dot{W}\cap X|<\omega\ \mbox{ or }\ |\dot{W}\setminus
X|<\omega).\]
Thus forcing with $\q^*_{{\rm s}\infty}(K^*_{\ref{blsh}},\Sigma^*_{
\ref{blsh}})$ makes ground model reals null too.
\end{proof}

Now we will present an application of forcing notions determined by
omittory creating pairs with the weak Halving Properties to questions coming
from localizing subsets of $\omega$. These problems were studied in
\cite{RoSh:501} and our example is built in a manner  similar to that of the
forcing notion constructed in \cite[2.4]{RoSh:501}. Moreover, all these
examples are relatives of the forcing notion presented in \cite{Sh:207}. The
creating pair constructed there can be build like $(K_{\ref{locex}}^\psi,
\Sigma^\psi_{\ref{locex}})$ for $\psi\equiv 1$. 

\begin{example}
\label{locex}
Let $\psi\in\baire$ be a non-decreasing function, $\psi(0)>0$. We construct a
creating pair $(K^\psi_{\ref{locex}},\Sigma^\psi_{\ref{locex}})$ which: 
\begin{enumerate}
\item[($\alpha$)] is strongly finitary, forgetful and omittory,
\item[($\beta$)]  has the weak Halving Property,
\item[($\gamma$)] is saturated with respect to nice pre-norms with values in
$\omega$,
\item[($\delta$)] is condensed and satisfies the demands ($\otimes_1$),
($\otimes_2$) and ($\oplus_3$) of \ref{thmalmbound}

\noindent [thus, by \ref{thmalmbound}, $(K^\psi_{\ref{locex}},\Sigma^\psi_{
\ref{locex}})$ is of the $\AB$--type],
\item[($\varepsilon$)] is anti-big and meagering.
\end{enumerate}
\end{example}

\begin{proof}[Construction] Let $\bH(m)=2$ for $m\in\omega$. First we describe
which creatures $t\in\CR[\bH]$ are taken to be in $K^\psi_{\ref{locex}}$. So,
$t=(\nor[t],\val[t],\dis[t])\in K^\psi_{\ref{locex}}$ if: 
\begin{itemize}
\item $\dis[t]=(T[t],L[t],R[t],D[t],\NOR[t],m^t_{\dn},m^t_{\up})$, where,
letting $T=T[t]$, $L=L[t]$, $R=R[t]$, $D=D[t]$ and $\NOR=\NOR[t]$: 
\begin{enumerate}
\item[(a)] $T$ is a finite tree, $D\subseteq\{\nu\in T: \suc_{T}(\nu)\neq
\emptyset\}$ and
\begin{enumerate}
\item[(i)] \ $(\forall\nu\in T\setminus D)(\suc_T(\nu)=\emptyset\ \mbox{ or }\
|\suc_T(\nu)|=\psi(L(\nu)))$,
\item[(ii)]  if $\nu\in T\setminus D$ and $\eta\in\suc_T(\nu)$ then either
$\eta\in D$ or $\suc_T(\eta)=\emptyset$,
\end{enumerate}
\item[(b)] $L,R: T\longrightarrow [m^t_{\dn},m^t_{\up})$ are such that for
each $\nu\in T$
\begin{enumerate}
\item[(i)] \ \ $L(\nu)\leq R(\nu)$,
\item[(ii)] \  if $\eta\in\suc_T(\nu)$ then $L(\nu)\leq L(\eta)\leq R(\eta)
\leq R(\nu)$,
\item[(iii)]   if $\suc_T(\nu)=\emptyset$ then $L(\nu)=R(\nu)$,
\item[(iv)] \  if $\eta_0,\eta_1\in\suc_T(\nu)$, $\eta_0\neq\eta_1$ then
\[[L(\eta_0),R(\eta_0)]\cap [L(\eta_1),R(\eta_1)]=\emptyset,\]
\end{enumerate}
\item[(c)] $\NOR$ is a function on $D$ such that for each $\nu\in D$,
$\NOR(\nu)$ is a nice pre-norm on $\suc_T(\nu)$ with values in $\omega$,
\end{enumerate}
\item if $\langle\rangle\notin D$ and $(\exists\nu\in\suc_T(\langle\rangle))(
\suc_T(\nu)=\emptyset)$, or $D=\emptyset$ then $\nor[t]=0$,\\
otherwise\quad $\nor[t]=\min\{\NOR(\nu)(\suc_T(\nu)): \nu\in D\}$,
\item $\val[t]=\{\langle u,v\rangle\in 2^{\textstyle m^t_{\dn}}\times
2^{\textstyle m^t_{\up}}: u\vartriangleleft v\ \&\ \{i\in [m^t_{\dn},m^t_{\up}
): v(i)=1\}\subseteq\{L(\nu): \nu\in T\ \&\ \suc_T(\nu)=\emptyset\}\}$.
\end{itemize}
For $t\in K^\psi_{\ref{locex}}$ we define 
\[\nor^0[t]=\left\{\begin{array}{ll}
\min\{\NOR[t](\nu)(\suc_{T[t]}(\nu)): \nu\in D[t]\}&\mbox{if } D[t]\neq
\emptyset,\\
0&\mbox{otherwise.} \end{array}\right.\]
Note that $\nor[t]\leq\nor^0[t]$ and in most cases they agree. One could use
$\nor^0[t]$ as the norm of $t$ and get the same forcing notion. We take
$\nor[t]$ for technical reasons only. Now we are going to describe a
composition operation $\Sigma^\psi_{\ref{locex}}$ on $K^\psi_{\ref{locex}}$
by giving basic operations which may be applied to creatures from
$K_{\ref{locex}}^\psi$. 

\noindent (1)\ \ \ For a creature $t\in K_{\ref{locex}}$ let $\uhalf(t)\in
K^\psi_{\ref{locex}}$ be such that 
\begin{itemize}
\item if $\nor[t]<2$ then $\uhalf(t)=t$,
\item if $\nor[t]\geq 2$ then $\val[\uhalf(t)]=\val[t]$, $T[\uhalf(t)]= T[t]$,
$L[\uhalf(t)]=L[t]$, $R[\uhalf(t)]=R[t]$, $D[\uhalf(t)]=D[t]$ and

if $\nu\in D[\uhalf(t)]$, $A\subseteq\suc_{T[\uhalf(t)]}(\nu)$ then
\[\NOR[\uhalf(t)](\nu)(A)=\max\{0,\NOR[t](\nu)(A)-\lfloor\frac{\nor[t]}{2}
\rfloor\}.\]
\end{itemize}
[Thus $m^{\uhalf(t)}_{\dn}=m^t_{\dn}$, $m^{\uhalf(t)}_{\up}=m^t_{\up}$ and
$\nor[\uhalf(t)]=\nor[t]-\lfloor\frac{\nor[t]}{2}\rfloor$ when $\nor[t]\geq
2$.] 

\noindent (2)\ \ \ For $t\in K^\psi_{\ref{locex}}$, $m_0\leq m^t_{\dn}$,
$m_1\geq m^t_{\up}$ let $s=S_{m_0,m_1}(t)\in K^\psi_{\ref{locex}}$ be a
creature such that $\dis[s]=(T[t],L[t],R[t],D[t],\NOR[t],m_0,m_1)$,
$\nor[s]=\nor[t]$ and 
\[\begin{array}{r}
\val[s]=\big\{\langle u,v\rangle\in 2^{\textstyle m_0}\times 2^{\textstyle
m_1}: u\vartriangleleft v\ \&\ v\rest [m_0,m^t_{\dn})={\bf 0}\ \&\ v\rest
[m^t_{\up},m_1)={\bf 0}\ \&\\
\&\ \langle v\rest m^t_{\dn},v\rest m^t_{\up}\rangle\in\val[t]\}.
  \end{array}\]
Thus, essentially, $S_{m_0,m_1}(t)=t\Rsh[m_0,m_1)$, the small difference in
the definition of $\dis$ is immaterial.

\noindent (3)\ \ \ For $t\in K_{\ref{locex}}$ let $\Sigma^*_{\ref{locex}}(t)$
consist of all $s\in K^\psi_{\ref{locex}}$ such that
\begin{enumerate}
\item[(i)]\ \ $m^s_{\dn}=m^t_{\dn}$, $m^s_{\up}=m^t_{\up}$, $T[s]\subseteq
T[t]$, $D[s]=D[t]\cap T[s]$,
\item[(ii)]\  $(\forall \nu\in T[s])(\suc_{T[s]}(\nu)=\emptyset\
\Leftrightarrow\ \suc_{T[t]}(\nu)=\emptyset)$
\noindent (thus if $\nu\in T[s]\setminus D[s]$ then $\suc_{T[s]}(\nu)=
\suc_{T[t]}(\nu)$),
\item[(iii)]  $L[s]=L[t]\rest T[s]$ and $R[s]=R[t]\rest T[s]$,
\item[(iv)]\  if $\nu\in D[s]$, $A\subseteq\suc_{T[s]}(\nu)$ then $\NOR[s](
\nu)(A)=\NOR[t](\nu)(A)$.
\end{enumerate}

\noindent (4)\ \ \ Suppose that $t_0,\ldots,t_{k-1}\in K^\psi_{\ref{locex}}$
are such that $k=\psi(L[t_0](\langle\rangle))$, $m^{t_i}_{\up}\leq m^{t_{i+1}
}_{\dn}$ (for $i<k-1$) and
\[(\forall i<k)(\langle\rangle\in D[t_i]\ \mbox{ or }\ \suc_{T[t_i]}(\langle
\rangle)=\emptyset).\]
Let $s=S^*(t_0,\ldots,t_{k-1})\in K^\psi_{\ref{locex}}$ be a creature such
that 
\begin{enumerate}
\item[(i)]\ \ $m^s_{\dn}=m^{t_0}_{\dn}$, $m^s_{\up}=m^{t_{k-1}}_{\up}$ and
$T[s]$ is a tree such that $|\suc_{T[s]}(\langle\rangle)|=k$ and for every
$\nu\in\suc_{T[s]}(\langle\rangle)$ there is a unique $i=i(\nu)<k$ such that
\[\hspace{-0.4cm}\{\eta\in T[s]\!:\nu\trianglelefteq\eta\}=\{\nu\conc\eta^*\!:
\eta^*\in T[t_i]\}\ \&\ L[s](\nu)=L[t_i](\langle\rangle)\ \&\ R[s](\nu)=R[t_i]
(\langle\rangle),\]
\item[(ii)]\ $D[s]=\{\nu\conc\eta^*: \nu\in\suc_{T[s]}(\langle\rangle)\ \&\
\eta^*\in D[t_{i(\nu)}]\}$,
\item[(iii)] $L[s](\langle\rangle)=L[t_0](\langle\rangle)$, $R[s](\langle
\rangle)=R[t_{k-1}](\langle\rangle)$ and for all $\nu\in\suc_{T[s]}(\langle
\rangle)$
\[(\forall\eta^*\in T[t_{i(\nu)}])(L[s](\nu\conc\eta^*)=L[t_{i(\nu)}](\eta^*)\
\&\ R[s](\nu\conc\eta^*)=R[t_{i(\nu)}](\eta^*)),\]
\item[(iv)]\ if $\nu\in\suc_{T[s]}(\langle\rangle)$, $\eta^*\in D[t_{i(
\nu)}]$ and $A\subseteq\suc_{T[t_{i(\nu)}]}(\eta^*)$ then 
\[\NOR[s](\nu\conc\eta^*)(\{\nu\conc\eta':\eta'\in A\})=\NOR[t_{i(\nu)}](
\eta^*)(A).\]
\end{enumerate}

\noindent (5)\ \ \ Suppose that $H:{\mathcal P}(m)\longrightarrow\omega$ is a
nice pre-norm and $t_0,\ldots,t_{m-1}\in K^\psi_{\ref{locex}}$ are such that
$m^{t_i}_{\up}\leq m^{t_{i+1}}_{\dn}$ for $i<m-1$. Let $s=S^{**}_H(t_0,\ldots,
t_{m-1})\in K^\psi_{\ref{locex}}$ be a creature such that
\begin{enumerate}
\item[(i)]\ \ $m^s_{\dn}=m^{t_0}_{\dn}$, $m^s_{\up}=m^{t_{m-1}}_{\up}$ and
$T[s]$ is a tree such that $|\suc_{T[s]}(\langle\rangle)|=m$ and for every
$\nu\in\suc_{T[s]}(\langle\rangle)$ there is a unique $j=j(\nu)<m$ such that
\[\hspace{-0.4cm}\{\eta\in T[s]\!:\nu\trianglelefteq\eta\}=\{\nu\conc\eta^*\!:
\eta^*\in T[t_j]\}\ \&\ L[s](\nu)=L[t_j](\langle\rangle)\ \&\ R[s](\nu)=R[t_j]
(\langle\rangle),\]
\item[(ii)]\ $D[s]=\{\langle\rangle\}\cup\{\nu\conc\eta^*: \nu\in\suc_{T[s]}(
\langle\rangle)\ \&\ \eta^*\in D[t_{j(\nu)}]\}$,
\item[(iii)] $L[s](\langle\rangle)=L[t_0](\langle\rangle)$, $R[s](\langle
\rangle)=R[t_{m-1}](\langle\rangle)$ and for every $\nu\in\suc_{T[s]}(\langle
\rangle)$
\[(\forall\eta^*\in T[t_{j(\nu)}])(L[s](\nu\conc\eta^*)=L[t_{j(\nu)}](\eta^*)\
\&\ R[s](\nu\conc\eta^*)=R[t_{j(\nu)}](\eta^*)),\]
\item[(iv)]\  if $A\subseteq\suc_{T[s]}(\langle\rangle)$ then $\NOR[s](\langle
\rangle)(A)=H(\{j(\nu): \nu\in A\})$,
\item[(v)]\ \ if $\nu\in\suc_{T[s]}(\langle\rangle)$, $\eta^*\in D[t_{j(\nu)
}]$ and $A\subseteq\suc_{T[t_{j(\nu)}]}(\eta^*)$ then 
\[\NOR[s](\nu\conc\eta^*)(\{\nu\conc\eta':\eta'\in A\})=\NOR[t_{j(\nu)}](
\eta^*)(A).\]
\end{enumerate}
Note that, under the respective assumptions, the procedures described in
(1)--(5) above determine creatures in $K^\psi_{\ref{locex}}$, though (in cases
(4) and (5)) not uniquely: there is some freedom in defining $\suc_{T[s]}(
\langle\rangle)$. However, this freedom becomes irrelevant when we identify
creatures that look the same. The last operation ($\Sigma^{**}_{\ref{locex}}$
below) is a way to describe which creatures are identified.

\noindent (6)\ \ \ For $t\in K^\psi_{\ref{locex}}$, let $\Sigma^{**}_{
\ref{locex}}(t)$ consist of all creatures $s\in K^\psi_{\ref{locex}}$ such that
$m^s_{\dn}=m^t_{\dn}$, $m^s_{\up}=m^t_{\up}$ and there is an (order)
isomorphism $\pi:T[s]\longrightarrow T[t]$ which preserves $L,R,D$ and $\NOR$. 

Finally, if $t_0,\ldots,t_{m-1}\in K^\psi_{\ref{locex}}$ are such that
$m^{t_i}_{\up}=m^{t_{i+1}}_{\dn}$ (for $i<m-1$) then $\Sigma^\psi_{
\ref{locex}}(t_0,\ldots,t_{m-1})$ consists of all creatures $s\in K^\psi_{
\ref{locex}}$ such that $m^s_{\dn}=m^{t_0}_{\dn}$, $m^s_{\up}=m^{t_{m-1}}_{
\up}$ and $s$ may be obtained from $t_0,\ldots,t_{m-1}$ by use of the
operations $\uhalf$, $S_{m_0,m_1}$, $S^*$, $S^{**}_H$, $\Sigma^*_{
\ref{locex}}$ and $\Sigma^{**}_{\ref{locex}}$ (with suitable parameters). 
\medskip

Let us check that $(K^\psi_{\ref{locex}},\Sigma^\psi_{\ref{locex}})$ has the
required properties. It should be clear that $(K^\psi_{\ref{locex}},
\Sigma^\psi_{\ref{locex}})$ is a finitary, forgetful and omittory creating
pair. The relation $\sim_{\Sigma^\psi_{\ref{locex}}}$ (see \ref{subcom}(3)) is
an equivalence relation on $K^\psi_{\ref{locex}}$ and $\Sigma^\psi_{
\ref{locex}}$ depends on $\sim_{\Sigma^\psi_{\ref{locex}}}$--equivalence
classes only (remember the definition of $\Sigma^{**}_{\ref{locex}}$; note
that $\Sigma^\psi_{\ref{locex}}(t_0)=\Sigma^\psi_{\ref{locex}}(t_1)$ implies
that $t_0,t_1$ are the same up to the isomorphism of the trees $T[t_0]$,
$T[t_1]$). Thus the value of $\Sigma^\psi_{\ref{locex}}(t_0,\ldots,t_{m-1})$
does not depend on the particular representation of the trees $T[t_i]$ (for
$i<m$). Hence, if $t\in \Sigma^\psi_{\ref{locex}}(t_0,\ldots,t_{m-1})$ then we
may think that $T[t]$ is a tree built of $T[t_0],\ldots,T[t_{m-1}]$  in the
following sense. There are $s_i\in\Sigma^*_{\ref{locex}}(t_i)$ (for $i<m$), a
front $F$ of $T[t]$ and a one-to-one mapping $\varphi: F\longrightarrow m$
such that for $\nu\in F$ 
\[\{\eta\in T[t]:\nu\trianglelefteq\eta\}=\{\nu\conc\eta':\eta'\in T[s_{
\varphi(\nu)}]\}\]
and the $L[t],R[t],D[t]$ above $\nu$ in $T[t]$, $\nu\in F$, look like
$L[s_{\varphi(\nu)}],R[s_{\varphi(\nu)}],D[s_{\varphi(\nu)}]$ (but the norms
given by $\NOR$ may be substantially different, still their values may be only
smaller). Now it should be clear that $(K^\psi_{\ref{locex}},\Sigma^\psi_{
\ref{locex}})$ is strongly finitary. The pair $(K^\psi_{\ref{locex}},
\Sigma^\psi_{\ref{locex}})$ has the weak Halving Property as witnessed by the
function $\uhalf$ defined in (1) above. [Why? Note that if $t_0\in K^\psi_{
\ref{locex}}$, $\nor[t_0]\geq 2$, $t\in\Sigma^\psi_{\ref{locex}}(\uhalf(
t_0))$, $\nor[t]>0$ (so $\nor[t]\geq 1$) then $t$ is obtained from $t_0$ by
alternate applications of $\uhalf$ and shrinking (i.e.~$\Sigma^*_{\ref{locex}
}$). Look at the tree $T[t]$ with $L[t],R[t],D[t]$: necessarily the last three
objects are the restrictions of $L[t_0], R[t_0], D[t_0]$ to $T[t]$. Let
$s\in\Sigma^*_{\ref{locex}}(t_0)$ be such that $T[s]=T[t]$ (it should be clear
clear that there is one; actually the $s$ is uniquely determined by the tree
$T[s]$). Since in the process of building $t$ the norms were decreased only,
and we started with $\uhalf(t_0)$, we may conclude that $\nor[s]\geq\nor[t]+
\lfloor\frac{1}{2}\nor[t_0]\rfloor\geq\frac{1}{2}\nor[t_0]$. Clearly
$\val[s]=\val[t]$.]  

Suppose that $H_0:{\mathcal P}(\omega)\longrightarrow\omega$ is a nice
pre-norm, $d\in D_{H_0}$ (see \ref{saturated}). Note that $\Sigma^{\bsum}_{d,
u_d}(t_i:n_0\leq i<n_1)$ (for $\langle t_i: n_0\leq i<n_1\rangle\in\PFC(
K^\psi_{\ref{locex}},\Sigma^\psi_{\ref{locex}})$, $u_d\subseteq [n_0,n_1)$)
corresponds to a creature $t$ obtained from $t_{n_0},\ldots,t_{n_1-1}$ by
suitable applications of $S_{m_0,m_1}$ and $S^{**}_H$ (the last for the
pre-norm $H=H_0\rest{\mathcal P}(u_d)$, the first is to omit creatures $t_i$
for $i\in [n_0,n_1)\setminus u_d$). The difference is in $\dis$, but this
causes no real problem as we may read $H_0\rest {\mathcal P}(u_d)$ from the
resulting creature (and it is essentially the $\NOR[t](\langle\rangle)$). We
could have changed the definition of $(K^\psi_{\ref{locex}},\Sigma^\psi_{
\ref{locex}})$ to make this correspondence more literal, but that would result
in unnecessary complications in the definition. Now note that  
\begin{quotation}
\noindent if $t\in\Sigma^\psi_{\ref{locex}}(\Sigma^{\bsum}_{d,u_d}(t_0,\ldots,
t_{m-1}))$, $\nor[t]>0$ and $d^*\in  D_{H_0}$ is such that $u_{d^*}\approx
\suc_{T[t]}(\langle\rangle)$ and $s_i\in\Sigma^*_{\ref{locex}}(t_i)$ for $i<m$
are such that  
\[\nu\in\suc_{T[t]}(\langle\rangle)\quad\Rightarrow\quad\{\nu\conc\eta':\eta'
\in T[s_{j(\nu)}]\}=\{\eta\in T[t]:\nu\trianglelefteq\eta\},\]
then
$\val[t]=\val[\Sigma^{\bsum}_{d^*,u_{d^*}}(s_0,\ldots,s_{m-1})]$. Moreover, if
$\nor[t_i]>0$ for $i<m$ then $\nor[s_i]>0$ for $i\in u_{d^*}$.
\end{quotation}
This shows that $(K^\psi_{\ref{locex}},\Sigma^\psi_{\ref{locex}})$ is
saturated with respect to nice pre-norms with values in $\omega$.

Suppose that $t\in\Sigma^\psi_{\ref{locex}}(t_0,\ldots,t_{n-1})$, $\nor[t]>0$,
$\nor[t_i]>0$ (for $i<n$). Then for some $i<n$ and $\nu\in T[t]$ we find $s\in
\Sigma^*_{\ref{locex}}(t_i)$ such that 
\[\{\eta\in T[t]:\nu\trianglelefteq\eta\}=\{\nu\conc\eta':\eta'\in T[s]\}\]
and $L[t],R[t]$ and $D[t]$ above $\nu$ are like $L[s],R[s],D[s]$, but remember
that $\NOR[t]$ may have nothing in common with $\NOR[s]$: the operation
$\uhalf$ may be involved. However, by the definition of $\nor[t]$, we know
that if $\eta'\in D[s]$ then 
\[\NOR[t](\nu\conc\eta')(\suc_{T[t]}(\nu\conc\eta'))>0,\]
and this is enough to conclude that $\NOR[s](\eta')(\suc_{T[s]}(\eta'))>0$. 
Moreover $D[s]\neq\emptyset$ as $\nor[t_i]>0$. Consequently $\nor[s]>0$. Since  
\[v\rest [m^{t_0}_{\dn}, m^s_{\dn})={\bf 0}\ \&\ \langle v,u\rangle\in\val[s]
\quad\Rightarrow\quad (\exists u^*)(\langle v,u^*\rangle\in\val[t]\ \&\
u\trianglelefteq u^*)\]
we conclude that $(K^\psi_{\ref{locex}},\Sigma^\psi_{\ref{locex}})$ is
condensed. 

One easily checks that the demands ($\otimes_1$), ($\otimes_2$) and
($\oplus_3$) of \ref{thmalmbound} are satisfied (remember that $t\Rsh
[m_0,m_1)$ is, basically, $S_{m_0,m_1}(t)$). 

Finally, let us check that $(K^\psi_{\ref{locex}},\Sigma^\psi_{\ref{locex}})$
is meagering and anti-big. For the first note that if $t\in K^\psi_{
\ref{locex}}$, $\nor[t]>2$ and $\nu\in T[t]$ is a maximal node of $T[t]$
then there is $s\in\Sigma^*_{\ref{locex}}(t)$ such that $\nu\notin T[s]$ and
$\nor[s]\geq \nor[t]-1$. [Why? Take the shortest $\eta\vartriangleleft\nu$
such that $\eta\in D[t]$ -- there must be one as $\nor[t]>0$ -- and choose
$s\in \Sigma^*_{\ref{locex}}$ so that $T[s]=\{\rho\in T[t]:\eta
\vartriangleleft\rho\ \Rightarrow\ \neg(\nu^*\trianglelefteq\rho)\}$, where
$\nu^*\in\suc_{T[t]}(\eta)$ is such that $\nu^*\trianglelefteq\nu$. Since
$\NOR[t](\eta)$ is a nice pre-norm we have
\[\NOR[t](\eta)(\suc_{T[s]}(\eta))\geq \NOR[t](\eta)(\suc_{T[t]}(\eta))-1\]
and hence $s$ is as required.] Now suppose that $(t_0,\ldots,t_{n-1})\in \PFC(
K^\psi_{\ref{locex}},\Sigma^\psi_{\ref{locex}})$, $\nor[t_i]>2$, $k_i\in
[m^{t_i}_{\dn},m^{t_i}_{\up})$ and $t\in \Sigma^\psi_{\ref{locex}}(t_0,\ldots,
t_{n-1})$, $\nor[t]>2$. Then there is a front $F$ of $T[t]$ such that for
every $\nu\in F$, for a unique $i=i(\nu)<n$ and some $s\in \Sigma^*_{
\ref{locex}}(t_i)$:
\[\{\eta\in T[t]:\nu\trianglelefteq\eta\}=\{\nu\conc\eta^*:\eta^*\in
T[s_i]\}\]
and $D[t],L[t],R[t]$ above $\nu$ look like $D[s_i],L[s_i],R[s_i]$ (but
$\NOR[t]$ might be different than that of $s_i$: the values may be
smaller). Now apply the previous remark to choose $s^*_i\in\Sigma^*_{
\ref{locex}}(s_i)$ such that $\nor[s^*_i]\geq\nor[s_i]-1$ and $k_i\notin\{
L[s^*_i](\eta):\eta\in T[s^*_i]\ \&\ \suc_{T[s^*_i]}(\eta)=\emptyset\}$.
Finally let $s\in\Sigma^*_{\ref{locex}}(t)$ be such that $F\subseteq T[s]$ and
for each $\nu\in F$
\[\{\eta\in T[s]:\nu\trianglelefteq\eta\}=\{\nu\conc\eta^*:\eta^*\in
T[s^*_{i(\nu)}]\}.\]
It is easy to check that this determines a creature $s$ as required in
\ref{meagering}(1) for $\langle k_i,t_i: i<n\rangle$, $t$.

To verify that $(K^\psi_{\ref{locex}},\Sigma^\psi_{\ref{locex}})$ is anti-big
define colourings 
\[c_t:\bigcup\limits_{u\in\basis(t)}\pos(u,t)\longrightarrow 3\]
for $t\in K^\psi_{\ref{locex}}$ by:
\[c_t(v)=\left\{\begin{array}{ll}
0&\mbox{ if }\ |\{k\in [m^t_{\dn},m^t_{\up}):v(k)=1\}|\ \mbox{ is even }>0,\\
1&\mbox{ if }\ |\{k\in [m^t_{\dn},m^t_{\up}):v(k)=1\}|\ \mbox{ is odd },\\
2&\mbox{ if }\ v\rest [m^t_{\dn},m^t_{\up})={\bf 0}_{[m^t_{\dn},m^t_{\up})}.
                \end{array}\right.\]
Suppose $(t_0,\ldots,t_{n-1})\in\PFC(K^\psi_{\ref{locex}},\Sigma^\psi_{
\ref{locex}})$, $t\in\Sigma^\psi_{\ref{locex}}(t_0,\ldots,t_{n-1})$,
$\nor[t_i]>1$ and $\nor[t]>1$. Clearly for some $i<n$ 
\[|\{L[t](\eta):L[t_i](\langle\rangle)\leq L[t](\eta)\leq R[t_i](\langle
\rangle)\ \mbox{ and }\ \suc_{T[t]}(\eta)=\emptyset\}|\geq 2.\]
Let $u\in\basis(t_0)=\prod\limits_{i<m^t_{\dn}}2$. Take $v_0,v_1\in\prod
\limits_{i<m^t_{\up}}2$ such that for $\ell=0,1$
\[\begin{array}{l}
u\vartriangleleft v,\qquad |\{k\in[m^t_{\dn},m^t_{\up}):
v_\ell(k)=1\}|=\ell+1,\qquad\mbox{ and}\\
\{k\in [m^t_{\dn},m^t_{\up}):v_\ell(k)=1\}\quad\mbox{ is contained in}\\
\{L[t](\eta):L[t_i](\langle\rangle)\leq L[t](\eta)\leq R[t_i](\langle\rangle)\
\&\ \eta\in T[t]\ \&\ \suc_{T[t]}(\eta)=\emptyset\}. 
  \end{array}\]
Now check that the $v_0,v_1$ are as required in \ref{meagering}(2) for
$\langle t_i: i<n\rangle$, $t$. 
\end{proof}

\begin{corollary}
\label{whatitdoes}
$\q^*_{{\rm s}\infty}(K^\psi_{\ref{locex}},\Sigma^\psi_{\ref{locex}})$ is a
proper almost $\baire$--bounding forcing notion which makes ground reals
meager and adds a Cohen real. 
\end{corollary}

\begin{proof}
By \ref{locex}, \ref{thmalmbound} and \ref{makemeager}.
\end{proof}

\begin{definition}
\label{muX}
\begin{enumerate}
\item For an infinite set $X\in\iso$ let $\mu_X:\omega\longrightarrow X$ be
the increasing enumeration of $X$.
\item Let $\psi\in\baire$ be non-decreasing. We define relations
$S^\psi_{\dn}, S^\psi_{\up}\subseteq\iso\times\iso$ by
\[\begin{array}{l}
(X,Y)\in S^\psi_{\dn}\qquad\mbox{ if and only if}\\
(\exists^\infty n\in\omega)(\forall i<\psi(\mu_Y(n)))(|[\mu_Y(n+i),
\mu_Y(n+i+1))\cap X|\geq 2),\\
(X,Y)\in S^\psi_{\up}\qquad\mbox{ if and only if}\\
(\exists^\infty n\in\omega)(\forall i<\psi(\mu_Y(n+1)))(|[\mu_Y(n+i),
\mu_Y(n+i+1))\cap X|\geq 2).
  \end{array}\]
If $\psi$ is a constant function, say $\psi\equiv k$, then $S^\psi_*$ may be
called $S_k$. 
\end{enumerate}
\end{definition}

\begin{remark}
We will consider the notions of $S^\psi_{\dn}$-- and
$S^\psi_{\up}$--localizations as given by \ref{gendef} for these relations as
well as the corresponding dominating numbers $\dominating(S^\psi_{\dn})$ and
$\dominating(S^\psi_{\up})$. Note that for a non-decreasing $\psi\in\baire$,
$S^\psi_{\up}$--localization implies $S^\psi_{\dn}$--localization (and
$\dominating(S^\psi_{\dn})\leq \dominating(S^\psi_{\up})$). If
$(\forall^\infty n\in\omega)(\psi(n)\leq\varphi(n))$ then
$S^\varphi_*$--localization implies $S^\psi_*$--localization and for
eventually constant $\psi$, $S^\psi_{\up}$--localization is the same as
$S^\psi_{\dn}$--localization.  
\end{remark}

\begin{proposition}
\label{noloc}
Suppose $\varphi,\psi\in\baire$ are non-decreasing functions such that
$(\forall^\infty n\in\omega)(1\leq\psi(n)<\varphi(n))$. Then the forcing
notion $\q^*_{{\rm s}\infty}(K^\psi_{\ref{locex}},\Sigma^\psi_{\ref{locex}})$
does not have the $S^\varphi_{\up}$--localization property. 
\end{proposition}

\begin{proof}
This is parallel to \cite[2.4.3]{RoSh:501}. The main step is done by the
following claim, which is essentially a repetition of \cite[2.4.2]{RoSh:501}. 

\begin{claim}
\label{cl34}
Suppose that $t\in K^\psi_{\ref{locex}}$ is such that $\nor^0[t]\geq 11$. Let
$Y\in\iso$. Then there is $s\in\Sigma^*_{\ref{locex}}(t)$ such that $\nor^0[s]
\geq\nor^0[t]-10$ and for every $n\in\omega$ there is $i\leq
\psi(\mu_Y(n+1))$ such that
\[\{L(s)(\nu): \nu\in T[s]\ \&\ \suc_{T[s]}(\nu)=\emptyset\}\cap [\mu_Y(n+i),
\mu_Y(n+i+1))=\emptyset.\]
\end{claim}

\noindent{\em Proof of the claim:}\ \ \ The proof is by induction on the
height of the tree $T[t]$. One could try just to apply the inductive
hypothesis to creatures determined by $\{\eta\in T[t]: \nu\trianglelefteq\eta
\}$ for each $\nu\in\suc_{T[t]}(\langle\rangle)$. However, this would not be
enough. What we need to do is to shrink $t$ to separate the sets
\[\{L[t](\eta): \nu\trianglelefteq\eta\ \&\ \eta\in T[t]\ \&\ \suc_{T[t]}(
\eta)=\emptyset\}\]
for distinct $\nu\in\suc_{T[t]}(\langle\rangle)$ by intervals $[\mu_Y(m),
\mu_Y(m+1))$. This will prevent ``bad events'' occurring above distinct $\nu
\in\suc_{T[t]}(\langle\rangle)$ from accumulating. The shrinking procedure
depends on the character of $\langle\rangle$ in $T[t]$, so the arguments brake
into two cases.
\medskip

\noindent{\sc Case 1:}\qquad $\langle\rangle\in D[t]$.\\
Let 
\[\begin{array}{l}
A_0=\{\nu\in\suc_{T[t]}(\langle\rangle): (\exists m\in\omega)(\mu_Y(2m)\leq
L[t](\nu)\leq R[t](\nu)<\mu_Y(2m+1))\},\\
A_1=\{\nu\in\suc_{T[t]}(\langle\rangle): (\exists m\in\omega)(\mu_Y(2m-1)\leq
L[t](\nu)\leq R[t](\nu)<\mu_Y(2m))\},\\
A_2=\{\nu\in\suc_{T[t]}(\langle\rangle): (L[t](\nu),R[t](\nu)]\cap Y\neq
\emptyset\}.
\end{array}\]
Note that $A_0\cup A_1\cup A_2=\suc_{T[t]}(\langle\rangle)$ and $\NOR[t](
\langle\rangle)(\suc_{T[t]}(\langle\rangle))\geq\nor^0[t]$. As $\NOR[t](\langle
\rangle)$ is a nice pre-norm, at least one of the following holds:
\begin{enumerate}
\item[(0)] $\NOR[t](\langle\rangle)(A_0)\geq\nor^0[t]-2$,
\item[(1)] $\NOR[t](\langle\rangle)(A_1)\geq\nor^0[t]-2$,
\item[(2)] $\NOR[t](\langle\rangle)(A_2)\geq\nor^0[t]-2$.
\end{enumerate}
Suppose that (0) holds. Let $s\in \Sigma^*_{\ref{locex}}(t)$ be a creature
such that  
\[T[s]=\{\eta\in T[t]: (\exists\nu\in A_0)(\eta\trianglelefteq\nu\ \mbox{ or }
\nu\vartriangleleft\eta)\}.\]
[It should be clear that this uniquely defines the creature $s$ and $\nor^0[s]
\geq\nor^0[t]-2$.] Note that by the definition of the set $A_0$, if
$\nu_1,\nu_2\in A_0$ are distinct and $R[t](\nu_1)\leq L[t](\nu_2)$ then there
are $m_1<m_2<\omega$ such that
\[\begin{array}{l}
\mu_Y(2m_1)\leq L[s](\nu_1)\leq R[s](\nu_1)<\mu_Y(2m_1+1)<\mu_Y(2m_2)\leq\\
L[s](\nu_2)\leq R[s](\nu_2)<\mu_Y(2m_2+1).
  \end{array}\]
Hence we may conclude that for every $n\in\omega$ there is $i<2$ such that
\[\{L[s](\eta): \eta\in T[s]\ \&\ \suc_{T[s]}(\eta)=\emptyset\}\cap
[\mu_Y(n+i),\mu_Y(n+i+1))=\emptyset,\] 
and thus $s$ is as required in the assertion of the claim.\\
If (1) holds then we may proceed in the same manner (considering the set
$A_1$).\\
So suppose that (2) holds true. Divide the set $A_2$ into three sets $A^0_2,
A^1_2,A^2_2$ such that for each $i<3$ and $\nu_1,\nu_2\in A^i_2$ with
$R[t](\nu_1)<L[t](\nu_2)$, there is $m\in\omega$ such that 
\[R[t](\nu_1)<\mu_Y(m)<\mu_Y(m+1)<L[t](\nu_2)\]
[e.g.~each $A^i_2$ contains every third element of $A_2$ counting according to
the values of $L[t](\nu)$]. For some $i<3$ we have
\[\NOR[t](\langle\rangle)(A^i_2)\geq\NOR[t](\langle\rangle)(A_2)-2\geq
\nor^0[t]-4.\] 
Fix $\nu\in A^i_2$ for a moment and let $T_\nu=\{\eta\in T[t]:\nu
\vartriangleleft\eta\ \mbox{or}\ \eta\trianglelefteq\nu\}$. If $T_\nu\cap
D[t]=\{\langle\rangle\}$ then necessarily $T_\nu$ does not contain sequences
of length $>2$ (remember clause {\bf (a)(ii)} of the definition of
$K^\psi_{\ref{locex}}$) and 
\[|\{\eta\in T_\nu:\suc_{T[t]}(\eta)=\emptyset\}|\in \{1,\psi(L[t](\nu))\}.\]
Let $n^*\in\omega$ be maximal such that $\mu_Y(n^*)\leq L[t](\nu)$ (if there
is no such $n^*$ the we let $n^*=-1$). Then, for every $n\geq n^*$,
$\psi(\mu_Y(n+1))\geq \psi(\mu_Y(n^*+1))\geq\psi(L[t](\nu))$ (as $\psi$ is
non-decreasing). Hence, letting $T^*_\nu=T_\nu$ we have that for each
$n\in\omega$ there is $i\leq \psi(\mu_Y(n+1))$ such that
\[\{L[t](\eta):\eta\in T^*_\nu\ \&\ \suc_{T[t]}(\eta)=\emptyset\}\cap
[\mu_Y(n+i),\mu_Y(n+i+1))=\emptyset.\]  
If $T_\nu\cap D[t]\neq\{\langle\rangle\}$ then we may look at a creature
$t^*\in K^\psi_{\ref{locex}}$ such that $T[t^*]=\{\eta':\nu\conc\eta'\in
T[t]\}$ and $L[t^*],R[t^*],D[t^*],\NOR[t^*]$ are copied in a suitable manner
from $t$ (so they are restrictions of the corresponding objects to $T_\nu$)
and $m^{t^*}_{\dn}=m^t_{\dn}$, $m^{t^*}_{\up}=m^t_{\up}$. Clearly
$\nor^0[t^*]\geq\nor^0[t]$ and the height of $T[t^*]$ is smaller than that of
$T[t]$. Thus we may apply the inductive hypothesis and we find $s^*\in
\Sigma^*_{\ref{locex}}(t^*)$ such that $\nor^0[s^*]\geq\nor^0[t^*]-10\geq
\nor^0[t]-10$ and for all $n\in\omega$ there is $i\leq \psi(\mu_Y(n+1))$ such
that
\[\{L[s^*](\eta):\eta\in T[s^*]\ \&\ \suc_{T[s^*]}(\eta)=\emptyset\}\cap
[\mu_Y(n+i),\mu_Y(n+i+1))=\emptyset.\] 
Let $T^*_\nu=\{\nu\conc\eta':\eta'\in T[s^*]\}$.\\
Look at the tree $T^*=\bigcup\limits_{\nu\in A^i_2} T^*_\nu$. It determines a
creature $s\in\Sigma^*_{\ref{locex}}(t)$ (i.e.~$s$ is such that $T[s]=T^*$,
$D[s]=D[t]\cap T^*$ etc). Clearly $\nor^0[s]\geq\nor^0[t]-10$ and it is easy
to check that $s$ is as required (remember the choice of the $A^i_2$ and the
$T^*_\nu$'s). 
\medskip

\noindent{\sc Case 2:}\qquad $\langle\rangle\notin D[t]$.\\
Since $T[t]\neq\{\langle\rangle\}$ (as $\nor^0[t]>0$) we have $|\suc_{T[t]}(
\langle\rangle)|=\psi(L[t](\langle\rangle))$. Moreover, for each $\nu\in
\suc_{T[t]}(\langle\rangle)$
\[\mbox{either }\ \ \ \nu\in D[t]\ \ \ \mbox{ or }\ \ \ \suc_{T[t]}(\nu)=
\emptyset\]
(remember the demand {\bf (a)(ii)} of the definition of $K^\psi_{
\ref{locex}}$). Note that necessarily $D[t]\cap\suc_{T[t]}(\langle\rangle)
\neq\emptyset$ (as $\nor^0[t]>0$). Fix $\nu\in D[t]\cap\suc_{T[t]}(\langle
\rangle)$.

\noindent Choose a set $A_\nu\subseteq\suc_{T[t]}(\nu)$ such that $\NOR[t](
\nu)(A_\nu)\geq\nor^0[t]-5$ and one of the following holds:
\begin{itemize}
\item $(\exists m\in\omega)(\forall\eta\in A_\nu)(\mu_Y(m)\leq L[t](\eta)\leq
R[t](\eta)<\mu_Y(m+1))$,
\item there are $m_0<m_1<m_2<m_3<\omega$ such that
\[\begin{array}{l}
(\forall\eta\in A_\nu)(\mu_Y(m_1)\leq L[t](\eta)\leq R[t](\eta)<\mu_Y(m_2))
\quad\mbox{ and}\\
(\exists\eta\in\suc_{T[t]}(\nu))(L[t](\eta)\leq \mu_Y(m_0))\qquad\mbox{ and}\\
(\exists\eta\in\suc_{T[t]}(\nu))(R[t](\eta)\geq \mu_Y(m_3)).
  \end{array}\]
\end{itemize}
Why is the choice of $A_\nu$ possible? For $m\in\omega$ let 
\[B^m_0=\{\eta\in\suc_{T[t]}(\nu):\mu_Y(m)\leq L[t](\eta)\leq R[t](\eta)<\mu_Y
(m+1)\}.\]
If there is $m$ such that $\NOR[t](\langle\rangle)(B^m_0)\geq \nor^0[t]-5$
then we may take the respective $B^m_0$ as $A_\nu$. So suppose that 
\[(\forall m\in\omega)(\NOR[t](\langle\rangle)(B^m_0)<\nor^0[t]-5).\]
Let $B_0=\bigcup\limits_{m\in\omega}B^m_0$ and suppose that $\NOR[t](\nu)(B_0)
\geq\nor^0[t]-1$. Let $k_0,k_1$ be the two smallest elements of $\{m: B^m_0
\neq\emptyset\}$ and let $k_2,k_3$ be the two largest elements of this set
(note that $|\{m: B^m_0\neq\emptyset\}|\geq 6$; remember that $\NOR[t](\nu)$
is a nice pre-norm). We let 
\[A_\nu=\{\eta\in\suc_{T[t]}(\nu): \mu_Y(k_1+1)\leq L[t](\eta)\leq R[t](\eta)<
\mu_Y(k_2)\},\]
and $m_0=k_0+1$, $m_1=k_1+1$, $m_2=k_2$, $m_3=k_3$. Easily
$\NOR[t](\nu)(A_\nu)\geq\NOR[t](\nu)(B_0)-4\geq\nor^0[t]-5$ and since 
$B^{k_0}_0\neq\emptyset$, $B^{k_3}_0\neq\emptyset$ we see that $A_\nu$ is as
required. So we are left with the possibility that $\NOR[t](\nu)(B_0)<\nor^0[
t]-1$. In this case we have
\[\NOR[t](\nu)(\{\eta\in\suc_{T[t]}(\nu): (L[t](\eta),R[t](\eta)]\cap Y\neq
\emptyset\})\geq\nor^0[t]-1.\]
Let $\eta_0,\eta_1,\eta_2,\eta_3\in\suc_{T[t]}(\nu)\setminus B_0$ be such that
$L[t](\eta_0)<L[t](\eta_1)$ are the first two members of $\{L[t](\eta):\eta\in
\suc_{T[t]}(\nu)\setminus B_0\}$ and $L[t](\eta_2)<L[t](\eta_3)$ are the last
two members of this set. Let $m_i$ be such that $\mu_Y(m_i)\in (L[t](\eta_i),
R[t](\eta_i)]$ (for $i<4$) and let 
\[A_\nu=\{\eta\in\suc_{T[t]}(\nu)\setminus B_0:\mu_Y(m_1)\leq L[t](\eta)\leq
R[t](\eta)<\mu_Y(m_2)\}.\]
Since 
\[\NOR[t](\nu)(\suc_{T[t]}(\nu)\setminus B_0\setminus\{\eta_0,\eta_1,\eta_2,
\eta_3\})\geq\nor^0[t]-5,\]
one easily checks that $A_\nu$ is as required.
\medskip

Let $T^\nu=\{\eta^*: \nu\conc\eta^*\in T[t]\ \&\ (\exists \eta\in A_\nu)(\eta
\vartriangleleft\nu\conc\eta^*\mbox{ or }\nu\conc\eta^*\trianglelefteq\eta)\}
$. We would like to apply the inductive hypothesis to the creature determined
by $T^\nu$ (with $L,R,D$ and $\NOR$ copied in a suitable way from $t$). 
However, this creature may have too small norm: it may happen that
$\NOR[t](\nu)(A_\nu)<11$. But we may repeat the procedure of {\sc Case 1},
noticing that the inductive hypothesis was applied there above some elements
of $\suc(\langle\rangle)$. Here, this corresponds to applying the inductive
hypothesis to creatures determined by $\{\eta^*: \eta\conc\eta^*\in T[t]\}$
for some $\eta\in A_\nu$, and these creatures have norms not smaller than
$\nor^0[t]$. Consequently we will get a tree 
\[T^\nu_*\subseteq\{\eta': \nu\vartriangleleft\eta'\ \&\ (\exists\eta\in
A_\nu)(\eta\trianglelefteq\eta')\}\]
corresponding to a creature $s_\nu$ with $\nor^0[s_\nu]\geq\nor^0[t]-10$ and
such that for each $n\in\omega$, for some $i\leq \psi(\mu_Y(n+1))$ we have
\[\{L[t](\eta'):\eta'\in T^\nu_*\ \&\ \suc_{T[t]}(\eta')=\emptyset\}\cap[
\mu_Y(n+i),\mu_Y(n+i+1))=\emptyset.\]
[Note that the procedure of {\sc Case 1} may involve further shrinking of
$A_\nu$ and dropping the norm by 4. Still, $5+4<10$ so the norm of $s_\nu$ is
above $\nor^0[t]-10$.] Next let 
\[T^*=\bigcup\{T^\nu_*:\nu\in\suc_{T[t]}(\langle\rangle)\cap D[t]\}\cup\suc_{T
[t]}(\langle\rangle)\cup\{\langle\rangle\},\]
and let $s$ be the restriction of the creature $t$ to $T^*$. Check that $s\in
\Sigma^*_{\ref{locex}}(t)$ is as required. This finishes the proof of the
claim. 
\medskip

Now we may prove the proposition. We are going to show that, in the model  
$\V^{\q^*_{{\rm s}\infty}(K^\psi_{\ref{locex}},\Sigma^\psi_{\ref{locex}})}$,
the set $\{m\in\omega: \dot{W}(m)=1\}$ witnesses that the
$S^\varphi_{\up}$--localization fails. So suppose that $p=(w,t_0,t_1,t_2\ldots)
\in\q^*_{{\rm s}\infty}(K^\psi_{\ref{locex}},\Sigma^\psi_{\ref{locex}})$
and $Y\in\iso$. Since $(K^\psi_{\ref{locex}},\Sigma^\psi_{\ref{locex}})$ is
omittory we may assume that $\psi(k)<\varphi(k)$ for $k\geq\lh(w)$ and for
each $i\in\omega$ 
\begin{enumerate}
\item[($\alpha$)] $\nor[t_i]\geq 11+i+m^{t_i}_{\dn}$, $\langle\rangle\in
D[t_i]$  and
\item[($\beta$)]  $|(R[t_i](\langle\rangle),L[t_{i+1}](\langle\rangle))\cap
Y|>2$, $|(\lh(w),L[t_0])\cap Y|>2$.
\end{enumerate}
Apply \ref{cl34} to get $s_i\in\Sigma^*_{\ref{locex}}(t_i)$ such that
$\nor^0[s_i]\geq\nor^0[t_i]-10$ and for every $n\in\omega$ there is $j\leq
\psi(\mu_Y(n+1))$ such that
\[\{L[s_i](\nu):\nu\in T[s_i]\quad\&\quad\suc_{T[s_i]}(\nu)=\emptyset\}\cap[
\mu_Y(n+j),\mu_Y(n+j+1))=\emptyset.\]
Note that by ($\alpha$) and the definition of $\Sigma^*_{\ref{locex}}$ and
$\nor^0$ we have 
\[\nor^0[t_i]=\nor[t_i],\quad\nor^0[s_i]=\nor[s_i].\]
Hence, letting $q=(w,s_0,s_1,s_2,\ldots)$ we will have $p\leq q\in\q^*_{{\rm
s}\infty}(K^\psi_{\ref{locex}},\Sigma^\psi_{\ref{locex}})$ and for every
$n>\lh(w)$  
\[q\forces(\exists j<\varphi(\mu_Y(n+1)))([\mu_Y(n+j),\mu_Y(n+
j+1))\cap \{m\in\omega: \dot{W}(m)=1\}=\emptyset);\]
remember that $q$ forces that $\{m\in\omega:\dot{W}(m)=1\}$ is a subset of 
\[\{m<m^{s_0}_{\dn}: w(m)=1\}\cup\{L[s_i](\nu): i<\omega\ \&\ \nu\in T[s_i]\
\&\ \suc_{T[s_i]}(\nu)=\emptyset\}.\]
This finishes the proof.
\end{proof}

\begin{proposition}
\label{getlocal}
Let $\varphi,\psi\in\baire$ be non-decreasing, $\varphi(0),\psi(0)>0$ and
$\lim\limits_{n\to\infty}\varphi(n)\leq\lim\limits_{n\to\infty}\psi(n)$.
Suppose that $N\prec({\mathcal H}(\chi),\in,<^*_\chi)$ is countable, $p\in N
\cap\q^*_{{\rm s}\infty}(K^\psi_{\ref{locex}},\Sigma^\psi_{\ref{locex}})$, 
$\varphi,\psi\in N$ and $Y\in\iso$ is such that 
\[(\forall X\in\iso\cap N)(\exists^\infty n\in\omega)(\forall i<\varphi(
\mu_Y(n)))(|[\mu_Y(n+i),\mu_Y(n+i+1))\cap X|\geq 2).\]
Then there is an $(N,\q^*_{{\rm s}\infty}(K^\psi_{\ref{locex}},\Sigma^\psi_{
\ref{locex}}))$--generic condition $q$ stronger than $p$ and such that 
\[\begin{array}{ll}
q\forces&\mbox{``for every }X\in\iso\cap N[\Gamma_{\q^*_{{\rm s}\infty}(
K^\psi_{\ref{locex}},\Sigma^\psi_{\ref{locex}})}],\\
\ &(\exists^\infty n\in\omega)(\forall i<\varphi(\mu_Y(n)))(|[\mu_Y(n+i),
\mu_Y(n+i+1))\cap X|\geq 2)\mbox{''.}
  \end{array}\]
Consequently, if $\psi$ is unbounded then the forcing notion $\q^*_{{\rm
s}\infty}(K^\psi_{\ref{locex}},\Sigma^\psi_{\ref{locex}})$ has the
$S^\varphi_{\dn}$--localization property for every non-decreasing $\varphi$. 
If $\psi$ is bounded, say $\lim\limits_{n\to\infty}\psi(n)=k$, then
$\q^*_{{\rm s}\infty}(K^\psi_{\ref{locex}},\Sigma^\psi_{\ref{locex}})$ has the
$S_k$-localization property. 
\end{proposition}

\begin{proof}
This is like \cite[2.4.5]{RoSh:501}. We will deal with the case
$\lim\limits_{n\to\infty}\psi(n)=\infty$ (if $\psi$ is bounded then the
arguments are similar).  

Suppose that $\varphi,\psi,N,p,Y$ are as in the assumptions. Let $\langle\dot{
\sigma}_n: n<\omega\rangle$ enumerate all $\q^*_{{\rm s}\infty}(K^\psi_{
\ref{locex}},\Sigma^\psi_{\ref{locex}})$--names from $N$ for ordinals and let
$\langle\dot{X}_n:n<\omega\rangle$ list all names for infinite subsets of
$\omega$. Further, for $n\in\omega$, let $\dot{\tau}_n\in N$ be a name for a
function in $\baire$ such that
\[\forces_{\q^*_{{\rm s}\infty}(K^\psi_{\ref{locex}},\Sigma^\psi_{
\ref{locex}})}(\forall m\in\omega)(\forall i\leq n)(|[m,\dot{\tau}_n(m))\cap
\dot{X}_i|>2).\]
We inductively construct sequences $\langle p_n,p^*_n: n<\omega\rangle$,
$\langle g_n: n<\omega\rangle$ and $\langle i_n:n<\omega\rangle$ such that for
each $n\in\omega$:
\begin{enumerate}
\item[(a)] $p_n,p^*_n\in \q^*_{{\rm s}\infty}(K^\psi_{\ref{locex}},
\Sigma^\psi_{\ref{locex}})\cap N$, $g_n\in\baire\cap N$ is strictly
increasing, $i_n\in\omega$,
\item[(b)] $p\leq^{{\rm s}\infty}_0 p_n\leq^{{\rm s}\infty}_n p^*_n\leq^{{\rm
s}\infty}_{n+1} p_{n+1}\leq^{{\rm s}\infty}_{n+1} p^*_{n+1}$,\quad $i_n<i_{n+
1}$,
\item[(c)] $p_n$ approximates $\dot{\sigma}_n$ at each $t^p_{n+\ell}$ for
$\ell\in\omega$,
\item[(d)] for each $\ell\in\omega$ and $v\in\pos(w^p,t^{p_n}_0,\ldots,
t^{p_n}_{n+\ell-1})$, 

if $t\in\Sigma^\psi_{\ref{locex}}(t^{p_n}_{n+\ell})$, $\nor[t]>0$ then there
is $w\in\pos(v,t)$ such that 

$(w,t^{p_n}_{n+\ell+1},t^{p_n}_{n+\ell+2},\ldots)\forces$ ``$\dot{\tau}_n(
g_n(\ell))<g_n(\ell+1)$'',
\item[(e)] $\varphi(g_n(\ell))<\psi(L[t^{p_n}_{n+\ell}](\langle\rangle))$ for
every $\ell\in\omega$,
\item[(f)] $\langle\rangle\in D[t^{p_n}_m]$ for all $m\in\omega$,
\item[(g)] $i_n$ is such that $(\forall i<\varphi(\mu_Y(i_n)))(|[\mu_Y(i_n+i),
\mu_Y(i_n+i+1))\cap\rng(g_n)|\geq 2)$,
\item[(h)] $t^{p^*_n}_m=t^{p_n}_m$ for $m<n$,
\item[(i)] for $i<\varphi(\mu_Y(i_n))$ let $j^n(i)=\min\{j\in\omega:g_n(j)\in
[\mu_Y(i_n+i),\mu_Y(i_n+i+1))\}$ and for $i\in [\varphi(\mu_Y(i_n)),
\psi(L[t^{p_n}_{n+j^n(0)}](\langle\rangle)))$ let $j^n(i)=j^n(\varphi(\mu_Y(
i_n))-1)+i$;\\
then $t^{p^*_n}_n=S_{m_0,m_1}(S^*(t^{p_n}_{n+j^n(i)}:i<\psi(L[t^{p_n}_{n+
j^n(0)}](\langle\rangle))))$ (see clauses (2), (4) of the definition of
$(K^\psi_{\ref{locex}},\Sigma^\psi_{\ref{locex}})$ in \ref{locex}, remember
{\bf (f)} above), where $m_0=m^{t^{p_n}_{n-1}}_{\up}$, $m_1=m^{t^{p_n}_k}_{
\up}$, $k=n+j^n(\psi(L[t^{p_n}_{n+j^n(0)}](\langle\rangle))-1)$,
\item[(j)] $t^{p^*_n}_{n+1+m}=t^{p_n}_{k+1+m}$ for every $m\in\omega$ (where
$k$ is as in {\bf (i)} above).
\end{enumerate}
The construction is quite straightforward and essentially described by the
requirements {\bf (a)}--{\bf (j)} above. Having defined $p^*_n$ we first
choose a condition $p^\prime_{n+1}\in N$ such that it approximates
$\dot{\sigma}_{n+1}$ at each $m\geq n+1$ and $p^*_n\leq^{{\rm s}\infty}_{n+1}
p^\prime_{n+1}$ (by \ref{deciding}). Then we use \ref{cl35} {\em inside\/} $N$
to find a condition $p_{n+1}\geq_{n+1} p^\prime_{n+1}$, $p_{n+1}\in N$ and a
function $g_{n+1}\in N$ such that the demands {\bf (a)}, {\bf (d)}--{\bf (f)}
are satisfied. Note here, that the creatures $t^p_\ell$ constructed in the
proof of \ref{cl35} came from the application of $\circledast^1_{\AB}$ of
\ref{ab}. This condition, in turn, is exemplified in our case by the use of
the operation $S^{**}_H$ (see item (5) of the definition of
$(K^\psi_{\ref{locex}},\Sigma^\psi_{\ref{locex}})$, see \ref{getnorm},
\ref{omitalmbound}). Consequently we may ensure that {\bf (f)} holds. As far
as {\bf (e)} is concerned, note that in the inductive construction of the
condition $p$ in \ref{cl35}, when choosing $p^*_\ell$, we may require
additionally that $\varphi(g(\ell))<\psi(L[t^{p^*_\ell}_\ell](\langle
\rangle))$ (remember $(K^\psi_{\ref{locex}},\Sigma^\psi_{\ref{locex}})$ is
omittory, $\psi$ is unbounded). Thus we have defined $p_{n+1},g_{n+1}$
satisfying {\bf (a)}--{\bf (f)}. By the assumptions on $Y$ we know
\[(\exists^\infty m\in\omega)(\forall i<\varphi(\mu_Y(m)))(|[\mu_Y(m+i),
\mu_Y(m+i+1))\cap \rng(g_{n+1})|\geq 2)\]
(remember $\rng(g_{n+1})\in\iso\cap N$). Hence we may choose $i_{n+1}>i_n$ as
required in {\bf (g)}. Clauses {\bf (h)}--{\bf (j)}, {\bf (b)} fully describe
the condition 
\[p^*_{n+1}=(w^p,t^{p_{n+1}}_0,\ldots,t^{p_{n+1}}_n,t^{p^*_{n+1}}_{n+1},
t^{p^*_{n+1}}_{n+2},\ldots)\geq_{n+1}p_{n+1}.\]
Note here that 
\[\varphi(\mu_Y(i_{n+1}))\leq \varphi(g_{n+1}(j^{n+1}(0)))<\psi(L[t^{p_{n
+1}}_{n+1+j^{n+1}(0)}](\langle\rangle))\]
(by {\bf (e)}; remember $\phi$ is non-decreasing), so there are no problems
with the definition of $t^{p^*_{n+1}}_{n+1}$ in clause {\bf (i)}. Moreover,
$t^{p^*_{n+1}}_{n+1}\in\Sigma^\psi_{\ref{locex}}(t^{p_{n+1}}_{n+1},\ldots,
t^{p_{n+1}}_k)$, where $k=n+1+j^{n+1}(\psi(L[t^{p_{n+1}}_{n+1+j^{n+1}(0)}](
\langle\rangle))-1$, and $\nor[t^{p^*_{n+1}}_{n+1}]>m^{t^{p^*_{n+1}}_{n+
1}}_{\dn}$. Clearly $p^*_{n+1}\in N$.

Now let $q=\lim\limits_n p_n=\lim\limits_n p^*_n\in \q^*_{{\rm s}\infty}
(K^\psi_{\ref{locex}},\Sigma^\psi_{\ref{locex}})$ (see \ref{propord}). We
claim that $q$ is as required in the assertion of the proposition. Clearly it
is stronger than $p$ (by {\bf (b)}) and is $(N,\q^*_{{\rm s}\infty}
(K^\psi_{\ref{locex}},\Sigma^\psi_{\ref{locex}}))$--generic (by {\bf (c)} and
\ref{smoessdec}). The proof will be finished if we show that for each
$k\in\omega$ 
\[q\forces(\exists^\infty m\in\omega)(\forall i<\varphi(\mu_Y(m)))(|[\mu_Y(m
+i),\mu_Y(m+i+1))\cap\dot{X}_k|\geq 2).\]
To this end suppose that $k,\ell\in\omega$, $q'\in\q^*_{{\rm s}\infty}
(K^\psi_{\ref{locex}},\Sigma^\psi_{\ref{locex}})$, $q\leq q'$. Passing to an
extension of $q'$ we may assume that for some $n>k$ we have $w^{q'}\in\pos(
w^q,t^q_0,\ldots,t^q_{n-1})$ and $\ell<i_n$ (remember {\bf (b)}). Let $m\geq
n$ be such that
$t^{q'}_0\in\Sigma^\psi_{\ref{locex}}(t^q_n,\ldots,t^q_m)$. Since
$\nor[t^{q'}_0]\geq m^{t^{q'}_0}_{\dn}>0$ we find $n_0\in [n,m]$ and $s\in
\Sigma^\psi_{\ref{locex}}(t^p_{n_0})$ such that $\nor[s]>0$ and 
\[\{L[s](\eta)\!:\eta\in T[s]\ \&\ \suc_{T[s]}(\eta)=\emptyset\}\subseteq
\{L[t^{q'}_0](\eta)\!:\eta\in T[t^{q'}_0]\ \&\ \suc_{T[t^{q'}_0]}(\eta)=
\emptyset\}\]
(compare the proof that $(K^\psi_{\ref{locex}},\Sigma^\psi_{\ref{locex}})$ is
condensed). Look at the creature $t^q_{n_0}$: it is $t^{p^*_{n_0}}_{n_0}$. So
look at the clause {\bf (i)}: the creature $t^{p^*_{n_0}}_{n_0}$ was obtained
from $t^{p_{n_0}}_{n_0+j^{n_0}(i)}$ (for $i<\psi(L[t^{p_{n_0}}_{n_0+j^{n_0}
(0)}](\langle\rangle))$) by applying the operation $S^*$ and
$S_{m_0,m_1}$. Since $s\in\Sigma^\psi_{\ref{locex}}(t^q_{n_0})$ we have
\[|\suc_{T[s]}(\langle\rangle)|=\psi(L[t^{p_{n_0}}_{n_0+j^{n_0}(0)}](\langle
\rangle))=\psi(L[s](\langle\rangle))\]
and for each $\nu\in\suc_{T[s]}(\langle\rangle)$ for unique $i=i(\nu)<\psi(
L[t^{p_{n_0}}_{n_0+j^{n_0}(0)}](\langle\rangle))$, the tree $T[s]$ above $\nu$
and $L[s],R[s]$ and $D[s]$ look exactly like $T[s_i],L[s_i],R[s_i]$ and
$D[s_i]$ for some $s_i\in\Sigma^*_{\ref{locex}}(t^{p_{n_0}}_{n_0+j^{n_0}(i)})$
with $\nor[s_i]>0$. Now, applying successively clause {\bf (d)} to each
$t^{p_{n_0}}_{n_0+j^{n_0}(i)}$, $s_i$ (for $i<\psi(L[s](\langle\rangle))$) we
find $w\in\pos(w^{q'},t^{q'}_0)$ such that for every $\nu\in\suc_{T[s]}(
\langle\rangle)$
\[(w,t^{p_{n_0}}_j,t^{p_{n_0}}_{j+1},\ldots)\forces_{\q^*_{{\rm s}\infty}(
K^\psi_{\ref{locex}},\Sigma^\psi_{\ref{locex}})}\mbox{`` }\dot{\tau}_{n_0}(
g_{n_0} (j^{n_0}(i(\nu))))<g_{n_0}(j^{n_0}(i(\nu))+1)\mbox{ '',}\]
where $j$ is such that $m^{t^{p_{n_0}}_j}_{\dn}=m^{t^{q'}_0}_{\up}$. By the
definition of the name $\dot{\tau}_{n_0}$ and the fact that $k<n\leq n_0$ we
get that for each $\nu\in\suc_{T[s]}(\langle\rangle)$
\[(w,t^{p_{n_0}}_j,t^{p_{n_0}}_{j+1},\ldots)\forces\mbox{`` }|[g_{n_0}(
j^{n_0}(i(\nu))),g_{n_0}(j^{n_0}(i(\nu))+1)\cap\dot{X}_k|\geq 2\mbox{ ''.}\]
By the choice of $i_{n_0}$, $j^{n_0}(i)$ we know that 
\[\mu_Y(i_{n_0}+i)\leq g_{n_0}(j^{n_0}(i))<g_{n_0}(j^{n_0}(i)+1)<\mu_Y(i_{n_0}
+i+1)\]
for each $i<\varphi(\mu_Y(i_{n_0}))\leq\varphi(g_{n_0}(j^{n_0}(0)))\leq\psi(
L[s])(\langle\rangle)$. Therefore 
\[(w,t^{p_{n_0}}_j,t^{p_{n_0}}_{j+1},\ldots)\forces\mbox{``}(\forall
i\!<\!\varphi(\mu_Y(i_{n_0})))(|[\mu_Y(i_{n_0}+i),\mu_Y(i_{n_0}+i+1))\cap
\dot{X}_k|\geq 2)\mbox{''.}\]
We finish noticing that $(w,t^{p_{n_0}}_j,t^{p_{n_0}}_{j+1},\ldots)\leq
(w,t^{q'}_1,t^{q'}_2,\ldots)$ and $\ell<i_n\leq i_{n_0}$.
\end{proof}

\begin{corollary}
The following is consistent with ZFC:
\begin{quotation}
\noindent $\dominating=\cov(\M)=\non(\M)=\aleph_2$ +

\noindent for every non-decreasing unbounded $\psi\in\baire$, $\dominating(
S^\psi_{\up})=\aleph_2$ +

\noindent for every non-decreasing $\varphi\in\baire$, $\dominating(
S^\varphi_{\dn})=\aleph_1$.
\end{quotation}
\end{corollary}

\begin{proof}
Start with a model for CH and build a countable support iteration $\langle
\p_\alpha,\dot{\q}_\alpha:\alpha<\omega_2\rangle$ and a sequence
$\langle\dot{\psi}_\alpha:\alpha<\omega_2\rangle$ such that  
\begin{enumerate}
\item[(a)] $\dot{\psi}_\alpha$ is a $\p_\alpha$--name for a non-decreasing
unbounded function in $\baire$,
\item[(b)] $\forces_{\p_\alpha}$ `` $\dot{\q}_\alpha=\q^*_{{\rm s}\infty}(
K^{\dot{\psi}_\alpha}_{\ref{locex}},\Sigma^{\dot{\psi}_\alpha}_{\ref{locex}
})$ '',
\item[(c)] for each $\p_{\omega_2}$--name $\dot{\psi}$ for a non-decreasing
unbounded function in $\baire$, for $\omega_2$ many $\alpha<\omega_2$,
$\forces_{\p_\alpha}$ `` $\dot{\psi}=\dot{\psi}_\alpha$ ''.
\end{enumerate}
By \ref{whatitdoes} we have
\[\forces_{\p_{\omega_2}}\mbox{`` }\cov(\M)=\non(\M)=\aleph_2=\con\ \ \&\ \
\unbounded = \aleph_1\mbox{ ''}\]
and by \ref{noloc} we get 
\[\forces_{\p_{\omega_2}}\mbox{`` }\dominating=\aleph_2 + \mbox{for each
non-decreasing unbounded }\psi\in\baire,\ \dominating(S^\psi_{\up})=\aleph_2
\mbox{ ''.}\]
To show that 
\[\forces_{\p_{\omega_2}}\mbox{`` }\dominating(S^\varphi_{\dn})=\aleph_1\mbox{
for every non-decreasing }\varphi\in\baire\mbox{ ''}\]
we use \ref{getlocal} and \cite[Ch XVIII, 3.6]{Sh:f} and we show that the
property described in \ref{getlocal} is preserved in countable support
iterations. 

So suppose that $\varphi\in\baire$ is non-decreasing and define a context
$(R^\varphi,S^\varphi,{\bf g}^\varphi)$ (see \cite[Ch XVIII, 3.1]{Sh:f}) as
follows. First, for $\eta\in\baire$ let $X_\eta\in\iso$ be such that
$\mu_{X_\eta}(n)=\sum\limits_{i\leq n}(\eta(i)+1)$. Next we let (in the ground
model $\V$): 
\begin{itemize}
\item $S^\varphi$ is the collection of all $N\cap{\mathcal H}(\aleph_1)$ for
$N$ a countable elementary submodel of $({\mathcal H}(\chi),\in,<^*_\chi)$,
\item for each $a\in S^\varphi$,\quad $d[a]=c[a]=\omega=d'[a]=c'[a]$,
\item $\alpha^*=1$,
\item $R^\varphi=R^\varphi_0$ is the relation determined by the
$S^\varphi_{\dn}$--localization: 

$\eta\; R^\varphi\; g$ if and only if  ($\eta,g\in\baire$ and)

$(\exists^\infty n\in\omega)(\forall i<\varphi(\mu_{X_g}(n)))(|[\mu_{X_g}(
n+i),\mu_{X_g}(n+i+1))\cap X_\eta|\geq 2)$

(remember that $\alpha^*=1$ so we have $R^\varphi_0$ only),
\item ${\bf g}^\varphi=\langle {\bf g}_a: a\in S^\varphi\rangle\subseteq
\baire$ is such that for every $a\in S^\varphi$, for each $\eta\in a\cap
\baire$ we have $\eta\; R^\varphi\; {\bf g}_a$ (exists as each $a$ is
countable, e.g.~one may take as ${\bf g}_a$ any real dominating $a$).
\end{itemize}
By the choice of $(R^\varphi,S^\varphi,{\bf g}^\varphi)$ we know that it
covers in $\V$ (see \cite[Ch XVIII, 3.2]{Sh:f}). 

\begin{claim}
\label{cl36}
Let $\p$ be a proper forcing notion such that 
\[\forces_{\p}\mbox{`` }(R^\varphi,S^\varphi,{\bf g}^\varphi)\mbox{
covers''.}\]
Then:
\begin{enumerate}
\item $\forces_{\p}$`` $(R^\varphi,S^\varphi,{\bf g}^\varphi)$ strongly covers
by Possibility B'' (see \cite[Ch XVIII, 3.3]{Sh:f}),
\item for every $\p$-name $\dot{\psi}$ for a non-decreasing unbounded function
in $\baire$,
\[\forces_{\p}\mbox{`` the forcing notion $\q^*_{{\rm s}\infty}(
K^{\dot{\psi}}_{\ref{locex}},\Sigma^{\dot{\psi}}_{\ref{locex}})$ is
$(R^\varphi,S^\varphi,{\bf g}^\varphi)$--preserving''.}\]
\end{enumerate}
\end{claim}

\noindent{\em Proof of the claim:}\ \ \ 1)\ \ \ As $\alpha^*=1$ it is enough
to show that the second player has an absolute (for extensions by proper
forcing notions) winning strategy in the following game $G_a$ (for each $a\in
S^\varphi$).
\begin{quotation}
\noindent The play lasts $\omega$ moves.\\
Player I, in his $n^{\rm th}$ move chooses functions $f^n_0,\ldots,f^n_n
\in\baire$ such that 
\[f^n_\ell\rest b_{n-1}=f^{n-1}_\ell\rest b_{n-1}\ \mbox{ for $\ell<n$ and }\
\ f^n_\ell\;R^\varphi\;{\bf g}_a\mbox{ for each $\ell\leq n$.}\]
Player II answers choosing a finite set $b_n\subseteq\omega$, $b_{n-1}
\subseteq b_n$.

\noindent At the end the second player wins if and only if for each $\ell\in
\omega$
\[\bigcup_{n\geq\ell} f^n_\ell\rest b_n\; R^\varphi\; {\bf g}_a.\]
\end{quotation}
But it should be clear that Player II has a (nice) winning strategy in this
game. In his $n^{\rm th}$ move he chooses as $b_n$ a sufficiently long initial
segment of $\omega$ to provide new ``witnesses'' for the quantifier
$(\exists^\infty n\in\omega)$ in the definition of $R^\varphi$ (for all
$f^n_0,\ldots,f^n_n$). 

\noindent 2)\ \ \ Since $\alpha^*=1$ what we have to prove is exactly the
statement of \ref{getlocal} (see \cite[Ch XVIII, 3.4A]{Sh:f}), so we are done
with the claim. 
\medskip

Now we may use \cite[Ch XVIII, 3.5, 3.6]{Sh:f} (and \ref{cl36}) to conclude
that $\forces_{\p_{\omega_2}}$ ``$(R^\varphi,S^\varphi,{\bf g}^\varphi)$
covers'' and hence immediately 
\[\forces_{\p_{\omega_2}}\dominating(S^\varphi_{\dn})=\aleph_1.\]
As every function from $\baire\cap \V^{\p_{\omega_2}}$ appears in $\baire\cap
\V^{\p_\alpha}$ for some $\alpha<\omega_2$ we finish the proof.
\end{proof}

\chapter{Around not adding Cohen reals}
The starting point for this chapter was the following request of
Bartoszy\'nski (see \cite[Problem 4]{Ba94}): construct a proper forcing notion
$\p$ such that: 
\begin{enumerate}
\item $\p$ is $\baire$--bounding,
\item $\p$ preserves non--meager sets,
\item $\p$ makes ground reals to have measure zero,
\item $\p$ has the Laver property,
\item countable support iterations of $\p$ with Laver forcing, random real
forcing and Miller's rational perfect set forcing do not add Cohen reals.
\end{enumerate}
A forcing notion with these properties would correspond to the invariant ${\bf
non}({\mathcal N})$ (the minimal size of a non-null set; see
\cite[7.3C]{BaJu95}). Forcing notions with properties (1)--(4) were known. The
fourth property is a kind of technical assumption and might be replaced by
\begin{enumerate}
\item[$4^-$.] $\p$ is $(f,g)$--bounding for some $f,g\in\baire$ (with $g(n)\ll
f(n)$, of course).
\end{enumerate}
At least we believe that that was the intension (see an example presented in
\cite[7.3C]{BaJu95}, see \ref{tomek2} here too). However it was not clear 
how one should take care of the last required property. The problem comes from
the fact that we do not have any good (meaning: iterable and sufficiently
weak) condition for ``not adding Cohen reals''. The difficulty starts already
at the level of compositions of forcing notions: adding first a dominating
real and then ``infinitely often equal real below it'' one produces a Cohen
real. Various iterable properties implying ``no Cohen reals'' are in use, but
the point is to find one capturing as many of them as possible. The first
section deals with $(f,g)$--bounding property. We generalize this property
in the following section (a special case of the methods developed there is
presented in \cite[7.2E]{BaJu95}, however not fully). The ``$(\bar{t},
\bar{{\mathcal F}})$--bounding'' property seems to be still not weak 
enough to capture the measure algebra. So we weaken this further and we
present a good candidate for a property ``responsible'' for not adding Cohen
reals in the third part of this chapter (see \ref{superconc} too). The tools
developed in this section are very general and will be used later too.

\section{$(f,g)$-bounding}\label{sekfg}
Let us recall that a proper forcing notion $\p$ is $(f,g)$--bounding (for some
increasing $f,g\in\baire$) if
\[\forces_{\p}(\forall x\in\prod_{i\in\omega} f(i))(\exists S\in\V\cap
\prod_{i\in\omega} \fsuo)(\forall i\in\omega)(|S(i)|\leq g(i)\ \ \&\ \
x(i)\in S(i)).\]
It is almost obvious that $(f,g)$--bounding forcing notions add neither Cohen
reals nor random reals (see e.g.~\cite[7.2.15]{BaJu95}). For the treatment of
this property in countable support iterations see \cite[A2.5]{Sh:326} or
\cite[Ch VI, 2.11A-C]{Sh:f}.

\begin{definition}
\label{fbig}
Let $f\in\baire$. We say that a weak creating pair $(K,\Sigma)$ for $\bH$ is
{\em essentially $f$-big} if 
\begin{enumerate}
\item[$(\otimes_{\ref{fbig}}^f)$] for every weak creature $t\in K$ and $u\in
\basis(t)$ such that $0<\lh(u)$ and $\nor[t]>f(0)$ and each function
$h:\pos(u,t)\longrightarrow f(\lh(u))$ {\em there is} $s\in\Sigma(t)$ such
that $u\in\basis(s)$, $h\rest\pos(u,s)$ is constant and $\nor[s]\geq \nor[t]-
\frac{f(0)}{\lh(u)}$.
\end{enumerate}
\end{definition}

\begin{remark}
Definition \ref{fbig} may be thought of as a kind of strengthening of
\ref{bigetc} and \ref{kbig}. 
\end{remark}

\begin{lemma}
\label{fgfinitary}
Let $f\in\baire$ be increasing. Suppose that $(K,\Sigma)$ is a finitary
essentially $f$-big tree-creating pair, $p\in\q^{\tree}_1(K,\Sigma)$ and
$\dot{\tau}$ is a $\q^{\tree}_1(K,\Sigma)$-name such that  
\[p\forces_{\q^{\tree}_1(K,\Sigma)}\mbox{``}\dot{\tau}\in\baire\mbox{ is such
that }(\forall n\in\omega)(\dot{\tau}(n)<f(n))\mbox{''}.\] 
Then there is a condition $q\in\q^{\tree}_1(K,\Sigma)$ stronger than $p$ and
such that for every $\rho\in T^q$ the condition $q^{[\rho]}$ forces a value to
$\dot{\tau}\rest(\lh(\rho)+1)$.
\end{lemma}

\begin{proof}
First note that the essential $f$--bigness of $(K,\Sigma)$ implies that if
$t\in K$, $\nor[t]>f(0)$ and $\lh(\mrot(t))>f(0)+1$ then the tree creature $t$
is 2-big. This is more than enough to carry out the proofs of
\ref{bigfront}(2) and \ref{treedec}(2) and thus we find a condition $q_0\geq 
p$, $\lh(\mrot(q_0))>f(0)+1$, and fronts $F_0,F_1,F_2,\ldots$ of $T^{q_0}$
such that   
\[(\forall n\in\omega)(\forall \eta\in F_n)(\lh(\eta)>n\ \mbox{ and }\
q_0^{[\eta]}\mbox{ decides }\dot{\tau}(n))\] 
and $(\forall \eta\in T^{q_0})(\nor[t^{q_0}_\eta]>2\cdot f(0)+1)$. For each
$n\in\omega$ we have a function $h_n:F_n\longrightarrow f(n)$ such that
\[(\forall\eta\in F_n)(q_0^{[\eta]}\forces_{\q^{\tree}_1(K,\Sigma)}
\dot{\tau}(n)=h_n(\eta)).\]
Suppose that $\nu\in T^{q_0}$ is such that $\pos(t^{q_0}_\nu)\subseteq F_n$
and $\lh(\nu)\geq n$ (note that there are $\nu\in T^{q_0}$ such that
$\pos(t^{q_0}_\nu)\subseteq F_n$ as $F_n$ is finite). Then we may apply
$(\otimes^f_{\ref{fbig}})$ and we find $s\in\Sigma(t^{q_0}_\nu)$ such that
$h_n\rest\pos(s)$ is constant and $\nor[s]\geq\nor[t^{q_0}_\nu]-\frac{f(0)}
{\lh(\nu)}$. Repeating this process downward and for all $n\in\omega$ we find
a quasi tree $T^*\subseteq T^{q_0}$ and $s_\nu\in\Sigma(t^{q_0}_\nu)$ for
$\nu\in T^*$ such that 
\begin{enumerate}
\item[$(\alpha)$] $\mrot(T^*)=\mrot(q_0)$,
\item[$(\beta)$]  $\nor[s_\nu]\geq\nor[t^{q_0}_\nu]-(\lh(\nu)+1)\cdot
\frac{f(0)}{\lh(\nu)}\geq\nor[t^{q_0}_\nu]-(f(0)+1)$,  
\item[$(\gamma)$] $\pos(s_\nu)=\suc_{T^*}(\nu)$,
\item[$(\delta)$] if $\nu,\eta_0,\eta_1\in T^*$, $n\leq\lh(\nu)$, $\nu
\vartriangleleft\eta_0$, $\nu\vartriangleleft\eta_1$, $\eta_0,\eta_1\in F_n$
then $h_n(\eta_0)=h_n(\eta_1)$. 
\end{enumerate}
This defines a condition $q^*\in\q^{\tree}_1(K,\Sigma)$. Clearly it is
stronger than $q_0$ and, by $(\delta)$ above, it has the required property.
\end{proof}

\begin{remark}
If $(K,\Sigma)$ is a local tree-creating pair (see \ref{local}),
$p\in\q^{\tree}_e(K,\Sigma)$ then 
\[(\forall\nu\in\dcl(T^p))(\mrot(p)\trianglelefteq\nu\ \ \Rightarrow\ \ \nu\in
T^p).\]   
\end{remark}

\begin{conclusion}
\label{confgbound}
Suppose that $f\in\baire$ is increasing, $(K,\Sigma)$ is essentially $f$-big
finitary and local tree creating pair. Then the condition
$q\in\q^{\tree}_1(K,\Sigma)$ provided by the assertion of \ref{fgfinitary},
gives at most $|T^q\cap\prod\limits_{m<n}\bH(m)|$ possible values to
$\dot{\tau} \rest (n+1)$ (for each $n$). Hence, if $g(n)=\prod\limits_{m<n}
|\bH(m)|$ then $\q^{\tree}_1(K,\Sigma)$ is $(f,g)$-bounding.
\end{conclusion}

\begin{definition}
\label{reducible}
A weak creating pair is {\em reducible} if for each $t\in K$ with $\nor[t]>3$
there is $s\in\Sigma(t)$ such that
$\frac{\nor[t]}{2}\leq\nor[s]\leq\nor[t]-1$. 
\end{definition}

\begin{definition}
\label{hlimited}
Let $h:\omega\times\omega\longrightarrow\omega$. We say that a weak creating
pair $(K,\Sigma)$ is {\em $h$--limited} whenever
\begin{center}
if $t\in K$, $u\in\basis(t)$, $\lh(u)\leq m_0$ and $\nor[t]\leq m_1$ then
$|\pos(u,t)|\leq h(m_0,m_1)$. 
\end{center}
If the function $h$ does not depend on the first coordinate (i.e.~$h(m_0,m_1)
=h_0(m_1)$) then we say that $(K,\Sigma)$ is {\em $h$--norm-limited}. We may
say then that $(K,\Sigma)$ is $h_0$--norm-limited or just $h_0$--limited.
\end{definition}

\begin{theorem}
\label{fgbounding}
Suppose that $f,g\in\baire$ are increasing, $h:\omega\times\omega
\longrightarrow\omega$ and 
\[(\forall^\infty n\in\omega)(\prod_{m<n} h(m,m)<g(n)<f(n)).\]
Assume that $(K,\Sigma)$ is a reducible finitary tree creating pair which is
$h$-limited and essentially $f$-big. Then the forcing notion
$\q^{\tree}_1(K,\Sigma)$ is $(f,g)$-bounding. 
\end{theorem}

\begin{proof}
Let $N$ be such that $(\forall n\geq N)(\prod\limits_{m<n}h(m,m)<g(n)<f(n))$. 
Suppose that $p\forces_{\q^{\tree}_1(K,\Sigma)}\dot{\tau}\in\prod\limits_{n<
\omega}f(n)$. By \ref{fgfinitary} we find $q\geq p$ such that for every
$\eta\in T^q$ the condition $q^{[\eta]}$ decides $\dot{\tau}\rest(\lh(\eta)
+1)$. As $(K,\Sigma)$ is reducible we may assume that $(\forall\eta\in
T^q)(\nor[t^q_\eta]\leq\lh(\eta))$ and $\lh(\mrot(q))>N$. For $n\in\omega$ let
\[F^*_n\stackrel{\rm def}{=}\{\eta\in T^q:\lh(\eta)\geq n\ \mbox{ and
}\ (\forall\nu\in T^q)(\nu\vartriangleleft \eta\ \ \ \Rightarrow\ \ \
\lh(\nu)<n)\}.\] 
Clearly each $F^*_n$ is a front of $T^q$ and if $\eta\in F^*_n$ then
$q^{[\eta]}$ decides the value of $\dot{\tau}(n)$. Now note that $|F^*_n|=1$
for $n\leq\lh(\mrot(q))$ and $|F^*_n|\leq\prod\limits_{m<n}h(m,m)<g(n)$ for
all other $n$. This allows us to finish the proof. 
\end{proof}

\begin{theorem}
\label{fgsin}
Assume that $(K,\Sigma)$ is a finitary and reducible creating pair which is
$h$--limited for some function $h$. Further suppose that $(K,\Sigma)$ is {\em
either} growing and big {\em or} omittory and omittory--big. Then the forcing
notion $\q^*_{{\rm s}\infty}(K,\Sigma)$ is $(f,g)$--bounding for any strictly
increasing functions $f,g\in\baire$.  
\end{theorem}

\begin{proof}
Suppose that $\dot{\tau}$ is a $\q^*_{{\rm s}\infty}(K,\Sigma)$--name for a
function in $\prod\limits_{n\in\omega}f(n)$ and $p\in\q^*_{{\rm s}\infty}(K,
\Sigma)$. Applying repeatedly \ref{decbel} (or \ref{omitdecbel} in the second
case) we may construct inductively an increasing sequence $n_0<n_1<\ldots<
\omega$ and a condition $q=(w^q,t^q_0,t^q_1,\ldots)\in\q^*_{{\rm s}\infty}(K,
\Sigma)$ such that $p\leq^{{\rm s}\infty}_0 q$ (so $w^q=w^p$) and for all
$k\in\omega$:  
\begin{enumerate}
\item[($\oplus_0$)] $g(n_k)>\prod\limits_{i\leq k} h(m^{t^q_i}_{\dn}, 2\cdot
m^{t^q_i}_{\dn}+1)$, 
\item[($\oplus_1$)] $\nor[t^q_k]\leq 2\cdot m^{t^q_k}_{\dn}+1$, and
\item[($\oplus_2$)] if $w\in\pos(w^q,t^q_0,\ldots,t^q_{k-1})$ then the
condition $(w,t^q_k,t^q_{k+1},\ldots)$ decides the value of $\dot{\tau}\rest
(n_k+1)$. 
\end{enumerate}
This is straightforward; to get $(\oplus_1)$ we use the assumption that
$(K,\Sigma)$ is reducible. Now we note that for each $k$:
\[|\pos(w^q,t^q_0,\ldots,t^q_k)|\leq\prod\limits_{i\leq k} h(m^{t^q_i}_{\dn}, 
2\cdot m^{t^q_i}_{\dn}+1)\]
and so the condition $q$ allows less than $g(n_k)$ candidates for values of
$\dot{\tau}$ on the interval $(n_k, n_{k+1}]$.
\end{proof}

\begin{theorem}
\label{winfg}
Let $(K,\Sigma)$ be a reducible and finitary creating pair. Suppose that
increasing functions $f,g\in\baire$ and a function $h:\omega\times\omega
\longrightarrow\omega$ are such that 
\begin{enumerate}
\item $(K,\Sigma)$ is $h$--limited,
\item $(K,\Sigma)$ is essentially $f$--big,
\item $(\forall^\infty n)(\prod\limits_{m<n} h(m,m)<g(n)<f(n))$.
\end{enumerate}
Lastly assume that $(K,\Sigma)$ captures singletons. Then the forcing notion
$\q^*_{{\rm w}\infty}(K,\Sigma)$ is $(f,g)$--bounding.
\end{theorem}

\begin{proof}
Take $N\in\omega$ such that $\prod\limits_{m<n}h(m,m)<g(n)<f(n)$ for all
$n\geq N$. Let $\dot{\tau}$ be a $\q^*_{{\rm w}\infty}(K,\Sigma)$--name for a
function in $\prod\limits_{n\in\omega}f(n)$, and let $p\in\q^*_{{\rm w}
\infty}(K,\Sigma)$. First note that, as $(K,\Sigma)$ captures singletons,
for each $t\in K$ we may find $s\in\Sigma(t)$ such that for some
$u\in\basis(s)$ (equivalently: for each $u$, remember $(K,\Sigma)$ is
forgetful) we have $|\pos(u,s)|=1$. Using this remark and \ref{sinwin} we find
a condition $(w^p,s_0,s_1,s_2,\ldots)\in\q^*_{{\rm w}\infty}(K,\Sigma)$ and a
sequence $0\leq\ell_0<\ell_1<\ell_2<\ldots<\omega$ such that:
\begin{enumerate}
\item[($\alpha$)] $\nor[s_{\ell_0}]\geq 2\cdot f(0)+2$, $\nor[s_{\ell_{i+1}}]
\geq 2\cdot f(0)\cdot |\pos(w^p,s_0,\ldots,s_{\ell_i})|+2(i+1)$,
\item[($\beta$)]  if $n\in\omega\setminus\{\ell_0,\ell_1,\ell_2,\ldots\}$ then
for some $u\in\basis(s_n)$ we have $|\pos(u,s_n)|=1$,
\item[($\gamma$)] $N+4<m^{s_{\ell_0}}_{\dn}$,\qquad $p\leq(w^p,s_0,s_1,s_2,
\ldots)$,
\item[($\delta$)] for each $i\in\omega$, $u\in\pos(w^p,s_0,\ldots,s_{\ell_i})$
the condition $(u,s_{\ell_i+1},s_{\ell_i+2},\ldots)$ decides $\dot{\tau}\rest
(m^{s_{\ell_i}}_{\up}+1)$.
\end{enumerate}
Next we slightly correct creatures $s_{\ell_i}$ to ensure that the value of
$\dot{\tau}\rest (m^{s_{\ell_i}}_{\dn}+1)$ is decided by any $u\in\pos(w^p,
s_0,\ldots,s_{\ell_i-1})$. For this we use the procedure similar to that in
the proof of \ref{fbig} (and based on the assumption that $(K,\Sigma)$ is
essentially $f$--big). Thus we get creatures $t_{\ell_i}\in\Sigma(s_{\ell_i})$
such that (for $i\in\omega$):
\[(\forall u\in\pos(w^p,s_0,\ldots,s_{\ell_i-1}))\big((u,t_{\ell_i},s_{\ell_i+
1},s_{\ell_i+2},\ldots)\ \mbox{ decides }\ \dot{\tau}\rest
(m^{s_{\ell_i}}_{\dn}+1)\big)\quad\mbox{and}\]
\[\nor[t_{\ell_i}]\geq\nor[s_{\ell_i}]-(m^{s_{\ell_i}}_{\dn}+1)\cdot\frac{
f(0)}{m^{s_{\ell_i}}_{\dn}}\cdot |\pos(w^p,s_0,\ldots,s_{\ell_i-1})|.\]
Note that $|\pos(w^p,s_0,\ldots,s_{\ell_0-1})|=1$ and hence $\nor[t_{\ell_0}]
\geq 1$. Moreover, for each $i\in\omega$, $|\pos(w^p,s_0,\ldots,s_{\ell_i})|=
|\pos(w^p,s_0,\ldots,s_{\ell_i},\ldots,s_{\ell_{i+1}-1})|$ and therefore
$\nor[t_{\ell_{i+1}}]\geq 2(i+1)$. Finally, as $(K,\Sigma)$ is reducible, we
may choose $t^*_{\ell_i}\in\Sigma(t_{\ell_i})$ (for $i\in\omega$) such that
$\frac{i+2}{2}\leq \nor[t^*_{\ell_i}]\leq m^{t^*_{\ell_i}}_{\dn}$. Now we let 
\[w^q=w^p,\qquad t^p_m=\left\{\begin{array}{ll}
s_m          &\mbox{if }m\in\omega\setminus\{\ell_0,\ell_1,\ldots\},\\
t^*_{\ell_i} &\mbox{if }m=\ell_i,\ i\in\omega.
\end{array}\right.\]
This defines a condition $q\in\q^*_{{\rm w}\infty}(K,\Sigma)$ stronger than
$p$ and such that for each $i\in\omega$:
\begin{enumerate}
\item[(a)] $(\forall u\in\pos(w^q,t^q_0,\ldots,t^q_{\ell_i-1}))\big( (u,
t^q_{\ell_i},t^q_{\ell_i+1},\ldots)\ \mbox{ decides }\ \dot{\tau}\rest
(m^{t^q_{\ell_i}}_{\dn}+1)\big)$,
\item[(b)] $|\pos(w^q,t^q_0,\ldots,t^q_{\ell_{i+1}-1})|\leq\prod\limits_{j
\leq i}h(m^{t^q_{\ell_j}}_{\dn},m^{t^q_{\ell_j}}_{\dn})<g(m^{t^q_{
\ell_i}}_{\dn}+1)$,
\item[(c)] $|\pos(w^q,t^q_0,\ldots,t^q_{\ell_0-1})|=1$.
\end{enumerate}
Now we easily finish. 
\end{proof}

\begin{definition}
\label{FHfast}
Let $\bH$ be finitary, $F\in\baire$ be increasing. A function $f:\omega\times
\omega\longrightarrow\omega$ is called {\em $(\bH,F)$--fast} if it is
$\bH$--fast (see \ref{fast}) and additionally
\[(\forall n,\ell\in\omega)(f(n+1,\ell)>f(n,\ell)+F(\ell)\cdot\fH(\ell)\cdot
\ell).\]
\end{definition}

\begin{theorem}
\label{FHfg}
Suppose that $(K,\Sigma)$ is a finitary local and $\bar{2}$--big creating pair
for $\bH$ which has the (weak) Halving Property. Let $F\in\baire$ be
increasing and $f:\omega\times\omega\longrightarrow\omega$ be $(\bH,F)$--fast.
Then the forcing notion $\q^*_f(K,\Sigma)$ is $(F,\fH)$--bounding.
\end{theorem}

\begin{proof}
Suppose that $\dot{\tau}$ is a $\q^*_f(K,\Sigma)$--name for an element of
$\prod\limits_{m<\omega}F(m)$ and $p\in \q^*_f(K,\Sigma)$. By
\ref{halbigdec} we find a condition $q\geq p$ which essentially decides all
the values $\dot{\tau}(m)$ (for $m\in\omega$). We may assume that $(\forall
i\in\omega)(\nor[t^q_i]>f(2,m^{t^q_i}_{\dn}))$. Applying, in a standard by now
way, the bigness (like in \ref{decbel} or \ref{fbig}) we build a condition
$r\geq_0 q$ such that $t^r_i\in\Sigma(t^q_i)$, $\nor[t^r_i]\geq \nor[t^q_i]-
F(m^{t^r_i}_{\dn})\cdot\fH(m^{t^r_i}_{\dn})\cdot m^{t^r_i}_{\dn}$ and for each
$i\in\omega$ and $u\in\pos(w^p,t^r_0,\ldots,t^r_{i-1})$ the condition
$(u,t^r_i,t^r_{i+1},\ldots)$ decides the value of $\dot{\tau}\rest
(\lh(u)+1)$. Now we easily finish (remembering that $(K,\Sigma)$ is local).
\end{proof}

\section{$(\bar{t},\bar{{\mathcal F}})$--bounding}
Here we introduce and deal with a property which, in our context, is a natural
generalization of the notion of $(f,g)$--bounding forcing notions. This is a
first step toward handling ``not adding Cohen reals'' and, in some sense, it
will be developed in the next parts of this chapter. After we formulate and
prove some basic results we show how one may treat this property in countable
support iterations. 

A particular case of this machinery was presented in \cite[7.2E]{BaJu95}. 

\begin{definition}
\label{tFboundetc}
Let $(K,\Sigma)$ be a creating pair, $\bar{t}=\langle t_n:n\in\omega\rangle
\in\PC(K,\Sigma)$.
\begin{enumerate}
\item For a function $h\in\baire$ we define $U_h(\bar{t})$ as the set
\[\{\bar{s}=\langle s_n\!: n\in\omega\rangle\!\in\!\PC(K,\Sigma)\!:\bar{t}
\leq\bar{s}\ \&\ (\forall n\!\in\!\omega)(\nor[s_n]\leq h(m^{s_n}_{\dn}))\}.\]
For $n\in\omega$ and $h\in\baire$ we let
\[V^n_h(\bar{t})=\{s\in\Sigma(t_n): \nor[s]\leq h(m^s_{\dn})\}.\]
\item Let $h_1,h_2\in\baire$. We say that a forcing notion $\p$ is $(\bar{t},
h_1,h_2)$--bounding if 
\[\forces_{\p}(\forall\bar{s}\in U_{h_1}(\bar{t}))(\exists\bar{s}^*\in
U_{h_2}(\bar{t})\cap\V)(\bar{s}^*\leq\bar{s}).\]  
\end{enumerate}
\end{definition}

\begin{remark}
\label{reminf}
\begin{enumerate}
\item We will be interested in the notions introduced in \ref{tFboundetc} only
for $\bar{t}\in\PC_\infty(K,\Sigma)$ (i.e.~$\lim\limits_{n\to\infty}\nor[t_n]= 
\infty$) and $(\forall^\infty n)(h_1(n)<h_2(n))$, $\lim\limits_{n\to\infty}
h_1(n)=\lim\limits_{n\to\infty}h_2(n)=\infty$.
\item Note that if $(K,\Sigma)$ is nice and simple (see \ref{simpglui}) then 
\[U_h(\bar{t})=\prod_{n\in\omega} V^n_h(\bar{t}).\]
\end{enumerate}
\end{remark}

\begin{definition}
\label{monotonic}
For a creating pair $(K,\Sigma)$ on $\bH$ we say that:
\begin{enumerate}
\item $(K,\Sigma)$ is {\em monotonic} if for each $t\in K$, $s\in\Sigma(t)$
we have $\val[s]\subseteq\val[t]$.
\item $(K,\Sigma)$ is {\em strictly monotonic} if it is monotonic and for all
$n\in\omega$, $t_0,\ldots,t_n\in K$ and $s\in \Sigma(t_0,\ldots,t_n)$ such
that $\nor[s]\leq\max\{\nor[t_\ell]-1:\ell\leq n\}$ we have: 
\[(\forall u\in\basis(t_0))\big(\pos(u,s)\varsubsetneq \pos(u,t_0,\ldots,t_n)
\big).\] 
\item $(K,\Sigma)$ is {\em spread} if for each $t\in K$, $u\in\basis(t)$ and
$v\in\pos(u,t)$ there is $s\in\Sigma(t)$ such that 
\[\nor[s]\leq\frac{1}{2}\nor[t]\quad\mbox{ and }\quad v\in\pos(u,s).\]
\end{enumerate}
\end{definition}

\begin{proposition}
\label{nocoh}
Let $(K,\Sigma)\in{\mathcal H}(\aleph_1)$ be a strictly monotonic and spread
creating pair for $\bH$. Suppose that $\bar{t}=\langle t_n:n\in\omega\rangle
\in\PC_\infty(K,\Sigma)$ and $h_1,h_2\in\baire$ are such that 
\[(\forall n)(0<h_1(m^{t_n}_{\dn})\leq h_2(m^{t_n}_{\dn})\leq\nor[t_n])\quad
\mbox{and}\quad (\forall^\infty n)(h_2(m^{t_n}_{\dn})\leq\nor[t_n]-1).\]
Then every $(\bar{t},h_1,h_2)$--bounding forcing notion does not add Cohen
reals. 
\end{proposition}

\begin{proof}
Let $w\in\basis(t_0)$ be such that 
\[(\forall n\in\omega)\big(\pos(w,t_0,\ldots,t_{n-1})\subseteq\basis(t_n)
\big).\] 
Look at the space
\[{\mathcal X}=\{x\in\prod_{m\in\omega}\bH(m): (\forall n\in\omega)(x\rest
m^{t_n}_{\up}\in \pos(w,t_0,\ldots,t_n))\}\]
equipped with the natural (product) topology. It is a perfect Polish space
(note that as $(K,\Sigma)$ is strictly monotonic and spread, by
$\lim\limits_{n\to\infty}\nor[t_n]=\infty$, for sufficiently large $n$, for
each $u\in\basis(t_n)$ we find two distinct $v_0,v_1\in\pos(u,t_n)$). Thus, if
a forcing notion $\p$ adds a Cohen real then it adds a Cohen real $c\in
{\mathcal X}$. In $\V[c]$, choose (e.g.~inductively) a sequence $\bar{s}=
\langle s_n:n\in\omega\rangle\in U_{h_1}(\bar{t})$ such that 
\[(\forall n\in\omega)\big(s_n\in\Sigma(t_n)\ \ \&\ \ c\rest m^{s_n}_{\dn}\in
\pos(w,s_0,\ldots,s_{n-1})\subseteq\basis(s_n)\big)\]
(possible by \ref{monotonic}(3), remember that $(K,\Sigma)$ is nice). We claim
that there is no $\bar{s}^*\in U_{h_2}(\bar{t})\cap\V$ with $\bar{s}^*\leq
\bar{s}$. Why? Suppose that $\bar{s}^*\in U_{h_2}(\bar{t})\cap\V$. Working in
$\V$, consider the set 
\[{\mathcal O}\stackrel{\rm def}{=}\{x\in {\mathcal X}: (\exists n\in\omega)
(x\rest m^{s^*_n}_{\dn}\notin\pos(w,s^*_0,\ldots,s^*_{n-1}))\}.\]
This set is open dense in $\mathcal X$ (for the density use strict
monotonicity of $(K,\Sigma)$; remember that for sufficiently large $n\in
\omega$, $\nor[s^*_n]\leq h_2(m^{s^*_n}_{\dn})\leq\nor[t_m]-1$, where $m$ is
such that $m^{t_m}_{\dn}=m^{s^*_n}_{\dn}$). Consequently, in $\V[c]$, $c\in
{\mathcal O}$ and $\bar{s}^*\not\leq\bar{s}$.
\end{proof}

\begin{definition}
\label{additive}
Let $(K,\Sigma)$ be a weak creating pair. 
\begin{enumerate}
\item We say that a weak creature $t\in K$ is {\em $(n,m)$--additive} if for
all $t_0,\ldots,t_{n-1}\in \Sigma(t)$ such that $\nor[t_i]\leq m$ (for $i<n$)
there is $s\in\Sigma(t)$ such that  
\[t_0,\ldots,t_{n-1}\in\Sigma(s)\quad\mbox{ and }\quad\nor[s]\leq\max\{
\nor[t_\ell]: \ell<n\}+1.\]
\item {\em $m$--additivity of a weak creature $t\in K$} is defined as 
\[\add_{m}(t)=\sup\{k<\omega: t\mbox{ is $(k,m)$--additive}\}.\]
[Note that each $t$ is at least $(1,m)$--additive.]
\item We say that $(K,\Sigma)$ is $(g,h)$--additive (for $g,h\in\baire$) if
$\add_{h(m_{\dn}(t))}(t)\geq g(m_{\dn}(t))$ for all $t\in K$.\\
Similarly, if $(K,\Sigma)$ is a creating pair and $\bar{t}\in\PC(K,\Sigma)$
then we say that {\em $\bar t$ is $(g,h)$--additive} if $(\forall n\in\omega)
(\add_{h(m^{t_n}_{\dn})}(t_n)\geq g(m^{t_n}_{\dn}))$.
\item If the function $g$ is constant, say $g\equiv n$, then instead of
``$(g,h)$--additive'' we may say ``$(n,h)$--additive'' etc.
\end{enumerate}
\noindent [Note that for creatures we have $m_{\dn}(t)=m^t_{\dn}$.]
\end{definition}

\begin{remark}
The notion of additivity of a weak creature is very close to that of bigness:
in most applications they coincide. One can easily formulate conditions under
which $(k,m)$--additivity is equivalent to $k$--bigness.
\end{remark}

Let us recall that a forcing notion $\p$ has {\em the Laver property} if it is
$(f,g^*)$--bounding for every increasing function $f\in\baire$, where $g^*(n)
=2^n$ ($g^*$ may be replaced by any other fixed increasing function in
$\baire$).

\begin{proposition}
\label{fgLavimptF}
Assume that $(K,\Sigma)$ is a strongly finitary (see \ref{strfin}) and simple
(see \ref{simpglui}) creating pair, $\bar{t}=\langle t_n: n<\omega\rangle\in
\PC(K,\Sigma)$ and $h_1,h_2\in\baire$ are such that 
\[(\forall n\in\omega)(h_1(m^{t_n}_{\dn})\leq h_2(m^{t_n}_{\dn}))\quad\mbox{
and }\quad(\forall^\infty n)(h_1(m^{t_n}_{\dn})+1\leq h_2(m^{t_n}_{\dn})).\]
\begin{enumerate}
\item If $f\in\baire$ is such that $(\forall n\in\omega)(|V^n_{h_1}(\bar{t})|
\leq f(m^{t_n}_{\dn}))$, $g\in\baire$ is strictly increasing and $\bar{t}$ is
$(g,h_1)$--additive then every $(f,g)$--bounding forcing notion is
$(\bar{t},h_1,h_2)$--bounding. 
\item If $\bar{t}$ is $(g^*,h_1)$--additive (where $g^*(n)=2^n$) then every
forcing notion with Laver property is $(\bar{t},h_1,h_2)$--bounding. 
\end{enumerate}
\end{proposition}

\begin{proof}
1)\ \ \ Suppose that $\langle \dot{s}_n: n\in\omega\rangle$ is a $\p$-name for
an element of $U_{h_1}(\bar{t})$, $p\in\p$. Since $(K,\Sigma)$ is simple we
know that $p\forces\dot{s}_n\in\Sigma(t_n)$. Consequently, we may apply the
assumption that $\p$ is $(f,g)$--bounding (remember the property of $f$) and
we get a condition $p_0\geq p$ and a sequence $\langle s^+_{n,\ell}:\ell<
g^+(m^{t_n}_{\dn}),\ n<\omega\rangle$ such that
\[\begin{array}{l}
(\forall n<\omega)(\forall\ell<g(m^{t_n}_{\dn}))\big(s^+_{n,\ell}\in
V^n_{h_1}(\bar{t})\big)\qquad\mbox{ and}\\
p_0\forces_{\p}(\forall n\in\omega)(\dot{s}_n\in\{s^+_{n,\ell}:\ell<g^+(
m^{t_n}_{\dn})\}),
  \end{array}\]
where $g^+\in\baire$ is such that for each $n\in\omega$: 
\[g^+(m^{t_n}_{\dn})=\left\{\begin{array}{ll}
g(m^{t_n}_{\dn})&\mbox{if }h_1(m^{t_n}_{\dn})+1\leq h_2(m^{t_n}_{\dn}),\\ 
1               &\mbox{otherwise}.
\end{array}\right.\]
Since $\bar{t}$ is $(g,h_1)$-additive we find $\langle s^*_n:n\in\omega\rangle
\in U_{h_2}(\bar{t})$ such that for each $n\in\omega$
\[(\forall\ell<g^+(m^{t_n}_{\dn}))(s^+_{n,\ell}\in\Sigma(s^*_n)).\]
Clearly, $p_0\forces \langle s^*_n: n<\omega\rangle\leq \langle \dot{s}_n:
n<\omega\rangle$.
\medskip

\noindent 2)\ \ \ Similarly. 
\end{proof}

\begin{definition}
\label{tgoodetc}
Let $(K,\Sigma)$ be a creating pair and $\bar{t}=\langle t_n:n\in\omega\rangle
\in\PC_\infty(K,\Sigma)$. 
\begin{enumerate}
\item We say that a partial ordering $\bar{\F}=(\F,<^*_{\F})$ on $\F\subseteq
\baire$ is {\em $\bar{t}$--good} if:
\begin{enumerate}
\item[(a)] $\bar{\F}$ is a dense partial order with no maximal and minimal
elements, 
\item[(b)] for each $h\in\F$ 
\[(\forall n\in\omega)(1<h(m^{t_n}_{\dn})\leq\nor[t_n])\quad\mbox{ and }\quad
\lim_{n\to\infty} [\nor[t_n]-h(m^{t_n}_{\dn})]=\infty,\] 
\item[(c)] if $h_1,h_2\in \F$, $h_1<^*_{\F} h_2$ then 
\[(\forall n\in\omega)(h_1(m^{t_n}_{\dn})\leq h_2(m^{t_n}_{\dn}))\quad\mbox{
and }\quad\lim\limits_{n\to\infty} [h_2(m^{t_n}_{\dn})-h_1(m^{t_n}_{\dn})]=
\infty.\] 
\end{enumerate}
\item Let $\bar{\F}$ be a $\bar{t}$--good partial order (on $\F\subseteq
\baire$). We say that a forcing notion $\p$ is $(\bar{t},\bar{\F})$--bounding
if $\p$ is $(\bar{t},h_1,h_2)$--bounding for all $h_1,h_2\in\F$ such that $h_1
<^*_{\F}h_2$.
\end{enumerate}
\end{definition}

It should be clear that if $\bar{\F}$ is $\bar{t}$--good, $\bar{t}\in
\PC_\infty(K,\Sigma)$ then the composition of $(\bar{t},\bar{\F})$--bounding
forcing notions is $(\bar{t},\bar{\F})$--bounding. To deal with the limit
stages (in countable support iterations) we have to apply the technique of
\cite[Ch VI, \S1]{Sh:f}. Theorem \ref{pretFbound} below fulfills the promise
of \cite[7.2.29]{BaJu95}.

\begin{theorem}
\label{pretFbound}
Let $(K,\Sigma)\in\cH(\aleph_1)$ be a simple and reducible creating pair
and let $\bar{t}=\langle t_n:n<\omega\rangle\in\PC_\infty(K,\Sigma)$,
$\nor[t_n]\geq 2$. Suppose that $\bar{\F}=(\F,<^*_{\F})$ is a $\bar{t}$--good
partial order such that $\bar{t}$ is $(2,h)$--additive for all $h\in\F$.
Assume that $\langle\p_\alpha,\dot{\q}_\alpha:\alpha<\beta\rangle$ is a
countable support iteration of proper forcing notions such that for each
$\alpha<\beta$: 
\[\forces_{\p_\alpha}\mbox{``$\dot{\q}_\alpha$ is $(\bar{t},\bar{\F}
)$--bounding''.}\] 
Then $\p_\beta$ is $(\bar{t},\bar{\F})$--bounding.
\end{theorem}

\begin{proof}
We are going to use \cite[Ch VI, 1.13A]{Sh:f} and therefore we will closely
follow the notation and terminology of \cite[Ch VI, \S1]{Sh:f}, checking all
necessary assumptions. First we define a fine covering model $(D^{\bar{\F},
\bar{t}},R^{\bar{\F},\bar{t}},<^{\bar{\F},\bar{t}})$ (see \cite[Ch VI,
1.2]{Sh:f}).  

For $h\in\F$ and $n\in\omega$ we fix a mapping $\psi_n^h:\omega
\stackrel{\rm onto}{\longrightarrow}V^n_h(\bar{t})$. Next, for $h\in\F$ and
$\eta\in\baire$ we let $\psi^h(\eta)=\langle \psi^h_n(\eta(n)): n<\omega
\rangle$. Note that, as $(K,\Sigma)$ is nice, $\psi^h(\eta)\in U_h(\bar{t})$
(for all $\eta\in\baire$). Further, for $h\in\F$ and for $\bar{s}=\langle s_n:
n<\omega\rangle\geq\bar{t}$ we let  
\[T_{\bar{s},h}=\{\nu\in\fseo:(\forall n<\lh(\nu))(\psi^{h}_n(\nu(n))\in
\Sigma(s_n))\}.\] 
Clearly each $T_{\bar{s},h}$ is a subtree of $\fseo$ and any node in
$T_{\bar{s},h}$ has a proper extension in $T_{\bar{s},h}$ (remember that
$(K,\Sigma)$ is reducible). We define: 
\begin{itemize}
\item $D^{\bar{\F},\bar{t}}$ is $\cH(\aleph_1)=\cH(\aleph_1)^{\V}$ (we want to
underline here that -- in the iteration -- $D^{\bar{\F},\bar{t}}$ is fixed and
consists of elements of the ground universe), 
\item for $x,T\in D^{\bar{\F},\bar{t}}$ we say that \quad $x\;R^{\bar{\F},
\bar{t}}\;T$ \quad if and only if 

$x=\langle h^*,h\rangle$ and $T=T_{\bar{s},h^*}$ for some $h^*,h\in\F$,
$h^*<^*_{\F} h$, $\bar{s}\in U_h(\bar{t})\cap D^{\bar{\F},\bar{t}}$,
\item for $\langle h^*,h\rangle, \langle h^{**},h'\rangle\in\dom(R^{\bar{\F},
\bar{t}})$ we say that $\langle h^*,h\rangle<^{\bar{\F},\bar{t}} \langle
h^{**}, h'\rangle$ if and only if $h^*=h^{**}<^*_{\F} h<^*_{\F}h'$.
\end{itemize}

\begin{claim}
\label{cl21}
$(D^{\bar{\F},\bar{t}},R^{\bar{\F},\bar{t}})$ is a weak covering model in $\V$
(see \cite[Ch VI, 1.1]{Sh:f}).  
\end{claim}

\noindent{\em Proof of the claim:}\ \ \ The demand (a) of the definition
\cite[Ch VI, 1.1]{Sh:f} of weak covering models holds by the way we defined
$R^{\bar{\F},\bar{t}}$. The clause (b) there is satisfied as for each $\eta
\in\baire$ and $x=\langle h^*,h\rangle$ such that $h^*,h\in\F$, $h^*<^*_{\F}
h$, we have $\psi^{h^*}(\eta)\in U_{h^*}(\bar{t})\subseteq U_h(\bar{t})$,
$x\;R^{\bar{\F},\bar{t}}\;T_{\psi^{h^*}(\eta),h^*}$ and
$\eta\in\lim(T_{\psi^{h^*}(\eta),h^*})$. 

\begin{claim}
\label{cl22}
A forcing notion $\p$ is $(\bar{t},\bar{\F})$--bounding\quad if and only if 

\noindent it is $(D^{\bar{\F},\bar{t}},R^{\bar{\F},\bar{t}})$--preserving (see
\cite[Ch VI, 1.5]{Sh:f}). 
\end{claim}

\noindent{\em Proof of the claim:}\ \ \ Suppose $\p$ is
$(\bar{t},\bar{\F})$--bounding. Let $\dot{\eta}$ be a $\p$--name for an
element of $\baire$ and $x=\langle h^*,h\rangle\in\dom(R^{\bar{\F},\bar{t}})$
(so $h^*,h\in\F$ and $h^*<^*_{\F} h$). Then $\forces_{\p}\psi^{h^*}(
\dot{\eta})\in U_{h^*}(\bar{t})$ and, as by \ref{tgoodetc}(2) $\p$ is
$(\bar{t},h^*,h)$--bounding, 
\[\forces_{\p}(\exists\bar{s}\in U_h(\bar{t})\cap\V)(\bar{s}\leq\psi^{h^*}
(\dot{\eta}))\]
and hence
\[\forces_{\p}(\exists T\in D^{\bar{\F},\bar{t}})(x\;R^{\bar{\F},\bar{t}}\;
T\quad\&\quad\dot{\eta}\in\lim(T))\] 
(so $(D^{\bar{\F},\bar{t}},R^{\bar{\F},\bar{t}})$ covers in $\V^{\p}$). 

On the other hand suppose that $\p$ is $(D^{\bar{\F},\bar{t}},R^{\bar{\F},
\bar{t}})$--preserving. Take $h^*,h\in\F$ such that $h^*<^*_{\F} h$ and
let $\dot{\bar{s}}$ be a $\p$--name for an element of $U_{h^*}(\bar{t})$. Let
$\dot{\eta}$ be a $\p$--name for an element of $\baire$ such that
$\forces_{\p}\psi^{h^*}(\dot{\eta})=\dot{\bar{s}}$. As, in $\V^{\p}$,
$(D^{\bar{\F},\bar{t}},R^{\bar{\F},\bar{t}})$ still covers (and $x=\langle
h^*,h\rangle\in\dom(R^{\bar{\F},\bar{t}})$) we have
\[\forces_{\p}(\exists\bar{s}\in D^{\bar{\F},\bar{t}})(\bar{s}\in
U_h(\bar{t})\ \&\ \dot{\eta}\in\lim(T_{\bar{s},h^*})).\]
Hence we may conclude that $\p$ is $(\bar{t},h^*,h)$--bounding.

\begin{claim}
\label{cl23}
$(D^{\bar{\F},\bar{t}},R^{\bar{\F},\bar{t}},<^{\bar{\F},\bar{t}})$ is a fine
covering model (see \cite[Ch VI, 1.2]{Sh:f}). 
\end{claim}

\noindent{\em Proof of the claim:}\ \ \ We have to check the requirements of
\cite[Ch VI, 1.2(1)]{Sh:f}. We will comment on each of them, referring to the
enumeration there. 
\begin{enumerate}
\item[($\alpha$)] $(D^{\bar{\F},\bar{t}},R^{\bar{\F},\bar{t}})$ is a weak
covering model (by \ref{cl21}). 
\item[($\beta$)]  $<^{\bar{\F},\bar{t}}$ is a partial order on $\dom(R^{
\bar{\F},\bar{t}})=\{\langle h^*,h\rangle\in\F\times\F: h^*<^*_{\F} h\}$ such
that: 
\begin{enumerate}
\item[(i)]   there is no minimal element in $<^{\bar{\F},\bar{t}}$,
\item[(ii)]  $<^{\bar{\F},\bar{t}}$ is dense as $<^*_{\F}$ is such,
\item[(iii)] if $x_1<^{\bar{\F},\bar{t}} x_2$ (so $x_1=\langle h^*,h_1
\rangle$, $x_2=\langle h^*,h_2\rangle$ and $h^*<^*_{\F}h_1<^*_{\F}h_2$) and
$x_1\;R^{\bar{\F},\bar{t}}\;T$ then there is $T^*\in D^{\bar{\F},\bar{t}}$
such that $T\subseteq T^*$ and $x_2\;R^{\bar{\F},\bar{t}}\;T$ -- namely $T$ 
itself may serve as $T^*$,
\item[(iv)]  if $x_1,x_2\in\dom(R^{\bar{\F},\bar{t}})$, $x_1<^{\bar{\F},
\bar{t}} x_2$ and $x_1\;R^{\bar{\F},\bar{t}}\;T_1$,\ \ $x_1\;R^{\bar{\F},
\bar{t}}\; T_2$\\
then there is $T\in D^{\bar{\F},\bar{t}}$ such that 
\[x_2\;R^{\bar{\F},\bar{t}}\; T,\quad T_1\subseteq T,\quad\mbox{ and
}\quad (\exists n\in\omega)(\forall\nu\in T_2)(\nu\rest n\in T_1\ \Rightarrow\
\nu\in T).\] 
\end{enumerate}
\end{enumerate}
For $(\beta){\bf (iv)}$ we use the assumption that $\bar{t}$ is
$(2,h)$--additive for $h\in\F$. Let $x_1=\langle h^*,h_1\rangle$, $x_2=\langle
h^*,h_2\rangle$ (so $h^*<^*_{\F}h_1<^*_{\F}h_2$) and let $\bar{s}_\ell=\langle
s_{\ell,m}:m<\omega\rangle\in U_{h_1}(\bar{t})$ be such that $T_\ell=
T_{\bar{s}_\ell,h^*}$ (for $\ell=1,2$). By \ref{tgoodetc}(1c) we find $n<
\omega$ such that 
\[(\forall m\geq n)(h_1(m^{t_n}_{\dn})+1<h_2(m^{t_n}_{\dn})).\]
For each $m\geq n$ we choose $s_m\in\Sigma(t_m)$ such that 
\[s_{1,m},s_{2,m}\in\Sigma(s_m)\quad\mbox{ and }\quad\nor[s_m]\leq\max\{\nor[
s_{1,m}],\nor[s_{2,m}]\}+1<h_2(m^{t_n}_{\dn})\]
(remember that $t_m$ is $(2,h_1(m^{t_m}_{\dn}))$--additive). For $m<n$ we let
$s_m=s_{1,m}\in V^m_{h_1}(\bar{t})\subseteq V^m_{h_2}(\bar{t})$. Finally let
$\bar{s}=\langle s_m: m<\omega\rangle$. As $(K,\Sigma)$ is nice, and by the
choice of the $s_m$'s we have $\bar{s}\in U_{h_2}(\bar{t})$. Look at
$T\stackrel{\rm def}{=} T_{\bar{s},h^*}\in D^{\bar{\F},\bar{t}}$. By
definitions, $x_2\;R^{\bar{\F},\bar{t}}\; T$, $T_1\subseteq T$ (remember that
$\Sigma(s_{1,m})\subseteq \Sigma(s_m)$ for all $m\in\omega$) and if $\nu\in
T_2$, $\nu\rest n\in T_1$ then $\nu\in T$ (as $\Sigma(s_{2,m})\subseteq
\Sigma(s_m)$ for $m\geq n$ and $s_m=s_{1,m}$ for $m<n$). 
\medskip

Now comes the main part: conditions $(\gamma)$ and $(\delta)$ of \cite[Ch VI,
1.2(1)]{Sh:f}. As the second one is stronger, we will verify it only. Let us
state what we have to show.
\begin{enumerate}
\item[$(\delta)$] If $\V^*$ is a generic extension of $\V$ and, in $\V^*$,
$(D^{\bar{\F},\bar{t}},R^{\bar{\F},\bar{t}})$ is a weak covering model
(i.e.~it still covers) then the following two requirements are satisfied (in
$\V^*$).    
\begin{enumerate}
\item[(a)] If $x,x^+,x_n\in\dom(R^{\bar{\F},\bar{t}})$ and $T_n\in
D^{\bar{\F},\bar{t}}$ are such that for each $n\in\omega$: 
\[x_n<^{\bar{\F},\bar{t}} x_{n+1}<^{\bar{\F},\bar{t}} x^+<^{\bar{\F},\bar{t}}
x\quad\mbox{ and }\quad x_n\; R^{\bar{\F},\bar{t}}\; T_n\]
then there are $T^*\in D^{\bar{\F},\bar{t}}$ and $W\in\iso$ such that\ \
$x\;R^{\bar{\F},\bar{t}}\; T^*$\ \ and
\[\{\eta\in\baire\!: (\forall i\!\in\! W)(\eta\rest \min(W\setminus(i{+}1))\in
\bigcup_{\scriptstyle j<i,\atop \scriptstyle j\in W} T_j\cup T_0)\}\subseteq
\lim(T^*).\] 
\item[(b)] If $\eta,\eta_n\in\baire$ are such that $\eta\rest n=\eta_n\rest n$
for every $n\in\omega$ and $x\in\dom(R^{\bar{\F},\bar{t}})$ then
\[(\exists T\in D^{\bar{\F},\bar{t}})(\exists^\infty n)\big(x\; R^{\bar{\F},
\bar{t}}\; T\quad\&\quad\eta_n\in\lim(T)\big).\]
\end{enumerate}
\end{enumerate}
So suppose that $\V\subseteq\V^*$ is a generic extension and
$\V^*\models$``$(D^{\bar{\F},\bar{t}},R^{\bar{\F},\bar{t}})$ covers''. We work
in $\V^*$. 
\medskip

\noindent $(\delta)${\bf (a)}\hspace{0.15in} Let $h^*,h_n,h^+,h\in\F$ be such
that $x_n=\langle h^*,h_n\rangle$, $x=\langle h^*,h\rangle$, $x^+=\langle h^*,
h^+\rangle$ and $h^*<^*_{\F}h_n<^*_{\F}h_{n+1}<^*_{\F}h^+<^*_{\F}h$. Choose
inductively an increasing sequence $0<n_0<n_1<\ldots<\omega$ such that 
\[(\forall m\geq n_0)(h_0(m^{t_m}_{\dn})+1<h_1(m^{t_m}_{\dn})\leq
h^+(m^{t_m}_{\dn}))\quad\mbox{ and}\]
\[(\forall i\in\omega)(\forall m\geq n_{i+1})(h_{n_i}(m^{t_m}_{\dn})+i+2<
h_{n_i+1}(m^{t_m}_{\dn})\leq h^+(m^{t_m}_{\dn}))\]
(possible by \ref{tgoodetc}(1c)). Let $\bar{s}_n=\langle s_{n,m}\!:m\!<\!
\omega\rangle\in U_{h_n}(\bar{t})\cap D^{\bar{\F},\bar{t}}$ be such
that $T_n=T_{\bar{s}_n,h^*}$. (Note: each $\bar{s}_n$, $h_n$ is in $\V$ but
the sequences $\langle\bar{s}_n: n<\omega\rangle$, $\langle h_n: n<\omega
\rangle$ do not have to be there.) Using $(2,h^+)$--additivity of $\bar{t}$ we
choose $s_m^+\in\Sigma(t_m)$ such that for all $m\in\omega$ we have
$\nor[s^+_m]\leq h^+(m^{s_m}_{\dn})$ and
\begin{quotation}
\noindent if $m<n_1$ then $s^+_m=s_{0,m}$ (so $s_{0,m}\in\Sigma(s^+_m)$),

\noindent if $n_{i+1}\leq m<n_{i+2}$, $i<\omega$ then $s_{0,m},s_{n_0,m},
\ldots,s_{n_i,m}\in\Sigma(s^+_m)$ 
\end{quotation}
(remember the choice of the $n_i$'s). As $(K,\Sigma)$ is nice we have
$\bar{s}^+=\langle s^+_m: m<\omega\rangle\in U_{h^+}(\bar{t})$. Since $h^+
<^*_{\F} h$ and $\V^*\models$``$(D^{\bar{\F},\bar{t}},R^{\bar{\F},\bar{t}})$
covers'' we find $\bar{s}\in U_h(\bar{t})\cap D^{\bar{\F},\bar{t}}$ such that
$\bar{s}\leq\bar{s}^+$ (compare the proof of \ref{cl22}). Look at $T^*
\stackrel{\rm def}{=} T_{\bar{s},h^*}\in D^{\bar{\F},\bar{t}}$. By the
definitions, $x\; R^{\bar{\F},\bar{t}}\; T^*$. Let $W\stackrel{\rm def}{=}\{
n_i: i\in\omega\}$. Suppose that $\eta\in\baire$ is such that 
\[(\forall i\in\omega)(\eta\rest n_{i+1}\in\bigcup_{j<i} T_{n_j}\cup T_0).\]
This means that for each $m\in\omega$:
\begin{quotation}
\noindent if $m<n_1$ then $\eta\rest (m+1)\in T_0=T_{\bar{s}_0,h^*}$ and so 
\[\psi^{h^*}_m (\eta(m))\in\Sigma(s_{0,m})\subseteq\Sigma(s^+_m)\subseteq
\Sigma(s_m),\]
if $n_{i+1}\leq m<n_{i+2}$, $i\in\omega$ then $\eta\rest(m+1)\in T_0\cup
T_{n_0}\cup\ldots\cup T_{n_i}$ and so
\[\psi^{h^*}_m (\eta(m))\in\Sigma(s_{0,m})\cup\Sigma(s_{n_0,m})\cup\ldots\cup
\Sigma(s_{n_i,m})\subseteq\Sigma(s^+_m)\subseteq\Sigma(s_m).\]
\end{quotation}
Hence $(\forall m\in\omega)(\psi^{h^*}_m(\eta(m))\in\Sigma(s_m))$ what implies
that $\eta\in\lim(T_{\bar{s},h^*})=\lim(T^*)$, finishing the proof of
$(\delta)${\bf (a)}.
\medskip

\noindent $(\delta)${\bf (b)}\hspace{0.15in} This is somewhat similar to
$(\delta)${\bf (a)}. Let $h^*,h\in\F$ be such that $x=\langle h^*,h\rangle$
(so $h^*<^*_{\F} h$). As $\bar{\F}$ is a dense partial order we may take
$h^+\in\F$ such that $h^*<^*_{\F} h^+<^*_{\F} h$. Choose $0<n_0<n_1<\ldots<
\omega$ such that
\[(\forall i\in\omega)(\forall m\geq n_i)(h^*(m^{t_m}_{\dn})+i+2<h^+(
m^{t_m}_{\dn})).\]
Now take $s^+_m\in\Sigma(t_m)$ such that for $m\in\omega$ we have $\nor[s^+_m]
\leq h^+(m^{s_m}_{\dn})$ and 
\begin{quotation}
\noindent if $m<n_0$ then $s^+_m=\psi^{h^*}_m(\eta_{n_0}(m))$ (so
$\psi^{h^*}_m(\eta_{n_0}(m))\in\Sigma(s^+_m)$), and

\noindent if $n_i\leq m<n_{i+1}$, $i<\omega$ then $\psi^{h^*}_m(\eta_{n_0}
(m)),\ldots,\psi^{h^*}_m(\eta_{n_{i+1}}(m))\in\Sigma(s^+_m)$ 
\end{quotation}
(possible as $\bar{t}$ is $(2,h^+)$--additive; remember the choice of the
$n_i$'s). Note that, since $\eta\rest n=\eta_n\rest n$ for all $n\in\omega$,
we have that
\[(\forall i\in\omega)(\forall m\in\omega)(\psi^{h^*}_m(\eta_{n_i}(m))\in
\Sigma(s^+_m)).\]
Let $\bar{s}^+=\langle s^+_m: m<\omega\rangle$. Thus $\bar{s}^+\in U_{h^+}(
\bar{t})$ and, as $(D^{\bar{\F},\bar{t}},R^{\bar{\F},\bar{t}})$ covers in
$\V^*$ and $h^+<^*_{\F} h$, we find $\bar{s}=\langle s_m: m<\omega\rangle\in
U_h(\bar{t})\cap D^{\bar{\F},\bar{t}}$ such that $\bar{s}\leq\bar{s}^+$. Now
we have 
\[(\forall i\in\omega)(\forall m\in\omega)(\psi^{h^*}_m(\eta_{n_i}(m))\in
\Sigma(s_m)),\]
and therefore $(\forall i\in\omega)(\eta_{n_i}\in\lim(T_{\bar{s},h^*}))$. As
$x\; R^{\bar{\F},\bar{t}}\; T_{\bar{s},h^*}$, we finish the proof of the
claim.
\medskip

Now, to finish the proof of the theorem we put together \ref{cl22}, \ref{cl23}
and \cite[Ch VI, 1.13A]{Sh:f}.
\end{proof}

\section{Quasi-generic $\Gamma$ and preserving them}
Here we will develop the technique announced in \cite[Ch XVIII, 3.14, 3.15]
{Sh:f}, putting it in a slightly more general setting, more suitable for our
context. We will get a reasonably weak, but still easily iterable, condition
for not adding Cohen reals -- this will be used in \ref{tomek2},
\ref{tomek2conc}. But the general schema presented here will be applied in the
next chapter too (to preserve some ultrafilters on $\omega$). 

\begin{definition}
\label{quasi}
Suppose that $(K,\Sigma)$ is a creating pair, $\bar{t}=\langle t_k: k<\omega
\rangle\in\PC(K,\Sigma)$. 
\begin{enumerate}
\item A function $W:\omega\times\omega\times\pfs\longrightarrow {\mathcal
P}(K)$ is called a {\em $\bar{t}$--system for $(K,\Sigma)$} if:
\begin{enumerate}
\item[(a)] if $k\leq\ell<\omega$ and $\sigma:[m^{t_k}_{\dn},m^{t_\ell}_{\up})
\longrightarrow\omega$ then $W(m^{t_k}_{\dn},m^{t_\ell}_{\up},\sigma)
\subseteq\Sigma(t_k,\ldots,t_\ell)$, in all other instances
$W(m_0,m_1,\sigma)$ is empty,
\item[(b)] if $s\in\Sigma(t_k,\ldots,t_\ell)$, $k\leq\ell<\omega$ then there
is $n=n_W(s)\in [k,\ell]$ such that for each $\sigma_0,\sigma_1:
[m^{t_k}_{\dn},m^{t_\ell}_{\up})\longrightarrow\omega$ 

if \qquad $\sigma_0\rest [m^{t_n}_{\dn},m^{t_n}_{\up})=\sigma_1\rest
[m^{t_n}_{\dn},m^{t_n}_{\up})$ 

then \quad $s\in W(m^s_{\dn},m^s_{\up},\sigma_0)\quad\Leftrightarrow\quad
s\in W(m^s_{\dn},m^s_{\up},\sigma_1)$,
\item[(c)] if $k_0<k_1<\ldots<k_i$, $s_j\in\Sigma(t_{k_j},\ldots,t_{k_{j+1}-
1})$ for $j<i$, $s\in\Sigma(s_0,\ldots,s_{i-1})$ and $j_0<i$ is such that
$n_W(s)\in [k_{j_0},k_{j_0+1})$ (see (b) above) and $\sigma: [m^s_{\dn},
m^s_{\up})\longrightarrow\omega$ then 
\[s_{j_0}\in W(m^{s_{j_0}}_{\dn},m^{s_{j_0}}_{\up},\sigma\rest [m^{
s_{j_0}}_{\dn},m^{s_{j_0}}_{\up}))\quad\Rightarrow\quad s\in W(m^s_{\dn},
m^s_{\up},\sigma),\]
\item[(d)] for some unbounded non-decreasing function $G:(1,\infty)
\longrightarrow\mbR^{\geq 0}$ (called sometimes {\em the weight of $W$}), for
every $s\in\Sigma(t_k,\ldots,t_\ell)$, $k\leq\ell<\omega$, $\nor[s]>1$, and
each $\sigma:[m^{t_k}_{\dn},m^{t_\ell}_{\dn})\longrightarrow\omega$ there is
$t\in W(m^{t_k}_{\dn},m^{t_\ell}_{\up},\sigma)$ such that   
\[t\in\Sigma(s)\quad\mbox{ and }\quad\nor[t]\geq G(\nor[s]).\]
If the function $G$ might be $G(x)=x-1$ then we call the $\bar{t}$--system $W$
{\em regular}.
\end{enumerate}
\item For a norm condition $\C(\nor)$ we let 
\[\p^*_{\C(\nor)}(\bar{t},(K,\Sigma))\stackrel{\rm def}{=}\{\bar{s}\in\PC_{\C(
\nor)}(K,\Sigma): \bar{t}\leq\bar{s}\}.\]
It is equipped with the partial order $\leq$ inherited from $\PC_{\C(\nor)}(K,
\Sigma)$. We introduce another relation $\preceq_{\C(\nor)}$ on
$\p^*_{\C(\nor)}(\bar{t},(K,\Sigma))$ letting 
\begin{quotation}
\noindent $\bar{s}_0\preceq_{\C(\nor)}\bar{s}_1$\quad\quad if and only if

\noindent there is $\bar{s}_2\in\p^*_{\C(\nor)}(\bar{t},(K,\Sigma))$ such that
$\bar{s}_0\leq \bar{s}_2$ and the sequence $\bar{s}_2$ is eventually equal to
$\bar{s}_1$. 
\end{quotation}
If the norm condition $\C(\nor)$ is clear we may omit the index to $\preceq$. 
\item For a $\bar{t}$--system $W$ (for $(K,\Sigma)$) and $\Gamma\subseteq
\p^*_{\C(\nor)}(\bar{t},(K,\Sigma))$ we say that $\Gamma$ is {\em
quasi-$W$-generic in $\p^*_{\C(\nor)}(\bar{t},(K,\Sigma))$} if
\begin{enumerate}
\item[(a)] $(\Gamma,\preceq)$ is directed (i.e.\ $(\forall \bar{s}_0,\bar{s}_1
\in\Gamma)(\exists \bar{s}\in\Gamma)(\bar{s}_0\preceq\bar{s}\ \&\ \bar{s}_1
\preceq\bar{s})$) and countably closed (i.e.~if $\langle \bar{s}_n:n<\omega
\rangle\subseteq\Gamma$ is $\preceq$--increasing then there is $\bar{s}\in
\Gamma$ such that $(\forall n\in\omega)(\bar{s}_n\preceq\bar{s})$),
\item[(b)] for every function $\eta\in\baire$ there is $\bar{s}=\langle s_m:
m<\omega\rangle\in\Gamma$ such that  
\[(\forall^\infty m)(s_m\in W(m^{s_m}_{\dn},m^{s_m}_{\up},\eta\rest [m^{s_m}_{
\dn},m^{s_m}_{\up}))).\]
\end{enumerate}
\end{enumerate}
\end{definition}

\begin{remark}
The demand \ref{quasi}(1d) is to ensure the existence of quasi-$W$-generic
sets $\Gamma$ (see \ref{exists}(2) below). Conditions \ref{quasi}(1b) and
\ref{quasi}(1c) are to preserve quasi-$W$-genericity in countable support
iterations. As formulated, they will be crucial in the proof of
\ref{cl25}(2). 

Natural applications of the notions introduced in \ref{quasi} will be when
$(K,\Sigma)$ is simple or ``simple plus at most omitting''. In both cases it
will be easy to check demands \ref{quasi}(1b,c). In the first case they are
trivial, see \ref{better} below.
\end{remark}

\begin{proposition}
\label{better}
Suppose that $(K,\Sigma)$ is a creating pair and $\bar{t}\in\PC(K,\Sigma)$.  
\begin{enumerate}
\item If $(K,\Sigma)$ is simple then the condition \ref{quasi}(1b) is empty
(so may be omitted) and the condition \ref{quasi}(1c) is equivalent to  
\begin{enumerate}
\item[(c)$^-$] if $t\in W(m^{t_n}_{\dn},m^{t_n}_{\up},\sigma)$, $\sigma:[m^{
t_n}_{\dn},m^{t_n}_{\up})\longrightarrow\omega$ and $s\in\Sigma(t)$

then $s\in W(m^{t_n}_{\dn},m^{t_n}_{\up},\sigma)$.
\end{enumerate}
\item If $W$ is a $\bar{t}$--system for $(K,\Sigma)$, $\eta\in\baire$ and
$\bar{s}_\ell=\langle s_{\ell,m}:m<\omega\rangle\in\p^*_\emptyset(\bar{t},(K,
\Sigma))$ (for $\ell<2$) are such that 
\[s_0\preceq s_1\quad\mbox{ and }\quad (\forall^\infty m)(s_{0,m}\in W(m^{s_{
0,m}}_{\dn},m^{s_{0,m}}_{\up},\eta\rest [m^{s_{0,m}}_{\dn},m^{s_{0,m}}_{\up}
)))\]
then 
\[(\forall^\infty m)(s_{1,m}\in W(m^{s_{1,m}}_{\dn},m^{s_{1,m}}_{\up},\eta
\rest [m^{s_{1,m}}_{\dn},m^{s_{1,m}}_{\up}))).\]
\end{enumerate}
\end{proposition}

\begin{proposition}
\label{exists}
Suppose that $(K,\Sigma)$ is a creating pair, $\C(\nor)$ is one of the norm
conditions introduced in \ref{conditions} (i.e.~it is one of $({\rm
s}\infty)$, $(\infty)$, $({\rm w}\infty)$ or $(f)$ for some fast function $f$)
and $\bar{t}\in\PC_{\C(\nor)}(K,\Sigma)$.
\begin{enumerate}
\item $(\p^*_{\C(\nor)}(\bar{t},(K,\Sigma)),\preceq_{\C(\nor)})$ is a
countably closed partial ordering.
\item Assume CH. Further suppose that if $\C(\nor)\in\{({\rm s}\infty), ({\rm
w}\infty)\}$ then $(K,\Sigma)$ is growing. Let $W:\omega\times\omega\times\pfs
\longrightarrow{\mathcal P}(K)$ be a $\bar{t}$--system which is regular if
$\C(\nor)=(f)$ (for some fast function $f$). Then there exists a
quasi-$W$-generic $\Gamma$ in $\p^*_{\C(\nor)}(\bar{t},(K,\Sigma))$.
\end{enumerate}
\end{proposition}

\begin{proof}
We will show this for $\C(\nor)=(\infty)$. In other instances the proof
is similar and requires very small changes only.
\medskip

\noindent 1)\ \ \ It should be clear that $\preceq$ is a partial order on
$\p^*_{\infty}(\bar{t},(K,\Sigma))$. To show that it is countably closed
suppose that $\bar{s}_n=\langle s_{n,m}:m<\omega\rangle\in\p^*_{\infty}(
\bar{t},(K,\Sigma))$ are such that $\bar{s}_n\preceq\bar{s}_{n+1}$ for all
$n\in\omega$. Choose an increasing sequence $m_0<m_1<\ldots<\omega$ such that
$(\forall i\in\omega)(\exists m\in\omega)(m_i=m^{s_{i,m}}_{\dn})$ and for each
$n\leq i<\omega$:  

if $m<\omega$, $m_i\leq m^{s_{n,m}}_{\dn}$ then $\nor[s_{n,m}]\geq i$ and  

if $m<\omega$, $m_i\leq m^{s_{i,m}}_{\dn}$ then $(\exists m'<m''<\omega)(
s_{i,m}\in\Sigma(s_{n,m'},\ldots,s_{n,m''-1}))$.

\noindent Now choose $\bar{s}=\langle s_m: m<\omega\rangle\in\PC(K,\Sigma)$
such that  

if $m\in\omega$, $m^{s_0,m}_{\dn}<m_1$ then $s_m=s_{0,m}$,

if $m\in\omega$, $m^{s_{k+1},m}_{\dn}\in [m_{k+1},m_{k+2})$ then
$s_{k+1,m}=s_{m^*}$ for some $m^*\in\omega$.

\noindent Clearly the choice is possible (and uniquely determined) and, by the
niceness, we are sure that $\bar{s}\in\PC(K,\Sigma)$. Moreover, by the choice
of $m_i$'s, we have that $\bar{s}\in\PC_\infty(K,\Sigma)$ and so $\bar{s}\in 
\p^*_{\infty}(\bar{t},(K,\Sigma))$. Plainly, $\bar{s}_n\preceq\bar{s}$ for all
$n\in\omega$.   
\medskip

\noindent 2)\ \ \ First note that if $\bar{s}\in\p^*_{\infty}(\bar{t},(K,
\Sigma))$,  $\eta\in\baire$ then there is $\bar{s}^*=\langle s^*_m\!:m<\omega
\rangle\in\p^*_{\infty}(\bar{t},(K,\Sigma))$ such that $\bar{s}\leq\bar{s}^*$
and $(\forall^\infty m)(s^*_m\in W(m^{s^*_m}_{\dn},m^{s^*_m}_{\up}, \eta\rest
[m^{s^*_m}_{\dn},m^{s^*_m}_{\up})))$ (by \ref{quasi}(1d), remember that
the weight of $W$ is unbounded and non-decreasing).    

Using this remark, (1) above, and the assumption of CH we may build a
$\preceq$-increasing sequence $\langle\bar{s}_\alpha: \alpha<\omega_1\rangle
\subseteq\p^*_{\infty}(\bar{t},(K,\Sigma))$ such that 
\[(\forall \eta\in\baire)(\exists\alpha<\omega_1)(\forall^\infty n)(s_{
\alpha,n}\in W(m^{s_{\alpha,n}}_{\dn},m^{s_{\alpha,n}}_{\up}, \eta\rest [
m^{s_{\alpha,n}}_{\dn},m^{s_{\alpha,n}}_{\up}))).\] 
This sequence gives a quasi-$W$-generic $\Gamma=\{\bar{s}_\alpha: \alpha<
\omega_1\}$.\\
Note that proving (2) for $\C(\nor)\in\{({\rm s}\infty),({\rm w}\infty)\}$ we
have to assume something about the creating pair $(K,\Sigma)$. The assumption
that it is growing is the most natural one (in our context). It allows us to
obtain the respective version of the first sentence of the proof of (2) for
$(\infty)$. Similarly, if $\C(\nor)$ is $(f)$ then we need too assume
something about the weight of the system $W$. The assumption that $W$ is
regular is much more than really needed. 
\end{proof}

\begin{remark}
\label{exirem}
If $W_i$ are $\bar{t}$--systems (for $i\in\omega_1$) then we may construct in
a similar way (under CH) $\Gamma$ which is quasi-$W_i$-generic for all
$i<\omega_1$. 
\end{remark}

\begin{definition}
\label{Gampres}
Let $(K,\Sigma)$ be a creating pair, $\bar{t}\in\PC_{\C(\nor)}(K,\Sigma)$
(where $\C(\nor)$ is a norm condition). Suppose that $W:\omega\times\omega
\times\pfs\longrightarrow{\mathcal P}(K)$ is a $\bar{t}$--system for
$(K,\Sigma)$, and $\Gamma\subseteq\p^*_{\C(\nor)}(\bar{t},(K,\Sigma))$ is
quasi-$W$-generic. We say that a proper forcing notion $\p$ is {\em
$(\Gamma,W)$--genericity preserving} (or {\em $\Gamma$--genericity preserving}
if $W$ is clear) if \quad $\forces_{\p}\mbox{`` }\Gamma$ is
quasi-$W$-generic''. 
\end{definition}

\begin{remark}
\begin{enumerate}
\item Note that if $\p$ is a proper forcing notion and $\Gamma\subseteq
\p^*_{\C(\nor)}(\bar{t},(K,\Sigma))$ is quasi-$W$-generic then
\[\forces_{\p}\mbox{``$(\Gamma,\preceq)$ is directed and countably
closed''}.\]
Which may fail after the forcing is condition \ref{quasi}(3b), so the real
meaning of \ref{Gampres} is that this condition is preserved.
\item The composition of $\Gamma$--genericity preserving forcing notions is
clearly $\Gamma$--genericity preserving. To handle the limit stages in
countable support iterations we use the main result of \cite[Ch XVIII,
\S3]{Sh:f}, see \ref{Gamiter} below. 
\end{enumerate}
\end{remark}

\begin{definition}
\label{senetc}
Let $(K,\Sigma)$ be a creating pair, $\bar{t}\in\PC(K,\Sigma)$. We say that a
$\bar{t}$--system $W:\omega\times\omega\times\pfs\longrightarrow{\mathcal
P}(K)$ is
\begin{enumerate}
\item {\em Cohen--sensitive} if for all sufficiently large $m'\leq m''<\omega$
\[(\forall s\in\Sigma(t_{m'},\ldots,t_{m''}))(\exists\sigma:[m^s_{\dn},
m^s_{\up})\longrightarrow\omega)(s\notin W(m^s_{\dn},m^s_{\up},\sigma)),\]
\item {\em directed} if for every $m'\leq m''<\omega$, $\sigma_0,\ldots,
\sigma_m: [m^{t_{m'}}_{\dn},m^{t_{m''}}_{\up})\longrightarrow\omega$
($m<\omega$) there is $\sigma: [m^{t_{m'}}_{\dn},m^{t_{m''}}_{\up})
\longrightarrow\omega$ such that
\[W(m^{t_{m'}}_{\dn},m^{t_{m''}}_{\up},\sigma)\subseteq\bigcap_{\ell\leq m}
W(m^{t_{m'}}_{\dn},m^{t_{m''}}_{\up},\sigma_\ell).\] 
\end{enumerate}
\end{definition}

\begin{proposition}
\label{ganoco}
Suppose that $(K,\Sigma)$ is a creating pair, $\bar{t}\in\PC_{\C(\nor)}(K,
\Sigma)$, $W:\omega\times\omega\times\pfs\longrightarrow {\mathcal P}(K)$ is a
$\bar{t}$--system and $\Gamma\subseteq\p^*_{\C(\nor)}(\bar{t},(K,\Sigma))$ is
quasi-$W$-generic. Let $\p$ be a proper forcing notion.
\begin{enumerate}
\item If $W$ is Cohen--sensitive and $\p$ is $\Gamma$--genericity preserving
then $\p$ does not add Cohen reals.
\item If $W$ is directed, $(K,\Sigma)$ is simple and $\p$ is
$\baire$--bounding then $\p$ is $\Gamma$--genericity preserving.
\end{enumerate}
\end{proposition}

\begin{proof}
1)\ \ Note that if $\eta\in\baire$ is a Cohen real over $\V$, $W$ is
Cohen--sensitive and $\bar{s}=\langle s_m: m<\omega\rangle\in\PC(K,\Sigma)\cap
\V$ then
\[(\exists^\infty m)(s_m\notin W(m^{s_m}_{\dn},m^{s_m}_{\up},\eta\rest
[m^{s_m}_{\dn},m^{s_m}_{\up}))).\]
2)\ \ Suppose that $\dot{\eta}$ is a $\p$--name for a real in $\baire$, $p\in
\p$. As $\p$ is $\baire$--bounding we find a function $\eta\in\baire$ and a
condition $q\geq p$ such that $q\forces_{\p}(\forall n<\omega)(\dot{\eta}(n)<
\eta(n))$. Since $W$ is directed we can find $\eta^*\in\baire$ such that for
each $m<\omega$
\[W(m^{t_m}_{\dn},m^{t_m}_{\up},\eta^*\rest [m^{t_m}_{\dn},m^{t_m}_{\up}))
\subseteq\bigcap\{W(m^{t_m}_{\dn},m^{t_m}_{\up},\sigma): \sigma\in\prod_{
m^{t_m}_{\dn}\leq k<m^{t_m}_{\up}}\eta(k)\}.\]
Next, as $\Gamma$ is quasi-$W$-generic, we find $\bar{s}=\langle s_m: m<\omega
\rangle\in\Gamma$ such that
\[(\forall^\infty m)(s_m\in W(m^{s_m}_{\dn},m^{s_m}_{\up},\eta^*\rest
[m^{s_m}_{\dn}, m^{s_m}_{\up}))).\]
We finish noting that $(\forall m<\omega)(s_m\in\Sigma(t_m))$, as $(K,\Sigma)$
is simple.
\end{proof}

\begin{definition}
\label{coherent}
Suppose that $(K_\ell,\Sigma_\ell)$ are simple creating pairs and
$\bar{t}_\ell=\langle t_{\ell,m}:m<\omega\rangle\in\PC(K_\ell,\Sigma_\ell)$
(for $\ell<2$) are such that $(\forall m<\omega)(m^{t_{0,m}}_{\dn}=m^{
t_{1,m}}_{\dn})$. 
\begin{enumerate}
\item Let $h_0,h_1\in\baire$. We say that $\bar{t}_1$--systems
$W_0,W_1:\omega\times\omega\times\pfs\longrightarrow{\mathcal P}(K_1)$ (for
$(K_1,\Sigma_1)$) are {\em $(\bar{t}_0,h_0,h_1)$--coherent} if there are
functions $\rho_0,\rho_1$ (called {\em $(\bar{t}_0,h_0,h_1)$--coherence
witnesses}) such that
\begin{enumerate}
\item[(a)] $\rho_0,\rho_1:\pfs\longrightarrow\bigcup\{\Sigma_0(t_{0,m}): m<
\omega\}$,
\item[(b)] if $\sigma:[m^{t_{1,m}}_{\dn},m^{t_{1,m}}_{\up})\longrightarrow
\omega$, $m<\omega$ then $\rho_\ell(\sigma)\in\Sigma_0(t_{0,m})$ and $\nor[
\rho_\ell(\sigma)]\leq h_\ell(m^{t_0,m}_{\dn})$ (for $\ell=0,1$),
\item[(c)] for sufficiently large $m<\omega$, for every $\sigma_0:[m^{t_{
1,m}}_{\dn},m^{t_{1,m}}_{\up})\longrightarrow \omega$, there is $\sigma^*_0:
[m^{t_{1,m}}_{\dn},m^{t_{1,m}}_{\up})\longrightarrow\omega$ such that  
\[W_1(m^{t_{1,m}}_{\dn},m^{t_{1,m}}_{\up},\sigma_1)\subseteq W_0(m^{t_{1,
m}}_{\dn},m^{t_{1,m}}_{\up},\sigma_0)\]
whenever $\sigma_1:[m^{t_{1,m}}_{\dn},m^{t_{1,m}}_{\up})\longrightarrow\omega$
is such that $\rho_0(\sigma^*_0)\in\Sigma_0(\rho_1(\sigma_1))$;\\
the sequence $\sigma^*_0$, as well as $\rho_0(\sigma^*_0)$, will be called
{\em a $\rho_0$--cover for} $\sigma_0$, 
\item[(d)] for each $s\in\Sigma_0(t_{0,m})$, $m<\omega$, if $\nor[s]\leq
h_1(m^{t_{0,m}}_{\dn})$ then there is $\sigma: [m^{t_{1,m}}_{\dn},m^{t_{1,m
}}_{\up})\longrightarrow\omega$ such that $\rho_1(\sigma)=s$.
\end{enumerate}
\item Suppose that $\bar{\F}=(\F,<^*_{\F})$ is a $\bar{t}_0$--good partial
order (see \ref{tgoodetc}). We say that a family $\langle W_h^k: h\in\F,
k\in\omega\rangle$ of $\bar{t}_1$--systems for $(K_1,\Sigma_1)$ is {\em
$(\bar{t}_0,\bar{\F})$--coherent} if for every $k\in\omega$ and $h_0\in\F$
there is $h_1\in\F$ such that $h_0 <^*_{\F} h_1$ and the systems $W_{h_0}^k,
W_{h_1}^{k+1}$ are $(\bar{t}_0,h_0,h_1)$--coherent. 
\end{enumerate}
\end{definition}

\begin{theorem}
\label{gamcoher}
Let $(K_0,\Sigma_0),(K_1,\Sigma_1)$ be simple creating pairs and let 
\[\bar{t}_0\in\PC_\infty(K_0,\Sigma_0),\qquad \bar{t}_1\in\PC_{\C(\nor)}
(K_1,\Sigma_1)\]
be such that $(\forall m<\omega)(m^{t_{0,m}}_{\dn}=m^{t_{1,m}}_{\dn})$. Assume
that $\bar{\F}=(\F,<^*_{\F})$ is a $\bar{t}_0$--good partial order and
$\langle W_h^k: h\in\F, k\in\omega\rangle$ is a $(\bar{t}_0,\bar{\F}
)$--coherent family of $\bar{t}_1$--systems for $(K_1,\Sigma_1)$. Further
suppose that $\Gamma\subseteq\p^*_{\C(\nor)}(\bar{t}_1,(K_1,\Sigma_1))$ is
quasi-$W_h^k$-generic for all $h\in\F$, $k\in\omega$. {\em Then} every
$(\bar{t}_0,\bar{\F})$--bounding proper forcing notion is
$(\Gamma,W_h^k)$--genericity preserving for all $h\in\F$, $k\in\omega$. 
\end{theorem}

\begin{proof}
Let $\p$ be a proper $(\bar{t}_0,\bar{\F})$--bounding forcing notion. We have
to show that for all $h\in\F$, $k\in\omega$ 
\[\forces_{\p}\mbox{``$\Gamma$ is quasi-$W_h^k$-generic''}.\]
For this suppose that $\dot{\eta}$ is a $\p$--name for a real in $\baire$,
$p\in\p$, $k\in\omega$ and $h_0\in\F$. Take $h_1\in\F$ such that $h_0<^*_{\F}
h_1$ and the systems $W^k_{h_0}$, $W^{k+1}_{h_1}$ are $(\bar{t}_0,h_0,
h_1)$--coherent and let functions $\rho_0,\rho_1:\pfs\longrightarrow\bigcup\{
\Sigma_0(t_{0,m}): m\in\omega\}$ witness this fact. Let $\dot{\bar{s}}=
\langle\dot{s}_m:m<\omega\rangle$ be a $\p$--name for an element of $U_{h_0}
(\bar{t}_0)$ such that for some $N\in\omega$
\[p\forces_{\p}\mbox{`` }(\forall m\geq N)(\dot{s}_m\mbox{ is a
$\rho_0$--cover for }\dot{\eta}\rest[m^{t_{1,m}}_{\dn},m^{t_{1,m}}_{\up}))
\mbox{ ''}\qquad\mbox{ (see \ref{coherent}(1c))}\] 
(remember clause \ref{coherent}(1b)). Now, as $\p$ is
$(\bar{t}_0,\bar{\F})$--bounding and $h_0<^*_{\F} h_1$, we find a condition
$q\geq p$ and $\bar{s}^*=\langle s^*_m: m<\omega\rangle\in U_{h_1}(\bar{t}_0)$
such that $q\forces_{\p} \bar{s}^*\leq\dot{\bar{s}}$. Since $(K_0,\Sigma_0)$
is simple this means that $q\forces_{\p}(\forall m\in\omega)(\dot{s}_m\in
\Sigma_0(s^*_m))$. Let $\eta\in\baire$ be such that for each $m\in\omega$ we
have $\rho_1(\eta\rest[m^{t_{1,m}}_{\dn},m^{t_{1,m}}_{\up}))=s^*_m$ (remember
$\nor[s^*_m]\leq h_1(m^{s^*_m}_{\dn})$; see clause \ref{coherent}(1d)). We
know that $\Gamma$ is quasi-$W^{k+1}_{h_1}$-generic, so there is $\bar{s}=
\langle s_m:m<\omega\rangle\in\Gamma$ such that
\[(\forall^\infty m)(s_m\in W^{k+1}_{h_1}(m^{s_m}_{\dn},m^{s_m}_{\up},\eta
\rest [m^{s_m}_{\dn},m^{s_m}_{\up}))).\]
But $(K_1,\Sigma_1)$ is simple too, so $m^{s_m}_{\dn}=m^{s^*_m}_{\dn}$,
$m^{s_m}_{\up}=m^{s^*_m}_{\up}$. Thus we may apply \ref{coherent}(1c) and
conclude that for sufficiently large $m$
\[q\forces_{\p} s_m\in W^{k+1}_{h_1}(m^{s_m}_{\dn},m^{s_m}_{\up},\eta\rest
[m^{s_m}_{\dn},m^{s_m}_{\up}))\subseteq W^k_{h_0}(m^{s_m}_{\dn},m^{s_m}_{\up},
\dot{\eta}\rest [m^{s_m}_{\dn},m^{s_m}_{\up})),\]
so we are done.
\end{proof}

\begin{theorem}
\label{Gamiter}
Suppose that $(K,\Sigma)\in\cH(\aleph_1)$ is a creating pair and $\C(\nor)$ is
a norm condition. Assume that 
\begin{quotation}
\noindent $\bar{t}=\langle t_k: k<\omega\rangle\in\PC_{\C(\nor)}(K,\Sigma)$,\\
$W:\omega\times\omega\times\pfs\longrightarrow{\mathcal P}(K)$ is a
$\bar{t}$--system and\\
$\Gamma\subseteq\p^*_{\C(\nor)}(\bar{t},(K,\Sigma))$ is a quasi-$W$-generic
for $(K,\Sigma)$.
\end{quotation}
Let $\langle\p_\alpha,\dot{\q}_\alpha:\alpha<\beta\rangle$ be a countable
support iteration of proper forcing notions such that for each $\alpha<\beta$:
\[\forces_{\p_\alpha}\mbox{`` }\dot{\q}_\alpha\mbox{ is $\Gamma$--genericity
preserving ''}.\] 
Then $\p_\beta$ is $\Gamma$--genericity preserving.
\end{theorem}

\begin{proof} We will use the preservation theorem \cite[Ch XVIII, 3.6]{Sh:f}
and therefore we will follow the notation and terminology of \cite[Ch XVIII,
\S3]{Sh:f}, checking all necessary details. First we have to define our
context $(\bar{R}^{\Gamma,W},S^{\Gamma,W},{\bf g}^{\Gamma,W})$ (see \cite[Ch
XVIII, 3.1]{Sh:f}).  

For each $m\in\{m^{t_k}_{\dn}: k\in\omega\}$ we fix a mapping 
\[\psi^m:\omega\stackrel{\rm onto}{\longrightarrow}\{s\in K: (\exists k\leq
\ell<\omega)(s\in\Sigma(t_k,\ldots,t_\ell)\quad\&\quad m^{t_k}_{\dn}=m)\}.\]
Next, for $\eta\in\baire$ we let $\psi(\eta)=\langle\psi^{m_0}(\eta(0)),
\psi^{m_1}(\eta(1)),\ldots\rangle$, where $m^{t_0}_{\dn}=m_0<m_1<m_2<\ldots<
\omega$ are chosen in such a way that $m^{\psi^{m_k}(\eta(k))}_{\up}=m_{k+1}$.
Note that $\psi(\eta)\in \p^*_\emptyset(\bar{t},(K,\Sigma))$, though it does
not have to be in $\p^*_{\C(\nor)}(\bar{t},(K,\Sigma))$ (we do not control the
norms). 

Now we choose  $(\bar{R}^{\Gamma,W},S^{\Gamma,W},{\bf g}^{\Gamma,W})$ such
that:
\begin{itemize}
\item $S^{\Gamma,W}$ is, in the ground model $\V$, the collection of all
intersections $N\cap\cH(\aleph_1)$, where $N$ is a countable elementary
submodel of $(\cH(\chi),{\in},{<^*_\chi})$; so $S^{\Gamma,W}\subseteq
([\cH(\aleph_1)^{\V}]^{\textstyle\leq\aleph_0})^{\V}$ is stationary (we could
replace $S^{\Gamma,W}$ by any stationary subset), 
\item for each $a\in S^{\Gamma,W}$ we let $d[a]=c[a]=\omega$ (so $d'[a]=c'[a]
=\omega$),
\item $\alpha^*=\omega$,
\item for $n<\alpha^*$ and $\eta,g\in\baire$ we let 

$\eta\; R_n\; g$ \quad\quad  if and only if 

if $\psi(g)=\langle s^g_m:m<\omega\rangle$ and $m\in\omega$ is such that
$m^{s^g_m}_{\up}\geq n$ then 
\[s^g_m\in W(m^{s^g_m}_{\dn},m^{s^g_m}_{\up}, \eta\rest [m^{s^g_m}_{\dn},m^{
s^g_m}_{\up})),\]
\item $\bar{R}^{\Gamma,W}$ is a three place relation such that $(\eta,n,g)\in
\bar{R}^{\Gamma,W}$ if and only if $\eta,g\in\baire$, $n\in\omega$ and $\eta\;
R_n\;g$

\noindent (note: this is a definition of a relation, not a fixed object from
$\V$), 
\item ${\bf g}^{\Gamma,W}=\langle {\bf g}_a: a\in S^{\Gamma,W}\rangle\subseteq
\baire$ is such that for every $a,a'\in S^{\Gamma,W}$:
\begin{enumerate}
\item[($\alpha$)] $\psi({\bf g}_a)=\langle{\bf g}_{a,m}^\psi:m<\omega\rangle
\in\Gamma$,
\item[($\beta$)]  $(\forall \eta\in a\cap\baire)(\forall^\infty m)({\bf
g}^\psi_{a,m}\in W(m^{{\bf g}^\psi_{a,m}}_{\dn},m^{{\bf g}^\psi_{a,m}}_{\up},
\eta\rest [m^{{\bf g}^\psi_{a,m}}_{\dn},m^{{\bf g}^\psi_{a,m}}_{\up})))$,
\item[($\gamma$)] if $a'\in a\cap S^{\Gamma,W}$, $\bar{s}\in a\cap\Gamma$ then
$\psi({\bf g}_{a'})\preceq\psi({\bf g}_a)$ and $\bar{s}\preceq\psi({\bf
g}_a)$.
\end{enumerate}
\end{itemize}
Note that we may choose ${\bf g}_a$ by $\in$--induction for $a\in
S^{\Gamma,W}$ considering all $a'\in a\cap S^{\Gamma,W}$, $\bar{s}\in a\cap
\Gamma$ and $\eta\in a\cap\baire$. So, before we choose ${\bf g}_a$, we first
take $\bar{s}_\eta=\langle s_{\eta,m}: m<\omega\rangle\in\Gamma$ for $\eta\in
a\cap\baire$ such that  
\begin{enumerate}
\item[($\beta^*$)]  \qquad$(\forall^\infty m)(s_{\eta,m}\in W(m^{s_{\eta,m}}_{
\dn},m^{s_{\eta,m}}_{\up},\eta\rest [m^{s_{\eta,m}}_{\dn},m^{s_{\eta,m}}_{
\up})))$
\end{enumerate}
(possible by \ref{quasi}(3b)). Next, as $(\Gamma,\preceq)$ is directed and
countably closed (by \ref{quasi}(3a)) we may find $\bar{s}_0\in\Gamma$ such
that 
\[(\forall\eta\in a\cap\baire)(\forall a'\in a\cap S^{\Gamma,W})(\forall
\bar{s}\in a\cap\Gamma)(\psi({\bf g}_{a'})\preceq\bar{s}_0\ \&\ \bar{s}
\preceq\bar{s}_0\ \&\ \bar{s}_\eta \preceq\bar{s}_0).\]  
Let ${\bf g}_a\in\baire$ be such that $\psi({\bf g}_a)=\bar{s}_0$. It is easy
to check that ${\bf g}_a$ is as required (in $(\alpha)$--$(\gamma)$ above; for
$(\beta)$ we use \ref{better}(2)). 

\begin{claim}
\label{cl24}
\begin{enumerate}
\item $(\bar{R}^{\Gamma,W},S^{\Gamma,W},{\bf g}^{\Gamma,W})$ covers in $\V$
(see \cite[Ch XVIII, 3.2]{Sh:f}), i.e.:

if $x\in\V$ then there is a countable elementary submodel $N$ of 

$(\cH(\chi),{\in},{<^*_\chi})$ such that $a\stackrel{\rm def}{=}N\cap\cH(
\aleph_1)^{\V}\in S^{\Gamma,W}$, $(\bar{R}^{\Gamma,W},S^{\Gamma,W},
{\bf g}^{\Gamma,W})$, 

$x\in N$ and $(\forall\eta\in N\cap\baire)(\exists n<\omega)(\eta\; R_n\;{\bf
g}_a)$.  
\item Let $\p$ be a proper forcing notion. Then the following conditions are
equivalent:
\begin{enumerate}
\item[$(\oplus)_1$] $\forces_{\p}$``$(\bar{R}^{\Gamma,W},S^{\Gamma,W},{\bf
g}^{\Gamma,W})$ covers'', 
\item[$(\oplus)_2$] $\p$ is $\Gamma$-genericity preserving,
\item[$(\oplus)_3$] {\em if}  $p\in\p$ and $N$ is a countable elementary
submodel of $(\cH(\chi),{\in},{<^*_\chi})$ such that $p,\p,S^{\Gamma,W},{\bf
g}^{\Gamma,W}\in N$, $a\stackrel{\rm def}{=}N\cap\cH(\aleph_1)\in S^{\Gamma,
W}$ {\em then} there is an $(N,\p)$--generic condition $q\in\p$ stronger than
$p$ and such that  
\[q\forces_{\p}\mbox{``}(\forall\eta\in\baire\cap N[\Gamma_{\p}])(\exists
n<\omega)(\eta\; R_n\; {\bf g}_a)\mbox{''}.\]
\end{enumerate} 
\end{enumerate}
\end{claim}

\noindent{\em Proof of the claim:}\ \ \ 1)\ \ By the choice of $(\bar{R}^{
\Gamma,W},S^{\Gamma,W},{\bf g}^{\Gamma,W})$ (see condition $(\beta)$ of the
choice of ${\bf g}^{\Gamma,W}$).  
\medskip

\noindent 2)\ \ Assume $(\oplus)_1$. Let $\dot{\eta}$ be a $\p$-name for a
real in $\baire$, $p\in\p$. By the assumption we find $q\geq p$, $a\in
S^{\Gamma,W}$ and a $\p$-name $\dot{N}$ for an elementary submodel such that 
\[q\forces_{\p}\mbox{``}\dot{N}\cap\cH(\aleph_1)^{\V}=a\quad\&\quad\dot{\eta}
\in\dot{N}\quad\&\quad(\forall\eta\in\baire\cap\dot{N})(\exists n<\omega)(
\eta\; R_n\;{\bf g}_a)\mbox{''}.\]
But, as $\psi({\bf g}_a)\in\Gamma$, this is enough to conclude $(\oplus)_2$
(see the definitions of $R_n$, $\bar{R}^{\Gamma,W}$). 
\smallskip

\noindent Now, suppose that $(\oplus)_2$ holds true. Let $N,p$ be as in the
assumptions of $(\oplus)_3$ (so $a\stackrel{\rm def}{=}N\cap\cH(\aleph_1)\in
S^{\Gamma,W}$). Let $q\in\p$ be any $(N,\p)$--generic condition stronger than
$p$. Then, by $(\oplus)_2$, the condition $q$ forces in $\p$ that
\[(\forall\eta\!\in\!\baire\!\cap\! N[\Gamma_{\p}])(\exists\langle s_m\!:m{<}
\omega\rangle\!\in\!\Gamma\!\cap\! N)(\forall^\infty m)(s_m\in
W(m^{s_m}_{\dn},m^{s_m}_{\up},\eta\rest [m^{s_m}_{\dn},m^{s_m}_{\up})))\]
(note that $\rng({\bf g}^{\Gamma,W})\in N$ is a cofinal subset of $\Gamma$). 
But now, using \ref{better}(2) and clause $(\gamma)$ of the choice of ${\bf 
g}^{\Gamma,W}$, we conclude
\[q\forces_{\p}(\forall\eta\in\baire\cap N[\Gamma_{\p}])(\exists
n<\omega)(\eta\; R_n\; {\bf g}_a).\]

\noindent The implication \qquad $(\oplus)_3\ \ \Rightarrow\ \ (\oplus)_1$
\qquad is straightforward.  

\begin{claim}
\label{cl25}
Suppose that $\p$ is a proper forcing notion such that 
\[\forces_{\p}\mbox{``}(\bar{R}^{\Gamma,W},S^{\Gamma,W},{\bf g}^{\Gamma,W})
\mbox{ covers''}.\]
Then:
\begin{enumerate}
\item If $\dot{\q}$ is a $\p$-name for a proper $\Gamma$-genericity preserving
forcing notion then 
\[\forces_{\p}\mbox{``}\dot{\q}\mbox{ is $(\bar{R}^{\Gamma,W},S^{\Gamma,W},
{\bf g}^{\Gamma,W})$--preserving (for Possibility $A^*$)''}\]
(see \cite[Ch XVIII, 3.4]{Sh:f}).
\item $\forces_{\p}$``$(\bar{R}^{\Gamma,W},S^{\Gamma,W},{\bf g}^{\Gamma,W})$
strongly covers in the sense of Possibility $A^*$'' (see \cite[Ch XVIII,
3.3]{Sh:f}).
\end{enumerate}
\end{claim}

\noindent{\em Proof of the claim:}\ \ \ 1)\ \ We have to show that the
following condition holds true in 
$\V^{\p}$: 
\begin{enumerate}
\item[($*$)] Assume
\begin{enumerate}
\item[(i)]   $\chi_1$ is large enough, $\chi>2^{\chi_1}$,
\item[(ii)]  $N$ is a countable elementary submodel of $(\cH(\chi),{\in},
{<^*_\chi})$, $a\stackrel{\rm def}{=} N\cap\cH(\aleph_1)^{\V}\in S^{\Gamma,W}$,
and $\dot{\q}, S^{\Gamma,W},{\bf g}^{\Gamma,W},\chi_1,\ldots\in N$,
\item[(iii)] $(\forall\eta\in\baire\cap N)(\exists n<\omega)(\eta\; R_n\; {\bf
g}_a)$, 
\item[(iv)]  $\dot{\eta}_0\in N$ is a $\dot{\q}$-name for a real in $\baire$,
\item[(v)]   $\eta^*_0\in\baire$,
\item[(vi)]  $p,p_n\in\dot{\q}\cap N$ are such that $p\leq_{\dot{\q}}p_n
\leq_{\dot{\q}} p_{n+1}$ for all $n\in\omega$,
\item[(vii)] $\eta^*_0,\langle p_n:n<\omega\rangle\in N$,
\item[(viii)]$(\forall x\in\omega)(\forall^\infty n)(p_n\forces_{\dot{\q}}
\dot{\eta}_0(x)=\eta^*_0(x))$, 
\item[(ix)]  $n_0<\omega$ is such that $\eta^*_0\; R_{n_0}\;{\bf g}_a$,
\item[(xi)]  there is a countable elementary submodel $N_1$ of $(\cH(\chi_1),
{\in},{<^*_{\chi_1}})$ such that $\dot{\q},S^{\Gamma,W},{\bf g}^{\Gamma,W}\in
N_1\in N$, and for every open dense subset ${\mathcal I}$ of $\dot{\q}$,
${\mathcal I}\in N_1$ for some $n\in\omega$ we have $p_n\in {\mathcal I}\cap
N_1$ (i.e.\ $\langle p_n: n\in\omega\rangle$ is a generic sequence over
$N_1$).  
\end{enumerate}
Then there is an $(N,\dot{\q})$--generic condition $q\in\dot{\q}$ stronger
than $p$ and such that
\begin{enumerate}
\item[(a)]   $q\forces_{\dot{\q}}$``$\dot{\eta}_0\; R_{n_0}\; {\bf g}_a$''
\qquad and 
\item[(b)]   $q\forces_{\dot{\q}}$``$(\forall\eta\in\baire\cap N[\Gamma_{
\dot{\q}}])(\exists n<\omega)(\eta\; R_n\;{\bf g}_a)$''. 
\end{enumerate}
\end{enumerate}
So suppose that $\chi_1,\chi,N,N_1,a,\dot{\eta}_0,\eta^*_0,n_0,p$ and
$p_n$ (for $n\in\omega$) are as in the assumptions of $(*)$. Passing to a
subsequence (in $N$) we may assume that $p_n\forces_{\dot{\q}}$``$\dot{\eta_0}
\rest n=\eta^*_0\rest n$''. Remember that we work in $\V^{\p}$.

So, as $(\bar{R}^{\Gamma,W},S^{\Gamma,W},{\bf g}^{\Gamma,W})$ covers and
$N\prec (\cH(\chi),{\in},{<^*_\chi})$, we find a countable elementary submodel
$N_2$ of $(\cH(\chi_1),{\in},{<^*_{\chi_1}})$ such that 
\[\dot{\q},S^{\Gamma,W},{\bf g}^{\Gamma,W},\dot{\eta}_0,\eta^*_0,\langle
p_n:n<\omega\rangle,N_1,\ldots\in N_2\in N,\quad a_2\stackrel{\rm def}{=}
N_2\cap\cH(\aleph_1)^{\V}\in S^{\Gamma,W}\] 
\[\mbox{and }\qquad (\forall\eta\in\baire\cap N_2)(\exists n<\omega)(\eta\;
R_n\;{\bf g}_{a_2}).\]
By the choice of ${\bf g}^{\Gamma,W}$ we know that $\psi({\bf g}_{a_2})
\preceq\psi({\bf g}_a)$ (as $a_2\in a$) and hence we find $m^*\in [n_0,
\omega)$ such that
\begin{quotation}
\noindent if $m<\omega$, $m^{{\bf g}^\psi_{a,m}}_{\up}\geq m^*$ then for some
$m'\leq m''<\omega$ we have
\[{\bf g}^\psi_{a,m}\in\Sigma({\bf g}^\psi_{a_2,m'},\ldots,{\bf g}^\psi_{a_2,
m''}).\]
\end{quotation}
Now, working in $N$, we inductively choose sequences $\langle n_\ell:\ell<
\omega\rangle$, $\langle k_\ell:\ell<\omega\rangle$, $\langle m_\ell:\ell<
\omega\rangle$, $\langle q_\ell:\ell<\omega\rangle$ and $\langle\sigma_\ell:
\ell<\omega\rangle$, all from $N$.
\smallskip

\noindent{\sc Step}\qquad $\ell=0$.\\
The $n_0$ is given already. Let $m_0$ be the first such that $m^*\leq m^{{\bf
g}^\psi_{a_2,m_0}}_{\dn}$ and let $k_0=n_0$, $q_0=p_{n_0}$, $\sigma_0=
\eta^*_0\rest m^{{\bf g}^\psi_{a_2,m_0}}_{\up}$. 
\smallskip

\noindent{\sc Step}\qquad $\ell+1$.\\
Suppose we have defined $n_\ell,k_\ell,m_\ell,q_\ell,\sigma_\ell$. We let
$n_{\ell+1}=m^{{\bf g}^\psi_{a_2,m_\ell}}_{\up}$ and we choose an
$(N_2,\dot{\q})$--generic condition $q_{\ell+1}\geq p_{n_{\ell+1}}$, integers
$k_{\ell+1}\in [n_{\ell+1},\omega)$ and $m_{\ell+1}\in (m_\ell,\omega)$, and a
finite function $\sigma_{\ell+1}: [n_{\ell+1},m^{{\bf g}^\psi_{a_2,m_{\ell+
1}}}_{\up})\longrightarrow\omega$ such that:
\begin{enumerate}
\item[($\alpha$)] $q_{\ell+1}\forces_{\dot{\q}}$``$(\forall\eta\in\baire\cap
N_2[\Gamma_{\dot{\q}}])(\exists n<\omega)(\eta\; R_n\;{\bf g}_{a_2})$'',
\item[($\beta$)]  $q_{\ell+1}\forces_{\dot{\q}}$``$\dot{\eta}_0\;
R_{k_{\ell+1}}\; {\bf g}_{a_2}$'',
\item[($\gamma$)] $m_{\ell+1}$ is the first such that $k_{\ell+1}\leq m^{{\bf
g}^\psi_{a_2,m_{\ell+1}}}_{\dn}$, 
\item[($\delta$)] $q_{\ell+1}\forces_{\dot{\q}}$``$\dot{\eta}_0\rest [n_{\ell+
1}, m^{{\bf g}^\psi_{a_2,m_{\ell+1}}}_{\up})=\sigma_{\ell+1}$''
\end{enumerate}
(possible as, in $\V^{\p}$, $\dot{\q}$ is a proper $\Gamma$-genericity
preserving forcing notion, remember \ref{cl24}(2)). Note that all parameters
needed for the construction are in $N$. After it we have  
\[n_0=k_0\leq m^{{\bf g}^\psi_{a_2,m_0}}_{\dn}<m^{{\bf g}^\psi_{a_2,m_0}}_{
\up}=n_1\leq k_1\leq m^{{\bf g}^\psi_{a_2,m_1}}_{\dn}<m^{{\bf g}^\psi_{a_2,
m_1}}_{\up}=n_2\leq k_2\leq\ldots\]
and $\sigma_0:[0,n_1)\longrightarrow\omega$, $\sigma_{\ell+1}: [n_{\ell+1},
n_{\ell+2})\longrightarrow\omega$. Let $\eta\stackrel{\rm def}{=}\sigma_0
\conc\sigma_1\conc\sigma_2\conc\ldots\in N\cap\baire$. By {\bf (iii)}, we find
$\ell>0$ such that $\eta\; R_{n_\ell}\;{\bf g}_a$. We claim that $q_\ell
\forces_{\dot{\q}}$``$\dot{\eta}_0\; R_{n_0}\; {\bf g}_a$''. If not, then we
find a condition $q\geq q_\ell$ and $m<\omega$ such that $m^{{\bf g}^\psi_{a,
m}}_{\up}\geq n_0$ and 
\[q\forces_{\dot{\q}}\mbox{``}{\bf g}^\psi_{a,m}\notin W(m^{{\bf
g}^\psi_{a,m}}_{\dn},m^{{\bf g}^\psi_{a,m}}_{\up},\dot{\eta}_0\rest [m^{{\bf
g}^\psi_{a,m}}_{\dn},m^{{\bf g}^\psi_{a,m}}_{\up}))\mbox{''}.\] 
Let $n=n_W({\bf g}^\psi_{a,m})$ (see \ref{quasi}(1b)) and let $m'<\omega$ be
such that  $m^{{\bf g}^\psi_{a_2,m'}}_{\dn}\leq m^{t_n}_{\dn}<m^{t_n}_{\up}
\leq m^{{\bf g}^\psi_{a_2,m'}}_{\up}$. Consider the following three
possibilities.  
\smallskip

\noindent{\sc Case 1:}\qquad $m^{{\bf g}^\psi_{a_2,m'}}_{\dn}\geq k_\ell$.\\
By the choice of $m_0$ and $m^*$ (remember $\ell>0$, so $m^{{\bf g}^\psi_{a,
m}}_{\up}>m^{{\bf g}^\psi_{a_2,m'}}_{\dn}\geq k_1>m^*$) we know that for some
$m''\leq m'''<\omega$ 
\[{\bf g}^\psi_{a,m}\in\Sigma({\bf g}^\psi_{a_2,m''},\ldots,{\bf g}^\psi_{a_2,
m'''})\]
and by the choice of $k_\ell$ we know that 
\[q_\ell\forces_{\dot{\q}}\mbox{``}{\bf g}^\psi_{a_2,m'}\in W(m^{{\bf
g}^\psi_{a_2,m'}}_{\dn}, m^{{\bf g}^\psi_{a_2,m'}}_{\up},\dot{\eta}_0\rest
[m^{{\bf g}^\psi_{a_2,m'}}_{\dn},m^{{\bf g}^\psi_{a_2,m'}}_{\up}))\mbox{''}.\]
By \ref{quasi}(1c) we conclude that
\[q_\ell\forces_{\dot{\q}}\mbox{``}{\bf g}^\psi_{a,m}\in W(m^{{\bf
g}^\psi_{a,m}}_{\dn},m^{{\bf g}^\psi_{a,m}}_{\up},\dot{\eta}_0\rest [m^{{\bf
g}^\psi_{a,m}}_{\dn},m^{{\bf g}^\psi_{a,m}}_{\up}))\mbox{''}\]
(remember the choice of $m'$, note $m''\leq m'\leq m'''$), a contradiction.
\smallskip

\noindent{\sc Case 2:}\qquad $m^{{\bf g}^\psi_{a_2,m'}}_{\dn}<n_\ell$.\\
Then, by the choice of $n_\ell$ (remember $\ell>0$), we have $m^{{\bf
g}^\psi_{a_2,m'}}_{\up}\leq n_\ell$ (and $m'\leq m_\ell$). As $q_\ell\forces_{
\dot{\q}}$``$\dot{\eta}_0\rest n_\ell=\eta^*_0\rest n_\ell$'' and $\eta^*_0\;
R_{n_0}\; {\bf g}_a$ we immediately get a contradiction (remember
\ref{quasi}(1b) and the choice of $m'$). So we are left with the following
possibility.  
\smallskip

\noindent{\sc Case 3:}\qquad $n_\ell\leq m^{{\bf g}^\psi_{a_2,m'}}_{\dn}<
k_\ell$ (so $m'<m_\ell$).\\
Now the choice of $\eta$, $n_\ell$ and clause $(\delta)$ of the choice of
$q_\ell$ work: we know that 
\[\begin{array}{l}
m^{{\bf g}^\psi_{a_2,m'}}_{\up}\leq m^{{\bf g}^\psi_{a_2,m_\ell}}_{\dn}<
m^{{\bf g}^\psi_{a_2,m_\ell}}_{\up}=n_{\ell+1},\\
q_\ell\forces_{\dot{\q}}\mbox{`` }\dot{\eta}_0\rest [n_\ell,n_{\ell+1})=
\eta\rest [n_\ell,n_{\ell+1})\mbox{ ''}\qquad\mbox{ and }\quad \eta\;
R_{n_\ell}\; {\bf g}_a.\end{array}\] 
Consequently we get a contradiction like in the previous cases.

Now, choosing an $(N,\dot{\q})$--generic condition $q\geq q_\ell$ such that 
\[q\forces_{\dot{\q}}\mbox{``}(\forall\eta\in\baire\cap N[\Gamma_{\dot{\q}}])
(\exists n<\omega)(\eta\; R_n\;{\bf g}_a)\mbox{''}\]
(possible by \ref{cl24}(2)) we finish.
\medskip

\noindent 2)\ \ Work in $\V^{\p}$. We know that $(\bar{R}^{\Gamma,W},S^{
\Gamma,W},{\bf g}^{\Gamma,W})$ covers. Clearly each $R_n$ (for $n<\omega$) is
(a definition of) a closed relation on $\baire$. So what is left are the
following two requirements: 
\begin{enumerate}
\item[$\otimes$] if $a_1,a_2\in S^{\Gamma,W}$, $a_1\in a_2$ then for every
$\eta\in\baire$ we have
\[(\exists n<\omega)(\eta\; R_n\; {\bf g}_{a_1})\quad\Rightarrow\quad (\exists
n<\omega)(\eta\; R_n\;{\bf g}_{a_2});\]
\item[$\oplus_1$] if $\dot{\q}$, $\dot{\eta}_0$, $p$, $N$, $N_1$, $G_1$, $n_0$
are such that (in $\V^{\p}$):
\begin{enumerate}
\item[(a)] $\dot{\q}$ is a proper forcing notion,
\item[(b)] $N\prec (\cH(\chi),{\in},{<^*_\chi})$ is countable, $a\stackrel{\rm
def}{=} N\cap \cH(\aleph_1)^\V\in S^{\Gamma,W}$, 
\[(\forall \eta\in\baire\cap N)(\exists n<\omega)(\eta\; R_n\; {\bf g}_a),\]
$\dot{\q}, S^{\Gamma,W}, {\bf g}^{\Gamma,W},\chi_1,\ldots\in N$, $p\in
\dot{\q}\cap N$, 
\item[(c)] $\dot{\eta}_0\in N$ is a $\dot{\q}$--name for a function in
$\baire$, 
\item[(d)] $\chi_1<\chi$ ($\chi_1$ large enough), $N_1\in N$, $N_1\prec (\cH(
\chi_1),{\in},{<^*_{\chi_1}})$ is countable, $\dot{\q}, p, S^{\Gamma,W}, {\bf
g}^{\Gamma,W},\dot{\eta}_0,\ldots\in N_1$, and $G_1\in N$ is a
$\dot{\q}$--generic filter over $N_1$, $p\in G_1$,
\item[(e)] $\dot{\eta}_0[G_1]\; R_{n_0}\; {\bf g}_a$,
\end{enumerate}
then for every $y\in N\cap \cH(\chi_1)$ there are $N_2, G_2$ satisfying the
parallel of clause (d) and such that $y\in N_2$ and $\dot{\eta}_0[G_2]\;
R_{n_0}\; {\bf g}_a$.
\end{enumerate}
Now, concerning $\otimes$, look at our choice of ${\bf g}_a$'s: by $(\gamma)$
we know that 
\[a_1\in a_2\quad\Rightarrow\quad\psi({\bf g}_{a_1})\preceq\psi({\bf g}_{
a_2}).\]
Thus we may use \ref{better}(2) and the definition of $\bar{R}^{\Gamma,W}$ to
get $\otimes$. 

\noindent To show $\oplus_1$ we proceed similarly as in the proof of (1)
above. So suppose that $\dot{\q}$, $\dot{\eta}_0$, $p$, $N$, $N_1$, $G_1$,
$n_0$ are as in the assumptions of $\oplus_1$ and $y\in N\cap\cH(\chi_1)$. We
work in the universe $\V^{\p}$. Choose a $\leq_{\dot{\q}}$--increasing
sequence $\langle p_n:n\in\omega\rangle\in N$, $\dot{\q}$--generic over $N_1$,
such that $\{p_n:n\in\omega\}\subseteq G_1$ and $p_n$ decides the value of
$\dot{\eta}_0\rest n$. Let $N_2\in N$ be a countable elementary submodel of
$(\cH(\chi_1),{\in},{<^*_{\chi_1}})$ such that $N_1,y,\ldots\in N_2$ and then
choose a countable $N_2^+\prec (\cH(\chi_1),{\in},{<^*_{\chi_1}})$ such that
$N_2\in N_2^+\in N$, $a_2\stackrel{\rm def}{=} N_2^+\cap \cH(\aleph_1)\in
S^{\Gamma,W}$ and 
\[(\forall\eta\in N_2^+\cap\baire)(\exists n\in\omega)(\eta\; R_n\; {\bf
g}_{a_2})\]
(remember that $(\bar{R}^{\Gamma,W},S^{\Gamma,W},{\bf g}^{\Gamma,W})$ covers
in $\V^{\p}$). Let $m^*\geq n_0$ be such that 
\[(\forall m\in \omega)\big (m^{{\bf g}^\psi_{a,m}}_{\up}\geq m^*\
\Rightarrow\ (\exists m'\leq m''<\omega)({\bf g}^\psi_{a,m}\in\Sigma({\bf
g}^\psi_{a_2,m'},\ldots,{\bf g}^\psi_{a_2,m''}))\big).\]
Next, working in $N$, construct inductively sequences $\langle n_\ell:\ell<
\omega\rangle$, $\langle\sigma_\ell:\ell<\omega\rangle$, $\langle G^\ell:
\ell<\omega\rangle$ (all from $N$) such that 
\begin{enumerate}
\item[$(\alpha)$] $G^\ell$ is $\dot{\q}$--generic over $N_2$, $p_{n_\ell}\in
G^\ell\in N^+_2$ (so $\dot{\eta}_0[G_\ell]\in N_2^+$), 
\item[$(\beta)$]  $\dot{\eta}_0[G^\ell]\; R_{n_{\ell+1}}\; {\bf g}_{a_2}$,
$n_{\ell+1}>n_\ell+m^*$ and $n_{\ell+1}= m^{{\bf g}^\psi_{a_2,k}}_{\dn}$ for
some $k\in\omega$,
\item[$(\gamma)$] $\sigma_\ell=\dot{\eta}_0[G^\ell]\rest [n_\ell,n_{\ell+1})$.
\end{enumerate}
Finally let $\eta=\dot{\eta}_0[G_1]\rest n_0\conc\sigma_0\conc\sigma_1\conc
\sigma_2\ldots\in\baire\cap N$ and let $\ell>0$ be such that $\eta\;
R_{n_\ell}\; {\bf g}_a$. As in (1), one shows now that $\dot{\eta}_0[G^\ell]\;
R_{n_0}\;{\bf g}_a$.

\begin{claim}
\label{cl26}
The forcing notion $\p_\beta$ is $(\bar{R}^{\Gamma,W},S^{\Gamma,W},{\bf
g}^{\Gamma,W})$--preserving and hence
\[\forces_{\p_\beta}\mbox{``}(\bar{R}^{\Gamma,W},S^{\Gamma,W},{\bf g}^{\Gamma,
W})\mbox{ covers''}.\]
\end{claim}

\noindent{\em Proof of the claim:}\ \ \ Due to \ref{cl25}, we may apply
\cite[Ch XVIII, 3.6(1)]{Sh:f} to get the conclusion.
\medskip

\noindent Putting together \ref{cl26} and \ref{cl24}(2) we finish the proof of
the theorem.
\end{proof}

\section{Examples}

\begin{example}
\label{supersystem}
Let $F\in\baire$ be strictly increasing.\\
There are increasing functions $f^F=f, g^F=g\in\baire$, and $(K^F_{
\ref{supersystem}},\Sigma^F_{\ref{supersystem}})=(K_{\ref{supersystem}},
\Sigma_{\ref{supersystem}})$, $(K^{\ell,F}_{\ref{supersystem}},\Sigma^{\ell,
F}_{\ref{supersystem}})=(K^\ell_{\ref{supersystem}},\Sigma^\ell_{
\ref{supersystem}})$,  $\bar{t}^F=\bar{t}$, $\bar{t}^F_\ell=\bar{t}_\ell$, 
$\bar{\F}^F_\ell=\bar{\F}_\ell=(\F_\ell,<^*_\ell)$ (for $\ell<\omega$) and
$\langle W^k_{\ell,h}: h\in\F_\ell,\ k,\ell<\omega\rangle$ such that for every
$\ell<\omega$: 
\begin{enumerate}
\item $(K_{\ref{supersystem}},\Sigma_{\ref{supersystem}})$,
$(K^\ell_{\ref{supersystem}},\Sigma^\ell_{\ref{supersystem}})$ are simple,
strongly finitary and forgetful creating pairs for $\bH_{\ref{supersystem}}$,
$\bH^\ell_{\ref{supersystem}}$, respectively; 
\item $\bar{t}\in\PC_\infty(K_{\ref{supersystem}},\Sigma_{\ref{supersystem}
})$, $\bar{t}_\ell\in\PC_\infty(K^\ell_{\ref{supersystem}},\Sigma^\ell_{
\ref{supersystem}})$; 
\item $\bar{\F}_\ell$ is a countable $\bar{t}_\ell$--good partial order;
\item $\bar{t}_\ell$ is $(2^g,h)$--additive for each $h\in\F_\ell$;
\item $\langle W^k_{\ell,h}: h\in\F_\ell,\ k\in\omega\rangle$ is a $(\bar{t
}_\ell,\bar{\F}_\ell)$--coherent sequence of regular $\bar{t}$--systems; 
\item each $W^k_{\ell,h}$ (for $h\in\F_\ell$, $k\in\omega$) is Cohen
sensitive;
\item $(\forall h\in\F_\ell)(\forall^\infty m)(|V^m_h(\bar{t}_\ell)|<f(m))$;
\item $(\forall m\in\omega)(g(m+1)=F(f(m)))$.
\end{enumerate}
Moreover, the sequence $\langle W^k_{\ell,h}: h\in\F_\ell,\ k,\ell<\omega
\rangle$ has the following property:
\begin{enumerate}
\item[$(\circledast)_{\ref{supersystem}}$] if $\Gamma\subseteq\p^*_\infty
(K_{\ref{supersystem}},\Sigma_{\ref{supersystem}})$ is
quasi-$W^k_{\ell,h}$-generic for every $k,\ell<\omega$, $h\in\F_\ell$ and
$\random$ is a measure algebra (i.e.~adding a number of random reals) then
$\random$ is $(\Gamma,W^k_{\ell,h})$--genericity preserving for all $k,\ell<
\omega$, $h\in\F_\ell$.
\end{enumerate}
\end{example}

\begin{proof}[Construction] Let $F\in\baire$ be a strictly increasing
function. We inductively define $f,g\in\baire$, $\langle n^*_i,\ell^*_i,k^*_i:
i<\omega\rangle\subseteq\omega$ and a function $\psi:\omega\times\omega
\longrightarrow\omega$ such that: 
\begin{enumerate}
\item[($\alpha$)]   $g(0)=F(1)$, $n^*_0>2^{2g(0)}$ satisfies $2^{2g(0)+1}<
\frac{(n^*_0)^{n^*_0}}{(n^*_0)!}$,
\item[($\beta$)]    $\ell^*_i$ is such that $n^*_i<\ell^*_i$ and 
\[(\ell^*_i\cdot n^*_i)^{n^*_i}\cdot\frac{(\ell^*_i\cdot n^*_i-n^*_i)!}{(
\ell^*_i\cdot n^*_i)!}\leq 2,\]
\item[($\gamma$)]   $k^*_i=n^*_i(\ell^*_i)^{n^*_i} + n^*_i +1$,
\item[($\delta_0$)] $\psi(i,0)=2^{g(i)(i+1)}\cdot (n^*_i)^{k^*_i}$,
\item[($\delta_1$)] $\psi(i,\ell+1)=[2^{2g(i)(i+1)+\psi(i,\ell)}\cdot
(\psi(i,\ell)!)]^{g(i)+1}\cdot (n^*_i)^{k^*_i}$, 
\item[($\epsilon$)] $f(i)=2^{\psi(i,i)}$, $g(i+1)=F(f(i))$,
\item[($\zeta$)]    $n^*_{i+1}$ is such that $2^{2g(i+1)}\cdot k^*_i<
n^*_{i+1}$ and
\[2^{2g(i+1)(i+1)^2+1}<\frac{(n^*_{i+1})^{n^*_{i+1}}}{(n^*_{i+1})!}.\]
\end{enumerate}
Why is the choice possible? For clauses $(\alpha)$, $(\zeta)$ remember that
$\lim\limits_{n\to\infty} \frac{n^n}{n!}=\infty$. For clause $(\beta)$ note
that
\[N^{n^*_i}\cdot\frac{(N-n^*_i)!}{N!}=\frac{N^{n^*_i}}{N\cdot (N-1)\cdot\ldots
\cdot(N-(n^*_i-1))} \stackrel{N\to\infty}{\longrightarrow} 1.\]
Now, we define $\bH_{\ref{supersystem}},\bH^\ell_{\ref{supersystem}}$. For
$i\in\omega$ we let 
\[\begin{array}{l}
\bH_{\ref{supersystem}}(i)=\\
\{e\subseteq{\mathcal P}([(n^*_i,k^*_i)]^{\textstyle n^*_i}): \bigcup e=(
n^*_i,k^*_i)\ \&\ (\forall u,u'\in e)(u\neq u'\ \Rightarrow\ u\cap u'=
\emptyset)\},
  \end{array}\]
\[\bH^\ell_{\ref{supersystem}}(i)=\{(j,x):j<\psi(i,\ell)\ \&\ x:(n^*_i,k^*_i)
\longrightarrow n^*_i\}.\]
A creature $t\in\CR[\bH_{\ref{supersystem}}]$ is in $K_{\ref{supersystem}}$ if
for some $m\in\omega$, $B\subseteq{\mathcal P}([(n^*_m,k^*_m)]^{\textstyle
n^*_m})$: 
\begin{itemize}
\item $\val[t]=\{\langle u,w\rangle\in\prod\limits_{i<m}\bH_{
\ref{supersystem}}(i)\times\prod\limits_{i\leq m}\bH_{\ref{supersystem}}(i):
u\vartriangleleft w\ \&\ B\subseteq w(m)\}$,
\item $\nor[t]=\frac{\log_2\left(1+\frac{k^*_m - n^*_m -1}{n^*_m}-|B|\right)}{
\log_2(\ell^*_m)}=\frac{\log_2\left(1+(\ell^*_m)^{n^*_m}-|B|\right)}{\log_2(
\ell^*_m)}$,
\item $\dis[t]=B$.
\end{itemize}
The composition operation $\Sigma_{\ref{supersystem}}$ on
$K_{\ref{supersystem}}$ is given by  
\[\Sigma_{\ref{supersystem}}(t)=\{s\in K_{\ref{supersystem}}: m^s_{\dn}=
m^t_{\dn}\ \&\ \dis[t]\subseteq\dis[s]\}\qquad\mbox{ for }t\in
K_{\ref{supersystem}}.\]
Now we define $(K^\ell_{\ref{supersystem}},\Sigma^\ell_{\ref{supersystem}})$
for $\ell<\omega$. The family $K^\ell_{\ref{supersystem}}$ consists of all
creatures $t\in\CR[\bH^\ell_{\ref{supersystem}}]$ such that for some $m\in 
\omega$ and a nonempty set $C\subseteq\bH^\ell_{\ref{supersystem}}(m)$ such
that $(\forall (j,f),(j'f')\in C)(j=j'\ \Rightarrow f=f')$ we have:
\begin{itemize}
\item $\val[t]=\{\langle u,w\rangle\in\prod\limits_{i<m}\bH^\ell_{
\ref{supersystem}}(i)\times\prod\limits_{i\leq m}\bH^\ell_{\ref{supersystem}}
(i):u\vartriangleleft w\ \&\ w(m)\in C\}$,
\item $\nor[t]=\frac{\log_2(|C|)}{g(m)}$, 
\item $\dis[t]=C$.
\end{itemize}
The composition operation $\Sigma^\ell_{\ref{supersystem}}$ on
$K^\ell_{\ref{supersystem}}$ is such that for $t\in
K^\ell_{\ref{supersystem}}$:
\[\Sigma^\ell_{\ref{supersystem}}(t)=\{s\in K^\ell_{\ref{supersystem}}:
m^s_{\dn}=m^t_{\dn}\ \&\ \dis[s]\hookrightarrow\dis[t]\},\]
where $\dis[s]\hookrightarrow\dis[t]$ means that there is an embedding
\[{\bf i}:\{j<\psi(m,\ell): (\exists x)((j,x)\in\dis[s])\}\stackrel
{1-1}{\longrightarrow}\{j<\psi(m,\ell): (\exists x)((j,x)\in\dis[t])\}\]
such that $(\forall j\in\dom({\bf i}))(\forall x)((j,x)\in\dis[s]\
\Rightarrow\ ({\bf i}(j),x)\in\dis[t])$. Later we may identify elements
$s_0,s_1\in\Sigma^\ell_{\ref{supersystem}}(t)$ such that $\dis[s_0]
\hookrightarrow\dis[s_1]$ and $\dis[s_1]\hookrightarrow\dis[s_0]$. Therefore
we may think that we have the following inequality:
\[|\Sigma^\ell_{\ref{supersystem}}(t)|\leq 2^{\dis[t]}.\]
It should be clear that $(K_{\ref{supersystem}},\Sigma_{\ref{supersystem}})$,
$(K^\ell_{\ref{supersystem}},\Sigma^\ell_{\ref{supersystem}})$ are strongly
finitary, simple and forgetful creating pairs. Now we have to define
$\bar{t}$, $\bar{t}_\ell$. The first is the minimal member of
$\PC(K_{\ref{supersystem}},\Sigma_{\ref{supersystem}})$: 
\[\bar{t}=\langle t_m: m<\omega\rangle\quad\mbox{ is such that }\quad
(\forall m\in\omega)(m^{t_m}_{\dn}=m\ \ \&\ \ \dis[t_m]=\emptyset).\]
Next, for each $m\in\omega$ (and $\ell\in\omega$), we choose $t_{\ell,m}\in
(K^\ell_{\ref{supersystem}},\Sigma^\ell_{\ref{supersystem}})$ such that
$m^{t_{\ell,m}}_{\dn}=m$ and 
\[(\forall x:(n^*_m,k^*_m)\longrightarrow n^*_m)\big(|\{j<\psi(m,\ell):
(j,x)\in\dis[t_{\ell,m}]\}|=\frac{\psi(m,\ell)}{(n^*_m)^{k^*_m-n^*_m-1}}
\big).\]
Then we let $\bar{t}_\ell=\langle t_{\ell,m}: m<\omega\rangle$. Note that 
\[\nor[t_m]=\frac{\log_2\left(1+(\ell^*_m)^{n^*_m}\right)}{\log_2(\ell^*_m)}
\longrightarrow \infty\quad\mbox{ and }\quad\nor[t_{\ell,m}]=\frac{\log_2(
\psi(m,\ell))}{g(m)}\longrightarrow\infty\]
(when $m$ goes to $\infty$, $\ell$ is fixed). Moreover, if $n$ is such that 
\[2^{g(m)(n+1)}<\frac{\psi(m,\ell)}{(n^*_m)^{k^*_m-n^*_m-1}}\qquad\qquad\mbox{
(where $m,\ell<\omega$)}\]
then the creature $t_{\ell,m}$ is $(2^{g(m)},n)$--additive. Why? Note that if
$\nor[s_i]\leq n$, $s_i\in\Sigma^\ell_{\ref{supersystem}}(t_{\ell,m})$ then
$\sum\limits_{i<2^{g(m)}}|\dis[s_i]|\leq 2^{g(m)(n+1)}$ and it is smaller than
$\frac{\psi(m,\ell)}{(n^*_m)^{k^*_m-n^*_m-1}}$, which is the number of
repetitions of each function from $(n^*_m)^{\textstyle (n^*_m,k^*_m)}$ in
$\dis[t_{\ell,m}]$.  

For each $\ell\in\omega$ we choose a countable $\bar{t}_\ell$--good partial
order $\bar{\F}_\ell=(\F_\ell,<^*_\ell)$ such that for every $h\in\F_\ell$:
\begin{enumerate}
\item[(i)]   $g(m)(h(m)+1)<\log_2\left(\frac{\psi(m,\ell)}{(n^*_m)^{k^*_m-
n^*_m-1}}\right)$ for all $m\in\omega$,
\item[(ii)]  there is $h^*\in\F_\ell$ such that 
\[h<^*_\ell h^*\quad\mbox{ and }\quad (\forall m\in\omega)(h^*(m)\leq h(m)+
m),\]
\item[(iii)] there is a function $h^\otimes_\ell\in\F_{\ell+1}$ such that 
\[h^\otimes_\ell(m)\geq g(m)(m+1)+\psi(m,\ell)+\log_2(\psi(m,\ell)!)\quad
\mbox{ (for all $m\in\omega$)}.\]
\end{enumerate}
There should be no problems in carrying the construction of the $\F_\ell$.
Note that we may do this inductively, building a linear order (and so it will
be isomorphic to rationals). The clause (iii) is not an obstacle (in the
presence of (i)) as $\psi(m,\cdot)$ is increasing fast enough:
\[\begin{array}{l}
(n^*_m+1)\log_2(n^*_m)+g(m)(g(m)(m+1)+\psi(m,\ell)+\log_2(\psi(m,\ell)!)+2)<\\
(n^*_m+1)\log_2(n^*_m)+(g(m)+1)(2g(m)(m+1)+\psi(m,\ell)+\log_2(\psi(m,\ell)!))
=\\
\log_2\left(\frac{\psi(m,\ell+1)}{(n^*_m)^{k^*_m-n^*_m-1}}\right).\\
\end{array}\]
Note that the clause (i) and the previous remark imply that $\bar{t}_\ell$
is $(2^g,h)$--additive for each $h\in\F_\ell$. Moreover, by the choice of the
function $f$ we have that for every $\ell\leq m<\omega$ and $h\in\F_\ell$
\[|V^m_h(\bar{t}_\ell)|<|\Sigma(t_{\ell,m})|\leq 2^{|\dis[t_{\ell,m}]|}=
2^{\psi(m,\ell)}\leq 2^{\psi(m,m)}=f(m).\]
Finally, we are going to define $\bar{t}$--systems $W^k_{\ell,h}$ for $k,\ell
\in\omega$ and $h\in\F_\ell$. First, for each $\ell\in\omega$, $h\in\F_\ell$
we fix a function $\rho^\ell_h:\pfs\longrightarrow\bigcup\limits_{m\in\omega}
V^m_h(\bar{t}_\ell)$ such that for $m\in\omega$:
\[\rho^\ell_h\rest \omega^{\textstyle [m,m+1)}:\omega^{\textstyle [m,m+1)}
\stackrel{\rm onto}{\longrightarrow}V^m_h(\bar{t}_\ell).\]
Next, for $m\in\omega$ and $\sigma:[m,m+1)\longrightarrow\omega$ (and $\ell,k
\in\omega$, $h\in\F_\ell$) we let:\\
if $m\leq k$ then $W^k_{\ell,h}(m,m+1,\sigma)=\Sigma_{\ref{supersystem}}
(t_m)$,\\
if $m>k$ then
\[\begin{array}{ll}
W^k_{\ell,h}(m,m+1,\sigma)=&\big\{s\in\Sigma_{\ref{supersystem}}(t_m):\mbox{
for some $u\in\dis[s]$ we have}\\
\ &\qquad |\{(j,x)\in\dis[\rho^\ell_h(\sigma)]: x[u]=n^*_m\}|<\frac{|\dis[
\rho^\ell_h(\sigma)]|}{2^{g(m)(m+1)(k+1)}}\big\}\\
  \end{array}\]
(in all other instances we let $W^k_{\ell,h}(m',m'',\sigma)=\emptyset$).

\begin{claim}
\label{cl27}
For each $k,\ell\in\omega$ and $h\in\F_\ell$, the function $W^k_{\ell,h}$ is a
Cohen--sensitive regular $\bar{t}$--system.
\end{claim}

\noindent{\em Proof of the claim:}\ \ \ First we have to check that
$W^k_{\ell,h}$ is a $\bar{t}$--system. Immediately by its definition we have
that \ref{quasi}(1a-c) are satisfied (remember $(K_{\ref{supersystem}},
\Sigma_{\ref{supersystem}})$ is simple; see \ref{better}(1)). What might be
problematic is \ref{quasi}(1d). So suppose that $k,\ell,m\in\omega$, $m>k$
(otherwise trivial), $h\in\F_\ell$, $\sigma:[m,m+1)\longrightarrow\omega$,
$s\in\Sigma_{\ref{supersystem}}(t_m)$, $\nor[s]>1$. The last means that 
\[\frac{k^*_m-n^*_m-1}{n^*_m}-|\dis[s]|>\ell^*_m-1.\]
Let $N=\ell^*_m\cdot n^*_m$. Choose a set $X\subseteq (n^*_m,k^*_m)$ such that
$|X|=N$ and $X\cap\bigcup\dis[s]=\emptyset$. Note that for each $(j,x)\in\dis[
\rho^\ell_h(\sigma)]$ we have
\[\frac{|\{u\in [X]^{\textstyle n^*_m}: x[u]=n^*_m\}|}{|[X]^{\textstyle
n^*_m}|}=\frac{|x^{-1}[\{0\}]\cap X|\cdot\ldots\cdot |x^{-1}[\{n^*_m-1\}]\cap
X|}{{N\choose n^*_m}}\leq\]
\[\left(\sum_{i<n^*_m}\frac{|x^{-1}[\{i\}]\cap X|}{n^*_m}\right)^{n^*_m}\cdot
\frac{(n^*_m)!\cdot (N-n^*_m)!}{N!}= N^{n^*_m}\cdot \frac{(N-n^*_m)!}{N!}
\cdot \frac{(n^*_m)!}{(n^*_m)^{n^*_m}}.\]
Now look at clause $(\beta)$ of the choice of $\ell^*_m$. It implies that
$N^{n^*_m}\cdot \frac{(N-n^*_m)!}{N!}\leq 2$ and hence
\[\frac{|\{u\in [X]^{n^*_m}: x[u]=n^*_m\}|}{|[X]^{\textstyle n^*_m}|}\leq
2\frac{(n^*_m)!}{(n^*_m)^{n^*_m}}< \frac{1}{2^{2g(m)m^2}}\leq\frac{1}{
2^{g(m)(m+1)(k+1)}}\] 
(the second inequality follows from clause $(\zeta)$ of the choice of $n^*_m$,
remember $m>k$). Consequently, applying the Fubini theorem, we find $u_0\in
[X]^{\textstyle n^*_m}$ such that
\[\frac{|\{(j,x)\in\dis[\rho^\ell_h(\sigma)]: x[u_0]=n^*_m\}|}{|\dis[
\rho^\ell_h(\sigma)]|}<\frac{1}{2^{g(m)(m+1)(k+1)}}.\]
Thus, choosing $s^*\in\Sigma_{\ref{supersystem}}(s)$ such that $\dis[s^*]=
\dis[s]\cup\{u_0\}$ we will have $\nor[s^*]\geq \nor[s]-1$ and $s^*\in W^k_{
\ell,h}(m,m+1,\sigma)$. This shows that each $W^k_{\ell,h}$ is a regular 
$\bar{t}$--system. 

Finally we show that $W^k_{\ell,h}$ is Cohen--sensitive. Suppose that
$s\in \Sigma_{\ref{supersystem}}(t_m)$, $m>k$. Choose a function $x^*:
(n^*_m,k^*_m)\longrightarrow n^*_m$ such that $(\forall u\in\dis[s])(x^*[u]=
n^*_m)$. Next take $\sigma:[m,m+1)\longrightarrow\omega$ such that
$\dis[\rho^\ell_h(\sigma)]=\{(j_0,f^*)\}$ for some $j_0<\psi(m,\ell)$. It
should be clear that $s\notin W^k_{\ell,h}(m,m+1,\sigma)$. 

\begin{claim}
\label{cl28}
For each $\ell\in\omega$, the sequence $\langle W^k_{\ell,h}: h\in\F_\ell,\
k\in\omega\rangle$ is $(\bar{t}_\ell,\bar{\F}_\ell)$--coherent.
\end{claim}

\noindent{\em Proof of the claim:}\ \ \ Let $h_0\in\F_\ell$, $k\in\omega$. By
the demand (ii) of the choice of $\F_\ell$ we find $h_1\in\F_\ell$ such that 
\[h_0<^*_\ell h_1\quad\mbox{ and }\quad (\forall m\in\omega)(h_1(m)\leq h_0(m)
+m).\]
We want to show that the systems $W^k_{\ell,h_0}$ and $W^{k+1}_{\ell,h_1}$ are
$(\bar{t}_\ell,h_0,h_1)$--coherent and that this is witnessed by the functions
$\rho^\ell_{h_0},\rho^\ell_{h_1}$. Clearly these functions satisfy the demands
(1a), (1b) and (1d) of \ref{coherent}, so what we have to check is
\ref{coherent}(1c) only.

Suppose that $m>k+1$, $\sigma_0:[m,m+1)\longrightarrow \omega$. Look at the
creature $s=\rho^\ell_{h_0}(\sigma_0)$. We know that $\nor[s]\leq h_0(m)$ and
hence $|\dis[s]|\leq 2^{g(m)h_0(m)}$. Since $2^{g(m)h_0(m)}<\frac{\psi(m,
\ell)}{(n^*_m)^{k^*_m-n^*_m-1}}$ (remember clause (i) of the choice of
$\F_\ell$) we find a creature $s^*\in\Sigma^\ell_{\ref{supersystem}}(t_{\ell,
m})$ such that $\dis[s]\hookrightarrow\dis[s^*]$, $\nor[s^*]=h_0(m)$
(i.e.~$|\dis[s^*]|=2^{g(m)h_0(m)}$) and for each function $x^*:(n^*_m,k^*_m)
\longrightarrow n^*_m$ we have
\[\frac{|\{(j,x)\in\dis[s]: x=x^*\}|}{|\dis[s]|}\leq 2\frac{|\{(j,x)\in\dis
[s^*]: x=x^*\}|}{|\dis[s^*]|}.\]
How? We just ``repeat'' each $(j,x)$ from $\dis[s]$ successively, till we get
the required size. We have enough space for this as the number of the required
repetitions for each function from $(n^*_m,k^*_m)$ to $n^*_m$ is less than
$2^{g(m)h_0(m)}$.\\ 
Take $\sigma^*_0:[m,m+1)\longrightarrow\omega$ such that $\rho^\ell_{h_0}(
\sigma^*_0)=s^*$. We want to show that this $\sigma^*_0$ is a
$\rho^\ell_{h_0}$--cover for $\sigma_0$. So suppose that $\sigma_1:[m,m+1)
\longrightarrow\omega$ is such that $\rho^\ell_{h_0}(\sigma^*_0)\in
\Sigma^\ell_{\ref{supersystem}}(\rho^\ell_{h_1}(\sigma_1))$. Let $t\in
W^{k+1}_{\ell,h_1}(m,m+1,\sigma_1)$. This means that we can find $u\in\dis[t]$
such that 
\[|\{(j,x)\in\dis[\rho^\ell_{h_1}(\sigma_1)]: x[u]=n^*_m\}|<\frac{|\dis[
\rho^\ell_{h_1}(\sigma_1)]|}{2^{g(m)(m+1)(k+2)}}\leq \frac{2^{g(m)h_1(m)}}{
2^{g(m)(m+1)(k+2)}}.\]
Hence, remembering that $\dis[s^*]\hookrightarrow\dis[\rho^\ell_{h_1}
(\sigma_1)]$ and the choice of $h_1$, we get
\[\begin{array}{l}
|\{(j,x)\in\dis[s^*]: x[u]=n^*_m\}|<\\
\displaystyle <\frac{2^{g(m)h_0(m)}}{2^{g(m)(m+1)(k+1)}}\cdot
\frac{2^{mg(m)}}{2^{g(m)(m+1)}}=\frac{1}{2^{g(m)}}\cdot
\frac{|\dis[s^*]|}{2^{g(m)(m+1)(k+1)}}.
  \end{array}\]
But we are interested in $s$. By the choice of $s^*$ we have
\[\frac{|\{(j,x)\in\dis[s]\!: x[u]\!=\!n^*_m\}|}{|\dis[s]|}\leq 2
\frac{|\{(j,x)\in\dis[s^*]\!: x[u]\!=\!n^*_m\}|}{|\dis[s^*]|}<
\frac{1}{2^{g(m)(m+1)(k+1)}}\]
and therefore $t\in W^k_{\ell,h_0}(m,m+1,\sigma_0)$. Thus we have proved
\[W^{k+1}_{\ell,h_1}(m,m+1,\sigma_1)\subseteq W^k_{\ell,h_0}(m,m+1,\sigma_0)\]
whenever $\rho^\ell_{h_0}(\sigma^*_0)\in\Sigma^\ell_{\ref{supersystem}}(
\rho^\ell_{h_1}(\sigma_1))$. This finishes the claim.

\begin{claim}
\label{cl29}
Suppose that $\Gamma\subseteq\p^*_\infty(K_{\ref{supersystem}},
\Sigma_{\ref{supersystem}})$ is quasi-$W^k_{\ell,h}$-generic for all
$k,\ell\in\omega$, $h\in\F_\ell$. Let $\random$ be a measure algebra. Then
\[\forces_{\random}\mbox{``$\Gamma$ is quasi-$W^k_{\ell,h}$-generic for all
$k,\ell\in\omega$, $h\in\F_\ell$''.}\]
\end{claim}

\noindent{\em Proof of the claim:}\ \ \ Let $\mu$ be a $\sigma$-additive
measure on the complete Boolean algebra $\random$. Let $k,\ell\in\omega$,
$h\in\F_\ell$. Suppose that $\dot{\eta}$ is a $\random$--name for a real in
$\baire$, $b\in\random^+$ (i.e.~$\mu(b)>0$). To simplify notation let us
define, for $m\in\omega$, 
\[K_m=|V_h^m(\bar{t}_\ell)|,\quad M_m=2^{g(m)(m+1)}\cdot K_m,\quad N_m=
\psi(m,\ell)!\]
and let $h^\otimes_\ell\in\F_{\ell+1}$ be the function given by the clause
(iii) of the choice of $\bar{\F}_{\ell+1}$. 

Fix $m\in\omega$ for a moment.

For $s\in V^m_h(\bar{t}_\ell)$ choose a creature $t(s)\in
\Sigma^{\ell+1}_{\ref{supersystem}}(t_{\ell+1,m})$ such that 
\begin{itemize}
\item $|\dis[t(s)]|=N_m$,
\item for each $x^*:(n^*_m,k^*_m)\longrightarrow n^*_m$
\[\frac{|\{(j,x)\in\dis[s]: x=x^*\}|}{|\dis[s]|}=\frac{|\{(j,x)\in\dis[t(s)]:
x=x^*\}|}{|\dis[t(s)]|}\] 
\end{itemize}
(possible be the choice of $N_m$ and the fact that each $x^*$ is repeated more
than $N_m$ times in $\dis[t_{\ell+1,m}]$). Further, we choose integers
$g(m,s)<M_m$ for $s\in V^m_h(\bar{t}_\ell)$ such that
\[\sum_{s\in V^m_h(\bar{t}_\ell)} g(m,s)=M_m\quad\mbox{ and }\quad |\frac{\mu
\big(\lbv \rho^\ell_h(\dot{\eta}\rest [m,m+1))=s\rbv_{\random}\cdot
b\big)}{\mu(b)}-\frac{g(m,s)}{M_m}|\leq\frac{1}{M_m}\]
(where $\lbv\cdot\rbv_{\random}$ stands for the Boolean value). Take a
creature $t^*_m\in\Sigma^{\ell+1}_{\ref{supersystem}}(t_{\ell+1,m})$ such that
for some sequence $\langle A_s: s\in V^m_h(\bar{t}_\ell)\rangle$ of disjoint
subsets of $\psi(m,\ell+1)$ we have
\begin{itemize}
\item $|A_s|=g(m,s)\cdot N_m$,
\item $\bigcup\{A_s: s\in V^m_h(\bar{t}_\ell)\}=\{j<\psi(m,\ell+1): (\exists
x)((j,x)\in\dis[t^*_m])\}$,
\item for some bijection $\pi_s:g(m,s)\times\dis[t(s)]\longrightarrow A_s$ we
have
\[(\forall k<g(m,s))(\forall (j,x)\in\dis[t(s)])\big((\pi_s(k,(j,x)),x)\in
\dis[t^*_m]\big).\]
\end{itemize}
Why is the choice of the $t^*_m$ possible? Note that our requirements imply
that 
\[|\dis[t^*_m]|=M_m\cdot N_m=2^{g(m)(m+1)}\cdot |V^m_h(\bar{t}_\ell)|\cdot
\psi(m,\ell)!<\frac{\psi(m,\ell+1)}{(n^*_m)^{k^*_m-n^*_m-1}}\]
and the last number says how often each function is repeated in $\dis[t_{\ell
+1,m}]$. Moreover
\[\begin{array}{l}
\nor[t^*_m]=\frac{g(m)(m+1)+\log_2(|V^m_h(\bar{t}_\ell)|)+\log_2(\psi(m,
\ell)!)}{g(m)}<\\
\frac{g(m)(m+1)+\psi(m,\ell)+\log_2(\psi(m,\ell)!)}{g(m)}\leq
\frac{h^\otimes_\ell(m)}{g(m)}<h^\otimes_\ell(m).
\end{array}\]
Thus $t^*_m\in V^m_{h^\otimes_\ell}(\bar{t}_{\ell+1})$.
\medskip

Let $\eta\in\baire$ be such that $\rho^{\ell+1}_{h^\otimes_\ell}(\eta\rest
[m,m+1))=t^*_m$. Since $\Gamma$ is
quasi-$W^{k+1}_{\ell+1,h^\otimes_\ell}$-generic, we find $\bar{s}=\langle s_m:
m<\omega\rangle\in \Gamma$ such that for some $m^*>k+4$
\[(\forall m\geq m^*)(s_m\in W^{k+1}_{\ell+1,h^\otimes_\ell}(m,m+1,\eta\rest
[m,m+1)).\]

Fix $m\geq m^*$ for a moment. We know that for some $u\in\dis[s_m]$ we have
\[|\{(j,x)\in\dis[t^*_m]: x[u]=n^*_m\}|<\frac{|\dis[t^*_m]|}{2^{g(m)(m+1)
(k+2)}}.\]
For $s\in V^m_h(\bar{t}_\ell)$ let
\[Y_m(s)\stackrel{\rm def}{=}\frac{|\{(j,x)\in\dis[s]: x[u]=n^*_m\}|}{|\dis[s]
|}=\frac{|\{(j,x)\in\dis[t(s)]: x[u]=n^*_m\}|}{|\dis[t(s)]|}\]
and note that $\sum\limits_s Y_m(s)\cdot g(m,s)\cdot N_m=|\{(j,x)\in\dis[
t^*_m]: x[u]=n^*_m\}|$. Let 
\[{\mathcal X}_m=\{s\in V^m_h(\bar{t}_\ell): Y_m(s)\geq \frac{1}{2^{g(m)(m+1)
(k+1)}}\}\]
(so $\rho^\ell_h(\sigma)\notin {\mathcal X}_m$ implies $s_m\in W^k_{\ell,h}
(m,m+1,\sigma)$). Note that 
\[\sum_{s\in {\mathcal X}_m}\frac{g(m,s)\cdot N_m}{2^{g(m)(m+1)(k+1)}}\leq
\sum_{s\in {\mathcal X}_m}Y_m(s)\cdot g(m,s)\cdot N_m\leq \frac{M_m\cdot N_m}
{2^{g(m)(m+1)(k+2)}}\]
and therefore $\sum\limits_{s\in {\mathcal X}_m}\frac{g(m,s)}{M_m}\leq
\frac{1}{2^{g(m)(m+1)}}$. Hence
\[\sum_{s\in {\mathcal X}_m}\frac{\mu\big(\lbv\rho^\ell_h(\dot{\eta}\rest
[m,m+1))=s\rbv_{\random}\cdot b\big)}{\mu(b)}\leq\sum_{s\in {\mathcal X}_m}
\left(\frac{g(m,s)}{M_m}+\frac{1}{M_m}\right)\leq\frac{1}{2^{g(m)m+g(m)-1}}.\]
Let $b_m=\sum_{s\in {\mathcal X}_m}^{\random}\lbv\rho^\ell_h(\dot{\eta}
\rest [m,m+1))=s\rbv_{\random}\cdot b$. By the above estimations we have
$\mu(b_m)\leq \frac{\mu(b)}{2^{m+4}}$ (remember $m\geq m^*>k+4\geq 4$). Look
at the condition $b^*=b-(\sum_{m\geq m^*}^{\random} b_m)$. Clearly
$\mu(b^*)>0$ and
\[b^*\forces_{\random}(\forall m\geq m^*)(\rho^\ell_h(\dot{\eta}\rest [m,m+1)
\notin {\mathcal X}_m))\]
and therefore 
\[b^*\forces_{\random}(\forall m\geq m^*)(s_m\in W^k_{\ell,h}(m,m+1,\dot{\eta}
\rest [m,m+1))).\]
This finishes the proof of the claim and thus checking that the construction
is as required.
\end{proof}

\begin{conclusion}
\label{superconc}
Let $F\in\baire$ be strictly increasing and $f,g$, $(K_{\ref{supersystem}},
\Sigma_{\ref{supersystem}})$, $(K^{\ell}_{\ref{supersystem}},
\Sigma^{\ell}_{\ref{supersystem}})$, $\bar{t}$, $\bar{t}_\ell$, $\bar{\F}_\ell
=(\F_\ell,<^*_\ell)$ (for $\ell<\omega$) and $\langle W^k_{\ell,h}\!: h\in
\F_\ell,k,\ell<\omega\rangle$ be given by \ref{supersystem} (for $F$). Suppose
that $\Gamma\subseteq\p^*_\infty(K_{\ref{supersystem}},\Sigma_{
\ref{supersystem}})$ is quasi-$W^k_{\ell,h}$-generic for every
$k,\ell<\omega$, $h\in\F_\ell$ (exists e.g.~under CH, see \ref{exists}(2),
\ref{exirem}).
\begin{enumerate}
\item Countable support iterations of proper forcing notions which are
$(\Gamma,W^k_{\ell,h})$--genericity preserving for all $k,\ell\in\omega$,
$h\in\F_\ell$ is $(\Gamma,W^k_{\ell,h})$--genericity preserving (for $k,\ell
\in\omega$, $h\in\F$) and hence does not add Cohen reals.
\item Every $(f,g)$-bounding proper forcing notion (this includes proper
forcing notions with the Laver property) and random real forcing are
$(\Gamma,W^k_{\ell,h})$--genericity preserving (for $k,\ell\in\omega$,
$h\in\F_\ell$),  
\item Assume CH. Then any countable support iteration of proper forcing
notions of one of the following types:
\begin{quotation}
\noindent $(f,g)$--bounding, Laver property, random forcing 
\end{quotation}
does not add Cohen reals.
\end{enumerate}
\end{conclusion}

\begin{proof}
 By \ref{supersystem}, \ref{ganoco}(1), \ref{gamcoher}, \ref{fgLavimptF}
and \ref{Gamiter}.
\end{proof}

\begin{example}
\label{tomek2}
We define a strictly increasing function $F\in\baire$ and a creating pair
$(K_{\ref{tomek2}}, \Sigma_{\ref{tomek2}})$ which: captures singletons, is
strongly finitary, reducible, forgetful, simple, essentially $f^F$--big and
$h$--limited for some function $h:\omega\times\omega \longrightarrow\omega$
such that\quad $(\forall^\infty n)(\prod\limits_{m<n}h(m,m)<g^F(n)<f^F(n))$,
\quad where $f^F,g^F$ are given by \ref{supersystem} for $F$.
\end{example}

\begin{proof}[Construction] For $N<\omega$ we define a nice pre-norm $H_N$ on
${\mathcal P}({\mathcal P}(N)\setminus\{N\})$ by: 
\[H_N(A)=\log_2(1+\max\{k<\omega: (\forall x\in [N]^{\textstyle k})(\exists
a\in A)(x \subseteq a)\})\]
(for $A\subseteq{\mathcal P}(N)\setminus\{N\}$). Note that if $B\subseteq
C\subseteq{\mathcal P}(N)\setminus\{N\}$, $H_N(C)>1$ then:
\begin{enumerate}
\item $\max\{k<\omega: (\forall x\in [N]^{\textstyle k})(\exists a\in C)(x
\subseteq a)\}\geq 2$,
\item $H_N(B)\leq H_N(C)$,
\item $\max\{H_N(B),H_N(C\setminus B)\}\geq H_N(C)-1$.
\end{enumerate}
For 3) above note that if it fails then we find $k_0\in\omega$ such that
\[1+\max\{H_N(B),H_N(C\setminus B)\}\leq \log_2(2k_0)<H_N(C).\]
By the first inequality we find $x_0,x_1\in [N]^{\textstyle k_0}$ such that 
\[(\forall a\in B)(x_0\not\subseteq a)\quad\mbox{ and }\quad(\forall a\in
C\setminus B)(x_1\not\subseteq a).\]
But the second inequality implies that there is $a\in C$ such that $x_0\cup
x_1\subseteq a$, what gives a contradiction.\\
As clearly $H_N(\{a\})=0$ for $a\in{\mathcal P}(N)\setminus\{N\}$, $H_N$ is
really a nice pre--norm.
\medskip

\noindent Let $F\in\baire$ be defined by $F(m)=2^{(m+1)^3\cdot 2^{2^m}}$ (for
$m\in\omega$) and let $f^F,g^F$ be from \ref{supersystem} (for $F$). To
simplify notation let $M_n=2^{f^F(n)}$, $N_n=(n+1)^2\cdot 2^{M_n}$ (for
$n\in\omega$).

Let $\bH:\omega\longrightarrow {\mathcal P}(\fsuo)$ be given by $\bH(n)=
[N_n]^{\textstyle 2^{M_n}}$.

Let $K_{\ref{tomek2}}$ consist of creatures $t\in\CR[\bH]$ such that, letting
$m=m^t_{\dn}$, 
\begin{enumerate}
\item[(a)] $\dis[t]$ is a subset of $\bH(m)$,
\item[(b)] $\val[t]=\{\langle w,u\rangle\in\prod\limits_{i<m}\bH(i)\times
\prod\limits_{i\leq m}\bH(i): w\vartriangleleft u\ \&\ u(m)\in\dis[t]\}$,
\item[(c)] $\nor[t]=\frac{H_{N_m}(\dis[t])}{(m+1)\cdot f^F(m)}$.
\end{enumerate}
The composition operation $\Sigma_{\ref{tomek2}}$ is the trivial one: it gives
a nonempty result for singletons only and then $\Sigma_{\ref{tomek2}}(t)=\{s
\in K_{\ref{tomek2}}: m^s_{\dn}=m^t_{\dn}\ \&\ \dis[s]\subseteq\dis[t]\}$.
Now, we have to check that $(K_{\ref{tomek2}},\Sigma_{\ref{tomek2}})$ has the
required properties. Clearly it is a strongly finitary, simple and forgetful
creating pair. Plainly, it captures singletons. Note that if $t_m\in 
K_{\ref{tomek2}}$ (for $m\in\omega$) is such that $\dis[t_m]=\bH(m)$ then
\[\nor[t_m]=\frac{H_{N_m}(\bH(m))}{(m+1)\cdot f^F(m)}=\frac{\log_2(1+2^{M_m})}
{(m+1)\cdot f^F(m)}>\frac{2^{f^F(m)}}{(m+1)\cdot f^F(m)}\stackrel{m\to\infty
}{\longrightarrow}\infty.\]
Consequently the forcing notions $\q^*_\infty(K_{\ref{tomek2}},
\Sigma_{\ref{tomek2}})$, $\q^*_{{\rm w}\infty}(K_{\ref{tomek2}},
\Sigma_{\ref{tomek2}})$ etc will be non-trivial.
 
\noindent To verify that $(K_{\ref{tomek2}},\Sigma_{\ref{tomek2}})$ is
reducible use the fact that $H_N$ is a nice pre-norm on ${\mathcal P}({
\mathcal P}(N)\setminus\{N\})$: if $a\in C\in {\mathcal P}({\mathcal P}(N)
\setminus\{N\})$, $H_N(C)>1$ then $H_N(C)-1\leq H_N(C\setminus\{a\})$. For
similar reasons $(K_{\ref{tomek2}},\Sigma_{\ref{tomek2}})$ is essentially
$f^F$--big, remember that we divide the respective value of $H_{N_m}$ by
$(m+1)\cdot f^F(m)$ (where $m=m^t_{\dn}$). 

\noindent Finally, let $h(m,k)=2^{N_m}$ for $m,k\in\omega$. Then $|\bH(m)|=
{N_m\choose 2^{M_m}}<2^{N_m}$ for all $m\in\omega$, so $(K_{\ref{tomek2}},
\Sigma_{\ref{tomek2}})$ is $h$--limited (actually much more). Moreover, for
every $m\in\omega$:
\[g^F(m\!+\!1)=F(f^F(m))>2^{(m+1)^3\cdot 2^{2^{f^F(m)}}}= 2^{(m+1)^3\cdot
2^{M_m}}>\prod_{k\leq m} 2^{N_k}=\prod_{k\leq m}h(k,k).\]
Thus $(K_{\ref{tomek2}},\Sigma_{\ref{tomek2}})$ is as required. 
\end{proof}

\begin{conclusion}
\label{tomek2conc}
The forcing notion $\q^*_{{\rm w}\infty}(K_{\ref{tomek2}},\Sigma_{
\ref{tomek2}})$:
\begin{enumerate}
\item is proper and $\baire$--bounding,
\item preserves non-meager sets,
\item is $(f^F,g^F)$--bounding,
\item makes the ground model reals have measure zero.
\end{enumerate}
Moreover, assuming CH, countable support iterations of $\q^*_{{\rm w}\infty}
(K_{\ref{tomek2}},\Sigma_{\ref{tomek2}})$ with Laver's forcing notion,
Miller's forcing notion and random forcing do not add Cohen reals. 
\end{conclusion}

\begin{proof} The first required property is a consequence of \ref{sinwin} and
\ref{winbound}. The second follows from \ref{winnonmea} and the third
property is a consequence of \ref{winfg}. The ``moreover'' part holds true
by \ref{superconc}. 

Let ${\mathcal X}=\prod\limits_{m\in\omega} N_m$ be equipped with the product
measure. Of course, this space is measure isomorphic to the reals, so what we
have to show is the following claim. 
\begin{claim}
\label{cl31}
$\forces_{\q^*_{{\rm w}\infty}(K_{\ref{tomek2}},\Sigma_{\ref{tomek2}})}$ 
`` $\V\cap {\mathcal X}$ is a null set''.
\end{claim}

\noindent{\em Proof of the claim:}\ \ \ Let $\dot{W}$ be the $\q^*_{{\rm w}
\infty}(K_{\ref{tomek2}},\Sigma_{\ref{tomek2}})$--name for the generic real
(see \ref{thereal}) and let $\dot{A}$ be a $\q^*_{{\rm w}\infty}(
K_{\ref{tomek2}},\Sigma_{\ref{tomek2}})$--name for a subset of ${\mathcal X}$
such that 
\[\forces_{\q^*_{{\rm w}\infty}(K_{\ref{tomek2}},\Sigma_{\ref{tomek2}})}
\mbox{``}\dot{A}=\{x\in {\mathcal X}: (\exists^\infty n)(x(n)\in \dot{W}(n))\}
\mbox{''}.\]
Note that $\forces$`` $(\forall m\in\omega)(\dot{W}(m)\in [N_m]^{2^{M_m}})$''
and $\frac{2^{M_m}}{N_m}=\frac{1}{(m+1)^2}$. Consequently, $\forces$``
$\dot{A}$ is a null subset of ${\mathcal X}$''. By the definition of
$(K_{\ref{tomek2}},\Sigma_{\ref{tomek2}})$ (remember the definition of $H_N$)
one easily shows that $\forces$`` $\V\cap {\mathcal X}\subseteq\dot{A}$'',
what finishes the proof of the claim and the conclusion.
\end{proof}

One may consider tree versions of \ref{tomek2}.

\begin{example}
\label{tomektree}
Let $F\in\baire$, $h:\omega\times\omega\longrightarrow\omega$ and $\bH:\omega
\longrightarrow{\mathcal P}(\fsuo)$ be as defined in \ref{tomek2}. There are
finitary tree--creating pairs $(K^\ell_{\ref{tomektree}},\Sigma^\ell_{
\ref{tomektree}})$ (for $\ell<3$) for $\bH$ such that
\begin{enumerate}
\item[(a)] $(K^0_{\ref{tomektree}},\Sigma^0_{\ref{tomektree}})$ is local,
reducible, $h$--limited and essentially $f^F$--big (where $f^F$ is given by
\ref{supersystem}),
\item[(b)] $(K^1_{\ref{tomektree}},\Sigma^1_{\ref{tomektree}})$ is t-omittory,
reducible and essentially $f^F$--big,
\item[(c)] $(K^2_{\ref{tomektree}},\Sigma^2_{\ref{tomektree}})$ is reducible,
$h$--limited and essentially $f^F$--big,
\item[(d)] the forcing notions $\q^{\tree}_0(K^0_{\ref{tomektree}},\Sigma^0_{
\ref{tomektree}})$, $\q^{\tree}_1(K^1_{\ref{tomektree}},\Sigma^1_{
\ref{tomektree}})$ and\\
$\q^{\tree}_1(K^2_{\ref{tomektree}},\Sigma^2_{\ref{tomektree}})$ are
equivalent, 
\item[(e)] the forcing notion $\q^{\tree}_1(K^0_{\ref{tomektree}},\Sigma^0_{
\ref{tomektree}})$ is proper, $\baire$-bounding, makes ground model reals
meager and null (even more), is $(f^F,g^+)$--bounding, where $g^+(n)=
\prod\limits_{i<n}|\bH(i)|$,
\item[(f)] the forcing notions $\q^{\tree}_1(K^1_{\ref{tomektree}},\Sigma^1_{
\ref{tomektree}})$ and $\q^{\tree}_1(K^2_{\ref{tomektree}},\Sigma^2_{
\ref{tomektree}})$ are proper, $\baire$--bounding, preserve non-meager sets,
make ground model reals null and are $(f^F,g^F)$--bounding.
\end{enumerate}
\end{example}

\begin{proof}[Construction] Let $F,N_m,M_m,\bH,h$ be as in \ref{tomek2}.

\noindent A tree creature $t\in\TCR_\eta[\bH]$ is in $K^0_{\ref{tomektree}}$
if
\begin{itemize}
\item $\dis[t]\subseteq\bH(\lh(\eta))$,
\item $\val[t]=\{\langle\eta,\nu\rangle: \eta\vartriangleleft\nu\ \&\ \lh(\nu)
=\lh(\eta)+1\ \&\ \nu(\lh(\eta))\in\dis[t]\}$,
\item $\nor[t]=\frac{H_{N_{\lh(\eta)}}(\dis[t])}{(\lh(\eta)+1)\cdot f^F(\lh(
\eta))}$.
\end{itemize}
The tree composition $\Sigma^0_{\ref{tomektree}}$ on $K^0_{\ref{tomektree}}$
is trivial: $\Sigma^0_{\ref{tomektree}}(t)=\{s\in K^0_{\ref{tomektree}}:
\val[s]\subseteq\val[t]\}$. 
\smallskip

\noindent The family $K^1_{\ref{tomektree}}$ consists of these tree--creatures
$t\in\TCR[\bH]$ that for some $\eta\trianglelefteq\eta^*\in\bigcup\limits_{
n<\omega}\prod\limits_{i<n}\bH(i)$ we have
\begin{itemize}
\item $\dis[t]\subseteq\bH(\lh(\eta^*))$,
\item $\val[t]=\{\langle\eta,\nu\rangle: \eta^*\vartriangleleft\nu\ \&\ \lh(
\nu)=\lh(\eta^*)+1\ \&\ \nu(\lh(\eta^*))\in\dis[t]\}$,
\item $\nor[t]=\frac{H_{N_{\lh(\eta^*)}}(\dis[t])}{(\lh(\eta^*)+1)\cdot f^F(
\lh(\eta^*))}$.
\end{itemize}
The tree composition $\Sigma^1_{\ref{tomektree}}$ is such that 
\[\Sigma^1_{\ref{tomektree}}(t_\nu\!:\nu\in\hat{T})=\{t\in
K^1_{\ref{tomektree}}\!: \dom(\val[t])\!=\!\{\mrot(T)\}\ \&\ \rng(\val[t])
\subseteq\max(T)\}.\]
In a similar manner we define $(K^2_{\ref{tomektree}},\Sigma^2_{
\ref{tomektree}})$. A tree creature $t\in\TCR_\eta[\bH]$ is in
$K^2_{\ref{tomektree}}$ if 
\begin{itemize}
\item $\dis[t]=(A_t,\langle\nu^t_x:x\in A_t\rangle)$, where $A_t\subseteq\bH(
\lh(\eta))$ and $\nu^t_x\in\bigcup\limits_{n<\omega}\prod\limits_{i<n}\bH(i)$
(for $x\in A_t$) are such that $\eta\vartriangleleft\eta\conc\langle x\rangle
\trianglelefteq \nu^t_x$,
\item $\val[t]=\{\langle\eta,\nu^t_x\rangle: x\in A_t\}$,
\item $\nor[t]=\frac{H_{N_{\lh(\eta)}}(A_t)}{(\lh(\eta)+1)\cdot f^F(\lh(
\eta))}$,
\end{itemize}
and the tree composition $\Sigma^2_{\ref{tomektree}}$ is defined like
$\Sigma^1_{\ref{tomektree}}$. 

\noindent Checking that $(K^\ell_{\ref{tomektree}},\Sigma^\ell_{
\ref{tomektree}})$ have the desired properties is straightforward and similar
to \ref{tomek2} (remember \ref{confgbound}, \ref{fgbounding},
\ref{tominonmea}, \ref{trmamea}). 
\end{proof}

\begin{corollary}
\label{fgforhal}
Let $F\in\baire$ be an increasing function, $f,\bH$ be as defined in
\ref{bighalex} (for $F$) and let $F^*(n)=F(\fH(n))$. Then the forcing notion
$\q^*_f(K_{\ref{bighalex}},\Sigma_{\ref{bighalex}})$ (defined as in
\ref{bighalex} for $F$) is proper, $\baire$--bounding, $(F^*,\fH)$--bounding
and makes ground model reals meager.
\end{corollary}

\begin{proof} 
By \ref{gfbound}, \ref{FHfg} and \ref{trmamea}(2).
\end{proof}

\chapter{Playing with ultrafilters}
This chapter originated in the following question of Matet and Pawlikowski.

Are the cardinals ${\mathfrak m}_1$, ${\mathfrak m}_2$ equal, where
\begin{enumerate}
\item[${\mathfrak m}_1$\ \ ] is the least cardinality of a set $Z\subseteq
\bigcup\limits_{E\subseteq\omega} \omega^E$ such that 
\begin{enumerate}
\item[(i)]   $(\exists F\in\baire)(\forall f\in Z)(\forall n\in\dom(f))(f(n)<
F(n))$,
\item[(ii)]  the family $\{\dom(f): f\in Z\}$ has the finite intersection
property, and
\item[(iii)] $(\forall x\in\prod\limits_{n\in\omega} [\omega]^{\textstyle \leq
n})(\exists f\in Z)(\forall^\infty n\in\dom(f))(f(n)\notin x(n))$;
\end{enumerate}
\item[${\mathfrak m}_2$\ \ ] is defined in a similar manner, but {\bf (iii)}
is replaced by
\begin{enumerate}
\item[(iii)$^-$] $(\forall g\in\baire)(\exists f\in Z)(\forall^\infty n\in
\dom(f))(f(n)\neq g(n))$?
\end{enumerate}
\end{enumerate}
It was known that ${\mathfrak m}_2\leq\lambda\leq{\mathfrak m}_1$, where
$\lambda$ is the least size of a basis of an ideal on $\omega$ which is not a
weak $q$-point (see Matet Pawlikowski \cite{MaPa9x}). We answer the Matet --
Pawlikowski question in \ref{conmatet}, showing that it is consistent that
$\lambda=\aleph_1$ (and so ${\mathfrak m}_2=\aleph_1$) and ${\mathfrak m}_1=
\aleph_2$. On the way to this result we have to deal with preserving some
special ultrafilters on $\omega$. The technology developed in the previous
section is very useful for this (both to describe the required properties and
to preserve them at limit stages of countable support iterations).

In the first part of the chapter we present the framework: ultrafilters
generated by quasi-$W$-generic $\Gamma$. Then we introduce several properties
of ultrafilters and discuss relations between them. The third section shows
how forcing notions constructed according to our schema may preserve some
special ultrafilters. Finally, in the last part, we apply all these tools to
answer the Matet -- Pawlikowski question. 

\section{Generating an ultrafilter}

\begin{definition}
\label{genult}
We say that a creating pair $(K,\Sigma)$ {\em generates an ultrafilter} if 
\begin{enumerate}
\item[$(\otimes_{\ref{genult}})$] for every $k<\omega$ there is $k^*<\omega$
such that 

{\em if} $t\in K$, $\nor[t]\geq k^*$ and $c:[m^t_{\dn},m^t_{\up})
\longrightarrow 2$, 

{\em then} for some $s\in\Sigma(t)$ and $i^*<2$ we have $\nor[s]\geq k$ and
\[[u\in\basis(s)\ \&\ v\in\pos(u,s)\ \&\ m^t_{\dn}\leq m<m^t_{\up}\ \&\ v(m)
\neq 0]\ \ \Rightarrow\ \ c(m)=i^*\]   
\end{enumerate}
and if $\nor[t]>0$, $u\in\basis(t)$ then there is $v\in\pos(u,t)$ such that
for some $m\in [m^t_{\dn},m^t_{\up})$ we have $v(m)\neq 0$.
\end{definition}

\begin{proposition}
\label{ultrafilter}
Assume that $(K,\Sigma)$ is an omittory and monotonic (see \ref{monotonic})
creating pair which generates an ultrafilter. Then 
\[\forces_{\q^*_{{\rm s}\infty}(K,\Sigma)}\mbox{``}\{m:\Wtil(m)\neq 0\}\mbox{
induces an ultrafilter on the algebra } \cP(\omega)^V/\Fin\mbox{''}.\]
\end{proposition}

\begin{proof}
For $k\in\omega$, let $g(k)$ be the $k^*$ given by $(\otimes_{\ref{genult}})$
(for $k$). Let $\dot{Z}$ be a $\q^*_{{\rm s}\infty} (K,\Sigma)$-name such that 
\[\forces_{\q^*_{{\rm s}\infty}(K,\Sigma)}\dot{Z}=\{m\in\omega:\dot{W}(m)\neq
0\}.\] 
Clearly $\forces_{\q^*_{{\rm s}\infty}(K,\Sigma)}\dot{Z}\in\iso$ (by the
second requirement of \ref{genult}). We have to show that for each $A\in\iso$ 
\[\forces_{\q^*_{{\rm s}\infty}(K,\Sigma)}\mbox{``either }\dot{Z}\cap A\mbox{
or } \dot{Z}\setminus A\mbox{ is finite''}.\]
So suppose that $p\in\q^*_{{\rm s}\infty}(K,\Sigma)$, $A\in\iso$ are such that 
\[p\forces_{\q^*_{{\rm s}\infty}(K,\Sigma)}\mbox{``both }\dot{Z}\cap A\mbox{
and } \dot{Z}\setminus A\mbox{ are infinite''}.\]
Since $(K,\Sigma)$ is omittory we may assume that
$p=(w^p,t_0,t_1,\ldots)$ where $\nor[t_\ell]\geq g(\ell+ m^{t_\ell}_{\dn})$
(see \ref{sinfty}.(2)). Let $c_\ell:[m^{t_\ell}_{\dn},m^{t_\ell}_{\up})
\longrightarrow 2$ be the characteristic function of $A\cap [m^{t_\ell}_{\dn},
m^{t_\ell}_{\up})$. Applying the condition $(\otimes_{\ref{genult}})$ and the
choice of $g$ for each $\ell<\omega$ we find $s^\prime_\ell\in\Sigma(t_\ell)$
such that $\nor[s^\prime_\ell]\geq\ell+m^{s^\prime_\ell}_{\dn}$ and for some
$i^\ell<2$ we have  
\begin{enumerate}
\item[$(\otimes^*)$] {\em if\/} $u\in\basis(s^\prime_\ell)$,
$v\in\pos(u,s^\prime_\ell)$, $m^{s^\prime_\ell}_{\dn}\leq m<
m^{s^\prime_\ell}_{\up}$ and $v(m)\neq 0$

\noindent{\em then \/} $c_\ell(m)=i^\ell$.
\end{enumerate}
Now choose $i^*<2$ and an increasing sequence $\langle \ell_k:k<\omega\rangle
\subseteq \omega$ such that $i^{\ell_k}=i^*$ (for $k<\omega$) and take $s_0=
s_{\ell_0}^\prime\Rsh [m^{s_0^\prime}_{\dn}, m^{s_{\ell_0}^\prime}_{\up})$,
$s_{k+1}= s^\prime_{\ell_{k+1}}\Rsh [m^{s^\prime_{\ell_k}}_{\up},
m^{s^\prime_{\ell_{k+1}}}_{\up})$. Once again, since $(K,\Sigma)$ is omittory
we get $q\stackrel{\rm def}{=}(w^p,s_0,s_1,s_2,\ldots)\in\q^*_{{\rm
s}\infty}(K,\Sigma)$ and it is stronger than $p$. Suppose $i^*=0$. We claim
that in this case $q\forces_{\q^*_{{\rm s}\infty}(K,\Sigma)}\dot{Z}\cap
A\subseteq m^{s_0}_{\dn}$ (contradicting the choice of $p$, $A$). If not 
then we find $k<\omega$, $w\in\pos(w^p,s_0,\ldots,s_k)$ and $n\in [m^{s_0}_{
\dn},m^{s_k}_{\up}) \cap A$ such that $w(n)\neq 0$. By smoothness we may
additionally demand that, if $k>0$ then $n\in[m^{s_k}_{\dn},m^{s_k}_{\up}) 
=[m^{s^\prime_{\ell_{k-1}}}_{\up},m^{s^\prime_{\ell_k}}_{\up})$. By the
smoothness and monotonicity of $(K,\Sigma)$ we have 
\[w\in\pos(w\rest m^{s_k}_{\dn},s_k)=\{v:\langle w\rest m^{s_k}_{\dn},v\rangle
\in\val[s_k]\}.\] 
Now, by the choice of $s_k$ we may conclude that 
\[w\rest [m^{s_k}_{\dn}, m^{s_{\ell_k}^\prime}_{\dn})={\bf 0}_{[m^{s_k}_{\dn},
m^{s_{\ell_k}^\prime}_{\dn})}\quad\mbox{ and thus }\  
n\in [m^{s_{\ell_k}^\prime}_{\dn}, m^{s_{\ell_k}^\prime}_{\up}).\]
Moreover $\langle w\rest m^{s^\prime_{\ell_k}}_{\dn},w\rangle\in
\val[s^\prime_{\ell_k}]$ so we may apply $(\otimes^*)$ (to $\ell_k$) and
conclude that $c_{\ell_k}(n)=i^*=0$, a contradiction. Similarly one shows
$q\forces\dot{Z}\setminus A\subseteq m^{s_0}_{\dn}$ if $i^*=1$.
\end{proof}

\begin{definition}
\label{Gamfilter}
Suppose that $(K,\Sigma)$ is a creating pair, $\bar{t}\in\PC_{\C(\nor)}(K,
\Sigma)$ and $W:\omega\times\omega\times\pfs\longrightarrow {\mathcal P}(K)$
is a $\bar{t}$--system (see \ref{quasi}). Let $\Gamma\subseteq\p^*_{\C(\nor)}
(K,\Sigma)$ be quasi-$W$-generic.
\begin{enumerate}
\item We define $\D(\Gamma)$ as the family of all sets $A\subseteq\omega$ such
that:\quad for some $\bar{s}=\langle s_n: n<\omega\rangle\in\Gamma$ and
$N<\omega$, for every $w\in \basis(s_0)$ and $u\in\pos(w,s_0,\ldots,s_N)$ we
have 
\[(\forall m>N)(\forall v\in\pos(u,s_{N+1},\ldots,s_m))(\{k\in [\lh(u),\lh(v))
\!: v(k)\neq 0\}\subseteq A).\]
\item We say that $\Gamma$ {\em generates a filter (an ultrafilter,
respectively) on $\omega$} if $\D(\Gamma)$ is a filter (an ultrafilter,
resp.). 
\end{enumerate}
\end{definition}

\begin{remark}
Note the close relation of \ref{Gamfilter} and \ref{ultrafilter}. Below it
becomes even closer. 
\end{remark}

\begin{definition}
\label{intcrea}
A creating pair $(K,\Sigma)$ is {\em interesting} if for each creature $t\in
K$ such that $\nor[t]>0$ and every $u\in\basis(t)$ we have
\[|\{m\in [m^t_{\dn},m^t_{\up}): (\exists v\in\pos(u,t))(v(m)\neq 0)\}|>1.\]
\end{definition}

\begin{proposition}
\label{facgamfil}
Let $(K,\Sigma)$ be a forgetful creating pair, $\bar{t}\in\PC_\infty(K,
\Sigma)$ and $W:\omega\times\omega\times\pfs\longrightarrow {\mathcal P}(K)$
be a $\bar{t}$--system.  
\begin{enumerate}
\item If $\Gamma\subseteq\p^*_\infty(\bar{t},(K,\Sigma))$ is quasi-$W$-generic
then $\D(\Gamma)$ is a filter on $\omega$ containing all co-finite sets. If,
additionally, $(K,\Sigma)$ is interesting and condensed (see \ref{ab}(3)) and
$A\subseteq\omega$ is such that
\[(\forall^\infty n)(|A\cap [m^{t_n}_{\dn},m^{t_n}_{\up})|\leq 1)\]
then $A\notin\D(\Gamma)$.
\item If $\Gamma$ is quasi-$W$-generic in $\p^*_\infty(\bar{t},(K,\Sigma))$
then $\D(\Gamma)$ is a $p$-point (see \ref{ramppoint}). 
\item Assume CH. Suppose that, additionally, $(K,\Sigma)$ is finitary,
omittory, monotonic and generates an ultrafilter. Then there exists a
quasi-$W$-generic $\Gamma\subseteq\p^*_\infty(K,\Sigma)$ such that $\D(
\Gamma)$ is an ultrafilter on $\omega$. 
\end{enumerate}
\end{proposition}

\begin{proof} 
1), 2) Should be obvious. 

\noindent 3) Modify the proof of \ref{exists}(2), noting that if $\bar{s}\in
\PC_\infty(K,\Sigma)$, $A\subseteq\omega$ then there is $\bar{s}^*=\langle
s^*_n: n<\omega\rangle\in\PC_\infty(K,\Sigma)$ such that $\bar{s}\leq
\bar{s^*}$ and either
\[(\forall u\!\in\!\basis(s_0^*))(\forall n\!<\!\omega)(\forall v\!\in\!
\pos(u,s_0^*,\ldots,s^*_n))(\{k\!\in\![\lh(u),\lh(v))\!: v(k)\!\neq\! 0\}
\subseteq A)\]
or a similar requirement with $\omega\setminus A$ instead of $A$ holds.
(Compare the proof of \ref{ultrafilter}, remember $(K,\Sigma)$ is forgetful.)
\end{proof}

\begin{conclusion}
\label{presextra}
Suppose that $(K,\Sigma)$ is a forgetful and monotonic creating pair,
$\bar{t}\in\PC_\infty(K,\Sigma)$ and $W$ is a $\bar{t}$--system. Assume that
$\Gamma\subseteq \p^*_\infty(\bar{t},(K,\Sigma))$ is quasi-$W$-generic and
$\D(\Gamma)$ is an ultrafilter. Let $\delta$ be a limit ordinal and $\langle
\p_\alpha,\dot{\q}_\alpha: \alpha<\delta\rangle$ be a countable support
iteration of proper forcing notions such that for each $\alpha<\delta$:
\begin{enumerate}
\item[(a)] $\forces_{\p_\alpha}$``$\Gamma$ is quasi-$W$-generic'',
\item[(b)] $\forces_{\p_\alpha}$``$\D(\Gamma)$ is an ultrafilter''.
\end{enumerate}
Then
\begin{enumerate}
\item $\forces_{\p_\delta}$``$\Gamma$ is quasi-$W$-generic'',
\item $\forces_{\p_\delta}$``$\D(\Gamma)$ is an ultrafilter''.
\end{enumerate}
\end{conclusion}

\begin{proof}
1)\ \ It follows from \ref{Gamiter}.

\noindent 2)\ \ Since, by \ref{facgamfil}(2), $\D(\Gamma)$ is a $p$-point we
may use \cite[Ch VI, 5.2]{Sh:f} (another presentation of this result might be
found in \cite[6.2]{BaJu95}).
\end{proof}

\begin{remark}
\begin{enumerate}
\item If $\bar{t}=\langle t_n:n<\omega\rangle\in\PC_\infty(K,\Sigma)$, $W$ is
a $\bar{t}$--system  and $\Gamma$ is quasi-$W$-generic generating a filter
then we make think of $\D(\Gamma)$ as a filter on $\bigcup\limits_{i\in\omega}
\{i\}\times (m^{t_i}_{\up}-m^{t_i}_{\dn})$ (just putting the intervals
$[m^{t_i}_{\dn}, m^{t_i}_{\up})$ vertically). This will be our approach in the
further part, where we will consider ultrafilters on $\bigcup\limits_{i\in
\omega}\{i\}\times (i+1)$. 
\item We may treat $\D(\Gamma)$ as a canonical filter on $\omega$ with a
property described by $\Gamma$ (or, more accurately, by the $\bar{t}$--system
$W$). This is the way we are going to use \ref{presextra} later: it will allow
us to claim that the additional property of an ultrafilter is preserved at
limit stages of an iteration.
\end{enumerate}
\end{remark}

\section{Between Ramsey and p-points}
Here we recall some definitions of special properties of ultrafilters on
$\omega$ and we introduce more of them. Then we comment on relations between
these notions.

\begin{definition}
\label{ramppoint}
Let $\D$ be a filter on $\omega$. We say that:
\begin{enumerate}
\item {\em $\D$ is Ramsey} if for each colouring
$F:[\omega]^{\textstyle 2}\longrightarrow 2$ there is a set $A\in\D$
homogeneous for $F$.
\item {\em $\D$ is a $p$-point} if for every partition $\langle A_n:
n\in\omega\rangle$ of $\omega$ into sets from the dual ideal (i.e.
$\omega\setminus A_n\in \D$) we find a set $A\in \D$ with
\[(\forall n\in\omega)(|A_n\cap A|<\omega).\]
\item {\em $\D$ is a $q$-point} if for every partition $\langle A_n:
n\in\omega\rangle$ of $\omega$ into finite sets there is a set $A\in\D$
with
\[(\forall n\in\omega)(|A_n\cap A|\leq 1).\]
\item {\em $\D$ is a weak $q$-point} if for each set $B\subseteq\omega$ such
that $\omega\setminus B\notin\D$ and a partition $\langle A_n:n\in\omega
\rangle$ of $B$ into finite sets there is a set $A\subseteq B$ such that
\[\omega\setminus A\notin\D\qquad\mbox{ and}\qquad(\forall n\in\omega)(|A_n
\cap A|\leq 1).\]
\end{enumerate}
\end{definition}

\begin{remark}
Clearly, if $\D$ is an ultrafilter on $\omega$ which is a weak $q$-point then
$\D$ is a $q$-point. (So the two notions coincide for ultrafilters).
\end{remark}

\begin{definition}
\begin{enumerate}
\item For a filter $\D$ on $\omega$ let $G^R(\D)$ be the game of two players,
I and II, in which Player I in his $n^{\rm th}$ move plays a set $A_n\in\D$
and Player II answers choosing a point $k_n\in A_n$. Thus a result of a play
is a pair of sequences $\langle\langle A_n: n\in\omega\rangle,\langle k_n:
n\in\omega\rangle\rangle$ such that $k_n\in A_n\in\D$.\\ 
Player I wins the play of the game $G^R(\D)$ if and only if\\ 
the result $\langle\langle A_n:n\in\omega\rangle,\langle k_n: n\in\omega
\rangle\rangle$ satisfies:\ \ \ $\{k_n:n\in\omega\}\notin\D$.

\item Similarly we define the game $G^p(\D)$ allowing the second player to
play finite sets $a_n\subseteq A_n$ (instead of points $k_n\in A_n$).\\
Player I wins if $\bigcup\limits_{n\in\omega} a_n\notin \D$
\end{enumerate}
\end{definition}

\begin{remark}
Let us recall that if $\D$ is an ultrafilter on $\omega$ then the
following conditions are equivalent (see \cite[Ch VI, 5.6]{Sh:f} or
\cite[4.5]{BaJu95}): 
\begin{enumerate}
\item[(a)] $\D$ is Ramsey,
\item[(b)] $\D$ is both a $p$-point and a $q$-point,
\item[(c)] Player I does not have a winning strategy in the game $G^R(\D)$.
\end{enumerate}
Similarly, an ultrafilter $\D$ is a $p$-point if and only if Player I does not
have a winning strategy in $G^p(\D)$.
\end{remark}
As we are interested in ultrafilters which are not $q$-points (see the
discussion of the Matet --- Pawlikowski problem at the beginning of this
chapter) it is natural to fix a partition of $\omega$ which witnesses this. 
Thus, after renaming, we may consider ultrafilters on $\bigcup\limits_{i\in
\omega}\{i\}\times (i+1)$ instead (compare with the last remark of the
previous section). 

\begin{definition}
\label{almramsey}
Let $\D$ be a filter on $\bigcup\limits_{i<\omega} \{i\}\times (i+1)$.
\begin{enumerate}
\item We say that the filter $\D$ is {\em interesting} if for each function
$h\in\prod\limits_{i\in\omega}(i+1)$ the set $\{(i,h(i)): i\in\omega\}$ is not
in $\D$. 
\item Let $G^{sR}(\D)$ be the game of two players in which Player I in
his $n^{th}$ move plays a set $A_n\in\D$ and the second player answers
choosing an integer $i_n$ and a set $a_n\in [A_n\cap\{i_n\}\times (i_n+1)]^{
\leq n}$.\\
Finally, Player I wins the play if $\bigcup\limits_{n\in\omega}a_n\notin\D$.

\item The game $G^{aR}(\D)$ is a modification of $G^{sR}(\D)$ such that now,
the first Player (in his $n^{th}$ move) chooses a set $A_n\in\D$, $L_n<\omega$
and a function 
\[f_n:\bigcup_{i<\omega}[\{i\}\times(i+1)]^{\textstyle {\leq n}}
\longrightarrow L_n.\]
Player II answers playing $i_n\in\omega$ and a set $a_n\in [A_n\cap\{i_n\}
\times(i_n+1)]^{\leq n}$ homogeneous for $f_n$ (i.e.~such that $f_n\rest
[a_n]^{\textstyle k}$ is constant for $k\leq n$).\\
Player I wins the play if $\bigcup\limits_{n\in\omega}a_n\notin\D$.
\item We say that the filter $\D$ is {\em semi--Ramsey} if the first
player has no winning strategy in the game $G^{sR}(\D)$.
\item The filter $\D$ is {\em almost Ramsey} if it is semi--Ramsey and for
every colouring 
\[f:\bigcup_{i<\omega}[\{i\}\times(i+1)]^{\textstyle {\leq n}}\longrightarrow
L,\qquad (n,L<\omega)\] 
there is a set $A\in \D$ which is almost homogeneous for $f$ in the following
sense:
\[(\forall i\in\omega)(f\rest [A\cap(\{i\}\times (i+1))]^{\textstyle k}\quad
\mbox{ is constant for each }k\leq n).\]
\end{enumerate}
\end{definition}

\begin{proposition}
\label{ultimp}
Suppose $\D$ is a non-principal ultrafilter on $\bigcup\limits_{i<\omega}\{i
\}\times (i+1)$. 
\begin{enumerate}
\item If $\D$ is interesting then it is not a $q$-point.
\item If $\D$ is Ramsey then $\D$ is almost Ramsey. 
\item If $\D$ is semi--Ramsey then it is a $p$-point.
\item If $\D$ is semi--Ramsey then
\[\D^*\stackrel{\rm def}{=}\{A\subseteq\omega: \bigcup\limits_{i\in A} \{i\}
\times (i+1)\in\D\}\]
is a Ramsey ultrafilter.
\end{enumerate}
\end{proposition}

\begin{proof}
Compare the games and definitions.
\end{proof}

\begin{theorem}
\label{semnotalm}
Assume CH. There exists an ultrafilter $\D$ on $\bigcup\limits_{i<\omega}\{i\}
\times(i+1)$ which is semi--Ramsey but not almost Ramsey.
\end{theorem}

\begin{proof} 
For $i\in\omega$ let $k_i$ be the integer part of the square root of $i+1$. 
Choose partitions $\langle e^m_i: m<k_i\rangle$ of $\{i\}\times (i+1)$ such
that $(\forall m<k_i)(|e^m_i|\geq k_i)$. Let $f:\bigcup\limits_{i\in\omega}
[\{i\}\times (i+1)]^{\textstyle 2}\longrightarrow 2$ be such that 
\[f((i,\ell_0),(i,\ell_1))=1\quad\mbox{ if and only if }\quad (\forall m<k_i)
((i,\ell_0)\in e^m_i\ \Leftrightarrow\ (i,\ell_1)\in e^m_i).\]
Assuming CH, we will construct a semi--Ramsey ultrafilter containing no almost
homogeneous set for $f$. To this end we choose an enumeration
$\{\varphi_\alpha: \alpha<\omega_1\}$ of all functions from $\omega$ to
$[\bigcup\limits_{i\in\omega}\{i\}\times (i+1)]^{\textstyle \omega}$. By
induction on $\alpha<\omega_1$ define sequences $\langle i^\alpha_n: n<\omega
\rangle$ and $\langle a^\alpha_n: n<\omega\rangle$ such that for $\alpha<\beta
<\omega_1$:
\begin{enumerate}
\item[(a)] $i^\alpha_0<i^\alpha_1<i^\alpha_2<\ldots<\omega$,
\item[(b)] $a^\alpha_n\in [\{i^\alpha_n\}\times (i^\alpha_n+1)]^{\textstyle
\leq n}$ (for $n\in\omega$),
\item[(c)] either $(\exists m\in\omega)(\forall n\in\omega)(a^\alpha_n\cap
\varphi_\alpha(m)=\emptyset)$ or $(\forall n\in\omega)(a^\alpha_{n+1}\subseteq
\varphi_\alpha(i^\alpha_n))$, 
\item[(d)] $\lim\limits_{n\to\infty}|\{m<k_{i^\alpha_n}: e^m_{i^\alpha_n}\cap
a^\alpha_n\neq\emptyset\}|=\infty$,\\ 
$\lim\limits_{n\to\infty}\min\{|a^\alpha_n\cap e^m_{i^\alpha_n}|:
m<k_{i^\alpha_n}\ \&\ a^\alpha_n\cap e^m_{i^\alpha_n}\neq \emptyset\}=\infty$,
\item[(e)] $|\bigcup\limits_{n\in\omega} a^\beta_n\setminus \bigcup\limits_{
n\in\omega} a^\alpha_n|<\omega$.
\end{enumerate}
There should be no problems with carrying out the construction. Let $A_\alpha
=\bigcup\limits_{n\in\omega} a^\alpha_n$. Clearly the sequence $\langle
A_\alpha: \alpha<\omega_1\rangle$ generates a non-principal ultrafilter $\D$
on $\bigcup\limits_{i\in\omega}\{i\}\times (i+1)$ (remember that the constant
functions are among the $\varphi_\alpha$'s). By the demand (d), no set
$A_\alpha$ is almost homogeneous for $f$, so $\D$ is not almost Ramsey. To 
show that $\D$ is semi--Ramsey suppose that $\sigma$ is a winning strategy for
the first player in the game $G^{sR}(\D)$. Then $\sigma$ is a function defined
on finite sequences $\bar{x}=\langle(i_0,a_0),\ldots,(i_{n-1},a_{n-1})\rangle$
such that $i_0<\ldots<i_{n_1}$ and 
\[(\forall \ell<n)(a_\ell\in [\{i_\ell\}\times(i_\ell+1)]^{\textstyle
\leq\ell})\]
and with values in $\D$. For $j<\omega$ put
\[\begin{array}{ll}
\varphi(j)=\bigcap\big\{\sigma((i_0,a_0),\ldots,(i_{n-1},a_{n-1})):& i_0<
\ldots<i_{n-1}\leq j\ \mbox{ and}\\ 
\ &a_\ell\in [\{i_\ell\}\times (i_\ell+1)]^{\textstyle \leq\ell}\big\}.
  \end{array}\]
Thus $\varphi:\omega\longrightarrow\D$, so for some $\alpha<\omega_1$ we have
$\varphi=\varphi_\alpha$. But now look at the sequence $\bar{a}=\langle
(i^\alpha_0,a^\alpha_0),(i^\alpha_1,a^\alpha_1),(i^\alpha_2,a^\alpha_2),\ldots
\rangle$. Since $A_\alpha\in\D$ and $\varphi_\alpha(m)\in\D$ for all $m\in
\omega$ we necessarily have
\[(\forall n\in\omega)(a^\alpha_{n+1}\subseteq\varphi_\alpha(i^\alpha_n)
\subseteq\sigma((i^\alpha_0,a^\alpha_0),\ldots,(i^\alpha_n,a^\alpha_n))).\]
This means that the sequence $\bar{a}$ is a result of a legal play of the
second player against the strategy $\sigma$. Hence $A_\alpha=\bigcup\limits_{
n\in\omega} a^\alpha_n\notin\D$, a contradiction.
\end{proof}

\begin{theorem}
Suppose that $\D$ is a semi--Ramsey ultrafilter on $\bigcup\limits_{i<\omega}
\{i\}\times(i+1)$. Then the following conditions are equivalent.
\begin{enumerate}
\item[(a)] $\D$ is almost Ramsey,
\item[(b)] the first player has no winning strategy in the game $G^{aR}(\D)$,
\item[(c)] for each $m,L\in\omega$ and a colouring $f:\bigcup\limits_{i<
\omega}[\{i\}\times(i+1)]^{\textstyle {\leq m}}\longrightarrow L$, the first
player has no winning strategy in the following modification $G^{sR}_f(\D)$ of
the game $G^{sR}(\D)$: rules are like in $G^{sR}(\D)$ but the sets $a_n$
chosen by the second player have to be homogeneous for $f$.
\end{enumerate}
\end{theorem}

\begin{proof}
The implications (b)\ $\Rightarrow$\ (c)\ $\Rightarrow$\ (a) are immediate by
the definitions. 

The implication (a)\ $\Rightarrow$\ (b) is easy too: suppose that $\sigma$ is
a strategy for the first player in the game $G^{aR}(\D)$. Let $\sigma^*$ be a
strategy for Player I in $G^{sR}(\D)$ such that if $\sigma((i_0,a_0),\ldots,
(i_{n-1},a_{n-1}))=(f_n,A_n)$ then $\sigma^*((i_0,a_0),\ldots,(i_{n-1},a_{n-
1}))\in\D$ is an almost $f_n$--homogeneous subset of $A_n$ (exists by the
assumption {\bf (a)}). Now, $\sigma^*$ cannot be the winning strategy for
Player I as $\D$ is semi--Ramsey. But then the play witnessing this shows that
$\sigma$ is not winning in $G^{aR}(\D)$.
\end{proof}

\begin{theorem}
\label{almboupre}
Assume that $\D$ is a semi--Ramsey ultrafilter on $\bigcup\limits_{i<\omega}
\{i\}\times(i+1)$. Suppose that $\p$ is a proper $\baire$--bounding forcing
notion such that 
\[\forces_{\p}\mbox{`` }\D\mbox{ generates an ultrafilter ''}.\]
Then
\[\forces_{\p}\mbox{`` }\D\mbox{ generates a semi--Ramsey ultrafilter ''}.\]
\end{theorem}

\begin{proof}
This is very similar to \cite[Ch VI, 5.1]{Sh:f}. We know that, in $\V^{\p}$,
$\D$ generates an ultrafilter. What we have to show is that Player I has no
winning strategy in the game $G^{sR}(\D)$ (in $\V^{\p}$). So suppose that
$\sigma\in\V^{\p}$ is a winning strategy of the first player in $G^{sR}(\D)$.
We may assume that the values of $\sigma$ are elements of $\D$ (so from the
ground model). But now, as $\p$ is proper and $\baire$--bounding, we find a
function $\sigma^+\in\V$ such that $\dom(\sigma^+)=\dom(\sigma)$, $\rng(
\sigma^+)\subseteq [\D]^{\textstyle <\omega}$ and $\sigma(\bar{a})\in\sigma^+
(\bar{a})$ for all $\bar{a}\in\dom(\sigma)$. Letting $\sigma^*(\bar{a})=
\bigcap\sigma^+(\bar{a})$ for $\bar{a}\in\dom(\sigma)$ we will get, in $\V$, a
winning strategy for Player I in $G^{sR}(\D)$, a contradiction.
\end{proof}

\begin{remark}
One may note that we did not mention anything about the existence of almost
Ramsey ultrafilters. Of course it is done like \ref{semnotalm}, under
CH. However we want to have an explicit representation of the ultrafilter as
$\D(\Gamma)$ for some quasi-generic $\Gamma$. This will give us the
preservation  of the ``colouring'' part of the definition of almost Ramsey
ultrafilters at limit stages. As the representation is very specific we
postpone it for a moment and we will present this in Examples (see
\ref{represent}, \ref{repexists}).
\end{remark}

\section{Preserving ultrafilters}
In this section we show when forcing notions of the type $\q^{\tree}_1$
preserve ultrafilters introduced in the previous part. The key property of a
tree--creating pair needed for this is formulated in the following definition.

\begin{definition}
\label{up}
Let $\D$ be a filter on $\omega$. We say that a tree creating pair
{\em $(K,\Sigma)$ is of the $\UP(\D)^{\tree}$--type} if the following
condition is satisfied:
\begin{enumerate}
\item[$(\circledast)^{\tree}_{\UP(\D)}$]\hspace{0.15in} Assume that $1\leq
m<\omega$, $p\in\q^{\tree}_\emptyset(K,\Sigma)$, $\nor[t^p_\nu]>m+1$ for each
$\nu\in T^p$, and $F_0, F_1,\ldots$ are fronts of $T^p$ such that 
\[(\forall n\in\omega)(\forall\nu\in F_{n+1})(\exists\eta\in F_n)(\eta
\vartriangleleft\nu).\]
Further suppose that $u_n\subseteq F_n$ (for $n\in\omega$) are sets such that
there is {\bf no} system $\langle s_\nu:\nu\in\hat{T}\rangle\subseteq K$ with:
\begin{enumerate}
\item $T\subseteq\{\nu\in T^p: (\exists\eta\in F_n)(\nu\trianglelefteq\eta)\}$
is a (well founded) quasi tree with $\max(T)\subseteq u_n$ and $\mrot(T)=\mrot
(T^p)$, 
\item for each $\nu\in\hat{T}$:

$\mrot(s_\nu)=\nu$, $\pos(s_\nu)=\suc_T(\nu)$, $\nor[s_\nu]\geq m$ and

$s_\nu\in\Sigma(t^p_\eta:\eta\in\hat{S}_\nu)$ for some (well founded) quasi
tree $S_\nu\subseteq T^p$.
\end{enumerate}
{\em Then} there is a condition $q\in\q^{\tree}_\emptyset(K,\Sigma)$ such that
\begin{enumerate}
\item[$(\alpha)$] $p\leq q$, 
\item[$(\beta)$]  the set 
\[Z\stackrel{\rm def}{=}\{n\in\omega:u_n\cap\dcl(T^q)=\emptyset\}\]
is not in the ideal $\D^c$ dual to $\D$, 
\item[$(\gamma)$] $(\forall\nu\in T^q)(\nor[t^q_\nu]\geq\min\{\nor[t^p_{\eta}]
-m\!: \nu\trianglelefteq\eta\in T^p\})$.
\end{enumerate}
\end{enumerate}
If we may additionally demand that the condition $q\in\q^{\tree}_\emptyset(K,
\Sigma)$ above satisfies
\begin{enumerate}
\item[$(\delta)$] $\mrot(q)=\mrot(p)$ and $(\forall n\in Z)(F_n\cap T^q$
is a front of $T^q)$
\end{enumerate}
then we say that $(K,\Sigma)$ is of the $\sUP(\D)^{\tree}$--type.

\noindent If $\D$ is the filter of all co-finite subsets of $\omega$ then we
say that $(K,\Sigma)$ is of the $\UP^{\tree}$--type ($\sUP^{\tree}$--type,
respectively) instead of $\UP(\D)^{\tree}$--type ($\sUP(\D)^{\tree}$--type,
resp.). 
\end{definition}

\begin{theorem}
\label{presramsey}
Let $\D$ be a Ramsey ultrafilter on $\omega$. Suppose that $(K,\Sigma)$ is a
finitary $2$--big tree--creating pair of the $\UP(\D)^{\tree}$--type. Then:
\[\forces_{\q^{\tree}_1(K,\Sigma)}\mbox{`` $\D$ generates an ultrafilter on
$\omega$ ''.}\]
Consequently, $\forces_{\q^{\tree}_1(K,\Sigma)}$`` $\D$ generates a Ramsey
ultrafilter on $\omega$ ''. 
\end{theorem}

\begin{proof}
Let $\dot{X}$ be a $\q^{\tree}_1(K,\Sigma)$--name for a subset of $\omega$,
$p_0\in\q^{\tree}_1(K,\Sigma)$. Consider the following strategy for Player I
in the game $G^R(\D)$: 
\begin{quotation}
\noindent in the $n^{\rm th}$ move he chooses a condition $p_{n+1}\in
\q^{\tree}_1(K,\Sigma)$ such that $p_n\leq_n^1 p_{n+1}$ and $p_{n+1}\forces
k_n\in\dot{X}$ (where $k_n$ is the last point played so far by Player II). 
Then he plays the set
\end{quotation}
\[B(p_{n+1},n+1)\stackrel{\rm def}{=}\{k\in\omega: (\exists q\in \q^{\tree}_1
(K,\Sigma))(p_{n+1}\leq^1_{n+1} q\ \ \&\ \ q\forces k\in\dot{X})\}.\]  
As $\D$ is Ramsey, this strategy cannot be the winning one for Player I and
therefore there is a play (determined by $k_0,k_1,k_2,\ldots$) according to
this strategy in which Player I looses. This means that one of the following
two possibilities holds:
\medskip

\noindent{\sc Case A}\hspace{0.15in}{\em In the course of the play all sets
$B(p_{n+1},n+1)$ are in the ultrafilter $\D$ and $\{k_0,k_1,\ldots\}\in\D$.}

\noindent In this situation we look at the sequence $\langle p_n:
n\in\omega\rangle$. By \ref{fusAxA} it has the limit
$p^*=\lim\limits_{n\in\omega} p_n$. Clearly $p^*\forces\{k_0,k_1,\ldots\}
\subseteq\dot{X}$ and we are done.
\medskip

\noindent{\sc Case B}\hspace{0.15in}{\em In the course of the play it occurs
that for some $n\in\omega$ the set $B(p_{n+1},n+1)$ is not in the
ultrafilter.}

\noindent Take $q\in\q^{\tree}_1(K,\Sigma)$, $p_{n+1}\leq^1_{n+1} q$ and
fronts $F^*, F_k$ of $T^q$ such that the condition $q^{[\eta]}$ decides the
truth value of ``$k\in\dot{X}$'' for each $\eta\in F_k$ and $k\in\omega$ and 
\[(\forall \nu\in F^*)(\forall\eta\in T^q)(\nu\trianglelefteq\eta\ \
\Rightarrow\ \ \nor[t^q_\eta]>n+3),\]
\[(\forall k\in\omega)(\forall\nu\in F_k)(\forall\eta\in T^q)(\nu
\trianglelefteq\eta\quad\Rightarrow\quad\nor[t^q_\eta]>2^k+n+3)\] 
(possible by \ref{treedec}(2) and \ref{fronor1}). Of course we may assume that 
\[(\forall k\in\omega)(\forall\nu\in F_{k+1})(\exists\eta\in F_k)(\eta
\vartriangleleft\nu)\]
and that the fronts $F_k$ are ``above'' $F^*$ and $F^*$ is ``above''
$F^1_{n+1}(q)$. Further let 
\[u_k=\{\eta\in F_k: q^{[\eta]}\forces_{\q^{\tree}_1(K,\Sigma)}
k\in\dot{X}\}.\]
Look at the set $C=\omega\setminus B(p_{n+1},n+1)\in\D$. If $k\in C$ then
necessarily for some $\rho\in F^*$ there is {\bf no} $\langle
s_\nu:\nu\in\hat{T}\rangle\subseteq K$ with:   
\begin{enumerate}
\item $T\subseteq\{\nu\in T^{q^{[\rho]}}: (\exists\eta\in
F_k)(\nu\trianglelefteq\eta)\}$ is a (well founded) quasi tree with
$\max(T)\subseteq u_k$ and $\mrot(T)=\rho$ 
\item for each $\nu\in\hat{T}$:

$\mrot(s_\nu)=\nu$, $\pos(s_\nu)=\suc_T(\nu)$, $\nor[s_\nu]\geq n+2$ and

$s_\nu\in\Sigma(t^q_\eta:\eta\in S_\nu)$ for some (well founded) quasi tree
$S_\nu\subseteq T^q$.
\end{enumerate}
(Otherwise we could build a condition $\leq^1_{n+1}$--stronger than $q$ (and
thus than $p_{n+1}$) and forcing that $k\in\dot{X}$, contradicting
$k\notin B(p_{n+1},n+1)$.)  As $F^*$ is finite we have one $\rho\in F^*$ for
which the above holds for all $k\in C_0$ for some $C_0\in\D$, $C_0\subseteq
C$. 

But now we may apply the fact that $(K,\Sigma)$ is of the
$\UP(\D)^{\tree}$--type (see \ref{up}, remember the choice of $F^*$): we get
$q^*\in\q^{\tree}_\emptyset (K,\Sigma)$ such that $q^{[\rho]}\leq q^*$,
$Z\stackrel{\rm def}{=}\{k\in C_0: u_k\cap\dcl(T^{q^*})=\emptyset\}\in\D$ and 
\[\nor[t^{q^*}_\nu]\geq \min\{\nor[t^q_\eta]-(n+2): \nu\trianglelefteq\eta\in
T^q\}.\] 
Now we easily see that in fact $q^*\in\q^{\tree}_1(K,\Sigma)$ (by the norm
requirement and the choice of $q$ and $F_k$'s) and
$q^*\forces_{\q^{\tree}_1(K,\Sigma)} Z\cap\dot{X}=\emptyset$. 

For the ``consequently'' part, note that, by \ref{treebound}, the forcing
notion $\q^{\tree}_1(K,\Sigma)$ is proper and $\baire$--bounding. Therefore we
may apply \cite[Ch VI, 5.1]{Sh:f}. This finishes the proof of the
theorem.
\end{proof}

\begin{definition}
\label{rich}
We say that a tree creating pair $(K,\Sigma)$ for $\bH$ is {\em rich} if:\\
for every system $\langle s_\nu:\nu\in\hat{T}\rangle\subseteq K$, $n\in\omega$
and $u$ such that
\begin{enumerate}
\item $T\subseteq\bigcup\limits_{k\in\omega}\prod\limits_{m<k}\bH(m)$ is a
well founded quasi tree, $u\subseteq \max(T)$,
\item $\mrot(s_\nu)=\nu$, $\pos(s_\nu)=\suc_T(\nu)$, $\nor[s_\nu]>n+3$,
\item there is {\bf no} $\langle s^*_\nu: \nu\in\hat{T}^*\rangle\subseteq K$
such that
\[\begin{array}{l}
T^*\subseteq T,\quad \max(T^*)\subseteq u,\quad \mrot(T^*)=\mrot(T),
\quad\mrot(s^*_\nu)=\nu,\\
\pos(s^*_\nu)=\suc_{T^*}(\nu),\quad\nor[s^*_\nu]>n+1,\quad\mbox{ and}\\
s^*_\nu\in\Sigma(s_\eta:\eta\in\hat{T}_\nu)\mbox{ for some }T_\nu\subseteq T
  \end{array}\]
\end{enumerate}
there is $\langle s^+_\nu:\nu\in \hat{T}^+\rangle\subseteq K$ such that
\[\begin{array}{l}
T^+\subseteq T,\quad \max(T^+)\subseteq\max(T)\setminus u,\quad\mrot(T^+)=
\mrot(T),\quad\mrot(s^+_\nu)=\nu,\\ 
\pos(s^+_\nu)=\suc_{T^+}(\nu), \quad
\nor[s^+_\nu] \geq \min\{\nor[s_\eta]: \eta\in T\}-(n+2),\\
\mbox{and }\quad s^+_\nu\in\Sigma(s_\eta:\eta\in\hat{T}_\nu)\mbox{ for some
}T_\nu\subseteq T.\\
  \end{array}\]
\end{definition}

\begin{theorem}
\label{prealmult}
Assume $(K,\Sigma)$ is a finitary $2$-big rich tree--creating pair of the
$\sUP^{\tree}$--type. Let $\D$ be an almost Ramsey interesting ultrafilter on
$\bigcup\limits_{i\in\omega}\{i\}\times (i+1)$. Then 
\[\forces_{\q^{\tree}_1(K,\Sigma)}\mbox{`` $\D$ generates an interesting
ultrafilter on $\bigcup\limits_{i\in\omega}\{i\}\times (i+1)$ ''.}\]
\end{theorem}

\begin{proof}
First note that if we show that, in $\V^{\q^{\tree}_1(K,\Sigma)}$, $\D$
generates an ultrafilter then the ultrafilter has to be interesting (remember
that $\D$ is interesting).

Let $\dot{X}$ be a $\q^{\tree}_1(K,\Sigma)$--name for a subset of $\bigcup
\limits_{i\in\omega}\{i\}\times (i+1)$.

We say that a condition $p\in\q^{\tree}_1(K,\Sigma)$ is {\em
$(\dot{X},n)$--special} if:
\begin{quotation}
\noindent there is a set $C\in\D$ such that:\\
for every $i\in\omega$ and $a\in [C\cap\{i\}\times (i+1)]^{\textstyle \leq n}$
there are a condition $p'\geq p$ and a front $F$ of $T^{p'}$ such that
\[\mrot(p)=\mrot(p'),\quad (\forall\nu\in T^{p'})(\nor[t^{p'}_\nu]>n+1),\] 
\[(\forall\nu\!\in\! T^{p'})(\forall\eta\!\in\! F)(\eta\trianglelefteq\nu\ \
\Rightarrow\ \ t^p_\nu=t^{p'}_\nu)\quad\mbox{ and }\quad p'
\forces_{\q^{\tree}_1(K,\Sigma)} a\subseteq\dot{X}.\]
\end{quotation}
The condition $p$ is {\em $n$--special} if it is either $(\dot{X},n)$--special
or $(\omega\setminus\dot{X}, n)$--special.

Note that the part of the definition of special conditions concerning the
existence of a front $F$ (of $T^{p'}$) is purely technical and usually easy to
get (once we have the rest):
\begin{quotation}
\noindent{\em If\/} the condition $p$ is such that
$(\forall\nu\in T^p)(\nor[t^p_\nu]>n+1)$ and the values of
$\dot{X}\cap\{i\}\times (i+1)$ are decided on some fronts $F_i$ of $T^p$

\noindent{\em then} if we have a condition $p'\geq_0 p$ such that
\[(\forall\nu\in T^{p'})(\nor[t^{p'}_\nu]>n+1)\]
and $p'\forces a\subseteq\dot{X}$, then we may find one (weaker than $p'$)
which has this property and a front $F$ as there. Moreover, in this situation,
if $p\leq_0 p_1$ and $p_1$ is $(\dot{X},n)$--special then $p$ is
$(\dot{X},n)$--special.   
\end{quotation}
Note that if $n\geq m$ and $p$ is $(\dot{X},n)$--special then it is
$(\dot{X},m)$-- special.

\begin{claim}
\label{cl14}
Let $n<\omega$. Suppose that $p\in\q^{\tree}_1(K,\Sigma)$ is such that
\[(\forall\nu\in T^p)(\nor[t^p_\nu]> (2^{2n}+1)(n+3))\]
and there are fronts $F_i$ of $T^p$ (for $i\in\omega$) with
\[(\forall i\in\omega)(\forall\nu\in F_i)(p^{[\nu]}\mbox{ decides }
\dot{X}\cap\{i\}\times (i+1)).\]
Then $p$ is $n$-special.
\end{claim}

\noindent{\em Proof of the claim:}\ \ \ Let $f^+,f^-:
\bigcup\limits_{i\in\omega} [\{i\}\times (i+1)]^{\textstyle n}\longrightarrow
2$ be such that
\begin{quotation}
\noindent $f^+(v)=1$\quad if and only if\quad there are $q\geq_0 p$ and a
front $F$ of $T^q$ such that $q\forces_{\q^{\tree}_1(K,\Sigma)} v\subseteq
\dot{X}$ and
\[(\forall\nu\!\in\! T^q)(\forall\eta\!\in\! F)(\eta\trianglelefteq\nu\
\Rightarrow\ t^q_\nu=t^p_\nu)\ \ \&\ \ (\forall\nu\!\in\! T^q)(\nor[t^q_\nu]>
n+1),\]
\end{quotation}
and $f^-$ is defined similarly replacing ``$\dot{X}$'' by ``$\omega\setminus
\dot{X}$''. 

As $\D$ is almost Ramsey and interesting we find $j^+,j^-<2$ and a set $C\in
\D$ such that for each $i\in\omega$ :
\[\begin{array}{l}
\mbox{if }\ C\cap\{i\}\times (i+1)\neq\emptyset\\
\mbox{then }\ |C\cap\{i\}\times (i+1)|\geq 2n,\quad f^+\rest [C\cap\{i\}\times
(i+1)]^{\textstyle n}=j^+\ \mbox{ and}\\
f^-\rest [C\cap\{i\}\times (i+1)]^{\textstyle n}=j^-.
\end{array}\]
If either $j^+=1$ or $j^-=1$ then plainly $p$ is $n$-special. So suppose that
$j^+=j^-=0$ (and we want to get a contradiction).

Take $i\in\omega$ such that $|C\cap\{i\}\times (i+1)|\geq 2n$ (remember the
choice of $C$) and fix $v\in [C\cap\{i\}\times (i+1)]^{\textstyle 2n}$. For
each $v_1\subseteq v$ let 
\[u^i_{v_1,v}=\{\nu\in F_i: p^{[\nu]}\forces_{\q^{\tree}_1(K,\Sigma)}
v\cap\dot{X}=v_1\}.\]
Since $(K,\Sigma)$ is rich we find $v_1$ and $\langle s^+_\nu: \nu\in\hat{T}^+
\rangle\subseteq K$ such that 
\[\begin{array}{l}
T^+\subseteq T^p,\quad \max(T^+)\subseteq u^i_{v_1,v},\quad
\mrot(T^+)=\mrot(p),\quad\pos(s^+_\nu)=\suc_{T^+}(\nu),\\
\mrot(s^+_\nu)=\nu,\quad\nor[s^+_\nu]>n+1\quad\mbox{ and }
s^+_\nu\in\Sigma(s_\eta: \eta\in\hat{T}_\nu)\ \ \mbox{ for some }T_\nu
\subseteq T^p.
  \end{array}\]
[How? We try successively each $v_1\subseteq v$. If we fail with one, we use
\ref{rich} to pass to a subtree with the minimum of norms dropping down
by at most $n+2$ and we try next candidate. For some $v_1\subseteq v$ we have
to succeed.]\\
Now look at this $v_1$ (and suitable $\langle s^+_\nu: \nu\in
\hat{T}^+\rangle$). Since $j^+=0$ we necessarily have $|v_1|<n$: if not then
we may take $v_2\in [v_1]^{\textstyle n}$ and then $f^+(v_2)=1$ as witnessed
by the condition $q$ starting with $\langle s^+_\nu: \nu\in\hat{T}^+\rangle$.
Similarly, by $j^-=0$, we have $|v\setminus v_1|<n$. Together contradiction to
$|v|=2n$. Thus the claim is proved.
\medskip

Now, let $p\in\q^{\tree}_1(K,\Sigma)$. By
\ref{fronor2} and \ref{treedec}(2) we find $p_1\geq p$ and fronts $F_n$ of
$T^{p_1}$ such that for $n\in\omega$:
\begin{enumerate}
\item $(\forall\nu\in F_{n+1})(\exists\eta\in F_n)(\eta\vartriangleleft\nu)$,
\item $(\forall\nu\in F_n)(\forall\eta\in T^{p_1})(\nu\trianglelefteq\eta\ \
\Rightarrow\ \ \nor[t^{p_1}_\eta]>(2^{2n+1}+5)(n+3))$,
\item $(\forall\nu\in F_n)(p^{[\nu]}_1\mbox{ decides }\dot{X}\cap\{n\}\times
(n+1))$. 
\end{enumerate}
Then, by \ref{cl14}, for each $\nu\in F_n$ the condition $p_1^{[\nu]}$ is
$n$--special. Let
\[u_n=\{\nu\in F_n: p^{[\nu]}_1\mbox{ is $(\dot{X},n)$--special}\}.\]
Now we consider two cases.
\medskip

\noindent{\sc Case A:}\ \ \ {\em There are $n\in\omega$, $\nu\in F_n$ such
that for each $m>n$ there is {\bf no} system $\langle s_\eta:\eta\in\hat{T}
\rangle\subseteq K$ with:
\[\begin{array}{l}
T\subseteq T^{p_1},\quad \max(T)\subseteq u_m,\quad \mrot(T)=\nu,\quad
\mrot(s_\eta)=\eta, \quad \pos(s_\eta)=\suc_{T}(\eta),\\
\nor[s_\eta]\geq (2^{2n}+2)(n+3)\quad\mbox{ and }s_\eta\in
\Sigma(t^{p_1}_\rho:\rho\in \hat{T}_\eta)\ \ \mbox{ for some }T_\eta
\subseteq T^{p_1}.
  \end{array}\]
}
\medskip

\noindent Since $(K,\Sigma)$ is of the $\sUP^{\tree}$--type and 
\[(\forall \eta\in T^{p_1})(\nu\trianglelefteq\eta\ \ \Rightarrow\ \
\nor[t^{p_1}_\eta]> (2^{2n+1}+5)(n+3)),\]
we find a condition $q\in\q^{\tree}_1(K,\Sigma)$ such that
\begin{enumerate}
\item[$(\alpha)$] $p_1^{[\nu]}\leq^1_0 q$,
\item[$(\beta)$]  $Z\stackrel{\rm def}{=}\{m>n:u_m\cap\dcl(T^q)=\emptyset\}\in
\iso,$ 
\item[$(\gamma)$] $(\forall\eta\in T^q)(\nor[t^q_\eta]\geq\min\{\nor[
t^{p_1}_\rho]: \eta\trianglelefteq\rho\in T^{p_1}\}-(2^{2n}+2)(n+3))$, 
\item[$(\delta)$] $(\forall m\in Z)(F_m\cap T^q$ is a front of $T^q)$.
\end{enumerate}
Let $m\in Z$ (so then $F_m\cap T^q$ is a front of $T^q$ and $u_m\cap
T^q=\emptyset$) and let $\eta\in F_m\cap T^q$. By $(\gamma)$ and $(2)$
above we know that
\[(\forall\rho\in T^q)(\eta\trianglelefteq\rho\ \ \Rightarrow\ \
\nor[t^q_\rho]>(2^{2m}+1)(m+3)).\]
Consequently we may use \ref{cl14} to conclude that $q^{[\eta]}$ is
$m$--special. It cannot be $(\dot{X},m)$--special as then the condition
$p^{[\eta]}_1$ would be $(\dot{X},m)$--special (compare the remark after the
definition of special conditions) contradicting $\eta\notin u_m$. Thus
$q^{[\eta]}$ is $(\omega\setminus\dot{X},m)$--special. 
\medskip

For $m\in Z$ and $\eta\in F_m\cap T^q$ fix a set $C^m_\eta\in\D$ witnessing
the fact that `` $q^{[\eta]}$ is $(\omega\setminus\dot{X},m)$--special ''. Let
$B_m=\bigcap\limits_{\eta\in F_m\cap T^q} C^m_\eta\in\D$ (for $m\in Z$).

Consider the following strategy for Player I in the game $G^{sR}(\D)$:
\medskip

\noindent{\em at position number $0$}:\\
Player I writes down to the side: $m_0=\min Z$, $q_0=q$\qquad and he plays:
$B_{m_0}$. 
\smallskip

\noindent [Note that he is (trivially) sure that if $(i_0,a_0)$ is an
answer of Player II then he may find $q_1\geq_0 q_0$ such that
$q_1\forces a_0\subseteq \omega\setminus\dot{X}$ and for some front $F$
of $q_1$, $(\forall\rho\in F)(\forall\eta\in
T^{q_1})(\rho\trianglelefteq\eta\ \Rightarrow\ t^{q_1}_\eta=t^q_\eta)$.]
\medskip

\noindent{\em at position number $k+1$}:\\
Player I looks at the last move $(i_k,a_k)$ of his opponent. He chooses a
condition $q_{k+1}\geq_k q_k$ and a front $F$ of $T^{q_{k+1}}$ such that
\[(\forall\rho\in F)(\forall\eta\in T^{q_{k+1}})(\rho\trianglelefteq\eta\ \
\Rightarrow\ \ t^{q_{k+1}}_\eta=t^q_{\eta})\quad\mbox{ and }\quad q_{k+1}
\forces a_k\subseteq \omega\setminus\dot{X}.\] 
Now he takes $m_{k+1}\in Z$ so large that $m_{k+1}>m_k$ and the front
$F_{m_{k+1}}$ is ``above'' both $F$ and $F^1_{k+2}(q_{k+1})$. Finally:\\
Player I writes down to the side: $m_{k+1}$, $q_{k+1}$\qquad and he plays:
$B_{m_{k+1}}$. 
\smallskip

\noindent [After this move he is sure that if $(i_{k+1},a_{k+1})$ is a
legal answer of the second player then he may find a condition
$q_{k+2}\geq _{k+1} q_{k+1}$ such that $q_{k+2}\forces a_{k+1}\subseteq
\omega\setminus\dot{X}$ and for some front $F$ of $T^{q_{k+2}}$,
$\rho\trianglelefteq\eta\ \ \Rightarrow\ \ t^{q_{k+2}}_\eta=t^q_\eta$
whenever $\rho\in F$, $\eta\in T^{q_{k+2}}$. Why? Remember the choice of
$Z,m_{k+1}, B_{m_{k+1}}$ and $F_{m_{k+1}}$; see \ref{fusfront}, clearly
$m_{k+1}\geq k+1$.]
\medskip

\noindent The strategy described above cannot be the winning one.
Consequently there is a sequence $\langle (i_0,a_0),(i_1,a_1),\ldots\rangle$
such that $(i_n,a_n)$ are legitimative moves of Player II against the strategy
and $\bigcup\limits_{k\in\omega}a_k\in\D$. But in this play, Player I
constructs (on a side) a sequence $p\leq q=q_0\leq_0 q_1\leq_1
q_2\leq_2\ldots$ of conditions such that $q_k\forces a_k\subseteq\omega
\setminus\dot{X}$. Take the limit condition $q_\infty=\lim\limits_{k\in\omega}
q_k$; it forces that $\bigcup\limits_{k\in\omega} a_k\subseteq\omega\setminus
\dot{X}$, finishing the proof of the theorem in Case A. 
\medskip

\noindent{\sc Case B:}\ \ \ {\em Not {\sc Case A}.}
\medskip

\noindent Thus for every $\nu\in F_n$, $n\in\omega$ we find $m>n$ and $\langle
s_\eta\!:\eta\in\hat{T}\rangle\subseteq K$ such that:
\[\begin{array}{l}
T\subseteq T^{p_1},\ \max(T)\subseteq u_m,\ \mrot(T)=\nu,\quad
\mrot(s_\eta)=\eta, \ \pos(s_\eta)=\suc_{T}(\eta),\\
\nor[s_\eta]\geq (2^{2n}+2)(n+3)\quad\mbox{ and }s_\eta\in
\Sigma(t^{p_1}_\rho:\rho\in\hat{T}_\eta)\ \ \mbox{ for some 
}T_\eta\subseteq T^{p_1}.  
  \end{array}\]
If $\langle s_\eta:\eta\in\hat{T}\rangle$, $m>n$ are as above then we 
will say that $\langle s_\eta:\eta\in\hat{T}\rangle$ {\em is $(m,n)$--good for
$\nu$}. For each $n\in\omega$ and $\nu\in u_n$ fix a set $C^n_\nu\in\D$
witnessing the fact that ``the condition $p_1^{[\nu]}$ is
$(\dot{X},n)$--special''. Now consider the following strategy for Player I in
the game $G^{sR}(\D)$:   
\medskip

\noindent{\em at position number $0$}:\\
For each $\nu\in F_0$, Player I chooses $n(\nu)>0$ and $\langle s^\nu_\eta:
\eta\in\hat{T}_\nu\rangle\subseteq K$ which is $(n(\nu),0)$--good. He
builds a condition $q_0$ which starts like these sequences. Thus $q_0$
is such that:
\begin{enumerate}
\item[(a)] $\mrot(q_0)=\mrot(p_1)$,
\item[(b)] if $\eta\in T^{p_1}$ is below the front $F_0$ then $\eta\in
T^{q_0}$, $t^{q_0}_\eta=t^{p_1}_\eta$,
\item[(c)] if $\eta\in\hat{T}_\nu$ for some $\nu\in F_0$ then $\eta\in
T^{q_0}$, $t^{q_0}_\eta=s^{\nu}_\eta$,
\item[(d)] if $\eta\in T^{p_1}$ and there are $\nu\in F_0$ and $\rho\in
\max(T_\nu)$ such that $\rho\trianglelefteq\eta$ then $\eta\in T^{q_0}$,
$t^{q_0}_\eta=t^{p_1}_\eta$.
\end{enumerate}
Note that $F^*_0\stackrel{\rm def}{=}\bigcup\limits_{\nu\in F_0}
\max(T_\nu)$ is a front of $T^{q_0}$ above $F^0_2(q_0)$.

Now, Player I writes down to the side: $q_0$, $F^*_0$\qquad and he plays: 
\[A_0=\bigcap\{C^{n(\nu)}_\eta: \nu\in F_0,\ \eta\in\max(T_\nu)\}.\]
[Thus Player I knows that $F^*_0$ is a front of $T^{q_0}$ above $F^0_2(q_0)$,
and for each $\nu\in F^*_0$ the condition $p^{[\nu]}_1=q_0^{[\nu]}$ is
$(\dot{X},1)$--special and the set $A_0$ witnesses this fact.] 
\medskip

\noindent{\em at the position number $k+1$}:\\
Player I looks at the last move $(i_k,a_k)$ of his opponent. He chooses a
condition $q_{k+1}\geq_k q_k$ and a front $F^*_{k+1}$ of $T^{q_{k+1}}$
such that $q_{k+1}\forces a_k\subseteq\dot{X}$, the front $F^*_{k+1}$ is
above $F^0_{k+3}(q_{k+1})$ and for each $\nu\in F^*_{k+1}$ the condition
$q_{k+1}^{[\nu]}=p_1^{[\nu]}$ is $(\dot{X},k+2)$--special. How does he
find $q_{k+1}$ and $F^*_{k+1}$? He has $q_k$ and $F^*_k$ and he knows
that $F^*_k$ is a front of $T^{q_k}$ above $F^0_{k+2}(q_k)$ and for each
$\nu\in F^*_k$, the condition $q_k^{[\nu]}=p_1^{[\nu]}$ is
$(\dot{X},k+1)$-special and the set $A_k$ (played by Player I before)
witnesses this fact (this is our inductive hypothesis). Now, as $a_k\subseteq
A_k$, for each $\nu\in F^*_k$ the first player may choose a condition
$q^+_\nu\geq_0 q^{[\nu]}_k=p_1^{[\nu]}$ which forces that
$a_k\subseteq\dot{X}$ and such that $(\forall\eta\in
T^{q^+_\nu})(\nor[t^{q^+_\nu}_\eta]>k+1)$ and for some front $F^+_\nu$ 
of $q^+_\nu$, if $\eta\trianglelefteq\rho\in T^{p_1}$, $\eta\in F^+_\nu$
then $\rho\in T^{q^+_\nu}$, $t^{q^+_\nu}_\rho=t^{p_1}_\rho$. We may
assume that the fronts $F^+_\nu$ are such that $F^+_\nu\subseteq
F_{m^+(\nu)}$ for some $m^+(\nu)>k+2$ and $F^+_\nu$ is above
$F^0_{k+3}(q^+_\nu)$. For each $\rho\in F^+_\nu$, $\nu\in F^*_k$ we may
choose $m(\rho)>m^+(\nu)$ and $\langle s^\rho_\eta: \eta\in\hat{T}_\rho
\rangle\subseteq K$ such that $\langle s^\rho_\eta: \eta\in \hat{T}_\rho
\rangle$ is $(m(\rho),m^+(\nu))$--good. Let
\[F^*_{k+1}=\bigcup\{\max(T_\rho): (\exists\nu\in F^*_k)(\nu\vartriangleleft
\rho\in F^+_\nu)\}.\]
The condition $q_{k+1}$ is such that:
\begin{enumerate}
\item[(a)] below $F^*_k$ it agrees with $q_k$,
\item[(b)] if $\nu\trianglelefteq\eta\vartriangleleft\rho \in F^+_\nu$,
$\nu\in F^*_k$, $\eta\in T^{q^+_\nu}$ then $\eta\in T^{q_{k+1}}$,
$t^{q_{k+1}}_\eta=t^{q^+_\nu}_\eta$,
\item[(c)] if $\nu\vartriangleleft\rho\trianglelefteq\eta\in
\hat{T}_\rho$, $\nu\in F^*_k$, $\rho\in F^+_\nu$ then $\eta\in
T^{q_{k+1}}$, $t^{q_{k+1}}_\eta=s^\rho_\eta$,
\item[(d)] if $\nu\vartriangleleft\rho\vartriangleleft\eta_0
\trianglelefteq\eta_1$, $\nu\in F^*_k$, $\rho\in F^+_\nu$,
$\eta_0\in\max(T_\rho)$, $\eta_1\in T^{p_1}$ then $\eta_1\in
T^{q_{k+1}}$ and  $t^{q_{k+1}}_{\eta_1}=t^{p_1}_{\eta_1}$.
\end{enumerate}
Thus $q_{k+1}\geq_k q_k$, $F^*_{k+1}$ is a front of $T^{q_{k+1}}$, and
if $\nu\vartriangleleft\rho\vartriangleleft\eta\in F^*_{k+1}$, $\nu\in
F^*_k$, $\rho\in F^+_\nu$ then $ q^{[\eta]}_{k+1}=p_1^{[\eta]}$ is
$(\dot{X},m(\rho))$--special, $m(\rho)>m^+(\nu)>k+2$ (so
$m(\rho)> k+3$).\\
Now the first player writes down to the side $q_{k+1}$, $F^*_{k+1}$\qquad and
he plays: 
\[A_{k+1}=\bigcap\{C^{m(\rho)}_\eta: (\exists\nu\in F^*_k)(\rho\in F^+_\nu\ \
\&\ \ \nu\vartriangleleft\rho\vartriangleleft\eta\in u_{m(\rho)})\}.\]
[Note that for each $\nu\in F^*_{k+1}$, the set $A_{k+1}$ witnesses that the
condition $q_{k+1}^{[\nu]}=p_1^{[\nu]}$ is $(\dot{X},k+2)$-special; the front
$F^*_{k+1}$ is above $F^0_{k+3}(q_{k+1})$.] 
\medskip

The strategy described above cannot be the winning one. Consequently, there is
a play according to this strategy in which Player I loses. Thus we have moves
$(i_0,a_0), (i_1,a_1),\ldots$ of Player II (legal in this play) for which
$\bigcup\limits_{k\in\omega} a_k\in\D$. But in the course of the play the
first player constructs conditions $p_1\leq q_0\leq_0 q_1\leq_1 q_2\leq_2
\ldots$ such that $q_k\forces a_k\subseteq\dot{X}$. Then the limit condition
$q_\infty=\lim\limits_{k\in\omega} q_k$ forces that $\bigcup\limits_{k\in
\omega} a_k\subseteq \dot{X}$. This finishes the proof of the theorem. 
\end{proof}

\begin{definition}
\label{ramnum}
For $n,k,m\in\omega$ let $R_n(k,m)$ be the smallest integer such
that for every colouring $f:[R_n(k,m)]^{\textstyle {\leq}n}\longrightarrow k$
there is $a\in [R_n(k,m)]^{\textstyle m}$ homogeneous for $f$ (so this is the
respective Ramsey number).  
\end{definition}

\begin{theorem}
\label{prealmram}
Let $\D$ be an almost Ramsey ultrafilter on $\bigcup\limits_{i\in\omega}\{i\}
\times (i+1)$. Suppose that $(K,\Sigma)$ is a finitary $2$-big rich
tree--creating pair such that 
\[\forces_{\q^{\tree}_1(K,\Sigma)}\mbox{`` $\D$ generates an ultrafilter on
$\bigcup\limits_{i\in\omega}\{i\}\times (i+1)$ ''.}\] 
Then $\forces_{\q^{\tree}_1(K,\Sigma)}$`` the ultrafilter $\tilde{\D}$
generated by $\D$ is almost Ramsey ''.
\end{theorem}

\begin{proof}
Due to \ref{treebound}(1) we may apply \ref{almboupre} and get that
\[\forces_{\q^{\tree}_1(K,\Sigma)}\mbox{`` $\D$ generates a semi--Ramsey
ultrafilter on $\bigcup\limits_{i\in\omega}\{i\}\times (i+1)$ ''.}\] 
So, what we have to do is to show that, in $\V^{\q^{\tree}_1(K,\Sigma)}$, the
ultrafilter $\tilde{\D}$ generated by $\D$ has the colouring property of
\ref{almramsey}(5). Suppose that a condition $p\in\q^{\tree}_1(K,\Sigma)$,
$n,L<\omega$ and a $\q^{\tree}_1(K,\Sigma)$--name $\dot{\varphi}$ are such
that $p\forces$``$\dot{\varphi}:\bigcup\limits_{i\in\omega}[\{i\}\times (i+
1)]^{\textstyle n}\longrightarrow L$''. Note that if $\D$ is not interesting
and it is witnessed by $h\in\prod\limits_{i\in\omega} (i+1)$ then the set
$\{(i,h(i)): i\in\omega\}$ is almost homogeneous for $\dot{\varphi}$ and it is
in $\D$. Consequently, we may assume that $\D$ is interesting. 

We say that a condition $q\in\q^{\tree}_1(K,\Sigma)$ is {\em $m$-beautiful}
(for $m \in\omega$) if:
\begin{quotation}
\noindent there is a set $C\in\D$ such that:\\
for every $i\in\omega$ and $a\in [C\cap\{i\}\times (i+1)]^{\textstyle m}$
there are a condition $q'\geq q$ and a front $F$ of $T^{q'}$ such that
\[\mrot(q)=\mrot(q'),\quad (\forall\nu\in T^{q'})(\nor[t^{q'}_\nu]>m+1),\]
\[\hspace{-0.3cm}q'\forces\mbox{``}a \mbox{ is }\dot{\varphi}
\mbox{--homogeneous''\quad and\quad }(\forall\nu\!\in\! T^{q'})(\forall\eta\!
\in\! F)(\eta\trianglelefteq\nu\ \ \Rightarrow\ \ t^q_\nu=t^{q'}_\nu).\]
\end{quotation}
Clearly, if $q$ is $m$-beautiful and $k\leq m$ then $q$ is $k$-beautiful.

Let $R^*(m)=R_n(L,m)$ (see \ref{ramnum}). 

\begin{claim}
\label{cl15}
Let $m<\omega$, $q\in\q^{\tree}_1(K,\Sigma)$. Assume that
\[(\forall \nu\in T^q)(\nor[t^q_\nu] > ({R^*(m)\choose m}+1)(m+3))\]
and there are fronts $F_i$ of $T^q$ (for $i\in\omega$) such that
conditions $q^{[\nu]}$ (for $\nu\in F_i$, $i\in\omega$) decide the value of
$\dot{\varphi}\rest [\{i\}\times(i+1)]^{\textstyle n}$. Then the condition $q$
is $m$--beautiful. 
\end{claim}

\noindent{\em Proof of the claim:}\ \ \ Look at the following colouring
$f:\bigcup\limits_{i\in\omega}[\{i\}\times (i+1)]^{\textstyle m}
\longrightarrow 2$:
\begin{quotation}
\noindent $f(v)=1$\quad if and only if\quad there are $q'\geq_0 q$ and a front
$F$ of $T^{q'}$ such that 
\[q'\forces_{\q^{\tree}_1(K,\Sigma)} \mbox{``}v\mbox{ is homogeneous for
}\dot{\varphi}\mbox{''}\quad\quad\mbox{ and}\]
\[(\forall\nu\in T^{q'})(\forall\eta\in F)(\eta\trianglelefteq\nu\ \Rightarrow
\ t^{q'}_\nu=t^q_\nu)\ \ \&\ \ (\forall\nu\in T^{q'})(\nor[t^{q'}_\nu]>m+1).\]
\end{quotation}
Since $\D$ is almost Ramsey and interesting we find $j<2$ and a set $C\in\D$
such that for each $i\in\omega$: 
\[\begin{array}{l}
\mbox{if }\ C\cap\{i\}\times (i+1)\neq\emptyset\\
\mbox{then }\ |C\cap\{i\}\times (i+1)|\geq R^*(m),\quad f\rest [C\cap\{i\}
\times (i+1)]^{\textstyle m}=j.\\
\end{array}\]
If $j=1$ then easily the condition  $q$ is $m$-beautiful. Thus we have
to exclude the other possibility. So assume $j=0$. Take $i\in\omega$ such that
$C\cap\{i\}\times (i+1)\neq\emptyset$ and choose $v\in [C\cap\{i\}\times
(i+1)]^{\textstyle R^*(m)}$. For $v_1\in [v]^{\textstyle m}$ put
\[u^i_{v_1,v}=\{\nu\in F_i: q^{[\nu]}\forces_{\q^{\tree}_1(K,\Sigma)}
\mbox{``}v_1\mbox{ is }\dot{\varphi}\mbox{--homogeneous''}\}.\]
Note that, by the definition of $R^*(m)$, for each $\nu\in F_i$ we find $v_1
\in [v]^{\textstyle m}$ such that $\nu\in u^i_{v_1,v}$. Consequently we may
apply the assumption that $(K,\Sigma)$ is rich and we find $v_1\in [v]^{
\textstyle m}$ and $\langle s^+_\nu: \nu\in\hat{T}^+\rangle\subseteq K$ such
that  
\[\begin{array}{l}
T^+\subseteq T^q,\ \max(T^+)\subseteq u^i_{v_1,v},\ \mrot(T^+)=\mrot(q),\quad
\pos(s^+_\nu)=\suc_{T^+}(\nu),\\ 
\mrot(s^+_\nu)=\nu,\ \nor[s^+_\nu]>m+1\quad\mbox{ and }s^+_\nu\in\Sigma(
s_\eta: \eta\in \hat{T}_\nu)\ \mbox{ for some }T_\nu\subseteq T^q. 
 \end{array}\]
[Exactly like in \ref{cl14}.] But with this in hands we easily conclude
that $f(v_1)=1$, contradicting $v_1\in[C\cap\{i\}\times(i+1)]^{\textstyle m}$
and the choice of $j,C$.  
\medskip

Choose a condition $q\geq p$ such that for some fronts $F_m$ of $T^q$ (for
$m\in\omega$) we have  
\begin{enumerate}
\item $(\forall\nu\in F_{m+1})(\exists\eta\in F_m)(\eta\vartriangleleft\nu)$,
\item $(\forall\nu\in F_m)(\forall\eta\in T^{q})(\nu\trianglelefteq\eta\ \
\Rightarrow\ \ \nor[t^{q}_\eta]>({R^*(m)\choose m}+1)(m+3))$,
\item $(\forall\nu\in F_m)(q^{[\nu]}\mbox{ decides }\dot{\varphi}\rest
[\{m\}\times(m+1)]^{\textstyle n})$ 
\end{enumerate}
(possible by \ref{fronor2} and \ref{treedec}(2)). By \ref{cl15} we know that
for each $\nu\in F_m$, $m\in\omega$ the condition $q^{[\nu]}$ is
$m$--beautiful. So for every $m<\omega$ and $\nu\in F_m$ we may fix a set
$C^m_\nu\in\D$ witnessing ``$q^{[\nu]}$ is $m$--beautiful''.

Consider the following strategy of the first player in the game $G^{sR}(\D)$:
\medskip

\noindent{\em at position number $0$}:\\
Player I writes down to the side $q_0=q$, $F^*_0=F_0$.\qquad He plays $\bigcap
\{C^0_\nu: \nu\in F_0\}\in\D$.
\medskip

\noindent{\em arriving at position $k+1$}:\\
Player I has a condition $q_k$ and a front $F^*_k$ of $T^{q_k}$ such that
$q_k$ above $F^*_k$ agrees with $q$. Moreover, the set played by him before
witnesses that each $q^{[\nu]}$ is $k$-beautiful (for $\nu\in F^*_k$). He
looks at the last move $(i_k,a_k)$ of his opponent. For each $\nu\in F^*_k$,
Player I can find a condition $q_\nu\geq_0 q^{[\nu]}_k=q^{[\nu]}$ such that
\begin{enumerate}
\item[(a)] $(\forall\eta\in T^{q_\nu})(\nor[t^{q_\nu}_\eta]>k+1)$,
\item[(b)] for some fronts $F_\nu$ of $T^{q_\nu}$, $q_\nu$ above $F_\nu$
agrees with $q$ and
\item[(c)] $q_\nu\forces_{\q^{\tree}_1(K,\Sigma)}$``$a_k$ is
$\dot{\varphi}$--homogeneous''. 
\end{enumerate}
We may think that for some $m>i_k$ the fronts $F_\nu$ are contained in $F_m$
(for $\nu\in F^*_k$). Now let $q_{k+1}$ be such that 
\begin{enumerate}
\item[$(\alpha)$] below $F^*_k$ it agrees with $q_k$,
\item[$(\beta)$]  if $\nu\trianglelefteq\eta\vartriangleleft\rho \in F_\nu$,
$\nu\in F^*_k$, $\eta\in T^{q_\nu}$ then $\eta\in T^{q_{k+1}}$,
$t^{q_{k+1}}_\eta=t^{q_\nu}_\eta$, 
\item[$(\gamma)$] if $\nu\vartriangleleft\rho\trianglelefteq\eta\in
T^{q_\nu}$, $\nu\in F^*_k$, $\rho\in F_\nu$ then $\eta\in T^{q_{k+1}}$,
$t^{q_{k+1}}_\eta=t^{q_\nu}_\eta=t^q_\eta$. 
\end{enumerate}
Let $F^*_{k+1}=\bigcup\{F_\nu: \nu\in F^*_k\}$. Clearly $F^*_{k+1}$ is a front
of $T^{q_{k+1}}$ contained in $F_m$, $q_{k+1}\geq_k q_k$ and $q_{k+1}$ forces
that $a_k$ is $\dot{\varphi}$-homogeneous. Now:\\
Player I writes down to the side $q_{k+1},F^*_{k+1}$\ \ and he plays the set 
\[\bigcap\{C^m_\eta:\eta\in F^*_{k+1}\}\in\D.\]
[Note that for every $\nu\in F^*_{k+1}$ the condition $q^{[\nu]}_{k+1}=
q^{[\nu]}$ is $k+1$-beautiful and the set played by Player I witnesses it.]
\medskip

This strategy cannot be winning for the first player. Consequently he loses
some play according to it. Let $\langle (i_0,a_0),(i_1,a_1),\ldots\rangle$ be
the sequence of the respective moves of the second player (so $\bigcup
\limits_{k\in\omega} a_k\in\D$) and let $q_0,q_1,q_2,\ldots$ be the sequence
of conditions written down by the first player during the play. Let $q_\infty$
be the limit condition $\lim\limits_{k\in\omega} q_k$. Then we have: 
\[q_\infty\forces_{\q^{\tree}_1(K,\Sigma)}(\forall k\in\omega)(a_k\mbox{ is }
\dot{\varphi}\mbox{--homogeneous }).\]
This finishes the proof of the theorem. 
\end{proof}

Let us finish this section with a theorem showing how several types of forcing
notions built according to our scheme may preserve $\Gamma$--genericity in the
context of ultrafilters. The proof of theorem \ref{presGpoint} below resembles
the proof that Blass--Shelah forcing notion preserves $p$--points (see
\cite[3.3]{BsSh:242}). 

\begin{definition}
\label{simomi}
\begin{enumerate}
\item A creating pair $(K,\Sigma)$ is {\em simple except omitting} if it is
omittory, $|\pos(u,t)|>1$ whenever $t\in K$, $\nor[t]>0$ and $u\in\basis(t)$,
and for every $(t_0,\ldots,t_{n-1})\in\PFC(K,\Sigma)$ and $s\in\Sigma(t_0,
\ldots,t_{n-1})$ there is $k<n$ such that $s\in\Sigma(t_k\Rsh [m^s_{\dn},
m^s_{\up}))$.
\item Suppose that $(K,\Sigma)$ is an omittory creating pair, $\bar{t}=\langle
t_k: k<\omega\rangle\in\PC_\infty(K,\Sigma)$ and $W:\omega\times\omega\times
\pfs\longrightarrow K$ is a $\bar{t}$--system. We say that $W$ is {\em
omittory compatible} if
\begin{enumerate}
\item[$(\alpha)$] $k\leq\ell<n<\omega$ and $s\in\Sigma(t_\ell\Rsh [m^{t_k}_{
\dn},m^{t_n}_{\dn}))$ imply $n_W(s)=\ell$ (where $n_W(s)$ is as in
\ref{quasi}(1b,c)), and 
\item[$(\beta)$]  $k\leq\ell<\omega$, $\sigma:[m^{t_k}_{\dn},m^{t_\ell}_{\up})
\longrightarrow\omega$ imply 
\end{enumerate}
\[W(m^{t_k}_{\dn},m^{t_\ell}_{\up},\sigma)=\bigcup\{s\!\Rsh\! [m^{t_k}_{\dn},
m^{t_\ell}_{\up})\!:s\in W(m^{t_n}_{\dn},m^{t_n}_{\up},\sigma\rest[m^{t_n}_{
\dn},m^{t_n}_{\up})),\ k\!\leq\! n\!\leq\!\ell\}.\]
\end{enumerate}
\end{definition}

\begin{theorem}
\label{presGpoint}
Suppose that a creating pair $(K_0,\Sigma_0)$ is simple except omitting,
forgetful and monotonic. Let $\bar{t}\in\PC_\infty(K,\Sigma)$ and $W$ be an
omittory--compatible $\bar{t}$--system. Assume that $\Gamma\subseteq\PC_\infty
(K_0,\Sigma_0)$ is quasi-$W$-generic and generates an ultrafilter.
\begin{enumerate}
\item If $(K_0,\Sigma_0)$ is strongly finitary, $(K,\Sigma)$ is a finitary,
monotonic, omittory and omittory--big creating pair then the forcing notion
$\q^*_{{\rm s}\infty}(K,\Sigma)$ is $\Gamma$--genericity preserving.
\item If $(K,\Sigma)$ is a finitary creating pair which captures singletons 
then the forcing notion $\q^*_{{\rm w}\infty}(K,\Sigma)$ is
$\Gamma$--genericity preserving.
\item If $(K,\Sigma)$ is a finitary t--omittory tree creating pair then the
forcing notion $\q^{\tree}_1(K,\Sigma)$ is $\Gamma$--genericity preserving.
\end{enumerate}
\end{theorem}

\begin{proof}
1)\ \ \ We have to show that, in $\V^{\q^*_{{\rm s}\infty}(K,\Sigma)}$,
the demand \ref{quasi}(3b) is satisfied (as we know that $\q^*_{{\rm s}\infty}
(K,\Sigma)$ is proper). So suppose that $\dot{\eta}$ is a $\q^*_{{\rm
s}\infty}(K,\Sigma)$--name for an element of $\baire$ and $p\in \q^*_{{\rm s}
\infty}(K,\Sigma)$. Since $(K_0,\Sigma_0)$ is strongly finitary we may
repeatedly use \ref{omitdecbel} and we get a condition $q\geq p$ such that for
each $n<\omega$ and $v\in\pos(w^q,t^q_0,\ldots,t^q_{n-1})$ the condition
$(v,t^q_n,t^q_{n+1},\ldots)$ decides the value of $W(m^{t_n}_{\dn},m^{t_n}_{
\up},\dot{\eta}\rest [m^{t_n}_{\dn},m^{t_n}_{\up}))$, where $\bar{t}=\langle
t_0,t_1,t_2,\ldots\rangle$. Fix $v\in\pos(w^q,t^q_0,\ldots\!,t^q_{n-1})$,
$n<\omega$ for a moment. For $k\geq n$ let $v^k=v\conc {\bf 0}_{[\lh(v),
m^{t^q_k}_{\up})}\in\pos(v,t^q_n,\ldots,t^q_k)$ (remember that $(K,\Sigma)$ is
omittory). Note that, for each $k\geq n$, the condition $(v^k,t^q_{k+1},
t^q_{k+2},\ldots)$ decides the value of $W(m^{t_{k+1}}_{\dn},m^{t_{k+1}}_{
\up},\dot{\eta}\rest [m^{t_{k+1}}_{\dn},m^{t_{k+1}}_{\up}))$. Thus we find a
function $\eta(v)\in\baire$ such that, for each $\ell\leq n$, 
\[(v,t^q_n,t^q_{n+1},\ldots)\forces W(m^{t_{\ell}}_{\dn},m^{t_{\ell}}_{\up},
\dot{\eta}\rest [m^{t_{\ell}}_{\dn},m^{t_{\ell}}_{\up}))=W(m^{t_{\ell}}_{\dn},
m^{t_{\ell}}_{\up},\eta(v)\rest [m^{t_{\ell}}_{\dn},m^{t_{\ell}}_{\up}))\]
and for each $k\geq n$ the condition $(v^k,t^q_{k+1},t^q_{k+2},\ldots)$ forces
that
\[\mbox{`` }W(m^{t_{k+1}}_{\dn},m^{t_{k+1}}_{\up},\dot{\eta}\rest [m^{t_{k+
1}}_{\dn},m^{t_{k+1}}_{\up}))=W(m^{t_{k+1}}_{\dn},m^{t_{k+1}}_{\up},\eta(v)
\rest [m^{t_{k+1}}_{\dn},m^{t_{k+1}}_{\up}))\mbox{ ''.}\] 
Since $\Gamma$ is quasi-$W$-generic we find $\bar{s}=\langle s_m: m<\omega
\rangle\in\Gamma$ such that for every $n<\omega$ and $v\in\pos(w^q,t^q_0,
\ldots,t^q_{n-1})$ we have
\[(\forall^\infty m\in\omega)(s_m\in W(m^{s_m}_{\dn},m^{s_m}_{\up},\eta(v)
\rest [m^{s_m}_{\dn},m^{s_m}_{\up})))\]
(remember $(\Gamma,\preceq)$ is directed countably closed). Next we choose
inductively an increasing sequence $\ell(0)<\ell(1)<\ell(2)<\ldots<\omega$
such that 
\begin{quotation}
\noindent if $v\in\pos(w^q,t^q_0,\ldots,t^q_{n-1})$,\ \ $m^{t_n}_{\up}\leq
m^{s_{\ell(i)}}_{\up}$,\ \ $n,i\in\omega$

\noindent then $(\forall m\geq\ell(i+1))(s_m\in W(m^{s_m}_{\dn},m^{s_m}_{\up},
\eta(v)\rest [m^{s_m}_{\dn},m^{s_m}_{\up})))$.
\end{quotation}
For $j<4$ let $Y_j=\bigcup\limits_{i\in\omega}[m^{s_{\ell(4i+j)}}_{\dn},
m^{s_{\ell(4i+j+1)}}_{\dn})$. Since $\Gamma$ generates an ultrafilter, exactly
one of these sets is in $\D(\Gamma)$. Without loss of generality we may assume
that $Y_2\in\D(\Gamma)$ (otherwise start the sequence of the $\ell(i)$'s from
$\ell(1)$, or $\ell(2)$ or $\ell(3)$). This means that we find $\bar{s}^*=
\langle s^*_m: m<\omega\rangle\in\Gamma$ such that $\bar{s}\leq\bar{s}^*$
and for sufficiently large $m$ (say $m\geq m^*$), there are $i=i(m)$ and
$k=k(m)$ such that 
\[m^{s_{\ell(4i+2)}}_{\dn}\leq m^{t_k}_{\dn}<m^{t_k}_{\up}\leq m^{s_{\ell(
4i+3)}}_{\dn}\quad \mbox{ and }\quad s^*_m\in\Sigma(t_k\Rsh [m^{s^*_m}_{\dn},
m^{s^*_m}_{\up}))\]
(remember $(K_0,\Sigma_0)$ is simple--except--omitting and monotonic). Let
$x(i)\in\omega$ be such that $m^{t_{x(i)}}_{\dn}=m^{s_{\ell(4i+3)}}_{\dn}$
(for $i\in\omega$). Now we define a condition $q^*=(w^{q^*},t^{q^*}_0,t^{
q^*}_1,\ldots)$:
\[w^{q^*}=w^q,\quad t^{q^*}_0=t^q_{x(0)}\Rsh [\lh(w^q),m^{t^q_{x(0)}}_{\up}),
\quad\mbox{and}\quad t^{q^*}_{i+1}=t^q_{x(i+1)}\Rsh [m^{t^q_{x(i)}}_{\up},
m^{t^q_{x(i+1)}}_{\up}).\]
Plainly, this defines a condition in $\q^*_{{\rm s}\infty}(K,\Sigma)$ stronger
than $q$. We claim that 
\[q^*\forces (\forall m>m^*)(s^*_m\in W(m^{s^*_m}_{\dn},m^{s^*_m}_{\up},
\dot{\eta}\rest [m^{s^*_m}_{\dn},m^{s^*_m}_{\up}))).\]
Assume not. Then we find $q'\geq q^*$ and $m>m^*$ such that 
\[q'\forces s^*_m\notin W(m^{s^*_m}_{\dn},m^{s^*_m}_{\up},\dot{\eta}\rest
[m^{s^*_m}_{\dn},m^{s^*_m}_{\up}))).\] 
Of course we may assume that $\lh(w^{q'})\geq m^{t^{q^*}_{i(m)}}_{\up}$. Let
$k^0,k^1,k^2$ be such that $m^{t_{k^\ell}}_{\dn}=m^{s_{\ell(4i(m)+\ell)}}_{
\dn}$ and let $v=w^{q'}\rest m^{t^q_{k^0}}_{\dn}$. Clearly
\[m^{t^{q^*}_{i(m)}}_{\dn}<m^{t^q_{k^0}}_{\up}<m^{t^q_{k^1}}_{\up}<m^{t^q_{
k^2}}_{\up}\leq m^{t^q_{k(m)}}_{\up}<m^{t^q_{x(i(m))}}_{\up}=m^{t^{q^*}_{i(m)
}}_{\up},\]
$v\in\pos(w^q,t^q_0,\ldots,t^q_{k^0-1})$ (by smoothness) and $w^{q'}\rest
[m^{t^q_{k^0}}_{\dn},m^{t^q_{x(i(m))}}_{\dn})\equiv{\bf 0}$ (by the definition
of $t^{q^*}_{i(m)}$, remember $(K,\Sigma)$ is monotonic). Consequently for
each $k\in [k^0,x(i(m)))$ we have $v^k\trianglelefteq w^{q'}$. Since $k^0<k(m)
<x(i(m))$ we conclude that the condition $q'$ forces 
\[\mbox{`` }W(m^{t_{k(m)}}_{\dn}\!,m^{t_{k(m)}}_{\up}\!,\dot{\eta}\rest
[m^{t_{k(m)}}_{\dn}\!,m^{t_{k(m)}}_{\up}))=W(m^{t_{k(m)}}_{\dn}\!,m^{t_{k(m)
}}_{\up}\!,\eta(v)\rest [m^{t_{k(m)}}_{\dn}\!,m^{t_{k(m)}}_{\up}))\mbox{ ''.}\] 
Since $m^{t_{k^0}}_{\up}<m^{s_{\ell(4i(m)+1)}}_{\up}$, for each $i\in
[\ell(4i(m)+2),\ell(4i(m)+3))$ we have $s_i\in W(m^{s_i}_{\dn},m^{s_i}_{\up},
\eta(v)\rest [m^{s_i}_{\dn},m^{s_i}_{\up}))$ (remember that choice of the
sequence of the $\ell(i)$'s). Now look at the condition \ref{quasi}(1c). Since
$W$ is omittory--compatible we have $n_W(s^*_m)\!=\!k(m)\!\in\! [k^2,x(i(m)))$
and therefore $s^*_m\in W(m^{s^*_m}_{\dn},m^{s^*_m}_{\up},\eta(v)\rest
[m^{s^*_m}_{\dn},m^{s^*_m}_{\up}))$. Hence, by \ref{simomi}(2$\beta$), we
easily get 
\[q'\forces \mbox{`` }s^*_m\in W(m^{s^*_m}_{\dn},m^{s^*_m}_{\up},\dot{\eta}
\rest [m^{s^*_m}_{\dn},m^{s^*_m}_{\up}))\mbox{ '',}\]
a contradiction. 
\medskip

\noindent 2)\ \ \ Repeat the proof of 1) noting that we do not have to assume
that $(K_0,\Sigma_0)$ is strongly finitary as we use \ref{sinwin} instead of
\ref{omitdecbel}. (Defining $v^k$ we use fixed sequences $u_n$ (for
$n\in\omega$) such that $u_n\in\pos(u,t^q_n)$ for each $u\in\basis(t^q_n)$.)
\medskip

\noindent 3)\ \ \ Similarly (remember that $(K,\Sigma)$ is finitary). 
\end{proof}

\begin{remark}
In \ref{presGpoint}(3) we need the assumption that the tree creating pair
$(K,\Sigma)$ is finitary. The forcing notion ${\mathbb D}_\omega$ of
\cite{NeRo93}, in which conditions are trees $\subseteq\fseo$ such that each
node has an extension which has all possible successors in the tree, adds a
Cohen real (see \cite[2.1]{NeRo93}). This forcing may be represented as
$\q^{\tree}_1(K',\Sigma')$ for some t--omittory (not finitary) tree creating
pair $(K',\Sigma')$. 
\end{remark}

\section{Examples}

\begin{example}
\label{represent}
Let $\bH(m)=2$ for $m<\omega$. We construct a creating pair $(K_{
\ref{represent}},\Sigma_{\ref{represent}})$ for $\bH$, $\bar{t}\in\PC_\infty(
K_{\ref{represent}}, \Sigma_{\ref{represent}})$ and $\bar{t}$-systems $W^n_L$
(for $n,L<\omega$) such that
\begin{enumerate}
\item $(K_{\ref{represent}},\Sigma_{\ref{represent}})$ is simple except
omitting, finitary, forgetful, monotonic, interesting, condensed and generates
an ultrafilter, the systems $W^n_L$ are omittory--compatible,  
\item if $\Gamma\subseteq\p^*_\infty(\bar{t},(K_{\ref{represent}},
\Sigma_{\ref{represent}}))$ is quasi-$W^n_L$-generic then $\D(\Gamma)$ is a
filter on $\bigcup\limits_{i\in\omega}\{i\}\times (i+1)$ such that for every
colouring $f:\bigcup\limits_{i\in\omega}[\{i\}\times(i+1)]^{\textstyle n}
\longrightarrow L$ there is a set $A\in\D(\Gamma)$ almost homogeneous for $f$,
\item if $\Gamma\subseteq\p^*_\infty(\bar{t},(K_{\ref{represent}},
\Sigma_{\ref{represent}}))$ is $\preceq$-directed and countably closed and
$\D(\Gamma)$ is a filter such that for every colouring $f:\bigcup\limits_{i\in
\omega}[\{i\}\times(i+1)]^{\textstyle n}\longrightarrow L$ ($n,L<\omega$)
there is a set $A\in\D(\Gamma)$ almost homogeneous for $f$ then $\Gamma$ is
quasi-$W^n_L$-generic for all $n,L$. 
\end{enumerate}
\end{example}

\begin{proof}[Construction] A creature $t\in \CR[\bH]$ is in
$K_{\ref{represent}}$ if for some non-empty subset $a_t$ of
$[m^t_{\dn},m^t_{\up})$ we have $\nor[t]=\log_2( |a_t|)$ and 
\[\val[t]=\{\langle u,v\rangle\in 2^{\textstyle m^t_{\dn}}\!\times\!
2^{\textstyle m^t_{\up}}\!: (\forall n\!\in\![m^t_{\dn},m^t_{\up}))(v(n)=1\
\Rightarrow\ n\in a_t)\}.\] 
We define $\Sigma_{\ref{represent}}$ by:

if $t_0,\ldots,t_n\in K_{\ref{represent}}$ are such that $m^{t_\ell}_{\up}=
m^{t_{\ell+1}}_{\dn}$ (for $\ell<n$) then 
\[\Sigma_{\ref{represent}}(t_0,\ldots,t_n)=\{s\in K_{\ref{represent}}:
m^{t_0}_{\dn}=m^s_{\dn}\ \&\ m^{t_n}_{\up}=m^s_{\up}\ \&\ (\exists\ell\leq n)
(a_s\subseteq a_{t_\ell})\}.\] 
It should be clear that $(K_{\ref{represent}},\Sigma_{\ref{represent}})$ is a
finitary, forgetful, monotonic and simple except omitting creating pair. It is
interesting as $\nor[t]>0$ implies $|a_t|>1$ (for $t\in K_{\ref{represent}}$).
By the definition of $\Sigma_{\ref{represent}}$, one easily shows that
$(K_{\ref{represent}},\Sigma_{\ref{represent}})$ is condensed and generates an
ultrafilter. Moreover, $(K_{\ref{represent}},\Sigma_{\ref{represent}})$ is
strongly finitary modulo $\sim_{\Sigma_{\ref{represent}}}$. Now, let
$\bar{t}=\langle t_i: i<\omega\rangle\in\PC_\infty(K_{\ref{represent}},
\Sigma_{\ref{represent}})$ be such that for $i\in\omega$:
\[m^{t_i}_{\dn}=\frac{i(i+1)}{2},\quad m^{t_i}_{\up}=\frac{(i+1)(i+2)}{2},
\quad a_{t_i}=[m^{t_i}_{\dn}, m^{t_i}_{\up})\]
(so $\nor[t_i]=\log_2(i+1)$). We will identify the interval $[m^{t_i}_{\dn},
m^{t_i}_{\up})$ with $\{i\}\times (i+1)$ (for $i\in\omega$). Thus, if $\Gamma
\subseteq\p^*_\infty(\bar{t},(K_{\ref{represent}},\Sigma_{\ref{represent}}))$
is quasi-$W$-generic for some $\bar{t}$--system $W$, then we may think of
$\D(\Gamma)$ as a filter on $\bigcup\limits_{i\in\omega}\{i\}\times(i+1)$. The
filter $\D(\Gamma)$ is interesting by \ref{facgamfil}(1). 

Fix $n,L<\omega$. For each $i\in\omega$ choose a mapping 
\[\psi_i: \omega^{\textstyle [m^{t_i}_{\dn},m^{t_i}_{\up})}\stackrel{\rm onto}
{\longrightarrow}\{f:\mbox{ $f$ is a function from $[[m^{t_i}_{\dn},
m^{t_i}_{\up})]^{\textstyle n}$ to $L\}$}.\]
Next, for $i\leq j<\omega$ and $\sigma: [m^{t_i}_{\dn}, m^{t_j}_{\up})
\longrightarrow\omega$ define
\[\begin{array}{ll}
W^n_L(m^{t_i}_{\dn}, m^{t_j}_{\up},\sigma)=&\{t\in\Sigma_{\ref{represent}}(t_i,
\ldots,t_j):\mbox{if }i\leq k\leq j,\ a_t\subseteq a_{t_k}\\
\ &\qquad\mbox{ then $a_t$ is homogeneous for }\psi_k(\sigma\rest[m^{t_k}_{
\dn},m^{t_k}_{\up}))\}\end{array}\]
(in all other instances we let $W^n_L(m',m'',\sigma)=\emptyset$).

\begin{claim}
\label{cl32}
$W^n_L$ is an omittory--compatible $\bar{t}$--system.
\end{claim}

\noindent{\em Proof of the claim:}\ \ \ The requirement \ref{quasi}(1a) is
immediate by the definition of $W^n_L$. For \ref{quasi}(1b,c) remember the way
we defined the composition operation $\Sigma_{\ref{represent}}$: if $s\in
\Sigma_{\ref{represent}}(s_0,\ldots,s_k)$ then $a_s\subseteq a_{s_\ell}$ for
some $0\leq\ell\leq k$. Finally note that if $s\in\Sigma_{\ref{represent}}(
t_k,\ldots,t_\ell)$, $k\leq\ell<\omega$, $\nor[s]>\log_2(R_n(L,m))$ (see
\ref{ramnum}) and $\sigma: [m^s_{\dn},m^s_{\up})\longrightarrow\omega$ then
there is $t\in W^n_L(m^s_{\dn},m^s_{\up},\sigma)$ such that $t\in
\Sigma_{\ref{represent}}(s)$ and $\nor[t]\geq\log_2(m)$ (by the definition of
$R_n(L,m)$). This gives the suitable for \ref{quasi}(1d) function $G$. Thus we
have verified that $W^n_L$ is a $\bar{t}$--system. It should be clear that $W$
is omittory--compatible.
\medskip

Note that if $\Gamma\subseteq\p^*_\infty(\bar{t},(K_{\ref{represent}},
\Sigma_{\ref{represent}}))$ is $\preceq$--directed then $\D(\Gamma)$ is the
filter generated by all sets $A^N_{\bar{s}}\stackrel{\rm def}{=}\bigcup
\limits_{N\leq n}a_{s_n}$, for $\bar{s}=\langle s_n: n<\omega\rangle\in\Gamma$
and $N<\omega$. Therefore we easily check that the systems $W^n_L$ (for $n,L<
\omega$) are as required in 2, 3 of \ref{represent}, noting that each
colouring $f:\bigcup\limits_{i\in\omega}[\{i\}\times(i+1)]^{\textstyle n}
\longrightarrow L$ corresponds via $\langle \psi_i: i<\omega\rangle$ to some
function $\eta\in\baire$. 
\end{proof}
 
\begin{proposition}
\label{repexists}
Assume CH. Then there exist $\Gamma\subseteq\p^*_\infty(\bar{t},(K_{
\ref{represent}},\Sigma_{\ref{represent}}))$ which is quasi-$W^n_L$-generic
for all $n,L<\omega$ and such that $\D(\Gamma)$ is a semi--Ramsey ultrafilter
on $\bigcup\limits_{i\in\omega}\{i\}\times (i+1)$. Consequently, $\D(\Gamma)$
is an interesting almost Ramsey ultrafilter on $\bigcup\limits_{i\in\omega}
\{i\}\times (i+1)$. 
\end{proposition}

\begin{proof}
This is somewhat similar to \ref{facgamfil}(3) (and so to \ref{exists}(2)),
but we have to be more careful to ensure that $\D(\Gamma)$ is semi--Ramsey. 
For this, as a basic step of the inductive construction of $\Gamma$, we use
the following observation. 

\begin{claim}
\label{cl33}
Suppose that $\bar{s}\in\p^*_\infty(\bar{t},(K_{\ref{represent}},
\Sigma_{\ref{represent}}))$, $n,L<\omega$, $\eta\in\baire$ and $\varphi:
\omega\longrightarrow [\bigcup\limits_{i\in\omega}\{i\}\times(i+1)]^{
\textstyle \omega}$. Then there exists $\bar{s}^*=\langle s^*_m:m<\omega
\rangle\in\p^*_\infty(\bar{t},(K_{\ref{represent}},\Sigma_{\ref{represent}}))$
such that $\bar{s}\leq\bar{s}^*$ and
\begin{enumerate}
\item $(\forall^\infty m)(s^*_m\in W^n_L(m^{s^*_m}_{\dn}, m^{s^*_m}_{\up},
\eta\rest [m^{s^*_m}_{\dn},m^{s^*_m}_{\up})))$,
\item $(\forall m\in\omega)(|a_{s^*_m}|\leq m+1)$,
\item if $i^*_m<\omega$ (for $m\in\omega$) are such that $a_{s^*_m}\subseteq
\{i^*_m\}\times (i^*_m+1)\cong [m^{t_{i^*_m}}_{\dn},m^{t_{i^*_m}}_{\up})$ then
\[\mbox{either }\ (\exists k\in\omega)(\forall m\in\omega)(a_{s^*_m}\cap
\varphi(k)=\emptyset)\quad\mbox{ or }\quad (\forall m\in\omega)(a_{s^*_{m+1}}
\subseteq\varphi(i^*_m)).\]
\end{enumerate}
\end{claim}

\noindent{\em Proof of the claim:}\ \ \ Let $\bar{s}=\langle s_m:
m<\omega\rangle$ and let $i_m$ be such that $a_{s_m}\subseteq \{i_m\}\times
(i_m+1)$ (remember $\bar{t}\leq\bar{s}$, see the definition of
$(K_{\ref{represent}},\Sigma_{\ref{represent}})$). We know that
$\lim\limits_{m\to\infty} |a_{s_m}|=\infty$, so we may choose $m_0<m_1<m_2<
\ldots<\omega$ such that $(\forall k\in\omega)(|a_{s_{m_k}}|\geq R_n(L,k))$
(see \ref{ramnum}). Choose $b_k\in [a_{s_{m_k}}]^{\textstyle k}$ homogeneous
for the colouring of $[\{i_{m_k}\}\times (i_{m_k}+1)]^{\textstyle n}$ (with
values in $L$) coded by $\eta\rest [m^{t_{i_{m_k}}}_{\dn},m^{t_{i_{m_k}}}_{
\up})$. Next choose $k(0)<k(1)<\ldots<\omega$ and $c_\ell\in [b_{k(\ell)}]^{
\textstyle \ell+1}$ (for $\ell<\omega$) such that
\[\mbox{either }\ (\exists k\in\omega)(\forall \ell\in\omega)(c_\ell\cap
\varphi(k)=\emptyset)\quad\mbox{ or }\quad (\forall m\in\omega)(c_{\ell+1}
\subseteq\varphi(i_{m_{k(\ell)}})).\]
Let $s^*_0\in K_{\ref{represent}}$ be such that $m^{s^*_0}_{\dn}=m^{s_0}_{
\dn}$, $m^{s^*_0}_{\up}=m^{s_{m_{k(0)}}}_{\up}$, $a_{s^*_0}=c_0$ and let
$s^*_{\ell+1}\in K_{\ref{represent}}$ (for $\ell\in\omega$) be such that
$m^{s^*_{\ell+1}}_{\dn}=m^{s_{m_{k(\ell)}}}_{\up}$, $m^{s^*_{\ell+1}}_{\up}=
m^{s_{m_{k(\ell+1)}}}_{\up}$ and $a_{s^*_{\ell+1}}=c_{\ell+1}$. Easily, the
sequence $\bar{s}^*=\langle s^*_n: n<\omega\rangle$ is as required.
\medskip

Assume CH. Using \ref{cl33}, we may construct a sequence $\langle
\bar{s}_\alpha: \alpha<\omega_1\rangle\subseteq\p^*_\infty(\bar{t},(
K_{\ref{represent}},\Sigma_{\ref{represent}}))$ such that
\begin{enumerate}
\item[($\alpha$)]  $\alpha<\beta<\omega_1\qquad\Rightarrow\qquad\bar{s}_\alpha
\preceq\bar{s}_\beta$,
\item[($\beta$)]   for each $n,L<\omega$ we have
\[(\forall\eta\in\baire)(\exists\alpha<\omega_1)(\forall^\infty m)(s_{\alpha,
m}\in W^n_L(m^{s_{\alpha,m}}_{\dn},m^{s_{\alpha,m}}_{\up},\eta\rest
[m^{s_{\alpha,m}}_{\dn},m^{s_{\alpha,m}}_{\up}))),\]
\item[($\gamma$)]  for each function $\varphi:\omega\longrightarrow
[\bigcup\limits_{i\in\omega}\{i\}\times (i+1)]^{\textstyle \omega}$ there is
$\alpha<\omega_1$ such that\\
if $a_{s_{\alpha,m}}\subseteq\{i_{\alpha,m}\}\times (i_{\alpha,m}+1)$ (for
$m\in\omega$) then
\[\mbox{either }\ (\exists k\in\omega)(\forall m\in\omega)(a_{s_{\alpha,m}}
\cap\varphi(k)=\emptyset)\quad\mbox{ or }\quad (\forall m\in\omega)(
a_{s_{\alpha,m+1}}\subseteq\varphi(i_{\alpha,m})).\]
\end{enumerate}
Like in \ref{semnotalm} and \ref{exists}(2) we check that $\Gamma\stackrel{\rm
def}{=}\{\bar{s}_\alpha:\alpha<\omega_1\}$ is quasi-$W^n_L$-generic for all
$n,L<\omega$ and $\D(\Gamma)$ is a semi--Ramsey ultrafilter on
$\bigcup\limits_{i\in\omega}\{i\}\times(i+1)$. 
\end{proof}

\begin{conclusion}
\label{limprealm}
Assume CH. Let $(K_{\ref{represent}},\Sigma_{\ref{represent}})$, $\bar{t}$,
$W^n_L$ be given by \ref{represent}, and let $\Gamma\subseteq\p^*_\infty(
\bar{t},(K_{\ref{represent}},\Sigma_{\ref{represent}}))$ be
quasi-$W^n_L$-generic for all $n,L<\omega$ such that $\D(\Gamma)$ is a
semi--Ramsey ultrafilter on $\bigcup\limits_{i\in\omega}\{i\}\times (i+1)$
(see \ref{repexists}). Suppose that $\delta$ is a limit ordinal and $\langle
\p_\alpha,\dot{\q}_\alpha: \alpha<\delta\rangle$ is a countable support
iteration of proper $\baire$--bounding forcing notions such that for each
$\alpha<\delta$:  
\[\forces_{\p_\alpha}\mbox{``$\Gamma$ generates an interesting almost Ramsey
ultrafilter''.}\] 
Then $\forces_{\p_\delta}$``$\Gamma$ generates an interesting almost Ramsey
ultrafilter''.  
\end{conclusion}

\begin{proof}
By \ref{represent}(3) we have
\[\forces_{\p_\alpha}\mbox{``$\Gamma$ is quasi-$W^n_L$-generic for all
$n,L<\omega$''}\] 
(for each $\alpha<\delta$). Hence, by \ref{presextra}, we get 
\[\forces_{\p_\delta}\mbox{``$\Gamma$ is quasi-$W^n_L$-generic for each 
$n,L<\omega$ and generates an ultrafilter''.}\] 
As $\p_\delta$ is $\baire$--bounding (by \cite[Ch VI, 2.3, 2.8]{Sh:f}) we may
apply \ref{almboupre} to conclude that 
\[\forces_{\p_\delta}\mbox{``$\D(\Gamma)$ is a semi--Ramsey ultrafilter''.}\]
Consequently, by \ref{represent}(2), we have
\[\forces_{\p_\delta}\mbox{``$\D(\Gamma)$ is an interesting almost Ramsey
ultrafilter''.}\] 
\end{proof}

\begin{example}
\label{matet}
Let $\psi\in\baire$ be such that $(\forall n\in\omega)(\psi(n)>(n+1)^2)$.\\
We build a tree creating pair $(K_{\ref{matet}}^\psi,\Sigma_{\ref{matet}}^\psi
)$ which is: finitary, $2$-big, rich (see \ref{rich}) and of the $\sUP(\D)^{
\tree}$-type (see \ref{up}) for every Ramsey ultrafilter $\D$ on $\omega$.  
\end{example}

\begin{proof}[Construction] 
Let $\bH(n)=[\psi(n)]^{\textstyle n+1}$.

For $\nu\in\prod\limits_{m<n}\bH(m)$ ($n\in\omega$) and $A\subseteq
\bigcup\limits_{m<\omega}\{m\}\times \psi(m)$ we will write $A\prec\nu$ if
\[(\forall m<n)(\forall k< \psi(m))((m,k)\in A\ \ \ \Rightarrow\ \ \
k\in\nu(m)).\]

Now we define $(K_{\ref{matet}}^\psi,\Sigma_{\ref{matet}}^\psi)$. A tree--like
creature $t\in\TCR_\eta[\bH]$ is in $K_{\ref{matet}}^\psi$ if: 
\begin{enumerate}
\item $\val[t]$ is finite and
\item $\nor[t]=\log_2(\min\{|A|\!: A\subseteq\bigcup\limits_{m\geq\lh(\eta)}
\{m\}\times\psi(m)\ \ \&\ \ (\forall\nu\in\pos(t))(A\not\prec\nu)\})$. 
\end{enumerate}
By the definition, $K_{\ref{matet}}^\psi$ is finitary. The tree composition
$\Sigma_{\ref{matet}}^\psi$ is generated similarly to  $\Sigma^{\tsum}$ of
\ref{treesum} but with norms as above. Thus, if $\langle t_\nu: \nu\in
\hat{T}\rangle\subseteq K^\psi_{\ref{matet}}$ is a system of tree-creatures
such that $T$ is a well founded quasi tree, $\mrot(t_\nu)=\nu$, and
$\rng(\val[t_\nu])=\suc_T(\nu)$ (for $\nu\in \hat{T}$) then we define
$S^*(t_\nu:\nu\in \hat{T})$ as the unique creature $t^*$ in $K^\psi_{
\ref{matet}}$ with $\rng(\val[t^*])=\max(T)$, $\dom(\val[t^*])=\{\mrot(T)\}$
and $\dis[t^*]=\langle\dis[t_\nu]:\nu\in\hat{T}\rangle$. Now we let   
\[\Sigma^\psi_{\ref{matet}}(t_\nu:\nu\in\hat{T})=\{t\in K^\psi_{\ref{matet}}:
\val[t]\subseteq\val[S^*(t_\nu:\nu\in\hat{T})]\}.\]
Clearly, $\Sigma_{\ref{matet}}^\psi$ is a tree-composition on
$K_{\ref{matet}}^\psi$. Note that if $t\in K^\psi_{\ref{matet}}$ then we may
identify elements of $\Sigma^\psi_{\ref{matet}}(t)$ with subsets of $\pos(t)$:
for each non-empty $u\subseteq\pos(t)$, $t^u$ is the unique creature in
$K^\psi_{\ref{matet}}$ with $\pos(t^u)=u$ and $\mrot(t^u)=\mrot(t)$. More
general, if $p\in\q^{\tree}_\emptyset(K^\psi_{\ref{matet}},\Sigma^\psi_{
\ref{matet}})$, $F$ is a front of $T^p$, $u\subseteq F$ then there is a unique
creature $t(p,u)\in K^\psi_{\ref{matet}}$ such that
\[\pos(t(p,u))\!=\! u,\ \ \mrot(t(p,u))\!=\!\mrot(p)\ \mbox{ and }\ t(p,u)\in
\Sigma^\psi_{\ref{matet}}(t^p_\nu\!: (\exists\eta\!\in\! F)(\nu
\vartriangleleft\eta)).\]

\begin{claim}
\label{cl16}
$\nor[S^*(t_\nu:\nu\in\hat{T})]\geq\min\{\nor[t_\nu]:\nu\in\hat{T}\}$.
\end{claim}

\noindent{\em Proof of the claim:}\ \ \ Suppose that
$m<2^{\min\{\nor[t_\nu]: \nu\in\hat{T}\}}$, but there is a set
$A\subseteq\bigcup\limits_{k\geq n_0}\{k\}\times \psi(k)$ (where
$n_0=\lh(\mrot(T))$) such that $|A|=m$ and $A\not\prec\nu$ for each
$\nu\in\max(T)$. Now we build inductively a bad $\nu\in\max(T)$: since
$m<2^{\nor[t_{\mrot(T)}]}$ we find $\nu_0\in\pos(t_{\mrot(T)})$ such that
$A\prec\nu_0$. Next we look at $A^1=A\cap\bigcup\limits_{k\geq\lh(\nu_0)}
\{k\}\times \psi(k)$. Since $m<2^{\nor[t_{\nu_0}]}$ we find
$\nu_1\in\pos(t_{\nu_0})$ such that $A^1\prec\nu_1$. Continuing in this
fashion, after finitely many steps, we get $\nu_k\in\max(T)$ such that $A\prec
\nu_k$, a contradiction.

\begin{claim}
\label{cl12}
$(K^\psi_{\ref{matet}},\Sigma^\psi_{\ref{matet}})$ is $2$-big.
\end{claim}

\noindent{\em Proof of the claim:}\ \ \ Let $t\in K^\psi_{\ref{matet}}$,
$\nor[t]>0$ and let $\pos(t)=u_0\cup u_1$. Take sets $A_0,A_1\subseteq
\bigcup\limits_{m\geq n_0}\{m\}\times \psi(m)$ (where $n_0=\lh(\mrot(t))$)
such that for $i=0,1$:
\[(\forall\nu\in u_i)(A_i\not\prec\nu)\quad\mbox{ and }\quad \log_2(|A_i|)
=\nor[t^{u_i}].\] 
Look at $A=A_0\cup A_1$. Clearly $(\forall\nu\in\pos(t))(A\not\prec\nu)$.
Hence, for some $i<2$:
\[\nor[t]\leq\log_2(|A|)\leq\log_2(|A_i|)+1=\nor[t^{u_i}]+1.\]

\begin{claim}
\label{cl17}
$(K^\psi_{\ref{matet}},\Sigma^\psi_{\ref{matet}})$ is rich.
\end{claim}

\noindent{\em Proof of the claim:}\ \ \ Suppose that
$\langle s_\nu:\nu\in\hat{T}\rangle\subseteq K^\psi_{\ref{matet}}$,
$n\in\omega$ and $u$ are such that
\begin{enumerate}
\item $T\subseteq\bigcup\limits_{k\in\omega}\prod\limits_{m<k}\bH(m)$ is a
well founded quasi tree, $u\subseteq \max(T)$,
\item $\mrot(s_\nu)=\nu$, $\pos(s_\nu)=\suc_T(\nu)$, $\nor[s_\nu]>n+3$,
\item there is no system $\langle s^*_\nu: \nu\in\hat{T}^*\rangle\subseteq
K^\psi_{\ref{matet}}$ such that
\[\begin{array}{l}
T^*\subseteq T,\quad \max(T^*)\subseteq u,\quad \mrot(T^*)=\mrot(T),
\quad\pos(s^*_\nu)=\suc_{T^*}(\nu)\\
\mrot(s^*_\nu)=\nu,\quad\nor[s^*_\nu]>n+1,\quad\mbox{ and}\\
s^*_\nu\in\Sigma^\psi_{\ref{matet}}(s_\eta:\eta\in\hat{T}_\nu)\quad\mbox{ for
some }T_\nu\subseteq T.
  \end{array}\]
\end{enumerate}
Let $t=S^*(s_\nu: \nu\in\hat{T})$. By \ref{cl16} we have $\nor[t]>n+3$
(remember 2. above). But now, considering $t^u$ (defined as before) we
note that necessarily $\nor[t^u]\leq n+1<\nor[t]-1$. Let $v=\pos(t)\setminus
u$. By the bigness (see \ref{cl17}) we have $\nor[t^v]\geq\nor[t]-1$,
finishing the claim (remember \ref{cl16}). 

\begin{claim}
\label{cl13}
Let $\D$ be a Ramsey ultrafilter on $\omega$. Then $(K^\psi_{\ref{matet}},
\Sigma^\psi_{\ref{matet}})$ is of the $\sUP(\D)^{\tree}$--type.
\end{claim}

\noindent{\em Proof of the claim:}\ \ \ Assume that $1\leq m<\omega$, $p\in
\q^*_\emptyset(K^\psi_{\ref{matet}},\Sigma^\psi_{\ref{matet}})$,
$\nor[t^p_\nu]>m+1$ for each $\nu\in T^p$ and $F_0, F_1,\ldots$ are fronts of
$T^p$ such that 
\[(\forall n\in\omega)(\forall\nu\in F_{n+1})(\exists\eta\in F_n)(\eta
\vartriangleleft\nu).\]
Further suppose that $u_n\subseteq F_n$ are such that there is no system
$\langle s_\nu: \nu\in\hat{T}\rangle$ with 
\[\pos(S^*(s_\nu:\nu\in \hat{T}))\subseteq u_n,\quad \mrot(T)=\mrot(T^p),\quad
\mbox{ and }\nor[s_\nu]\geq m.\]
In particular, this means that $\nor[t(p,u_n)]<m$ for each $n\in\omega$
($t(p,u_n)$ is as defined earlier). Thus, for each $n\in\omega$, we find a set
\[A_n\subseteq\bigcup\limits_{k\geq\lh(\mrot(p))}\{k\}\times \psi(k)\]
such that
\[|A_n|=2^m\quad\mbox{ and }\quad (\forall \nu\in u_n)(A_n\not\prec \nu).\]
Now we use the assumption that $\D$ is Ramsey:\quad we find sets $Z_0\in\D$
and $A^*\subseteq\bigcup\limits_{i\in\omega}\{i\}\times \psi(i)$ such that 
\[(\forall n_0,n_1\in Z_0)(n_0<n_1\quad\Rightarrow\quad A_{n_0}\cap A_{n_1}
=A^*).\]
[How? Just consider the colouring $f$ of $[\omega]^{\textstyle 2}$ such that
$f(n_0,n_1)$ codes the trace of $A_{n_0}$ on $A_{n_1}$ (for $n_0<n_1$, in the
canonical enumerations of $A_n$'s) and take an $f$-homogeneous set.] Now
choose $Z\subseteq Z_0$, $Z\in\D$ such that if $n_0<n_1$, $n_0,n_1\in Z$ then
\[\min\{i: (A_{n_1}\setminus A^*)\cap\{i\}\times \psi(i)\neq\emptyset\}
>\max\{\lh(\nu): \nu\in F_{n_0}\}.\]
[How? Consider the following strategy for Player I in the game $G^R(\D)$:\\
at stage $k+1$ of the game he looks at the last move $i_k$ of the second
player and he chooses $N$ such that if $n\geq N$, $n\in Z_0$ then
\[(A_n\setminus A^*)\cap \bigcup\big\{\{i\}\times\psi(i): i\leq\max\{\lh(\nu):
\nu\in F_{i_k}\}\big\}=\emptyset.\] Now he plays $Z_0\cap (N,\omega)$. This
strategy cannot be the winning one.] Using the set $Z$ we build the suitable
condition $q\in\q^{\tree}_\emptyset (K^\psi_{\ref{matet}},\Sigma^\psi_{
\ref{matet}})$. It will be constructed in such a way that $p\leq_0^1 q$,
$T^q\subseteq\{\mrot(q)\}\cup\bigcup\limits_{\ell\in Z} F_\ell$ and each
$F_\ell\cap T^q$ will be a front of $T^q$ (for $\ell\in Z$). Let $\ell_0=\min
Z$ and $v_{\ell_0}=\{\nu\in F_{\ell_0}: A_{\ell_0}\prec\nu\}$. By the choice
of $A_{\ell_0}$ we know that $v_{\ell_0}\cap u_{\ell_0}=\emptyset$. Let
$t^q_{\mrot(q)}=t(p,v_{\ell_0})$. Note that $2^{\nor[t(p,v_{\ell_0})]}+2^m\geq
2^{\nor[t(p,F_{\ell_0})]}$, and hence, as $\nor[t(p,F_{\ell_0})]>m+1$ (see
\ref{cl16}), we have 
\[\nor[t(p,v_{\ell_0})]\geq \nor[t(p,F_{\ell_0})]-m\geq\min\{\nor[t^p_\nu]:
\mrot(p)\trianglelefteq\nu\in T^p\}-m.\]
We put $t^q_{\mrot(q)}=t(p,v_{\ell_0})$. Suppose that we have defined $T^q$ up
to the level of $F_\ell$, $\ell\in Z$ (thus we know $T^q\cap F_\ell$ already).
Let $\eta\in T^q\cap F_\ell$ and let $\ell'=\min(Z\setminus(\ell+1))$. By the
first step of the construction we know that $A^*\prec \eta$. By the choice of
$Z$ we have that 
\[(A_{\ell'}\setminus A^*)\cap\bigcup_{i<\lh(\eta)}\{i\}\times
\psi(i)=\emptyset.\]
We take $v_{\eta,\ell'}=\{\nu\in F_{\ell'}: \eta\vartriangleleft\nu\ \&\
A_{\ell'}\prec \nu\}$. As before, $v_{\eta,\ell'}\cap u_{\ell'}=\emptyset$ and
\[\nor[t(p^{[\eta]}, v_{\eta,\ell'})]\geq \nor[t(p^{[\eta]},F_{\ell_0}\cap
T^{p^{[\eta]}})]-m\geq\min\{\nor[t^p_\nu]: \eta\trianglelefteq\nu\in T^p\}-m\]
(remember \ref{cl16}). 

Now we easily check that the condition $q$ constructed above is as required by
\ref{up} to show that $(K^\psi_{\ref{matet}},\Sigma^\psi_{\ref{matet}})$ is of
the $\sUP(\D)^{\tree}$--type.
\end{proof}

One easily checks that the forcing notion $\q^{\tree}_1(K^\psi_{\ref{matet}},
\Sigma^\psi_{\ref{matet}})$ is non-trivial (i.e.~$K^\psi_{\ref{matet}}$
contains enough tree like creatures with arbitrarily large norms). Let us show
another property of the tree creating pair $(K^\psi_{\ref{matet}},
\Sigma^\psi_{\ref{matet}})$. 

\begin{proposition}
\label{addpro}
Let $(K^\psi_{\ref{matet}},\Sigma^\psi_{\ref{matet}})$ be the tree creating
pair of \ref{matet}. Then:
\begin{enumerate}
\item $\forces_{\q^{\tree}_1(K^\psi_{\ref{matet}},\Sigma^\psi_{\ref{matet}})}
(\forall m\in\omega)(\dot{W}(m)\in [\psi(m)]^{\textstyle m+1})$ and
\item if $h$ is a partial function, $\dom(h)\in\iso$ and $(\forall m\in
\dom(h))(h(m)<\psi(m))$ then
\[\forces_{\q^{\tree}_1(K^\psi_{\ref{matet}},\Sigma^\psi_{\ref{matet}})}
(\exists^\infty m\in\dom(h))(h(m)\in\dot{W}(m))\]
\end{enumerate}
(where $\dot{W}$ is the generic real added by
$\q^{\tree}_1(K^\psi_{\ref{matet}},\Sigma^\psi_{\ref{matet}})$, see
\ref{thereal}). 
\end{proposition}

\begin{proof}
1)\ \ \ Should be clear. 

\noindent 2)\ \ \ Let
$p\in\q^{\tree}_1(K^\psi_{\ref{matet}},\Sigma^\psi_{\ref{matet}})$ and let $h$
be as in the assumptions. We may assume that $(\forall \eta\in T^p)(\nor[
t^p_\eta]>5)$. Let $m_0=\min\big(\dom(h)\setminus\lh(\mrot(p))\big)$. Take 
a front $F$ of $T^p$ such that $(\forall\eta\in F)(\lh(\eta)>m_0)$ and look at
the set $u=\{\eta\in F: h(m_0)\in\eta(m_0)\}$. Plainly, $\nor[t(p,u)]\geq
\nor[t(p,F)]-1\geq 4$. Let $q$ be a condition in $\q^{\tree}_1(K^\psi_{
\ref{matet}},\Sigma^\psi_{\ref{matet}})$ such that $\mrot(q)=\mrot(p)$, $T^q
\subseteq T^p$, $t^q_{\mrot(q)}=t(p,u)$ and if $\eta\in u$, $\eta
\trianglelefteq\nu\in T^p$ then $\nu\in T^q$, $t^q_\nu=t^p_\nu$. Clearly
$q\geq p$ and $q\forces h(m_0)\in\dot{W}(m_0)$. Now we may easily finish.
\end{proof}

\begin{conclusion}
\label{conmatet}
The following is consistent with ZFC:
\begin{enumerate}
\item there is an almost Ramsey interesting ultrafilter on $\bigcup\limits_{i
\in\omega} \{i\}\times (i+1)$ which is generated by $\aleph_1$ elements (so
${\mathfrak m}_2=\lambda=\aleph_1$, see the introduction to this chapter) and  
\item ${\mathfrak m}_1=\aleph_2$, and even more: for each function $\psi\in
\baire$ and a family ${\mathcal F}$ of $\aleph_1$ partial infinite functions
$h:\dom(h)\longrightarrow\omega$ such that 
\[(\forall m\in\dom(h))(h(m)<\psi(m))\]
{\em there is} $W\in \prod\limits_{m\in\omega}[\psi(m)]^{\textstyle m+1}$ such
that  
\[(\forall h\in{\mathcal F})(\exists^\infty m\in\dom(h))(h(m)\in W(m)).\]
\end{enumerate}
\end{conclusion}

\begin{proof}
Start with $\V\models$ CH. By \ref{repexists} we have an interesting almost
Ramsey ultrafilter $\D=\D(\Gamma)$ on $\bigcup\limits_{i\in\omega}\{i\}\times
(i+1)$ generated by a quasi generic $\Gamma\subseteq\p^*_\infty(\bar{t},(
K_{\ref{represent}},\Sigma_{\ref{represent}}))$ as there. Build a countable
support iteration $\langle\p_\alpha,\dot{\q}_\alpha:\alpha<\omega_2\rangle$
and a list $\langle\dot{\psi}_\alpha: \alpha<\omega_2\rangle$ such that for
each $\alpha<\omega_2$: 
\begin{enumerate}
\item $\dot{\psi}_\alpha$ is a $\p_\alpha$--name for a function in $\baire$
such that $(\forall n\in\omega)(\dot{\psi}_\alpha(n)>(n+1)^2)$,
\item $\dot{\q}_\alpha$ is the $\p_\alpha$--name for the forcing notion
$\q^{\tree}_1(K^{\dot{\psi}_\alpha}_{\ref{matet}},
\Sigma^{\dot{\psi}_\alpha}_{\ref{matet}})$,
\item $\langle\dot{\psi}_\beta:\beta<\omega_2\rangle$ lists with
$\omega_2$--repetitions all (canonical) $\p_{\omega_2}$--names for functions
$\psi\in\baire$ such that $(\forall n\in\omega)(\psi(n)>(n+1)^2)$. 
\end{enumerate}
We claim that if $G\subseteq \p_{\omega_2}$ is a generic filter over $\V$,
then, in $\V[G]$, the two sentences of the conclusion hold true. Why? One can
inductively show that for each $\alpha\leq\omega_2$:
\[\forces_{\p_\alpha}\mbox{``$\Gamma$ generates an almost Ramsey interesting
ultrafilter on }\bigcup\limits_{i\in\omega}\{i\}\times(i+1)\mbox{''}\]
(at successor stages use \ref{prealmult}, \ref{prealmram} and \ref{matet}; at
limit stages use \ref{limprealm}). Hence, in  $\V[G]$, the ultrafilter
$\D(\Gamma)$ witnesses the first property. For the second assertion use
\ref{addpro}.
\end{proof}

\chapter{Friends and relatives of PP} 
In this chapter we answer a question of Balcerzak and Plewik, showing that the
cardinal number $\kappa_{\rm BP}$ (see \ref{BPnumber}) may be smaller than the
continuum (\ref{BPbelcon}) and that it may be larger than the dominating
number (\ref{BPabod}). As this cardinal turns out to be bounded by a cardinal
number related to the strong PP--property, we take this opportunity to have a
look at several properties close to the PP--property.

\section{Balcerzak--Plewik number}
For an ideal $\J$ of subsets of $\can$ it is natural to ask if it has the
following property {\bf (P)}:
\begin{enumerate}
\item[({\bf P})$_{\J}$] \quad every perfect subset of $\can$ contains a
perfect set from $\J$.
\end{enumerate}
The property {\bf (P)} has numerous consequences and applications (see e.g.\
Balcerzak \cite{Ba91}, some related results and references may be found in
Balcerzak Ros{\l}anowski \cite{BaRo95}) and it is usually easy to decide if
${\bf (P)}_{\J}$ holds. However, that was not clear for some of Mycielski's
ideals $\M^*_{2,\K}$.  

Suppose that $\K\subseteq\iso$ is a non-empty family such that
\begin{enumerate}
\item[$(\oplus)$] \quad $(\forall X\in\K)(\exists X_0,X_1\in\K)(X_0,X_1
\subseteq X\ \ \&\ \ X_0\cap X_1=\emptyset)$.
\end{enumerate}
Let $\M^*_{2,\K}$ consist of these sets $A\subseteq\can$ that
\[(\forall X\in\K)(\exists f:X\longrightarrow 2)(\forall g\in A)(\neg f
\subseteq g).\]
It is easy to check that $\M^*_{2,\K}$ is a $\sigma$--ideal of subsets of
$\can$. These ideals are relatives of the ideals from Mycielski \cite{My69}
and were studied e.g.\ in Cicho\'n Ros{\l}anowski Stepr\={a}ns W{\c e}glorz
\cite{CRSW93} and \cite{Ro94}. If the family $\K$ is countable than easily the
ideal $\M^*_{2,\K}$ determined by it has the property {\bf (P)}. D{\c e}bski,
Kleszcz and Plewik \cite{DKP92} showed that 
the ideal $\M^*_{2,\iso}$ does not satisfy {\bf (P)}. Then Balcerzak and
Plewik defined the following cardinal number $\kappa_{\rm BP}$ (see
Balcerzak Plewik \cite{BaPl96}). 

\begin{definition}
\label{BPnumber}
The Balcerzak--Plewik number $\kappa_{\rm BP}$ is the minimal size of a family
$\K\subseteq\iso$ such that
\begin{quotation}
\noindent for some perfect set $Q\subseteq\can$, for every perfect subset $P$
of $Q$ there is $X\in\K$ such that
\[P\rest X\stackrel{\rm def}{=}\{f\in 2^{\textstyle X}: (\exists g\in P)(f
\subseteq g)\}=2^{\textstyle X}.\]
\end{quotation}
[Note that $\kappa_{\rm BP}$ is the minimal size of $\K\subseteq\iso$
satisfying $(\oplus)$ for which the ideal $\M^*_{2,\K}$ does not have the
property {\bf (P)}.]
\end{definition}
They proved that $\dominating\leq\kappa_{\rm BP}$ and asked if
\begin{itemize}
\item it is consistent that $\kappa_{\rm BP}<\con$,
\item it is consistent that $\dominating<\kappa_{\rm BP}$.
\end{itemize}
A full answer to these questions will be given in the final part of this
chapter. Now we want to give an upper bound to $\kappa_{\rm BP}$.

\begin{definition}
\label{sPPrel}
Let ${\mathcal X}$ be the space of all sequences $\bar{w}=\langle w_i: i\in
B\rangle$ such that $B\in\iso$ and $(\forall i\in B)(w_i\in
[\omega]^{\textstyle i})$. We define a relation $R^{\rm sPP}\subseteq\baire
\times{\mathcal X}$ by
\[(\eta,\bar{w})\in R^{\rm sPP}\quad\mbox{ if and only if }\quad \eta\in
\baire,\ \bar{w}\in {\mathcal X},\mbox{ and }(\forall i\in\dom(\bar{w}))
(\eta(i)\in w_i).\]
\end{definition}
Note that the space $\mathcal X$ carries a natural Polish topology (inherited
from the product space of all $\bar{w}=\langle w_i: i<\omega\rangle$ such that
for each $i\in\omega$, either $w_i=\emptyset$ or $|w_i|=i$). The relation
$R^{\rm sPP}$ describes the strong PP--property of \cite[Ch VI, 2.12E]{Sh:f}:
a proper forcing notion $\p$ has the strong PP--property if and only if it has
the $R^{\rm sPP}$--localization property (see \ref{gendef}). 

\begin{theorem}
\label{kbelPP}
$\kappa_{\rm BP}\leq\dominating(R^{\rm sPP})$.
\end{theorem}

\begin{proof}
Construct inductively a perfect tree $T\subseteq\fs$ and an increasing
sequence $0=k_0<k_1<k_2<\ldots<\omega$ such that for every $i\in\omega$: 
\begin{enumerate}
\item[$(\alpha)$] $(\forall\nu\in T\cap 2^{\textstyle k_i})(\exists \eta_0,
\eta_1\in T\cap 2^{\textstyle k_{i+1}})(\nu\vartriangleleft\eta_0\ \&\
\nu\vartriangleleft\eta_1\ \&\ \eta_0\neq\eta_1)$,
\item[$(\beta)$] for each colouring $f:T\cap 2^{\textstyle k_i}\longrightarrow
2$ there is $n\in [k_i,k_{i+1})$ such that 
\[(\forall \eta\in T\cap 2^{\textstyle k_{i+1}})(\eta(n)=f(\eta\rest k_i)).\]
\end{enumerate}
The construction is straightforward. It is not difficult to check that if
$P\subseteq [T]$ is a perfect set then there is $X\in\iso$ such that $P\rest
X=2^{\textstyle X}$ (or see \cite{DKP92}). 

Let $D\subseteq{\mathcal X}$ be such that $|D|=\dominating(R^{\rm sPP})$ and
\[(\forall\eta\in\baire)(\exists\bar{w}\in D)((\eta,\bar{w})\in R^{\rm
sPP}).\] 
Let $N$ be an elementary submodel of ${\mathcal H}(\chi)$ such that $D,T,
\langle k_i: i<\omega\rangle\in N$, $D\subseteq N$ and $|N|=|D|$. We are going
to show that 

for each perfect set $P\subseteq [T]$ there is $X\in N\cap\iso$ such that
$P\rest X=2^{\textstyle X}$ 

\noindent (what will finish the proof of the theorem). To this end suppose
that $T^*\subseteq T$ is a perfect tree. Since $N\cap\baire$ is a dominating
family (as the strong PP--property implies $\baire$--bounding) we may choose
an increasing sequence $\langle n_i: i<\omega\rangle\in N\cap\baire$ such that
\[(\forall i\in\omega)(\forall\nu\in T^*\cap 2^{\textstyle k_{n_i}})(|\{\eta
\in T^*\cap 2^{\textstyle k_{n_{i+1}}}: \nu\vartriangleleft\eta\}|> 2(i+1)).\]
As we may encode (in a canonical way) subsets of $2^{\textstyle n}$ as
integers, we may use the choice of $D$ and $N$ and find a sequence $\langle
w_i: i\in B\rangle\in N$ such that for each $i\in B$:
\begin{enumerate}
\item[(i)] \ \ $0<|w_i|\leq i$,
\item[(ii)]\   $a\in w_i\quad\Rightarrow\quad a\subseteq T\cap 2^{\textstyle
k_{n_i}}$, 
\item[(iii)] $T^*\cap 2^{\textstyle k_{n_i}}\in w_i$.
\end{enumerate}
By shrinking each $w_i$ if necessary, we may additionally demand that $|a|\geq
2\cdot\min(B)$ for each $a\in w_{\min(B)}$ and if $i<j$ both are from $B$,
$a\in w_j$ then   
\begin{enumerate}
\item[(iv)]\  for each $\nu\in a$ the set $\{\eta\in a:\nu\rest k_{n_i}=\eta
\rest k_{n_i}\}$ has at least $2j$ elements
\end{enumerate}
(remember the choice of the $n_i$'s). Now, working in $N$, we inductively
build a perfect tree $T^+\subseteq T$. First for each $a\in w_{\min(B)}$ we
choose $\eta^{\langle\rangle}_{a,0},\eta^{\langle\rangle}_{a,1}\in a$ such
that there are no repetitions in $\langle\eta^{\langle\rangle}_{a,\ell}: a\in
w_{\min(B)},\;\ell<2\rangle$ (possible as $|w_{\min(B)}|\leq\min(B)$ and $|a|
\geq 2\cdot\min(B)$ for $a\in w_{\min(B)}$). We declare that 
\[T^+\cap 2^{\textstyle k_{n_{\min(B)}}}=\{\eta^{\langle\rangle}_{a,\ell}:
a\in w_{\min(B)},\; \ell<2\}.\]
Suppose that we have defined $T^+\cap 2^{\textstyle k_{n_i}}$, $i\in B$ and
$j=\min(B\setminus (i+1))$. For each $\nu\in T^+\cap 2^{\textstyle k_{n_i}}$
and $a\in w_j$ we choose $\eta^{\nu}_{a,0},\eta^{\nu}_{a,1}\in T\cap
2^{\textstyle k_{n_j}}$ such that 
\begin{enumerate}
\item there are no repetitions in $\langle\eta^{\nu}_{a,0},\eta^{\nu}_{a,1}:
a\in w_j\rangle$,
\item if $a\in w_j$ is such that $(\exists\eta\in a)(\nu\vartriangleleft\eta)$
then $\eta^{\nu}_{a,0},\eta^{\nu}_{a,1}\in a$.
\end{enumerate}
Again, the choice is possible by {\bf (iv)}. We declare that 
\[T^+\cap 2^{\textstyle k_{n_j}}=\{\eta^{\nu}_{a,\ell}: \nu\in T^+\cap
2^{\textstyle k_{n_i}},\ a\in w_j,\ \ell<2\}.\]
This fully describes the construction (in $N$) of the tree $T^+\subseteq T$. 
Next, working still in $N$, we choose integers $m_j\in [k_{n_j},k_{n_j+1})$ 
such that for each $j\in B$:
\begin{enumerate}
\item[(+)] {\em if\/} $\nu\in T^+\cap 2^{\textstyle k_{n_i}}$, $i\in B$ is
such that $j=\min(B\setminus (i+1))$ (or $\nu=\langle\rangle$ and $j=\min(B)$)
and $a\in w_j$ and $\eta^{\nu}_{a,\ell}\vartriangleleft\eta\in T\cap
2^{\textstyle k_{n_j+1}}$  {\em then\/} $\eta(m_j)=\ell$
\end{enumerate}
(possible by $(\beta)$ of the choice of the tree $T$ and the first demand of
the choice of the $\eta^{\nu}_{a,\ell}$'s). Let $X=\{m_j\!: j\in B\}\in\iso
\cap N$. Suppose $f\!:X\longrightarrow 2$. By clause {\bf (iii)} we know
that for each $j\in B$, $b_j\stackrel{\rm def}{=}T^*\cap 2^{\textstyle
k_{n_j}}\in w_j$. Consequently we may build inductively an infinite branch
$\eta\in [T^*]$ such that $\eta\rest k_{n_{\min(B)}}=\eta^{\langle\rangle}_{
b_{\min(B)},f(m_{\min(B)})}$ and for each $i\in B$, if $j=\min(B\setminus(i+
1))$, $\nu_i=\eta\rest k_{n_i}$ then $\eta\rest k_{n_j}=\eta^{\nu_i}_{b_j,
f(m_j)}$. It follows from {\bf (+)} that $f\subseteq\eta$.
\end{proof}

\begin{remark}
Note that a sequence $\langle w_i: i\in B\rangle$, where $w_i\in [\omega]^{
\textstyle i}$, may be interpreted as a sequence $\langle v_i: i\in B\rangle$,
where $|v_i|=i$ and members of $v_i$ are functions from the interval
$[\frac{i(i+1)}{2},\frac{(i+1)(i+2)}{2})$ to $\omega$. Consequently,
considering suitable diagonals, we may use Bartoszy\'nski--Miller's
characterization of $\non(\M)$ (see \cite[2.4.7]{BaJu95}) and show that
$\non(\M)\leq\dominating(R^{\rm sPP})$. As clearly $\dominating\leq\dominating
(R^{\rm sPP})$ we conclude that $\cof(\M)\leq\dominating(R^{\rm sPP})$ (by
\cite[2.2.11]{BaJu95}). On the other hand, it follows from Bartoszy\'nski's
characterization of $\cof(\N)$ (see \cite[2.3.9]{BaJu95}) that $\dominating(
R^{\rm sPP})\leq\cof(\N)$ (just note that the Sacks property implies the
strong PP--property).
\end{remark}

\section{An iterable friend of the strong PP--property}
The PP--property is preserved in countably support iterations of proper
forcing notions (see \cite[Ch VI, 2.12]{Sh:f}). However, it is not clear if
the strong PP--property is preserved. Here, we introduce a property stronger
then the strong PP which is preserved in countable support iterations. 

\begin{definition}
\label{DsPP}
Let $\D$ be a filter on $\omega$, $x\in\baire$ be a non-decreasing function and
let $\bar{\F}=(\F,<_{\F})$ be a partial order on $\F\subseteq\baire$.
\begin{enumerate}
\item The filter $\D$ is {\em weakly non-reducible} if it is non-principal and
for every partition $\langle X_n: n<\omega\rangle$ of $\omega$ into finite
sets there exists a set $Y\in\iso$ such that $\bigcup\limits_{n\in\omega
\setminus Y} X_n\in\D$.

If above we allow partitions into sets from the dual ideal $\D^c$ then we say
that $\D$ is {\em non-reducible}.
\item The partial order $\bar{\F}$ is {\em PP--ok} if it is dense, has no
maximal and minimal elements, each member of $\F$ is non-decreasing, the
identity function belongs to $\F$ and
\begin{enumerate}
\item[($\otimes$)] if $h_0,h_1\in \F$, $h_0<_{\F}h_1$ then $(\forall n\in
\omega)(1\leq h_0(n)\leq h_1(n))$ and 

$\lim\limits_{n\to\infty}\frac{h_1(n)}{h_0(n)}=\infty$ and there is $h\in\F$
such that 
\[(\forall^\infty n\in\omega)(h(n)\leq\frac{h_1(n)}{h_0(n)}).\]  
\end{enumerate}
\item We say that a proper forcing notion $\p$ has {\em the $(\D,x)$--strong
PP--property} if 
\[\begin{array}{ll}
\forces_{\p}&\mbox{`` for every }\eta\in\baire\mbox{ there are } B\in\D\cap\V
\mbox{ and }\langle w_i:i\in B\rangle\in\V\mbox{ such that}\\
\ &\ \ (\forall i\in B)(|w_i|\leq x(i)\ \ \&\ \ \eta(i)\in w_i)\ \mbox{''.}
  \end{array}\]
\item A proper forcing notion $\p$ has {\em the $(\D,\bar{\F})$--strong
PP--property} if it has the $(\D,x)$--strong PP--property for every $x\in\F$.
\end{enumerate}
\end{definition}

\begin{remark}
\begin{enumerate}
\item Each non-principal ultrafilter on $\omega$ is non-reducible. Clearly
non-reducible filters are weakly non-reducible. 
\item One can easily construct a countable partial order $\bar{\F}$ which is
PP--ok and such that $x\in\F$ for any pregiven non-decreasing unbounded
function $x\in\baire$.   
\end{enumerate}
\end{remark}

\begin{theorem}
\label{ppiniter}
Suppose that $\D$ is a non-principal p-filter on $\omega$ (see
\ref{ramppoint}(2)) and $\bar{\F}=(\F,<_{\F})$ is a PP--ok partial order.  
Then:  
\begin{enumerate}
\item every proper forcing notion which has the $(\D,\bar{\F})$--strong
PP--property has the strong PP--property,
\item the $(\D,\bar{\F})$--strong PP--property is preserved in countable
support iterations of proper forcing notions.
\end{enumerate}
\end{theorem}

\begin{proof}
1) Should be clear.

\noindent 2)\ \ \ We will apply \cite[Ch VI, 1.13A]{Sh:f}, so we will follow
the terminology of \cite[Ch VI, \S1]{Sh:f}. However, we will not quote the
conditions which we have to check, as that was done in the proof of
\ref{pretFbound} (and the proof here is parallel to the one there). We will
present the proof in a slightly more complicated way than needed, but later we
will be able to refer to it in a bounded context (in \ref{presRthm}).
Moreover, in this way the analogy to \ref{pretFbound} will be more clear. 

For each $m\geq 1$ we fix a function $\psi^m:\omega\stackrel{\rm onto}{
\longrightarrow} [\omega]^{\textstyle {\leq}m}$ and for $h\in\F$ we define
$\psi_h:\baire\longrightarrow\prod\limits_{n\in\omega}\fsuo$ by $\psi_h(\eta)
(n)=\psi^{h(n)}(\eta(n))$. Further, for $h^*,h\in\F$, $B\in\D$ and $\bar{w}=
\langle w_i: i\in B\rangle$ such that $h^*<_{\F}h$ and $(\forall i\in B)(w_i
\in [\omega]^{\textstyle{\leq} h(i)})$ we put 
\[T_{\bar{w},h^*}=\{\nu\in\fseo: (\forall i\in B\cap\lh(\nu))(\psi^{h^*(i)}(
\nu(i))\subseteq w_i)\}.\]
Each $T_{\bar{w},h^*}$ is a perfect subtree of $\fseo$. Now we define:
\begin{itemize}
\item $D_{\D,\bar{\F}}$ is $\cH(\aleph_1)^{\V}$,
\item for $x,T\in D_{\D,\bar{\F}}$ we say that\qquad $x\;R_{\D,\bar{\F}}\; T$
\qquad if and only if\\
$x=\langle h^*,h\rangle$ and $T=T_{\bar{w},h^*}$ for some $h^*,h\in\F$ and
$\bar{w}=\langle w_i: i\in B\rangle\in D_{\D,\bar{\F}}$ such that\qquad
$h^*<_{\F}h$,\quad $B\in\D$,\quad and\quad $(\forall i\in B)(w_i\in [
\omega]^{\textstyle {\leq}h(i)})$,
\item for $\langle h^*,h\rangle,\langle h^{**},h'\rangle\in \dom(R_{\D,\bar{
\F}})$ we say that  $\langle h^*,h\rangle<_{\D,\bar{\F}}\langle h^{**},h'
\rangle$ if and only if $h^*=h^{**}<_{\F}h<_{\F}h'$.
\end{itemize}

\begin{claim}
\label{cl41}
\begin{enumerate}
\item $(D_{\D,\bar{\F}},R_{\D,\bar{\F}})$ is a weak covering model in $\V$.
\item In any generic extension $\V^*$ of $\V$ in which $(D_{\D,\bar{\F}},R_{
\D,\bar{\F}})$ covers, a forcing notion $\p$ is $(D_{\D,\bar{\F}},R_{\D,
\bar{\F}})$--preserving if and only if it has the $(\D,\bar{\F})$--strong
PP--property.
\end{enumerate}
[Compare \ref{cl21}, \ref{cl22}.]
\end{claim}

\noindent{\em Proof of the claim:}\ \ \ 1)\ \ \ Check.

\noindent 2)\ \ \ Suppose that $(D_{\D,\bar{\F}},R_{\D,\bar{\F}})$ covers in
$\V^*$ and $\p\in\V^*$ is a forcing notion with the $(\D,\bar{\F})$--strong
PP--property. Let $\langle h^*,h\rangle\in\dom(R_{\D,\bar{\F}})$ and $\eta\in
\baire\cap (\V^*)^{\p}$. Choose $h_0,h_1\in\F$ such that $h_0<_{\F}h_1$ and
$(\forall^\infty n\in\omega)(h_1(n)\leq\frac{h(n)}{h^*(n)})$ (possible by
\ref{DsPP}(2)). Take $\bar{w}^*=\langle w^*_i: i\in B^*\rangle\in \V^*$ such
that $B^*\in\D\cap\V$ and
\[(\forall i\in B^*)(w^*_i\in [\omega]^{\textstyle {\leq}h_0(i)}\ \ \&\ \
\eta(i)\in w^*_i).\]
Let $\eta^*\in\baire\cap\V^*$ be such that $\psi^{h_0(i)}(\eta^*(i))=w^*_i$
(for $i\in B^*$). Since $(D_{\D,\bar{\F}},R_{\D,\bar{\F}})$ covers in $\V^*$
and $\langle h_0,h_1\rangle\in\dom(D_{\D,\bar{\F}})$, we find $\bar{w}=\langle
w_i: i\in B\rangle\in D_{\D,\bar{\F}}$ such that $\eta^*\in\lim(T_{\bar{w},
h_0})$ and $(\forall i\in B)(w_i\in [\omega]^{\textstyle {\leq}h_1(i)})$. Let
$B^+=B\cap B^*\in\D\cap\V$, $w^+_i=\bigcup\limits_{k\in w_i} \psi^{h^*(i)}(k)$
for $i\in B^+$. Note that for sufficiently large $i\in B^+$
\[|w^+_i|\leq|w_i|\cdot h^*(i)\leq h_1(i)\cdot h^*(i)\leq h(i)\]
and we may assume that this holds for all $i\in B^+$. Letting $\bar{w}^+=
\langle w^+_i: i\in B^+\rangle$ we will have $\langle h^*,h\rangle\;R_{\D,
\bar{\F}}\; T_{\bar{w}^+,h^*}$ and $\eta\in\lim(T_{\bar{w}^+,h^*})$, and hence
$(D_{\D,\bar{\F}},R_{\D,\bar{\F}})$ covers in $(\V^*)^{\p}$. 
\smallskip

The converse implication is even simpler. 

\begin{claim}
\label{cl42}
$(D_{\D,\bar{\F}},R_{\D,\bar{\F}},<_{\D,\bar{\F}})$ is a fine covering model.

[Compare \ref{cl23}.]
\end{claim}

\noindent{\em Proof of the claim:}\ \ \ Immediately by the definition of
$(D_{\D,\bar{\F}}, R_{\D,\bar{\F}},<_{\D,\bar{\F}})$ one sees that the demands
($\alpha$), ($\beta$)(i)--(iii) of \cite[Ch VI, 1.2(1)]{Sh:f} are satisfied. 
To verify the condition ($\beta$)(iv) of \cite[Ch VI, 1.2(1)]{Sh:f} suppose
that $\langle h^*,h\rangle<_{\D,\bar{\F}} \langle h^*,h'\rangle$ and $\langle
h^*,h\rangle\;R_{\D,\bar{\F}}\;T_{\bar{w}^\ell,h^*}$, $\bar{w}^\ell=\langle
w^\ell_i: i\in B_\ell\rangle$, $B_\ell\in\D$ (for $\ell=1,2$). Take $n\in
\omega$ such that $(\forall m\geq n)(2\cdot h(m)<h'(m))$ (possible by
\ref{DsPP}(2)) and let $B=B_1\cap B_2\cap [n,\omega)\in\D$. Put $w_i=w^1_i\cup
w^2_i$ for $i\in B$ and look at the tree $T_{\bar{w},h^*}$. Clearly $\langle
h^*,h'\rangle\; R_{\D,\bar{\F}}\; T_{\bar{w},h^*}$ and $T_{\bar{w}^1,h^*}\cup
T_{\bar{w}^2,h^*}\subseteq T_{\bar{w},h^*}$ (so more than needed). 

Checking clauses ($\gamma$) and ($\delta$) of \cite[Ch VI, 1.2(1)]{Sh:f} we
restrict ourselves to the stronger condition ($\delta$). So suppose that
$\V^*$ is a generic extension (via a proper forcing notion) of $\V$ such that
$\V^*\models$ ``$(D_{\D,\bar{\F}},R_{\D,\bar{\F}})$ covers''.  

\noindent{\bf (a)}\ \ \ Assume that $x,x^+,x_n\in\dom(R_{\D,\bar{\F}})$, $T_n
\in D_{\D,\bar{\F}}$ are such that for $n\in\omega$
\[x_n<_{\D,\bar{\F}}x_{n+1}<_{\D,\bar{\F}}x^+<_{\D,\bar{\F}}x\quad \mbox{ and
}\quad x_n\;R_{\D,\bar{\F}}\; T_n.\]
Let $x=\langle h^*,h\rangle$, $x_n=\langle h^*,h_n\rangle$, $x^+=\langle h^*,
h^+\rangle$ (so $h^*<_{\F}h_n<_{\F}h_{n+1}<_{\F}h^+<_{\F}h$) and let
$\bar{w}^n=\langle w^n_i: i\in B_n\rangle\in D_{\D,\bar{\F}}$ be such that
$T_n=T_{\bar{w}^n,h^*}$ (so $B_n\in\D\cap\V$ and $|w^n_i|\leq h_n(i)$). Look
at the sequence $\langle B_n: n\in\omega\rangle\subseteq\D\cap\V$. It does not
have to belong to $\V$, but it may be covered by a countable set from $\V$ (as
$\V^*$ is a proper forcing extension of $\V$). Hence, as $\D$ is a $p$-filter,
we find a set $B\in\D$ such that $(\forall n\in\omega)(|B\setminus B_n|<
\omega)$. Take $h_0^-,h^-_1\in\F$ such that $h^-_0<_{\F}h^-_1$ and
$(\forall^\infty n\in\omega)(h_1^-(n)\leq\frac{h(n)}{h^+(n)})$ (remember
$(\otimes)$ of \ref{DsPP}(2)) and choose an increasing sequence $m_0<m_1<m_2<
\ldots<\omega$ such that $(\forall i\geq m_0)(h_1^-(i)\leq\frac{h(i)}{h^+(i)}
)$, $B\setminus B_0\subseteq m_0$ and for $n\in\omega$:  
\[B\setminus B_{m_n}\subseteq m_{n+1}\quad\mbox{ and }\quad(\forall i\geq
m_{n+1})((n+2)\cdot h_{m_n}(i)<h^+(i)).\]
Let $\eta\in\baire\cap\V^*$ be such that if $i\in B\cap [m_n,m_{n+1})$ then
$\psi^{h^+(i)}(\eta(i))=w^0_i\cup\bigcup\limits_{k<n} w^{m_k}_i$, and let
$\eta^*\in \baire\cap\V^*$ be such that $\psi^{h^-_0(i)}(\eta^*(i))=\{\eta(i)
\}$ for each $i\in\omega$. Since $(D_{\D,\bar{\F}},R_{\D,\bar{\F}})$ covers in
$\V^*$ we find $T_{\bar{w}^-,h_0^-}\in D_{\D,\bar{\F}}$ such that $\eta^*\in
\lim(T_{\bar{w}^-,h_0^-})$ and $\langle h^-_0,h^-_1\rangle\;R_{\D,\bar{\F}}\;
T_{\bar{w}^-,h_0^-}$. Let $B^*\stackrel{\rm def}{=}B\cap \dom(\bar{w}^-)
\setminus m_0\in\D$. By the choice of $\eta$, $T_{\bar{w}^-,h_0^-}$ and
$h_1^-$ we find $\bar{w}^*=\langle w^*_i: i\in B^*\rangle\in D_{\D,\bar{\F}}$
such that $|w^*_i|\leq h^-_1(i)\cdot h^+(i)\leq h(i)$ and $w^0_i\cup
\bigcup\limits_{k<n} w^{m_k}_i\subseteq w^*_i$ whenever $i\in B^*\cap [m_n,
m_{n+1})$, $n\in\omega$. Clearly $T_{\bar{w}^*,h^*}\in D_{\D,\bar{\F}}$ and
$\langle h^*,h\rangle\; R_{\D,\bar{\F}}\; T_{\bar{w}^*,h^*}$. Put
$W\stackrel{\rm def}{=}\{m_0,m_1,m_2,\ldots\}$ and suppose that $\rho\in
\baire\cap\V^*$ is such that for every $n\in\omega$, $\rho\rest m_{n+1}\in
\bigcup\limits_{k<n}T_{m_k}\cup T_0$. If $i\in B^*$, $m_n\leq i<m_{n+1}$ then,
by the assumptions on $\rho$, $\psi^{h^*(i)}(\rho(i))\subseteq w^0_i\cup
\bigcup\limits_{k<n} w^{m_k}_i\subseteq w^*_i$. Hence $\rho\in\lim(T_{
\bar{w}^*,h^*})$ (remember $B^*\subseteq [m_0,\omega)$).

\noindent{\bf (b)}\ \ \ Assume that $x=\langle h^*,h\rangle\in\dom(R_{\D,
\bar{\F}})$, $\eta_n,\eta\in\baire$ are such that $\eta\rest n=\eta_n\rest n$
for $n\in\omega$. Take $h'\in\F$ such that $h^*<_{\F} h'<_{\F} h$ and choose
an increasing sequence $0=m_0<m_1<m_2<\ldots<\omega$ such that for each $n\in
\omega$
\[(\forall m\geq m_{n+1})((n+2)\cdot h^*(m)<h'(m)).\]
Let $\eta^*\in\baire$ (in $\V^*$) be such that
\begin{quotation}
\noindent if $m\in [m_n,m_{n+1})$, $n\in\omega$, $0<k<n+1$

\noindent then $\psi^{h^*(m)}(\eta_{m_k}(m))\subseteq\psi^{h'(m)}(\eta^*(m))$
\end{quotation}
(remember the choice of the $m_n$'s). Since $(D_{\D,\bar{\F}}, R_{\D,
\bar{\F}})$ covers in $\V^*$ we find $\bar{w}=\langle w_i: i\in B\rangle$ such
that $\langle h',h\rangle\; R_{\D,\bar{\F}}\; T_{\bar{w},h'}$ (so in
particular $B\in\D$ and $|w_i|\leq h(i)$) and $\eta^*\in\lim(T_{\bar{w},h'})$.
But now look at the tree $T_{\bar{w},h^*}$. Clearly it satisfies $\langle h^*,
h\rangle\; R_{\D,\bar{\F}}\; T_{\bar{w},h^*}$. Moreover, one can inductively
show that $\eta_{m_{n+1}}\in\lim(T_{\bar{w},h^*})$ for each $n\in\omega$. 
[Why? Plainly for each $m\in B$ we have $\psi^{h^*(m)}(\eta_{m_1}(m))\subseteq
\psi^{h'(m)}(\eta^*(m))\subseteq w_m$, so $\eta_{m_1}\in\lim(T_{\bar{w},h^*}
)$. Looking at $\eta_{m_{n+1}}$, $n>0$, note that $\eta_{m_{n+1}}\rest m_n=
\eta_{m_n}\rest m_n$ and for each $m\geq m_n$ we have $\psi^{h^*(m)}(\eta_{
m_{n+1}}(m))\subseteq \psi^{h'(m)}(\eta^*(m))$.] Thus $T_{\bar{w},h^*}$ is as
required, finishing the proof of the claim. 
\medskip

Finally, due to \ref{cl41}, \ref{cl42} we may apply \cite[Ch VI, 1.13A]{Sh:f}
to conclude the theorem.
\end{proof}

\begin{theorem}
\label{getDsPP}
Let $\D$ be a weakly non-reducible filter on $\omega$, $x\in\baire$ be an
unbounded non-decreasing function. 
\begin{enumerate}
\item If $(K,\Sigma)$ is a finitary t-omittory tree--creating pair then the
forcing notion $\q^{\tree}_1(K,\Sigma)$ has the $(\D,x)$--strong PP--property.
\item If $(K,\Sigma)$ is a finitary creating pair which captures singletons
then the forcing notion $\q^*_{{\rm w}\infty}(K,\Sigma)$ has the
$(\D,x)$--strong PP--property.  
\end{enumerate}
\end{theorem}

\begin{proof}
1)\ \ \ Suppose that $\dot{\eta}$ is a $\q^{\tree}_1(K,\Sigma)$--name for an
element of $\baire$, $p\in\q^{\tree}_1(K,\Sigma)$. Choose a condition $q\geq
p$ and fronts $F_n$ of $T^q$ such that for each $n\in\omega$ 
\begin{enumerate}
\item if $\nu\in F_n$ then the condition $q^{[\nu]}$ decides the value of
$\dot{\eta}(n)$,
\item $(\forall \nu\in F_n)(\nor[t^q_\nu]>n+1)$,
\item $(\forall \nu\in F_{n+1})(\exists \nu'\in F_n)(\nu'\vartriangleleft\nu)$
\end{enumerate}
(possible by \ref{treedec}(2), \ref{remtomit}). Next choose an increasing
sequence $0=n_0<n_1<n_2<\ldots<\omega$ such that for each $k\in\omega$
\[|\bigcup\{\pos(t^q_\nu):\nu\in F_{n_k}\}|<x(n_{k+1}).\]
Since $\D$ is weakly non-reducible we find $Y\in\iso$ such that $B
\stackrel{\rm def}{=}\bigcup\limits_{k\in\omega\setminus Y} [n_k,n_{k+1})\in
\D$. Now construct inductively a condition $q^*\geq q$ such that $\mrot(q^*)=
\mrot(q)$ and
\begin{enumerate}
\item[(a)] $T^{q^*}\subseteq\{\mrot(q)\}\cup \bigcup\{\pos(t^q_\nu): \nu\in
F_{n_k}\ \&\ k\in Y\}$,
\item[(b)] if $\nu\in T^{q^*}$ then $\pos(t^{q^*}_\nu)\subseteq\pos(t^q_{
\nu^*})$ and $\nor[t^{q^*}_\nu]\geq \nor[t^q_{\nu^*}]-1$ for some $\nu^*\in
\bigcup\limits_{k\in Y} F_{n_k}$.
\end{enumerate}
It should be clear that one can build such $q^*$ (remember $(K,\Sigma)$ is
t-omittory). Note that if $k_0,k_1\in Y$, $k_0<k_1$ and $(k_0,k_1)\subseteq
\omega\setminus Y$ then for each $n\in (n_{k_0},n_{k_1}]$ we have
\[|\dcl(T^{q^*})\cap F_n|=|\dcl(T^{q^*})\!\cap\! F_{n_{k_0+1}}|\leq |\bigcup
\{\pos(t^q_\nu)\!:\nu\!\in\! F_{n_{k_0}}\}|<x(n_{k_0+1})\leq x(n).\]
For $n\in B$ let $w_n=\{m\in\omega: (\exists\nu\in F_n\cap\dcl(T^{q^*}))(q^{[
\nu]}\forces\dot{\eta}(n)=m)\}$. By the above remark we have $|w_n|\leq x(n)$
(for $n\in B$) and clearly $q^*\forces(\forall n\in B)(\dot{\eta}(n)\in w_n)$. 
\medskip

\noindent 2)\ \ \ Similar.
\end{proof}

\section{Bounded relatives of PP}
In the following definition we introduce relations which determine
localization properties (see \ref{gendef}(2)) close to the PP--property when
restricted to functions from $\prod\limits_{n\in\omega}f(n)$. Not surprisingly
they include (the relation responsible for) the $(f,g)$--bounding property
too.  

\begin{definition}
\label{fgPP}
Let $f,g\in\baire$ be non-decreasing functions such that $(\forall n\in\omega)
(0<g(n)<f(n))$. Define:
\begin{enumerate}
\item $S_{f,g}=\prod\limits_{n\in\omega}[f(n)]^{\textstyle g(n)}$,\quad
$S^*_{f,g}=S_{f,g}\times\iso$, 
\item $R^\exists_{f,g},R^\forall_{f,g}\subseteq\prod\limits_{n\in\omega}f(n)
\times S_{f,g}$ are given by

\noindent $\eta\; R^\exists_{f,g}\; \bar{A}$\qquad if and only if 

\noindent $\eta\in\prod\limits_{n\in\omega}f(n)$, $\bar{A}=\langle A_n: n\in
\omega\rangle\in S_{f,g}$ and $(\exists^\infty n\in\omega)(\eta(n)\in A_n)$,

\noindent $\eta\; R^\forall_{f,g}\; \bar{A}$\qquad if and only if 

\noindent $\eta\in\prod\limits_{n\in\omega}f(n)$, $\bar{A}=\langle A_n: n\in
\omega\rangle\in S_{f,g}$ and $(\forall^\infty n\in\omega)(\eta(n)\in A_n)$,
\item a relation $R^*_f\subseteq\prod\limits_{n\in\omega}f(n)\times
\prod\limits_{n\in\omega}f(n)$ is defined by

\noindent $\eta_0\; R^*_f\; \eta_1$\quad if and only if\quad $\eta_0,\eta_1
\in\prod\limits_{n\in\omega}f(n)$ and $(\exists^\infty n\in\omega)(\eta_0(n)=
\eta_1(n))$,
\item a relation $R^{**}_{f,g}\subseteq \prod\limits_{n\in\omega}f(n)\times
S^*_{f,g}$ is such that 

\noindent $\eta\; R^{**}_{f,g}\; (\bar{A},K)$\qquad if and only if 

\noindent $\eta\in\prod\limits_{n\in\omega}f(n)$, $\bar{A}=\langle A_n: n\in
\omega\rangle\in S_{f,g}$, $K=\{k_0,k_1,k_2,\ldots\}\in\iso$ (the increasing
enumeration) and $(\forall m\in\omega)(\exists n\in [k_m,k_{m+1}))(\eta(n)\in
A_n)$. 
\end{enumerate}
\end{definition}

\begin{remark}
\begin{enumerate}
\item The spaces $S_{f,g}$, $S^*_{f,g}$ and $\prod\limits_{n\in\omega}f(n)$
carry natural (product) Polish topologies. 
\item The relation $R^\forall_{f,g}$ corresponds to the $(f,g)$--bounding
property, of course. The cardinal number $\dominating(R^\forall_{f,g})$ is the
$c(f,g)$ of \cite{GoSh:448} (see there for various ZFC dependencies between
the cardinals determined by different functions as well as for consistency
results).  
\item Note that the relation $R^{**}_{f,g}$ (or actually the corresponding
localization property) is really very close to the PP--property. The cardinal
numbers $\dominating(R^*_f)$ and $\dominating(R^{**}_{f,g})$ appear naturally
in \cite{BRSh:616}.
\item There are other natural variants of relations introduced in \ref{fgPP}. 
We will deal with them (and the corresponding cardinal invariants) in the
continuation of this paper. 
\end{enumerate}
\end{remark}
Below we list some obvious relations between the localization properties
introduced in \ref{fgPP} and the corresponding cardinal numbers.

\begin{proposition}
\label{depPPfg}
Let $f,g,h\in\baire$ be non-decreasing functions such that $0<g(n)<f(n)$ for
each $n\in\omega$. Then:
\begin{enumerate}
\item The $R^\forall_{f,g}$--localization implies the
$R^\exists_{f,g}$--localization and $\dominating(R^\exists_{f,g})\leq
\dominating(R^\forall_{f,g})$. 
\item Suppose that for some increasing sequence $m_0<m_1<m_2<\ldots<\omega$ we
have 
\[(\forall n\in\omega)(g(n)\leq m_{n+1}-m_n\ \ \&\ \ f(n)\geq\prod_{k\in
[m_n,m_{n+1})} h(k)).\]
Then the $R^\exists_{f,g}$--localization implies the
$R^*_h$--localization 
and $\dominating(R^*_h)\leq\dominating(R^\exists_{f,g})$.
\item The $R^*_f$--localization implies the $R^\exists_{f,g}$ localization
and $\dominating(R^\exists_{f,g})\leq\dominating(R^*_f)$. 
\item The $R^\exists_{f,g}$--localization plus $\baire$--bounding imply the
$R^{**}_{f,g}$--localization. The $R^{**}_{f,g}$--localization implies the
$R^\exists_{f,g}$--localization.  Hence
$\dominating(R^\exists_{f,g})\leq\dominating(R^{**}_{f,g})\leq\max\{ 
\dominating,\dominating(R^\exists_{f,g})\}$.
\item If $g$ is unbounded then the strong PP--property implies the
$R^\exists_{f,g}$--localization, and $\dominating(R^\exists_{f,g})\leq
\dominating(R^{\rm sPP})$.
\item Assume ground model reals are not meager. Then the extension has the
$R^*_f$--localization property and thus $\dominating(R^*_f)\leq\non(\M)$.
\item The $R^{**}_{f,g}$--localization implies that there is no Cohen real
over the ground model, and thus $\cov(\M)\leq\dominating(R^{**}_{f,g})$.
\end{enumerate}
\end{proposition}
For getting the $R^\forall_{f,g}$--localization (i.e.\ $(f,g)$--bounding
property) for forcing notions built according to our schema see
\ref{sekfg}. Let us note that the other properties appear naturally too.

\begin{proposition}
\label{xxpp}
Let $f,g\in\baire$. Suppose that $\p$ is a forcing notion one of the following
type
\begin{itemize}
\item $\q^*_{{\rm s}\infty}(K,\Sigma)$ for a finitary creating pair
$(K,\Sigma)$ which is either growing and big or omittory and omittory--big,
\item $\q^*_{{\rm w}\infty}(K,\Sigma)$ for a finitary creating pair which
captures singletons,
\item $\q^{\tree}_0(K,\Sigma)$ for a finitary t-omittory tree--creating pair
$(K,\Sigma)$. 
\end{itemize}
Then $\p$ has the $R^{**}_{f,g}$--localization property.
\end{proposition}

\begin{proof}
It should be clear, so we will sketch the proof for the first case only. Let
$\dot{\eta}$ be a $\q^*_{{\rm s}\infty}(K,\Sigma)$--name for a function in
$\prod\limits_{n\in\omega}f(n)$ and let $p\in\q^*_{{\rm s}\infty}(K,
\Sigma)$. Using \ref{decbel} or \ref{omitdecbel} construct a condition $q\in
\q^*_{{\rm s}\infty}(K,\Sigma)$, an enumeration $\langle u_k:k\in\omega
\rangle$ of $\bigcup\limits_{n\in\omega}\pos(w^q,t^q_0,\ldots,t^q_{n-1})$ and 
a sequence $\langle m_k: k<\omega\rangle$ such that
\begin{enumerate}
\item $p\leq_0 q$,\quad $m_0<m_1<\ldots<\omega$,
\item if $u_k\in\pos(w^q,t^q_0,\ldots,t^q_{n-1})$ then the condition $(u_k,
t^q_n,t^q_{n+1},\ldots)$ decides the value of $\dot{\eta}(m_k)$.
\end{enumerate}
Plainly the construction is possible and easily it finishes the proof.
\end{proof}

It may be not clear how one can preserve (in countable support iterations) the
localization properties introduced in \ref{fgPP}. To deal with the
$R^\exists_{f,g}$--localization property we may adopt the approach of 
\ref{ppiniter}. It slightly changes the meaning of this notion but the change
is not serious and makes dealing with compositions much easier.

\begin{definition}
\label{presRdef}
Let $h_0,h_1,f\in\baire$ be non-decreasing unbounded functions, $\D$ be a
filter on $\omega$ and $\bar{\F}=(\F,<_{\F})$ be a partial order on $\F
\subseteq\prod\limits_{n\in\omega}f(n)$. We say that a proper forcing notion
$\p$: 
\begin{enumerate}
\item {\em has the $\D$-$R^\exists_{f,h_0,h_1}$--localization property} if 
\[\begin{array}{ll}
\forces_{\p}&\mbox{`` for every }\eta\in\prod\limits_{n\in\omega}[f(n)]^{
\textstyle {\leq}h_0(n)}\mbox{ there are } B\in\D\cap\V \mbox{ and }\langle
w_i:i\in B\rangle\in\V\\   
\ &\ \mbox{ such that }\ (\forall i\in B)(|w_i|\leq h_1(i)\ \ \&\ \ \eta(i)
\subseteq w_i)\ \mbox{'',} 
  \end{array}\]
\item {\em has the $(\D,\bar{\F})$--$R^\exists_f$--localization property} if
it has the $\D$-$R^\exists_{f,h_0,h_1}$--localization property for every
$h_0,h_1\in\F$ such that $h_0<_{\F} h_1$. 
\end{enumerate}
\end{definition}

\begin{proposition}
\label{presRthm}
Suppose that $\D$ is a non-principal p-filter on $\omega$, $f\in\baire$ is
non-decreasing unbounded and $\bar{\F}=(\F,<_{\F})$ is a PP--ok partial order
on $\F\subseteq\prod\limits_{n\in\omega}f(n)$ (see \ref{DsPP}(2), except that
it does not have to contain the identity function). Let $\langle\p_\alpha,
\dot{\q}_\alpha: \alpha<\delta\rangle$ be a countable support iteration such
that for each $\alpha<\delta$  
\[\begin{array}{ll}
\forces_{\p_\alpha}&\mbox{`` }\dot{\q}_\alpha\mbox{ is a proper forcing notion
which has}\\
\ &\ \mbox{ the $(\D,\bar{\F})$--$R^\exists_f$--localization property''.}
  \end{array}\]
Then the forcing notion $\p_\alpha$ has the $(\D,\bar{\F})$--$R^\exists_f
$--localization property.
\end{proposition}

\begin{proof}
Repeat the proof of \ref{ppiniter} making suitable adjustments to the fact
that we are ``below the function $f$''. No real changes are required. 
\end{proof}

\begin{remark}
Note that with no serious changes we may formulate and prove a variant of
\ref{presRthm} which would be an exact reformulation of \ref{pretFbound} for
the current context.
\end{remark}

\begin{proposition}
\label{omRex}
Let $f,h_0,h_1\in\baire$ be non-decreasing unbounded functions such that
$(\forall n\in\omega)(1\leq h_0(n)\leq h_1(n)\leq f(n))$ and $\lim\limits_{
n\to\infty}\frac{h_1(n)}{h_0(n)}=\infty$. Assume $\D$ is a non-reducible
p-filter on $\omega$. Suppose that $(K,\Sigma)$ is a finitary, omittory and
omittory--big creating pair. Then the forcing notion $\q^*_{{\rm s}\infty}(K,
\Sigma)$ has the $\D$-$R^\exists_{f,h_0,h_1}$--localization property.
\end{proposition}

\begin{proof}
Like \ref{getDsPP} plus \ref{xxpp}.
\end{proof}

\begin{remark}
Note that if $x(n)\leq\frac{h_1(n)}{h_0(n)}$ for $n\in\omega$ then the 
$(\D,x)$--strong PP--property implies the $\D$-$R^\exists_{f,h_0,
h_1}$--localization property. Consequently we may use \ref{getDsPP} to get the
conclusion of \ref{omRex} for the two types of forcing notions specified in
\ref{getDsPP}.  
\end{remark}

\section{Weakly non-reducible p-filters in iterations}
One could get an impression that \ref{ppiniter}, \ref{presRthm} together with
\ref{getDsPP} and \ref{omRex} are everything we need: the properties involved
are iterable and we may get them for various forcing notions. However, to be
able to make a real use of \ref{getDsPP} or \ref{omRex} we have to know that
if we start with a weakly non-reducible $p$-filter and then we iterate
suitable forcing notions, the filter remains weakly non-reducible. One could
start with a $p$-point and consider forcing notions which are $p$-point
preserving only. However this is much too restrictive: we may iterate forcing
notions mentioned in \ref{getDsPP} and \ref{omRex} and the iterations will
preserve the fact that the filter is weakly non-reducible. The first step in
proving this is the following observation.

\begin{proposition}
\label{preweak}
Suppose that $\D$ is a weakly non-reducible filter on $\omega$. Let $\p$ be an
almost $\baire$-bounding forcing notion. Then
\[\forces_{\p}\mbox{`` (The filer generated by) $\D$ is weakly non-reducible
''.}\]
[Note that this covers $\baire$-bounding forcing notions.]
\end{proposition}

\begin{proof}
Suppose that $\langle\dot{X}_n:n\in\omega\rangle$ is a $\p$--name for a
partition of $\omega$ into finite sets. Let $\dot{f}$ be a $\p$--name for a 
function in $\baire$ such that 
\[\forces_{\p}\mbox{`` }(\forall n\in\omega)(\exists m\in\omega)(\dot{X}_m
\subseteq [n,\dot{f}(n)))\mbox{ ''}.\]
Suppose $p\in\p$. Since $\p$ is almost $\baire$-bounding we find an increasing
function $g\in\baire$ such that
\[(\forall A\in\iso)(\exists q\geq p)(q\forces_{\p}\mbox{`` }(\exists^\infty
n\in A)(\dot{f}(n)<g(n))\mbox{ ''}).\]
Let $0=n_0<n_1<n_2<\ldots<\omega$ be such that $g(n_k)<n_{k+1}$. As the filter
$\D$ is weakly non-reducible, we find $Y\in\iso$ such that $Z\stackrel{\rm
def}{=}\bigcup\limits_{k\in\omega\setminus Y} [n_k,n_{k+1})\in\D$. Let $A=\{
n_k: k\in Y\}$. By the choice of the function $g$, there is a condition $q\geq
p$ such that 
\[q\forces_{\p}\mbox{`` }(\exists^\infty n\in A)(\dot{f}(n)<g(n))\mbox{ ''}.\]
Now look at the choice of $\dot{f}$ -- necessarily
\[q\forces_{\p}\mbox{`` }(\exists^\infty m\in\omega)(\dot{X}_m\cap Z=\emptyset
)\mbox{ ''},\]
which is enough to conclude the proposition.
\end{proof}

Note that \ref{preweak} captures {\em almost} all forcing notions mentioned in
\ref{getDsPP}, \ref{omRex}. So what is needed more is that ``$\D$ is weakly
non-reducible'' is preserved at limit stages of countable support iterations
of proper forcing notions. This is done like preserving unbounded families
(i.e.\ by \cite[Ch VI, \S3]{Sh:f}). 

\begin{theorem}
\label{presnon}
Let $\D$ be a weakly non-reducible p-filter on $\omega$. Suppose that $\langle
\p_\alpha,\dot{\q}_\alpha:\alpha<\delta\rangle$ is a countable support
iteration of proper forcing notions such that $\delta$ is limit and for each
$\alpha<\delta$
\[\forces_{\p_\alpha}\mbox{`` (The filter generated by) $\D$ is weakly
non-reducible ''.}\]
Then $\forces_{\p_\delta}$`` (The filter generated by) $\D$ is weakly
non-reducible ''. 
\end{theorem}

\begin{proof}
We will use \cite[Ch VI, 3.13]{Sh:f} and thus we will follow the notation
there. Let $F\subseteq\baire$ be the family of all increasing enumerations of
elements of $\D$ (i.e.\ $F=\{\mu_X: X\in\D\}$, see \ref{muX}). Let $R$ be (a
definition of) the following two--place relation on $\baire$: 

$g\; R\; f$\quad if and only if\quad ($g,f\in\baire$ and)

$(\exists^\infty k)([n^g_k,n^g_{k+1})\cap\rng(f)=\emptyset)$, where
$n^g_0=0$, $n^g_{k+1}=n^g_k+g(k)+1$ for $k\in\omega$. 

\noindent As $\D$ is weakly non-reducible, the family $F$ is $R$--bounding
(i.e.\ $(\forall g\in\baire)(\exists f\in F)(g\; R\; f)$).

\begin{claim}
$(F,R)$ is $S$--nice (see \cite[Ch VI, 3.2]{Sh:f}; here $S\subseteq
[F]^{\textstyle \omega}$ is arbitrary).
\end{claim}

\noindent{\em Proof of the claim:}\ \ \ We have to show that for each $N\in S$
there is $g\in F$ such that for each $m_0\in\omega$ (the $n_0$ of \cite[Ch VI,
3.2.3($\beta$)]{Sh:f} is irrelevant here) the second player has an absolute
winning strategy in the following game.
\begin{quotation}
\noindent At the stage $k$ of the game, Player I chooses $f_k\in\baire$ and
$g_k\in F\cap N$ such that $f_k\rest m_{\ell+1}=f_\ell\rest m_{\ell+1}$ for
all $0\leq \ell<k$ and $f_k\; R\; g_k$. Then Player II answers playing an
integer $m_{k+1}>m_k$.\\
Player II wins the game if $(\bigcup\limits_{k\in\omega} f_k\rest m_k)\; R\;
g$. 
\end{quotation}
But this is easy: let $g\in F$ be such that $\rng(g)\subseteq^*\rng(f)$ for
all $f\in F\cap N$ (remember $\D$ is a $p$-filter). Then 
\[f_k\; R\; g_k\quad\mbox{ implies }\quad (\exists^\infty\ell\in\omega)([
n^{f_k}_\ell,n^{f_k}_{\ell+1})\cap\rng(g)=\emptyset).\]
Thus, at stage $k$ of the game, the second player may choose $m_{k+1}>m_k$
such that $[n^{f_k}_\ell,n^{f_k}_{\ell+1})\cap \rng(g)=\emptyset$ for some
$\ell\in (m_k,m_{k+1})$.
\medskip

As we iterate proper forcing notions, countable subsets of $F$ from
$\V^{\p_\alpha}$ can be covered by countable subsets of $F$ from $\V$. By our
assumptions, $F$ is $R$--bounding in each $\V^{\p_\alpha}$ (for
$\alpha<\delta$) and it is nice there (like in the claim above). Consequently
we may apply \cite[Ch VI, 3.13(3)]{Sh:f} and we conclude that 
\[\forces_{\p_\delta}\mbox{`` $F$ is $R$--bounding ''.}\]
But this is exactly what we need. 
\end{proof}

\section{Examples}

\begin{example}
\label{BPex}
Let $P\subseteq\can$ be a perfect set. We construct a finitary function
$\bH^P$, an $\bH^P$--fast function $f^P:\omega\times\omega\longrightarrow
\omega$ and a $\bar{2}$--big trivially meagering simple creating pair
$(K^P_{\ref{BPex}},\Sigma^P_{\ref{BPex}})$ for $\bH^P$ with the (weak) Halving
Property such that 
\[\begin{array}{ll}
\forces_{\q^*_{f^P}(K^P_{\ref{BPex}},\Sigma^P_{\ref{BPex}})}&\mbox{`` there is
a perfect set $Q\subseteq P$ such that}\\
\ &\ \ (\forall K\in\iso\cap \V)(Q\rest K\neq 2^{\textstyle K})\mbox{ ''.}
  \end{array}\]
\end{example}

\begin{proof}[Construction]
The creating pair $(K^P_{\ref{BPex}},\Sigma^P_{\ref{BPex}})$ will be
constructed in a way slightly similar to $(K^1_{\ref{two}},\Sigma^1_{ 
\ref{two}})$. For positive integers $i,m$ let $R^i_2(2,m)$ be the minimal
integer $k$ such that for every function $\varphi:\prod\limits_{\ell<i}
[k]^{\textstyle 2}\longrightarrow 2$ there are sets $a_0,\ldots,a_{i-1}\in
[k]^{\textstyle m}$ such that $\varphi\rest [a_0]^{\textstyle 2}\times\ldots
\times [a_{i-1}]^{\textstyle 2}$ is constant. [Thus this is the Ramsey number
for polarized partition relations; $R^1_2(2,m)$ is essentially $R_2(2,m)$ of
\ref{ramnum}.]  

Define inductively $Z^i_2(k)$ for $k\!\in\!\omega$ by $Z^i_2(0)\!=\!
R^i_2(2,4)$, $Z^i_2(k\!+\!1)\!=\! R^i_2(2,2\cdot Z^i_2(k))$, and for a finite
set $X$ let  
\[H_i(X)\stackrel{\rm def}{=}\min\{k\in\omega: |X|\leq Z^i_2(k)\}.\]
Let $T\subseteq\fs$ be a perfect tree such that $P=[T]$. Now construct
inductively functions $\bH^P=\bH$ and $f^P=f$ and an increasing sequence
$\bar{n}=\langle n_i: i\in\omega\rangle$ such that
\begin{enumerate}
\item[(i)] \ \ $f(0,\ell)=\ell+1$, $f(k+1,\ell)=2^{\fH(\ell)+1}\cdot
(f(k,\ell)+\fH(\ell)+2)$\quad (compare \ref{fast}),
\item[(ii)]\   $n_0=0$, $n_{i+1}$ is the first such that for every $\nu\in
T\cap 2^{\textstyle n_i}$
\[H_{2^i}(\{\eta\in T\cap 2^{\textstyle n_{i+1}}: \nu\vartriangleleft\eta\})
\geq 2^{f(i,i)},\]
\item[(iii)]   $\bH(i)$ is the family of all non-empty subsets of $T\cap 2^{
\textstyle n_{i+1}}$.
\end{enumerate}
It should be clear that the clauses (i)--(iii) uniquely determine $\bH$, $f$
and $\bar{n}$. 

Call a sequence $u\in\prod\limits_{i<m}\bH(i)$\ {\em acceptable}\ if for each
$i_0<i_1<m$ 
\[\begin{array}{l}
|u(0)|=2,\quad u(i_0)=\{\eta\rest n_{i_0+1}: \eta\in u(i_1)\},\quad\mbox{
  and}\\
(\forall\nu\in u(i_0))(|\{\eta\in u(i_0+1):\nu \vartriangleleft\eta\}|=2).
\end{array}\]
Note that if $W\in\prod\limits_{m\in\omega}\bH(m)$ is such that each $W\rest
m$ is acceptable then the sequence $W$ determines a perfect tree 
\[T^W\stackrel{\rm def}{=}\{\nu\in T:(\exists m\in\omega)(\exists\eta\in W(m)
)(\nu\trianglelefteq\eta)\}\subseteq T\]
with the property that $|T^W\cap 2^{\textstyle n_i}|=2^i$ and each node 
from $T^W\cap 2^{\textstyle n_i}$ has a ramification below $n_{i+1}$.

A creature $t\in\CR[\bH]$ is in $K^P_{\ref{BPex}}$ if $m^t_{\up}=m^t_{\dn}+1 
=i+1$ and 
\begin{itemize}
\item $\dis[t]=\langle m^t_{\dn}, \langle B^t_\nu: \nu\in T\cap 2^{\textstyle
n_i}\rangle, r^t\rangle$, where $r^t$ is a non-negative real and $B^t_\nu
\subseteq\{\eta\in T\cap 2^{\textstyle n_{i+1}}: \nu\vartriangleleft\eta\}$
(for all $\nu\in T\cap 2^{\textstyle n_i}$; remember $i=m^t_{\dn}$),
\item $\val[t]=\{\langle u,v\rangle\in\prod\limits_{k<m^t_{\dn}}\bH(k)\times
\prod\limits_{k\leq m^t_{\dn}}\bH(k): u\vartriangleleft v$ both are acceptable
and 

\qquad\qquad\qquad\qquad if $\eta\in v(m^t_{\dn})$ then $\eta\in B^t_{\eta
\rest n_i}\}$,
\item $\nor[t]=\max\{0,\;\min\{H_{2^i}(B^t_\nu): \nu\in T\cap 2^{\textstyle
n_i}\}-r^t\}$.
\end{itemize}
The operation $\Sigma^P_{\ref{BPex}}$ is defined by 
\[\Sigma^P_{\ref{BPex}}(t)=\{s\in K^P_{\ref{BPex}}: m^t_{\dn}=m^s_{\dn}\ \&\
(\forall\nu\in T\cap 2^{\textstyle n_i})(B^s_\nu\subseteq B^t_\nu)\ \&\
r^s\geq r^t\}.\]
It should be clear that $(K^P_{\ref{BPex}},\Sigma^P_{\ref{BPex}})$ is a simple
finitary creating pair and the forcing notion $\q^*_f(K^P_{\ref{BPex}},
\Sigma^P_{\ref{BPex}})$ is not trivial. 

To check that $(K^P_{\ref{BPex}},\Sigma^P_{\ref{BPex}})$ is $\bar{2}$--big
suppose that $t\in K^P_{\ref{BPex}}$, $\nor[t]>1$, $u\in\basis(t)$ and
$c:\pos(u,t)\longrightarrow 2$ (note that $\basis(t)=\dom(\val[t])$ and
$\pos(u,t)=\{v\in\rng(\val[t]): u\vartriangleleft v\}$). Then $u$ is
acceptable and, if $\lh(u)>0$, $|u(\lh(u)-1)|=2^{\lh(u)}$. Let $i=\lh(u)=
m^t_{\dn}$. Let $k=\min\{H_{2^i}(B^t_\nu): \nu\in u(i-1)\}$. By the definition
of the norm of $t$ we know that $k\geq\nor[t]+r^t>1$, so necessarily $k\geq
2$. Note that under natural interpretation $\prod\limits_{\nu\in u(i-1)}
[B^t_\nu]^{\textstyle 2}\subseteq\pos(u,t)$, so we may restrict our colouring
$c$ to this set and use the definition of $H_{2^i}$ (and the choice of $k$). 
Thus we find sets $B^*_\nu\subseteq B^t_\nu$ (for $\nu\in u(i-1)$) such that 
\[(\forall\nu\in u(i-1))(|B^*_\nu|=2\cdot Z^{2^i}_2(k-2))\quad\mbox{ and
}\quad c\rest (\prod_{\nu\in u(i-1)} B^*_\nu)\ \mbox{ is constant.}\]
Note that $H_{2^i}(B^*_\nu)\geq k-1\geq\nor[t]-1+r^t>r^t$. Let $s\in K^P_{
\ref{BPex}}$ be a creature determined by
\[m^s_{\dn}=m^t_{\dn},\quad r^s=r^t,\quad B^s_\nu=B^*_\nu\ \mbox{ if }\nu\in
u(i-1),\quad \mbox{ and } B^s_\nu=B^t_\nu \mbox{ otherwise}.\]
Clearly $s\in\Sigma^P_{\ref{BPex}}(t)$, $\nor[s]\geq\nor[t]-1$ and $c\rest
\pos(u,s)$ is constant.

Plainly $(K^P_{\ref{BPex}},\Sigma^P_{\ref{BPex}})$ is trivially meagering, as
if $x,y\in X$ then 
\[H_i(X\setminus\{x,y\})\geq H_i(X)-1.\]

Let $\uhalf:K^P_{\ref{BPex}}\longrightarrow K^P_{\ref{BPex}}$ be such that if
$\nor[t]\geq 2$ then
\[\dis(\uhalf(t))=\langle m^t_{\dn},\;\langle B^t_\nu: \nu\in T,\; \lh(\nu)=
n_{m^t_{\dn}}\rangle,\; r^t+\frac{1}{2}\nor[t]\rangle,\]
and $\uhalf(t)=t$ otherwise. Exactly like in \ref{two} one checks that the
function $\uhalf$ witnesses the fact that $(K^P_{\ref{BPex}},\Sigma^P_{
\ref{BPex}})$ has the (weak) Halving Property. 

To show the last assertion of \ref{BPex} we prove that
\[\forces_{\q^*_f(K^P_{\ref{BPex}},\Sigma^P_{\ref{BPex}})}\mbox{`` }(\forall
K\in\iso\cap\V)([T^{\dot{W}}]\rest K\neq 2^{\textstyle K})\mbox{ '',}\]
where $\dot{W}$ is the name for the generic real (see \ref{thereal}) and
$T^{\dot{W}}$ is the tree defined before. To this end suppose that $p\in
\q^*_f(K^P_{\ref{BPex}},\Sigma^P_{\ref{BPex}})$ and $K\in\iso$. We may assume
that $\lh(w^p)=j_0>0$ and $(\forall i\in\omega)(\nor[t^p_i]>f^P(0,m^{t^p_i}_{
\dn})\geq 2)$. Choose $j_1\in\omega$ such that $|[n_{j_0},n_{j_1})\cap K|\geq
2^{j_0}$ and fix one-to-one mapping 
\[k:w^p(j_0-1)\longrightarrow [n_{j_0},n_{j_1})\cap K:\nu\mapsto k(\nu)\]
(remember that $w^p$ is acceptable, so $|w^p(j_0-1)|=2^{j_0}$ and
$w^p(j_0-1)\subseteq T\cap 2^{\textstyle n_{j_0}}$).

Fix $\nu\in w^p(j_0-1)$ for a moment. Let $i(\nu)=i<j_1-j_0$ be such that
$k(\nu)\in [n_{j_0+i},n_{j_0+i+1})$. For each $\eta\in T\cap 2^{\textstyle
n_{j_0+i}}$ such that $\nu\trianglelefteq\eta$ choose $c_{j_0+i}(\eta)\in 2$
such that the set $B^\nu_\eta\stackrel{\rm def}{=}\{\rho\in B^{t^p_i}_\eta:
\rho(k(\nu))=c_{j_0+i}(\nu)\}$ has at least $\frac{1}{2}|B^{t^p_i}_\eta|$
elements (so then $H_{2^{j_0+i}}(B^\nu_\eta)\geq H_{2^{j_0+i}}(B^{t^p_i}_\eta)
-1$). If $\eta=\nu$ then we finish the procedure. Otherwise, for each $\rho\in
T\cap 2^{\textstyle n_{j_0+i-1}}$ such that $\nu\trianglelefteq\rho$ we choose
$c_{j_0+i-1}(\rho)\in 2$ such that the set $B^\nu_\rho\stackrel{\rm def}{=}\{
\eta\in B^{t^p_{i-1}}_\rho: c_{j_0+i-1}(\rho)=c_{j_0+i}(\eta)\}$ has at least
$\frac{1}{2}|B^{t^p_{i-1}}_\rho|$ elements (and so $H_{2^{j_0+i-1}}(
B^\nu_\rho)\geq H_{2^{j_0+i-1}}(B^{t^p_{i-1}}_\rho)-1$). Continuing this
procedure downward till we arrive to $\nu$ we determine sets $\langle
B^\nu_\eta: \nu\trianglelefteq\eta\in T\cap 2^{\textstyle n_j},\; j_0\leq j
\leq j_0+i(\nu)\rangle$ and $c_{j_0}(\nu)\in 2$ such that
\begin{enumerate}
\item[$(\alpha)_\nu$] if $\nu\trianglelefteq\eta\in T\cap 2^{\textstyle n_j}$,
$j_0\leq j\leq j_0+i(\nu)$\\
then $B^\nu_\eta\subseteq B^{t^p_{j-j_0}}_\eta$ and $H_{2^j}(B^\nu_\eta)\geq
H_{2^j}(B^{t^p_{j-j_0}}_\eta)-1$, 
\item[$(\beta)_\nu$]  if $\nu\trianglelefteq\eta\in T\cap 2^{\textstyle
n_{j_1}}$ is such that $(\forall j\in [j_0,j_0+i(\nu)))(\eta\rest n_{j+1}\in
B^\nu_{\eta\rest n_j})$\\
then $\eta(k(\nu))=c_{j_0}(\nu)$.
\end{enumerate}

For each $i<j_1-j_0$ choose a creature $s_i\in\Sigma^P_{\ref{BPex}}(t^p_i)$
such that $r^{s_i}=r^{t^p_i}$ and for every $\eta\in T\cap 2^{\textstyle
n_{j_0+i}}$
\begin{quotation}
\noindent if $i(\eta\rest n_{j_0})\leq i$ then $B^{s_i}_\eta=B^{\eta\rest
n_{j_0}}_\eta$, otherwise $B^{s_i}_\eta=B^{t^p_i}_\eta$.
\end{quotation}
Clearly $\nor[s_i]\geq \nor[t^p_i]-1$ and thus $q=(w^p,s_0,\ldots,
s_{j_1-j_0-1},t^p_{j_1-j_0},\ldots)$ is a condition in $\q^*_{f^P}(K^P_{
\ref{BPex}},\Sigma^P_{\ref{BPex}})$ stronger than $p$. Let $\sigma:K\cap
[n_{j_0},n_{j_1})\longrightarrow 2$ be such that $\sigma(k(\nu))=1-c_{j_0}(
\nu)$ for each $\nu\in w^p(j_0-1)$. Note that
\[u\in\pos(w^p,s_0,\ldots,s_{j_1-j_0-1})\quad\Rightarrow\quad (\forall \eta
\in u(j_1-1))(\eta\rest (K\cap [n_{j_0},n_{j_1})\neq \sigma)),\]
what finishes the proof.
\end{proof}

\begin{conclusion}
\label{BPabod}
It is consistent that $\kappa_{\rm BP}=\non(\M)=\dominating(R^{\rm sPP})=\con
=\aleph_2$ and $\dominating=\aleph_1$.
\end{conclusion}

\begin{proof}
Start with a model for CH and build inductively (with a suitable bookkeeping)
a countable support iteration $\langle\p_\alpha,\dot{\q}_\alpha:\alpha<
\omega_2\rangle$ and a sequence $\langle\dot{P}_\alpha:\alpha<\omega_2\rangle$
such that 
\begin{enumerate}
\item[($\alpha$)] $\langle \dot{P}_\alpha: \alpha<\omega_2\rangle$ lists with
$\omega_2$--repetitions all $\p_{\omega_2}$--names for perfect subsets of
$\can$; each $\dot{P}_\alpha$ is a $\p_\alpha$--name,
\item[($\beta$)]  $\dot{\q}_\alpha$ is a $\p_\alpha$--name for the forcing
notion $\q^*_{f^{\dot{P}_\alpha}}(K^{\dot{P}_\alpha}_{\ref{BPex}},\Sigma^{
\dot{P}_\alpha}_{\ref{BPex}})$.
\end{enumerate}
By \ref{halbigprop} and \ref{gfbound} we know that each $\dot{\q}_\alpha$ is
a (name for) proper $\baire$--bounding forcing notion and hence
$\forces_{\p_{\omega_2}}$ ``$\dominating=\aleph_1$''. By \ref{trmamea}(2) we
easily conclude that $\forces_{\p_{\omega_2}}$ ``$\non(\M)=\aleph_2$'' and
finally we note that by the last property of $\q^*_{f^P}(K^P_{\ref{BPex}},
\Sigma^P_{\ref{BPex}})$ stated in \ref{BPex}, and by the choice of $\langle
\dot{P}_\alpha: \alpha<\omega_2\rangle$, we have $\forces_{\p_{\omega_2}}$
``$\kappa_{\rm BP}=\aleph_2$''. To finish remember \ref{kbelPP}.
\end{proof}

\begin{conclusion}
\label{BPbelcon}
It is consistent that $\kappa_{\rm BP}=\dominating(R^{\rm sPP})=\aleph_1$ and
$\non(\N)=\con=\aleph_2$. 
\end{conclusion}

\begin{proof} Force over a model of CH with countable support iteration,
$\omega_2$ in length, of forcing notions $\q^*_{{\rm w}\infty}(K_{
\ref{tomek2}},\Sigma_{\ref{tomek2}})$. 

By \ref{getDsPP}(2) we know that the forcing notion $\q^*_{{\rm w}\infty}
(K_{\ref{tomek2}},\Sigma_{\ref{tomek2}})$ has the $(\D,x)$-strong PP--property
for any weakly non-reducible filter $\D$ on $\omega$ and an unbounded
non-decreasing $x\in\baire$. Consequently, if (in $\V$) we take a $p$-point
$\D$ and a PP--ok partial order $\bar{\F}$ then the iteration will have the
$(\D,\bar{\F})$-strong PP--property, so in particular the strong PP--property
(by \ref{ppiniter}; remember that by \ref{preweak}+\ref{presnon} the filter
generated by $\D$ in the intermediate universes is weakly non-reducible). 
Hence, in the resulting model we have $\dominating(R^{\rm   sPP})=\aleph_1$
and thus $\kappa_{\rm BP}=\aleph_1$ (by \ref{kbelPP}). Finally, it follows
from \ref{tomek2conc} that in this model $\non(\N)=\con=\aleph_2$.

Note that one can use the forcing notion $\q^{\tree}_1(K^1_{\ref{tomektree}},
\Sigma^1_{\ref{tomektree}})$ of \ref{tomektree} as well. 
\end{proof}

Let us recall the following notions from \cite{Sh:326}.
\begin{definition}
Let $\F\subseteq\baire$ and $g\in\baire$.
\begin{enumerate}
\item We say that the family $\F$ is {\em $g$--closed} if
\[(\forall f\in\F)(\exists f^+,f^*\in\F)(\forall^\infty n\in\omega)(f(n)^{
g(n)}\leq f^*(n)\ \ \&\ \ \prod_{m<n}(f(m)+1)\leq f^+(n)).\]
\item We say that a proper forcing notion $\p$ has the $(\F,g)$--bounding
property if it has the $(f,g^\varepsilon)$--bounding property for each
$\varepsilon>0$ and $f\in\F$.
\end{enumerate}
\end{definition}
These notions are important when we want to iterate $(f,g)$--bounding forcing
notions: if $\F$ is a $g$--closed family then each countable support iteration
of proper $(\F,g)$--bounding forcing notions is $(\F,g)$--bounding (see
\cite[A2.5]{Sh:326}, compare to \ref{pretFbound}).

\begin{proposition}
\label{newFgclo}
Suppose that $\F\subseteq\baire$ is a $g$--closed family and $\varphi\in
\baire$ is an increasing function. Let $g_\varphi=g\comp\varphi$ and let
$\F_\varphi=\{f\comp\varphi: f\in\F\}$. Then $\F_\varphi$ is
$g_\varphi$--closed. 
\end{proposition}

\begin{proof}
Check.
\end{proof}

\begin{conclusion}
Let $g(n)=n^n$ for $n\in\omega$ and let $\F\subseteq\baire$ be a countable
$g$--closed family. Suppose that $F\in\baire$ is an increasing function which
dominates all elements of $\F$ (i.e.~$(\forall f\in\F)(\forall^\infty n\in
\omega)(f(n)<F(n))$) and let $\bH=\bH^F$, $f=f^F$ (and $\fH$) be as defined in
\ref{bighalex} for $F$. Next, let $f_0\in\F$ and $\langle m_k: k\in\omega
\rangle\subseteq\omega$ and $h\in\baire$ be such that $m_0=0$, $m_{k+1}=m_k+
\fH(k)^{\fH(k)}$, $h$ is non-decreasing and $h(m_{k+1})\leq f_0(\fH(k))$ (for
$k\in\omega$). Assume that $\p_{\omega_2}$ is the countable support iteration
of the forcing notions $\q^*_f(K_{\ref{bighalex}},\Sigma_{\ref{bighalex}})$. 
Then 
\[\forces_{\p_{\omega_2}}\mbox{`` }\dominating(R^*_{\bH})=\aleph_2\ \ \&\ \
\dominating=\dominating(R^*_h)=\dominating(R^\forall_{f_0\comp\fH,g\comp\fH})
=\aleph_1\mbox{ ''.}\]
\end{conclusion}

\begin{proof}
We know that $\F_{\fH}=\{f'\comp\fH:f'\in\F\}$ is $g_{\fH}$--closed,
where $g_{\fH}=g\comp\fH$. By \ref{fgforhal} and the choice of $g$, $F$ we
have that the forcing notion $\q^*_f(K_{\ref{bighalex}},\Sigma_{\ref{bighalex}
})$ is proper, $\baire$--bounding, $(\F_{\fH},g_{\fH})$--bounding and so is
the iteration. Hence, for each $f_1\in\F$, $\forces_{\p_{\omega_2}}$
``$\dominating(R^\forall_{f_1\comp\fH,g\comp\fH})=\dominating=\aleph_1$''. 
Next note that for the function $h$ defined in the assumptions and for
sufficiently large $k$ we have
\[\prod_{n\in [m_k,m_{k+1})}h(n)\leq h(m_{k+1})^{\fH(k)}\leq ((f_0\comp\fH)(
k))^{\fH(k)}\leq (f^*_0\comp\fH)(k),\]
where $f^*_0\in\F$ is such that $(\forall^\infty n\in\omega)(f_0(n)^{g(n)}\leq
f^*_0(n))$. Use \ref{depPPfg}(2) to conclude that $\forces_{\p_{\omega_2}}$
``$\dominating(R^*_h)=\aleph_1$''. Finally, note that if $\dot{W}$ is the name
for the generic real then 
\[\forces_{\q^*_f(K_{\ref{bighalex}},\Sigma_{\ref{bighalex}})} (\forall x\in
\prod_{n\in\omega}\bH(n)\cap\V)(\forall^\infty n\in\omega)(\dot{W}(n)\neq
x(n))\] 
and therefore $\forces_{\p_{\omega_2}}$ ``$\dominating(R^*_{\bH})=\aleph_2$''. 
\end{proof}

\begin{conclusion}
It is consistent that $\non(\M)=\dominating=\aleph_2$ and $\cov(\M)=\unbounded
=\aleph_1=\dominating(R^\exists_{f,g})$ for every non-decreasing unbounded
$g\in\baire$ and any $f\in\baire$ such that $\lim\limits_{n\to\omega}
\frac{f(n)}{g(n)}=\infty$ . 
\end{conclusion}
 
\begin{proof}
Start with a model of CH and iterate $\omega_2$ times with countable
support the Blass-Shelah forcing notion $\q^*_{{\rm s}\infty}(K^*_{
\ref{blsh}},\Sigma^*_{\ref{blsh}})$. By \ref{blshrev} we immediately conclude
that the iteration forces ``$\non(\M)=\dominating=\aleph_2\ \ \&\ \
\unbounded=\aleph_1$''. As each function in $\baire$ appears in an 
intermediate model we may restrict our attention to $f,g\in\baire\cap\V$. By
\ref{omRex} and \ref{presRthm} we conclude that the iteration has the
$R^\exists_{f,g}$--localization property (just build a suitable PP--ok partial
order $\bar{\F}$ and take any $p$-point $\D\in\V$; by \ref{presnon},
\ref{preweak} and \ref{blshrev} we know that $\D$ generates a non-reducible
$p$-filter in the intermediate universes). Hence we get that in the resulting
model $\cov(\M)=\dominating(R^\exists_{f,g})=\aleph_1$.
\end{proof} 

\begin{example}
\label{notsPP}
We construct a finitary 2--big tree--creating pair $(K_{\ref{notsPP}},\Sigma_{
\ref{notsPP}})$ of the $\NMP$--type (see \ref{nmp}(2)) such that the forcing
notion $\q^{\tree}_1(K_{\ref{notsPP}},\Sigma_{\ref{notsPP}})$ does not have
the strong PP--property.
\end{example}

\begin{proof}[Construction]
This example is similar to that of \ref{matet} (what is not surprising
if you notice some kind of duality between ${\mathfrak m}_1$ and
$\dominating(R^{\rm sPP})$). 

Let $\bH(n)=n^n$. Let ${\mathcal A}$ be the family of all pairs $(n,x)$ such
that $x\in [\bH(n)]^{\textstyle n}$. For $\nu\in\prod\limits_{k<m}\bH(k)$,
$m_0<m$ and $A\subseteq{\mathcal A}$ we will write $\nu\prec^*_{m_0} A$ if
\[(\exists (n,x)\in A)(m_0\leq n<m\ \ \&\ \ \nu(n)\in x).\]
Now we define $(K_{\ref{notsPP}},\Sigma_{\ref{notsPP}})$. A tree--like
creature $t\in\TCR_\eta[\bH]$ is taken to be in $K_{\ref{notsPP}}$ if:
\begin{itemize}
\item $\val[t]$ is finite,\quad and
\item $\nor[t]=\log_2(\min\{|A|: A\subseteq{\mathcal A}\ \ \&\ \
(\forall\nu\in\rng(\val[t]))(\nu\prec^*_{\lh(\eta)}A)\})$.
\end{itemize}
The tree composition $\Sigma_{\ref{notsPP}}$ is defined like in \ref{matet}: 
\qquad if $\langle t_\nu: \nu\in\hat{T}\rangle\subseteq K_{\ref{notsPP}}$ is a
system such that $T$ is a well founded quasi tree, $\mrot(t_\nu)=\nu$, and
$\rng(\val[t_\nu])=\suc_T(\nu)$ (for $\nu\in \hat{T}$) then we define
$S^*(t_\nu:\nu\in \hat{T})$ as the unique creature $t^*$ in $K_{\ref{notsPP}}$
with $\rng(\val[t^*])=\max(T)$, $\dom(\val[t^*])=\{\mrot(T)\}$ and $\dis[t^*]=
\langle\dis[t_\nu]:\nu\in\hat{T}\rangle$. Next we put   
\[\Sigma_{\ref{notsPP}}(t_\nu:\nu\in\hat{T})=\{t\in K_{\ref{notsPP}}:
\val[t]\subseteq\val[S^*(t_\nu:\nu\in\hat{T})]\}.\]
It should be clear that $(K_{\ref{notsPP}},\Sigma_{\ref{notsPP}})$ is a
finitary tree--creating pair and the forcing notion $\q^{\tree}_1(K_{
\ref{notsPP}}, \Sigma_{\ref{notsPP}})$ is non-trivial. 

\begin{claim}
\label{cl43}
$(K_{\ref{notsPP}},\Sigma_{\ref{notsPP}})$ is 2--big.
\end{claim}

\noindent{\em Proof of the claim:}\ \ \ Let $t\in K_{\ref{notsPP}}$, $\nor[t]>
0$ and suppose that $\pos(t)=u_0\cup u_1$. Let $s_\ell\in\Sigma_{\ref{notsPP}}
(t)$ be such that $\pos(s_\ell)=u_\ell$ (for $\ell=0,1$). Take $A_\ell
\subseteq{\mathcal A}$ such that
\[|A_\ell|=2^{\nor[s_\ell]}\quad\mbox{ and }\quad (\forall \nu\in\pos(s_\ell))
(\nu\prec^*_{m_0} A_\ell),\]
where $m_0=\lh(\mrot(s_\ell))=\lh(\mrot(t))$. Clearly $(\forall \nu\in\pos(t))
(\nu\prec^*_{m_0} A_0\cup A_1)$ and thus
\[\nor[t]\leq\log_2(|A_0|+|A_1|)\leq 1+\max\{\nor[s_0],\nor[s_1]\}.\]

\begin{claim}
\label{cl44}
$(K_{\ref{notsPP}},\Sigma_{\ref{notsPP}})$ is of the $\NMP$--type (see
\ref{nmp}(2)). 
\end{claim}

\noindent{\em Proof of the claim:}\ \ \ Suppose that $\langle t_\eta:\eta\in
T\rangle\in\q^{\tree}_\emptyset(K_{\ref{notsPP}},\Sigma_{\ref{notsPP}})$ is
such that $(\forall \eta\in T)(\nor[t_\eta]>1)$ and $F_0,F_1,F_2,\ldots$ are
fronts of $T$ such that
\[(\forall\nu\in F_{i+1})(\exists\nu'\in F_i)(\nu'\vartriangleleft\nu).\]
Clearly these fronts are finite (as $(K_{\ref{notsPP}},\Sigma_{\ref{notsPP}})$
is finitary). Further suppose that $g:\bigcup\limits_{i\in\omega} F_i
\longrightarrow\bigcup\limits_{i\in\omega} F_{i+1}$ is such that $\nu
\vartriangleleft g(\nu)\in F_{i+1}$ for $\nu\in F_i$. 

\noindent Let $r=\min\{2^{\nor[t_\eta]}:\eta\in T\}$. Choose increasing
sequences $\langle n_k: k\leq r^r+r\rangle$, $\langle i_k,j_k: k<r^r+r\rangle$
such that for each $k<r^r+r$:
\begin{enumerate}
\item[(i)]\ \ \ $i_k<j_k<i_{k+1}$ and if $k\in [r, r+r^r)$ then $j_k=i_k+1$,
\item[(ii)]\ \  if $k<r$, $\nu\in F_{i_k}$ then $|\{\rho\in F_{j_k}:\nu
\vartriangleleft\rho\}|\geq r$,
\item[(iii)]\   $n_k<\min\{\lh(\nu):\nu\in F_{i_k}\}<\max\{\lh(\nu):\nu\in
F_{j_k}\}<n_{k+1}$. 
\end{enumerate}
Choose a mapping $\pi: F_{j_{r-1}}\longrightarrow r^{\textstyle r}$ such that
\begin{enumerate}
\item[$(*)_0$] if $\nu\vartriangleleft\eta_0\in F_{j_{r-1}}$, $\nu
\vartriangleleft\eta_1\in F_{j_{r-1}}$, $\nu\in F_{j_k}$, $k<r$\\
then $\pi(\eta_0)(k)=\pi(\eta_1)(k)$,
\item[$(*)_1$] for each $\nu\in F_{i_k}$, $k<r$ and $\ell<r$ there is $\eta\in
F_{j_{r-1}}$ such that $\nu\vartriangleleft\eta$ and $\pi(\eta)(k)=\ell$. 
\end{enumerate}
(It is easy to define such a mapping if you remember clause (ii) above.)
Let $\pi^*:r^{\textstyle r}\longrightarrow r^r$ be the isomorphism of
$r^{\textstyle r}$ equipped with the lexicographical order and $r^r$ with the
natural order of integers. Take a tree--creature $s\in\Sigma_{\ref{notsPP}}(
t_\eta: (\exists\nu\in F_{j_{r^r+r-1}})(\eta\vartriangleleft\nu))$ such that 
\[\begin{array}{ll}
\rng(\val[s])=\{\nu\in F_{j_{r^r+r-1}}:&\mbox{if }\nu\rest m_0\in F_{j_{r-1}},
\ m_0\in\omega\mbox{ and }k=\pi^*(\pi(\nu\rest m_0))\\
\ &\mbox{and }\nu\rest m_1\in F_{i_{r+k}},\ m_1\in\omega\ \ \mbox{ then
  }g(\nu\rest m_1)\vartriangleleft\nu\}. 
  \end{array}\]
(Note that we may find a suitable $s$ by the definition of $\Sigma_{
\ref{notsPP}}$.) By the choice, this $s$ satisfies the demand $(\beta)^\tree$
of \ref{nmp}(2). But why does it have large enough norm? Suppose that
$A\subseteq{\mathcal A}$ is such that $|A|<r$. Let $k_0<r$ be such that 
\[(n,x)\in A\quad\Rightarrow\quad n\notin [n_{k_0},n_{k_0+1}).\]
Since $|A|<2^{\nor[t_\eta]}$ for each $\eta\in T$ we may inductively build a
sequence $\nu_0\in F_{i_{k_0}}$ such that 
\[(\forall (n,x)\in A)(\lh(\mrot(T))\leq n<\lh(\nu_0)\ \Rightarrow\ \nu_0(n)
\notin x).\]
Let $\sigma_0:k_0\longrightarrow r$ be such that $\sigma_0(k)$ is the value of
$\pi(\eta)(k)$ for each $\eta\in F_{j_{r-1}}$, $\nu_0\vartriangleleft\eta$. 
Take $\ell<r$ such that 
\begin{quotation}
\noindent if $\sigma\in r^{\textstyle r}$, $\sigma_0\conc \langle\ell\rangle
\trianglelefteq\sigma$

\noindent then there is no $(n,x)\in A$ with $n_{\pi^*(\sigma)+r}\leq n<
n_{\pi^*(\sigma)+r+1}$ 
\end{quotation}
(remember the choice of $\pi^*$ and that $|A|<r$). Now take $\nu_1\in
F_{j_{k_0}}$ such that $\nu_0\vartriangleleft\nu_1$ and 
\[(\forall\eta\in F_{j_{r-1}})(\nu_1\trianglelefteq\eta\ \Rightarrow\
\pi(\eta)(k_0)=\ell)\]
(possible by $(*)_0$ + $(*)_1$). By the choice of $k_0$ we know that
\[(\forall (n,x)\in A)(\lh(\mrot(T))\leq n<\lh(\nu_1)\ \Rightarrow\ \nu_1(n)
\notin x)\]
(look at (iii)). Next continue like at the beginning to get $\eta\in
F_{j_{r-1}}$ such that $\nu_1\trianglelefteq\eta$ and $\neg (\eta\prec^*_{\lh(
\mrot(T))} A)$. We are sure that $\sigma_0\conc\langle\ell\rangle
\trianglelefteq\pi(\eta)$ and therefore there is no $(n,x)\in A$ with
$n_{\pi^*(\pi(\eta))+r}\leq n<n_{\pi^*(\pi(\eta))+r+1}$. Consequently we may
continue the procedure applied to build $\eta$ and we construct $\eta^*\in
F_{j_{r^r+r-1}}$ such that

$\eta\vartriangleleft\eta^*$, if $\eta^*\rest m\in F_{i_{r+\pi^*(\pi
(\eta))}}$ then $g(\eta^*\rest m)\vartriangleleft\eta^*$, and $\neg (\eta
\prec^*_{\lh(\mrot(T))} A)$.

\noindent Since, by its construction, the sequence $\eta^*$ is in
$\rng(\val[s])$, it exemplifies that $A$ cannot witness the minimum in the
definition of $\nor[s]$. Consequently, $\nor[s]\geq\log_2(r)$ and thus the
tree--creature $s$ satisfies the demand $(\alpha)^\tree$ of \ref{nmp}(2).

\begin{claim}
The forcing notion $\q^{\tree}_1(K_{\ref{notsPP}},\Sigma_{\ref{notsPP}})$ does
not have the strong PP--property.
\end{claim}

\noindent{\em Proof of the claim:}\ \ \ We will show that the generic real
$\dot{W}$ shows that the strong PP--property fails for $\q^{\tree}_1(
K_{\ref{notsPP}},\Sigma_{\ref{notsPP}})$. So suppose that $\langle w_i: i\in
B\rangle$ is such that $B\in\iso$ and $(\forall i\in B)(w_i\in [\omega]^{
\textstyle i})$ and let $p\in\q^{\tree}_1(K_{\ref{notsPP}},\Sigma_{
\ref{notsPP}})$. We may assume that $\nor[t^p_\eta]>1$ for each $\eta\in
T^p$. Take $i_0\in B\setminus\lh(\mrot(p))$ and build inductively a condition
$q\geq p$ such that for each $\eta\in T^q$ and $\nu\in \pos(t^q_\eta)$
\[t^q_\eta\in\Sigma_{\ref{notsPP}}(t^p_\eta)\quad\mbox{and}\quad
\nor[t^q_\eta]\geq\nor[t^p_\eta]-1\quad\mbox{and}\quad (i_0<\lh(\nu)\
\Rightarrow\ \nu(i_0)\notin w_{i_0})\] 
(remember the definition of the norm of elements of $K_{\ref{notsPP}}$). Now
clearly $q\forces$ `` $\dot{W}(i_0)\notin w_{i_0}$ '', finishing the proof of
the claim and the construction. 
\end{proof}

\begin{conclusion}
It is consistent that $\cof(\M)<\dominating(R^{\rm sPP})$.
\end{conclusion}

\begin{proof}
Start with a model of CH and force with countable support iteration of
length $\omega_2$ of forcing notions $\q^{\tree}_1(K_{\ref{notsPP}},
\Sigma_{\ref{notsPP}})$. We know that $\q^{\tree}_1(K_{\ref{notsPP}},
\Sigma_{\ref{notsPP}})$ is proper, $\baire$--bounding and Cohen--preserving
(by \ref{treenonmea} + \ref{treebound}). Consequently the iteration is of the 
same type (see \cite[6.3.21, 6.3.22]{BaJu95}) and, by standard arguments, in
the final model we have $\non(\M)=\dominating=\aleph_1$. But this implies that
$\cof(\M)=\aleph_1$ too (see \cite[2.2.11]{BaJu95}). Finally, as
$\q^{\tree}_1(K_{\ref{notsPP}},\Sigma_{\ref{notsPP}})$ does not have the
strong PP--property we easily conclude that the iteration forces that
$\dominating(R^{\rm sPP})=\aleph_2$.
\end{proof}

\backmatter
\chapter*{List of definitions}
\begin{enumerate}

\item[\ref{weakcreat}]\ \ \  weak creatures, $\WCR[\bH]$;
\item[\ref{morecreat}]\ \ \  finitary $\bH$, finitary $K$;
\item[\ref{subcom}]\ \ \     sub-composition operation, weak creating pair,
the relation $\sim_\Sigma$; 
\item[\ref{basis}]\ \ \      basis $\basis(t)$, possibilities $\pos(w,
{\mathcal S})$;
\item[\ref{maindef}]\ \ \    forcing notion $\q_{{\mathcal
C}(\nor)}(K,\Sigma)$ (for a weak creating pair $(K,\Sigma)$ and a norm
condition ${\mathcal C}(\nor)$); 
\item[\ref{conditions}]\ \   $m_{\dn}(t)$, norm conditions and corresponding
forcing notions $\q_{{\rm s}\infty}(K,\Sigma)$, $\q_{\infty}(K,\Sigma)$,
$\q_{{\rm w}\infty}(K,\Sigma)$, $\q_f(K,\Sigma)$,
$\q_\emptyset(K,\Sigma)$;
\item[\ref{fast}]\ \         fast function, $\bH$-fast function
$f:\omega\times\omega\longrightarrow\omega$;
\item[\ref{thereal}]\ \      name for the generic real $\dot{W}$;
\item[\ref{creatures}]\ \ \  $m^t_{\dn}$, $m^t_{\up}$, creatures, $\CR[\bH]$;
\item[\ref{compos}]\ \ \     composition operation on $K$, creating pairs
$(K,\Sigma)$;
\item[\ref{candidates}]\ \ \   finite candidates $\FC(K,\Sigma)$, pure finite
candidates $\PFC(K,\Sigma)$, pure candidates $\PC(K,\Sigma)$,
$\C(\nor)$--normed pure candidates $\PC_{\C(\nor)}(K,\Sigma)$ and partial
orders on them;
\item[\ref{niceandsmo}]\ \ \  creating pairs which are: nice, smooth,
forgetful, full; 
\item[\ref{forcingstar}]\ \ \  forcing notions $\q^*_{{\mathcal
C}(\nor)}(K,\Sigma)$ for creating pairs $(K,\Sigma)$;
\item[\ref{essapprox}]\ \ \   when a condition $p$ essentially decides a name
$\dot{\tau}$, approximates $\dot{\tau}$;
\item[\ref{orders}]\ \        partial orders $\leq_{\apr}$, $\leq^{{\rm
s}\infty}_n$, $\leq^{\infty}_n$, $\leq^{{\rm w}\infty}_n$, $\leq^f_n$;
\item[\ref{quasitree}]\ \ \   quasi trees, well founded quasi trees, downward
closure $\dcl(T)$, successors $\suc_T(\eta)$ of $\eta$ in $T$, $T^{[\eta]}$,
$\spliting(T)$, $\max(T)$, $\hat{T}$, $\lim(T)$, fronts of a quasi tree $T$;
\item[\ref{treecreature}]\ \ \  tree--creatures, $\TCR[\bH]$,
tree--composition, bounded tree--composition;
\item[\ref{treeforcing}]\ \ \  forcing notions $\q^{\tree}_e(K,\Sigma)$ for
$e<5$, $\q^{\tree}_\emptyset(K,\Sigma)$, condition $p^{[\eta]}$ for
$p\in\q^{\tree}_e(K,\Sigma)$, $\eta\in T^p$;
\item[\ref{thick}]\ \ \       $e$-thick antichains in $T^p$ for
$p\in\q^{\tree}_e(K,\Sigma)$;
\item[\ref{AxA}]\ \            partial orders $\leq^e_n$ (for $e<3$);
\item[\ref{local}]\ \ \       local weak creating pairs;

\item[\ref{omitoryetc}]\ \ \  creature $t\Rsh [m_0,m_1)$, for a creating pair
$(K,\Sigma)$ we say when it is omittory, growing;
\item[\ref{simpglui}]\ \ \    creating pairs which are: gluing, simple;
\item[\ref{singleton}]\ \     creating pairs which capture singletons;
\item[\ref{bigetc}]\ \ \      big creating pairs;
\item[\ref{omitbig}]\ \ \     omittory--big creating pairs;
\item[\ref{halving}]\ \ \     Halving Property and weak Halving Property;
\item[\ref{kbig}]\ \ \        big tree--creating pairs;
\item[\ref{tomit}]\ \ \       t--omittory tree creating pairs;
\item[\ref{prenorms}]\ \ \    pre--norm on ${\mathcal P}(A)$, nice pre-norm;
\item[\ref{exprno}]\ \ \      pre--norms $\dpt^i$, $\dpt^i_n$ (for $i<3$,
$n\in\omega$); 

\item[\ref{meapres}]\ \ \     Cohen--preserving proper forcing notions;
\item[\ref{nmp}]\ \ \         creating pairs of the $\NMP$--type, tree creating
pairs of the $\NMP^{\tree}$--type;
\item[\ref{trmea}]\ \ \       trivially meagering weak creating pairs;
\item[\ref{tnnp}]\ \ \        weak creating pairs of the $\NNP$--type;
\item[\ref{gluing}]\ \ \      gluing and weakly gluing tree creating pairs;
\item[\ref{strfin}]\ \ \      strongly finitary creating pairs;
\item[\ref{srsp}]\ \ \        when a weak creating pair $(K,\Sigma)$ strongly
refuses Sacks property;

\item[\ref{meagering}]\ \ \   creating pairs which are: meagering, anti--big;
\item[\ref{sum}]\ \ \         $\Sigma^{\bsum}$;
\item[\ref{dsum}]\ \ \        $(d,u)$--sum $\Sigma^{\bsum}_{d,u}$; 
\item[\ref{saturated}]\ \ \   when a creating pair is saturated with respect
to a family of pre--norms;
\item[\ref{treesum}]\ \ \     $\Sigma^{\tsum}$;
\item[\ref{ab}]\ \ \          decision functions, creating pairs of the
$\AB$-type, condensed creating pairs;
\item[\ref{abplus}]\ \ \      creating pairs of the $\AB^+$--type;

\item[\ref{fbig}]\ \ \        essentially $f$-big weak creating pairs;
\item[\ref{reducible}]\ \ \   reducible weak creating pairs;
\item[\ref{hlimited}]\ \ \    $h$--limited weak creating pairs;
\item[\ref{FHfast}]\ \        $(\bH,F)$--fast function;
\item[\ref{tFboundetc}]\ \ \  $U_h(\bar{t})$, $V^n_h(\bar{t})$, $(\bar{t},h_1,
h_2)$--bounding forcing notions;
\item[\ref{monotonic}]\ \ \   creating pairs which are monotonic, strictly
monotonic, spread;
\item[\ref{additive}]\ \ \    m--additivity $\add_m(t)$ of a weak creature $t$,
$(g,h)$-additive weak creating pairs; 
\item[\ref{tgoodetc}]\ \ \    $\bar{t}$--good families of functions,
$(\bar{t},\F)$--bounding forcing notions; 
\item[\ref{quasi}]\ \ \       $\bar{t}$--systems, regular $\bar{t}$--systems,
$\p^*_{\C(\nor)}(\bar{t},(K,\Sigma))$ and the partial order $\preceq$ on it,
quasi-$W$-generic $\Gamma$;
\item[\ref{Gampres}]\ \ \     $(\Gamma,W)$--genericity preserving forcing
notions;
\item[\ref{senetc}]\ \ \      Cohen sensitive $\bar{t}$--systems, directed
$\bar{t}$--systems;
\item[\ref{coherent}]\ \      $(\bar{t}_0,h_1,h_2)$--coherent
$\bar{t}$--systems, $(\bar{t}_0,\bar{\F})$--coherent sequences of
$\bar{t}$--systems;

\item[\ref{genult}]\ \ \      creating pair which generates an ultrafilter;
\item[\ref{Gamfilter}]\ \ \   when $\Gamma$ generates a filter (ultrafilter),
$\D(\Gamma)$;
\item[\ref{intcrea}]\ \ \     interesting creating pair;
\item[\ref{ramppoint}]\ \ \   Ramsey filter, p--point, q--point, weak
q--point;
\item[\ref{almramsey}]\ \ \   interesting ultrafilters, games $G^{sR}(\D)$,
$G^{aR}(\D)$, semi--Ramsey ultrafilters, almost Ramsey ultrafilters;
\item[\ref{up}]\ \ \          tree creating pairs of the $\UP(\D)^{\tree}$,
$\sUP(\D)^{\tree}$ --types;
\item[\ref{rich}]\ \ \        rich tree creating pairs;
\item[\ref{ramnum}]\ \ \      $R_n(k,m)$;
\item[\ref{simomi}]\ \ \      simple except omitting creating pairs,
omittory--compatible $\bar{t}$--systems; 

\item[\ref{BPnumber}]\ \ \    $\kappa_{BP}$;
\item[\ref{sPPrel}]\ \ \      $R^{\rm sPP}$;
\item[\ref{DsPP}]\ \ \        non-reducible filters, PP--ok partial orders,
forcing notions with $(\D,x)$-strong PP--property, $(\D,\bar{\F})$-strong
PP--property;
\item[\ref{fgPP}]\ \ \        $R^\exists_{f,g}$, $R^\forall_{f,g}$, $R^*_f$,
$R^{**}_{f,g}$;
\item[\ref{presRdef}]\ \ \    $(\D,\bar{\F})$-$R^\exists_f$--localization
property.
\end{enumerate}


\end{document}